\documentclass{article}

\usepackage[latin1]{inputenc}
\usepackage[T1]{fontenc}
\usepackage[a4paper]{geometry}
\usepackage{supertabular}
\usepackage{amsmath}
\usepackage{amssymb}
\usepackage{theorem}
\usepackage{amscd}
\usepackage{makeidx}

\newtheorem{thm}{Theorem}[subsection]
\newtheorem{defn}{Definition}[subsection]
\newtheorem{lemma}{Lemma}[subsection]
\newtheorem{prop}{Proposition}[subsection]
\newtheorem{exam}{Example}[subsection]
\newtheorem{cory}{Corollary}[subsection]
\newtheorem{rmk}{Remark}[subsection]

\makeatletter

\@addtoreset{equation}{section}
\makeatother

\makeindex

\newcommand{\oct}{\mathbb{O}}
\newcommand{\R}{\mathbb{R}}
\newcommand{\C}{\mathbb{C}}
\newcommand{\h}{\mathbb{H}}
\newcommand{\g}{\mathfrak{g}}

\newcommand{\lm}{\lambda}
\newcommand{\Lm}{\Lambda}
\newcommand{\iif}{if and only if }
\newcommand{\Gt}{G_{\tau}}
\newcommand{\Gtc}{G_{\tau}^{\C}}
\newcommand{\gtc}{\g_{\tau}^{\C}}
\newcommand{\gtau}{\g_{\tau}}

\newcommand{\im}{\text{\normalfont Im\,}}
\newcommand{\eps}{\varepsilon}

\newcommand{\dl}[2]{\frac{\partial{#1}}{\partial{#2}}}
\newcommand{\mk}{\mathfrak{m}}
\newcommand{\hk}{\mathfrak{h}}
\newcommand{\lie}{\mathrm{Lie}\,}
\newcommand{\Is}{\mathrm{Is}}
\newcommand{\Ad}{\mathrm{Ad}}
\newcommand{\ad}{\mathrm{ad}}
\newcommand{\so}{\mathfrak{so}}
\newcommand{\mrm}[1]{\mathrm{#1}}

\newcommand{\tm}{\tau_{|\mathfrak{m}}}
\newcommand{\Int}{\mathrm{Int}}

\newcommand{\tl}[1]{\tilde{#1}}

\newcommand{\taum}{\tau_{\mathfrak{m}}}

\newcommand{\tmh}{\tau_{\mid\mathfrak{h}}}
\newcommand{\und}[1]{\underline{#1}}
\newcommand{\ovr}[1]{\overline{#1}}
\newcommand{\mak}[1]{\mathfrak{#1}}
\newcommand{\Adm}{\mathrm{Ad}_{\mk}}
\newcommand{\adm}{\mathrm{ad}_{\mk}}
\newcommand{\Id}{\mathrm{Id}}
\newcommand{\Z}{\mathbb{Z}}

\newcommand{\mal}[1]{\mathcal{#1}}
\newcommand{\syst}{\mathrm{Syst}}
\newcommand{\systm}{\mathrm{Syst}(m)}
\newcommand{\systmt}{\mathrm{Syst}(m,\tau)}
\newcommand{\zdk}{\mathcal{Z}_{2k}}
\newcommand{\rdn}{\R^{2n}}
\newcommand{\U}{\mathbb{U}}
\newcommand{\ul}{\mathfrak{u}}
\newcommand{\com}{\mathrm{Com}}
\newcommand{\End}{\mathrm{End}}
\newcommand{\Um}{\mathfrak{U}}
\newcommand{\ver}{\mathcal{V}}
\newcommand{\hor}{\mathcal{H}}
\newcommand{\pk}{\mathfrak{p}}
\newcommand{\kk}{\mathfrak{k}}
\newcommand{\nk}{\mathfrak{n}}
\newcommand{\mU}{\mathrm{U}}
\newcommand{\mB}{\mathrm{B}}
\newcommand{\mR}{\mathrm{R}}
\newcommand{\maB}{\mathcal{B}}
\newcommand{\maU}{\mathcal{U}}
\newcommand{\mJ}{\mathcal{J}}
\newcommand{\mZ}{\mathcal{Z}}
\newcommand{\rabla}[2]{\nabla_{#2}^{[#1]}}
\newcommand{\mijo}{\mathfrak{I}_{J_0}}
\newcommand{\mai}{\mathfrak{I}}
\newcommand{\onabla}[3]{\overset{#1}{\nabla_{#3}^{#2}}}
\newcommand{\nablamet}[2]{\overset{\mathrm{met}}{\nabla_{#1}}{}^{#2}}
\newcommand{\mv}{\mathrm{v}}
\newcommand{\mH}{\mathrm{H}}
\newcommand{\maH}{\mathsf H}
\newcommand{\undj}{\underline{\mathrm{J}}}
\newcommand{\wrt}{w.r.t.\ }
\newcommand{\circact}{\circlearrowleft \negthickspace\negthickspace\negthickspace\cdot \,}
\newcommand{\circacti}{\circlearrowleft\,\negthickspace\negthickspace\!\cdot \,}
\newcommand{\proof}{\textbf{Proof. }}
\newcommand{\comprf}{This completes the proof. }
\newcommand{\hsq}{\hfill$\square$}
\newcommand{\con}{\mathrm{Con}}

\author{Idrisse Khemar}
\title{Elliptic Integrable Systems: a Comprehensive Geometric Interpretation}
\date{}

\begin{document}
\maketitle

\null\hfill  \textbf{Abstract.}  \hfill\null\\
{\small
In this paper, we study all the elliptic integrable systems, in the sense of C.L. Terng \cite{tern}. That is to say the family  of all the $m$-th elliptic integrable systems
associated to a $k'$-symmetric space $N=G/G_0$. Here  $m\in \mathbb N$ and $k'\in \mathbb N^*$ are integers. For example, it is known that the first elliptic integrable system associated to a symmetric space (resp. to a Lie group) is the equation for harmonic maps into this symmetric space (resp. this Lie group). Indeed it is well known that this harmonic maps equation can be  written as a zero curvature equation:\index{harmonic map}
$
d\alpha_{\lm} + \frac{1}{2}[\alpha_{\lm}\wedge\alpha_{\lm}]=0, \quad \forall\lm\in \C^*,
$
where $\alpha_\lm = \lm^{-1}\alpha_1' + \alpha_0 + \lm \alpha_1''$ is a 1-form on a Riemann surface $L$ taking values in the Lie algebra $\g$. This 1-form $\alpha_\lm$ is obtained as follows. Let $f\colon L\to N=G/G_0$ be a map from the Riemann surface $L$ into the symmetric space $G/G_0$. Then let $F\colon L\to G$ be a  lift of $f$, and consider $\alpha=F^{-1}.dF$ its Maurer-Cartan form. Then decompose $\alpha$ according to the symmetric decomposition
$ \g=\g_0\oplus\g_1$ of $\g$ : $\alpha=\alpha_0 + \alpha_1$. Finally, we define $\alpha_\lm:= \lm^{-1}\alpha_1' + \alpha_0 + \lm \alpha_1''$, $\forall\lm\in \C^*$, where $\alpha_1',\alpha_1''$ are the resp. the $(1,0)$ and $(0,1)$ parts of $\alpha_1$. Then the zero curvature equation for this $\alpha_\lm$, for all $\lm\in \C^*$, is equivalent to the harmonic maps equation for $f\colon L\to N=G/G_0$, and is by definition  the first elliptic integrable system associated to the symmetric space $G/G_0$. Thus the methods of integrable system theory  apply to give generalised Weierstrass representations, algebro-geometric solutions, spectral deformations and so on. In particular, we can apply the DPW method \cite{DPW} to obtain a generalised Weierstrass representation. More precisely, we have a Maurer-Cartan equation in some loop Lie algebra $\Lambda \g_\tau=\{\xi\colon S^1\to \g |\xi(-\lm)=\tau(\xi(\lm))\}$, then we can integrate it in the corresponding loop group and finally apply some factorizations theorems in loop groups to obtain a generalised Weierstrass representation: this is the DPW method. Moreover, these methods of integrable system theory hold for all the systems written in the forms of a zero curvature equation for some $\alpha_\lm=\lm^{-m}\hat{\alpha}_{-m} + \cdots +\hat{\alpha}_0 + \cdots + \lm^m\hat{\alpha}_m$. Namely, these method apply to construct the solutions of all the $m$-th elliptic integrable systems. So it is natural to ask what is the geometric interpretation of these systems. Do they correspond to some generalisations of harmonic maps? This is the problem that we solve in this paper: to describe the geometry behind this family of integrable systems whose we know how to construct (at least locally) all the solutions. The introduction below gives an overview of all the main  results, as well as some related subjects and works, and some additional motivations.
}

\tableofcontents
%
%
%
%%%%%%%%%%%%%%%%%%%%%%%%%%%%%%%%%%%%%%%%%%%%%%%%%%%%%%%%%%%%%%%%%%%%%%%%%

%     Introduction

%%%%%%%%%%%%%%%%%%%%%%%%%%%%%%%%%%%%%%%%%%%%%%%%%%%%%%%%%%%%%%%%%%%%%%%
%

%
%\renewcommand\thesection{\Roman{\section}}
%
\section*{Introduction}
\addcontentsline{toc}{section}{\protect\numberline{$\bullet$}Introduction}
\index{m elliptic integrable system@$m$-th elliptic integrable system}\index{elliptic integrable system|see{$m$-th elliptic integrable system}}\index{k symmetric space@$k'$-symmetric space|(}
In this paper, we give a geometric interpretation of all the $m$-th elliptic integrable systems
associated to a $k'$-symmetric space $N=G/G_0$ (in the sense of C.L. Terng \cite{tern}).\medskip\\ 
\index{eigenspace decomposition of the automorphism $\tau$}
Let $\g$ be a real Lie algebra and $\tau\colon \g \to \g$ be an automorphism of finite order $k'$. This automorphism admits an eigenspace decomposition $\g^\C=\oplus_{j\in\Z_{k'}} \g_j^\C$, where $\g_j^\C$ is the eigenspace  of  $\tau $ with respect to the eigenvalue $\omega_{k'}^{j}$. We denote by $\omega_{k'}$ a $k'$-th primitive root of  unity. Moreover  the automorphism $\tau$ defines a $k'$-symmetric space $N=G/G_0$. Furthermore let $L$ be a Riemann surface.\\
The $m$-th elliptic integrable system associated to $(\g, \tau)$ can be written as a zero curvature equation 
$$
d\alpha_{\lm} + \frac{1}{2}[\alpha_{\lm}\wedge\alpha_{\lm}]=0, \quad \forall\lm\in \C^*,
$$
where $\alpha_\lm= \sum_{j=0}^m \lm^{-j}u_j + \lm^j \bar u_j =\sum_{j=-m}^m \lm^j\hat{\alpha}_j$ is a 1-form on the Riemann surface $L$ taking values in the Lie algebra $\g$. The "coefficient" $u_j$ is a $(1,0)$-type 1-form on $L$ with values in the eigenspace $\g_{-j}^\C$.\\
\index{order m@order $m$ of the system}
Moreover, we call the integer $m$ the \emph{order} of the system. Let us make precise that the order $m$ has nothing to do with the order of any PDE, but this is only the maximal power on $\lm$ in the (finite) Fourier decomposition of $\alpha_\lm$ \wrt $\lm$.
\\
First, we remark that any solution of the  system of order $m$ is a solution of the system of order $m'$, if $m\leq m'$ (and the automorphism $\tau$ is fixed). In other words, the system of order $m$ is a reduction of the system of order $m'$, if $m\leq m'$.\\
Moreover, it turns out that we have to introduce the integer $m_{k'}$ defined by 
$$
m_{k'}=\displaystyle \left[\frac{k'+1}{2}\right]=\begin{cases} k  \text{ if } k'=2k\\ k +1 \text{ if } k'=2k+1 \end{cases} \text{ if }k'>1, \text{ and }m_{1}=0.
$$
Then the general problem splits into three cases : the primitive case ($m<m_{k'}$), the
determined case\index{determined@determined, case or system} ($m_{k'}\leq m \leq k'-1$) and the underdetermined case ($m \geq k'$).\index{underdetermined@underdetermined, case or system}\index{primitive@primitive, case or system|(}
\subsection{The primitive systems}\label{primitive_intro}
The primitive systems have an interpretation in terms of $F$-holomorphic maps,\index{f holomorphic@$f$-holomorphic} with respect to an $f$-struture \index{f structure@$f$-structure} $F$ on the target space $N=G/G_0$ (i.e. an endomorphism $F$ satisfying $F^3 + F= 0$). More precisely:\medskip\\
$\bullet$\index{canonical!$f$-structure}\index{horizontal subbundle} \index{vertical subbundle}\index{even case|(} \emph{In the even case} ($k'=2k$), we have a fibration $\pi\colon G/G_0\to G/H$ over a $k$-symmetric space $M=G/H$ (defined by the square $\tau^2$ of the automorphism $\tau$ of order $k'$ defining $N=G/G_0$). We also have  a $G$-invariant splitting $TN=\hor \oplus \ver$ corresponding to this fibration\footnote{i.e. $\hor$ is a connection  on this fibration.}, where $\ver=\ker d\pi$.  Moreover $N$ is naturally endowed with an $f$-structure $F$ which defines a complex structure on the horizontal subbundle $\hor$ and vanishes on the vertical subbundle $\ver$. Furthermore, the eigenspace decomposition of the order $k'$ automorphism $\tau$ gives us some $G$-invariant decomposition $\hor=\oplus_{j=1}^{k-1}[\mk_j]$, where $\mk_j\subset\g$ is defined by $\mk_j^\C=\g_{-j}^\C\oplus\g_j^\C$, and $[\mk_j]\subset TN$ is the corresponding $G$-invariant subbundle. This allows to define, by multiplying  $F$ on the left by the projections on the subbundles $\hor^m=\oplus_{j=1}^m[\mk_j]$, a family of $f$-structures $F^{[m]}, 1\leq m\leq k-1$. Then the primitive system of order $m$ ($m<m_{k'}=k$) associated to $G/G_0$ is exactly the equation for $F^{[m]}$-holomorphic maps. Therefore the solutions of  the primitive systems are exactly the $F$-holomorphic maps.\medskip\\
$\bullet$\index{canonical!almost complex structure}\index{J holomorphic@$J$-holomorphic}\index{odd case|(}
\emph{In the odd case} ($k'=2k+1$), $N=G/G_0$ is naturally endowed with an almost complex structure $\undj$. Then the solutions of the primitive systems are exactly the $\undj$-holomorphic curves. Moreover, in the same way as for the even case, the eigenspace decomposition of $\tau$ provides a $G$-invariant decomposition $TN=\oplus_{j=1}^{k}[\mk_j]$, which allows to define a family of $f$-structures $F^{[m]}, 1\leq m\leq k$, with $F^{[k]}=\undj $. Then the primitive system of order $m$ ($m<m_{k'}=k+1$) associated to $G/G_0$ is exactly the equation for $F^{[m]}$-holomorphic maps. In other words, the solutions of the primitive system of order $m$ are exactly the integral holomorphic curves of the  complex Pfaffian system $\oplus_{j=1}^m[\mk_j]\subset TN$ in the almost complex manifold $(N,\undj)$.
\index{primitive@primitive, case or system|)}
\subsection{The determined case}\index{determined@determined, case or system} 
We  call "the minimal determined system"\index{determined minimal@determined, minimal|(} the determined system  of minimal order $m_{k'}$, and "the maximal determined system"\index{determined maximal@determined, maximal|(} the determined system of maximal order $k'-1$.\\
 Any solution of a determined system is solution of the corresponding maximal determined system. More precisely, a map $f\colon L\to G/G_0$ is solution of a determined system (associated to $G/G_0$) \iif it is solution of the maximal determined system (associated to $G/G_0$)  and satisfies an additional holomorphicity condition.\index{holomorphicity condition}  When this holomorphicity condition is maximal, then we obtain the minimal determined system.%\medskip\\
\subsubsection{The minimal determined  system}
\index{determined minimal@determined, minimal} \index{vertically harmonic|(}
The minimal determined  system has an interpretation in terms of horizontally holomorphic and vertically harmonic maps $f\colon L\to N=G/G_0$. It also has  an equivalent interpretation in terms of vertically harmonic twistor lifts in some twistor space. Let us make precise this point.
\paragraph{In the even case.}\index{canonical!$f$-structure}\index{horizontal subbundle}
\index{vertical subbundle}
As we have seen in the subsection \ref{primitive_intro} below, the homogeneous space $N=G/G_0$ admits a $G$-invariant splitting $TN=\hor \oplus \ver$ corresponding to the fibration $\pi\colon N\to M$ and  $N$ is naturally endowed with an $f$-structure $F$ which defines a complex structure on the horizontal subbundle $\hor$ and vanishes on the vertical subbundle $\ver$. Then we say that a map $f\colon L\to N$ is horizontally holomorphic if \index{horizontally holomorphic}
$$
(df\circ j_L)^{\mal H}=F\circ df.
$$
Then we prove that \emph{the even minimal determined system $(\syst(k,\tau)$) means that the geometric map $f$ is horizontally holomorphic and vertically harmonic}, i.e.\index{vertical tension field}
$$
\tau^v(f):=\mrm{Tr}_g(\nabla^v d^v f)=0
$$
(for any hermitian metric $g$ on the Riemann surface $L$). Here $\nabla^v$ is the vertical component of the Levi-Civita connection $\nabla$ (of some $G$-invariant metric on $N$). Moreover,  vertically harmonic maps are exactly the critical points (\wrt vertical variations) of the vertical energy functional:
$$
E^v(u)=\dfrac{1}{2}\int_L |d^v u|^2 d\mrm{vol}_g .
$$
%\medskip\\
%
\index{twistor}
We prove also that this system  also has  an equivalent interpretation in terms of vertically harmonic twistor lifts in some twistor space $\mal Z_{2k,j}(M,J_2)$ which is a subbundle of $\zdk(M)$, where 
$$
\mal Z_{k'}(M)=\{J\in SO(TM)|\, J^{k'}=\Id,J^p\neq\Id\text{ if } p<k', \ker(J\pm\Id)=\{0\}\}
$$
is the bundle of isometric endomorphisms of $TM$ with finite order $k'$ and with no eigenvalues $=\pm1$. More precisely denoting by $J_2$  the section  of $SO(TM)$, of  order $k$, defined by $\taum^2$, then we define $\mal Z_{2k,j}(M,J_2)=\{J\in \zdk(M)|J^2=J_2\}$. Then we prove that $N=G/G_0$ can be embedded into the twistor space $\mal Z_{2k,j}(M,J_2)$ via a natural morphism of bundle over $M=G/H$. We prove that $f\colon L\to N$ is solution of the system \iif the corresponding map $J^f\colon L\to \mal Z_{2k,j}(M,J_2)$ is a
vertically harmonic twistor lift.
\paragraph{In the odd case.}
We obtain an analogous interpretation as in the even case. An interpretation in terms  of horizontally holomorphic and vertically harmonic maps $f\colon L\to N=G/G_0$. Let us make precise that in the odd case, the action functional in the variational interpretation  has a Wess-Zumino term in addition to the vertical energy (See below in this introduction).\\
Moreover  by embedding $G/G_0$ into the twistor space $\mal{Z}_{2k+1}(N)$ of order $2k+1$ isometric endomorphisms in $TN$, we obtain an interpretation in terms of vertically harmonic twistor lift.
\index{vertically harmonic|)}\index{determined minimal@determined, minimal|)}
\subsubsection{The general structure of the maximal determined case.}\index{model@model case or system|(}
First, the maximal determined system has 3 model cases. This means that  we can distinguish 3 maximal determined systems, namely the three  maximal determined systems with lowest order of symmetry (2,3,4). Their corresponding geometric equations (when put all together) contain already all the structure terms - in a simple form- that will appear in the further maximal determined systems in a more complex and general form due to the more complex geometric structure in the further maximal determined systems. That is in this sense that we can say that all the further determined systems associated to target spaces $N$ with higher order of symmetry will be modeled on these model systems.
\subsubsection{The model system in the even case}\index{harmonic map}
In the even case, this model is the first elliptic integrable system associated to a symmetric space ($m=1, k'=2$) which is - as it is well known - exactly the equation of harmonic maps from the Riemann surface $L$ into the   symmetric space under consideration. This is the "smallest" determined system, i.e. with lowest order of symmetry in the target space $N=G/G_0$. In this case -$N$ is symmetric- the determined case is reduced to one system, the one of order 1.
\subsubsection{The model system in the odd case}\index{canonical!almost complex structure|(}
In the odd case, this model is the second elliptic integrable system associated to a 3-symmetric space. This is the "smallest" determined system in the odd case, i.e. with lowest odd order of symmetry in the target space $N=G/G_0$. We prove that this system is exactly the equation for \emph{holomorphically harmonic maps} into the almost complex manifold $(N,\undj)$ with respect to the anticanonical connection $\nabla^1=\nabla^0 +[\ ,\ ]_{[\mk]}$,\index{anticanonical connection} where $\nabla^0$ is the canonical connection. Or equivalently this is  the equation for holomorphically harmonic maps into the almost complex manifold $(N,-\undj)$ with respect to the canonical connection $\nabla^0$.
\paragraph{Holomorphically harmonic maps.}\index{holomorphically harmonic|(}
Given a general  almost complex manifold $(N,J)$ with a connection $\nabla$, we define holomorphically harmonic maps $f\colon L\to N$ as the solutions of the equation
\begin{equation}\label{eq-def-hol-harm}
\left[ \bar\partial^\nabla\partial f\right]^{1,0}=0
\end{equation} 
where $[\ ]^{1,0}$ denotes the $(1,0)$-component according to the splitting $TN^\C=T^{1,0}N\oplus T^{0,1}N$  defined by $J$. This equation is equivalent to 
$$
d^\nabla df+ Jd^\nabla * df=0
$$
or, equivalently, using any Hermitian metric $g$ on $L$
$$
T_g(f) + J\tau_g(f)=0
$$
where $T_g(f)=*f^*T=f^*T(e_1,e_2)$, with $(e_1,e_2)$ an orthonormal basis of $TL$, and $\tau_g(f)= *d^\nabla *df=\mrm{Tr}_g(\nabla df)$ is the tension field of $f$.\index{tension field} Of course $\mrm{Tr}_g$ denotes the trace with respect to $g$, and the expression $\nabla df$ denotes the covariant derivative of $df$ with respect to the connection induced in $T^*L\otimes f^*TN$ by $\nabla$ and the Levi-Civita connection in $L$.\\
In particular, we see that if $\nabla$ is torsion free or more generally if $f$ is torsion free, i.e. $f^*T=0$, then holomorphic harmonicity is equivalent to (affine) harmonicity. Therefore, this new notion is interesting only in the case of a non torsion free connection $\nabla$.
\paragraph{The vanishing of some $\bar\partial\partial$-derivative.}\index{dd derivative@$\bar\partial\partial$-derivative}
Now, let us suppose that the connection $\nabla$ on $N$ is \emph{almost complex}, i.e. $\nabla J=0$. Then, according to equation~(\ref{eq-def-hol-harm}), we see that any holomorphic curve $f\colon (L,j_L)\to (N,J)$ is anti-holomorphically harmonic, i.e. holomorphically harmonic with respect to $-J$. In particular, this allows to recover that a 1-primitive solution (i.e. of order $m=1$) of the elliptic system associated to a 3-symmetric space is also solution of the  second elliptic system associated to this space.\medskip\\
Moreover, the holomorphically harmonic maps admit a formulation very analogous to that of harmonic maps in terms of the vanishing of some $\bar\partial\partial$-derivative, which implies a well kown characterisation in terms of holomorphic 1-forms. Indeed we prove that $f\colon (L,j_L) \to (N,J,\nabla)$  is holomorphically harmonic \iif 
\begin{equation}\label{d-bar-d-f=0}
\overline{\hat\partial}^{\widehat\nabla}\hat\partial f =0,
\end{equation}
i.e. $\hat\partial f$ is a holomorphic section of $T_{1,0}^*L\otimes_\C f^*TN$. Here the hat " $\hat{} $ " means that we extend a 1-form on $TL$, like $d$ or $\nabla$, by $\C$-linearity as a linear map from $TL^\C$ into the complex bundle $(TN,J)$. In other words instead of extending these 1-forms as $\C$-linear maps from $TL^\C$ into $TN^\C$ as it is usual, we use the already existing structure of complex vector bundle in $(TN,J)$  and extend these very naturally as  $\C$-linear map from $TL^\C$ into the complex bundle $(TN,J)$. Therefore we can conclude that holomorphically harmonic maps have the same formulation as harmonic maps with the difference that instead of working in the complex vector bundle  $TN^\C$, we stay in $TN$ which is already a complex vector bundle in which we work.
\paragraph{The sigma model with a Wess-Zumino term.}
\index{Wess-Zumino|(} \index{sigma model with a Wess-Zumino term|see{Wess-Zumino}}
Finally, let us suppose that $N$ is endowed with a $\nabla$-parallel Hermitian metric $h$. Therefore $(N,J,h)$ is an almost Hermitian manifold with a Hermitian connection $\nabla$. Suppose also that $J$ anticommutes with the torsion $T$ of $\nabla$ i.e.
$$
T(X,JY)=-JT(X,Y)
$$
which is equivalent to 
$$
T=\dfrac{1}{4}N_J 
$$
where $N_J$ denotes the torsion of $J$ i.e its Nijenhuis  tensor.\\
Suppose also that the torsion of $\nabla$ is totally skew-symmetric i.e. the trilinear map\index{skew-symmetric torsion|(}
$$
T^*(X,Y,Z)=\langle T(X,Y),Z\rangle
$$
is a 3-form. Lastly, we suppose that the torsion is $\nabla$-parallel, i.e. $\nabla T^*=0$ which is equivalent to $\nabla T=0$.  Then we prove that this implies that the 3-form 
$$
H(X,Y,Z)=-T^*(X,Y,JZ)=\langle JT(X,Y),Z\rangle
$$ 
is closed $dH=0$.\footnote{Let us point out that in general $T^*$ is not closed even if it is $\nabla$-parallel. For example, in a Riemannian naturally reductive homogeneous space $G/H$,  endowed with its canonical connection $\nabla^0$, we have $\nabla^0 T=0$ but $dT^*(X,Y,Z,V)= -2\langle \mrm{Jac}_\mk(X,Y,Z),V\rangle$ where $\mrm{Jac}_\mk$ is the $\mk$-component of the Jacobi identity (i.e. the sum of the circular permutations of $[X,[Y,Z]_\mk]_\mk$). Of course $\mk$ denotes the $\Ad H$-invariant summand in the reductive decomposition $\g=\hk\oplus \mk$.}\\
Then \emph{the equation for holomorphically harmonic maps $f\colon L\to N$ is the  equation of motion (i.e. Euler-Lagrange equation) for the sigma model in $N$ with the Wess-Zumino term defined by the closed 3-form $H$.}
The action functional is given by
$$
S(f)=E(f)  + S^{WZ}(f)= \dfrac{1}{2}\int_L|df|^2 d\mathrm{vol}_g + \int_B H,
$$
where $B$ is 3-submanifold (or indeed a 3-chain) in $N$ whose boundary is $f(L)$.\\
Then since $dH=0$, the variation of the Wess-Zumino term is a boundary term 
$$
\delta S^{WZ}=\int_B L_{\delta f}H=\int_Bd\imath_{\delta f}H=\int_{f(L)}\imath_{\delta f}H,
$$
whence its contribution  to the Euler-Lagrange equation  involves only the original map $f\colon L\to N$.\\
In particular, applying this result to the case we are interested in, i.e. $N$ is 3-symmetric, we obtain:\\
\emph{The second elliptic system associated to a 3-symmetric space $N=G/G_0$ is the equation of motion for the sigma model in $N$ with the Wess-Zumino term defined by the closed 3-form $H(X,Y,Z):=T^*(X,Y,\undj Z)$, where $T$ is the torsion of the canonical connection $\nabla^0$ and $\undj$ is the canonical almost complex structure.}%
\footnote{In fact, we need a naturally reductive metric on $N$ to ensure that $T^*$ is a 3-form. But if we allow Pseudo-Riemannian metrics and if $\g$ is semisimple then the metric defined by the Killing form is naturally reductive. In fact, the elliptic integrable system is a priori written in an affine context, i.e. its natural - in the sense of initial- geometric interpretation takes place in the context of affine geometry in terms of the linear connections $\nabla^t=\nabla^0 + t[\ ,\ ]_{[\mk]}$. If we want that this interpretation takes place in the context of Riemannian geometry we need, of course, to add some hypothesis of compactness, like the compactness of $\Ad_\mk G_0$ and the natural reductivity. But we do not need  these hypothesis if we work in the Pseudo-Riemannian context.}
\paragraph{The good geometric context/setting}
In the previous variational interpretation, we need to make 3 hypothesis on the torsion of the  Hermitian connection: $T$ anticommutes with $J$, is totally skew-symmetric and $\nabla$-parallel. It is natural to ask ourself what do these hypothesis mean geometrically and  what is the good geometric context in which these take place. It  turns out that the good geometric context is the one of \emph{Nearly K\"{a}hler manifold}.\\
An almost Hermitian manifold $(N,J,h)$ is Nearly K\"{a}hler \iif $(\nabla_X^h J)X=0$, for all $X\in TN$, where $\nabla^h$ is the Levi-Civita connection. Then we prove that  the almost Hermitian manifolds for which there exists an Hermitian connection satisfiying the three hypothesis above are exactly Nearly K\"{a}hler manifolds, and that this Hermitian connection is then unique and coincides with the canonical Hermitian connection. Then the variational interpretation can be rewritten as follows:
\begin{thm}
Let $(N,h,J)$ be a nearly K\"{a}hler manifold then the equation of holomorphic harmonicity, \wrt the canonical Hermitian connection, for maps $f\colon L\to N$ is exactly the Euler-Lagrange equation for the sigma model in $N$ with a Wess-Zumino term defined by the 3-form:
$$
H=\dfrac{1}{3}d\Omega_J
$$
where $\Omega_J=\langle J\cdot, \cdot\rangle$ is the K\"{a}hler form. 
\end{thm}
Therefore: \emph{the second elliptic system associated to a 3-symmetric space $N=G/G_0$, endowed with its  canonical almost complex structure $\undj$, is the equation of motion for the sigma model in $N$ with the Wess-Zumino term defined by the closed 3-form $H=-\dfrac{1}{3}d\Omega_{\undj}$.}\index{Wess-Zumino|)}
\index{skew-symmetric torsion|)}\index{totally skew-symmetric torsion|see{skew-symmetric torsion}}
\paragraph{$J$-twisted harmonic maps.}
We prove that we can also interpret the holomorphic harmonicity in terms of \emph{$J$-twisted harmonic maps} (\wrt the Levi-Civita connection). Let us define this notion. Let $(E,J)$ be a complex vector bundle over an almost complex manifold $(M,j_M)$. Then let $\overline\nabla$ be a connection on $E$. Then we can decompose it in an unique way as the sum of a $J$-commuting and a $J$-anticommuting part, i.e. in  the form
$$
\overline\nabla =\nabla^0 + A
$$
where $\nabla^0J=0$ and $A\in \mathcal C(T^*M\otimes \End(E))$, $AJ=-JA$. More precisely, we have $A=\dfrac{1}{2}J\overline\nabla J$. Then we set
$$
\overline\nabla^{J}=\nabla^0 - (A\circ j_M)J=\overline\nabla - \dfrac{1}{2}J\overline\nabla J -\dfrac{1}{2}\overline\nabla J \circ j_M.
$$
Now let $f\colon (L,j_L)\to (N,J)$ be a map from a Riemann surface into  the almost complex manifold $(N,J)$ endowed with a connection $\nabla$. Then let us take in what precede $(M,j_M)=(L,j_L)$ and $(E, \overline\nabla) = (f^*TN, f^*\nabla)$. Then we say that the map $f\colon (L,j_L)\to (N,J,\nabla)$ is $J$-twisted harmonic \iif
$$
\mrm{Tr}_g(\overline{\nabla}^J df)=0
$$
(for any hermitian metric $g$ on the Riemann surface $L$).\index{holomorphically harmonic|)}
\index{canonical!almost complex structure|)}
\subsubsection{The coupled model system}
This is the third elliptic integrable system associated to a 4-symmetric space. The corresponding geometric equation in the 4-symmetric space $G/G_0$ can be viewed as a coupling between the equation of  harmonic maps into the symmetric fibre $H/G_0$ and the equation for harmonic maps into the symmetric space $G/H$. In other words, it can be viewed as a coupling between the even model  system associated to the symmetric space $G/H$ and the even model system associated to the symmetric space $H/G_0$. In  particular, if harmonicity is heuristically replaced by holomorphicity "in the basis $G/H$", then we recover the horizontally holomorphic and vertically harmonic maps into $G/G_0$, that is to say the solutions of the second elliptic integrable system associated to $G/G_0$. We will come back to this in the end of  section~\ref{gene-max-deter-even-intro} in this introduction.
\index{model@model case or system|)}
\subsubsection{The General maximal determined odd system. ($k'=2k+1,m=2k$)}
The maximal determined odd system has a  geometric interpretation in terms of \emph{stringy harmonic}\index{stringy harmonic|(} maps $f\colon L\to (G/G_0,\undj)$, with respect to the canonical connection and the canonical almost complex structure.\index{canonical!almost complex structure}
\paragraph{Stringy harmonic maps.}
Let $(N,J)$ be an almost complex manifold with $\nabla$ an linear connection then we will say that a map $f\colon L\to N$ from a Riemann surface into $N$ is \emph{stringy harmonic} if it is a solution of \emph{the harmonic map equation with a $JT$-term}:
$$
-\tau_g(f) + (J\cdot T)_g(f)=0.
$$
We have used the notation $J\cdot B=-JB(J\cdot,\cdot)$,  $\forall B\in \mal C(\Lambda^2 T^*N\otimes TN)$. This action of $J$ on  $ B\in \mal C(\Lambda^2 T^*N \otimes TN)$ can be written more naturally if $(N,J)$ is endowed with a  Hermitian metric $h$. Indeed, in this case, we have an identification, $ \mal C(\Lambda^2 T^*N\otimes TN) = \mal C(\Lambda^2 T^*N\otimes T^*N)\subset \mal C(\otimes^3 T^*N)$, between $TN$-valued 2-forms on $N$ and trilinear forms on $N$ skew-symmetric \wrt the 2 first variables: $B(X,Y,Z):=\langle B(X,Y),Z\rangle$. Then $J\cdot B$ is written:
$$
J\cdot B= B(J\cdot,J\cdot,J\cdot)=:-B^c .
$$
We remark that  if $T$ anticommutes with $J$ then stringy harmoniciy coincides with holomorphic harmonicity (since in this case $J\cdot T= JT$). More particulary, if $T=0$, then the stringy harmonicity coincides with the harmonicity.\\
Furthermore, we look for a general geometric setting  in which the stringy harmonicity has an interesting interpretation. First of all, let us remark that in the context of homogeneous reductive space, in which our system takes place, we have a canonical connection, with respect to which the stringy harmonicity can be written "canonically". But in general we do not have a "special" connection with respect to which one can consider the stringy harmonicity. Therefore, if one wants to place stringy harmonicity in a more meaningfull, interesting and fruitful  context (than the general context of almost complex manifolds endowed with some linear connection) and, in so doing, obtain a better understanding of our elliptic integrable system by writting its geometric interpretation in the best geometric context, a first problem - that we solved - is to find a general class of (almost complex) manifold in which there exists some unique "canonical" connection, with respect to which we then could consider the stringy harmonicity. This   provides us, firstly, some special connection (in the same sense that the Levi-Civita connection is special in Riemannian geometry), which solves the problem of the choice of the connection, but secondly it turns out that it provides also a variationnal interpretation of the stringy harmonicity. 
\paragraph{Best geometric setting}
\index{connection!characteristic|(}\index{g1 manifold@$\mal G_1$-manifold}\index{skew-symmetric torsion|(}
It  turns out that the more rich geometric context in which stringy harmonicity admits interesting properties  is the one of $\mal G_1$-manifolds, more precisely $\mal G_1$-manifolds whose the characteristic connection has a parallel torsion.\index{Gray-Hervella classification} Making systematic use of the covariant derivative of the K\"{a}hler form, A. Gray and L. M. Hervella, in the late seventies, classified almost Hermitian structures into sixteen classes \cite{Gray-Hervella}. Denote by $\mal W$ the space of all trilinear forms (on some Hermitian vector space, say $T_{y_0}N$ for some reference point $y_0\in N$) having the same algebraic properties as $\nabla^h\Omega_J$. Then they proved that we have a $U(n)$-irreducible decomposition $\mal W=\mal W_1\oplus \mal W_2\oplus\mal W_3\oplus\mal W_4$. The sixteen classes are then respectively  the classes of almost Hermitian manifolds for which $\nabla^h\Omega_J$ `lies in' the $U(n)$-invariant subspaces $\mal W_I=\oplus_{i\in I}\mal W_i$, $I\subset \{1,\ldots,4\}$, respectively. In particular, if we take as invariant subspace $\{0\}$, we obtain the K\"{a}hler manifolds, if we take $\mal W_1$, we obtain the class of nearly K\"ahler manifolds. Moreover the class of $\mal G_1$-manifolds is the one defined by $\mal G_1= \mal W_1\oplus\mal W_3\oplus\mal W_4$. It is characterised by : $(N,J,h)$ is of type $\mal G_1$ \iif the Nijenhuis tensor $N_J$ is totally skew-symmetric (i.e. a 3-form).\index{Nijenhuis tensor!} \\
In this paper, we prove the following theorem\footnote{After the fact, we  realized that it has already been proved by Friedrich-Ivanov\cite{Friedrich-ivanov}. However we give a different and completely written proof. See remark~\ref{rmk-frie-ivan}}:
\begin{thm}
An almost  Hermitian manifold $(N,J,h)$ admits a Hermitian connection with totally skew-symmetric  torsion \iif the Nijenhuis tensor $N_J$ is itself totally skew-symmetric. In this case, the connection is unique and determined by its torsion which is given by 
$$
T=-d^c\Omega_J + N_J.
$$
The characteristic connection  is then given by $\nabla =\nabla^h-\frac{1}{2}T$.
\end{thm}
Then we prove:  
\begin{prop}
Let us suppose that the almost Hermitian manifold $(N,J,h)$ is a $\mal G_1$-manifold. Let us suppose that its  characteristic connection $\nabla$ has a parallel torsion  $\nabla T=0$. Then the 3-form 
$$
H(X,Y,Z)=T(JX,JY,JZ)=\langle (J\cdot T)(X,Y),Z\rangle
$$
is closed $dH=0$.
\end{prop}
Which then gives us  the following variational interpretation
\begin{thm}\index{Wess-Zumino|(}
Let us suppose that the almost Hermitian manifold $(N,J,h)$ is a $\mal G_1$-manifold. Let us suppose that its  characteristic connection $\nabla$ has a parallel torsion  $\nabla T=0$.\\
Then the equation for stringy harmonic maps $f\colon L\to N$ is exactly the Euler-Lagrange equation for the sigma model in $N$ with a Wess-Zumino term defined by the closed 3-form
$$
H=-d\Omega_J + JN_J.
$$ 
\end{thm}
\index{canonical!almost complex structure}
Moreover  any $(2k+1)$-symmetric space $(G/G_0,\undj,h)$ endowed with its canonical almost complex structure and a naturally reductive $G$-invariant metric $h$ (for which $\undj$ is orthogonal) is a $\mal G_1$-manifold and moreover its characteristic connection coincides with its canonical connection $\nabla^0$. Finally, the torsion of the canonical connection is obviously parallel. Therefore we obtain an interpretation of the maximal determined  system associated to a $(2k+1)$-symmetric space in terms of a sigma model with a Wess-Zumino term.
\begin{rmk}\em
Let us add  about stringy harmonicity that we  also prove, in this paper, that the stringy harmonicity \wrt an almost complex connection $\nabla$ is equivalent to the holomorphic harmonicity \wrt a new almost complex connection $\nabla^\star$. \index{holomorphically harmonic}
\end{rmk}
\index{g1 manifold@$\mal G_1$-manifold|)}\index{skew-symmetric torsion|)}
\subsubsection{General maximal determined even system. ($k'=2k,m=2k-1$)}\label{gene-max-deter-even-intro}
\index{canonical!$f$-structure}\index{horizontal subbundle}\index{vertical subbundle}
In the even case, the geometric structure of the target space $G/G_0$ is more complex (and more rich): as we have already seen, there is a fibration $\pi\colon N=G/G_0\to M=G/H$,  a splitting $TN=\hor\oplus\ver$ with $\ver=\ker\pi$, and an $f$-structure $F$ such that $\ker F=\ver$ and $\im F=\hor$ (in particular $\Bar J:=F_\hor$ is a complex structure on $\hor$). Moreover the geometric PDE obtained from our elliptic integrable system uses this geometric structure. In particular, this geometric PDE splits into its horizontal and vertical parts and can be viewed as a coupling between the equation of $\Bar J$-stringy harmonicity and the equation of vertical harmonicity, the coupling terms calling out the curvature of $\hor$.\medskip\\
The maximal determined even system has a  geometric interpretation in terms of \emph{stringy harmonic} maps $f\colon L\to (G/G_0,F)$, $F$ being the canonical $f$-structure on $G/G_0$.
\paragraph{Stringy harmonic maps \wrt an $f$-structure.}
Let  $(N,F)$ be an $f$-manifold with $\nabla$ a linear connection. Then we will say that a map $f\colon L\to N$ from a Riemann surface into $N$ is \emph{stringy harmonic} if it is solution of \emph{the stringy harmonic maps equation}:
$$
-\tau_g(f) + (F\bullet T)_g(f)=0.
$$
where $F\bullet B$, for $B\in \mal C(\Lambda^2T^*N\otimes TN)$, denotes some natural (linear) action of $F$ on $\mal C(\Lambda^2T^*N\otimes TN)$. For more simplicity, let us write it in the case where $(N,F)$ is endowed with a compatible metric $h$ (i.e. $\ver\perp\hor$ and $\Bar J=F_{|\hor}$ is orthogonal with respect to $h_{|\hor\times\hor}$):
\begin{eqnarray*}
F\bullet B & =  & B(F\cdot,F\cdot,F\cdot) + \dfrac{1}{2}F\circact (B-B_{\hor^3})\\
F\circact A & = &  A(F\cdot, \cdot,\cdot) + A(\cdot,F\cdot,\cdot) + A(\cdot,\cdot,F\cdot)
\end{eqnarray*}
for all $B,A\in\mal C(\Lambda^2 T^*N\otimes TN)$.\\
Now, we want to proceed as in the case of stringy harmonicity with respect to an almost complex structure. That is to say: to find a class of $f$-manifolds for which there exists some unique characteristic connection which preserves the structure and then to look for a variational interpretation of the stringy harmonicity with respect to this connection.
\paragraph{Best Geometric context.}\index{global type $\mal G_1$|(}\index{skew-symmetric torsion|(}
We look for metric $f$-manifolds $(N,F,h)$ for which there exists a metric $f$-connection $\nabla$ (i.e. $\nabla F=0$ and $\nabla h=0$) with skew-symmetric torsion $T$. In a first step, we consider metric connections which preserve the splitting $TN=\ver\oplus\hor$ (i.e. $\nabla q=0$, where $q$ is  the projection on $\ver$) and we characterize the manifolds $(N,h,q)$ for which there exists such a connection with skew-symmetric torsion, and call these \emph{reductive} metric $f$-manifolds.\\
Then, saying about a metric $f$-manifolds $(N,F,h)$ that it is of global type $\mal G_1$ if its extended Nijenhuis tensor $\tl N_F$\index{Nijenhuis tensor!extended} is skew-symmetric, we prove the following theorem:
\begin{thm}\label{characteristic-equiv-G1} 
A metric $f$-manifold $(N,F,h)$ admits a metric $f$-connection $\nabla$ with skew-symmetric torsion \iif it is reductive and of global type $\mal G_1$. Moreover, in this case, for any $\alpha\in \mal C(\Lambda^3\ver^*)$, there exists a unique metric connection $\nabla$ with skew-symmetric torsion such that $T_{|\Lambda^3\ver}=\alpha$. This unique connection is given by 
$$
T=(-d^c\Omega_F + N_{F|\hor^3}) + \mrm{Skew}(\Phi) + \mrm{Skew}(\mR_\ver) + \alpha.
$$
where $\Omega_F=\langle F\cdot, \cdot\rangle$, $\Phi$ and $\mR_\ver$ are resp. the curvature of $\hor$ and $\ver$ resp., and $\mrm{Skew}$ is the sum of all the circular permutations on the three variables.
\end{thm}
On a metric $f$-manifold $(N,F,h)$, a metric $f$-connection $\nabla$ with skew-symmetric torsion will be called  a \emph{characteristic connection}.\medskip\\ 
Moreover, we prove that for any reductive metric $f$-manifold of global type $\mal G_1$, the closure of $H=F\bullet T$ is equivalent to the closure of the  3-form $F\cdot N_F -\dfrac{1}{2}F{\circact} \left(\mrm{Skew}(\Phi) + \mrm{Skew}(\mR_\ver)\right)$, so that:
\begin{thm}
Let $(N,F,h)$ be a reductive metric $f$-manifold of global type $\mal G_1$. Let us suppose that  the  3-form $F\cdot N_F -\dfrac{1}{2}F{\circact} \left(\mrm{Skew}(\Phi) + \mrm{Skew}(\mR_\ver)\right)$ is closed, where $\Phi$ and $\mR_\ver$ are respectively the curvatures of the horizontal and vertical subbundles. Let $\nabla$ be one characteristic connection.\\
Then the equation for stringy harmonic maps,  $f\colon L\to N$, (\wrt $\nabla$) is exactly the Euler-Lagrange equation for the sigma model in $N$ with a Wess-Zumino term defined by the closed 3-form
$$
H=-d\Omega_F + F\cdot N_F -\dfrac{1}{2}F{\circact} \left(\mrm{Skew}(\Phi) + \mrm{Skew}(\mR_\ver)\right).
$$ 
\end{thm}
Contrary to the case of stringy harmonic maps into an almost Hermitian $\mal G_1$-manifolds, in the present case, the hypothesis that the torsion of one characteristic connection is parallel $\nabla T=0$ does not imply the closure of the 3-form $H=F\bullet T$. However, we characterize this closure under the hypothesis $\nabla T=0$ and $\mR_\ver=0$, by some 2 conditions that we will not explain in this introduction (see section~\ref{dH=0}): the horizontal complex structure $\Bar J$ is a cyclic permutation of the horizontal curvature,  and  the 2-forms $N_{\Bar J}$ and $\Phi$ have orthogonal supports.
\begin{thm}
Let $(N,F,h)$ be a reductive metric $f$-manifold of global type $\mal G_1$. Let us suppose that one of its  characteristic connections, $\nabla$,  has a parallel torsion $\nabla T=0$. Let us supppose that $\mR_\ver=0$ and that the horizontal curvature $\Phi$ is pure. The following statement are equivalent:\medskip\\
$\bullet$ The horizontal 3-form $F\cdot N_F$ is closed.\smallskip\\
$\bullet$ $F\cdot N_F$ and $ F{\circact} \left(\mrm{Skew}(\Phi) + \mrm{Skew}(\mR_\ver)\right)$ are closed.\smallskip\\
$\bullet$ The horizontal complex structure $\Bar J$ is a cyclic permutation of the horizontal curvature,  and  the 2-forms $N_{\Bar J}$ and $\Phi$ have orthogonal supports.\smallskip\\
In this case, $H=F\bullet T$ is closed.
\end{thm}
\index{canonical!$f$-structure|(}
Moreover any $2k$-symmetric space $(G/G_0,F,h)$ endowed with its canonical $f$-structure and a naturally reductive $G$-invariant metric $h$ (compatible with $F$) is reductive and of global type $\mal G_1$, and moreover its canonical connection $\nabla^0$ is a characteristic connection. Furthermore, the torsion of the canonical connection is obviously parallel. Finally, we prove that any $2k$-symmetric space $(G/G_0,F,h)$ satisfies the two hypothesis above. Therefore we obtain an interpretation of the maximal determined even system associated to a $2k$-symmetric space in terms of a sigma model with a Wess-Zumino term. 

\paragraph{A particular case: Horizontally K\"{a}hler $f$-manifolds.} \index{horizontally K\"{a}hler}
Let $(N,F,h)$ be a metric $f$-manifold. We will say that $(N,F,h)$ is  horizontally K\"ahler if $DF_{|\hor^3}=0$, where $D$ is the Levi-Civita connection of $h$.\\
Then, we prove that in this case, the two hypothesis above (which characterise the closure of the 3-form $H$) are satisfied. Moreover, any characteristic connection $\nabla$ satisfies $T_{\hor^3}=0$, which leads to special properties.\\
A example of this situation is given by any 4-symmetric space, endowed with its canonical $f$-structure and a naturally reductive $G$-invariant metric (compatible with $F$). 
\index{determined maximal@determined, maximal|)}\index{connection!characteristic|)}\index{global type $\mal G_1$|)}
\index{skew-symmetric torsion|)}
\subsubsection{The intermediate determined systems}
For the intermediate determined systems ($m_{k'}<m<k'-1$), these are obtained from the maximal determined case by adding holomorphicity in the subbundle $ \hor^{ m}=\oplus_{j=1}^{\und m}[\mk_j]\subset \hor$, where $\und m=k'-1-m$. It means that the $m$-th determined system has a geometric interpretation in terms of  stringy harmonic maps\index{stringy harmonic|)} which are $\hor^{ m}$-holomorphic: 
$$
(df\circ j_L)^{\hor^{m}}=F^{[\und m]}\circ df.
$$
We have seen that the maximal determined system has an interpretation in terms of a sigma model with a Wess-Zumino term defined by a 3-form $H$. In fact, more generally, let $m_{k'}\leq m \leq k' -1$,  and let us consider the splitting $TN=\hor^m \oplus\ver^m$ defined above. Then one can prove that  any $m$-th determined system is the Euler-Lagrange equation \wrt vertical variations (i.e. in $\ver^m$) of the following functional
$$
E^{\bar v}(f)= \dfrac{1}{2}\int_L |d^{\bar v } f|^2 d\mrm{vol}_g  + \int_B H^{\bar v }
$$
where $d^{\bar v } f= [df]^{\ver^m}$, $H^{\bar v }=H-\ovr H= H_{|\mal S(\ver^m,\hor^m)}$, $\ovr H=H_{|(\hor^m)^3}$, and $B$ is a 3-submanifold of $N$ with boundary $\partial B= f(L)$.\index{Wess-Zumino|)}
\index{canonical!$f$-structure|)}
\subsection{The underdetermined case} \index{underdetermined@underdetermined, case or system}
We prove that the $m$-th underdetermined system ($m> k'-1$) is in fact equivalent to some $m$-th determined or primitive system associated to some new automorphism $\tl\tau$ defined in a product $\g^{q+1}$, of the initial Lie algebra $\g$. More precisely, we write
$$ 
m=qk' + r,\quad 0\leq r \leq k'-1
$$
the Euclidean division of $m$ by $k'$. Then we consider the automorphism in $\g^{q+1}$ defined by
$$
\tl\tau \colon (a_0,a_1,\ldots,a_{q})\in\g^{q+1}\longmapsto (a_{1}, \ldots,a_{q},\tau(a_0))\in\g^{q+1}.
$$
Then $\tl\tau$ is of order $(q+1)k'$. We prove that the initial $m$-th (underdetermined) system associated to $(\g,\tau)$ is in fact equivalent to the $m$-th (determined) system associated to $(\g^{q+1},\tl\tau)$.
\subsection{In the twistor space.}\index{twistor|(}
For each previous geometric interpretation in the target space $N=G/G_0$, there is a corresponding geometric interpretation  in the twistor space.\\
$\bullet$ These previous geometric interpretations take place in some  manifolds endowed with some particular structure. This could simply be, for example,  a structure of almost complex manifold in the case of the interpretation of the stringy harmonicity in terms of the vanishing of some $\bar\partial\partial$-derivative  but it could also be  the more strong structure of $\mal G_1$-manifolds whose characteristic connection has a parallel torsion). Moreover, our $k'$-symmetric spaces are very particular examples of this kind of manifolds. Therefore it is natural to try to make these interpretations more universal by writting them in a more general setting. More precisely we want to find some universal prototype of these "special" manifolds, which can be endowed canonically with the needed geometric structure and such that any of our special manifolds can be embedded in this prototype.\\
Indeed, as concerns $k'$-symmetric spaces, we know that they can be embedded canonically into some twistor bundles. In the even case we have an injective morphim of bundle over $M=G/H$ defined by the embedding\index{canonical!embedding}
$$
G/G_0\hookrightarrow \zdk(G/H),
$$
whereas in the odd case we have a section defined by the embedding
$$
G/G_0\hookrightarrow \mal Z_{2k+1}(G/G_0).
$$
$\bullet$ In the even case, the fibration $\pi\colon G/G_0\to G/H$  imposes to view canonically any $2k$-symmetric space as a subbundle of $\zdk(G/H)$ so that the twistorial interpretation is in some sense dictated by the structure of the $2k$-symmetric space.
The geometric interpretations in the twistor spaces are universal since these twistor spaces are defined for any Riemannian manifold and are endowed canonically with the different geometric structures that we need. That is to say the geometric structures  we need to endow the target space $N$ with, in our previous geometric interpretations.\\
More generally, suppose that we want to study stringy harmonicity in  metric $f$-manifolds $(N,F,h)$. It is then natural to consider the particular case where the vertical subbundle $\ver$ is the tangent space to the fibre of a Riemannian submersion $\pi\colon (N,h)\to (M,g)$, i.e. $\ver=\ker d\pi$. For example, we have  seen that among the list of sufficient conditions in our variationnal interpretation of stringy harmonicity there is the condition $\mR_\ver=0$.  In the particular case of a Riemannian submersion, the $f$-structure $F$ defines a complex structure $\Bar J$ on $\pi^*TM=\hor$ which itself gives rise to a morphism of submersion $\mal I \colon N\to \Sigma (M)$, $y\mapsto (\pi(y),\Bar J(y))$. This shows that the twistor bundle $\Sigma(M)$ appears naturally in the general context - even though the morphism $\mal I$ is not injective in general.\\
Furthermore, an interesting class of Riemannian submersion $\pi\colon (N,h)\to (M,g)$ is the one of homogeneous fibre bundles, of which the twistor bundles $\mal Z_p(M)$ are particular examples ($p\in \mathbb N^*$). For example, vertically harmonic sections of homogeneous fibre bundles have been investigated by C.M. Wood \cite{cmw1,cmw2}.
\index{twistor|)}\index{even case|)}\index{odd case|)}\index{homogeneous fibre bundle}
\subsection{Related subjects and works, and motivations}
\subsubsection{Relations with surface theory.}
The theory of harmonic maps of surfaces has been greatly enriched by ideas and methods
from integrable systems \cite{bpw,BFPP,harmonicTori, 12,DPW,H2,H3,Uhlenbeck}. In particular, these ideas have revolutionised the theory of harmonic maps from a surface
into a symmetric space and so, via an appropriate Gauss map, the theory of constant
mean curvature surfaces and Willmore surfaces among others. For example, Pinkall and
Sterling \cite{Pink-Ster} were able to give an algebro-geometric construction of all constant mean curvature
tori while Dorfmeister-Pedit-Wu \cite{DPW}  gave a Weierstrass formula for all constant mean
curvature immersions of any (simply connected) surface in terms of holomorphic data.
These advances were taken up by H\'{e}lein and Romon \cite{HR1,HR2,HR3} who showed that similar ideas could
be applied to the study of Hamiltonian stationary Lagrangian surfaces in a 4-dimensional
Hermitian symmetric space. It was the first example of second elliptic integrable system associated to a 4-symmetric space. In \cite{ki1}, we presented a new class of isotropic surfaces in the Euclidean space of dimension 8 by identifying $\R^8$ with the set of octonions $\oct$, and we proved that these surfaces are solutions of a second elliptic integrable system associated to a  4-symmetric space. By restriction to $\R^4=\h$, we obtained the Hamiltonian stationary Lagrangian surfaces and by restriction to $\R^3=\im \h$, we obtained the CMC surfaces.  Furthermore, in \cite{ki3}, we presented a geometric interpretation of all the second elliptic integrable systems associated to a 4-symmetric space in terms of vertically harmonic twistor lifts of conformal immersions  into the Riemannian symmetric space  (associated to our 4-symmetric space) (see also \cite{bk}). When the previous Riemannian symmetric space   is 4-dimensional, then any conformal immersion admits an unique twistor lift and the vertical harmonicity of this twistor lift is equivalent to the holomorphicity of the mean curvature vector of the conformal immersion (see \cite{bk}). In particular, when the Riemannian symmetric space is Hermitian, one obtains a conceptual explanation of the result of Hélein-Romon.

\subsubsection{Relations with mathematical physics.}
\paragraph{Metric connections  with totally skew-symmetric torsion.}
We refer the reader to \cite{srni,Friedrich-ivanov} about connections with  skew-symmetric torsions and their relations to physics and more particulary  string theory. Linear metric connections  with totally skew-symmetric torsion recently became a subject of interest in theoretical and mathematical physics. Let us give here some examples (taken from \cite{Friedrich-ivanov}). \medskip\\
$\bullet$ The target space of supersymmetric sigma models with Wess-Zumino term carries a geometry of a metric connection with skew-symmetric torsion.\\
$\bullet$ In supergravity theories, the geometry of the moduli space of a class of black holes is
carried out by a metric connection with skew-symmetric torsion.\\
$\bullet$ The geometry of NS-5 brane solutions of type II supergravity theories is generated by a metric connection with skew-symmetric torsion.\\
$\bullet$ The existence of parallel spinors with respect to a metric connection with skew-symmetric torsion on a Riemannian spin manifold is of importance in string theory, since they are associated with some string solitons (BPS solitons).
\paragraph{The sigma-models} 
Nonlinear sigma-models provide a much-studied class of field theories of both
phenomenological and theoretical interest. The chiral model for example summarizes
many low energy QCD interactions while 2-dimensional sigma-models may
possess nontrivial classical field configurations and have analogies with 4-dimensional
Yang-Mills equations but are simpler to handle. It was in the study of the chiral model that Wess and Zumino
introduced their parity violating term satisfying anomalous Ward-identities. 
This term now has a far broader interpretation in terms of anomalies. Sigma-models
also have connections with string theories.\\
Two dimensional sigma-models have already proven a fertile arena for the interplay of
topology, geometry, and physics. 
The  (supersymmetric) sigma-models have been generalized by introducing a Wess-Zumino term into the Lagrangian. This term may be interpreted as adding torsion to the canonical Levi-Civita connection of the
earlier models. The addition of such torsion can have a marked effect and imposes constraints on the possible geometries of the target \cite{braden}.\smallskip\\
\emph{Our contribution.} We give new examples of integrable two-dimensional non linear sigma models. These new examples take place in some homogeneous spaces, namely $k'$-symmetric spaces, which are not symmetric spaces. At our knowledge, all the already known integrable two-dimensional non linear sigma models take place in symmetric spaces or (equivalently) in Lie groups.

\subsubsection{Relations of $F$-stringy harmonicity and supersymmetry}
The P.D.E of $F$-stringy harmonicity splits following the splitting $TN=\hor\oplus\ver$ defined by the $f$-structure $F$.\\ More precisely, this equation is a coupling between the equation of $\Bar J$-stringy harmonicity and the equation of vertical harmonicity, the coupling terms calling out the curvatures of $\hor$ and $\ver$, and the component $N_{F|\ver\times\hor}$ of the Nijenhuis tensor.\smallskip\\
Moreover, let us suppose that we have  we have a fibration $\pi \colon N \to M$. Then we have a supersymmetric interpretation of the $F$-stringy harmonicity: $F$-stringy harmonicity can be viewed as a supersymmetric extention of the $J$-stringy harmonicity. In the splitting $TN=\hor\oplus\ver$, the  horizontal subbundle played the role of the odd part and the vertical subbundle plays the role of the even part. In other words,  the bosonic equation is a harmonic map equation (the vertical harmonicity) and  the fermionic equation is the $\Bar J$-stringy harmonic map equation.\smallskip\\
Let us also mention that, in \cite{ki2}, we obtained a supersymmetric interpretation of all the second elliptic integrable systems associated to a 4-symmetric space in terms of superharmonic maps.\index{k symmetric space@$k'$-symmetric space|)}
\medskip\\

\noindent
\textbf{Aknowledgements} The author wishes to thank Josef Dorfmeister for his useful comments on the first parts of this paper. He is also  grateful to him for his interest in the present work, his encouragements as well as his support during the preparation of this paper.

%%%%
%%%%
%with an almost complex connection $\nabla$ (i.e. $\nabla J=0$)
%or equivalently $d^\nabla T=0$. This hypothsis is in particular satisfied if the torsion is $\nabla$-parallel, %$\nabla T=0$

%
%%%%%%%%%%%%%%%%%%%%%%%%%%%%%%%%%%%%%%%%%%%%%%%%%%%%%%%%%%%%%%%%%%%%%%%%%

%     Index of notations

%%%%%%%%%%%%%%%%%%%%%%%%%%%%%%%%%%%%%%%%%%%%%%%%%%%%%%%%%%%%%%%%%%%%%%%

\section*{Notations, Conventions and general definitions.}
\addcontentsline{toc}{section}{\protect\numberline{0}Notations, Conventions and general definitions}
\subsection{List of notational conventions and organisation of the paper.}\label{generalities}
$\bullet$ Let $k\in \mathbb N^*$. Then we will often confuse - when it is convenient to do it- an element in $\Z_k$ with one of its representants. For example, let $(a_i)_{i\in\Z_{k}}$ be a family of elements in some vector space $E$, and $0\leq m < k/2$ an integer. Then we will write 
$$
a_i=a_{-i}\quad 1\leq i\leq m
$$
to say that this equality holds for all $i\in \{1+k\Z,\ldots,m+k\Z\}\subset \Z_{k}$.
\\
$\bullet$ Let us suppose that a vector space $E$ admits some decomposition $E=\oplus_{i\in I} E_i$. Then, for any vector $v\in E$ we denote by $[v]_{E_i}$ its component in $E_i$. It could also happen sometimes that we write this component $[v]^{E_i}$.
\\
$\bullet$ We denote $E^\C:=E\otimes \C$ the complexification of a real vector space. 
\\
$\bullet$ Furthermore let $A\in \End(E)$ be an endomorphism of a finite dimensional real vector space that we suppose to be diagonalizable. Then if its complex spectrum is $\{\lm_i, i\in I\}$, we will write its eigenspace decomposition in the form $E^\C=\oplus_{i\in I} E_i^\C $, where $E_i^\C=\ker (A - \lm_i \Id)$. Moreover, if $\lm_i\in \R$ for some $i\in I$, then we set $E_i:=E_i^\C\cap E$ so that $E_i^\C=(E_i)^\C$, for these particular $i\in I$. 
\\
$\bullet$ We denote by $\omega_p$ a $p$-th primitive root of unity, which will be often chosen equal to $e^{2i\pi/p}$.
\\
$\bullet$ We denote by $C^\infty(M,N)$ the set of smooth maps from a manifold $M$ into a manifold $N$. Now, let $\pi\colon N\to M$ be a surjection, then we denote by $\mal C(\pi)$ the set of sections of $\pi$, i.e. the maps $s$ such that $\pi\circ s=\Id_M$. If there is no risk of confusion we will also use the notation $\mal C(N)$. Furthermore, let $p\colon E\to M$ be a vector bundle. Then we will denote by the same letter $p$ its tensorial extensions: $p\colon \End(E)\to M$ and so on.
\\
$\bullet$ Let $(M,g)$ be a Riemannian manifold. Then we denote by $*$ its Hodge operator (and $*_g$ if the metric need to be precised). Let $(E,h)\to M$ be a  Riemannian vector bundle over a manifold $M$, we denote by $\Sigma(E)$ the bundle of orthogonal almost complex structures in $E$. In particular if $E=TM$ then we set $\Sigma(M):=\Sigma(TM)$.\\
$\bullet$ More generally, let $\mZ(\R^n)\subset\End(\R^n)$ (resp. $ O(n)$) be some submanifold defined for any $n\in\mathbb N^*$. This allows us to define $\mZ(E)$ for any (Riemannian) vector bundle $E$ and we set in particular $\mal Z(M)=\mZ (TM)$ for any (Riemannian) manifold.
\\
$\bullet$ Moreover it will happen that we will write "$SO(E)$" without precising that the Riemannian vector bundle $E$ is  supposed to be oriented. We will consider that this hypothesis is implicit once we write this kind of symbol.
\\
$\bullet$ Let $(A,+,\times)$ be an associative $\mrm K$-algebra over the field $\mrm K$. Then for any  $a\in A$, we set $\com(a)=\{b\in A| ab=ba\}$ and $\mrm{Ant}(a)=\{b\in A| ab=-ba\}$.\\
$\bullet$ \wrt : with respect to.\medskip\\
\textbf{A list of notations and an index} are available at the end of the paper.\medskip\\
\textbf{Bibliographic remarks and a summary of our own contributions} are avalaible at the end of each section from  section~\ref{melliptic} to section~\ref{affineharmhom}.\medskip\\
\textbf{The numbering} of the equations is made by section: e.g. equation (2.13) is the 13th equation of section 2. Moreover the numbering of theorems, definitions and such, is made by subsection: e.g. Theorem~2.2.1 is the first theorem in subsection 2.2 (but it is not in subsubsection 2.2.1).

\subsection{Almost complex geometry}
Let $E$ be a real vector space endowed with a complex structure: $J\in\End(E)$, $J^2=-\Id$. Then we denote by $E^{1,0}$ and  $E^{0,1}$ respectively the eigenspaces of $J$ associated to the eigenvalues $\pm i$ respectively. Then we have the  following eigenspace decomposition
\begin{equation}\label{decEJ}
E^\C= E^{1,0} \oplus E^{0,1}
\end{equation}
and the following equalities
\begin{equation}\label{E(1,0)}
\begin{array}{c}
E^{1,0}=\ker(J-i\Id)=(J + i\Id)E^\C\\
E^{0,1}=\ker(J+i\Id)=(J - i\Id)E^\C
\end{array}
\end{equation}
so that  remarking that $(J \pm i\Id)iE=(\Id \mp iJ)E=(\Id \mp iJ)JE=(J \pm i\Id)E$, we can also write
\begin{equation}\label{E(0,1)}
\begin{array}{c}
E^{1,0}=(J + i\Id)E=(\Id - iJ)E= \{X -iJX, X\in E\} \\
E^{0,1}=(J - i\Id)E=(\Id + iJ)E=\{X +iJX, X\in E\}
\end{array}
\end{equation}
In the same way we denote by 
$$
(E^*)^\C=E_{1,0}^*\oplus E_{0,1}^*
$$
the decomposition induced on the dual $E^*$ by the complex structure $J^*\colon \eta\in E^*\to \eta J\in E^*$. Besides, given a vector $Z\in E^\C$, we denote by 
$$
Z= [Z]^{1,0} + [Z]^{0,1}
$$
its decomposition according to (\ref{decEJ}). Let us remark that 
$$
[Z]^{1,0}=(\Id -iJ)Z \quad \text{and} \quad [Z]^{0,1}=(\Id + iJ)Z.
$$
Moreover, given $\eta$ a $n$-form  on $E$, we denote by $\eta^{(p,q)}$ its component  in $\Lm^{p,q} E^*$ according to the decomposition 
$$
\Lm^n E^* =\oplus_{p+q=n}\Lm^{p,q} E^*,
$$ 
where $\Lm^{p,q} E^*=\left( \Lm^p E_{1,0}^*\right) \wedge \left( \Lm^q E_{0,1}^*\right)$. However for 1-forms, we will often prefer the notation $\eta=\eta' + \eta''$, where $\eta'$ and $\eta''$ denote respectively $\eta^{(1,0)}$ and $\eta^{(0,1)}$.\\ 
More generally, all what precedes holds naturally when $E$ is a real vector bundle over a manifold $M$, endowed with a complex structure $J$. \\
We will write 
$$
d=\partial + \bar\partial
$$
the decomposition of the exterior derivative of differential forms on an almost complex manifold $(M,J)$, according to the decomposition $TM^\C=T^{1,0}M\oplus T^{0,1}M$.\medskip\\
We will denote by $\mrm{Hol}((M,J^M),(N,J^N)):=\{f\in C^\infty(M,N)|\, df\circ J^M=J^N\circ df\}$ the set of holomorphic maps between two almost complex manifolds $(M,J^M)$ and $(N,J^N)$.\index{J holomorphic@$J$-holomorphic}\medskip\\
In this paper, we will use the following definitions.
\begin{defn}\label{f-structure}\index{f structure@$f$-structure}
Let $E$ be a real vector bundle. An $f$-structure in $E$ is an endomorphism $F\in \mal C(\End E)$ such that $F^3 + F=0$. An $f$-structure on a manifold $M$ is an $f$-structure in $TM$. A manifold $(M,F)$ endowed with an $f$-structure is called an $f$-manifold\index{f manifold@$f$-manifold}.
\end{defn}
An $f$-structure $F$ in a vector bundle $E$ is determined by its eigenspaces decomposition that we will denote by
$$
E^\C=E^+\oplus E^-\oplus E^0
$$
where $E^\pm=\ker(F \mp i\Id)$ and $E^0=\ker F$. In particular if $E=TM$, then we will set $T^iM= (TM)^i$, $\forall i\in \{0,\pm 1\}$.
\begin{defn}\index{f holomorphic@$f$-holomorphic}
Let $(M,F^M)$  and $(N,F^N)$ be  $f$-manifolds. Then a map $f\colon (M,F^M)\to (N,F^N)$ is said to be $f$\textbf{-holomorphic} if it satisfies
$$
df\circ F^M= F^N\circ df
$$
\end{defn}
\begin{defn}\index{horizontally holomorphic}
Let $(M,J^M)$ be an almost complex manifold and $N$ a manifold with a splliting $TN=\hor\oplus\ver$. Let us suppose that the subbundle $\hor$ is naturally endowed with a complex structure $J^\hor$. Then we will say that a map $f\colon (M,J^M)\to N$ is $\hor$\textbf{-holomorphic} if it satisfies the equation
$$
\left[ df\right]^\hor\circ J^M = J^\hor \circ [df]^\hor,
$$
where $[df]^\hor$ is the projection of $df$ on $\hor$ along $\ver$. 
Moreover, if for some reason, $\hor$ inherits the name of horizontal subbundle, then we will say that $f$ \textbf{is horizontally holomorphic}.
\end{defn}
\begin{rmk} \em
Let us remark that an $f$-structure in a manifold $N$ is equivalent to a splitting $TN=\ver\oplus \hor$ together with a complex structure $J^\hor$ on $\hor$.
\end{rmk}
\begin{defn}
An affine manifold $(N,\nabla)$ is a manifold endowed with a linear connection.\\
An almost complex affine manifold $(N,J,\nabla)$ is an almost complex manifold endowed with an almost complex connection: $\nabla J=0$.
\end{defn}
%In other words, $f\colon (M,J^M)\to (N,F)$ is $f$-holomorphic, when  $N$ is endowed with the $f$-structure $F=J^\hor\oplus 0_\ver$.
%This situation occurs for example if $N$ is endowed with an $f$-structure $F$ which leaves invariant $\hor$ and $F_{|%\hor}$ is a complex structure (i.e. $T^0N\cap \hor =\{0\}$).

%
%%%%%%%%%%%%%%%%%%%%%%%%%%%%%%%%%%%%%%%%%%%%%%%%%%%%%%%%%%%%%%%%%%%%%%%%%%

%   Invariant connections on reductive homogeneous spaces

%%%%%%%%%%%%%%%%%%%%%%%%%%%%%%%%%%%%%%%%%%%%%%%%%%%%%%%%%%%%%%%%%%%%%%%%%%
  
\section{Invariant connections on reductive homogeneous spaces}\label{1}
The references for this section where we recall some results that we will need in this paper, are \cite{KN}, \cite{pham}, \cite{BuRaw}, and to a lesser extent \cite{agricola} and \cite{higaki}.
\subsection{Linear isotropy representation}\label{linearisotropy}
Let $M=G/H$ be a homogeneous space with $G$ a real Lie group and $H$ a closed subgroup of $G$. $G$ acts transitively on $M$ in a natural manner which defines a natural representation: $\phi\colon g\in G\mapsto(\phi_g\colon x\in M\mapsto g.x)\in\mrm{Diff}(M)$. Then $\ker\phi$ is the maximal normal subgroup of $G$ contained in $H$. Further, let us consider the linear isotropy representation:
$$
\rho\colon h\in H\mapsto d\phi_h(x_0)\in GL(T_{x_0}M)
$$
where $x_0=1.H$ is the reference point in $M$. Then we have $\ker\rho\supset\ker\phi$. Moreover the linear isotropy representation is faithful (i.e. $\rho$ is injective) \iif $G$ acts freely on the bundle of linear frames $L(M)$.\\
We can always suppose without loss of generality that the action of $G$ on $M$ is effective (i.e. $\ker\phi=\{1\}$) but it does not imply in general that the linear isotropy representation is faithful. However if there exists on $M$ a $G$-invariant linear connection, then the linear isotropy representation is faithful provided that $G$ acts effectively on $M$. (Indeed, given a manifold $M$ with a linear connection, and $x\in M$, an affine transformation $f$ of $M$ is determined by $(f(x),df(x))$, i.e. $f$ is the identity \iif it leaves one linear frame fixed).
\subsection{Reductive homogeneous space}\label{reductivehomspaces}
Let us suppose now that $G/H$ is reductive, i.e. there exists a decomposition $\g=\hk\oplus\mk$ such that $\mk$ is $\Ad H$-invariant: $\forall h\in H, \Ad h(\mk)=\mk$. Then the surjective map $\xi\in\g\mapsto \xi.x_0\in T_{x_0}M$ has $\hk$ as kernel and  so its restriction to $\mk$ is an isomorphism $\mk\cong T_{x_0}M$. This provides an isomorphism of the associated bundle $G\times_H\mk$ with $TM$ by:
\begin{equation}\label{can}
[g,\xi]\mapsto g.(\xi.x_0)=\Ad g(\xi).x
\end{equation}
where $x=\pi(g)=g.x_0$.\\
Moreover, we have a natural inclusion $G\times_H\mk\mapsto G\times_H\g$ and the associated bundle $G\times_H\g$ is canonically identified with the trivial bundle $M\times\g$ via
\begin{equation}\label{io}
[g,\xi]\mapsto (\pi(g),\Ad g(\xi)).
\end{equation}
Thus we have an identification of $TM$ with a subbundle $[\mk]$ of $M\times\g$, which we may view as a $\g$-valued 1-form $\beta$ on $M$ given by: 
$$
\beta_x(\xi.x)=\Ad g[\Ad g^{-1}(\xi)]_\mk,
$$
where $\pi(g)=x,\xi\in\g$ and $[\ ]_\mk$ is the projection on $\mk$ along $\hk$. Equivalently, for all
$X\in T_xM$, $\beta(X)$ is the unique element $\xi\in [\mk]_x$ ($=\Ad g(\mk)$, with $\pi(g)=x$) such that
$X=\xi.x$, in other words $\beta(X)$ is caracterized by 
$$ 
\beta(X)\in [\mk]_x\subset\g \quad \text{and} \quad X=\beta(X).x\ .
$$
In fact, $\beta$ is nothing but the projection on
$M$ of the $H$-equivariant 1-form, $\theta_\mk$, on $G$, i.e.
$\theta_\mk$ is the $H$-equivariant lift of $\beta$. Here,  $\theta_\mk$ is defined as the $\mk$-component of the left invariant Maurer-Cartan form $\theta$ of $G$.  This can be written as follows
\begin{equation}\label{beta}
(\pi^*\beta)_g=\Ad g(\theta_\mk)\  \forall g\in G 
\end{equation}
with $\theta_g(\xi_g)=g^{-1}.\xi_g$ for all $g\in G$, $\xi_g\in T_gG$.\\[1mm]
\textbf{Notation} For any $\Ad H$-invariant subspace $\mak l\subset \mk$, we will denote by $[\mak l]$ the subundle of $[\mk]\subset M\times \g$ defined by $[\mak l]_{g.x_0}=\Ad g(\mak l)$.
\subsection{The (canonical) invariant connection}\label{invariant}\index{canonical!connection, $G$-invariant|(}
On a reductive homogeneous space $M=G/H$, the $\Ad(H)$-invariant summand $\mk$ provides by left translation in
$G$, a $G$-invariant distribution $\mal H(\mk)$, given by $\mal
H(\mk)_g=g.\mk$ which is horizontal for $\pi\colon G\to M$ and right
$H$-invariant and thus defines a $G$-invariant connection in the
principal bundle $\pi\colon G\to M$. In fact this procedure defines
a bijective correspondance between  reductive summands $\mk$ and
$G$-invariant connections in $\pi\colon G\to M$ (see \cite{KN},
chap. 2, Th 11.1). Then the corresponding $\hk$-valued connection 1-form $\omega$ on $G$ (of this $G$-invariant connection) is the $\hk$-component of the left invariant Maurer-Cartan form of $G$:
$$
\omega=\theta_\hk.
$$ 
\subsection{Associated covariant derivative}\label{assocovarder}
The connection 1-form $\omega$ induces a covariant derivative in the associated bundle ${G\times_H\mk}\cong TM$ and thus a $G$-invariant covariant derivative $\nabla^0$ in the tangent bundle $TM$. In particular, we can conclude according to section~\ref{linearisotropy} that if $G/H$ is reductive then the linear isotropy representation is faithful (provided that $G$ acts effectively) or equivalently that $\ker\Adm=\ker\rho=\ker\phi$. We will suppose in the following that,  without  explicit  reference to the contrary, the  action of $G$ is effective and (thus) the linear isotropy representation is faithful.\\
One can compute explicitly $\nabla^0$.
\begin{lemma}\label{br}\emph{\cite{BuRaw}}
$$
\beta(\nabla_X^0 Y)=X.\beta(Y)-[\beta(X),\beta(Y)],\quad  X,Y\in \Gamma(TM).
$$
\end{lemma}
Let us write (locally) $\beta(X)=\Ad U(X_\mk)$, $\beta(Y)=\Ad U(Y_\mk)$ where $U$ is a (local) section of $\pi$ and $X_\mk,Y_\mk\in C^\infty(M,\mk)$ then we have (using the previous lemma)
\begin{eqnarray*}
\beta(\nabla_X^0 Y) & = & \Ad U \left(dY_\mk(X) + [\alpha(X), Y_\mk] - [X_\mk,Y_\mk]\right)\\
                    & = & \Ad U \left(dY_\mk(X) + [\alpha_\hk(X), Y_\mk] + [\alpha_\mk(X) - X_\mk,Y_\mk]\right)
\end{eqnarray*}
where $\alpha=U^{-1}.dU$. Besides since $U$ is a section of $\pi$ ($\pi\circ U=\Id$), then pulling back (\ref{beta}) by $U$, we obtain $\beta=\Ad U(\alpha_\mk)$ and  then $\alpha_\mk(X)=X_\mk$, so that
\begin{equation}\label{nabla-0}
\beta(\nabla_X^0 Y)=\Ad U \left(dY_\mk(X) + [\alpha_\hk(X), Y_\mk]\right)
\end{equation}
\begin{rmk}\label{lift}
\em We could also say that $X_\mk,Y_\mk$ are respectively the pullback by $U$ of the $H$-equivariant lifts $\tl X,\tl Y$ of $X,Y$ (given by $\beta(X_{\pi(g)})=\Ad g(\tl X(g))$).\\
Then $\nabla^0_{X} Y$ lifts as the $\mk$-valued $H$-equivariant map on $G$:
$$
\widetilde{\nabla^0_{ X} Y}= d\tl Y(\tl X) + [\theta_\hk(\tl X),\tl Y]
$$
and then taking the $U$-pullback we obtain the previous result (without using lemma~\ref{br}).\\
Moreover, we can express $\nabla^0$ in terms of the flat differentiation in the trivial bundle $M\times\g$ ($\supset [\mk]$). Let us differentiate the equation $Y=\Ad U(Y_\mk)$ (we use the identification $TM=[\mk]\subset M\times\g$)
$$
dY=\Ad U \left(dY_\mk + [\alpha,Y_\mk]\right) = \Ad U \left(dY_\mk +[\alpha_\hk,Y_\mk]\right) + \Ad U\left(([\alpha_\mk,Y_\mk]\right))=\nabla^0 Y + [\beta, Y].
$$
Finally, we obtain
\begin{equation}\label{flat-diff}
dY=\nabla^0 Y + [\beta,Y]
\end{equation}
and we recover lemma~\ref{br}.
\end{rmk}
\index{canonical!connection, $G$-invariant|)}
\subsection{$G$-invariant linear connections in terms of equivariant bilinear maps}\label{g-invariantaffin}
Now let us recall the following results about invariant connections on reductive homogeneous spaces.
\begin{thm}\label{correspondance}\emph{\cite{KN}}
Let $\pi_P\colon P\to M$ be a $K$-principal bundle over the reductive homogeneous space $M=G/H$ and suppose that $G$ acts on $P$ as a group of automorphisms and let $u_0\in P$ be a fixed point in the fibre of $x_0 \in M$ ($\pi_P(u_0)=x_0$). There is a bijective correspondance between the set of $G$-invariant connections $\omega$ in $P$ and the set of linear maps $\Lambda_\mk\colon\mk\to \mak k$ such that
\begin{equation}\label{lambda-m}
\Lm_\mk(hXh^{-1})=\lm(h)\Lm_\mk(X)\lm(h)^{-1}\quad \text{for } X\in\mk \text{ and } h\in H
\end{equation}
where $\lm\colon H\to K$ is the morphism defined by $hu_0 = u_0\lm(h)$ ($H$ stabilizes the fibre $P_{x_0}=u_0.K$). The correspondance is given by
\begin{equation}\label{omega-u}
\Lm(X)=\omega_{u_0}(\tl X),\quad \forall X\in \g
\end{equation}
where $\tl X$ is the vector field on $P$ induced by $X$ (i.e. $\forall u\in P$, $\tl X(u)=\frac{d}{dt}_{|t=0}\exp(tX).u$) and $\Lm\colon\g\to\mak k$ is defined by $\Lm_{|\mk}=\Lm_\mk$ and $\Lm_{|\hk}=\lm$ (hence completely determined by $\Lm_\mk$).
\end{thm}
\begin{cory}\label{correspondance'}
In the previous theorem, let us suppose that $P$ is a $K$-structure on $M=G/H$, i.e. $P$ is a subbundle of the bundle $L(M)$ of linear frames on $M$ with structure group $K\subset GL(n,\R)=GL(\mk)$ (we identify as usual $\mk$ with $T_{x_0}M$ by $\xi\mapsto \xi.x_0$, and $T_{x_0}M$ to $\R^n$ via the linear frame $u_0\in P\subset L(M)$). Then in terms of  the $G$-invariant covariant derivative $\nabla$ corresponding to $\omega$, the $G$-invariant linear connection in $P$,  the previous bijective correspondance may be given by
$$
\Lm(X)(Y)=\nabla_{\tl X} \tl Y
$$
where $\tl X,\tl Y$ are any (local) left $G$-invariant  vector fields extending $X,Y$ i.e. there exists a local section, $g\colon U\subset M\to G$,  of $\pi\colon G\to M$, such that $\tl X_x=\Ad g_{(x)}(X).p$, for all $x\in M$.
\end{cory}
\begin{rmk}\label{rmk-can-con-restr}\em
In theorem~\ref{correspondance}, the $G$-invariant connection in $P$ defined by $\Lm_{\mk}=0$ is called the canonical connection (with respect to the decomposition $\g=\hk + \mk$). If we set $P(M,K)=G(G/H,H)$ with group of automorphisms $G$, the $G$-invariant connection defined by the horizontal distribution $\mal H(\mk)$ is the canonical connection.\\
Now, let $P$ be a $G$-invariant $K$-structure on $ M=G/H$ as in corollary~\ref{correspondance'}. Let $P'$ be a $G$-invariant subbundle of $P$ with structure group $K'\subset K$, then the canonical connection in $P'$ defined by $\Lm_\mk=0$ is (the restriction of ) the canonical connection in $P$ which is itself the restriction to $P$ of the canonical connection in $L(M)$. In particular, if we set $P'=G.u_0$, this is a subbundle of $P$ with group $H$, which is isomorphic to the bundle $G(G/H,H)$. Then the canonical linear connection in $P'$  corresponds to the invariant connection in $G(G/H,H)$ defined by the distribution $\mal H(\mk)$.
\end{rmk}
\begin{thm}
Let $P\subset L(M)$ be a $K$-structure on $M=G/H$. Then the canonical linear connection ($\Lm_\mk=0$) in $P$  defines the covariant derivative $\nabla^0$ in $TM$ (obtained from $\mal H(\mk)$ in the associated bundle $G\times_H\mk\cong TM$). Moreover there is a bijective correpondence between the set of  of $G$-invariant linear connections $\nabla$, on $M$,  determined by a connection in $P$, and the set of linear maps $\Lm_\mk\colon\mk\to \mak k\subset \mak{gl}(\mk)$ such that  
\begin{equation}\label{equi}
\Lm_\mk(hXh^{-1})=\Adm( h) \Lm_\mk(X) \Adm (h)^{-1}\quad \forall X\in \mk,\forall h\in H,
\end{equation}
given by
$$
\nabla=\nabla^0 + \Bar\Lm_\mk
$$
i.e. $\nabla_X Y = \nabla_X^0 Y + \Bar\Lm_\mk(X)Y$ for any vector fields $X,Y$ on $M$, where  with the help of (\ref{equi}) we extended  the $\Ad(H)$-equivariant map $\Lm_\mk\colon \mk\times\mk \to\mk$ to the bundle $G\times_H\mk=TM$ to obtain a map $\Bar\Lm_\mk \colon TM\times TM\to TM$.
\end{thm}
\begin{exam}
Let us suppose that $M$ is Riemannian (i.e. $\Adm H$ is compact and $\mk$ is endowed with an $\Ad H$ invariant inner product which defines a $G$-invariant metric on $M$)  and let us take $P=O(M)$ the bundle of orthonormal frames on $M$, the previous  correspondance is between the set  of $G$-invariant metric linear connections and the set of $\Ad(H)$-equivariant linear maps $\Lm_\mk\colon\mk \to \so(\mk)$.\\
In particular the canonical connection $\nabla^0$ is metric (for any $G$-invariant metric on $M$).
\end{exam}
\begin{thm}\label{geodesic}
\begin{description}
\item[$\bullet$]
$G$-invariant tensors on the reductive homogeneous space $M=G/H$ (or more generally $G$-invariant sections of associated bundles) are parallel with respect to the canonical connection.
\item[$\bullet$]
The canonical connection is complete (the geodesics are exactly the curves $x_t=\exp(tX).x_0$, for $X\in \mk$).
\item[$\bullet$]
Let $P$ be a $G$-invariant $K$-structure on $M=G/H$, then the $G$-invariant connection defined by $\Lm\colon\mk\to \mak k$ has the same geodesics as  the canonical connection \iif 
$$
\Lm_\mk(X)X=0,\quad\forall X\in \mk
$$
\end{description}
\end{thm}
\begin{thm}\label{T-R}
The torsion tensor $T$ and the curvature tensor $R$ of the $G$-invariant connection corresponding to $\Lm_\mk$ is given  at the origin point $x_0$ as follows:
\begin{enumerate}
\item
$T(X,Y)= \Lm_\mk(X)Y -\Lm_\mk(Y)X - [X,Y]_\mk,$
\item $R(X,Y)=[\Lm_\mk(X),\Lm_\mk(Y)]-\Lm_\mk([X,Y]_\mk)-\adm([X,Y]_\hk)$,
\end{enumerate}
for  $X,Y\in \mk$.\\
In particular, for the canonical connection we have $T(X,Y)=-[X,Y]_\mk$ and $R(X,Y)=-\adm([X,Y]_\hk)$, for $X,Y\in \mk$; moreover we have $\nabla T =0$, $\nabla R=0$.
\end{thm}

\subsection{A Family of connections on the reductive space $M$}\label{family}
We take in what precede (i.e. in section~\ref{g-invariantaffin}) $P=L(M)$. Then let us consider the one parameter family of connections $\nabla^t$, $0\leq t\leq 1$ defined by
$$
\Lm_\mk^t(X)Y=t[X,Y]_\mk,\quad 0\leq t\leq 1.
$$
For $t=0$, we obtain the canonical connnection $\nabla^0$. Since for any $t\in [0,1]$,  $\Lm_\mk^t(X)X=0$, $\forall X\in \mk$, $\nabla^t$ has the same geodesics as $\nabla^0$ and in particular is complete. The torsion tensor is given (at $x_0$) by 
\begin{equation}\label{Tt}
T^t(X,Y)=(2t-1)[X,Y]_\mk .
\end{equation}
In particular $\nabla^{\frac{1}{2}}$ is the unique torsion free $G$-invariant linear connection having the same geodesics as the canonical connection (according to theorems~\ref{geodesic} and \ref{T-R}).\\[1mm]
Assume now that $M$ is Riemannian, and  let us take $P=O(M)$. Then $\nabla^t$ is metric \iif $\Lm_\mk^t$ takes values in $\mak k=\so(\mk)$ which is equivalent (for $t\neq 0$) to say that $M$ is naturally reductive (which means by definition that $\forall X\in \mk$, $[X,\cdot]_\mk$ is skew symmetric). Now (still in the Riemannian case) let us construct a family of linear connections, $\overset{\mrm{met}}{\nabla^t}$, $0\leq t\leq 1$, which are always metric: 
$$
\overset{\mrm{met}}{\nabla^t}=\nabla^0 + t\left([\ ,\ ]_{[\mk]} +  \mU^M\right)
$$
where $ \mU^M\colon TM\oplus TM\to TM$ is the "natural reductivity term" which is the symmetric bilinear map defined by\footnote{$\mU^M$ is the  $G$-invariant extension of $\mU^\mk\colon \mk\oplus\mk \to \mk$, its restriction to $\mk\oplus\mk$.} 
\begin{equation}\label{defofU}
\langle  \mU^M(X,Y),Z\rangle=\langle[Z,X]_{[\mk]},Y\rangle + \langle X,[Z,Y]_{[\mk]}\rangle
\end{equation}
for all $X,Y,Z\in [\mk]$.
Since $\mU^M$ is symmetric, the torsion of $\overset{\mrm{met}}{\nabla^t}$ is once again given by 
$$
T^t(X,Y)=(2t-1)[X,Y]_{[\mk]}
$$
and thus $\overset{\mrm{met}}{\nabla^\frac{1}{2}}$ is torsion free and metric and we recover that $\overset{\mrm{met}}{\nabla^\frac{1}{2}}$ is the Levi-Civita connection
$$
\overset{\mrm{met}}{\nabla^\frac{1}{2}}=\nabla^{\mrm{L.C.}}. 
$$
Obviously if $M$ is naturally reductive then $\overset{\mrm{met}}{\nabla^t}=\nabla^t$, $\forall t\in [0,1]$. Moreover if $M$ is (locally) symmetric, i.e. $[\mk,\mk]\subset\hk$, then all these connections coincide and are equal to the Levi-Civita connection: $\overset{\mrm{met}}{\nabla^t}=\nabla^t=\nabla^0=\nabla^{\mrm{L.C.}}. $
\begin{rmk}\label{nabla1}\index{anticanonical connection}
\em $\nabla^1$ is interesting since it is nothing but the flat differentiation in the trivial bundle $M\times \g$ followed by the projection onto $[\mk]$ (along $[\hk]$) (see remark~\ref{lift}). So this connection is very natural and following \cite{agricola}, we will call it \textbf{the anticanonical connection}.
\end{rmk}
\subsection{Differentiation in $\mrm{End}(T(G/H))$}\label{End}
According to section~\ref{reductivehomspaces}, we have 
$$ 
\mrm{End}(T(G/H))=G\times_H\End(\mk)\subset (G/H)\times\End(\g),
$$
the previous inclusion being given by $[g,A]\mapsto(\pi(g),\Ad g\circ A\circ\Ad g^{-1})$ and we embedd $\End(\mk)$ in $\g$ by extending  any endomorphism of $\mk$ to the corresponding endomorphism of $\g$ which vanishes on $\hk$. In other words $\End(T(G/H))$ can be identified to the subbundle $[\End(\mk)]$ of the trivial bundle $(G/H)\times\End(\g)$, with fibers $[\End(\mk)]_{g.x_0}=\End(\Ad g(\mk))=\Ad g(\End(\mk))\Ad g^{-1}=\Ad g(\End(\mk)\oplus\{0\})\Ad g^{-1}$.\\
Now, let us compute in terms of  the Lie algebra setting, the derivative of the inclusion map $\mak I \colon \End(T(G/H))\to M\times\End(\g)$ or more concretely the flat derivative in $M\times\End(\g)$ of any section of $\End(T(G/H))$. To do that, we compute the derivative of 
$$
\tl{\mak I} \colon (g,A_\mk)\in G\times \End(\mk)\longmapsto(g.x_0,\, \Ad g\circ A_\mk\circ \Ad g^{-1})\in M\times \End(\g).
$$
We obtain 
$$
d\tl{\mak I}(g,A_\mk)=\left( \Ad g(\theta_\mk)\, .\,\pi(g),\,\Ad g\left(d A_\mk +[\ad\theta, A_\mk]\right)\Ad g^{-1}\right).
$$
Next, let us decompose the endomorphisms in $\g$ into blocs according to the vector space decomposition $\g=\hk\oplus\mk$:
\begin{equation}\label{dec-endog}
\End(\g)=\begin{pmatrix}
\End(\hk) & \End(\mk,\hk)\\
\End(\hk,\mk) & \End(\mk)
\end{pmatrix}.
\end{equation}
By regrouping terms, we obtain the following splitting
$$
\End(\g)=\End(\mk)\oplus\left( \End(\mk,\hk)\oplus \End(\hk,\mk)\oplus \End(\hk)\right), 
$$
which applied to $d\tl{\mak I}(g,A_\mk)$, gives us the decomposition
\begin{eqnarray}\label{dec-endo}
d\tl{\mak I}(g,A_\mk) & = & \left(0,\,\Ad g\left(d A_\mk +[\adm\theta_\hk, A_\mk] + \left[[\adm\theta_\mk]_\mk, A_\mk\right]\right)\Ad g^{-1}\right)\\
& &  + \left( \Ad g(\theta_\mk)\, .\,\pi(g),\,\Ad g\left([\adm\theta_\mk ]_{\hk}\circ A_\mk -A_\mk\circ\ad_\hk\theta_\mk\right)\Ad g^{-1}\right).\nonumber
\end{eqnarray}
The first term is in  the vertical space $\mal V_{\tl{\mak I}(g,A_\mk)}=\Ad g(\End(\mk)) \Ad g^{-1}=\End(T_{\pi(g)}M)$ and the previous decomposition (\ref{dec-endo}) provides us with a splitting $T\End(M)=\mal V\oplus\mal H= \pi_M^*(\End(M)) \oplus \mal H$, i.e. a connection on $\End(M)$. Let us determine this connection: we see that the projection on the vertical space (along the horizontal space) corresponds to the projection on $[\End(\mk)]$ following (\ref{dec-endog}) so that according to remark~\ref{nabla1}, we can conclude that the horizontal distribution $\mal H$ defines the connection $\nabla^1$ on $\End(TM)= TM^*\otimes TM$.
\begin{rmk}
\em 
\textbf{a)} We can recover this  previous fact directly from the first term of (\ref{dec-endo}) and the definition of $\nabla^1$. Indeed, first recall that given two linear connections $\nabla$, $\nabla'$ on $M$, we can write $\nabla'=\nabla + F$, where $F$ is a section of $TM^*\otimes\End(TM)$, and then for any section $A$ in $\End(TM)$,
$$
\nabla' A= \nabla A + [F,A].
$$
Besides, $\nabla^1=\nabla^0 + [\ ,\ ]_{[\mk]}$, and moreover, if we write (locally) $A=(\pi(U),\, \Ad U\circ A_\mk\circ\Ad U^{-1})$ where $U$ is a local section of $\pi$ and $A_\mk\in C^\infty(M,\End(\mk))$, then according to (\ref{nabla-0}),
\begin{equation}\label{nabla0-endo}
\nabla^0 A=\Ad U\left(d A_\mk + [\ad_\mk\alpha_\hk,A_\mk]\right),
\end{equation}
so that we conclude that
$$
\nabla^1 A=\Ad U\left(d A_\mk + [\ad_\mk\alpha_\hk,A_\mk]+ \left[[\adm\alpha_\mk]_\mk, A_\mk\right] \right)\Ad U^{-1}
$$
which is the (pullback of) the first term of (\ref{dec-endo}).\\[1.5mm] 
\textbf{b)} Furthermore, if $G/H$ is (locally) symmetric (i.e. $[\mk,\mk]\subset\hk$), then $\nabla^{L.C.}=\nabla^0=\nabla^1$ and in particular 
\begin{equation}\label{L.C.sym}
\nabla^{L.C.}A=\Ad U\left(d A_\mk + [\ad_\mk\alpha_\hk,A_\mk]\right).
\end{equation}
\end{rmk}

%
%
%%%%%%%%%%%%%%%%%%%%%%%%%%%%%%%%%%%%%%%%%%%%%%%%%%%%%%%%%%%%%%%%%%%%%%%%%%%%%%%%%%%%%

%  $m$-th elliptic integrable system associated to a $k'$-symmetric space

%%%%%%%%%%%%%%%%%%%%%%%%%%%%%%%%%%%%%%%%%%%%%%%%%%%%%%%%%%%%%%%%%%%%%%%%%%%%%%%%%%%%%

\section{$m$-th elliptic integrable system associated to a $k'$-symmetric space}\label{melliptic}
\index{m elliptic integrable system@$m$-th elliptic integrable system|(}
\subsubsection{Definition of $G^\tau$ (even when $\tau$ does not integrate in $G$)}\label{def-g-tau}
Here, we will extend the notion of subgroup fixed by an automorphism of Lie group to the situation where only a Lie algebra automorphism is provided. Indeed, let $\tau\colon G \to G$ be a Lie group automorphism, then usually one can define $G^\tau=\{g\in G|\, \tau(g)=g\}$ the subgroup fixed by $\tau$. Now, we want to extend this definition to the situation where we only have a Lie algebra automorphism, and so that the two definitions coincide when the Lie algebra automorphism integrates in $G$.\\ 

\noindent
Let $\g$ be a real Lie algebra and $\tau\colon \g \to \g$ be an automorphism. Then let us denote by
\begin{equation}\label{g0def}
\g_0=\g^\tau:=\{\xi\in\g |\,\tau(\xi)=\xi\}
\end{equation}
the subalgebra of $\g$ fixed by $\tau$. Let us  assume that $\tau$ defines in $\g$ a $\tau$-invariant reductive decomposition 
$$
\g=\g_0\oplus \nk,\quad [\g_0,\nk]\subset \nk,\quad \tau(\nk)=\nk.
$$
Moreover we suppose that we have $\nk=\im(\Id-\tau)$, i.e. that the decomposition $\g=\g_0\oplus\nk$ coincides with the Fitting decomposition of $\Id -\tau$ (remark that this decomposition is then automatically reductive since $\ad\xi\circ(\Id-\tau)=(\Id-\tau)\circ\ad \xi$, $\forall \xi\in\g_0$).\\
Without loss of generality, we assume that the center of $\g$ is trivial. Moreover, we assume also, without loss of generality, that $\g_0$ does not contain non-trivial ideal of $\g$, i.e. that $\ad_\nk \colon \g_0\to \mak{gl}(\nk)$ is injective (the kernel is a $\tau$-invariant ideal of $\g$ that we factor out). We then have
\begin{equation}\label{g0cara}
\g_0=\{\xi\in \g|\,\tau\circ\ad\xi\circ\tau^{-1}=\ad\xi\}
\end{equation}
Let $G$ be a connected Lie group with Lie algebra $\g$. Then  let us consider the subgroup
$$
G_0=\{g\in G|\, \tau\circ\Ad g \circ \tau^{-1}=\Ad g\}.
$$
Then $G_0$ is a closed subgroup of $G$ and  $\lie G_0= \g_0$. \\
Moreover, without loss of generality, we will suppose  that $G_0$ does not contain non-trivial normal subgroup of $G$ (by factoring out, if needed, by some discrete subgroup of $G$), i.e. that $\Ad_\nk\colon G_0\to GL(\nk)$ is injective (see section~\ref{1}). Now, we want to prove that if $\tau$ integrates in $G$, then we have $G_0=G^\tau$, where $G^\tau$ is the subgroup fixed by $\tau\colon G\to G$. First, we have $\forall g\in G^\tau$, 
$\tau\circ\Ad g \circ \tau^{-1}= \Ad\tau(g)=\Ad g$, thus $G^\tau\subset G_0$. Conversely, we have $\forall g\in G_0$, $\Ad g(\nk)=\nk$  and $\Ad g=\tau\circ\Ad g \circ \tau^{-1}= \Ad\tau(g)$ so that $\Ad_\nk g=\Ad_\nk \tau(g)$ and thus $g=\tau(g)$ since $\Ad_\nk\colon G_0\to GL(\nk)$ is injective. We have proved $G^\tau=G_0$. This allows us to make the following:
% and 
\begin{defn}
Let $\g$ be a real Lie algebra and $\tau\colon \g \to \g$ be an automorphism, and $G$  a connected Lie group with Lie algebra $\g$. Let us  assume that $\tau$ defines in $\g$ a $\tau$-invariant reductive decomposition: $\g=\g_0\oplus \nk$ with $\nk=\im(\Id-\tau)$. Then we will set  
$$
G^\tau:= \{g\in G|\, \tau\circ\Ad g \circ \tau^{-1}=\Ad g\}.
$$
\end{defn}
Let us conclude this subsection by some notations:\\[1.5mm]
\textbf{Notations and conventions}
In all the paper, when a real Lie algebra $\g$ and an automorphism $\tau$ will be given, then we will suppose without loss of generality that the center of $\g$ is trivial,  $\g_0$ will denote the Lie subalgebra defined by (\ref{g0def}), $G$ will denote a connected Lie group with Lie algebra $\g$ and $G_0\subset G$ a closed subgroup such that $(G^\tau)^0\subset G_0\subset G^\tau$ (which implies that its Lie algebra is $\g_0$).\\
Moreover, without loss of generality, we will always suppose that $\g_0$ does not contain non-trivial ideal of $\g$ - we will then say that $(\g,\g_0)$ is effective - and also suppose that 
$G^\tau$ does not contain non-trivial normal subgroup of $G$ (by factoring out, if needed, by some discrete subgroup of $G$). Consequently, when $\tau$ can be integrated in $G$, then $G^\tau$ will coincide with the subgroup of $G$ fixed by $\tau\colon G\to G$.
\begin{rmk}\label{rmk-G0&G'}\em
$\bullet$ Let $\g$ be a Lie algebra and $\tau\colon \g \to \g$ be an automorphism. Let us consider $G'=\Ad G =:\Int(\g)$  the adjoint group of $\g$ (which does not depend on the choice of the connected group $G$), and let $C=\ker \Ad$ be the center of $G$. Then we can identify the morphism $\Ad$ with the covering $\pi\colon G\to G/C$ and $G'$ to $G/C$. Besides $\tau$ always integrates in $G'$ into $\tau'$ defined by $\tau'=\Int \tau\colon \Ad g\in G'\mapsto \tau\circ\Ad g\circ\tau^{-1}$ and we have $\tau'\circ\pi=\pi\circ\tau$ (See \cite[p. 127]{Hel}). Then, we see that $G_0=\Ad^{-1}({G'}^{\tau'})=\pi^{-1}({G'}^{\tau'})$.\\
$\bullet$ Now, let us come back to the hypothesis of trivialness of the center of $\g$. Let $\g$ be a real Lie algebra with center $\mak c$ and $\tau\colon\g\to\g$ an automorphism. Then we have $\tau(\mak c)=\mak c$. Therefore, if $\tau$ is of order $k'$, for any $m\in \mathbb N^*$, the $m$-th elliptic integrable system associated to $\tau$ splits into two independent systems: a general one on $\g'=\g/\mak c$ and a trivial one on $\mak c =\R^n$ : $d\alpha_\lm=0$, $\forall \lm\in\C^*$. Then the solutions of this trivial system are the maps $f_\lm\colon L\to \mak c$, $f_\lm=\sum_{j=-m}^{m} \lm^j f_j$, $\forall \lm\in\C^*$, such that $f_{-j}=\ovr{f_j}$, $j\geq  0$, and $f_{-j}$ is an holomorphic map into the eigenspace $\mak c^\C\cap\g_{-j}^\C$. (Of course  $\alpha_\lm=df_\lm$.) Therefore, the hypothesis of trivialness of the center can be done without loss of generality. 
\end{rmk}
\subsection{Finite order Lie algebra automorphisms}\label{finitorderauto}
\index{eigenspace decomposition of the automorphism $\tau$}\index{k symmetric space@$k'$-symmetric space|(}
Let $\g$ be a real Lie algebra and $\tau\colon \g \to \g$ be an  automorphism of order $k'$. Let $\omega_{k'}$ be a  $k'$-th primitive root of unity. Then we have the following eigenspace decomposition:
$$
\g^\C=\bigoplus_{j\in\Z/k'\Z}\g_j^\C,\qquad [\g_j^\C,\g_l^\C]\subset\g_{j+l}^\C
$$
where $\g_j^\C$ is the $\omega_{k'}^j$-eigenspace of $\tau$.
\\[1mm]
We then have to distinguish two cases.
\subsubsection{The even case: $k'=2k$}\label{evendetercase}\index{even case|(}
Then we have $\g_0^\C=(\g_0)^\C$. Moreover let us remark that 
\begin{equation}\label{gj-bar}
\overline{\g_j^\C} =\g_{-j}^\C,\quad \forall j\in\Z/k'\Z.
\end{equation}
Therefore $\overline{\g_k^\C}=\g_{-k}^\C=\g_k^\C$ so that we can set $\g_k^\C=(\g_k)^\C$ with
$$
\g_k=\{\xi\in\g|\,\tau(\xi)=-\xi\}.
$$
Moreover, owing to (\ref{gj-bar}), we can define $\mk_j$ as the unique real subspace in $\g$ such that its complexified is given by
$$
 \mk_j^\C=\g_j^\C\oplus\g_{-j}^\C\, \text{ for } j\neq 0,k, 
$$
and $\nk$ as the unique real subspace such that
$$ 
\nk^\C=\bigoplus_{j\in \Z_k'\setminus \{0\} }\g_j^\C ,
$$
that is $\nk=(\oplus_{j=1}^{k-1}\mk_j)\oplus \g_k$. In particular $\tau$ defines a $\tau$-invariant reductive decomposition $\g=\g_0\oplus\nk$.\\
 Hence the eigenspace decomposition is written:
$$
\g^\C=\left(\g_{-(k-1)}^\C\oplus\ldots\oplus\g_{-1}^\C\right)\oplus\g_0^\C\oplus\left(\g_1^\C\oplus\ldots\oplus \g_{k-1}^\C\right)\oplus\g_k^\C
$$
so that by grouping
\begin{eqnarray*}
\g^\C & = & \g_0^\C\oplus\g_k^\C\oplus\left[\oplus_{j=1}^{k-1}\mk_j^\C\right]\\
      & = & \hk^\C\oplus\mk^\C
\end{eqnarray*}
where $\hk=\g_0\oplus\g_k$ and $\mk=\oplus_{j=1}^{k-1}\mk_j$. Considering the automorphism $\sigma=\tau^2$, we have $\hk=\g^\sigma$ and $\g=\hk\oplus\mk$ is the reductive decomposition defined by the order $k$ automorphism $\sigma$. Without loss of generality, and according to our convention applied to  $\g$ and $\sigma$,   we will suppose in the following that $(\g,\hk)$ is effective i.e. $\hk$ does not contain non trivial ideal of $\g$. This implies in particular that $(\g,\g_0)$ is also effective.\medskip\\
 %\linebreak 
Now let us integrate our setting: let $G$ be a Lie group with Lie algebra $\g$ and we choose $G_0$ such that $(G^\tau)^0\subset G_0\subset G^\tau$. Then $G/G_0$ is a (locally) $2k$-symmetric space (it is globally $2k$-symmetric if $\tau$ integrates in $G$) and is in particular a reductive homogeneous space (reductive decomposition $\g=\g_0\oplus\mak n$).\\
 Moreover since $\sigma=\tau^2$ is an order $k$ automorphism, then for any subgroup $H$, such that $(G^\sigma)^0\subset H\subset G^\sigma$, $G/H$ is a (locally) $k$-symmetric space. In all the following we will always do this choice for $H$ and suppose that $H\supset G_0$ (it is already true up to covering since $\hk\supset\g_0$) so that $N=G/G_0$ has a structure of associated bundle over $M=G/H$ with fibre $H/G_0$: $G/G_0\cong G\times_H H/G_0$. We can add that on $\hk$, $\tau$ is an involution: $(\tmh)^2=\Id_\hk$, whose symmetric decomposition is $\hk=\g_0\oplus\g_k$, and gives rise to the (locally) symmetric space $H/G_0$. The fibre $H/G_0$ is thus (locally) symmetric (and globally symmetric if the inner automorphism $\Int\tm$ stabilizes $\Adm H$). Owing to the effectivity of $(\g,\hk)$, we have the following characterization:
\begin{eqnarray}\label{go-gk}
\g_0 & = & \{\xi\in\hk| [\adm\xi,\tm]=0\}\\
\g_k & = & \{\xi\in\hk| \{\adm\xi,\tm\}=0 \}
\end{eqnarray} 
$\{\}$ being the anticommutator.%\medskip\\
\paragraph{Two different types of $2k$-symmetric spaces.}
Since $(\tm)^k$ is an involution, there exist two invariant subspaces $\mk'$ and $\mk''$, of $\mk$, each sum of certain $\mk_j$'s,  such that
$$
(\tm)^k=-\Id_{\mk'}\oplus\Id_{\mk''}.
$$
These subspaces $\mk'$ and $\mk''$ can be computed easily : $\mk'=\oplus_{j=0}^{[\frac{k-2}{2}]}\mk_{2j+1}$ and $\mk''=\oplus_{j=1}^{[\frac{k-1}{2}]}\mk_{2j}$. In other words, their complexifications are given by 
$$
{\mk'}^\C=\bigoplus_{\begin{array}{c} z^k=-1\\ z\neq -1\end{array}}\ker (\tau- z\Id),\quad {\mk''}^\C=\bigoplus_{\begin{array}{c} z^k=1\\ z\neq \pm 1\end{array}}\ker (\tau- z\Id).
$$
At this stage, there are two possibilities:
\begin{description}
\item[$\bullet$] if $\mk''=0$ then $(\tm)^k=-\Id_\mk$ and $\tm$ admits eigenvalues only in the set $\{z^k=-1, z\neq -1\}$.
\item[$\bullet$] if $\mk''\neq 0$ then $(\tm)^k\neq -\Id_\mk$ and $\tm$ admits eigenvalues in each the sets $\{z^k=1, z\neq \pm 1\}$ and $\{z^k=-1, z\neq -1\}$.
\end{description} 
These two cases give rise to two different types of $2k$-symmetric spaces (see section~\ref{return})%\medskip.\\
\paragraph{$G$-invariant metrics.}
Now, let us suppose  that $M=G/H$ is Riemannian (i.e. $\Adm H$ is compact) then we can choose an $\Ad H$-invariant inner product on $\mk$ for which $\tm$ is an isometry\footnote{See the Appendix, theorem~\ref{AdmH,taum-compact}, for the proof of the existence of such a inner product}.  In the next of the paper, we will always choose this kind of inner product on $\mk$. Therefore, $\tm$ is  an order $2k$ isometry. We will study this kind of endomorphism in section~\ref{isometry-twistor}.\\
Moreover, let us remark that if $G/H$ is Riemannian then so is $G/G_0$. Further, since the elliptic system we will study in this paper is given in the Lie algebra setting it is useful to know how the fact that $G/H$ is Riemannian can be read in the Lie algebra setting. In fact, under our hypothesis of effectivity, $G/H$ is Riemannian \iif $\hk$ is compactly imbedded\footnote{See \cite[p. 130]{Hel} for a definition.} in $\g$ and $\Ad H/\Ad H^0$ is finite. Moreover, according to proposition~\ref{AdHsurAdH0finite} in the Appendix, $\Ad H/\Ad H^0$ is always finite so that $G/H$ is Riemannian \iif $\hk$ is compactly imbedded in $\g$.
\index{even case|)}
\subsubsection{The odd case: $k'=2k+1$}\label{2.1.2}
\index{odd case|(} 
 As in the even case we have $\g_0^\C=(\g_0)^\C$ and $\overline{\g_j^\C} =\g_{-j}^\C,\, \forall j\in\Z/k'\Z$. Then we obtain the following eigenspace decomposition:
\begin{equation}\label{decimpaire}
\g^\C=\left(\g_{-k}^\C\oplus\ldots\oplus\g_{-1}^\C\right)\oplus\g_0^\C \oplus\left(\g_{1}^\C \oplus\ldots\oplus\g_{k}^\C\right),
\end{equation}
which provides in particular the following reductive decomposition:
$$
\g=\g_0\oplus\mk
$$
with $\mk=\oplus_{j=1}^k\mk_j$ and $\mk_j$ is the real subspace whose the complexification is $\mk_j^\C=\g_{-j}^\C\oplus\g_j^\C$. According to our convention, we suppose  that $(\g,\g_0)$ is effective.\\
Then, as in the even case, integrating our setting and choosing $G_0$ such that $(G^\tau)^0\subset G_0\subset G^\tau$, we consider  $N=G/G_0$ which is a locally $(2k+1)$-symmetric space and in particular a reductive homogeneous space. Moreover, the decomposition (\ref{decimpaire}) gives rises to  a splitting $TN^\C=T^{1,0}N\oplus T^{0,1}N$ defined by  
\begin{equation}\label{def-odd}
\begin{array}{rcccc}
\displaystyle TN^\C & = & \dfrac{}{}\left(\oplus_{j=1}^k[\g_{-j}^\C]\right) & \oplus & \dfrac{}{} \left(\oplus_{j=1}^k[\g_{j}^\C]\right)\\
\displaystyle    & =  &    \dfrac{}{}T^{1,0}N                           &   \oplus     &  \dfrac{}{} T^{0,1}N
\end{array}
\end{equation}
This splitting defines a canonical almost complex structure on $G/G_0$, that we will denote by $\undj$.\index{canonical!almost complex structure}\medskip\\
Let us suppose that $N=G/G_0$ is Riemannian then the subgroup generated by $\Adm G_0$ and  $\tm$ is compact (because $\Adm G_0$ is compact and $\tm\Adm g\,\tm^{-1}=\Adm g$, $\forall g\in G_0$, and $\taum$ is of finite order). Therefore there exists an $\Ad G_0$-invariant inner product on $\mk$ for which $\tm$ is an isometry. In the next of the paper, we will always choose this kind of inner product on $\mk$ (when $N$ is Riemannian).
\index{odd case|)}

\subsection{Definitions and general properties of the $m$-th elliptic system.} 
\subsubsection{Definitions}\label{def-subsec}
Let  $\tau\colon \g \to \g$ be an order $k'$ automorphism with $k'\in\mathbb{N}^*$ (if $k'=1$ then $\tau=\Id$). We use the notations of \ref{finitorderauto}. Let us begin by defining some useful notations.\medskip\\
\textbf{Notation and convention}
Given $I\subset \mathbb N$, we denote by $\prod_{j\in I}\g_j^\C$, the product $\prod_{j\in I}\g_{j \,\mrm{mod}\,k'}^\C$. In the case $\sum_{j\in I}\g_{j \,\mrm{mod}\,k'}^\C$ is a direct sum in $\g^\C$, we will identify it with the previous product via the canonical isomorphism
\begin{equation}\label{via}
(a_j)_{j\in I}\longmapsto \sum_{j\in I}a_j,
\end{equation}
and we will denote these two subspaces by the same notation $\oplus_{j\in I}\g_j^\C$\medskip.\\
Now, let us define the $m$-th elliptic integrable system associated to a $k'$-symmetric space, in the sense of Terng \cite{tern}.
\begin{defn}\label{m-syst-}
Let $L$ be a Riemann surface.
The $m$-th $(\g,\tau)$-system (with the $(-)$-convention) on $L$ is the equation for $(u_0,\ldots,u_m)$, $(1,0)$-type 1-form on $L$ with values in $\prod_{j=0}^m\g_{-j}^\C$:
$$
\left\{
  \begin{array}{lcr}
\displaystyle\bar\partial u_j + \sum_{i=0}^{m-j} [\bar u_i\wedge u_{i+j}]=0 &  (S_j), & \text{  if  } 1\leq j\leq m ,\\
\displaystyle\bar\partial u_0 + \partial \bar u_0 + \sum_{j=0}^{m} [u_j\wedge \bar u_j]=0 & (S_0) &
\end{array}\right.\qquad \qquad\qquad\mrm{(Syst)}\\
$$
It is equivalent to say that the 1-form
\begin{equation}\label{alphalambda}
\alpha_\lm= \sum_{j=0}^m \lm^{-j}u_j + \lm^j \bar u_j =\sum_{j=-m}^m \lm^j\hat{\alpha}_j
\end{equation}
satisfies the zero curvature equation:
\begin{equation}\label{courbnul}
d\alpha_{\lm} + \frac{1}{2}[\alpha_{\lm}\wedge\alpha_{\lm}]=0, \quad \forall\lm\in \C^*.
\end{equation}
 \end{defn}
\begin{defn}\label{m-syst+}
Let $L$ be a Riemann surface.
The $m$-th $(\g,\tau)$-system (with the $(+)$-convention) on $L$ is the equation ($\syst$) as in definition~\ref{m-syst-} but for  $(u_0,\ldots,u_m)$, $(1,0)$-type 1-form on $L$ with values in $\prod_{j=0}^m\g_{j}^\C$:\footnote{instead of $\prod_{j=0}^m\g_{-j}^\C$.}\\
It is equivalent to say that the 1-form
\begin{equation}\label{alphalambda+}
\alpha_\lm= \sum_{j=0}^m \lm^{j}u_j + \lm^{-j} \bar u_j =\sum_{j=-m}^m \lm^j\hat{\alpha}_j
\end{equation}
satisfies the zero curvature equation (\ref{courbnul}).
\end{defn}
\begin{rmk}\label{rmk-convention}\em
The difference between the two conventions is that in the first one  $\alpha_\lm '=\sum_{j=0}^m \lm^{-j}u_j$ involves negative powers of $\lambda$ whereas in the second one $\alpha_\lm '$ involves positive powers of $\lambda$ (in other words $\hat\alpha_{-j}''=0$, for $j\geq 1$ in the first one whereas $\hat\alpha_{j}''=0$, for $j\geq 1$ in the second one).
In fact the second system is the first system associated to $\tau^{-1}$ and vice versa.\\
The first convention is the traditional one: it was used for harmonic maps into symmetric space (see \cite{DPW}) and by H\' elein-Romon \cite{HR1,HR2,HR3} for Hamiltonian stationary Lagrangian surfaces in Hermitian symmetric space -- first example of second elliptic integrable system associated to a 4-symmetric space. Then this convention was used  in \cite{ki1,ki2,ki3}. Terng \cite{tern}, herself,  in her definition  of the elliptic integrable system  uses also this convention. However in \cite{bk}, it is the second convention which is used.\\
The $(+)$-convention is in fact the most natural, as we will see, since  in the $(-)$-convention, there is a minus sign  which appears when we pass from the Lie algebra setting to the geometric setting. This is what happens for example when we will associate to the automorphism $\tau$, a almost complex complex structure in the target space (see also \cite{ki3}, remark 13). Namely, we will consider that the eigenspace spaces $\g_{-j}$ will define the  $(1,0)$-part and the subspaces $\g_j$  the $(0,1)$-part of the tangent space of the target. But the $(+)$-convention  leads to several changes in the traditional conventions, like for example in the DPW method \cite{DPW}, we must use the Iwasawa decomposition $\Lm\Gtc= \Lm\Gt .\Lm_{\mathcal B}^-\Gtc$ instead of
$\Lm\Gtc= \Lm\Gt .\Lm_{\mathcal B}^+\Gtc$ and in particular the holomorphic potential involves positive power of $\lambda$ instead of negative one as it is the case traditionally. We decided here to continue to perpetuate the tradition as  in \cite{ki3} and to use the first convention. So in the following \textbf{when we will speak about the $m$-th elliptic integrable system, it will be according to the definition~\ref{m-syst-}}.
\end{rmk}
\textbf{Notation} Sometimes, to avoid confusion we will denote (Syst) either by (Syst($m,\g,\tau$)), (Syst($m,\tau$)) or simply by (Syst($m$)) depending on the context and the needs.\\
For shortness we will also often say the $(m,\g,\tau)$-system instead of the $m$-th $(\g,\tau)$-system. We will also say the $m$-th elliptic (integrable) system associated to (the $k'$-symmetric space) $G/G_0$.\\
We will say that a family of 1-forms $(\alpha_\lm)_{\lm\in \C^*}$ (denoted by abuse of notation, simply by $\alpha_\lm$) is solution of the  $(m,\g,\tau)$-system  (or of (Syst)) if it corresponds to some solution $u$ of this system, according to (\ref{alphalambda}). Therefore $\alpha_\lm$ is solution of the  $(m,\g,\tau)$-system \iif it can be written in the form (\ref{alphalambda}), for some $(1,0)$-type 1-form $u$ on $L$ with values in $\prod_{j=0}^m\g_{-j}^\C$, and satisfies the zero curvature equation (\ref{courbnul}).\medskip\\
It will turn out that the $(m,\g,\tau)$-system has distinctively different behaviour if $k'$ is even and if $k'$ is odd. Moreover, for every $k'$ there are three different types of behaviour, according to the size of $m$ relatively to $k'$.  

\begin{defn}\label{defn-cases}\index{determined maximal@determined, maximal}\index{determined minimal@determined, minimal}\index{primitive@primitive, case or system|(}\index{underdetermined@underdetermined, case or system|(}
If $k'=1$, set  $m_{k'}=0$. If $k'>1$, we set  $m_{k'}=\displaystyle \left[\frac{k'+1}{2}\right]=\begin{cases} k  \text{ if } k'=2k\\ k +1 \text{ if } k'=2k+1 \end{cases}$.\\
We will say that the $m$-th $(\g,\tau)$-system is:\\
-- in the \textbf{primitive case} (or that  the system is primitive) if $0\leq m < m_{k'}$,\\
-- in the \textbf{determined case} (or that  the system is  determined) if $ m_{k'} \leq m \leq k' -1 $,\\
-- and in the \textbf{underdetermined case} (or that  the system is  underdetermined) if $ m > k'-1 $.\\
 Moreover, the determined system of minimal order $ m_{k'}$ will be called "the minimal determined system", and the one of maximal order $k'-1$ will be called "the maximal determined system".
\end{defn}
Let us consider the $\g$-valued 1-form $\alpha:=\alpha_{\lm=1}$. Then we have $\alpha = \sum_{j=0}^m u_j + \bar u_j$ according to (\ref{alphalambda}) which is equivalent to $\alpha'=\sum_{j=0}^m u_j$, since $\alpha$ is $\g$-valued\medskip.\\
$\bullet$ In the primitive and determined cases ($m\leq k'-1$), $\sum_{j=0}^m \g_{-j}^\C$ is a direct sum so that
$u=(u_0,\ldots, u_m)$  can be identified with $\sum_{j=0}^m u_j =\alpha'$  via  (\ref{via}) and according to our convention. We will then write simply $u=\alpha'$. In particular  we have 
$$
u_j=\alpha_{-j}' \quad \forall j, 0\leq j\leq m
$$ 
with $\alpha_j:=[\alpha]_{\g_j}$ $\forall j\in \Z/k'\Z$. Hence in the primitive and determined cases the $m$-th $(\g,\tau)$-system can be considered as a system on $\alpha$. Moreover, we can recover  $\alpha_\lm$ from $\alpha$ and we will speak about the "extended Maurer Cartan form" $\alpha_{\lm}$ which is then associated to $\alpha$ by
$$
\alpha_\lm= \sum_{j=1}^m \lm^{-j}\alpha_{-j}' + \alpha_0 + \sum_{j=1}^m \lm^j \alpha_{j}''
$$
according to (\ref{alphalambda})\medskip.\\
$\bullet$ In the underdetermined case, $\sum_{j=0}^m \g_{-j}^\C$ is not a direct sum so that to a given $\alpha$ (coming from some solution $\alpha_\lm$ of the $m$-th $(\g,\tau)$-system, according to $\alpha=\alpha_{\lm=1}$) there are a priori many (other) corresponding solutions $u=(u_0,\ldots, u_m)$ since
$$
\forall j\in \Z/k'\Z, \alpha_{-j}'=\sum_{i\equiv j[k']} u_i . 
$$
In fact, we will prove that there are effectively an infinity of other solutions $\alpha_\lm$ satisfying the condition $\alpha_{\lm=1}=\alpha$ (see theorem~\ref{undetermined-non-injective}; see also section~\ref{underdetermined} for a conceptual  explanation).
%
%
%%%%%%%%%%%%%%%%%%%%%%%%%%%%%%%%%%
%
\subsubsection{The geometric solution}\label{geom-sol}
\index{geometric!solution}
\textbf{Convention.} 
Our study, in the present paper, is local therefore we will suppose (when it is necessary to do so) either that $L$ is   simply connected or that all lifts (of maps defined on $L$) and integrations (of 1-forms on $L$) are made locally. We consider that these considerations are implicit and will not specify these most of the time.\medskip\\
The equations (\ref{courbnul}) and (\ref{alphalambda}) are  invariant by gauge transformations by the group $C^\infty(L,G_0)$:
$$
U_0\cdot \alpha_\lm =\Ad U_0 (\alpha_\lm) -dU_0.U_0^{-1}.
$$
where $U_0\in C^\infty(L,G_0)$. This means that  if $\alpha_\lm$ satisfies \eqref{courbnul}  then so is  $U_0\cdot \alpha_\lm$ and if $\alpha_\lm$ can be written in the form \eqref{alphalambda} then so is  $U_0\cdot \alpha_\lm$. Therefore if $\alpha_\lm$ is a solution of (Syst) then so is $U_0\cdot \alpha_\lm$.
This allows us to define a  \textbf{geometric solution} of (Syst) as follows:
\begin{defn}
A map $f\colon L\to G/G_0$ is a \textbf{geometric solution} of \emph{(Syst)} if for any (local) lift   $U$ of $f$,  into $G$,   there exists a (local) solution $\alpha_{\lm}$   of \emph{(Syst)} such that $U^{-1}.dU=\alpha_{\lm=1}$.
\end{defn}
In other words, we obtain the set of geometric solutions as follows: for each solution $\alpha_\lm$ of (Syst), consider the $\g$-valued 1-form $\alpha:=\alpha_{\lm=1}$, then integrate it by $U\colon L\to G$, $U^{-1}.dU=\alpha$, and finally project $U$ on $G/G_0$ to obtain the map $f\colon L\to G/G_0$.
\medskip\\
Now, to simplify the exposition, let us suppose that $L$ is simply connected (until the end of \ref{geom-sol}). Then  $(\alpha_{\lm})_{\lm\in \C^*}\mapsto \alpha=\alpha_{\lm=1}$ is a surjective map from the set of solutions of (Syst) to the set of Maurer-Cartan forms of  lifts of  geometric solutions. According to the discussion at the end of subsection \ref{def-subsec}, this map  is bijective in the primitive and determined case ($m\leq k'-1$) and not injective in the underdetermined case ($m>k'-1$). By quotienting by $C^\infty(L,G_0)$, we obtain a surjective map $\pi_m$ with the same properties, taking values in the set of geometric solutions\medskip.\\
Let us make more precise all that. We suppose, until the end of this subsection~\ref{geom-sol},  that the automorphism $\tau\colon\g\to \g$ is fixed (so that the only  data which varies in the $(m,\g,\tau)$-system is  the order $m$). First, let us give an explicit expression of the space $\mathcal S(m)$ of solutions $\alpha_\lm$ of the system (Syst($m$)), i.e. the solutions of  the zero curvature equation (\ref{courbnul}), which satisfies the equality (\ref{alphalambda}) for some  $(1,0)$-type 1-form $u$ on $L$ with values in $\prod_{j=0}^m\g_{-j}^\C$. To do that, we want to express the condition to be written in the form (\ref{alphalambda}) as a condition on  the loop $\alpha_\bullet\colon\lm\in S^1\mapsto \alpha_\lm\in \mal C(T^*L\otimes \g)$. The " $\bullet$ " means of course functions on the parameter\footnote{Remark that $\alpha_\bullet$ determines $(\alpha_{\lm})_{\lm\in \C^*}$, when this latter satisfies (\ref{alphalambda}).} $\lm\in S^1$.  More precisely, we will consider $\alpha_\bullet$ as a 1-form on $L$ with values  in the loop Lie algebra $C(S^1,\g)$. Then the condition to be written in the form (\ref{alphalambda}) means: 
\begin{equation}\label{express}
(\ref{alphalambda}) \Longleftrightarrow (\alpha_\bullet\in\Lm_m\gtau \quad \mrm{and}\quad \alpha_\bullet'\in \Lm^-\gtau^\C)
\end{equation}
where
\begin{eqnarray*}
\Lm\gtau & = & \{\eta_\bullet\colon S^1\to \g|\,\eta_{\omega\lm}=\tau(\eta_\lm),\forall \lm \in S^1\}\\
\Lm_m\gtau & = & \{\eta_\bullet\in\Lm\gtau|\,\eta_\lm=\sum_{|j|\leq m}\lm^j \hat{\eta}_j\}\\
\Lm^-\gtau^\C & = & \{\eta_\bullet\in\Lm\gtau^\C|\,\eta_\lm=\sum_{j\leq 0} \lm^j \hat{\eta}_j\}
\end{eqnarray*}
and $\omega$ is a $k'$-th primitive root of unity. We refer the reader to \cite{PS} for more details about loop groups ad their Lie algebras (in particular about the possible choices of topology which makes $\Lm\gtau$ be a Banach Lie algebra). Therefore the space of solutions of (Syst($m$)) is given by
\begin{equation}\label{Sm}
\mal{S}(m)=\{\alpha_\bullet\in\mal C(T^*L\otimes\Lm_m\gtau)|\, \alpha_\bullet'\in \Lm^-\gtau^\C \quad \mrm{and} \quad
d\alpha_{\bullet} + \frac{1}{2}[\alpha_{\bullet}\wedge\alpha_{\bullet}]=0\}.
\end{equation}
Let us remark that the condition $\alpha_\bullet'\in \Lm^-\gtau^\C$ can be interpretated as a condition of $\C$-linearity. Indeed, the Banach vector space $\Lm\gtau/\g_0 $ is naturally endowed with the complex structure defined by the following decomposition 
\begin{equation}\label{complex-struct-loop}
(\Lm\gtau/\g_0)^\C=\Lm\gtau^\C/\g_0^\C=\Lm_*^-\gtau \oplus \Lm_*^+\gtau,
\end{equation} 
where $\Lm_*^{\pm}\gtau=\{\eta_\bullet\in\Lm\gtau^\C|\,\eta_\lm=\sum_{j \gtrless 0} \lm^j \hat{\eta}_j\}$. Then the condition $\alpha_\bullet'\in \Lm^-\gtau^\C$ means that $[\alpha_\bullet']_*\colon TL\to (\Lm\gtau/\g_0)^\C$ is $\C$-linear, where $[\ ]_*$ denotes the component in $\Lm_*\gtau=\{\eta_\bullet\in\Lm\gtau| \eta_\lm=\sum_{j \neq 0} \lm^j \hat{\eta}_j\}\cong \Lm\gtau/\g_0$.\medskip\\
Now let us integrate our setting. Firstly, let us define the twisted loop group (\cite{PS})
$$
\Lm\Gt=\{U_{\bullet} \colon S^1\to G | U_{\omega\lm}=\tau\left(U_\lm\right)\}.
$$
Then, let us set
\begin{eqnarray*}
\mathcal{E}^m & = & \{ U_{\bullet}\colon L\to \Lm \Gt | U_{\lm}(0)=1, \forall \lm\in S^1; \alpha_\lm:= U_\lm^{-1}.d U_\lm \text{ is a solution of } (\mrm{Syst}(m))\}\\
\mathcal{E}_1^m & = & \{ U \colon L\to G|\exists U_{\bullet}\in \mathcal{E}^m, U=U_1\}\\
\mathcal{G}_1^m & = & \{f\colon L\to G/G_0\ \text{geometric  solution of}\ (\mrm{Syst}(m)), f(0)=1.G_0\}\\
\mathcal{G}^m & = & \{f_\bullet=\pi_{G/G_0}\circ U_\bullet,\, U_\bullet\in \mathcal{E}^m\}
\end{eqnarray*}
Remark that because of the gauge invariance: $\mal E(m).\mal K\subset \mal E(m)$ where $\mal K=C_*^\infty(L,G_0)=\{U\in C^\infty(L,G_0)| U(0)=1\}$, any lift\footnote{With initial condition $U(0)=1$.} $U_\bullet\colon L\to \Lm \Gt$   of an element $f_\bullet\in \mathcal{G}^m$  belongs to $\mathcal{E}^m$. 
\begin{defn}
An element $f_\bullet\in \mathcal{G}^m$ will be called an \emph{extended geometric solution} of $(\systm)$.
\end{defn}
The space of geometric solutions is obviously obtained from the space of extended geometric solutions $\mathcal{G}^m$ by the evaluation at $\lm=1$. Moreover $\mal S(m)\simeq\mal E^m$ is determined by $\mathcal{G}_\bullet^m$ because of the gauge invariance: $\mal E(m).\mal K\subset \mal E(m)$,  so that we can write $\mathcal{G}^m=\mal E(m)/\mal K$. Consequently, we have also $\mathcal{G}_1^m=\mathcal{E}_1^m/\mal K$.\\
Finally, we obtain the following diagram 
$$
\begin{CD}
\mal S(m)@> \mrm{int}>\cong > \mal E^m@>\pi_{\mal K}>> \mal E^m/\mal K @=\mathcal{G}_\bullet^m\\
@V\mrm{ev}_1 VV @VVV @VVV @V\pi_m VV\\
\mal S(m)_1@> \mrm{int}>\cong > \mal E_1^m@>>> \mal E_1^m/\mal K @=\mathcal{G}^m
\end{CD}.
$$
Therefore, the surjective map $\pi_m$ is bijective for $m\leq k'-1$ and not injective for $m>k'-1$ (because so is $\mrm{ev}_1$). We will need the following definition:
\begin{defn}\index{geometric!map} 
Given a $\g$-valued Maurer-Cartan 1-form $\alpha$ on $L$, we define \textbf{the geometric map corresponding to} $\alpha$, as $f=\pi_{G/G_0}\circ U$, where $U$ integrates $\alpha$: $U^{-1}.dU=\alpha$, $U(0)=1$.
\end{defn}
We have seen that in the primitive and determined cases, we can consider  $(\systm)$ as a system on the $\g$-valued 1-form $\alpha:=\alpha_{\lm=1}$. Since the Maurer-Cartan equation for $\alpha$ is always contained in $(\systm)$ according to (\ref{courbnul}), this systems on $\alpha$ is itself equivalent to a system on the geometric map $f$ corresponding to $\alpha$. This system on $f$ is then a $G$-invariant elliptic PDE on $f$ of order $\leq 2$.\\
Let us  summarize:
\begin{prop}
The natural surjective map $\pi_m \colon \mal G^m \to \mal G_1^m$ from the set of extended geometric solutions of the  $(m,\g,\tau)$-system into the set of geometric solutions  is bijective in the primitive and determined cases ($m\leq k'-1$) and not injective in the underdetermined case ($m> k'-1$). Moreover, in the primitive and determined cases, the $(m,\g,\tau)$-system - which is initially a system on the $\Lm_m\gtau$-valued 1-form $\alpha_\lm$ - is in fact a system on the 1-form $\alpha:=\alpha_{\lm=1}$,  itself  equivalent to an elliptic PDE of order $\leq 2$ on the corresponding geometric map $f\colon L\to G/G_0$.
\end{prop}
Furthermore, let us interpret the $\C$-linearity of $[\alpha_\bullet']_*\colon TL\to (\Lm\gtau/\g_0)^\C$ in terms of  the corresponding extended geometric solution $f_\bullet\colon L\to \Lm\Gt/G_0$, defined by $f_\bullet =\pi_{G/G_0}\circ U_\bullet$ where $U_\bullet$ integrates $\alpha_\bullet$. Firstly, the complex structure defined in $\Lm\gtau/\g_0 $ by (\ref{complex-struct-loop}) is $\Ad G_0$-invariant so that it defines a $\Lm \Gt$-invariant complex structure on the homogeneous space $\Lm\Gt/G_0$. Therefore the $\C$-linearity of $[\alpha_\bullet']_*$ means exactly that $f_\bullet\colon L\to \Lm\Gt/G_0$ is holomorphic. Now, let us interpret the condition $\alpha_\bullet \in \Lm_m\gtau$ in terms of the map $f_\bullet$. Let us consider the following $\Ad G_0$-invariant decomposition
$$
\Lm\gtau/\g_0=\Lm_{m*}\gtau\oplus \Lm_{>m}\gtau
$$
where $\Lm_{m*}\gtau=\Lm_m\gtau\cap \Lm_*\gtau$ and $\Lm_{>m}\gtau=\{\eta_\bullet\in\Lm\gtau|\,\eta_\lm=\sum_{|j| > m} \lm^j \hat{\eta}_j\}$, which gives rise respectively to some $\Lm\Gt$-invariant splitting
$$
T(\Lm\Gt/G_0)=\hor_m^\Lm\oplus\ver_m^\Lm.
$$
 Then $\hor_m^\Lm$ and $\ver_m^\Lm$ inherit respectively  the qualificatifs horizontal and vertical subbundle respectively. Therefore, in the same spirit as \cite{DPW} (remark~{2.5} and proposition~{2.6}), the equation~(\ref{Sm}) gives us the following familiar twistorial characterization
\begin{prop}\label{prop-holom-hori}
A map $f_\bullet\colon L\to \Lm\Gt/G_0$ is an extended geometric solution of the $(m,\g,\tau)$-system \iif it is holomorphic and horizontal.
\end{prop}
\subsubsection{The increasing sequence of spaces of solutions: $(\mal S(m))_{m\in \mathbb N}$}\label{increasing}
Again, we suppose in all \ref{increasing} that the automorphism $\tau$ is fixed and that $L$ is simply connected. Then according to the realisation of $(\systm)$ in the forms (\ref{courbnul}) and (\ref{alphalambda}), we see  that any solution of $(\systm)$ is solution of $(\syst(m'))$ for $m\leq m'$ (take $u_j=0$ for $m<j\leq m'$). More precisely, $(\systm)$ is a reduction of $(\syst(m'))$: $(\systm)$ is obtained from $(\syst(m'))$ by putting
$u_j=0$, $m<j\leq m'$, in $(\syst(m'))$. In particular, $\mal S(m)\subset\mal S(m')$ for $m\leq m'$; so that any solution in the primitive case ($m< m_{k'}$) is solution of any determined system ($m_{k'}\leq m\leq k'-1$), and any solution of a determined system is solution of any underdetermined system ($m>k'-1$). 
$$
\{\text{Primitive case}\} \subset \{\text{determined case}\} \subset \{\text{underdetermined case}\}.
$$
\begin{rmk}\em
We have ${\pi_{m'}}_{|\mal G^m}=\pi_m$ if $m\leq m'$. In particular, $\pi_{m'}(\mal G^m)=\mal G_1^m $. We can set  $\mal S(\infty)=\cup_{m\in\mathbb N}\mal S(m)$, $\mal E^\infty=\cup_{m\in\mathbb N}\mal E^m$ and $\mal G^\infty =\cup_{m\in \mathbb N} \mal G^m$. Then we have $\mal G^\infty = \mal E^\infty /\mal K $. Moreover we can define the surjective map $\pi_\infty\colon \mal G^\infty\to \mal G_1^\infty$ such that $ {\pi_\infty}_{|\mal G^m}=\pi_m$, $\forall m\in \mathbb N$. Then ${\pi_\infty}_{|\mal G^m}$ is a bijection onto $\mal G_1^m$ for each $m\leq k'-1$.\\
We can call $\mal S(\infty)$ the $(\g,\tau)$-system, and then we can speak about its subsystem of order $m$, namely
$\mal S(m)$. In particular, we have the following characterization:
$$
\mal S(\infty)=\{\alpha_\bullet\in \mal C(T^*L\otimes\Lm_{(\infty)}\gtau)|\, \alpha_\bullet'\in\Lm^-\gtau \quad \mrm{and} \quad
d\alpha_{\bullet} + \frac{1}{2}[\alpha_{\bullet}\wedge\alpha_{\bullet}]=0\}
$$
where $\Lm_{(\infty)}\gtau=\cup_{m\in \mathbb N}\Lm_m\gtau$.
\end{rmk}
\textbf{Important Remark}
It could happens that  the eigenspaces $\g_j$ vanishes for the first values of $j$, i.e. $j$ close to $0$. For example it is a priori possible that $\g_0=0$ (which implies that our $k'$-symmetric space is a group). Then for the values of $m\geq k'$ close to $k'$, the underdetermined systems $\syst(m)$ coincide trivially  with the determined system $\syst(k'-1)$, because then $u_j\in\g_{-j}=\{0\}$, for $k'\leq j\leq m$. This is why, when we say "underdetermined", we mean in fact "underdetermined but not determined" to exclude the (potential) formal underdetermined systems which are in fact trivially determined because of the (potential) vanishing of the first eigenspaces. Remark that this eventuality  to happen need that all the eigenspaces from $0$ to a small value, vanish (they can not all vanish otherwise $\g=0$), and in particular, this implies that $\g_0=0$ (the automorphism $\tau$ has no fixed point). 
\paragraph{The primitive and determined cases ($m\leq k'-1$)}
Now, let us apply the previous discussion (about the increasing sequence $(\mal S(m))_{m\in \mathbb N}$) to the study of the determined case. Let us recall that in this case, we can consider that the system ($\systm$) deals only with $\g$-valued 1-forms $\alpha$. Let us also keep  in mind  that, in this case, we have   $\alpha'=u\in \oplus_{j=0}^m\g_{-j}^\C$ (see the discussion in the end of \ref{def-subsec}, after definition~\ref{defn-cases}). Then we obtain immediately:
\begin{prop}\label{hol-hor-condition}\index{determined maximal@determined, maximal}\index{determined minimal@determined, minimal}\index{holomorphicity condition}
The solutions of a determined system $(\systm)$, $m_{k'}\leq m\leq k'-1$, are exactly the solutions of the maximal determined system, i.e. $(\syst(k'-1))$, which satisfy the \emph{holomorphicity conditions}:
$$
\alpha_{-j}''=0,\   1\leq j\leq k'-1-m.
$$
Moreover, the solutions of a primitive system $(\systm)$, $1\leq m\leq m_{k'}-1$, are the solutions of the minimal determined system, i.e. $(\syst(m_{k'}))$, which satisfy 
\begin{description}
\item[(i)] if $k'=2k$ is even, the \emph{horizontality conditions}:
$$
\alpha_{k}=\alpha_{\pm(k-1)}=\dots\alpha_{\pm(m+1)}=0
$$
\item[(ii)] if $k'=2k+1$ is odd, 
$$
\begin{array}{ll}
\bullet\emph{ the holomorphicity condition : }  \alpha_{-k}''= 0 & \text{ if } m=k,\\
\bullet\emph{ the horizontality conditions : } \alpha_{\pm k}=\alpha_{\pm(k-1)}=\dots\alpha_{\pm(m+1)}=0 & \text{ if } m\leq k-1.
\end{array}
$$
\end{description}
\end{prop}

\paragraph{The non injectivity of $\pi_m$ in the underdetermined case} Now, let us turn ourself to the underdetermined case. We want to study the surjective map $\pi_m\colon\mal G^m\to \mal G_1^m$, in this case.
%We want to prove our claim: the map $\alpha_\bullet\in \mal S(m)\to \alpha\in \mal S(m)_1$ is not injective.

\begin{thm}\label{undetermined-non-injective}
The surjective map $\pi_m\colon\mal G^m\to \mal G_1^m$, from the set of extended geometric solutions into the set of geometric solutions, is a principal bundle with as structure group some group of holomorphic curves into $\Lm\Gt/G_0$. In the determined and primitive cases, this group is trivial whereas in the underdetermined case it is of infinite dimension.\\
Therefore, in the underdetermined case,  the surjective map $\pi_m\colon\mal G^m\to \mal G_1^m$  is not injective and its  fibers are of an infinite dimension (even up to conformal transformations of $L$). A fortiori, so is for the map $\mrm{ev_1}\colon\alpha_\bullet\in \mal S(m)\to \alpha:=\alpha_{\lm=1}\in \mal S(m)_1$.\\

\end{thm}
\proof
The natural map $\Lm \Gt \to G$, $g\to g(1)$ is a morphism of group and therefore a principal fibre bundle, with as structure group the kernel
$$
\maH =\{g\in\Lm \Gt| g(1)=1_G\}.
$$
Then it induces  the fibration $\pi\colon \Lm\Gt/G_0 \to G/G_0$ which is also a  $\maH$-principal bundle (remark that $G_0\cap \maH = \{1\}$), $\maH$ acting by the left on $\Lm\Gt/G_0$.\\  
Let $\mrm{Hol}^m(L,\Lm\Gt/G_0)$ be the set  of holomorphic integrale curves $f_\bullet$   of the holomorphic $\Lm \Gt$-invariant distribution $\hor_m^\Lambda$, such that $f_\bullet(o)=1.G_0$, where $o$ is a reference point in $L$.
According to proposition~\ref{prop-holom-hori}, we have $ \mal G^m= \Lm\Gt . \mrm{Hol}^m(L,\Lm\Gt/G_0) $. Therefore, denoting $\left( \mal G_1^m\right)_o = \{g\in \mal G_1^m|g(o)=1.G_0\}$, we  want to prove that the following map 
\begin{equation}\label{the-map-proj-ev1}
\pi_m \colon f_\bullet \in \mrm{Hol}^m(L,\Lm\Gt/G_0)\mapsto f_1 \in \left( \mal G_1^m\right)_o
\end{equation}
is a principal fibre bundle.
To fix ideas, let us first consider the fibre  defined by the constant map $1.G_0$. This is nothing but the holomorphic integral curves of the holomorphic $\Lm \Gt$-invariant distribution $\left( \hor_m^\Lm\right)_*$ defined by the complex subspace $\{\eta\in \Lm_{m*}\gtau|\, \eta_1=(J\eta)_1=0\}$,  such that $f_\bullet(o)=1.G_0$. Of course $J$ is  the complex structure on $\Lm \Gt/G_0$. Furthermore, remark that $\left( \hor_m^\Lm\right)_*$ vanishes in the determined case and is non trivial in the undetermined case.\\
More generally, let us compute the  fibre of any $f_\bullet\in\mrm{Hol}^m(L,\Lm\Gt/G_0)$. Firstly, we remark that $\maH$ is immersed into $\Lm\Gt/G_0$ via the projection $\pi_{G_0}\colon \Lm \Gt \mapsto \Lm \Gt/G_0$. The image of $\maH$ in $\Lm\Gt/G_0$ is $\{h\in\Lm\Gt/G_0| h(1)=1.G_0\}$. Remark also that each element of this image has one and only one lift in $\maH\subset \Lm\Gt$.  Now, let $f, f'\in \mrm{Hol}^m(L,\Lm\Gt/G_0)$ be in the same fibre of $\pi_m$. Therefore, there exists $h\colon L\to \maH$ such that $f'=h.f$. Moreover it is not difficult to see that $h\colon L\to \Lm\Gt/G_0$ is holomorphic. Then we have $df'=dh.f + h.df$ but $df'$ and $h.df$ take values in $\hor_m^\Lm$ so that $dh.f $ takes values in $\hor_m^\Lm$ also. Therefore $dh$ takes values in $\hor_m^\Lm$ i.e. $h$ is a holomorphic integral curve of  $\hor_m^\Lm$ which takes values in $\maH$ and satisfies $h(o)=1.G_0$. This is  in fact equivalent to say that $h$ is a holomorphic integral curve of $\left( \hor_m^\Lm\right)_*$ such that $h(o)=1.G_0$ (then $h\in \maH$ automatically). Conversely, any holomorphic integral curve of $\left( \hor_m^\Lm\right)_*$ such that $h(o)=1.G_0$ satisfies that $d(h.f)$ takes values in $\hor_m^\Lm$, and thus $f':=h.f$ is in $\mrm{Hol}^m(L,\Lm\Gt/G_0)$ and in the same fibre as $f$.\\
We have proved that the fibres of $\pi_m$ are all isomorphic to the set of holomorphic integral curves $h$ of $\left( \hor_m^\Lm\right)_*$ such that $h(o)=1.G_0$. \\
Moreover, proceeding as above we prove that the previous space $\mrm{Hol}_*^m(L,\Lm\Gt/G_0)$ of holomorphic maps is a subgroup of $C^\infty(L,\maH)$ (take for $f$ an element of this space, i.e. in the fibre of the constant map $1.G_0$, and then applying what precedes proves that this space is stable by multiplication of two elements). \comprf\hsq 
\index{underdetermined@underdetermined, case or system|)}

\subsubsection{The decreasing sequence $(\syst(m,\tau^p))_{p/k'}$}
We will call $m$-th $\g$-system, the $m$-th $(\g,\Id)$-system  (i.e.  $u=(u_0,\ldots,u_m)$ takes values in $(\g^\C)^{m+1}$, in definition~\ref{m-syst-}).\\
Any solution of the $m$-th $(\g,\tau)$-system is solution of the $m$-th $\g$-system.% (take $\tau=\Id$, i.e. $u=(u_0,\ldots,u_m)\in(\g^\C)^{m+1})$.
 More precisely, the $m$-th $(\g,\tau)$-system is the restriction to $\oplus_{j=0}^m\g_{-j}^\C(\tau)$ of the $m$-th $\g$-system. \\
More generally, for any $p\in \mathbb N^*$ such that $p$ divides $k'$, the $m$-th $(\g,\tau)$-system is the $m$-th $(\g,\tau^p)$-system restricted to $\oplus_{j=0}^m\g_{-j}(\tau)$, or equivalently - in terms of $\alpha_\lm\in \Lm\g_{\tau^p}^\C$ - restricted to $\Lm\gtc$.
%
%
%\subsubsection{The primitive and determined cases ($m\leq k'-1$)}
%In this case, we can consider that the system ($\systm$) deals only with Maurer-cartan forms $\alpha$ and consequently also with geometric maps $f$. Then we have to translate the equations on $\alpha$ into geometric conditions on $f$. This is what we will begin to do now in the section~\ref{determined}
%
%
%%%%%%%%%%%%%%%%%%%%%%%%%%%%%%%%%%%%%%%%%%%%%%%%%%%%%%%%%%%%%%%%%%%%%%%%

\subsection{The minimal determined case}\label{determined}
\index{determined minimal@determined, minimal|(}\index{vertically harmonic|(}\index{holomorphicity condition|(}\index{vertical tension field|(}
We study here the elliptic system ($\systm$) in the minimal determined case and by the way its subcase the primitive case. Let us recall again that in this case, we can consider that the system ($\systm$) deals only with Maurer-cartan forms $\alpha$ and consequently also with geometric maps $f$. Then we have to translate the equations on $\alpha$ into geometric conditions on $f$. This is what we will begin to do now.\\
The minimal determined case splits into two cases.
\subsubsection{The even minimal determined case: $k'=2k$ and $m=k$}
\index{even case|(}\index{horizontally holomorphic|(}
Let us recall the following decomposition
$$
\g^\C=\left(\g_{-(k-1)}^\C\oplus\ldots\oplus\g_{-1}^\C\right)\oplus\g_0^\C\oplus\left(\g_1^\C\oplus\ldots\oplus \g_{k-1}^\C\right)\oplus\g_k^\C.
$$
It is useful for the following to keep in mind that $k=-k\, \mrm{mod}\, 2k$.
\begin{prop}\label{prop-evendeter}
The system $(\syst(k,\tau))$ can be written
\begin{equation}\label{evendeter}
\left\{
\begin{array}{lc}
\displaystyle\alpha_{-j}''=0, \ 1\leq j\leq k-1 & (H_j)\\
\displaystyle d\alpha + \frac{1}{2}[\alpha\wedge\alpha]=0  & (\mrm{MC})\\
\displaystyle\bar\partial\alpha_{-k}' + [\alpha_0''\wedge\alpha_{-k}']=0 & (S_k)
\end{array}
\right. .
\end{equation}
More precisely the equations $(S_j)$, $0\leq j\leq k-1$, of  $(\syst(k,\tau))$ are respectively  the projection on $\g_{-j}^\C$, $0\leq j\leq k-1$, of (MC) (owing to the holomorphicity conditions $(H_j)$ given by proposition~\ref{hol-hor-condition}).
\end{prop}
\proof
The equation $(S_j)$, for $1\leq j\leq k-1$, is written in terms of $\alpha$:
$$
\Bar\partial \alpha_{-j}' + \sum_{i=0}^{k-j} [\alpha_{i}''\wedge \alpha_{-j-i}']=0,
$$
according to definition~\ref{m-syst-}. Since we have $\alpha_{-j}''=0, \ 1\leq j\leq k-1$, according to proposition~\ref{hol-hor-condition}, this equation is equivalent to 
\begin{equation}\label{eq-MC-j}
d\alpha_{-j} + \sum_{i=0}^{k-j} [\alpha_{i}\wedge \alpha_{-j-i}]=0
\end{equation}
which is nothing but the projection of (MC) on $\g_{-j}^\C$. Indeed,  we have
\begin{multline*}
\frac{1}{2}[\alpha\wedge\alpha]_{\g_{-j}^\C} = \dfrac{1}{2}\left( \sum_{i+l=-j}[\alpha_{i}''\wedge \alpha_{l}'] + \sum_{i+l=-j}[\alpha_{i}'\wedge \alpha_{l}'']\right)
 = \sum_{i+l=-j}[\alpha_{i}''\wedge \alpha_{l}'] =  \sum_{i=-k +1 }^{k} [\alpha_{i}''\wedge \alpha_{-j-i}'] \\
=  \sum_{i=0}^{k-j} [\alpha_{i}''\wedge \alpha_{-j-i}']
\end{multline*}
where we have used,  in the last line,  the holomorphicity conditions: $\alpha_i''=0$ if $-k+1\leq i \leq -1$ and $\alpha_{-j-i}'=0$ if $k-j \leq i \leq k$.\\
Now, since $\alpha$ is real (i.e. is $\g$-valued), this equation \eqref{eq-MC-j} is also equivalent to its complex conjugate i.e. the projection of (MC) on $\g_{j}^\C$. Moreover, the equation $(S_0)$ is written in terms of $\alpha$:
$$
\Bar\partial \alpha_0' + \partial \alpha_0'' + \sum_{j=0}^k [\alpha_{-j}'\wedge\alpha_j'']=0
$$
and moreover, using the holomorphicity conditions, we have 
$$
\sum_{j=0}^k [\alpha_{-j}'\wedge\alpha_j'']= \sum_{j=0}^{k -1}[\alpha_{-j}\wedge\alpha_j] + \dfrac{1}{2}[\alpha_{-k}\wedge\alpha_k]=\dfrac{1}{2}[\alpha\wedge\alpha]_{\g_0}
$$
which proves that $(S_0)$ is equivalent to the projection of (MC) on $\g_0$. Finally, the equation $(S_k)$ is written in terms of $\alpha$ as   in \eqref{evendeter}. 
\comprf\hsq\medskip\\
Moreover, always owing to the holomorphicity conditions, the projection of (MC) on $\g_k$ gives us
$$
d\alpha_k + [\alpha_0\wedge\alpha_k]=0
$$
which is the real part of $(S_k)$. Hence the only new information (in addition to (MC) and ($H$)) given by the minimal determined elliptic integrable system in the even case is the imaginary part of $(S_k)$:
$$
  d(*\alpha_k) + [\alpha_0\wedge (*\alpha_k)]=0 \quad (E_k)
$$
which is as we will see the vertical part of a harmonic map equation. Hence $(\syst(k,\tau))$ is equivalent to
$$
\left\{
\begin{array}{lc}
\displaystyle \frac{}{} \alpha_{-j}''=0,\ 1\leq j\leq k-1 & (H_j)\\
\displaystyle  d(*\alpha_k) + [\alpha_0\wedge (*\alpha_k)]=0 & (E_k)\\
\displaystyle d\alpha + \frac{1}{2}[\alpha\wedge\alpha]=0  & (\mrm{MC})
\end{array}
\right. .
$$
$\bullet$ Now, let us contemplate the equation $(S_k)$. In terms of the geometric map $f\colon L\to G/G_0$, corresponding to $\alpha$, this  equation means
$$
\bar\partial^{( \nabla^0)^v}\partial^v f =0
$$ 
where $( \nabla^0)^v$ is the vertical part of the canonical connection $\nabla^0$ on the homogeneous space $N=G/G_0$. The vertical and horizontal spaces are defined by $\mal V=[\g_k]$ and $\mal H=[\mk]$ since we can do the splitting: $T(G/G_0)=[\mk]\oplus[\g_k]$.\index{horizontal subbundle}\index{vertical subbundle} Moreover,  we will see in section~\ref{vertically-harmonic} (Theorem~\ref{difference-tensor} and §\ref{homspacefibr}) that since the fibration $\pi\colon G/G_0\to G/H$ has a symmetric fibre,  the vertical component $\nabla^v$ of the Levi-Civita connection $\nabla$ on the Riemannian homogeneous space $G/G_0$, coincides in the vertical subbundle $\mal V$ with the  vertical component of the canonical connection $( \nabla^0)^v$. Therefore $(S_k)$ means
$$
\bar\partial^{\nabla^v}\partial^v f =0.
$$
Hence the equation $(E_k)$ is equivalent to  
$$
d^{\nabla^v}(*d^v f)=0 \Longleftrightarrow \tau^v(f):=\mrm{Tr}_g(\nabla^v d^v f)=0
$$
(for any hermitian metric $g$ on the Riemann surface $L$).We will say that $f$ is \emph{vertically harmonic}.\\
We will come back to this with more precisions and details in section \ref{vertically-harmonic} (see in particular theorem~\ref{evendeterhom}).

\paragraph{The primitive case}
The $m$-primitive case is obtained by putting $\alpha_k=\alpha_{\pm(k-1)}=\dots\alpha_{\pm(m+1)}=0$ in the minimal determined case (\ref{evendeter}). In particular $\alpha_k=0$ and $(S_k)$ is trivial so that the only  conditions on the geometric map $f\colon L\to G/G_0$ (whose existence is guaranted by (MC)) are the equations $(H_j)$: $\alpha_{-j}''=\alpha_j'=0$, $1\leq j\leq m$, and $\alpha_k=\alpha_{\pm(k-1)}=\dots\alpha_{\pm(m+1)}=0$, according  to proposition~\ref{hol-hor-condition} .
\begin{prop}\label{partialf}
Let $\tau\colon\g\to\g$ be an order $2k$ automorphism, and $m<k$ a positive integer  then the $m$-th elliptic integrable system $(\systmt)$ means that the geometric map $f\colon L\to G/G_0$ satisfies
$$
\partial f\in \oplus_{j=1}^m[\g_{-j}^\C]\subset T(G/G_0)^\C .
$$
\end{prop}
\textbf{Proof.} Let $f$ be the geometric map corresponding to the Maurer-Cartan form $\alpha$, that we integrate by $U\colon L\to G$, then we have $\partial f=\Ad U(\alpha_{\mak n}')$ and $\alpha$ is solution of $(\systmt$) \iif
$$
\alpha_{\mak n}'=\alpha_{-1}'+\ldots +\alpha_{-m}'\in\oplus_{j=1}^m\g_{-j}^\C \Longleftrightarrow \partial f \in \oplus_{j=1}^m[\g_{-j}^\C].
$$
This completes the proof.\hfill$\square$
\begin{rmk}\emph{
In particular, in the  primitive case $f$ is horizontal ($\alpha_k=0$). Therefore $(S_k)$ is trivial and (owing to the holomorphicity conditions $(H_j), 1\leq j\leq k-1$) the free curvature equation (\ref{courbnul}) is equivalent to (MC) in the primitive case.
} \end{rmk}

\begin{defn}
Let  $m<k$ be a positive integer. We will call  $m$-primitive map (into the locally $(2k)$-symmetric space $G/G_0$) a geometric solution of the system $(\systmt)$.
\end{defn}
%
%%%%%%%%%%%%
%
\paragraph{Geometric interpretation of the equations $(H_j)$.}\index{canonical!$f$-structure|(}\index{f structure@$f$-structure}
For $m< k$, let $F^{[m]}$ be the $f$-structure  on $N=G/G_0$ defined by the following (eigenspace) decomposition:
\begin{equation}\label{def-even}
\begin{array}{rcccccc}
\displaystyle TN^\C & = & \left(\displaystyle\oplus_{j=1}^m[\g_{-j}^\C]\right) & \oplus & \left(\displaystyle\oplus_{|j|>m}[\g_{j}^\C]\right) & \oplus &\left(\displaystyle\oplus_{j=1}^m[\g_{j}^\C]\right)\\
& = & T^+N & \oplus & T^0N  &  \oplus & T^-N
\end{array}.
\end{equation}
We will set $F=F^{[k-1]}$, and we will call $F$ \emph{the canonical $f$-structure on $N$}.
Then according to proposition~\ref{partialf} we have
\begin{thm}\index{f holomorphic@$f$-holomorphic}
A map $f\colon L \to G/G_0$ is $m$-primitive \iif it is $F^{[m]}$-holomorphic.
\end{thm}
\begin{rmk}\emph{
The equations $(H_j)$: $\alpha_{-j}''=0$, $1\leq j\leq m$, on a Maurer-Cartan 1-form $\alpha$ means that the corresponding geometric map $f\colon L\to G/G_0$ satisfies $\mrm{pr}_m\circ(df\circ j_L)=F^{[m]}\circ df$ where $j_L$ is the complex structure in $L$, and $\mrm{pr}_m\colon TN\to \oplus_{j=1}^m[\mk_j]$ is the projection on $\oplus_{j=1}^m[\mk_j]$ along $\left(\oplus_{j=m+1}^{k-1}[\mk_j]\right)\oplus[\g_k]$. This means that the projection $\mrm{pr}_m\circ df\colon TL\to \oplus_{j=1}^m[\mk_j]$ is a morphism of complex vector bundle. Let us denote $C_m^\infty(L,G/G_0)=\{f\in C^\infty(L,G/G_0)|\, df\in\oplus_{j=1}^m[\mk_j] \}$. Then we have the following equivalences between the Maurer-Cartan 1-form $\alpha$ and its geometric map:
\begin{equation}\label{equiv-m}
\begin{array}{rcl}
\displaystyle\alpha\in\g_0\oplus(\oplus_{j=1}^m\mk_j) & \Longleftrightarrow & f\in C_m^\infty(L,G/G_0)\\
\displaystyle (H_j):\ \alpha_{-j}''=0,\ 1\leq j\leq m  & \Longleftrightarrow & \mrm{pr}_m\circ(df\circ j_L)=F^{[m]}\circ df
\end{array}
\end{equation}
}\end{rmk}
\smallskip
Then additionning these two equivalences, we recover the equivalence: "$\alpha$ solves $(\systmt)$" $\Longleftrightarrow$ "$f$ is $F^{[m]}$-holomorphic". Moreover, the equations $\alpha_{-j}''=0$, $1\leq j\leq k-1$, mean that $f$ is horizontally holomorphic.
\begin{prop}
Let $\alpha$ be a $\g$-valued 1-form on $L$ and $f$ its geometric map. The following statements are equivalent:
\begin{description}
\item[(i)]  $\alpha_{-j}''=0$, $1\leq j\leq k-1$
\item[(ii)] $f$ is horizontally holomorphic: $(df\circ j_L)^{\mal H}=F^{[k-1]}\circ df$, $\mal H=[\mk]$ being the horizontal space and ${F^{[k-1]}}_{|\mal H}$ defining a complex structure on $\mal H$.
\end{description}
\end{prop}
So that we can conclude:
\begin{thm}\label{conclusion-even-mindeter}
 The even minimal determined system $(\syst(k,\tau)$) means that the geometric map $f$ is horizontally holomorphic and vertically harmonic.
\end{thm}
\begin{rmk}\emph{
We can express what precedes in terms of the projection map $\bar\pi_{G/G_0}\colon\alpha\to f$ defined as follows. Let $\mal{MC}$ be the set of $\g$-valued integrable 1-form on $L$ and for $m<k$, $\mal{MC}^m$ the subset of interable 1-form taking values in $\g_0\oplus(\oplus_{j=1}^m\mk_j)$, then $\bar\pi_{G/G_0}\colon \mal{MC}\to C^\infty(L,G/G_0)$ is defined by:
$$
\begin{CD}
\bar\pi_{G/G_0}\colon\alpha\in\mal{MC} @>\mrm{int}>> U\in C_*^\infty(L,G)@>\pi_{G/G_0}>> f=\pi_{G/G_0}\circ U\in C^\infty(L,G/G_0)
\end{CD}.
$$
The preceding results can be summarized as follows: for any $m<k$
$$
\bar\pi_{G/G_0}(\mal{MC}^m)=C_m^\infty(L,G/G_0)\quad \mrm{and}\quad \bar\pi_{G/G_0}(\mal S(m))=\mrm{Hol}((L,j),(G/G_0,F^{[m]})),
$$
the set of $F^{[m]}$-holomorphic maps; and the equations $(H_j)$, $1\leq j\leq m$, in $\mal{MC}$ are transformed by $\bar\pi_{G/G_0}$ into the equation $\mrm{pr}_m\circ(df\circ j_L)=F^{[m]}\circ df$ in $C^\infty(L,G/G_0)$.
}\end{rmk}
$\bullet$ \emph{The even minimal determined system will be studied with much more details, precisions and results in section~\ref{vertically-harmonic}.}
\index{canonical!$f$-structure|)}\index{even case|)}\index{horizontally holomorphic|)}
%
%%%%%%%%%%%%%%%%%%%

\subsubsection{The minimal determined odd case}\label{detodcase}
\index{odd case|(}
The order of the automorphism $\tau$ is odd: $k'=2k+1$, and $m=k+1$. Let us recall the following decomposition
$$
\g^\C=\left(\g_{-k}^\C\oplus\ldots\oplus\g_{-1}^\C\right)\oplus\g_0^\C\oplus\left(\g_{1}^\C\oplus\ldots\oplus\g_{k}^\C\right).
$$
\begin{prop}\label{prop-odd-mindeter}
The equations $(S_j)$, $0\leq j\leq k-1$, of $(\syst(k+1,\tau))$ are respectively the projection on $\g_{-j}^\C$, $0\leq j\leq k-1$, of the Maurer-Cartan equation (MC) (owing to the holomorphicity conditions given by proposition~\ref{hol-hor-condition}). Hence the elliptic system $(\syst(k+1,\tau))$ can be written:
\begin{equation}\label{odd-deter}
\left\{
\begin{array}{lc}
\displaystyle \frac{}{}\alpha_{-j}''=0, \ 1\leq j\leq k-1 & (H_j)\\
\displaystyle \frac{}{} \bar\partial\alpha_{k}' + [\alpha_0''\wedge\alpha_{k}']=0 & (S_{k+1})\\
\displaystyle\bar\partial\alpha_{-k}' + [\alpha_0''\wedge\alpha_{-k}'] + [\alpha_1''\wedge\alpha_k']=0 & (S_k)\\
\displaystyle d\alpha + \frac{1}{2}[\alpha\wedge\alpha]=0  & (\mrm{MC})
\end{array}
\right. .
\end{equation}
\end{prop}
\proof
Proceed as for proposition~\ref{prop-evendeter}. \hsq \bigskip\\
Then we see that - in the presence of $(H)$ - the projection on $\g_{-k}^\C$ of (MC):
\begin{equation}\label{MCgk}
d\alpha_{-k} + [\alpha_0\wedge\alpha_{-k}] + [\alpha_1\wedge\alpha_k]=0
\end{equation}
 is nothing but $(S_k)+ (\overline{S_{k+1}})$.\\[1mm]
Now we have to distinguish two cases.%\\[1mm]
\paragraph{$\bullet$ The strictly minimal determined case\protect\footnote{That is to say the  minimal, but non maximal, determined cases.}} Let us suppose that $k\geq 2$, then we have 
$$
(S_k)\vee(S_{k+1})\equiv (S_k)+ (S_{k+1})\Longleftrightarrow \left[\bar\partial (\Ad U(\alpha_{\mk_k}'))\right]_{[\mk_k]}=0
$$
where $U$ integrates $\alpha$. For the last equivalence, just do the computation:
\begin{eqnarray*}
\Ad U^{-1} \left[\bar\partial (\Ad U(\alpha_{\mk_k}'))\right]_{[\mk_k]} & = & \bar\partial\alpha_{\mk_k}' + [\alpha''\wedge\alpha_{\mk_k}']_{\mk_k}\\
& = & \bar\partial\alpha_{\mk_k}' + [\alpha_0''\wedge\alpha_{\mk_k}'] +
[\alpha_1''\wedge\alpha_k'] + [\alpha_{-1}''\wedge\alpha_{-k}']\\
& = & (S_k)+ (S_{k+1})
\end{eqnarray*}
since $\alpha_{-1}''=0$. Hence we obtain that
$$
\begin{array}{c}
(\syst)\Longleftrightarrow
\left\{
\begin{array}{lc}
\displaystyle \frac{}{} \alpha_{-j}''=0, \ 1\leq j\leq k-1 & (H_j)\\
\displaystyle \bar\partial\alpha_{k}' + [\alpha_0''\wedge\alpha_{k}']=0 & (S_{k+1})\\
\displaystyle d\alpha + \frac{1}{2}[\alpha\wedge\alpha]=0  & (\mrm{MC})
\end{array}
\right.
\Longleftrightarrow\\
\left\{
\begin{array}{lc}
\displaystyle \frac{}{} \alpha_{-j}''=0, \ 1\leq j\leq k-1 & (H_j)\\
\displaystyle\bar\partial\alpha_{-k}' + [\alpha_0''\wedge\alpha_{-k}'] + [\alpha_1''\wedge\alpha_k']=0 & (S_k)\\
\displaystyle d\alpha + \frac{1}{2}[\alpha\wedge\alpha]=0  & (\mrm{MC})
\end{array}
\right.
\Longleftrightarrow
\left\{
\begin{array}{lc}
\displaystyle \frac{}{} \alpha_{-j}''=0, \ 1\leq j\leq k-1 & (H_j)\\
\displaystyle\left[\bar\partial (\Ad U(\alpha_{\mk_k}'))\right]_{[\mk_k]}=0 & (S_{\mk_k})\\
\displaystyle d\alpha + \frac{1}{2}[\alpha\wedge\alpha]=0  & (\mrm{MC})
\end{array}
\right.
\end{array}
$$
In terms of the geometric map $f\colon L\to G/G_0$, we have according to remark~\ref{nabla1}  the following geometric interpretation:\index{anticanonical connection|(}
$$
(S_{\mk_k})\Longleftrightarrow \bar\partial^{(\nabla^1)^v}\partial^v f=0,
$$
where the splitting $TN=\mal H\oplus \mal V$ is defined by $\mal H= [\mk']$, $\mal V=[\mk_k]$ and $\mk'=\oplus_{j=1}^{k-1}\mk_j$.
Moreover since $2\,\mrm{Re}(S_{\mk_k})$ is 
$$
d\alpha_{\mk_k} + [\alpha_0\wedge \alpha_{\mk_k} ] + [\alpha_{\mk_1} \wedge \alpha_{\mk_k} ]_{\mk_k} =0
$$
which is nothing but $[\mrm{MC}]_{\mk_k}$ (in the presence of $(H)$), the projection of (MC) on $\mk_k $, then the only new information (in addition to (MC) and ($H$)) given by the determined elliptic integrable system in the odd case is the imaginary part of $(S_{\mk_k}$) which means that $f$ is vertically harmonic (with respect to $\nabla^1$):
$$
2\,\mrm{Im}(S_{\mk_k}): d*\alpha_{\mk_k} + [\alpha_0\wedge *\alpha_{\mk_k} ] + [\alpha_{\mk_1} \wedge *\alpha_{\mk_k} ]_{\mk_k} =0 \Longleftrightarrow \tau_1^v(f):=\mrm{Tr}((\nabla^1)^v d^v f)=0.
$$
\paragraph{$\bullet$ The model case.} \index{canonical!almost complex structure|(}\index{holomorphically harmonic|(} \index{model@model case or system} 
Let us suppose that $k=1$. Let us remark that, in this situation, the determined case reduces to the (model) system $(\syst(2,\tau))$ which is then simultaneously minimal and maximal. Furthermore, coming back to (\ref{odd-deter}), we have 
$$
(S_1)\Longleftrightarrow \left[\bar\partial (\Ad U(\alpha_{\mk}'))\right]_{[\mk]}^{1,0}=0
$$
where $U$ integrates $\alpha$ and $[\ ]^{(1,0)}$ denotes the $(1,0)$-component with respect to the canonical almost complex structure $\undj$ in $N$, defined let us recall it, by the decomposition (\ref{def-odd}), i.e. in our case $TN^\C=[\g_{-1}^\C]\oplus [\g_1^\C]$. Indeed we have 
\begin{eqnarray*}
\Ad U^{-1} \left[\bar\partial (\Ad U(\alpha_{\mk}'))\right]_{[\mk]} & = & \bar\partial\alpha_{\mk}' + [\alpha''\wedge\alpha_{\mk}']_{\mk}\\
& = & \bar\partial\alpha_{\mk}' + [\alpha_0''\wedge\alpha_{\mk}'] + [\alpha_{\mk}''\wedge\alpha_{\mk}']_{\mk}\\
& = & \bar\partial\alpha_{-1}' + [\alpha_0''\wedge\alpha_{-1}'] +
[\alpha_1''\wedge\alpha_1']\\
&  + & \bar\partial\alpha_{1}' + [\alpha_0''\wedge\alpha_{1}'] + [\alpha_{-1}''\wedge\alpha_{-1}']
\end{eqnarray*}
the up term being the $(1,0)$-component and the down one, the $(0,1)$-component. Then recalling that $(S_1) + \overline{(S_2)}$ is the projection on $\g_{-1}^\C$ of (MC), we obtain that
$$
(\syst)\Leftrightarrow\left\lbrace \begin{array}{l}  (S_2)\\ (\mrm{MC}) \end{array} \right.\Leftrightarrow \left\lbrace \begin{array}{l}  (S_1)\\ (\mrm{MC})\end{array}\right.\Leftrightarrow \left\lbrace \begin{array}{lc}
\displaystyle\left[\bar\partial (\Ad U(\alpha_{\mk}'))\right]_{[\mk]}^{1,0}=0 & (S_{\mk}^{(1,0)})\\
\displaystyle d\alpha + \frac{1}{2}[\alpha\wedge\alpha]=0  & (\mrm{MC})
\end{array}
\right.
$$
and the only new information (in addition to (MC)) given by the determined elliptic integrable system in this case is $(S_\mk^{(1,0)})$.\\
In terms of the geometric map, $f\colon L\to G/G_0$, we have according to remark~\ref{nabla1} the following geometric interpretation:
$$
(S_1)\equiv (S_{\mk}^{(1,0)})\Longleftrightarrow \left[\bar\partial^{\nabla^1}\partial f\right]^{1,0}=0,
$$  
we will say that $f$ is holomorphically harmonic \wrt $\nabla^1$ and $\undj$ (see section~\ref{affineharmonic}, definition~\ref{holoharmdef}, for a precise definition).\\
In the same way, we also have  the following (equivalent) geometric interpretation\index{canonical!connection, $G$-invariant}
$$
(S_2)\Longleftrightarrow \left[\bar\partial^{\nabla^0}\partial f\right]^{0,1}=0,
$$  
we will say that  $f$ is holomorphically harmonic\footnote{See also section~\ref{affineharmhom}, theorem~\ref{thm3sym}.} \wrt $\nabla^0$ and $-\undj$\medskip.\\
Moreover, let us write the equations of the system as a real equation and then write the corresponding geometric  equation (which will then take place in $TN$ and not in $TN^\C$). We have that  the real equation $(S_1) + \overline{(S_2)}$ is the projection on $\g_{-1}^\C$ of (MC). Now, let us write  the  equation 
$$
(E_{-1})\equiv  \dfrac{(S_1) - \overline{(S_2)}}{i}\equiv \left( d*\alpha_{-1} + [\alpha_0\wedge *\alpha_{-1}] + \dfrac{1}{2i}[\alpha_1\wedge\alpha_1]=0 \right), 
$$
then taking the sum with its complex conjugate, we obtain
$$
\begin{array}{rcl}
(E_{-1}) + \overline{(E_{-1})} & \equiv  & \left(  d*\alpha_{\mk} + [\alpha_0\wedge *\alpha_\mk] + \dfrac{1}{2i}( [\alpha_1\wedge\alpha_1] - [\alpha_{-1}\wedge\alpha_{-1}])=0\right)\\
& \equiv &  \left(  d*\alpha_{\mk} + [\alpha_0\wedge *\alpha_\mk] - \dfrac{1}{2}\undj_0[\alpha_\mk\wedge\alpha_\mk]=0\right),
\end{array} 
$$
since $ \frac{1}{2i}( [\alpha_1\wedge\alpha_1] - [\alpha_{-1}\wedge\alpha_{-1}])= \frac{1}{2}[\undj_0\alpha_\mk\wedge\alpha_\mk]=-\frac{1}{2}\undj_0[\alpha_\mk\wedge\alpha_\mk]$.
Then written in terms of the geometric maps $f\colon L\to N$, the last equation means\index{tension field}
$$
\tau^0(f) + \undj T^0(f)=0,
$$
where $\tau^0(f)=\mrm{Tr}_g(\nabla^0 d f)$ is the tension field of $f$ \wrt $\nabla^0$ (and $g$  is some Hermitian metric on $L$), and $T^0(f)=*(f^*T^0)$, $T^0$ being the torsion of $\nabla^0$. As we will see in section~\ref{holharmmap}, this  equation is in fact a general characterization for holomorphically harmonic maps. \index{holomorphically harmonic|)}
\paragraph{The primitive case.} The $m$-primitive case is obtained by putting, in the minimal determined case (\ref{odd-deter}), $\alpha_k'=0$, if $m=k$, and  $\alpha_j=0$, $m+1\leq |j|\leq k$, if $m\leq k-1$. As in the even case we obtain:
\begin{prop}
Let $\tau\colon\g\to\g$ be an order $2k+1$ automorphism, and $m\leq k$ a positive integer then the $m$-elliptic integrable system $(\systmt)$ means that the geometric map $f\colon L\to G/G_0$ satisfies:
$$
\partial f\in \oplus_{j=1}^m[\g_{-j}^\C]\subset T(G/G_0)^\C .
$$
\end{prop}
%
%%%%%%%%%
%
\paragraph{Geometric interpretation of the equations $(H_j)$.}
Recall that $\undj$ denotes the canonical complex structure on $N=G/G_0$ (see (\ref{def-odd})) and set  $F^{[m]}:=\mrm{pr}_m\circ \undj=\undj\circ\mrm{pr}_m$ for $m\leq k$, where $\mrm{pr}_m\colon TN\to \oplus_{j=1}^m[\mk_j]$ is the projection on $\oplus_{j=1}^m[\mk_j]$ along $\oplus_{j\geq m+1}^k[\mk_j]$ (remark that $\mrm{pr}_k=\Id$). Then $F^{[m]}$ is an $f$-structure on $N$ (remark that $F^{[k]}=\undj$ is a complex structure). Then we have:
\begin{thm}\index{J holomorphic@$J$-holomorphic}\index{f holomorphic@$f$-holomorphic}
A map $f\colon L\to G/G_0$ is $m$-primitive \iif it is $F^{[m]}$-holomorphic. In particular, $f$ is $k$-primitive \iif it is
holomorphic (with respect to the canonical almost complex structure on $G/G_0$), and thus any $m$-primitive map is in particular a holomorphic curve in $G/G_0$. More precisely,  $m$-primitive maps are exactly the integral holomorphic curves of the complex Pfaff system $\oplus_{j=1}^m[\mk_j]\subset TN$.
\end{thm}
\begin{rmk}\emph{
The equivalences (\ref{equiv-m}) hold also in the odd case. However for $m=k$, the first equivalence of (\ref{equiv-m}) is trivial: $\alpha\in\g\Longleftrightarrow f\in C^\infty(L,G/G_0)$. There is no restriction (in the form "$\alpha$ takes values in a subspace of $\g$") in the highest primitive case.
}\end{rmk}
\begin{prop}
Let $\alpha$ be a $\g$-valued 1-form on $L$ and $f$ its geometric map. The following statements are equivalent:
\begin{description}
\item[(i)]  $\alpha_{-j}''=0$, $1\leq j\leq k-1$
\item[(ii)] $f$ is horizontally holomorphic: $(df\circ j_L)^{\mal H}=F^{[k-1]}\circ df$, where $\mal H=[\mk']$ is the horizontal space, and ${F^{[k-1]}}_{|\mal H}=\undj_{|\mal H}$ defines a complex structure on $\mal H$.
\end{description}
\end{prop}
We can conclude:
\begin{thm}\label{conclusion-odd-mindeter}\index{holomorphically harmonic}
The odd minimal determined system $(\syst(k+1,\tau)$) means that the geometric map $f$ is horizontally holomorphic and vertically harmonic \wrt $\nabla^1$ if $k\geq 2$; and  it means that $f$ is holomorphically harmonic \wrt $\nabla^1$, or equivalently anti-holomorphically harmonic \wrt $\nabla^0$ if $k=1$.
\end{thm}
\begin{rmk} \em The odd minimal determined system will be studied with much more details, precisions and results in section~\ref{subsec-odd-vert-hol-harm-Pfaff}.
\end{rmk}
\index{determined minimal@determined, minimal|)}\index{vertically harmonic|)}\index{primitive@primitive, case or system|)}\index{vertical tension field|)}\index{anticanonical connection|)}
\subsection{The maximal determined case}\label{eq-maximal-deter}
\index{determined maximal@determined, maximal}
Let us see how can be rewritten the elliptic system in this case in more geometric terms. Let us recall that  the determined system can be considered as a system on the 1-form $\alpha$. 
\begin{thm}\label{thm-max-deter-odd}
Let $\tau\colon \g\to\g$ be an automorphism of odd order $k'=2k+1$. Let us set $J_0=\taum$ and let $\undj_0$ be the corresponding complex structure on $\mk$ i.e. the value of $\undj$ at the reference point $y_0=1.G_0\in N$ (see equation (\ref{def-odd})). Then the associated maximal determined system, $\syst(k'-1,\tau)$, is equivalent to
$$
\left\lbrace \begin{array}{l}
\displaystyle d*\alpha_\mk + [\alpha_0\wedge *\alpha_\mk] + \dfrac{1}{2}[\,\undj_0 \alpha _\mk\wedge\alpha_\mk]_\mk + \sum_{1\leq i< j\leq k}[\,\undj_0 \alpha _{\mk_i}\wedge\alpha_{\mk_j}]_{\mk_{j-i}}=0 \quad (E_\mk)\\
d\alpha + \dfrac{1}{2}[\alpha\wedge\alpha] =0 \quad (\mrm{MC})
\end{array}\right.
$$
\end{thm}
\index{canonical!almost complex structure|)}
\begin{thm}\label{thm-max-deter-even}\index{canonical!$f$-structure}\index{even case}
Let $\tau\colon \g\to\g$ be an automorphism of even order $k'=2k$. Let us set $J_0=\taum$ and let $\undj_0$ be the corresponding complex structure on $\mk$ i.e. the value of $F_{|\hor}$ at the reference point $y_0=1.G_0\in N$  (see equation (\ref{def-even})). Then the associated maximal determined system, $\syst(k'-1,\tau)$, is equivalent to
$$
\left\lbrace \begin{array}{l}
\displaystyle d*\alpha_k + [\alpha_0\wedge *\alpha_k] + \dfrac{1}{2}[\undj_0 \alpha_\mk\wedge\alpha_\mk]_{\g_k}=0 \quad (E_k)\\
\displaystyle d*\alpha_\mk + [\alpha_0\wedge *\alpha_\mk] + \dfrac{1}{2}[\,\undj_0 \alpha _\mk\wedge\alpha_\mk]_\mk + \sum_{1\leq i< j\leq k}[\,\undj_0 \alpha _{\mk_i}\wedge\alpha_{\mk_j}]_{\mk_{j-i}} + [\alpha_k\wedge \undj_0 \alpha_\mk]=0 \quad (E_\mk)\\
d\alpha + \dfrac{1}{2}[\alpha\wedge\alpha] =0 \quad (\mrm{MC})
\end{array}\right.
$$
\end{thm}
\textbf{Proof of theorems~\ref{thm-max-deter-odd} and \ref{thm-max-deter-even}.}
This is a straightforward computation. Indeed, it suffices to check that we have
$$
\left\lbrace\begin{array}{l}
\displaystyle \left[(\mrm{MC})\right]_{\g_0}= (S_0)\\
\displaystyle \left[(\mrm{MC}) \right]_{\g_{-j}}= (S_j)  + \ovr{(S_{k'-j})}, \quad 1\leq j \leq k'-1\\
\displaystyle \left[(E) \right]_{\g_{-j}}= \dfrac{(S_j)  - \ovr{(S_{k'-j})}}{i}, \quad 1\leq j \leq k'-1
\end{array}\right. 
$$
where $(E)=(E_\mk)$ in the odd case and $(E)=(E_\mk) + (E_k)$ in the even case. \comprf\hsq
\begin{rmk}\em
We will see, in sections~\ref{affineharmonic} and \ref{gene-Harm-f-structure} respectively, that the equation $(E)$ means that the geometric  map $f$ is stringy harmonic.
\end{rmk}
\subsubsection*{Adding holomorphicity conditions; the intermediate determined systems.}
According to proposition~\ref{hol-hor-condition},  we know that the determined system $(\systm)$, $m_{k'}\leq m\leq k'-1$, is obtained by adding the holomorphicity conditions, $\alpha_{-j}''=0$,   $1\leq j\leq k'-1-m$, to the maximal determined system, whose equations are given respectively by theorems~\ref{thm-max-deter-odd} and \ref{thm-max-deter-even}. Therefore, using the notations defined in the proof of theorems~\ref{thm-max-deter-odd} and \ref{thm-max-deter-even}, we have the following.
\begin{prop}\label{prop-intermediate-syst}
Let $\tau\colon \g\to\g$ be an automorphism of  order $k'$. Let $m$ be an integer such that $m_{k'}\leq m\leq k'-1$, and let us set $\und m= k'-1 - m$. Then the determined system $(\systm)$ is equivalent to
$$
\begin{array}{rcl}
\left\lbrace \begin{array}{l}
\displaystyle
\bar\partial \alpha_{-j}' + \sum_{i=0}^{m-j} [\alpha_i''\wedge\alpha_{-i-j}'] = 0 \quad (S_j), \quad \und m + 1 \leq j \leq m \\
d\alpha + \dfrac{1}{2}[\alpha\wedge\alpha] =0 \quad (\mrm{MC})\\
\alpha_{-j}''=0  \quad (H_j), \quad   1\leq j\leq \und m
\end{array}\right.
&  \Longleftrightarrow & 
\left\lbrace \begin{array}{l}
\displaystyle
[(E)]_{\pk} \\
d\alpha + \dfrac{1}{2}[\alpha\wedge\alpha] =0 \quad (\mrm{MC})\\
\alpha_{-j}''=0  \quad (H_j), \quad   1\leq j\leq \und m  
\end{array}\right.
\end{array}
$$
We have set $\mk' =\oplus_{j=1}^{\und m} \mk_j$ and $\pk=\oplus_{j=\und m + 1}^k \mk_j$, where in the even case $\mk_k:=\g_k$. Then $[(E)]_{\pk}$ denotes the projection of $(E)$ on $\pk$ according to the decomposition $\mk=\mk' \oplus\pk$ in the odd case, and $\nk=\mk'\oplus\pk$ in the even case.
\end{prop}
\proof Consider $(\systm)$  as a subsystem of $\syst(k'-1)$. Then we already have  seen in the proof of theorems~\ref{thm-max-deter-odd} and \ref{thm-max-deter-even} that $\left[(\mrm{MC})\right]_{\g_0}= (S_0)$, $\left[(\mrm{MC}) \right]_{\g_{-j}}= (S_j)  + \ovr{(S_{k'-j})}$,  $1\leq j \leq m$ and $ \left[(E) \right]_{\g_{-j}}= \dfrac{(S_j)  - \ovr{(S_{k'-j})}}{i}$, $1\leq j \leq m$. Then it suffices to check that  the holomorphicity conditions, i.e. $\alpha_{-j}''=0$,   $1\leq j\leq \und m$, imply that the equation $(S_j)$, for $1\leq j \leq \und m$, is nothing but the projection on $\g_{-j}$ of the Maurer-Cartan equation. This can be made by computation, with the same method as in the proof of proposition~\ref{prop-evendeter}. \comprf \hsq
\index{odd case|)}\index{holomorphicity condition|)} 
\subsection{The underdetermined case}\label{underdetermined}
\index{underdetermined@underdetermined, case or system}
Now, we prove that any underdetermined system can be written as a determined system in a new setting.
\begin{thm}
Let us consider an underdetermined system $(\syst(m,\g,\tau))$, $m \geq k'$. Let us write
$$
m=qk' + r,\quad 0\leq r \leq k'-1
$$
the Euclidean division of $m$ by $k'$. Then let us consider the automorphism in $\g^{q+1}$ defined by
$$
\tl\tau \colon (a_0,a_1,\ldots,a_{q})\in\g^{q+1}\longmapsto (a_{1}, \ldots,a_{q},\tau(a_0))\in\g^{q+1}.
$$
Then $\tl\tau$ is of order $(q+1)k'$. Moreover the  $m$-th system associated to $(\g,\tau)$ is in fact equivalent to the $m$-th system associated to $(\g^{q+1},\tl\tau)$. More precisely, denoting by $\tl\omega$ a $(q+1)k'$-th primitive root of  unity\footnote{Chosen such that $\tl\omega^{q+1}=\omega_{k'}$}, then the map 
$$
\alpha_\lm \longmapsto (\alpha_\lm,\alpha_{\tl \omega \lm},\ldots,\alpha_{\tl\omega^q\lm})
$$
is  a bijection from the set of solutions of the \textbf{underdetermined} $(m,\g,\tau)$-system into the set of solutions of the \textbf{determined}  $(m,\g^{q+1},\tl\tau)$-system. 
\end{thm}
\proof
It suffices to check that $\alpha_\lm$ satisfies the properties in equation~\eqref{Sm} \iif $\tl\alpha_\lm:=(\alpha_\lm,\alpha_{\tl \omega \lm},\ldots,\alpha_{\tl\omega^q\lm})$ do so (in the setting $(\g^{q+1},\tl\tau)$), which is immediate. \hsq \medskip\\ 
This theorem tells us that the study of the underdetermined case reduces to that of the determined case.\\
Moreover this theorem provides a new point of view about the projection map $\pi_m$. It tells us that by putting together several exemplar of $\pi_m$, to obtain some $p$-uplet taking values in a $p$-power of $\mal G_1^m$, then we construct a diffeomorphism (remarking that if $\alpha_\lm$ is a solution then $\tl\alpha_\lm:=\alpha_{t\lm}$ is also a solution, for any $t\in S^1$).
\subsection{Examples}\label{Examples}
\subsubsection{The trivial case: the 0-th elliptic system associated to a Lie group.}
We consider the determined system $(\syst(m,\tau))$ with $\tau=\Id$ and (thus) $k'=1$ so that $m_k'=m_1=0=k'-1$. Then the determined system $(\syst(0,\Id))$ is nothing but the Maurer-Cartan equation for $\g$-valued 1-form $\alpha$ (i.e. in other words the "equation" for the trivial geometric map $f\colon L\to G/G=\{1\}$).
\subsubsection{Even determined case}\index{even case}
\paragraph{The first elliptic system associated to a symmetric space \cite{DPW}.}\index{harmonic map}
We consider the even determined system $(\syst(k,\tau))$, with $k=1$ and $\tau$ an involution. Then the horizontal subbundle is trivial $\hor=[\mk]=\{0\}$ and $TN=[\ver]=[\g_1]$ so that the horizontal holomorphicity is trivial and vertical harmonicity means harmonicity. Hence the first elliptic system associated to a symmetric space, $(\syst(1,\tau))$, \textbf{is the equation for harmonic maps} $f\colon L\to G/G_0$.
\paragraph{The second elliptic system associated to a 4-symmetric space (\cite{ki3,bk}).}
\index{vertically harmonic}
Here $\tau$ is an order four automorphism and (thus) $k=2$. Then we consider the even determined system $(\syst(2,\tau))$. A geometric interpretation in terms of vertically harmonic twistor lifts of this system is given in \cite{ki3} (see also \cite{bk}). Let us give some examples.
\begin{description}
\item[$\bullet$]  Surfaces with holomorphic mean curvature vector in 4-dimensional symmetric spaces (\cite{bk}).
\item[$\bullet$]  Hamiltonian stationary Lagrangian surfaces in Hermitian symmetric spaces (\cite{HR3,HR1,HR2}).
\item[$\bullet$]  Surfaces with holomorphic mean curvature vector in 4-dimensional spaces forms (\cite{bk}).
\item[$\bullet$]  surfaces with anti-holomorphic mean curvature vector in $\C P^2$ (\cite{bk}).
\item[$\bullet$]  $\rho$-harmonic surfaces in $\oct$ \cite{ki1}.
\end{description}
The second, third and fourth examples are particular cases of the first one (see \cite{bk}).
\subsubsection{Primitive case}\index{primitive@primitive, case or system}
Let us give some examples where the automorphism $\tau$ is of order $k'\geq 3$, and  $m=1$.  
\begin{description}
\item[$\bullet$] Minimal surfaces in $\C P^2$. $(m=1,k'=3, N=SU(3)/S(U(1)^3)\,)$. (\cite{bpw},\cite{harmonicTori}).
\item[$\bullet$] Minimal Lagrangian surfaces in $\C P^2$. $(m=1,k'=6,N=SU(3)/S(U(1)^2)\, )$. (\cite{tern}).
\item[$\bullet$] Special Lagrangian surfaces in 4-dim Hermitian symmetric spaces. $(m=1,k'=4)$. (\cite{HR3,HR1,HR2})
\item[$\bullet$] Constant mean curvature (CMC) surfaces in $H^3$: their Gauss maps are primitive maps into the unit tangent bundle which is a 4-symmetric space. (\cite{dik})
\item[$\bullet$] Affine Toda fields or Toda lattice. $(m=1,k'\geq 3)$. (\cite{bpw},\cite{12}).
\item[$\bullet$] Non-superminimal (weakly) conformal harmonic maps into $S^n$. $(m=1, k'=2r+2, N=F^r(S^n)\,)$.  (\cite{harmonicTori}).
\item[$\bullet$] (Weakly) conformal non-isotropic harmonic maps into $\C P^n$.  $(m=1,k'=r+2,N=F^r(\C P^n)\,)$. (\cite{harmonicTori}).
\end{description}
$F^r(S^n)$ denotes some bundle of isotropic flags over $S^n$ , and $F^r(\C P^n)$ denotes some bundle of  flags over $\C P^n$ (see \cite{harmonicTori}).

\subsubsection{Underdetermined  case}
\index{underdetermined@underdetermined, case or system}
\paragraph{First elliptic integrable system associated to a Lie group (\cite{Uhlenbeck, tern}).}
We consider the system $(\syst(m,\tau))$ with $m=1$ and $\tau=\Id$. Therefore $k'=1$ and $m_1=0<m$, and thus this is an underdetermined system. Then $(\syst(1,\Id))$ \textbf{is the equation for harmonic maps into the Lie group} $G$, $f\colon L\to G$.
\paragraph{Second elliptic integrable system associated to the symmetric space $Gr_{3,1}(\R^{n+1,1})$.}
Constrained Willmore surfaces in $S^n$ corresponds to particular solutions of this system. (See \cite{BuPedPink}, \cite{quintito}).
\index{m elliptic integrable system@$m$-th elliptic integrable system|)}\index{k symmetric space@$k'$-symmetric space|)}
\subsection{Bibliographical remarks and summary of the results.}
The previous list of examples given in §\ref{Examples} is a complete list of all known cases. Consequently, all the results of the present section~\ref{melliptic} about the elliptic integrable systems are completely new. More generally all the results about the elliptic integrable systems contained in the present paper (sections~\ref{melliptic}--\ref{affineharmhom}) are completely new.\\
We would like to mention some bibliographic references about $k$-symmetric spaces: Kowalski \cite{kowalski},  Wolf and Gray \cite {wolf}.

%
%%%%%
%%%
%%%%%%%%%%%%%%%%%%%%%%%%%%%%%%%%%%%%%%%%%%%%%%%%%%%%%%%%%%%%%%%%%%%%%%
%
%
%         Finite order isometries and Twistor spaces
%
%
%
%
%%%%%%%%%%%%%%%%%%%%%%%%%%%%%%%%%%%%%%%%%%%%%%%%%%%%%%%%%%%%%%%%%%%%%%%%%
%

%
\section{Finite order isometries and Twistor spaces}\label{isometry-twistor}
\index{twistor|(}\index{even case|(}\index{odd case}
The aim of this section is to realise any $k'$-symmetric space as the subbundle of some bundle of endomorphisms. To make things more precise we need first to introduce some notations.\medskip\\
Let $E$ be an Euclidean space and let us define (for $p\in \mathbb N^*$)
$$
\begin{array}{c}
\mal{U}_p(E)=\{A\in SO(E),\, A^p=\Id,\, A^i\neq \Id \text{ if } 1\leq i < p\}\\
\mal{U}_p^*(E)=\{A\in\mal{U}_p(E)|1\notin\mrm{Spect}(A)\},\quad \mal{U}_p^{**}(E)=\{A\in\mal{U}_p(E)|\pm 1\notin\mrm{Spect}(A)\}.
\end{array}
$$
Then for $k\in\mathbb N^*$ we set
$$
\mal{Z}_{2k}(E)=\mal{U}_{2k}^{**}(E)\quad \text{and}\quad \mal{Z}_{2k+1}(E)=\mal{U}_{2k+1}^{*}(E)=\mal{U}_{2k+1}^{**}(E).
$$
Let us now explain the general spirit  of this section.
There are two main ideas which will allow us to provide a twistor interpretation of our elliptic integrable systems. These are the two ideas that we will follow in our exposition.
\paragraph{1st idea.}
A $2k$-symmetric  homogenous space has a very particular geometry. Then we want to write it more universally by embedding it in a more universal object. To do that, we prove that any $2k$-symmetric space $G/G_0$ can be embedded canonically in the twistor bundle  $\zdk(M)\to M$ of the associated $k$-symmetric $M=G/H$ via a morphism of bundles over $M$ (recall that $G/G_0=G\times_H(H/G_0)\to G/H$ is a bundle over $G/H$ with a symmetric fibre $H/G_0$).  Therefore  to any geometric map $f\colon L\to G/G_0$ corresponds a map in the twistor space $J\colon L\to \zdk(G/H)$.
\paragraph{2nd idea.}
Conversely, we write each connected component of $\zdk(\rdn)$ as a $2k$-symmetric homogeneous space $\zdk^\alpha(\rdn)=SO(2n)/K$, where $K$ is some Lie subgroup of $SO(2n)$. Indeed for all $J\in \zdk(\rdn)$, the inner automorphism $\Int J$ is of finite order $2k$ and hence  defines a $2k$-symmetric space which is nothing but the connected component of $J$ in $\zdk(\rdn)$. In particular the corresponding Lie algebra automorphism $\Ad J\colon \so(2n)\to \so(2n)$ gives rise to an eigenspace decomposition $\so(2n)^\C=\oplus_{j\in\Z_{2k}} \so_j^\C(2n)$ analogous to the decomposition $\g^\C=\oplus_{j\in\Z_{2k}}\g_j^\C$. Therefore, the geometric properties of a map $f\colon L\to G/G_0$ can be translated in terms of the corresponding twistor map $J\colon L\to \zdk(G/H)$.\medskip\\
We are then led to study the space of endomorphisms $\zdk(\rdn)$ which we do in the subsection~\ref{order 2k}. The guidelines of this study is to generalize the space $\Sigma(\rdn):=\mal Z_4(\rdn)$ of  (orthogonal) almost complex structures. The latter is a homogeneous symmetric space of which the tangent space at an element $J$ is the subspace of elements in $\so(2n)$ which anticommute with $J$. We will prove  analogous properties in $\zdk(\rdn)$.\\
In a second time (section~\ref{can-sect-zdk}), we will prove the embedding $G/G_0\hookrightarrow \zdk(G/H)$.\medskip\\
Even if we will use the twistor space and the embedding describred above especially in the even case, however we will sometimes need these in the odd case $k'=2k+1$. Moreover, it will also happen in our study that we will need to consider the spaces $\left( \zdk(M)\right)^j=\{ J^j,\, J\in \zdk(M)\}$ which are in  general twistor spaces $\mal U_p(M)$ for some $p\in\mathbb N$. Therefore, we will also study these spaces. Nevertheless, the even case, $\zdk(M)$, is the most rich and most complex. Indeed it satisfies all the properties which hold in the general case $\mal U_p(M)$ (in particular in the odd case $\mZ_{2k+1}(M)$), so that the study of $\zdk(M)$ contains the one of $\mal U_p(M)$. Moreover,  the even case $\zdk(M)$ admits special properties and  features which need to be studied separately. This is why we will concentrate our study on this case, in a first time. Then in a second time we will make precise which properties hold in  general and which ones are particular to the even case.
\subsection{Isometries of order $2k$ with no eigenvalues $=\pm 1$}\label{order 2k}
In this section~\ref{order 2k}, we will  study of $\mal{Z}_{2k}(E)$, where again $E$ is an Euclidean space.\\
For each $A\in\mal{Z}_{2k}(E)$, we have the following eigenspace decomposition:
$$
E^\C=\oplus_{j=1}^{k-1}\left(E_A(\omega_{2k}^j)\oplus E_A(\omega_{2k}^{-j})\right)
$$
with $E_A(\lm)=\ker(A-\lm\Id)$ and $\omega_{2k}=e^{i\pi/k}$.\\
Let us set $\mk_j^\C=E_A(\omega_{2k}^j)\oplus E_A(\omega_{2k}^{-j})$ for $j\geq 0$. Then we have
$\dim_\R\mk_j=\displaystyle\frac{1}{2}\left(\dim_\R E_A(\omega_{2k}^j) +  \dim_\R E_A(\omega_{2k}^{-j})\right)= \dim_\R E_A(\omega_{2k}^j)=2\dim_\C E_A(\omega_{2k}^j)$.
Hence $\dim_\R\mk_j$  is even and hence we will suppose now that $E=\R^{2n}$ (in all the section~\ref{order 2k}).
\begin{exam}\emph{
 We have $\mal{Z}_2(E)=\emptyset$, and $\mal{Z}_4(E)=\Sigma (E)$ the set of almost complex structures in $E$.}
\end{exam}
\paragraph{Situation in the Euclidean plane}
Here $E=\R^2$, and any element of $A\in\mal{Z}_{2k}(E)$ is written in the form $A=R(\frac{l\pi}{k})$, with $(l,2k)=1$, $R(\theta)$ being the  rotation of angle $\theta\in\R/2\pi\Z$. In other words, considered as a complex number, $A$ is a primitive $(2k)$-th root of the unity. Hence $\mrm{card}(\mal{Z}_{2k}(\R^2))=\phi(2k)$, $\phi$ being the Euler phi-function.
\subsubsection{The set of connected components in the general case}\label{setconnectcomp}
\begin{thm}\label{thm-connectcomp}
$\pi_0(\mal{Z}_{2k}(\R^{2n}))$, the set of connected components of $\mal{Z}_{2k}(\R^{2n})$, is (in one to one correspondance with):
$$
X_{2k}:=\left\{(\varepsilon,p)\in\Z_2\times\mathbb N^{k-1}\left |\,\sum_{j=1}^{k-1}p_j=n\ \mrm{and}\ \mrm{lcm}\left(\left\{ \frac{2k}{2k\wedge j}, p_j\neq 0 \right\}\right)=2k \right.\right\}
$$
\end{thm}
\textbf{Proof.} Let $A\in\mal{Z}_{2k}(\R^{2n})$, then $A_{|\mk_j^\C}= \omega_{2k}^j\Id_{E_A(\omega_{2k}^j)}\oplus\omega_{2k}^{-j} \Id_{E_A(\omega_{2k}^{-j})}$. We choose an orientation on each $\mk_j$ (such that the induced orientation on $\oplus_1^{k-1}\mk_j$ is the one of $\R^{2n}$). Then there exist oriented planes $P_j^l$ such that $\mk_j=\oplus_{l=1}^{p_j}P_j^l$ (sum of non oriented spaces), where $p_j=\frac{\dim\mk_j}{2}$, and
$$
 A_{|\mk_j}=\oplus_{l=1}^{p_j} R_{P_j^l}(\theta_j)
$$
where $R_{P_j^l}(\theta_j)$ is the rotation on $P_j^l$ of angle $\theta_j=\frac{j\pi}{k}$. Let $\varepsilon_j$ be the  sign of the orientation\footnote{see remark~\ref{orientation}} of $\oplus_{l=1}^{p_j}P_j^l$ (sum of oriented spaces) \wrt the one of $\mk_j$. Now let us consider the map
$$
f\colon A\in\mal Z_{2k}(\R^{2n})\mapsto (\Pi_{j=1}^{k-1}\varepsilon_j,(p_j)_{1\leq j\leq k-1})\in X_{2k}.
$$
Then it is a continuous\footnote{By lower semicontinuity of the rank, and the fact that $\sum_1^{k-1}p_j$ is constant.} surjection and $f^{-1}(\{(\eps,\und p)\})$ is an $SO(2n)$-orbit of the  $SO(2n)$-adjoint action on $\mal Z_{2k}(\R^{2n})$. This completes the proof.\hfill$\square$
\begin{rmk}\label{orientation}\em
Let $P$ be an oriented plane and $A=R_P(\theta)\setminus \{\pm \Id_P\}$ a rotation of this plane. Then we can define
$\eps_P(A):=\mrm{sign} \det_P(x,Ax)$ which is well defined because independent of $x \in P$. We have set $\det_P :=\det_{(e_1,e_2)}$, for any direct orthonormal basis $(e_1,e_2)$ of $P$.\\
 Now, let $E$ be an oriented Euclidean vector space of  dimension $2n$ and let $E=\oplus_{l=1}^q P_l^0$ be a decomposition of $E$ as a sum of oriented planes (endowed with the  orientation induced by $E$).  Let $A\in \oplus_{l=1}^q R_{P_l^0}(\theta_l)\in SO(E)\setminus \{ \oplus_{l=1}^q (\pm \Id_{P_l^0}) \}$. We can define
$$
\eps(A)=\eps_{P_1^0}(A_{|P_1^0})\times \cdots \times\eps_{P_q^0}(A_{|P_l})= \mrm{sign} \det_{\mal B}(x_1, Ax_1,\cdots, x_q,A x_q)
$$
where $\mal B$ is any  direct orthonormal basis of $E$, and $x_l\in P_l^0$.\\
In the other hand  let $ P_l$ be the plane $P_l^0$ endowed with the unique orientation such that $\theta_l\in ]0,\pi[$. Then we have
$$
\eps(A)=\mrm{orient}_E\left(\oplus_{l=1}^q P_l \right),
$$
 the right hand side meaning the sign of the orientation of $\oplus_{l=1}^q P_l$ \wrt to the one of $E$.

\end{rmk}
\begin{rmk}\label{continuous}\em
Each connected component is a $SO(2n)$-orbit and thus is compact. Moreover each connected component is an open and closed submanifold of $\zdk(\rdn)$ (which is itself a compact submanifold in $SO(2n)$).
\end{rmk}
\begin{defn} We will denote by $\zdk^\alpha(\rdn)$ (and sometimes only by $\zdk^\alpha$) the connected component $f^{-1}(\{\alpha\})$, for $\alpha=(\eps,\und p)\in X_{2k}$. We define
\begin{eqnarray*}
\zdk^0(\rdn) & = & \left\{ A\in\zdk(\rdn)| A^k = -\Id\right\}=\bigsqcup_{\{\alpha|\forall j, p_{2j}=0\}}\zdk^\alpha\\
\zdk^*(\rdn) & = & \left\{ A\in\zdk(\rdn)| A^k \neq -\Id\right\}=\bigsqcup_{\{\alpha|\exists j, p_{2j}\neq 0\}}\zdk^\alpha
\end{eqnarray*}
\end{defn}
$\zdk^0(\rdn)$ is the union of \textbf{order $k$ components} in $\zdk(\rdn)$, and $\zdk^*(\rdn)$ is the union of \textbf{order $2k$ components} in $\zdk(\rdn)$ (see below for the meaning of this terms).\\
In the following we will denote by $\zdk^a(\rdn)$, for $a\in\{0,*\}$, any of the two spaces $\zdk^0(\rdn)$ and $\zdk^*(\rdn)$, and $r$ the order of these two spaces i.e. $r=\begin{cases} 2k & \text{ in } \zdk^*(\rdn)\\
 k & \text{ in } \zdk^0(\rdn)\end{cases}$. $r$ is in fact the order of $\Ad J$, for $J\in\zdk(\rdn)$ (see \ref{AdJ} below). Let us compute the tangent space of $\zdk(\rdn)$: $\forall J\in\zdk(\rdn)$,
\begin{equation}\label{TJ}
T_J\zdk(\rdn)=\left\{A\in J.\so(2n)\left|\sum_{p+l=2k-1} J^pAJ^l=0\right.\right\}
\end{equation}
and for $J\in\zdk^0(\rdn)$, we have in addition
\begin{equation}\label{TJ0}
T_J\zdk(\rdn)=T_J \zdk^0(\rdn)=\left\{A\in J.\so(2n)\left|\sum_{p+l=k-1} J^pAJ^l=0\right.\right\}
\end{equation}
It could seem strange that the two expressions (\ref{TJ}) and (\ref{TJ0}) are equal for $J\in \zdk^0(\rdn)$, but as we will see below, it comes from the fact that the "even" eigenspaces of $\Ad J$ vanish, for $J\in\zdk^0(\rdn)$, which leads to this last equality (which is in general an inclusion "$\supset$")
\begin{exam}
\em If $k=2$, then $X_{2k}=\{\pm 1\}=\Z_2$ and $\mal Z_4(\rdn)=\mal Z_4^0(\rdn)=\Sigma(\rdn)=\{J\in SO(\rdn)|J^2=-\Id
\}=\Sigma^+(\rdn)\bigsqcup\Sigma^-(\rdn)$ (resp. the positive and negative components of $\Sigma(\rdn)$), whereas $\mal Z_4^*(\rdn)=\emptyset$.
\end{exam}
\subsubsection{Study of $\Ad J$, for $J\in\zdk^a(\rdn)$}\label{AdJ}
Let $J\in\zdk^a(\rdn)$. $\Ad J$ is then an order $r$ automorphism of $\mrm{End}(\rdn)$ (since $(\Ad J)^p=\Id\Leftrightarrow J^p=\pm \Id$) thus we have the following eigenspaces decomposition:
$$
\mrm{End}(\rdn)^\C=\bigoplus_{j\in\Z/r\Z}\ker(\Ad J-\omega_r^j\Id)
$$
with $\omega_r=e^{2i\pi/r}$. Let us set 
$$
\mal A_j^\C(J)=\ker(\Ad J -\omega_r^j\Id).
$$
Then $\mal A_0(J)=\mrm{Com}(J):=\{A\in\mrm{End}(\rdn)| [A,J]=0\}$ and for $j\neq 0$ we have: $\forall A\in\mal A_j^\C $,
$$
\sum_{l+p=r-1} J^l A J^p=\sum_{l+p=r-1}(\omega_r^j)^l A J^{p+l}=\left[\sum_0^{l-1}(\omega_r^j)^l\right]J^{r-1}=0.
$$
Hence
$$
\bigoplus_{j\in\Z/r\Z\setminus \{0\}}\mal A_j^\C(J)\subset\ker\left(\sum_{l+p=r-1}L(J^l)\circ R(J^p)\right)
$$
where $ L,R$ denotes respectively the left and right multiplications. This inclusion is in fact an equality. Indeed, let $A\in\mrm{End}(\rdn)^\C$, then $A=\sum_{j=0}^{r-1}A_j$, with $A_j\in\mal A_j^\C(J)$, thus $\sum_{j=0}^{r-1}J^l A J^{r-1-l}=rA_0J^{r-1} + 0=rA_0J^{r-1}$ which vanishes if and only if $A_0=0$. This proves:
\begin{prop}\label{sum-Aj}
The following equality holds 
$$
\bigoplus_{j\in\Z/r\Z\setminus \{0\}}\mal A_j^\C(J)=\ker\left(\sum_{l+p=r-1}L(J^l)\circ R(J^p)\right).
$$
\end{prop}
Now, let us restrict ourself to $J.\so(2n)$, resp. to $\so(2n)$, (which does not change the order of $\Ad J_{|J.\so(2n)}$, resp.  $\Ad J_{|\so(2n)}$) and set\footnote{$\mal B_j^\C(J)$ is stable by $\Ad J$ and we have $\mal B_j^\C(J)=J.\so_j^\C(J)=\so_j^\C(J).J$. Besides we have more generally $J.\mal A_j^\C(J)=\mal A_j^\C(J).J=\mal A_j^\C(J)$.} 
$$
\mal B_j^\C(J)=\mal A_j^\C(J)\cap(J.\so(2n))^\C
,\quad \text{resp.} \quad \so_j^\C(J)=\mal A_j^\C(J)\cap\so(2n)^\C . 
$$
Then we have, according to (\ref{TJ})-(\ref{TJ0}),
\begin{prop}\label{T_j=sum-B_j}
The tangent space of $\zdk^a(\rdn)$ at one element $J$ is given by
\begin{equation}\label{TJB}
T_J\zdk^a(\rdn)=\left(\oplus_{j=1}^{r-1}\mal{B}_j^\C(J)\right)\cap \End(\rdn).
\end{equation}
\end{prop}
The inner automorphism\footnote{The conjugaison by $J$ is denoted by $\Int J\colon GL_n(\R)\to GL_n(\R)$ when  the domain of definition is a Lie subgroup and by $\Ad J\colon \mak{gl}_n(\R)\to \mak{gl}_n(\R)$ when it is a Lie subalgebra.} $T=\Int J_{|SO(2n)}$ gives rise to the $r$-symmetric space $SO(2n)/\U_0(J)$, where $\U_0(J)=SO(2n)^T=\com(J)\cap SO(2n)$, which is nothing but the connected component $\zdk^\alpha$ of $J$ (which is also the orbit $SO(2n)\cdot J=\Int (SO(2n))(J)$\,):
$$
\zdk^\alpha(\rdn)=SO(2n)/\U_0(J).
$$
Consider now
$$
\U_{j-1}(J):=\com(J^j)\cap SO(2n)=SO(2n)^{T^j},
$$
Then $T$ is an order $j$ automorphism\footnote{We confuse $j\in \Z_r$ and its representant in $\{1,\ldots,r\}$.}   on $\U_{j-1}(J)$ and gives rises to the $j$-symmetric space $\U_{j-1}(J)/\U_0(J)$ which is in fact equal to
$$
\mal Z_{2k,j}^\alpha(\rdn,J^j):=\{J'\in\zdk^\alpha(\rdn)|(J')^j=J^j\}.
$$
Indeed let $J'\in\zdk^\alpha(\rdn)$, then there exists $g\in SO(2n)$ such that $J'=gJg^{-1}$, then $(J')^j=J^j$ if and only if $gJ^jg^{-1}=J^j$ i.e. $g\in\U_{j-1}(J)$, which proves that $\mal Z_{2k,j}^\alpha(\rdn,J^j)=\Int(\U_{j-1}(J))(J)$ i.e.
$$
\mal Z_{2k,j}^\alpha(\rdn,J^j)=\U_{j-1}(J)/\U_0(J).
$$
Let us recapitulate what precedes.
\begin{prop}\label{zdkj-k-sym-space}
Let $J\in\zdk(\rdn)$. The connected component $\zdk^\alpha(\rdn)$ of $J$ is a $SO(2n)$-orbit (via the adjoint action):
$$
\zdk^\alpha(\rdn)=SO(2n)/\U_0(J).
$$
Moreover, denoting by $r$ the order of this component (i.e. the order of $\Ad J$), then it is also a $r$-symmetric space defined by the inner automorphism $T=\Int J_{|SO(2n)}$.\\
Furthermore, for any $j\in\{1,\ldots,r\}$, let us consider the subgroup $\U_{j-1}(J):=\{g\in SO(2n), gJ^jg^{-1}=J^{j} \}$.
Then the submanifold of $\zdk^\alpha(\rdn)$, defined by 
$$
\mal Z_{2k,j}^\alpha(\rdn,J^j):=\{J'\in\zdk^\alpha(\rdn)|(J')^j=J^j\}
$$
is the $\U_{j-1}(J)$-orbit of $J$.  Moreover,  the restriction of $T$ to $\U_{j-1}(J)$ is an order $j$ automorphism of which  the associated $j$-symmetric space is 
$$
\mal Z_{2k,j}^\alpha(\rdn,J^j)=\U_{j-1}(J)/\U_0(J).
$$ 
We have then a increasing sequence $\left( N_j\right) $ indexed by $j\in\{1,\ldots,r\}$, of $j$-symmetric spaces, all included in $N_r=\zdk^\alpha(\rdn)$.
\end{prop}

\begin{rmk}\emph{
Obviously, in this equation $J$ can be replaced by any $J'\in\mal Z_{2k,j}^\alpha(\rdn,J^j)$.}
\end{rmk}
\begin{exam}
\em If $k=2$, then we have $\mal Z_{4,2}^\alpha(\rdn,J^2)=\mal Z_{4,2}^\alpha(\rdn,-\Id)=\mal Z_{4}^\alpha(\rdn)=\Sigma^\alpha(\rdn)= SO(2n)/U(n)$, and the other values of $j$ are trivial $\mal Z_{4,\pm 1}^\alpha(\rdn,J^{\pm 1})=\{J\}$.
\end{exam}
\begin{rmk}\label{nouvellenotation}\em
Sometimes, we will need to specify clearly what is the eigenvalues of the eigenspaces $\mal A_i^\C(J)$ and $\so_i^\C(J)$, then we will simply use the notation
$$
\mal A_{(\omega)}^\C(J)=\ker(\Ad J -\omega\Id)
$$
and idem for $\so_{(\omega)}^\C(J)$ and $\mal B_{(\omega)}^\C(J)$.\\
Besides, sometimes for the homogeneity of the equations, we will extend the notations  $\mal A_i^\C(J)$ for real index and set for $t\in \R$
$$
\mal A_{t}^\C(J)=\ker(\Ad J -\omega_r^t\Id).
$$
\end{rmk}
\subsubsection{Study of $\Ad J^j$}\label{AdJj}
\paragraph{Expression of the eigenspaces of $\Ad J^j$ in terms of those of $\Ad J$.} 
Let $j\in\Z^*$. Then we have
$$
\Ad J^j=(\Ad J)^j= \oplus_{l=0}^{r-1}\, \left(\omega_r^l\right)^j\,\Id_{\mal A_l^\C(J)},
$$
$\omega_r^j$ is of order $p=\displaystyle\frac{r}{(r,j)}$, i.e. it is in $\Hat U_p=\{z\in S^1|z^p=1\}=\exp\left((\Z/p\Z)\cdot \displaystyle\frac{2i\pi}{p}\right)$. Then we have
$$
\Ad J^j=  \bigoplus_{q=0}^{(r,j)-1}\left[\oplus_{l=0}^{p-1}\, \left(\omega_r^j\right)^l\,\Id_{\mal A_{qp+l}^\C(J)}\right],
$$
hence writing (that $\Ad J^j$ is of order $p$):
$$
\Ad J^j=\oplus_{l=0}^{p-1}\, \omega_p^l\,\Id_{\mal A_l^\C(J^j)},
$$
we obtain 
\begin{lemma}
We have the following relation between the eigenspaces of $\Ad J^j$ and those of $\Ad J$:
\begin{equation}\label{alj}
\mal A_l^\C(J^j)=\oplus_{q=0}^{(r,j)-1}\mal A_{qp+l'}^\C(J)
\end{equation}
where $l'=(j')^{-1}l$ in the ring $\Z/p\Z$, and $j'=\left[\displaystyle\frac{j}{(j,r)}\right]_{\text{mod }p}$ ($j'$ is inversible in the ring $\Z/p\Z$, since $(j',p)=1$ by definition of $(r,j)$).
\end{lemma}
In particular,
\begin{equation}\label{comJ}
\com(J^j)^\C=\mal A_0^\C(J^j)=\oplus_{q=0}^{(r,j)-1} \mal A_{qp}^\C(J).
\end{equation}
More particulary,
$$
\com (J^k)^\C=\begin{cases}
\oplus_{q=0}^{k-1}\mal A_q^\C(J)=\mrm{End}(\rdn)^\C & \text{ if } r=k\\
\oplus_{q=0}^{k-1}\mal A_{2q}^\C(J) & \text{ if } r=2k
\end{cases}
$$
and
$$
\com (J^2)=\begin{cases}
\left.
\begin{array}{ll}
\mal{A}_0(J)\oplus\mal A_k(J) & \text{ if } r=2k\\
\mal{A}_0(J)\oplus\mal A_{\frac{k}{2}}(J) & \text{ if } r=k\in 2\Z
\end{array}\right\}
= \mal{A}_0(J)\oplus\mal A_{\frac{r}{2}}(J) & \text{ if } r \text{ is even}\\
\mal A_0(J) & \text{ if } r=k \text{ is odd}
\end{cases}.
$$
\paragraph{Decomposition of the tangent space $T_J \mal Z_{2k,j}^\alpha(\rdn,J^j)$ in terms the eigenspaces $\mal B_{j}^\C(J)$ of $\Ad J$.}
We can rewrite all what precedes in $J.\so(2n)$ (resp. in $\so(2n)$) by replacing $\mal A_l^\C$ by $\mal B_l^\C$ (resp. $\so_l^\C$). In particular we have, according to (\ref{comJ}),
$$
\ul_{j-1}(J):=\lie (\U_{j-1}(J))=\so_0(J^j)=\left( \oplus_{q=0}^{(r,j)-1}\so_{qp}^\C(J)\right)\cap \so(2n), 
$$
this\footnote{We mean $\ul_{j-1}(J)^\C=\oplus_{q=0}^{(r,j)-1}\so_{qp}^\C(J)$}  is the eigenspace decomposition of the order $j$ automorphism  obtained by restricting $T=\Ad J$ to $\ul_{j-1}(J)$. Therefore we have
\begin{prop}
The tangent space of $\mal Z_{2k,j}^\alpha(\rdn,J^j)$ is given by
\begin{equation}\label{tangentspace}
T_J \mal Z_{2k,j}^\alpha(\rdn,J^j)=\left( \bigoplus_{q=1}^{(r,j)-1}\mal B_{qp}(J)\right) \cap\End(\rdn).
\end{equation}
\end{prop}
\textbf{Proof.}
Indeed, $g\in\U_{j-1}(J)\mapsto gJg^{-1}\in\mal Z_{2k,j}^\alpha(\rdn,J^j)$ is a surjective submersion whose the (surjective) derivative at $g=1$,
$$
A\in\ul_{j-1}(J)\mapsto [A,J]=\sum_{q=0}^{(r,j)-1}[A_{qp},J]=\sum_{q=0}^{(r,j)-1}(1-\omega_r^{qp}) A_{qp}J\in\ T_J\mal Z_{2k,j}^\alpha(\rdn,J^j)
$$
has $\left( \oplus_{q=1}^{(r,j)-1}\mal B_{qp}(J)\right)\cap \End(\rdn)$ as image, which proves the equality (\ref{tangentspace}).\hfill $\square$\medskip \\
More simply by differentiating the definition equation of $\mal Z_{2k,j}^\alpha(\rdn,J^j)$ we  obtain
\begin{equation}\label{diffdef}
T_J\mal Z_{2k,j}^\alpha(\rdn,J^j)=\left\{A\in J.\so(2n)\left|\sum_{p+l=j-1} J^pAJ^l=0\right.\right\}
\end{equation}
Now, let us apply (\ref{tangentspace}) for $j=2$ :
$$
T_J \mal Z_{2k,2}^\alpha(\rdn,J^2)=\begin{cases}
\mal B_{\frac{r}{2}}(J)=\{A\in J.\so(2n)|AJ+JA=0\} & \text{ if } r \text{ is even}\\
0 & \text{ if } r \text{ is odd.}
\end{cases}
$$
This can be recovered from (\ref{diffdef}) by remarking that if $r$ is odd then $-1$ is not a $r$-th root of unity (and thus not an eigenvalue of $\Ad J$).
\begin{rmk}
\em
If $(j,2k)=1$ (so that $(j,r)=1$ also) then $J^j\in\zdk(\rdn)$ and $T^j$ is of order $r$ and we have, according to (\ref{alj})
$$
\mal A_l^\C(J^j)=\mal A_{[j]_r^{-1}\cdot l}^\C(J),\quad \forall l\in\Z/r\Z,
$$
in others words $\mal A_l^\C(J)=\mal A_{j\cdot l}^\C(J^j)$, $\forall l\in\Z/r\Z$. In particular
$$
\U_{j-1}(J)=\mal A_0(J^j)\cap SO(2n)=\mal A_0(J)\cap SO(2n)=\U_0(J).
$$
Hence
$$
\mal Z_{2k,j}^\alpha(\rdn,J^j)=\{J\}.
$$
More generally, we have, according to (\ref{comJ}),  since $(jl,r)=(l,r)$,
$$
\com((J^j)^l)^\C=\com(J^{jl})^\C= \oplus_{q=0}^{(l,r)-1}\mal A_{qp}^\C(J)=\com(J^l)^\C
$$
 with $p=\dfrac{r}{(l,r)}$, and thus
$$
\com((J^j)^l)=\com(J^l)\quad \forall l\in\Z/r\Z.
$$
In particular, $\U_{l-1}(J^j)=\U_{l-1}(J)$ $\forall l\in\Z/r\Z$ and thus
$$
\mal Z_{2k,l}^\alpha(\rdn,(J^j)^l)=\mal Z_{2k,l}^{[j]_{2k}^{-1}\cdot\alpha}(\rdn,J^l)
$$
where $[j]_{2k}^{-1}\cdot\alpha$ is the action of $[j]_{2k}^{-1}$ on $\alpha\in X_{2k}$, the action of $l\in(\Z/r\Z)^*$ on $X_{2k}$ being defined by the bijective map
$$
J\in\zdk(\rdn)\mapsto J^l\in \left(\zdk(\rdn)\right)^l=\zdk(\rdn)
$$
which sends a connected component onto another one
$$
\left(\zdk^\alpha(\rdn)\right)^l=:\zdk^{l\cdot\alpha}(\rdn).
$$
In particular, for $j=1$, we have
$$
\U_{l-1}(J^{-1})=\U_{l-1}(J)\quad \forall l\in \Z/r\Z
$$
and thus
$$
\mal Z_{2k,l}^\alpha(\rdn,(J^{-1})^l)=\mal Z_{2k,l}^{-\alpha}(\rdn,J^l)
$$
where $-(\eps,\und p)=((-1)^n\eps,\und p)$ in $X_{2k}$. Hence $\mal Z_{2k,l}^\alpha(\rdn,(J^{-1})^l)=\mal Z_{2k,l}^\alpha(\rdn,J^l)$ \iif $n$ is even (i.e. $J$ and $J^{-1}$ are in the same connected component).
\end{rmk}
%
%%%%%%%%%%%%%%%%%%%%%%%%%%%%%%%%%%%%%%%%%%%%%%%%%%%
%
\subsection{Isometries of order $2k+1$ with no eigenvalue $=1$}
\index{odd case}
We can do exactly the same study for $\mal Z_{2k+1}(E)$ as we did for $\zdk(E)$, with however the following simplification: all the connected components have the same order $r=2k+1$ and we do not have to distinguish two types of orbits as previously. Let us review all the  results obtained in \ref{order 2k} to see those which hold  and for these latters, see  if some modifications are necessary.\medskip\\
$\bullet$ The space $\mal Z_{2k+1}(\R^{2n+1})$ is empty.\\
$\bullet$ Theorem~\ref{thm-connectcomp} holds: 
$$
\pi_0(\mal{Z}_{2k+1}(\R^{2n}))= X_{2k+1}:= \left\{(\varepsilon,p)\in\Z_2\times\mathbb N^{k}\left |\,\sum_{j=1}^{k}p_j=n\ \mrm{and}\ \mrm{lcm}\left(\left\{ \frac{2k+1}{(2k+1)\wedge j}, p_j\neq 0 \right\}\right)=2k+1\right.\right\}
$$
i.e. more generally
$$
\pi_0(\mal{Z}_{k'}(\R^{2n}))= X_{k'}:= \left\{(\varepsilon,p)\in\Z_2\times\mathbb N^{m_{k'}-1}\left |\,\sum_{j=1}^{m_{k'}-1}p_j=n\ \mrm{and}\ \mrm{lcm}\left(\left\{ \frac{k'}{k'\wedge j}, p_j\neq 0 \right\}\right)=k'\right.\right\}
$$
$\bullet$ Equation~\eqref{TJ} holds in general for any order. \\
$\bullet$ Propositions~\ref{sum-Aj} and \ref{T_j=sum-B_j} hold also for any order $k'$.\\
$\bullet$ Propositions~\ref{zdkj-k-sym-space} holds for  $\mZ_{2k+1}(\rdn)$ (replace everywhere $2k$ and $r$ by $2k+1$).\medskip\\
More generally, let us consider the family of spaces $\mal U_{k'}(\rdn)$. Then we could write explicitely
$\pi_0(\mal U_{k'})$. Moreover equation~\eqref{TJ} holds in general. Propositions \ref{sum-Aj} and \ref{T_j=sum-B_j} hold also in general. Finally, Propositions \ref{zdkj-k-sym-space} holds also in general. Let us make precise that then $r$ is always the order of $\Ad J$, in the considered connected component. 

%%%%%%%%%%%%%%%%%%%%%%%%%%%%%%%%%%%%%%%%%%%%%%%%%%%
%
\subsection{The effect of the power maps on the finite order isometries}\label{effect-power-map}
Let $J\in\mal U_{k'}(\rdn)$ then $J^j\in \mal U_p(\rdn)$ with $p=\frac{k'}{(k',j)}$. Moreover it is easy to see (from the diagonalisation) that the power map
$$
J\mapsto J^j
$$
is surjective from $\mal U_{k'}(\rdn)$ onto $\mal U_p(\rdn)$ (since $z\in \Hat U_{k'}(\rdn)\mapsto z^j\in \Hat U_p(\rdn)$ is surjective). Besides, since each connected component in $\mal U_{k'}(\rdn)$ (and in $\mal U_p(\rdn)$) is a $SO(2n)$-orbit, then the power map $J\mapsto J^j$ sends one component in $\mal U_{k'}(\rdn)$ onto another one in $\mal U_p(\rdn)$ so that it induces a map:
$$
\alpha\in \pi_0(\mal U_{k'}(\rdn))\longrightarrow j\cdot\alpha\in\pi_0(\mal U_p(\rdn))
$$
such that
\begin{equation}\label{def-j-alpha}
(\mal U_{k'}^\alpha)^j(\rdn)=:\mal U_p^{j\cdot\alpha}(\rdn),\quad  \forall \alpha\in \pi_0(\mal U_{k'}(\rdn)).
\end{equation}
\begin{rmk}
\em In general we have $(\zdk)^j(\rdn)\nsubseteq \mal Z_p(\rdn)$. For example, for $j=2$, we have
$$
(\zdk(\rdn))^2=\begin{cases} \mal U_k^*(\rdn) & \text{ if } k \text{ is even}\\
\mal Z_k(\rdn) & \text{ if } k \text{ is odd.}
\end{cases}
$$
Besides, given $J\in\zdk^\alpha(\rdn)$, then $\mal Z_{2k,j}^\alpha(\rdn,J^j)$ is the inverse image of $J^j$ by the map 
$$ 
J'\in \zdk^\alpha(\rdn)\longmapsto (J')^j\in (\zdk^\alpha(\rdn))^j=\mal U_p^{j\cdot \alpha}(\rdn). 
$$
Since $(J')^j$ is constant in $\mal Z_{2k,j}^\alpha(\rdn,J^j)$, we can denote it by $J_j$ and then
$$
\mal Z_{2k,j}^\alpha(\rdn,J^j)=\mal Z_{2k,j}^\alpha(\rdn,J_j).
$$
\null\hfill $\blacksquare$
\end{rmk}

\noindent Furthermore, we have also for any $J\in \mal U_{k'}^\alpha(\rdn)$, $\U_0(J^j)=\U_{j-1}(J)$ so that, according to \eqref{def-j-alpha}, the component of $J^j$ in  $\mal U_{p} (\rdn)$ satisfies
$$
\mal U_{p}^{j\cdot \alpha}(\rdn)=SO(2n)/\U_0(J^j)=SO(2n)/\U_{j-1}(J),
$$
hence 
\begin{equation}\label{powerj}
(\mal U_{k'}^\alpha(\rdn))^j=SO(2n)/\U_{j-1}(J)
\end{equation}
which we can recover directly by taking the power $j$ in the equality $\mal U_{k'}^\alpha(\rdn)=\{ gJg^{-1}, g\in SO(2n)\} $.\\[1mm]
\textbf{Convention}: for each $\alpha\in \pi_0(\zdk(\rdn))$, we will choose (and fix) a canonical representant in $\zdk^\alpha(\rdn)$. For example, let $(\epsilon_1,\ldots,\epsilon_{2n})$ be the canonical basis in $\rdn$, and 
\begin{eqnarray*}
e_{2l+1} & = & \frac{\epsilon_{2l+1} + i \epsilon_{2l+2}}{\sqrt{2}},\ 0\leq l \leq n-1,\\
e_{2l} & = &  \overline e_{2l-1}, \ 1\leq l \leq n.
\end{eqnarray*}
Then $e=(e_1,\ldots, e_{2n})$ is a hermitian basis  in $\C^{2n}$ and  we can take $J_0^\alpha$ such that
\begin{equation}\label{eq-def-J0alpha}
\mal{M}at_e(J_0^\alpha)=\mrm{Diag}\left( \begin{pmatrix} e^{i\theta_j}\Id_{p_j} & 0 \\ 0 & e^{-i\theta_j}\Id_{p_j}
\end{pmatrix}, 1\leq j \leq n-1  \right) 
\end{equation}
where $\und p=(p_1,\ldots, p_{k-1})$ is determined by $\alpha=(\varepsilon,\und p)\in\pi_0(\zdk(\rdn))$ (see section~\ref{setconnectcomp}), $\theta_j=\displaystyle\frac{j\pi}{k}$, and $\mal M at_e(\cdot)$ means "the matrice in the basis $e$ of ".%\\[1.5mm]
\subsection{The Twistor spaces of a Riemannian manifolds and its reductions}\label{3twistorspacesandreduction}
Let $M$ be  an oriented (even dimensional) Riemannian manifold and let us consider the bundle of order $2k$ isometries $\mal U_{2k}(M)$ as well as its subbundles $\mal U_{2k}^*(M)$ and $\zdk(\rdn)$. Let us fix $\alpha\in \pi_0(\zdk(\rdn))$ and consider the  component $\zdk^\alpha(M)$. Then denoting by $\mal{SO}(M)$ the $SO(2n)$-bundle of positively oriented orthonormal frames on $M$, we have\footnote{The canonical isomorphism which allows to write this  equality will be explicitly given in section~\ref{examhomfibrbund}.} 
$$
\zdk^\alpha(M)=\mal{SO}(M)/\U_0(J_0^\alpha).
$$
We want to ask the following question: does $\mal{SO}(M)$ admit  a $\U_{j-1}(J_0^\alpha)$-reduction for $1\leq j \leq r$.  We know (according to \cite{KN}) that $\mal{SO}(M)$ admits a $\U_{j-1}(J_0^\alpha)$-reduction \iif the associated bundle $\mal{SO}(M)/\U_{j-1}(J_0^\alpha)$ ($=\mal{SO}(M)\times_{SO(2n)}SO(2n)/\U_{j-1}(J_0^\alpha)$) admits a global section $J_j\colon M\to \mal{SO}(M)/\U_{j-1}(J_0^\alpha)$.\\
Besides,  according to (\ref{powerj}) and \eqref{def-j-alpha} applied to $\zdk^\alpha(\rdn)$, we have 
$$
\mal{SO}(M)/\U_{j-1}(J_0^\alpha)=(\zdk^\alpha(\rdn))^j=\mal U_p^{j\cdot\alpha}(M)
$$
with $p=\dfrac{2k}{(2k,j)}$. Hence $J_j$ (when it exists) is a global section of $(\zdk^\alpha(M))^j$ and then the $\U_{j-1}(J_0^\alpha)$-reduction of $\mal{SO}(M)$ is given in terms of $J_j$ by:
$$
\Um_{j-1}^\alpha(M):=\{e=(e_1,\ldots,e_{2n})\in \mal{SO}(M)|\,\mal Mat_e(J_j)=(J_0^\alpha)^j\}.
$$ 
Then we have  
$$
\Um_{j-1}^\alpha(M)/\U_0(J_0^\alpha)=\mal Z_{2k,j}^\alpha(M,J_j).
$$
In particular, since $(J_0^\alpha)^r=\pm \Id$, we have $\U_{r-1}(J_0^\alpha)=SO(2n)$ and $\mal{SO}(M)$ has always an (unique and trivial) $SO(2n)$-reduction for which $ J_r=\pm\Id_{TM}$ and thus $\Um_{r-1}^\alpha(M)=\Um_0^\alpha(M)= \mal{SO}(M)$ and $\mal Z_{2k,r}^\alpha(M,J_r)=\zdk^\alpha(M)$.
\begin{exam} 
\em If  $k=2$, and thus $r=2$, then $J_0=-\Id$ defines the trivial reduction. Moreover, for $j=1$, a global section $J_1$ (when it exists) defines on $M$ an almost complex  structure and $\Um_{j-1}^\alpha(M)=\Um_0^\alpha(M)$ is then the subbundle of hermitian frames on $M$ (with respect to this almost complex structure).
\end{exam}

%%%%%%%%%%%%%%%%%%%%%%%%%%%%%%%%%%%%%%%%%%%%%%%%%%%
%
\subsection{Return to  an order $2k$ automorphism $\tau\colon\g\to\g$.}\label{return}
\index{k symmetric space@$k'$-symmetric space|(}
We give ourself the same ingredients as in section~\ref{evendetercase} and we use the same notations. In particular, we suppose that the subgroup $H$ is chosen such that $(G^\sigma)^0\subset H\subset G^\sigma$. In addition to that we suppose $G/H$ Riemannian.
\subsubsection{Case $r=k$}
Suppose that we have $\tm^k=-\Id$ i.e. $\tm\in\zdk^0(\mk)$. Then $\g_{2j}^\C=0$ for all $2j\in\Z/(2k)\Z\setminus\{0,k\}$. Hence we have
$$
[\g_p^\C,\g_l^\C]=\{0\}\qquad \text{if } p+l\neq 0,k .
$$
Indeed, if $p$ or $l$ is even then the corresponding eigenspace vanishes. If $p$ and $l$ are odd then $[\g_p^\C,\g_l^\C]\subset \g_{p+l}^\C$ and $p+l$ is even, thus $\g_{p+l}^\C=\{0\}$ except if $p+l=0$ or $k$. Consequently, we have $[\mk,\mk]\subset\hk$ and thus $G/H$ is a (locally) symmetric space. Let us distinguish the following two cases.
\paragraph{$k$ is odd.}
Then $[\g_k,\g_j^\C]\subset\g_{k+j}^\C=\{0\}$ for all $j$ odd $\neq 0,k$. Hence $[\g_k,\mk]=\{0\}$ i.e. $\adm\g_k=0$ so that  $\g_k=0$ and thus this case is trivial because $H=G_0$ up to covering and thus the fibre $H/G_0$ is trivial (i.e. a discret set). Moreover we have $[\mk,\mk]\subset\hk=\g_0$ and $G/H=G/G_0$ (up to covering) is the (locally) symmetric space associated to the involution $\tau^k$.
\paragraph{$k$ is even.}
Then the symmetric decomposition $\g=\hk\oplus\mk$ is the eigenspace decomposition of $\tau^k$, and $G/H$ is the (locally) symmetric space corresponding to this involution $\tau^k$.
\paragraph{}In conclusion, if $r=k$, then $G/H$ is the (locally) symmetric space corresponding to $\tau^k$.

\begin{exam} \em If $2k=4$, then we  always have $r=k=2$ (since $\tm^2=-\Id$) and $G/H$ is the symmetric space corresponding to $\sigma=\tau^2$.
\end{exam}
\subsubsection{Action of $\Ad\tm$ on $\ad\g_j^\C$}\label{Action-of-Ad-taum}
We have $\tau\circ\ad X\circ\tau^{-1}=\ad\tau(X)$, $\forall X\in\g$. In particular, for all $j\in\Z/(2k)\Z $ we have
$$
\forall X_j\in \g_j^\C,\ \tau\circ\ad X_j\circ\tau^{-1}=\omega_{2k}^j\ad X_j .
$$
Hence by taking the restriction to $\mk$ and  projecting on $\mk$:
\begin{equation}\label{eq-Ad-taum}
\tm\circ[\adm X_j]_\mk\circ\tm^{-1}=\omega_{2k}^j[\adm X_j]_\mk
\end{equation}
so that\footnote{See remark~\ref{possible}.}
$$
[\adm\g_j^\C]_\mk\subset\mal A_j^\C(\tm)  \quad \forall j\in\Z/(2k)\Z.
$$
If $r=k$ then $[\mk,\mk]\subset\hk$, hence $[\adm\g_j^\C]_\mk=0$, for all $j\in\Z/(2k)\Z\setminus\{0,k\}$.\medskip\\
Moreover, let us recall that we always have ($r=2k$ or $k$)
\begin{eqnarray*}
[\adm\g_0]_\mk & = & \adm\g_0\subset \so_0(\tm)=\com(\tm)\cap\so(\mk)\\
{}[\adm\g_{k} ]_{\mk}  & = & \adm\g_k\subset
%\adm\g_{k}
%[\adm\g_0]_\mk & = & \adm\g_0\subset
\begin{cases}
\so_{(-1)}(\tm)=\mrm{Ant}(\tm)\cap\so(\mk) =\so_{\frac{r}{2}}(\tm) & \text {if } r \text{ is even}\\
0 & \text {if } r=k \text{ is odd (trivial case)}
\end{cases}
\end{eqnarray*}
where\footnote{According to the notation defined in \ref{generalities}.} $\mrm{Ant}(\tm)=\{A\in\End(\mk)|A\,\tm +\tm A=0\}$, and we used the notation defined in remark~\ref{nouvellenotation} for $\so_{(-1)}(\tm) $.
\begin{rmk}\label{possible}
\em In general, we do not  have $[\adm\g_j^\C]_\mk\subset\so_j^\C(\tm)$  and we also do not have   $[\adm\g_j^\C]_\mk\subset\mal B_j(\tm)$. However, if the metric is  naturally reductive, then these inclusions hold.
\end{rmk}
We can easily generalize \eqref{eq-Ad-taum} as follows. We keep in mind the conventions and notations defined  in the subsection~\ref{def-g-tau}. We use also the notations of the begining of subsection~\ref{finitorderauto}.
\begin{prop}\label{generali-def-g-tau}
Let $\tau\colon \g\to\g$ be an automorphism of Lie algebra of finite order $k'$. Let us consider some decomposition $k'=pq$, with $p,q\in \mathbb N^*$.  Then let us consider the order $p$ automorphism $\hat\sigma=\tau^q$, and set $\hat\hk=\g^{\hat{\sigma}}$. We denote by $\g=\hat\hk\oplus\hat\mk$ the corresponding $\hat\sigma$-invariant reductive decomposition of the Lie algebra $\g$. Then we have 
$$
\tau_{|\hat\mk}\circ[\ad_{\hat\mk} X_j]_{\hat\mk}\circ\tau_{|\hat\mk}^{-1}=\omega_{k'}^j[\ad_{\hat\mk} X_j]_{\hat\mk}, \quad \forall X_j\in \g_j^\C.
$$
\end{prop}
\begin{rmk}\em
In particular, we can apply the previous proposition in the case $k'=2k+1$ and $q=1$.
\end{rmk}
\subsection{The canonical section in $(\zdk(G/H))^2$, the canonical embedding, and the Twistor lifts}\label{can-sect-zdk}
\index{canonical!embedding|(}\index{embedding, canonical|see{canonical embedding}}
\subsubsection{The canonical embedding}\label{subsub-canonical-embedding}
Once more, we give ourself the same ingredients and  notations as in section~\ref{evendetercase}. We suppose that $(G^\sigma)^0\subset H\subset G^\sigma$ and that $G/H$ is Riemannian. We denote by $p_0:=1.H$ the reference point in $G/H$.  According to the definition of $H$, we have
\begin{lemma}\label{lemma2}
Let $J_0$ be the element in $\zdk^{\alpha_0}(T_{p_0}M)$ corresponding \footnote{$\alpha_0$ denotes off course the connected component of $J_0$ in $\zdk(T_{p_0}M)$.} to $\tm$  under the identification $T_{p_0}M=\mk$. Then we have $\forall g\in H$, $gJ_0^2g^{-1}=J_0^2$. Hence there exists a unique section
$$
 J_2\colon G/H \mapsto (\zdk^{\alpha_0})^2=\mal U_k^{\alpha_0^2}(G/H)
$$
defined by
$$
g.p_0\in G/H\mapsto gJ_0^2 g^{-1}\in (\zdk^{\alpha_0})^2.
$$
\end{lemma}
Proceeding as in \cite{ki3}, Theorem~3, we obtain:
\begin{thm}\label{embedding}
Let $\tau\colon\g\to\g$ be an order $2k$ automorphism and $M=G/H$ a (locally) Riemannian $k$-symmetric space corresponding to $\sigma=\tau^2$. Let us make $G$ acting on $\zdk(M)$: $g\cdot J=gJg^{-1}$. Let $J_0\in\zdk^{\alpha_0}(T_{p_0}M)$ be the finite order isometry corresponding to $\tm$ under the identification $T_{p_0}M=\mk$. Then the orbit of $J_0$ under the action of $G$ is an immersed submanifold in $\zdk^{\alpha_0}(M)$. Denoting by $G_0$ the stabilizer of $J_0$, then $G_0=G^\tau\cap H$ and thus $N=G/G_0$ is a locally $2k$-symmetric bundle over $M$ and the natural map:
$$
\begin{array}{crcl}
\mak I_{J_0}\colon & G/G_0 & \longrightarrow & \mal Z_{2k,2}^{\alpha_0}(G/H, J_2)\\
              & g.G_0 & \longmapsto    &    gJ_0g^{-1}
\end{array}
$$
is an injective immersion and a morphism of bundle. Moreover, if the image of $G$ in $\Is(M)$ (the group of isometry of $M$) is closed, then $\mak I_{J_0}$ is an embedding.
\end{thm}
\begin{rmk}\em
Remark that we could also have chosen $J_0=\tm^j$ with $(j,2k)=1$, since in this case $\tau^j$ is an order $2k$-isomorphism. See \cite{ki3}, theorem 3 and remark 13.
\end{rmk}
\textbf{Notation}
For a geometric map $f\colon L\to G/G_0$, we will denote by $J$ the corresponding map $\mak I_{J_0}\circ f\colon L\to Z_{2k,2}(G/H, J_2)$ under the previous inclusion $G/G_0 \hookrightarrow Z_{2k,2}(G/H,J_2)$.
\paragraph{A remark about the  canonical section}
In fact the lemma~\ref{lemma2} can be written more generally as follows:
\begin{lemma}\label{lemma-def-J_1}
Let $N=G/G_0$ be a (locally) $k'$-symmetric space. We denote by $\g=\g_0\oplus\nk$  the reductive decomposition of the Lie algebra. Let $y_0=1.G_0\in N$ and let $J_0$ be the element in $\mal{U}_p^*(T_{y_0}N)$ corresponding  to $\tau_{|\nk}$  under\footnote{Or more generally corresponding to $\tau_{|\nk}^j$ with $(j,k')=1$.} the identification $T_{y_0}N=\nk$. Then we have $\forall g\in G_0$, $gJ_0g^{-1}=J_0$. Hence there exists a canonical section $ J_1\colon G/G_0 \mapsto \mal{U}_p^*(T_{y_0}N)$ defined by
$$
g.y_0\in G/G_0\mapsto gJ_0 g^{-1}\in \mal{U}_p^*(T_{y_0}N).
$$
\end{lemma}

\subsubsection{The twistor lifts}\label{subsub-twistorlift}
\begin{defn}\index{structure@$2k$-structure}
An isometry $A\in SO(\rdn)$ will be called an $e^{i\theta}$-structure if \,$\mrm{Spect}(A)=\{e^{i\theta},e^{-i\theta}\}$. An isometry $A\in SO(\rdn)$ will be called a $2k$-structure if $A\in\zdk(\rdn)$.
\end{defn}
\begin{defn}\label{defn-j-undj}
Let $(E,h)\to M$ be a Riemannian vector bundle over a manifold $M$. Then for each $2k$-structure $J\in \zdk(E)$, we denote by $\undj$ the complex structure in $E$ defined by 
\begin{eqnarray*}
\ker(\undj - i\Id) &  = & \oplus_{j=1}^{k-1} \ker( J -\omega_{2k}^{-j}\Id)\\
\ker(\undj + i\Id) &  = & \oplus_{j=1}^{k-1} \ker( J -\omega_{2k}^{j}\Id)
\end{eqnarray*}
\end{defn}
\begin{rmk}\em
Let us remark that if $J\in \Sigma (E)$ is complex structure then $\undj=-J$. This sign is needed because of our convention chosen in remark~\ref{rmk-convention}. See also \cite{ki3},  remark 13.
\end{rmk}
\begin{rmk}\em
We see that definition~\ref{defn-j-undj} defines a $SO(2n)$-invariant map
$$
P\colon \zdk^\alpha (\rdn) \longrightarrow \Sigma(\rdn)
$$ 
which is nothing but the $SO(2n)$-invariant projection 
$$
SO(2n)/\U_0(J_0^\alpha)\longrightarrow SO(2n)/\U(\undj_0^\alpha),
$$ 
since $\U_0(J_0^\alpha)\subset \U(\undj_0^\alpha)$. Indeed by definition $\undj_0^\alpha$ stabilizes the eigenspaces of $J_0^\alpha$ : 
$$
\undj_0^\alpha= \left( \oplus_{j=1}^{k-1} i \Id_{\left[\g_{-j}^\C\right]}\right) \oplus \left( \oplus_{j=1}^{k-1} -i \Id_{\left[\g_{j}^\C\right]}\right).
$$
Moreover, the restriction of $P$ to $\mZ_{2k,j}(\rdn,(J_0^\alpha)^j)$ is the $\U_{j-1}(J_0^\alpha)$-invariant projection
$$
\U_{j-1}(J_0^\alpha)/\U_0(J_0^\alpha)\longrightarrow \U_{j-1}(J_0^\alpha)/\U(\undj_0^\alpha).
$$
\end{rmk}
\begin{defn}\label{def-undj-odd}
Let $(E,h)\to M$ be a Riemannian vector bundle over a manifold $M$. Then for each element $J\in \mal Z_{2k+1}(E)$, we denote by $\undj\in \Sigma (E)$ the complex structure in $E$ defined by 
\begin{eqnarray*}
\ker(\undj - i\Id) &  = & \oplus_{j=1}^{k} \ker( J -\omega_{2k+1}^{-j}\Id)\\
\ker(\undj + i\Id) &  = & \oplus_{j=1}^{k} \ker( J -\omega_{2k+1}^{j}\Id)
\end{eqnarray*}
\end{defn}
\begin{rmk}\index{canonical!almost complex structure}\em
We see then that the canonical almost complex structure $\undj$ defined in a $(2k+1)$-symmetric space $N=G/G_0$ (see \eqref{def-odd}), is the almost complex structure associated to the element $J_1\in \mal Z_{2k+1}(N)$ defined by lemma~\ref{lemma-def-J_1}, according to definition~\ref{def-undj-odd}.
\end{rmk}
\begin{defn}\label{def-twistor-lift}
Let $(L,i)$ be a complex manifold (of dimension $d\geq 1$), $M$ an oriented Riemannian manifold and $u\colon L\to M$ a immersion. Then an element $J\colon L\to u^*(\zdk(M))$ is an admissible twistor lift of $u$ if one of the following equivalent statements holds:
\begin{description}
\item[(i)] $[\Bar\partial u]_{\mk_j(J)^\C}\in \ker( J -\omega_{2k}^{j}\Id)$   for all $1\leq j\leq k-1$, where
$ \mk_j(J)\subset u^* TM $ is defined by 
$$
\mk_j(J)^\C = \ker( J -\omega_{2k}^{-j}\Id)\oplus \ker( J -\omega_{2k}^{j}\Id).
$$
\item[(ii)]  Let $\undj$ be the complex structure on $u^*(TM)$ defined  by the  $2k$-structure $J\in\mal C(u^*(\zdk(M)))$, then $u$ is $\undj$-holomorphic: $du\circ i=\undj\circ du$
\end{description}
In particular, if $(L,i)$ is a Riemann surface, then we can add that the existence of  an admissible twistor lift $J$ of $u$ implies in particular that $u$ is a conformal immersion.
\end{defn}
\begin{rmk}\em
Let $E_j$ be the orthogonal projection of the tangent subbundle $u_*(TL)$ on the subbundle $\mk_j(J)$. Then $J$ is an admissible twistor lift means also : for all $j\in \{ 1,\ldots,k-1\}$, $J$ stabilizes $E_j$ and $J_{|E_j}$ is a $\omega_{2k}^j$-structure, i.e. a rotation of the plan $E_j$, and moreover this rotation defined by $J_{|E_j}$ and the rotation defined by $i$ have opposite sens of rotation.
\end{rmk}
\begin{thm}\label{thm-twistorlift&horihol}
In the situation described in theorem~\ref{embedding}, let $\alpha$ be a $\g$-valued Maurer-Cartan 1-form on a Riemann surface $L$ and $f\colon L\to G/G_0$ its geometric map and $J=\mak I_{J_0}\circ f$. The the following statements are equivalent:
\begin{description}
\item[(i)] $\alpha_{-j}''=0$, $1\leq j\leq k-1$
\item[(ii)] $J\colon L\to \mal Z_{2k,2}(G/H, J_2)$ is an admisible twistor lift.
\end{description}
\end{thm}
\proof
This follows immediately from definition~\ref{def-twistor-lift}-(i).  \hsq
\index{twistor|)}\index{even case|)}\index{canonical!embedding|)}\index{k symmetric space@$k'$-symmetric space|)}

\subsection{Bibliographical remarks and summary of the results.}
To our knowledge, the "generalized" twistor spaces defined in this section have never been used before. Until now, among this family of twistor spaces, only the twistor bundle of orthogonal almost complex structure $\Sigma(M)$ has been used, in twistor theory. The twistor theory invented by R. Penrose \cite{Penrose} could be defined as a way to use the methods of complex (or pseudo-complex) geometry to solve problems of Riemannian (or pseudo-Riemannian) geometry. More particulary, there are a lost of works which have been down about twistor methods for harmonic maps \cite{rawnsley,Salomon,BuRaw}. In these works the method is to use twistor lifts of harmonic maps and  to prove that a map $f\colon (L,j)\to (M,g)$ is harmonic if and only if it has a holomorphic twistor lift $J\colon (L,j)\to (\Sigma(M),\check{J})$. Here $(L,j)$ is a Riemann surface, $(M,g)$ a Riemannian manifold and $\check{J}$ is some almost complex structure on the twistor space $\Sigma(M)$.\\
The definition~\ref{def-twistor-lift} above generalizes the notion of twistor lift to our "generalized" twistor spaces $\zdk(M)$. Moreover the theorem~\ref{thm-twistorlift&horihol} tells us that a map $f\colon L\to G/G_0$ is horizontally holomorphic \iif the corresponding map  $J=\mak I_{J_0}\circ f\colon L\to \mal Z_{2k}(G/H)$ is an admissible twistor lift. We will prove in section~\ref{vertically-harmonic} that the even minimal determined system has an interpretation in terms of vertically harmonic admissible twistor lifts.\\
Let us add that our theorem~\ref{embedding} about the canonical embedding $G/G_0 \hookrightarrow  \mal Z_{2k,2}^{\alpha_0}(G/H, J_2)$ is of course totally new and it generalises our previous result \cite[ Theorem~3]{ki3} which concerns the particular case of 4-symmetric spaces.

%%%%%%%%%%%%%%%%%%%%%%%%%%%%%%%%%%%%%%%%%%%%%%%%%%%%%%%%%%%%%%%%%%%%%%%%%%%%%%%%%%%%%%%%%%%%%%

%  Vertically Harmonic maps and Harmonic sections of submersions

%%%%%%%%%%%%%%%%%%%%%%%%%%%%%%%%%%%%%%%%%%%%%%%%%%%%%%%%%%%%%%%%%%%%%%%%%%%%%%%%%%%%%%%%%%%%%%%%%%%%%%%%%%%%%%%%%

%%%%%%%%%%%%%%%%%%%%%%%%%%%%%%%%%%%%%%%%%%%%%%%%%%%%%%%%%%%%%%%%%%%%%%%%%%%%%%%%%%%%%%%%%%%%%%

%  Vertically Harmonic maps and Harmonic sections of submersions

%%%%%%%%%%%%%%%%%%%%%%%%%%%%%%%%%%%%%%%%%%%%%%%%%%%%%%%%%%%%%%%%%%%%%%%%%%%%%%%%%%%%%%%%%%%%%%%%%%%%%%%%%%%%%%%%%

\section{Vertically Harmonic maps and Harmonic sections of submersions}\label{vertically-harmonic}
\index{vertically harmonic|(}
We will first recall  some definitions and properties about vertical harmonicity and homogeneous fibre bundles (sections~\ref{preliminaires} and \ref{Homfibrbund}). We refer to \cite{cmw1,cmw2} for proofs. Then in section~\ref{examhomfibrbund} we will apply this to the study of the examples which interest us, examples that we will already introduce and start to study in  \ref{preliminaires-examples} : homogeneous spaces and twistor spaces. Finally, we will conclude with a geometric interpretation of the even minimal determined elliptic integrable system in terms of vertically harmonic twistor lifts (section~\ref{geom-interpret}).
\subsection{Definitions, general properties and examples}\label{preliminaires}
\subsubsection{The vertical energy functional}\label{verticalenergyfunctional}
\index{horizontal subbundle|(}\index{vertical tension field|(}
\index{vertical subbundle|(}
Let $(M,g)$, $(N,h)$ be Riemannian manifolds and $\pi\colon N\to M$ a surjective submersion. We can do the splitting $TN=\ver\oplus\hor$, where the vertical and horizontal subbundles are defined by $\ver=\ker d\pi$ and  $\hor=(\ker d\pi)^\perp=\ver^\perp$.\\
For any  map $u \colon M\to N$, we  denote by $d^v u=(du)^v$ the vertical component of $du$. Following \cite{cmw1}, this allows us to define the vertical energy density of $u$, $e^v(u)=\dfrac{1}{2}|d^v u|^2$, and the associated vertical energy functional:
$$
E^v(u)=\dfrac{1}{2}\int_M |d^v u|^2 d\mrm{vol}_g .
$$
Let us define the vertical tension field of $u\colon M\to N$ by
$$
\tau^v(u)=\mrm{Tr}_g(\nabla^v d^v u)
$$
where $\nabla^v$ denotes the vertical component of the Levi-Civita connection (of $N$) in $TN$, and $\mrm{Tr}_g$ the trace with respect to $g$. Then we have 
\begin{thm}\cite[Theorem~2]{cmw1}
The map $u\colon M \to N$ is a critical point of $E^v$ with respect to vertical variations \iif $\tau^v(u)=0$. In particular, if $u$ is a section, i.e. $\pi\circ u=\Id_M$, then it is a critical point of $E^v$ with respect to variations through sections \iif  $\tau^v(u)=0$.
\end{thm}
\begin{defn} A map $u\colon M \to N$ is vertically harmonic  if $\tau^v(u)=0$. If moreover $u$ is a section we will say that it is a  harmonic section. 
\end{defn}
\subsubsection{Examples}\label{preliminaires-examples}
\begin{exam}\label{example1}\em
Let $\pi\colon N\to M$ be like above. Let $(L, b)$ be a Riemannian manifold and $f\colon L \to N$ a map. Then we can consider the projection $u=\pi\circ f\colon L\to M$ and the manifold 
$$
u^*N:=\{(z,n)\in L\times N,\, n\in\pi^{-1}(\{u(z)\}) \}.
$$
Then we have  the submersion $u^*\pi \colon (z,n)\in u^*N\mapsto z\in L$. Furthermore, $u^*N$ can be endowed canonically with a Riemannian  metric: take the metric induced by the product metric
$$
|(dz,dn)|^2=|dz|^2 + |dn|^2
$$
in $L\times N \supset u^*N$.
\begin{defn} 
We will say that $f\colon L\to N$ is vertically harmonic if
$$
\mrm{Tr}_b(\nabla^v d^v f)=0.
$$
\end{defn}
When $u=\pi\circ f$ is an isometry and $\pi$ a Riemannian submersion \textbf{this is equivalent to say that the corresponding section $\tl f\colon L \to u^*N$ is a harmonic section} (see the Appendix, theorem~\ref{same-vertical-Harm}).
\end{exam}
%
%%%%%%%%%%
%
\begin{exam}\label{sigmaE}\index{twistor|(}\em
\textbf{The twistor bundle of almost complex structures $\Sigma(E)$.}\\
Let $\mrm p\colon (E,\nabla, \langle\cdot,\cdot\rangle)\mapsto (M,g)$ be a Riemannian vector bundle of rank $2n$ (in particular $\langle\cdot,\cdot\rangle$ is $\nabla$-parallel). Then we consider the bundle of orhogonal almost complex structure: $N_\Sigma=\Sigma(E)=\{ (x,J_x), J_x\in \Sigma(E_x)\}$, where $\Sigma(E_x)=\{ J\in \so(E_x)| J^2=-\Id\}$. We have a fibration $\pi_\Sigma\colon N_\Sigma\to M$. The vertical space is given by: $\forall J \in N_\Sigma$,
$$
\ver_J:= T_J\Sigma(E_x)=\{ A \in \so(E_x)| AJ + JA=0\} 
$$
where $x=\pi_\Sigma(J)$.\\
The metric connection $\nabla $ gives us a splitting : $T\Sigma(E)=\ver^\Sigma\oplus \hor^\Sigma$. Indeed we have the following splitting (coming from $\nabla$)
\begin{equation}\label{splitting*}
T\so(E)=\mrm p^*\so(E)\oplus \hor
\end{equation}
where $\mrm p\colon \so(E)\to M$ is the natural fibration\footnote{we denote by the same letter the fibration $\mrm p\colon E\to M$ and all its "tensorial extensions": $\mrm p\colon \End(E)\to M$, $\mrm p\colon \so(E)\to M$, etc..}. Then for any (local) section $J\colon U\subset M\to  \Sigma(E)$, we have 
$$
0=\nabla J^2=(\nabla J)J + J(\nabla J)
$$
so that $\nabla J\in \ver^\Sigma$ and thus in the decomposition (\ref{splitting*}): $[dJ]_{so(E)}=\nabla J\in \ver^\Sigma$ and thus 
$[dJ]_\hor=dJ -\nabla J\in TN_\Sigma $ which allows  us to conclude that 
$$
T\Sigma(E)=\ver^\Sigma\oplus\hor_{|\Sigma(E)}.
$$
Then we can endow $N_\Sigma$ with the metric 
\begin{equation}\label{metric}
h=\pi^*g  + \langle\ ,\ \rangle_{\ver^\Sigma}
\end{equation}
where $\langle\ ,\ \rangle_{\ver^\Sigma}$ is the fibre metric in $\ver^\Sigma$ induced by the metric in $\so(E)$:
\begin{equation}\label{tracemetric}
\langle A,B\rangle = \mrm{Tr}(A^t.B).
\end{equation}
With this metric  we have obviously $\hor^\Sigma={\ver^\Sigma}^\perp$.\\[1mm]
Furthermore, let us remark that $T\Sigma(E)$ is a subbundle of $T\so(E)_{|\Sigma(E)}$ and that we have
\begin{eqnarray}\label{dec1}
T\so(E)_{|\Sigma(E)}=\pi_\Sigma^*\so(E)\oplus\hor_{|\Sigma(E)} & = & \so_-(\pi_\Sigma^*E) \oplus \so_+(\pi_\Sigma^*E) \oplus \hor_{|\Sigma(E)}\\
 & =  & T\Sigma(E)\oplus\so_+(\pi_\Sigma^*E)
\end{eqnarray}
with\footnote{using the notations defined in section~\ref{AdJ} (i.e. the definition of $\so_j^\C(J)$ for $j\in\Z/2\Z$).}
\begin{eqnarray*}
\so_+(\pi_\Sigma^*E)_J & = & \so_+(E_x,J):=\{ A\in \so(E_x)|\, [A,J]=0\}\\
\so_-(\pi_\Sigma^*E)_J & = & \so_-(E_x,J):=\{ A\in \so(E_x)|\, AJ + JA=0\}=\ver_J^\Sigma
\end{eqnarray*}
for all $J\in\Sigma(E)$ (and where $x=\pi_\Sigma(J)$). In other words, $\pi_\Sigma^*E$ is canonically endowed with a complex structure: $\mal{I}_J=J$, $\forall J\in N_\Sigma$, and this complex structure defines the two spaces $\so_{\pm}(\pi_\Sigma^*E)$ by
$$
\so_{\pm}(\pi_\Sigma^*E)=\so_{\pm}(\pi_\Sigma^*E,\mal{I}).
$$
Now given a section $J\in\mal C(\pi_\Sigma)$, then we consider the vertical part of the rough Laplacian $\nabla^*\nabla J$, in the decomposition~(\ref{dec1}): $(\nabla^*\nabla J)^{\ver^\Sigma}=\dfrac{1}{2}J[J,\nabla^*\nabla J]$. We will see in section~\ref{twistor2} that this is in fact exactly the vertical tension field of $J$ in $N_\Sigma$:
$$
\tau^v(J)=\dfrac{1}{2}J[J,\nabla^*\nabla J].
$$
In particular, we recover the definition of vertical harmonicity used in $\cite{ki3}$ and $\cite{bk}$.
\end{exam}
%
%%%%%%%%%%%%%%%%%%
%
\begin{exam}\label{zdkE}\em
\textbf{The twistor bundle $\zdk(E)$ of a Riemannian vector bundle.}\\
Let $\mrm p\colon (E,\nabla, \langle\cdot,\cdot\rangle)\mapsto (M,g)$ be a Riemannian vector bundle of rank $2n$. Then we consider more generally the bundle of order $2k$ isometries $\mal U_{2k}(E)$ as well as its subbundles $\mal U_{2k}^*(E)$ and $\zdk(E)$. Let us fix $\alpha\in \pi_0(\zdk(\rdn))$ and consider the component $\zdk^\alpha(E):=N_\mal{Z}$. We have a natural fibration $\pi_\mal{Z}\colon\zdk^\alpha(E) \to M$. The vertical space is given by  
\begin{equation}\label{vzj}
\forall J\in N_\mal{Z},\quad \ver_J^\mal Z =T\zdk^\alpha(E_x)=\left( \bigoplus_{j\in\Z/r\Z\setminus\{0\}} \negthickspace \negthickspace \mal B_j^\C(E_x,J)\right)  \bigcap\End(E_x) =J.\so_*(E_x,J)
\end{equation}
according to section~\ref{AdJ} (more particulary equation (\ref{TJB})) and where
$$
\so_*(E_x,J):=\left( \bigoplus_{j\in\Z/r\Z\setminus\{0\}} \negthickspace \negthickspace \so_j^\C(E_x,J)\right) \bigcap\so(E_x).
$$
The metric connection $\nabla$ gives us a splitting: $T\zdk^\alpha(E)=\ver^{\mal Z}\oplus\hor^{\mal Z}$. Indeed we have the following splitting (coming from $\nabla$)
\begin{equation}\label{splittingsoe}
TSO(E)= \ver^{SO(E)}\oplus \hor
\end{equation}
where $\ver_J^{SO(E)}=T_J SO(E_x) = J.\so(E_x)$ (since $0=\nabla( J^tJ)=(\nabla J)^t J + J^t(\nabla J)$ $\Longrightarrow$
 $\nabla J\in T_J SO(E_x)$). Then for all (local) section $J\colon U\subset M \to N_{\mal Z}$, we have
$$
0=\nabla J^{2k}=\sum_{p+l=2k-1}J^p(\nabla J) J^l
$$
so that according to (\ref{TJ}), $\nabla J\in\ver^{\mal Z}$ and thus in the decomposition~(\ref{splittingsoe}), we have $[dJ]_{\ver^{\so(E)}}\in\ver^{\mal Z}$ and hence $[dJ]_\hor=dJ -[dJ]_{\ver^{\so(E)}}\in TN_\mal{Z} $ which  leads to
\begin{equation}\label{splittingzdk}
T\zdk^\alpha(E)=\ver^{\mal Z}\oplus\hor_{|\zdk^\alpha(E)}.
\end{equation}
Then we can endow $N_\mal{Z}$ with the metric defined as in (\ref{metric}) and where the fibre $\langle\ ,\ \rangle_{\ver^\mal Z}$ is induced by  the trace metric (\ref{tracemetric}), for which we have  $\hor^{\mal Z}={\ver^{\mal Z}}^\perp$.\\[1mm]
Furthermore let us  remark that $T\zdk^\alpha(E)$ is a subbundle of $TSO(E)_{|\zdk^\alpha(E)}$ and that we have
\begin{eqnarray*}
TSO(E)_{|\zdk^\alpha(E)} & = & \maB_0(\pi_{\mal Z}^*E)\oplus\maB_*(\pi_{\mal Z}^*E)\oplus \hor_{|\zdk^\alpha(E)}\\
     & = & \maB_0(\pi_{\mal Z}^*E)\oplus T\zdk^\alpha(E)
\end{eqnarray*}
where \footnote{still with the notation defined in section~\ref{AdJ}}
\begin{eqnarray*}
\maB_0(\pi_{\mal Z}^*E)_J & = & \maB_0(E_x,J)\quad \text{and}\\
\maB_*(\pi_{\mal Z}^*E)_J & = & \maB_*(E_x,J):=\left( \bigoplus_{j\in\Z/r\Z\setminus\{0\}}\negthickspace \negthickspace \mal B_j^\C(E_x,J)\right)  \bigcap\End(E_x)=\ver_J^{\mal Z}
\end{eqnarray*}
for all $J\in\zdk^\alpha(E)$.  In other words, $\pi_{\mal Z}^*E$ is canonically endowed with a $2k$-structure: $\mal{I}_J=J$, $\forall J\in N_{\mal Z}$, and this $2k$-structure defines the spaces $\maB_j^\C(\pi_{\mal Z}^*E):=\maB_j^\C(\pi_{\mal Z}^*E,\mal I)$.\\
Now let us make precise the relation between $SO(E)$ and   $\so(E)$ and in particular the relation $T_JSO(E_x)=J.\so(E_x)$. For $J\in SO(E)$, let 
$$
L_J\colon A\in \End(E_x)\longmapsto J.A\in \End(E_x)
$$ 
be the left multiplication by  $J$ in $\End(E_x)$, with $x=\mrm p(J)$. Let us still denote by $\mal I$, the tautological section of $\mrm p^*SO(E)$ defined by  $\mal I_J=J$, $\forall J\in SO(E)$, and whose restriction to $N_{\mal Z}$ is our canonical $2k$-structure $\mal I$ on $\pi_{\mal Z}^* E$. Then let 
$L_{\mal I} \colon SO(E)\longrightarrow \mrm{Aut}(\End(\mrm p^*E))$ be the section  of the bundle  of linear automorphism of the vector bundle $\End(\mrm p^*E)$ defined by
$$
L_{\mal I}\colon J\in SO(E)\longmapsto L_J\in \mrm{Aut}(\End(E_{\mrm p(J)}))
$$ 
or more concretely
$$
L_{\mal I}\colon (J,A)\in \End(\mrm p^*E)\longmapsto (J,J.A)\in \End(\mrm p^*E).
$$ 
Then we have 
$$
\ver^{SO(E)}=L_{\mal I}(\so(\mrm p^*E))\quad \text{and}\quad \maB_j^\C(\pi_{\mal Z}^*E)=L_{\mal I}(\so_j^\C(\pi_{\mal Z}^*E))
$$
which we will denote more simply by
$$ 
\ver^{SO(E)}={\mal I}.\so(\mrm p^*E)\quad \text{and}\quad \maB_j^\C(\pi_{\mal Z}^*E)={\mal I}.\so_j^\C(\pi_{\mal Z}^*E).
$$
\end{exam}
\begin{exam}\label{zdkjE}\em
\textbf{The Twistor subbundle $\mZ_{2k,j}^\alpha(E)$.}\\
Let us consider the previous example and let us suppose that there exists a (global) section $J_j$ of $(\zdk(E))^j=\mal U_{p'}^{j\cdot\alpha}(E)$ for some $j\in\Z$ and $p'=\dfrac{2k}{(2k,j)}$. Let us consider the subbundle
$$
N_\mal{Z}^j:=\mal Z_{2k,j}^\alpha(E,J_j)=\{ J\in \zdk^\alpha(E)|\, J^j=J_j\}
$$ 
for which we have the natural fibration $\pi_\mZ^j\colon \mZ_{2k,j}^\alpha(E,J_j)\to M$.
The vertical space is given by  
$$
\forall J\in N_\mZ^j,\quad \ver_J^{\mZ,j}=T_J\mZ_{2k,j}^\alpha(E_x,J_j)=\left( \oplus_{q=1}^{(r,j)-1}\mal B_{qp}^\C(E_x,J)\right) \cap \End( E_x) =J.\ul_{j-1}^*(E_x,J)
$$ 
according to (\ref{tangentspace})\footnote{and with the notations of section~\ref{AdJj}, in particular $p=\dfrac{r}{(r,j)}$.}, where 
$$
\ul_{j-1}^*(E_x,J)=\left( \oplus_{q=1}^{(r,j)-1}\so_{pq}^\C(E_x,J)\right) \cap \so(E_x)=\ul_{j-1}(E_x,J)/\ul_0(E_x,J)=\so_0(E_x,J^j)/\so_0(E_x,J).
$$
Furthermore, differentiating the definition equation of $\mZ_{2k,j}(E,J_j)$:  $J^j=J_j$, we obtain: for all (local) section $J$ of $\pi_\mZ^j$,
\begin{equation}\label{nablaJj}
\nabla J^j=\sum_{l+q=j-1}J^l\nabla J\,J^q =\nabla J_j
\end{equation} 
so that 
$$
\nabla J\in \ver^{\mZ,j} \Longleftrightarrow \nabla J_j=0,
$$
therefore in general, we have $\nabla J\notin \ver^{\mZ,j}$. We will simply set
$$
\hor^{\mZ,j}={\ver^{\mZ,j}}^\perp\cap T\mal Z_{2k,j}^\alpha(E,J_j).
$$
Then the splitting 
$$
T\mal Z_{2k,j}^\alpha(E,J_j)=\ver^{\mZ,j}\oplus\hor^{\mZ,j}
$$
do not correspond to the splitting (\ref{splittingsoe}) or equivalently to (\ref{splittingzdk}), in general. In other words, the connection  in $\pi_\mZ$ defined by the horizontal distribution $\hor_{|N_\mZ}$ is not reducible to a connection in $\pi_\mZ^j$ (which could only be $\hor^{\mZ,j}$): it happens \iif $\hor$ is tangent to $N_\mZ^j$. Besides we have two different ways to decompose the orthogonal of $\ver^{\mZ,j}$ in $TN_\mZ$, using the decompositions $TN_\mZ=\ver^\mZ\oplus\hor^\mZ$ or ${TN_\mZ}_{|N_\mZ^j}= TN_\mZ^j\oplus {TN_\mZ^j}^\perp$:
$$
\begin{array}{rcl}
{TN_\mZ}_{|N_\mZ^j} & = & \ver^{\mZ,j}  \oplus  {\ver^{\mZ,j}}^\perp\\
       & = & \ver^{\mZ,j}  \oplus  {\ver^{\mZ,j}}^\perp \cap \ver^\mZ\oplus \hor_{|N_\mZ^j} \\
       & = & \underbrace{\ver^{\mZ,j}  \oplus   \hor_\mZ^j }_{TN_\mZ^j}\oplus {TN_\mZ^j}^\perp
\end{array}
$$    
In particular, we have for any (local) section $J\colon U\subset M\to N_\mZ^j$
$$
[dJ]_{\ver^{\mZ,j}}=[\nabla J]_{\ver^{\mZ,j}}=\mrm{pr}_{\ver^{\mZ,j}}^{\ver^\mZ}(\nabla J)
$$
where $[\ ]_{\ver^{\mZ,j}}\colon TN_\mZ\to \ver^{\mZ,j}$ and $\mrm{pr}_{\ver^{\mZ,j}}^{\ver^\mZ}\colon \ver_\mZ\to \ver^{\mZ,j}$ are resp. the orthogonal projections.
Moreover, let us decompose $T\mal{SO}(E)_{|\mZ_{2k,j}^\alpha(E)}$ into an orthogonal sum making appear the vertical subbundle $\ver^{\mZ,j}$ of $N_\mZ^j$:
$$
\begin{array}{rcllcc}
T\mal{SO}(E)_{|N_\mZ^j}  & = &  \mal B_0({\pi_\mZ^j}^* E, \mal I)  \oplus \dfrac{}{} & \ \mal B_*({\pi_\mZ^j}^* E, \mal I) & \oplus & \hor_{|N_\mZ^j}\\ 
 & = &  \underbrace{\mal B_0({\pi_\mZ^j}^* E, \mal I)  \oplus \dfrac{}{} \ver^{\mZ,j}}_{\mal I.\so_*({\pi_\mZ^j}^* E,\mal I^j)} & \oplus \dfrac{}{}\underbrace{\mal I.\so_*({\pi_\mZ^j}^* E,\mal I^j)}_{{\ver^{\mZ,j}}^\perp \cap \ver^\mZ} & \oplus  & \hor_{|N_\mZ^j}
\end{array}.
$$
Now let us see how we can determine $\hor^{\mZ,j}$ from the section $J_j$. First we remark that $\hor^{\mZ,j}\cap \ver^\mZ=\{0\}$ (indeed $\ker d\pi_\mZ\cap \hor^{\mZ,j}\subset \ker d\pi_\mZ\cap TN_\mZ^{j}=\ver^{\mZ,j}$ and of course $\ver^{\mZ,j}\cap \hor^{\mZ,j}=\{0\}$). Therefore $\hor^{\mZ,j}$ is a vector subbundle of $({\ver^{\mZ,j}}^\perp \cap \ver^\mZ)\oplus \hor$ which satisfies $({\ver^{\mZ,j}}^\perp \cap \ver^\mZ)\cap \hor^{\mZ,j}= \{0\}$. Thus $\hor^{\mZ,j}$ is the graph of some linear map\footnote{i.e. a morphism of vector bundle} 
$\Gamma\colon \hor\to {\ver^{\mZ,j}}^\perp \cap \ver^\mZ$, \footnote{In the following reasoning, we will forget the index "$|N_\mZ^{j}$" in $\hor_{|N_\mZ^{j}}$ to do not weigh down the equations. The right notation will reapear in the final equation.}
$$
\Id+ \Gamma\colon W\in \hor \mapsto W + \Gamma(W)\in \hor\oplus ({\ver^{\mZ,j}}^\perp \cap \ver^\mZ)
$$
has $\hor^{\mZ,j}$ as  image.\\
Let us concentrate ourself on (\ref{nablaJj}). $\nabla J$ is in $\ver^\mZ$ so that we can write it $\nabla J=\sum_{i=1}^{r-1} JA_i$ with $A_i\in\so_i^\C(E,J)$, according to  (\ref{vzj}). Then we have $\forall i\in \{1,\ldots,r-1\}$,
$$
\sum_{l + q=j-1}J^q (JA_i) J^l = \sum_{l=0}^{j-1}J^{j-1} \omega_r^l J A_i=\dfrac{1-(\omega_r^i)^j}{1-\omega_r^i} J^j A_i
$$
so that  
$$
\sum_{l + q=j-1}J^q (\nabla J) J^l = J^j\left( \sum_{\underset{i\notin p.\Z }{i=1}}^{r-1} \dfrac{1-(\omega_r^i)^j}{1-\omega_r^i}\, A_i\right) 
$$
where as usual $p=\dfrac{r}{(r,j)}$ is the order of $\omega_r^j$. In particular, we remark  that (with obvious notation)\footnote{remark that $\mal B_*(E,J_j)=J^j. \so_*(E,J^j)=J^j.\left( \left( \oplus_{i\in \Z_r\setminus p.\Z_r} \so_i^\C(E,J) \right) \cap\so(E,J)\right)  $}
\begin{equation}\label{pjJ}
\sum_{l + q=j-1} L(J^l)\circ R(J^q)^{-1}\colon \mal B_*(E,J)\longmapsto  \mal B_*(E,J_j)
\end{equation} 
is a surjective map with kernel\footnote{Obviously, since $J$ is a local section, everything is local here and $E$ must be replaced by $E_U:=\mrm p^{-1}(U)$, but we do not want to weigh down the notations.} 
$$
\left( \bigoplus_{i\in p.\Z_r\setminus\{0\}} \mal B_i^\C(E,J)\right) \bigcap \End (E) = J^*\ver^{\mZ,j},
$$ 
so that it induces an isomorphism from
$$
J.\so_*(E,J_j)=J^*({\ver^{\mZ,j}}^\perp\cap\ver^\mZ)
$$ 
onto 
$\mal B_*(E,J_j)$. Let us denote by $P^j(J)$ the surjective map (\ref{pjJ}) and by $P^j(J)^{-1}$ the inverse map  of the isomorphism induced on $J.\so_*(E,J_j)$. Then we have
$$
P^j(J)(\nabla J)=\nabla J_j
$$ 
so that $[\nabla J]^{{\ver^{\mZ,j}}^\perp\cap\ver^\mZ} =P^j(J)^{-1}(\nabla J_j)$, but we have $\nabla J = [dJ]^{\ver^\mZ}$, and therefore
\begin{equation}\label{dJvperp}
[d J]^{{\ver^{\mZ,j}}^\perp\cap\ver^\mZ} =P^j(J)^{-1}(\nabla J_j).
\end{equation}
On the other hand, $d\pi^\mZ\circ dJ=\Id_{TM}$ so that $d\pi^\mZ\circ [dJ]^\hor=\Id_{TM}$, which with (\ref{dJvperp}) allows to conclude that
$$
\Gamma =P^j(J)^{-1}\circ(\nabla J_j)\circ d\pi^\mZ_{|\hor}
$$
that is to say, for all $W\in \hor_{|N_\mZ^{j}}$
$$
\Gamma(W)=P^j(J_0)^{-1}\cdot \nabla_{(\pi^\mZ)_* W} J_j
$$
where $W=(J_0,W_{J_0})$, $J_0\in N_\mZ^{j}$, $W_{J_0}\in \hor_{J_0}$.
\end{exam}
\index{twistor|)}

\subsubsection{$\Psi$-torsion, $\Psi$-difference tensor, and curvature of a Pfaffian system}\label{phitorsion}
\paragraph{$\Psi$-torsion, $\Psi$-difference.}
Let us consider a vector bundle morphism
$$
\begin{CD}
TM @>\Psi>> (E,\nabla)\\
@VVV @VVV\\
M @>\psi>> N 
\end{CD},
$$
$\nabla$ being a connections on the vector bundle $E$. Then the $\Psi$-torsion of $\nabla$ is the $\psi^*E$-valued 2-form on $M$,
$$
T^\Psi(X,Y) =\nabla_X(\Psi Y)-\nabla_Y(\Psi X)-\Psi[X,Y]= d^\nabla\Psi(X,Y)\quad  \forall X,Y\in\mal C(TM).
$$
Let us give now some examples.
\begin{exam}\em
Let $N$ be a  manifold and suppose that we have a splitting $TN=\ver\oplus \hor$ and suppose also that the vertical bundle $\ver$ is endowed with a covariant derivative $\nabla^c$ and let $\phi\colon TN\to \ver$ be the projection (morphism) on $\ver$ along $\hor$, then we can speak about the $\phi$-torsion of $\nabla^c$,           $T^\phi=d^{\nabla^c}\phi$.
{\tiny
$$
\begin{CD}
TN @>\phi>> (\ver,\nabla^c)\\
@VVV @VVV\\
N @>\Id>> N
\end{CD},
$$ }
\end{exam}
\begin{exam}\em
Let $s\colon M\to (N,\nabla)$ be a map from a  manifold $M$ into an affine manifold $(N,\nabla)$ and suppose that  we have  a splitting $TN=\ver\oplus\hor$, then let us  consider the morphism of bundle
$$
\begin{array}{lr}
\begin{CD}
TM @>d^v s>> (\ver,\nabla^v)\\
@VVV @VVV\\
M @>s>> N
\end{CD}
& \begin{array}{c}
\text{where } \nabla^v \text{ is the vertical part}\\
\text{of the linear connection} \nabla. 
\end{array}
\end{array}
$$ 
Then the \emph{vertical $s$-torsion} of $N$ is  $T^s:=T^{d^vs}= d^{\nabla^v} d^v s$.
\end{exam}
\begin{exam}\em
In particular let us take $s=\Id_N$ (in the previous example) and thus $d^vs=\phi$ the projection on $\ver$ and then the $\phi$-torsion of $\nabla^v$ (or $\Id_N$-torsion of $N$) is the \emph{vertical torsion} in $\ver$: $T^v=d^{\nabla^v}\phi$.\\
Now for any map $s\colon M\to N$ we have 
$$
T^s=s^*T^v.
$$
We will say that $s$ is \emph{vertically torsion free} if $T^s=0$.
\end{exam}
\begin{exam}\em
Let $E\to N$ be a vector bundle and let us suppose that we have a isomorphism of bundle $\psi\colon TN\to E$ (over $\Id_N$). Let $\nabla$ be a connection  in $E$. Then we have $T^\psi=\psi\circ T$, where $T$ is the torsion of the linear connection $\psi^{-1}\circ\nabla\circ \psi$ on  $N$. 
$\blacksquare$
\end{exam}
Now, we define the $\psi$-difference and the $\psi$-equivalence.
\begin{defn}
Let $E\to N$ be a vector bundle and let us suppose that we have a morphism of bundle $\psi\colon TN\to E$ (over $\Id_N$). Given two connections $\nabla$ and $\nabla'$ in $E$ , the \textbf{$\psi$-difference tensor} $S^\psi$ for the pair  $(\nabla,\nabla')$ is defined by
$$
S^\psi(X,Y)=\nabla_X(\psi Y) - \nabla_X'(\psi Y)=(\nabla -\nabla')_X(\psi Y).
$$
Then $S^\psi$ is symmetric precisely when $\nabla$ and $\nabla'$ have the same $\psi$-torsion. On the other hand, if $S^\psi$ is skew-symmetric we will say  (following \cite{cmw2}) that $\nabla$ and $\nabla'$ are \textbf{$\psi$-equivalent}: it means  that these have the same $\psi$-geodesics, a $\psi$-geodesic of $\nabla$ being a path $y(t)$ in $N$ solution of the equation
$$
\nabla_{y'(t)} (\psi y'(t))=0
$$
(if $\psi$ is an isomorphism then $\psi$-geodesics of $\nabla$ are precisely the geodesics of $\psi^{-1}\circ\nabla\circ\psi$).
\end{defn}
\index{curvature, of a subundle|(}
\paragraph{Curvature of a Pfaffian system}
\begin{defn}\label{CurvaturPfaffian}
Let $\mal P$ be a Pfaffian system on the manifold $N$. Then for any local sections of $\mal P$, $X,Y\colon (N,a)\to \mal P$, defined in the neighbourhood of $a\in N$, the image $([X,Y]_a)_\ver$ of $[X,Y]_a$ by the canonical projection $T_aN\to \ver_a=T_aN/\mal P_a$, depends only on the values $X_a,Y_a$ at $a\in N$, of the vector fields $X,Y$. We define the curvature of $\mal P$ as the tensor $\mR\in \Lm^2\mal P^*\otimes \ver$,
$$
 \mR_a(X_a,Y_a):=-([X,Y]_a)_\ver.
$$
\end{defn}
\begin{defn}\label{courburedef}
Let $N$ be a manifold endowed with a Pfaffian system $\ver$ ("vertical subbundle") and let us suppose that $\ver$ admits a connection i.e. a complement $\hor$ ("horizontal subbundle"): $TN=\ver\oplus\hor$. Then $\ver$ is identified to $TN/\mal \hor$ so that the curvature of the connection $\hor$ becomes the tensor $\mR\in\Lm^2\hor^*\otimes \ver$ defined by
$$
\mR(X,Y)=-[X,Y]_\ver\quad \forall X,Y\in\mal C(\hor)
$$
the subscripts "$\ver$" denoting the $\ver$-component along $\hor$.
\end{defn}
\textbf{Convention} We will often extend $\mR$ to the corresponding horizontal 2-form on $N$, still denoted by $\mR$, i.e. $\mR\in\Lm^2T^*N\otimes\ver$ such that $\mR(X,Y)=0$ if $X$ or $Y\in\ver$.
\begin{prop}\label{Tc-Phi}
Let $N$ be a  manifold and suppose that we have a splitting $TN=\ver\oplus \hor$ and suppose also that the vertical bundle $\ver$ is endowed with a covariant derivative $\nabla^c$, then we have
$$
T_{|\hor\times\hor}^c=\mR_\hor,
$$
$\mR_\hor$ being the curvature of  $\hor$.
\end{prop}
\textbf{Proof.} For any $H_1,H_2\in\hor$, we have $T^c(H_1,H_2)= \nabla^c_{H_1} H_2^v - \nabla_{H_2} H_1^v -[H_1,H_2]^v= -[H_1,H_2]^v = \mR_\hor(H_1,H_2)$. \hfill$\square$
\begin{prop}
Let $\pi\colon Q\to Q/H=M$ be  a $H$-principal bundle endowed with a connection 1-form $\omega\colon TQ\to \hk$. Let $\hor=\ker\omega\subset TQ$ be the corresponding horizontal subbundle. Let be $\Omega=d\omega + \frac{1}{2}[\omega\wedge\omega]$ the curvature 2-form. Then we have 
$$
(\mR_\hor)_q(X,Y)=q.\Omega_q(X,Y)\quad \forall q\in Q,\,\forall X,Y\in \hor_q
$$
$\mR_\hor$ being the curvature of the connection $\hor$. In other words, we have 
$$
\mR_\hor=\Omega^*
$$
where $\Omega_q^*=q.\Omega_q$.
\end{prop}
\subsection{Harmonic sections of homogeneous fibre bundles}\label{Homfibrbund}
\index{homogeneous fibre bundle|(}
In this section, we study fibre bundles $\pi\colon N\to M$ for which the fibre is a homogeneous space $H/K$. To do that, we follow, in subsections~\ref{geometry} and \ref{vert-harm-eqn},  the exposition of \cite{cmw2}. Next, we present a generalisation of the results (of \cite{cmw2}) to non section maps in the end of \ref{vert-harm-eqn}. Finally we study the  homogeneous fibre bundle reductions in \ref{Reductions}.
\subsubsection{Definitions and Geometric properties}\label{geometry}
Let $\pi_M\colon Q\to M$ be a principal $H$-bundle, with $H$ a Lie group. Let $K$ be a Lie subgroup of $H$ and $N=Q/K$. Then the map $\pi_N\colon Q \to N$ is a principal $K$-bundle and we have $\pi_M=\pi\circ\pi_N$ where $\pi\colon N\to M$ is a fibre bundle with fibre $H/K$, which is naturally isomorphic to the associated bundle $Q\times_H H/K$. We assume the following hypothesis
\begin{description}
\item[(i)] $H/K$ is reductive: $\hk=\kk \oplus \pk$, and $\Ad K(\pk)\subset \pk$, where $\hk$ and $\kk$ are respectively the Lie algebras of $H$ and $K$.
\item[(ii)] $M$ is endowed with a Riemannian metric $g$
\item[(iii)] $H/K$ is Riemannian: there exists  a $H$-invariant Riemannian metric on $H/K$ (equivalently an $\Ad K$-invariant (positive definite) inner product on $\pk$). Equivalently $\Ad_\pk K$ is compact.
\item[(iv)] The principal $H$-bundle $\pi_M\colon Q\to M$ is endowed with a connection. We denote by $\omega$ the corresponding $\hk$-valued connection form on $Q$.

\end{description}
Then the splitting $TQ=\ver_0\oplus\hor_0$ defined  by $\omega$ ($\ver_0=\ker d\pi_M$, $\hor_0=\ker\omega $) gives rise by $d\pi_N$, to the following decomposition $TN=\ver\oplus\hor$, where $\ver=\ker d\pi=d\pi_N(\ver_0)$ and 
$\hor=d\pi_N(\hor_0)$. Let $\pk_Q:=Q\times_K\pk\to N$ be the vector bundle associated to $\pi_N\colon Q\to N$ with fibre $\pk$. Let us denote by $[q,a]\in \pk_Q$ the element defined by $(q,a)\in Q\times \pk$. Then  we have the following vector bundle isomorphism 
$$
\begin{array}{crcl}
I  \colon & \ver & \longrightarrow & \pk_Q\\
        & d\pi_N(q.a) & \longmapsto & [q,a]
\end{array}
$$
 where $q\in Q$, $a\in \pk$ and  as usual $q.a=\dfrac{d}{dt}_{|t=0}q.\exp(ta)\in T_q Q$. Decomposing $\omega =\omega_\hk + \omega_\pk$ following $\hk=\kk\oplus\pk$, then since $H/K$ is reductive, $\omega_\pk$ is a $K$-equivariant ($\omega_\pk(X.h)=\Ad h\, \omega_\pk(X)$) and $\pi_N$-horizontal (${\omega_\pk}_{|\ver_0}=0$) $\pk$-valued 1-form on $Q$ and hence projects to a $\pk$-valued  1-form $\phi$ on $N$:
$$
\phi( d\pi_N(X))=[q,\omega_\pk(X)].
$$
Then we have 
$$
\phi_{|\ver}=I \quad \text{and} \quad  \ker\phi =\hor.
$$
We can  now  construct  a Riemannian metric $h$ on $N$:
\begin{equation}\label{metrich}
h=\pi^*g + \langle \phi,\phi \rangle
\end{equation}
where $\langle\ ,\ \rangle$ is the fibre metric induced on $\pk_Q$ by the inner product  on $\pk$. \\[1mm]
In the same way, let $\Phi$ be the  $\pk_Q$-valued 2-form on $N$ defined by the component $\Omega_\pk$ of the curvature form $\Omega$ of $\omega$. Since $\Omega_\pk$ is $\pi_M$-horizontal ($\Omega(X,Y)=0$ if $X\in \ver_0$ or $Y\in \ver_0$),  then $\Phi$ is $\pi$-horizontal: $\Phi(X,Y)=0$, if $X\in \ver$ or $Y\in \ver$.
\begin{rmk}\em
In \cite{cmw2}, $\pk_Q$ is called the \emph{canonical bundle}, $I$ the \emph{canonical isomorphism}, $\phi$ the \emph{homogeneous  connection form}, and $\Phi$ the \emph{homogeneous curvature form}.
\end{rmk}
\index{curvature, of a subundle|)}
\begin{defn}
We will call  the following datas $(Q,H,K,\omega)$ a homogeneous fibre bundle structure on $\pi\colon N\to M$ (or on $N$, when the fibration $\pi$ is    considered as implicitely given). Moreover $\pi\colon N\to M$ will then be sayed to be a homogeneous fibre bundle.
\end{defn}
The 1-form $\omega_\kk$ (which is a connection form in $\pi_N$ because $H/K$  is reductive) defines a connection in $\pi_N$ called the\emph{ canonical connection}. This connection induces  a covariant derivative $\nabla^c$ in the  associated bundle $\pk_Q$, with respect  to which the fibre metric is parallel. $\nabla ^c$ defines a exterior derivative $d^c$ on the space of $\pk_Q$-valued differential forms on $N$. This allows us to define the \emph{canonical torsion} $T^c$ which is nothing but the $\phi$-torsion of $\nabla^c$ (see section~\ref{phitorsion})
\begin{equation}\label{Tc}
T^c(A,B)=d^c\phi(A,B)=\nabla_A^c(\phi B) -\nabla_B^c(\phi A) -\phi[A,B],\quad \forall A,B \in \mal C(TN)
\end{equation}
Let $\hk_Q:=Q\times_H\hk\to M$ be the vector bundle associated to $\pi_M$ with fibre $\hk$, and in the same way $\kk_Q:=Q\times_K \kk\to N$ the bundle associated to $\pi_N$ with fibre $\hk$. Then we have 
$$
\pi^*\hk_Q=\kk_Q\oplus \pk_Q.
$$
The Lie bracket of $\hk$ induces a bracket on the fibres of $\hk_Q$, and those of $\pi^*\hk_Q$, which we continue to denote  by $[\ ,\ ]$  , and we denote also its $\pk_Q$-component (when there is no risk of confusion) by $[\ ,\ ]_\pk$ (otherwise we denote it by $[\ ,\ ]_{\pk_Q}$). Taking the $\pk$-component of the structure equation $d\omega=\Omega -\dfrac{1}{2}[\omega\wedge\omega]$ and then projecting on $N$, we obtain  \emph{the homogeneous structure equation}:
\begin{equation}\label{hme}
T^c=\Phi-\dfrac{1}{2}[\phi\wedge\phi]_\pk
\end{equation}
and thus
$$ 
\begin{array}{lc}
T^c_{|\ver\times\ver}\dfrac{}{}  =  -[I\cdot, I\cdot]_\pk, &  T^c_{|\ver\times\hor}=0\\
T^c_{|\hor\times\hor}  =\dfrac{}{}  \Phi_{|\hor\times\hor}.  &   
\end{array}
$$
In particular, $T^c$ is horizontal \iif $H/K$ is a (locally) symmetric space, and in this case 
\begin{equation}\label{tcphi}
T^c=\Phi
\end{equation}
\begin{rmk}
\em According to (\ref{hme}) and (\ref{Tc}),  for all $X,Y\in\hor$, (extended to vector fields in $N$ denoted by the same letters), we have
$$
\Phi(X,Y)=T^c(X,Y)= -\phi([X ,Y])
$$
so that 
$$
\Phi=\mR_\hor,
$$
according to definition~\ref{courburedef}. The homogeneous curvature form is nothing but the curvature of the connection $\hor$. $\centerdot$
\end{rmk}
Now, let $\mU$ be the $\pk_Q$-valued symmetric bilinear form defined on $\pk_Q$ by:
\begin{equation}\label{defofmU}
\langle \mU(a,b), c\rangle =\langle[c,a]_\pk,b\rangle + \langle a, [c, b]_\pk\rangle
\end{equation}
where $\langle\ ,\ \rangle$ is the fibre metric, and $a,b,c\in\pk_Q$. Let us set 
\begin{equation}\label{eq-B=U+croch}
\mB=\mU+ [\ ,\ ]_\pk
\end{equation}
which is a $\pk_Q$-valued bilinear form on $\pk_Q$, whose the symmetric  and skew symmetric components are respectively $\mU$ and $[\ ,\ ]_\pk$. $\mU$ vanishes \iif $H/K$ is naturally reductive and $\mB$ \iif $H/K$ is (locally) symmetric . Then denoting by $\nabla$ the Levi-Civita connection on $N$, we have:
\begin{thm}\label{difference-tensor}\cite[p. 197]{cmw2}
Let us consider the difference tensor:
$$
S(A,B)=\phi(\nabla_A  B) - \nabla_A^c(\phi B)
$$
then we have 
$$
2S=\phi^*\mU - T^c = \phi^*\mB - \Phi.
$$
Consequently, $\forall V\in\mal C(\ver)$ 
\begin{equation}\label{egality}
I(\nabla_A^v V)=\nabla_A^c (I V) + \frac{1}{2}\mB(\phi A, I V).
\end{equation}
In particular, if $H/K$ is a (locally) symmetric space, we have 
$$ 
I \nabla_A^v V =\nabla_A^c (IV).
$$
\end{thm}
\begin{rmk}\em
If $H/K$ is a symmetric space, under the canonical identification $I\colon \ver \overset{\simeq}{\longrightarrow} \pk_Q$, we have $\nabla^v=\nabla^c$ on $\ver$. More generally the difference between $\nabla^v$ and $\nabla^c$ looks like to the difference between the Levi-Civita and canonical connections of a reductive Riemannian homogeneous space (see section~\ref{family}).\\
Moreover, $\nabla^v$ is $\phi$-equivalent to $\nabla^c$ when $H/K$ is naturally reductive, according to (\ref{egality}).
\end{rmk}
Let $\nabla^\omega$ be the covariant derivative in the vector bundle $\hk_Q$ (associated to $\pi_M$), defined by the connection form $\omega$. Let us decompose (the $\pi$-pullback of) $\nabla^\omega$ following the decomposition $\pi^*\hk_Q=\kk_Q\oplus \pk_Q$, and the $\pk_Q$-component gives us a connection $\nabla^\pk$ in $\pk_Q$.
\begin{thm}\label{nabla-p}\cite[p.~198]{cmw2}
For all $V\in\mal C(\pk_Q)$,
$$
\nabla^\pk V=\nabla^\omega V- [\phi,\alpha]_\hk
$$
and 
$$
\nabla^c V=\nabla^\pk V-[\phi, V]_\pk=\nabla^\omega V - [\phi, V].
$$
Consequently, $\nabla^\pk$ and $\nabla^c$ are $\phi$-equivalent (since their $\phi$-difference  is $[\phi,\phi]_\pk$). In particular $\nabla^c=\nabla^\pk$ if $H/K$ is a (locally) symmetric space.
\end{thm}
\begin{exam}\label{2.1}\em
Let us consider the situation described by example~\ref{example1} and suppose that $u\colon L\to M$ is an isometry. Then if $\pi\colon N\to M$ is a homogeneous fibre bundle like above then this is also the case for $u^*\pi\colon u^*N\to L$.\\
Indeed  let us set  
$$
u^*Q=\{(z,q)\in L\times Q,\ q\in \pi_M^{-1}(\{u(z)\})\,\}=\bigsqcup_{z\in L}\{z\}\times f(z).H
$$
then $u^*\pi_M\colon (z,q)\in u^*Q\mapsto z\in L$ is a principal $H$-bundle over $L$. Then we have $u^*N=u^*Q/K$, 
and $u^*\pi \colon u^*N \to L$ is a fibre bundle with fibre $H/K$.\\
Finally we have to define a connection on $u^*\pi_M\colon u^*Q\to L$. Let us extend the connection $\omega$, to a connection on $\Id_L\times \pi_M \colon L\times Q\to L\times M$ by $\tl\omega_{(z,q)}(dz + dq)=\omega_q(dq)$ and then let us set 
$$
u^*\omega:=\tl\omega_{|T(u^* Q)}.
$$
In the same way, the homogeneous connection and curvature forms on $u^*N$ are given respectively by 
$$
u^*\phi:=\tl\phi_{|T(u^*N)} \quad \text{and}\quad  u^*\Phi:=\tl\Phi_{|T(u^*N)\oplus T(u^*N)}.
$$
The canonical torsion $T^c$ on $u^*N$ is also given by $u^*T^c:=\tl T^c_{|T(u^*N)\oplus T(u^*N)}$.
\end{exam}
\begin{defn}
Let $(L,b)$ be a Riemannian manifold. We will then say that $f\colon L\to N$ is vertically geodesic if $\phi(\nabla df)=0$, horizontally geodesic if $\pi_*(\nabla df)=0$, and superflat if its \emph{vertical second fondamental form} vanishes $\Pi^v (f):=\nabla^v d^vf=0$.
\end{defn}
\subsubsection{Vertical harmonicity equation}\label{vert-harm-eqn}
We know that the structure group $H$  of $\pi_M\colon Q\to M$ is reducible to $K$ (i.e. there exists an $K$-bundle $\pi_M'\colon Q'\to M$) \iif the associated bundle $\pi \colon N\to M$ admits a (global) section $s\colon M\to N$ (see
\cite{KN}) so that there is a one to one correspondance between the $K$-reductions of $\pi_M$ and the space of sections $\mal C(\pi)$.\\
Let $\omega'={\omega_\hk}_{|TQ'}$. Then $\omega'$ is a connection in $\pi_M'$, and $\omega$ is reducible \iif $\omega_{|TQ'}=\omega'$ (see \cite{KN}). The reducibility of $\omega$ can be characterized as follows.
\begin{prop}
The following statements are equivalent:\smallskip\\
\textbf{\em (i)} $s$ is horizontal;\\
\textbf{\em (ii)} $s^*\phi=0$;\\
\textbf{\em (iii)} $s$ is an isometric immersion;\\
\textbf{\em (iv)} $\omega$ is reducible
\end{prop}
Now we have the following expression of the tension field for sections $s\colon M\to N$.
\begin{thm}\label{coderivative}\cite[Theorem~3.2]{cmw2}
For all $s\in\mal C(\pi)$,
$$
I(\tau^v(s))=-d^*(s^*\phi) + \frac{1}{2}\mrm{Tr}(s^*\phi^*\mU)
$$
where $d^*$ is the coderivative for $s^*\pk_Q$-valued differential forms on $M$ relative to the $s$-pullback of any connection in $\pk_Q$ which is $\phi$-equivalent to $\nabla^c$. In particular, if $H/K$ is naturally reductive then $s$ is an harmonic section \iif $s^*\phi$ is coclosed. 
\end{thm} 
\begin{rmk}\label{naturreductiv}\em
If $H/K$ is naturally reductive, to compute the vertical tension field $\tau^v(s)=\mrm{Tr}(\nabla^vd^v s)$, we can use instead of $\nabla^v$ any connection in $\ver\cong\pk_Q$ which is $\phi$-equivalent to $\nabla^c$.
\end{rmk}
From the homogeneous structure equation (\ref{hme}), we obtain
$$
s^*\Phi=d^c(s^*\phi) + \dfrac{1}{2}[s^*\phi\wedge s^*\phi]_\pk,
$$
hence every horizontal section is flat (i.e. $s^*\Phi=0$).%\\
Let us introduce the following 3-covariant tensor $\langle s^*\phi\otimes s^*\Phi\rangle$ on $M$: 
$$
\langle s^*\phi\otimes s^*\Phi\rangle(X,Y,Z)=\langle s^*\phi(X), s^*\Phi(Y,Z)\rangle,
$$
then we have 
\begin{thm}\label{phinablads}\cite[Theorem~3.4]{cmw2}
For all $s\in\mal C(\pi)$ we have 
\begin{description}
\item[(i)] $ \phi(\nabla ds)=\nabla^c(s^*\phi) + \dfrac{1}{2}s^*\phi^*\mB -\dfrac{1}{2}s^*\Phi$.\\
In particular, if $s$ is vertically geodesic then $s$ is a harmonic section. 
\item[(ii)] $2g(\pi_*\nabla ds (X,Y), Z)=\langle s^*\phi\otimes s^*\Phi\rangle(X,Y,Z) + \langle s^*\phi\otimes s^*\Phi\rangle(Y,X,Z)$.\\
Therefore $s$ is horizontally geodesic \iif $\langle s^*\phi\otimes s^*\Phi\rangle$ is a 3-form on $M$. In particular, if $s$ is flat then $s$ is horizontally geodesic.
\end{description}
\end{thm}
\begin{thm}\label{tauv-tau}\cite[Theorem~3.5]{cmw2}
\begin{description}
\item[(i)] The symmetric and skew symmetric components of $\Pi^v s:=\nabla^v d^v s$ are given by:
$$
I(\Pi^vs)=\phi\circ\nabla ds + \dfrac{1}{2}s^*\Phi.
$$
\item[(ii)] The section $s$ is superflat \iif $s$ is flat and totally geodesic. In particular, if $s$ is flat then  $s$ is totally geodesic \iif $s$ is super-flat.
\item[(iii)] Moreover $\tau^v(s)$ is the vertical component of the tension field $\tau(s)$. So if $s$ is  an harmonic map, then it is certainly a harmonic section.
\end{description}
\end{thm}
\begin{thm}\label{harmonicsection-map}\cite[Theorem~3.6]{cmw2}
An harmonic section $s$ is a harmonic map \iif $\langle s^*\phi, s^*\Phi\rangle=0$ where
$$
\langle s^*\phi, s^*\Phi\rangle(X)=\sum_i\langle s^*\phi\otimes s^*\Phi\rangle(E_i,E_i,X)
$$
for any orthonormal tangent frame $(E_i)$ of $M$.\\
In particular, if $s$ is flat ($s^*\Phi=0$) then $s$ is a harmonic map \iif $s$ is a harmonic section.
\end{thm}
\begin{rmk}\em
Let us consider the situation described by examples~\ref{example1} and \ref{2.1}. Then if $f^*\Phi=0$, $f\colon L\to N$ is vertically harmonic \iif $\tl f \colon L\to u^*N$ is an harmonic section \iif $\tl f\colon L\to u^*N$ is an harmonic map. But it does not imply that $f\colon L\to N$ is harmonic! (See the Appendix.) Indeed in the previous theorem it is essential that $s$ be a section: $\pi\circ s=\Id$.
\end{rmk}
In fact the previous theorems can be easily generalized for non section map . The proofs in \cite{cmw2} holds without any change for theorems \ref{coderivative}, \ref{phinablads}-(i), \ref{tauv-tau}-(i,\,iii), while for  theorems~\ref{phinablads}-(ii), \ref{tauv-tau}-(ii), \ref{harmonicsection-map}: follow the proof of \cite{cmw2}, just replace the starting equation $\pi\circ s=\Id$ by $\pi\circ s= u$. Then we obtain
\begin{thm}\label{m,n}
For all $s\in C^\infty(M,N)$, we have 
\begin{description}
\item[(i)] $I(\tau^v(s))=-d^*(s^*\phi) + \dfrac{1}{2}\mrm{Tr}(s^*\phi^*\mU)$
\item[(ii)] $ \phi(\nabla ds)=\nabla^c(s^*\phi) + \dfrac{1}{2}s^*\phi^*\mB - \dfrac{1}{2}s^*\Phi$.\\[1mm] 
In particular if $s$ is vertically geodesic then $s$ is a harmonic section.
\item[(iii)] $I(\Pi^vs)=\phi\circ\nabla ds + \dfrac{1}{2}s^*\Phi.$\\[1.5mm]
The map $s$ is superflat \iif $s$ is flat and vertically geodesic. Moreover $\tau^v(s)$ is the vertical component of the tension field $\tau(s)$. So if $s$ is  an harmonic map, then it is certainly vertically harmonic.
\item[(iv)] $ 2g(\pi_*\nabla ds (X,Y), u_*Z)=\langle s^*\phi\otimes s^*\Phi\rangle(X,Y,Z) + \langle s^*\phi\otimes s^*\Phi\rangle(Y,X,Z)+ 2g\left(\nabla du(X,Y),u_*Z\right)$.\\[1.5mm]
Let us suppose now that $u$ is an immersion, then this equation determines the horizontal part of $\nabla ds$. 
In particular, if $s$ is flat then $s$ is horizontally geodesic  \iif $u$ is totally geodesic;
  and $s$ is totally geodesic \iif $s$ is superflat and $u$ is totally geodesic.
\item[(v)] A vertically harmonic map $s$ is a harmonic map \iif
$$
g(\tau(u),\cdot ) + \langle s^*\phi, s^*\Phi\rangle=0.
$$
In particular if $s$ is flat, then $s$ is a harmonic map \iif $s$ is vertically harmonic and $u=\pi\circ s$ is harmonic.
\end{description}
\end{thm} 
We could also deduce this generalisation from the previous theorems~\ref{coderivative}-\ref{harmonicsection-map} themself. Indeed we can apply these to the section $\tl s\in \mal C(u^*N)$ corresponding to $s$ and use theorems~\ref{same-vertical-Harm} and \ref{same2form} in the Appendix, but we must suppose in addition that $u$ is an isometry. \\
Let us go further in the generalisation and consider maps $f\in\mal C^\infty(L,N)$ with$(L,b)$ a Riemannian manifold (see examples~\ref{example1} and \ref{2.1}). Then the proofs in \cite{cmw2} holds for theorems~\ref{coderivative}, \ref{phinablads}-(i), \ref{tauv-tau}-(i,\,iii), whereas theorems~\ref{phinablads}-(ii), \ref{tauv-tau}-(ii), \ref{harmonicsection-map}  are no longer valid. Indeed the equation in theorem~\ref{m,n}-(iv) holds, but it gives us only $[\pi_*\nabla df]_{u_*TL}$, the component of  $[\pi_*\nabla df]$ in the tangent bundle $u_*TL$. So if we want $[\pi_*\nabla df]_{(u_*TL)^\perp}$  we must introduce the 3-linear form $\langle f^*\phi\otimes \Phi(f_*\cdot,\cdot) \rangle \in \mal C(T^*L\otimes T^*L\otimes f^*\hor)$ defined by :
$$
\langle f^*\phi\otimes \Phi(f_*\cdot,\cdot) \rangle(a,b,Z)=\langle f^*\phi(a),\Phi(f_*a, Z)\rangle.
$$
Then we  have
\begin{thm}\label{thm-LtoNvertHarm-vs-Harm}
For all $f\in\mal C^\infty(L,N)$, we have 
\begin{description}
\item[(i)] $g(\pi_*\nabla df (a,b), \pi^*Z)=\langle f^*\phi\odot \Phi(f_*\cdot,\cdot)\rangle(a,b,Z) +
g(\nabla du(a,b), \pi_*Z)$.\\[1.5mm]
In particular, if $f$ is strongly flat i.e. ($\Phi(f_*\cdot,\cdot)=0$) then:\\
- $f$ is horizontaly geodesic \iif $u$ is totally geodesic.\\
- $f$ is totally geodesic \iif $f$ is superflat and $u$ is totally geodesic.
\item[(ii)] A vertically harmonic map $f$ is a harmonic map \iif 
$$
g(\tau(u),\cdot) + \langle f^*\phi,\Phi(f_*\cdot,\cdot)\rangle=0
$$
where $\langle f^*\phi,\Phi(f_*\cdot,\cdot)\rangle(X)=\sum_i\langle f^*\phi\otimes \Phi(f_*\cdot,\cdot)\rangle(e_i,e_i,X)
$ 
for any tangent frame $(e_i)$ of $L$. In particular if $f$ is strongly flat ($\Phi(f_*\cdot,\cdot)=0$) then $f$ is a harmonic map \iif $f$ is vertically harmonic and $u$ is harmonic. 
\end{description}
\end{thm}                           
Examples and illustrations of  the theorems~\ref{thm-LtoNvertHarm-vs-Harm} and \ref{m,n} will be given in section~\ref{sect-harm-vs-vertharm}.
%%%%%%%%%%%%%%%%%%%%%%%%%%%%%%%%%%%%%%%%%
%

\subsubsection{Reductions of homogeneous fibre bundles}\label{Reductions}
We  consider reductions  $(Q^\mv,H^\mv,K^\mv,\omega^\mv)$  of  $(Q,H,K,\omega)$. We distinguish two cases: first we consider the particular case $K^\mv\subset K$ and then, in a second time, we explain what does change in the general case $K^\mv\subset K$. Finally, we specify which notion depends or does not depend on the choice of the metric on the fibre. 
\paragraph{When the subgroup $K$ is the same.}
Let us suppose now that the structure group $H$ of $\pi_M\colon Q\to M$ is reducible to a (closed) subgroup $H^\mrm{v}\supset K$. That is to say, there exists a principal $H^\mrm{v}$-subbundle $\pi_M^{\mrm v}\colon Q^{\mrm v} \to M$. Since $H^{\mrm v}/K$ is Riemannian, then it is reductive. Moreover, remarking that $\dim(\hk^\mv\cap \pk)=\dim \hk^\mv +\dim\pk -\dim(\hk) =\dim \hk^\mv -\dim \kk$, we deduce that $\pk^\mv:=\hk^\mv\cap\pk $ satisfies $\hk^{\mrm v}=\kk\oplus\pk^{\mrm v}$ and $\Ad k(\pk^{\mrm v})=\pk^{\mrm v}$, $\forall k\in K$. The restriction to $\pk^{\mrm v}$ of the $\Ad K$-invariant inner product on $\pk$ defines a $H^{\mrm v}$-invariant metric on $H^{\mrm v}/K$ which is nothing but the metric induced by the $H$-invariant metric on $H/K$, so that the inclusion $H^{\mrm v}/K\to H/K$ is an isometric embedding.\\
Let $\pk'=(\pk^{\mrm v})^\perp$ in $\pk$ and let us suppose that $\pk'$ is $\Ad H^{\mrm v}$-invariant, therefore $\hk=\hk^{\mrm v}\oplus \pk'$ is a reductive decomposition and $H/H^{\mrm v}$ is reductive.\smallskip\\
Conversely if $H/H^{\mrm v}$ is reductive: $\hk=\hk^{\mrm v}\oplus \pk'$ with $\pk'$ $\Ad H^\mrm{v}$-invariant, then we can always complete any $\Ad K$-inner product  in $\pk'$ by  an $\Ad K$-invariant inner product in $\pk:=\pk^\mrm{v}\oplus\pk'$:
\begin{equation}\label{eq-inner-p'+pv}
\langle\ ,\ \rangle_\pk= \langle\ ,\ \rangle_{\pk^\mrm{v}} + \langle\ ,\ \rangle_{\pk'}
\end{equation}
for which $\pk'=(\pk^\mrm{v})^\perp$ in $\pk$. Of course $\pk^\mv$ denotes an $\Ad K$-invariant summand such that: $\hk^\mv=\kk\oplus \pk^\mv$ ($H^\mv/K$ is reductive) and $\langle\ ,\ \rangle_{\pk^\mrm{v}}$ is an $\Ad K$-invariant inner product.\medskip\\
$\bullet$ \textbf{In the following} we suppose that $H/H^{\mrm v}$ is reductive and that the inner product in $\pk$ is chosen as described above, i.e according to \eqref{eq-inner-p'+pv}. According to the previous discussion, this is equivalent to suppose that the $\Ad K$-invariant inner product chosen in $\pk$ is such that the subspace $\pk'=(\pk^{\mrm v})^\perp$ is $\Ad H^{\mrm v}$-invariant, where $\pk^\mv:=\hk^\mv\cap\pk $.%\smallskip\\
\begin{rmk}\label{rmk-natural-red-autom}\em
Let us now consider the particular case where $H/K$ is naturally reductive. Then in this case, according to \ref{appendix-natred} in the Appendix, any naturally reductive inner product  on $\pk$ can be extended into an $\Ad H$-invariant PseudoRiemannian product $B$ on $\hk$. Then in this case,  the hypothesis above are automatically satisfied since then $\pk':=(\hk^\mv)^\perp$, the orthogonal of $\hk^\mv$ \wrt $B$, is $\Ad H^\mv$-invariant.
In particular, this is  the case if $H$ is compact.
\end{rmk}

%
%\medskip\\ 
Now let us turn toward the connection 1-form $\omega$. Its restriction $\omega^\mv:=\omega_{\hk^\mv|TQ^\mv}$ defines a connection on $\pi_M^\mv\colon \to M$. We endow $Q^\mv$ with $\omega^\mv$ and $(Q^\mv,H^\mv,K,\omega^\mv)$ is then a homogeneous fibre bundle as defined in the begining of \ref{Homfibrbund}.
\begin{defn}
We will then say that $\pi\colon N^\mv\to M$ is a homogeneous fibre subbundle of $\pi\colon N\to M$.
\end{defn}
Moreover $\omega$ is reducible (to $\omega^\mv$) in $Q^\mv$ \iif one of the following equivalent statements holds (\cite{KN})
\begin{description}
\item[$\bullet$] $\forall q\in Q^\mv$, $(\hor_0)_q$ is tangent to $Q^\mv$.
\item[$\bullet$] $\omega_{|TQ^\mv}=\omega^\mv$ (i.e. $\omega_{|TQ^\mv}$ is $\hk^\mv$-valued).
\item[$\bullet$] The canonical cross section $s^\mv$ of the associated bundle $E^\mv:=Q/H^\mv=Q\times_H(H/H^\mv)$, which defines the $H^\mv$-reduction $Q^\mv$ is horizontal.
\end{description}
\begin{defn}\label{def-red-Hfb}
Denoting by   $i_\mv\colon N^\mv \to N$ the natural inclusion,  we will say that $i_\mv\colon N^\mv \to N$ is a reduction of homogeneous fibre bundles  if $\omega$ is reducible in $Q^v$ (being implicit that the homogeneous fibre structure are respectively $(Q,H,K,\omega)$ and $(Q^\mv,H^\mv,K,\omega^\mv)$).
\end{defn}
The vertical bundle (in $TQ$), $\ver_0$, splits as follows
$$
\ver_0=\ver_0'\oplus \ver_0^\mv
$$
where $(\ver_0^\mv)_q= q.\hk^\mv=T_q (q.H^\mv)$ and $(\ver_0')_q=q.\pk'$, and quotienting by $\kk$, i.e. by applying $d\pi_N$ we obtain the following decomposition of $\ver$:
$$
\ver=\ver'\oplus\ver^\mv
$$
with $\ver'=d\pi_N(\ver_0')$ and $\ver^\mv=d\pi_N(\ver_0^\mv)$.\\
Then the canonical isomorphism $I\colon \ver\to \pk_Q$ sends the previous decomposition onto the following $\pk_Q=\pk_Q'\oplus\pk_Q^\mv$ (i.e. $\ver'$, $\ver^\mv$ are sent resp. onto $\pk_Q'$ and $\pk_Q^\mv$).\\
Then the vertical space in $TN^\mv$ is $\ver_{|N^\mv}^\mv$ that we will also denote by $\ver^\mv$ when there is no possibilities of confusion. The splitting of $TN^\mv$ defined by $\omega^\mv$ is then
$$
TN^\mv=\ver^\mv\oplus\hor_{|N^\mv}^\mv
$$ 
where $\hor^\mv=d\pi_N(\hor_0^\mv)$ and $\hor_0^\mv=\ker \omega_{\hk^\mv}$. Let us remark that $\omega$ is reducible \iif one of the following equivalent statements holds
\begin{description}
\item[$\bullet$] $\hor_{|N^\mv}^\mv=\hor_{|N^\mv}$,
\item[$\bullet$] the inclusion $i_\mv\colon (N^\mv,h^\mv)\to (N,h)$ is an isometry, where $h^\mv$ is the metric defined by $\omega^\mv$ and $\langle\ ,\ \rangle_{\pk^\mrm{v}}$, according to \eqref{metrich}. 
\end{description}
The canonical bundle on $N^\mv$, $\pk_{Q^\mv}^\mv=Q^\mv\times_K \pk^\mv\to N^\mv$ is the restriction to $N^\mv$ of $\pk_Q^\mv\to N$, and the canonical isomorphism $I^\mv\colon \ver_{|N^\mv}^\mv\to \pk_{Q^\mv}^\mv$ is the restriction to $\ver_{|N^\mv}^\mv$ of $I\colon\ver\to \pk_Q$.\\
Since $\omega_\pk=\omega_{\pk^\mv} + \omega_{\pk'}$, the homogeneous connection form on $N^\mv$, $\phi^\mv$  (the $\pk_{Q^\mv}^\mv$-valued 1-form on $N^\mv$ defined by $\omega_{\pk^\mv}$) is the restriction to $N^\mv$ of the 
$\pk_{Q}^\mv$-component of $\phi$:
$$
\phi^\mv = (i_\mv)^*\left[\phi\,\right] _{\pk_Q^\mv}=[(i_\mv)^*\phi\,]_{\pk_{Q^\mv}^\mv},
$$
where $i_\mv\colon N^\mv \to N$  is the natural inclusion.
The homogeneous curvature form $\Phi^\mv$ (defined by $\Omega_{\pk^\mv}^\mv$, with $\Omega^\mv=d\omega^\mv +\dfrac{1}{2}[\omega^\mv\wedge\omega^\mv]$) is given by
$$
\Phi^\mv=(i_\mv)^*\left[\Phi\right] _{\pk_Q^\mv}=[(i_\mv)^*\Phi]_{\pk_{Q^\mv}^\mv}.
$$
Furthermore, $\pk_{Q^\mv}^\mv$ is $\nabla^c$-parallel: the covariant derivative on $\pk_{Q^\mv}^\mv$ defined by $\omega_{\kk|TQ^\mv}$ is the restriction of $\nabla^c$ to  $\pk_{Q^\mv}^\mv$. In other words $\nabla^c$ commutes with the projection on $\pk_Q^\mv$. The canonical torsion on $N^\mv$ is given by 
$$
T^{c,\mv} = d^c\phi^\mv=[(i_\mv)^*T^c\,]_{\pk_Q^\mv}.
$$
In particular, the results of what precedes is:
\begin{prop}\label{prop-tool-corresp}
Let us suppose that $\omega$ is reducible in $Q^\mv$. Then we have 
\begin{description}
\item[(i)] $\phi^\mv = (i_\mv)^*\phi$,
\item[(ii)] $\Phi^\mv = (i_\mv)^*\Phi$,
\item[(iii)] $T^{c,\mv} = (i_\mv)^* T^c$.
\end{description}
\end{prop}
Now, let us remark that 
\begin{lemma}\label{lemma-ivdeB}
We have the following identity
$$
\mB^\mv=(i_\mv)^*\mB
$$
\end{lemma}
\proof  Let $a_{\pk^\mv} \, , b_{\pk^\mv}\in \pk^\mv$ and $c_\pk \in \pk$. Then we have
$$
[ a_{\pk^\mv},b_{\pk^\mv} ]_\pk = [ a_{\pk^\mv}, b_{\pk^\mv} ]_{\pk^\mv}
$$ 
since $\hk^\mv=\kk\oplus\pk^\mv$ is a Lie subalgebra of $\hk=\kk\oplus\pk$. Moreover, we have also
\begin{eqnarray*}
\langle U(a_{\pk^\mv},b_{\pk^\mv}), c_\pk \rangle & = & \langle [c_\pk,a_{\pk^\mv}]_\pk, b_{\pk^\mv}\rangle + \langle [c_\pk,a_{\pk^\mv}]_\pk, b_{\pk^\mv}\rangle \\
& = &  \langle [c_{\pk^\mv},a_{\pk^\mv}]_{\pk^\mv}, b_{\pk^\mv}\rangle + \langle [c_{\pk^\mv},a_{\pk^\mv}]_{\pk^\mv}, b_{\pk^\mv}\rangle 
= \langle U^\mv(a_{\pk^\mv},b_{\pk^\mv}), c_{\pk^\mv} \rangle
\end{eqnarray*}
since $[\pk',\pk^\mv]\subset \pk'$ because the decomposition $\hk=\hk^\mv\oplus\pk'$ is $\Ad H^\mv$-invariant.
We have denoted by $U^\mv$  the map associated with the bundle $\pi\colon N^\mv\to M$ according to equation \eqref{eq-B=U+croch}. This completes the proof.\hfill$\square$\medskip\\
Therefore we  deduce from what precedes and theorem~\ref{difference-tensor}, the following.
\begin{prop}\label{prop-tau&Pi-corresp}
Let us suppose that $\omega$ is reducible in $Q^\mv$. Let $(L,b)$ be a Riemannian manifold, and $f\colon  L\to N^\mv$ be a map. Then we have 
\begin{description}
\item[(i)] $\tau^\mv (f)= \tau^v(f)$,
\item[(ii)] $\Pi^\mv(f) = \Pi^v(f)$,
\end{description}
where $\tau^\mv(f)$ is the vertical tension field of $f$ in $N^\mv$ and $\Pi^\mv(f)$ its vertical second fondamental form in $N^\mv$.
\end{prop}
\paragraph{When the subgroup $K^\mv=K\cap H^\mv\subset K$.}
Now, let us suppose more generally, that the structure group $H$ of $\pi_M\colon Q\to M$ is reducible to a (closed) subgroup $H^\mrm{v}$, i.e. there exists a principal $H^\mrm{v}$-subbundle $\pi_M^{\mrm v}\colon Q^{\mrm v} \to M$. Then this gives  rise to a fibration $\pi^\mv\colon N^\mv\to M$, with homogeneous fibre $H^\mv/K^\mv$, where $N^\mv=Q^\mv/K^\mv$ and $K^\mv=K\cap H^\mv$.\smallskip\\
We  suppose that the $\Ad K^\mv$-invariant subspace $\pk^\mv:=\hk^\mv\cap \pk$ satisfies $\hk^\mv=\kk^\mv\oplus \pk^\mv$. As previously,  we endow $\pk^\mv$ with the inner product induced by restriction to $\pk^{\mrm v}$ of the $\Ad K$-invariant inner product on $\pk$, i.e. $H^{\mrm v}/K^\mv$  with the metric induced by the  metric on $H/K$, so that the inclusion $H^{\mrm v}/K\to H/K$ is an isometric embedding. Moreover, we suppose that $\pk'=(\pk^{\mrm v})^\perp\subset\pk$  is $\Ad H^{\mrm v}$-invariant, so that in particular $\hk=\hk^{\mrm v}\oplus \pk'$ is a reductive decomposition and $H/H^{\mrm v}$ is reductive. \smallskip\\
Then we can check easily that remark~\ref{rmk-natural-red-autom}, proposition~\ref{prop-tool-corresp}, lemma~\ref{lemma-ivdeB}, proposition~\ref{prop-tau&Pi-corresp} still hold. 
%Moreover, if $\omega$ is reducible to $Q^\mv$, then $\nabla^\mv$ leaves invariant $\ver^\mv\simeq \pk_{Q^\mv}^\mv$.
%
%
\paragraph{About the different possible choices of metrics.} Sometimes, it could happen that the situation imposes a metric on $H^\mv/K^\mv$ different from the one induced by the metric of $H/K$. This leads us to see what is dependent or not of the choice of this metric.\\
Remark that in a homogeneous fibre bundle,  $\phi$, $\Phi$, $T^c$ and $\nabla^c$ do not depend on the $H$-invariant metric on the fibre $H/K$.
Moreover, according to theorem~\ref{difference-tensor}, when this previous metric is naturally reductive, then the difference tensor is independent of the choice of this naturally reductive metric. Therefore, proposition~\ref{prop-tau&Pi-corresp} still hold when $H^\mv/K^\mv$ and $H/K$ are endowed with any naturally reductive metrics, even if the former is not induced by the latter. These considerations leads us to the following.
\begin{defn}
In the situation of definition~\ref{def-red-Hfb}, we will say that  $i_\mv\colon (N^\mv, h^\mv)\to (N,h)$ is a \emph{metric} reduction of homogeneous fibre bundle if is a reduction of homogeneous fibre bundle and a isometry. The term \emph{metric} will be implicitely assumed when we make precise the metrics "$\, i_\mv\colon (N^\mv, h^\mv)\to (N,h)$ is a  reduction of homogeneous fibre bundle". Otherwise, when we only say that  "$\, i_\mv\colon N^\mv\to N$ is  reduction of homogeneous fibre bundle" we will consider that  $\langle\cdot,\cdot\rangle_{\pk^\mv}$ is not induced a priori by $\langle\cdot,\cdot\rangle_{\pk}$ (i.e. $h^\mv$ is not a priori induced by $h$).
\end{defn}

\subsection{Examples of Homogeneous fibre bundles}\label{examhomfibrbund}
\index{canonical!connection, $G$-invariant|(}
In this section, we give examples and applications for the theory developped in the previous sections whose we use here the same notations.
\subsubsection{Homogeneous spaces fibration}\label{homspacefibr}
Let us take $Q=G$ a Lie group, and $K\subset H\subset G$ subgroups of $G$, $(H,K)$ satisfying the hypothesis in the begining of section~\ref{geometry}. Let us suppose that $M=G/H$ is reductive and Riemannian: that is to say if $\g=\hk\oplus\mk$ is the reductive decomposition, then $\Adm H$ is compact and we choose an $\Ad H$-invariant inner product $\langle\ ,\ \rangle_\mk$ in $\mk$. For $\omega$, we take the canonical connection on $\pi_M\colon G\to G/H$ which is given, let us recall it, by $\omega=\theta_\hk$ where $\theta$ is the Maurer-Cartan form in $G$ (see section~\ref{invariant}). Then the corresponding decomposition $TQ=\ver_0\oplus\hor_0$ is given by
$$
T_gG = g.\g= \underbrace{g.\hk}_{\ver_0 } \oplus   \underbrace{g.\mk}_{\hor_0}.
$$
Since  $\mak n:=\pk \oplus\mk$ is $\Ad K$- invariant, then $\g=\kk\oplus \mak n$ is a reductive decomposition and $ N=G/K$ is reductive. Let us recall that we have the canonical identification $G\times_K\g \cong N\times \g$ given by (\ref{io}), which gives us an identification $\mak n_G=G\times_K \mak n\cong [\mak n]$. Then under this    last identification and under the one given by the Maurer-Cartan form of $G/K$, $\beta\colon TN\overset{\cong}{\longrightarrow} [\mak n]$ (see section~\ref{reductivehomspaces}), the splitting $TN=\ver\oplus \hor$ is 
$$
TN=[\pk]\oplus [\mk],
$$
the canonical isomorphism $I\colon \ver \to \pk_G$ is then the identity, and $\phi\colon TN\to \pk_G$ the projection  on $[\pk]$ along $[\mk]$. The metric $h$ on $G/K$ is then   defined by the $\Ad K$-invariant inner product:
\begin{equation}\label{eq-met-p+m}
\langle \ ,\ \rangle_\mak{n}=\langle\ ,\ \rangle_\pk + \langle\ ,\ \rangle_\mk.
\end{equation}
Furthermore, $\Omega= d\omega + \dfrac{1}{2}[\omega\wedge\omega]= d\theta_\hk  + \dfrac{1}{2} [\theta_\hk\wedge\theta_\hk]$ and thus 
$$
\Omega_\pk= d\theta_\pk+ [\theta_\kk\wedge\theta_\pk] +\dfrac{1}{2} [\theta_\pk\wedge\theta_\pk]_\pk.
$$
Since $d\theta +\dfrac{1}{2}[\theta \wedge\theta ]=0$, then (projecting on $\hk$) we have
$$
d\theta_\hk + \dfrac{1}{2}[\theta_\hk \wedge\theta_\hk ] + \dfrac{1}{2}[\theta_\mk\wedge\theta_\mk]_\hk=0
$$
thus 
$$
\Omega=-\dfrac{1}{2}[\theta_\mk\wedge\theta_\mk]_\hk
$$
so that
\begin{equation}\label{lift-Omega}
\Omega_\pk=-\dfrac{1}{2}[\theta_\mk\wedge\theta_\mk]_\pk
\end{equation}
therefore
\begin{equation}\label{formulePhi}
\Phi=-\dfrac{1}{2}[\psi\wedge\psi]_\pk
\end{equation}
where $\psi\colon TN\to  \hor=[\mk]$ is the projection on $\hor$ along $\ver=[\pk]$.\\
The covariant derivative $\nabla^c$, which lifts into $ d +\theta_\kk$ in $G$, is nothing but the canonical linear connection $\nabla^0$ in $N=G/K$ restricted to $[\pk]\subset TN$ (see section~\ref{assocovarder} and \ref{g-invariantaffin}).\\
The canonical torsion $T^c$, which lifts in $G$ into
\begin{equation}\label{lift-T}
d\theta_\pk+ [\theta_\kk\wedge\theta_\pk]= -\dfrac{1}{2} [\theta_\mk\wedge\theta_\mk]_\pk -\dfrac{1}{2} [\theta_\pk\wedge\theta_\pk]_\pk
\end{equation}
is given by \footnote{In particular, according to (\ref{formulePhi}\label{HSE}), we recover, for this example, the Homogeneous structure equation (\ref{hme}).} 
\begin{eqnarray*}
T^c & = &  -\dfrac{1}{2}[\psi\wedge\psi]_\pk -\dfrac{1}{2}[\phi\wedge\phi]_\pk \\
    &  = &  [\ ,\ ]_\mak{n} + [\phi\wedge\psi]_\pk.
\end{eqnarray*}
The associated bundle
$$
\hk_G=G\times_H\hk\cong [\hk]^M:=\{(g.x_0,\Ad g(a)), g\in G,a\in\hk\}\subset M\times \g
$$
can be embedded into $\so(TM)$ by\footnote{Or in other words $[g,a]\in\hk_G=G\times_H\hk\mapsto [g,\adm a]\in G\times_H\so(\mk)\cong\so(TM)$.}
\begin{equation}\label{embedd-h}
\xi=\Ad g(a)\in \Ad g(\hk)= [\hk]_{g.x_0}^M\longmapsto \ad \xi_{|\Ad g(\mk)}=\Ad g\circ\adm a \circ\Ad g^{-1} \in\so(\Ad g(\mk)).
\end{equation}
In the same way, $\kk_G=G\times_K\kk\cong [\kk]^N$ embedds in $\so(N)$. Moreover let us remark that we have 
$$
\pi^*\hk_G=G\times_K \hk=[\hk]^N=\{(g.n_0,\Ad g(a)), g\in G, a\in\hk\},
$$
and that $\pi^*\hk_G$ embedds into $\so(\pi^*TM)$.
\begin{rmk}\label{rmk-metric-inducedbyTr}\em
The subalgebra $\hk$ can be endowed with the $\Ad H$-invariant inner product induced by the $\Ad SO(\mk)$-invariant inner product in $\so(\mk)$: $\langle A, B\rangle= \mrm{Tr}(A^tB)$. Therefore, this inner product induced an $\Ad K$-invariant inner product on $\pk$, which could be a possible choice for $\langle \cdot,\cdot\rangle_\pk$ in \eqref{eq-met-p+m}.
\end{rmk}
As concerns the covariant derivative $\nabla^\omega$, defined in $\hk_G$, it lifts into $d + \theta_\hk$ in $G$. Moreover, under the embedding~(\ref{embedd-h}), it is nothing but the restriction to the subbundle $\so(TM)$ of  the endomorphism connection on $M$ (i.e. the tensor product connection in $T^*M\otimes T^*M$) defined by the canonical linear connection in $M$, $\overset{M}{\nabla^0}$. Indeed under the embedding~(\ref{embedd-h}), $\nabla^\omega$ lifts to the derivative $d + \ad_\mk\theta_\hk$ and equation~(\ref{nabla0-endo}) allows to conclude.\\
Therefore $\nabla^\pk$ is given by its lift \footnote{ $[\ ]_\pk$ denotes as usual the $\pk$-component.} 
\begin{equation}\label{ex-nabla-p}
[(d + \theta_\hk)_{|\pk}]_\pk= d + \ad_\pk\theta_\kk + [\theta_\pk,\cdot_{|\pk}]_\pk 
\end{equation}
that is to say $\nabla^\pk$ is the $[\pk]$-component of the linear connection $\nabla^1$ in $G/K$ (see section~\ref{family}) restricted to $[\pk]\subset TN$. (Indeed we have $[\theta_{\mak n}, \cdot_{|\pk}]_\pk= [\theta_\pk,\cdot_{|\pk}]_\pk + [\theta_\mk, \cdot_{|\pk}]_\pk$, but $[\mk, \pk]\subset[\mk,\hk]\subset \mk$ by reductivity.)\\[1mm]
Let us see to which corresponds $\nabla^c$ under the embedding~(\ref{embedd-h}). The canonical connection $\overset{N}{\nabla}{}^0$ gives rise by restriction to a linear connection on the horizontal subbundle $\hor=\pi^* TM$. Then $\nabla^c$ is nothing but the restriction to $\pk_G\subset \so(\pi^* TM)$ of the tensor product connection in $\End(\pi^* TM)$ defined by $\overset{N}{\nabla}{}^0$.\\[1mm]
Moreover the Levi-Civita connection in $N$ is given by (see section~\ref{family})
$$
\overset{N}{\nabla}=\overset{\mrm{met}}{\nabla^{\frac{1}{2}}} = \nabla^0 + \dfrac{1}{2}\mB^N
$$
where $\mB^N=[\ ,\ ]_{[\mak n]} + \mU^N$ and $\mU^N$ is defined by equation (\ref{defofU}). Then we have by taking the projection on the vertical subbundle $[\pk]$:
$$
\phi(\nabla_A^N V)=\nabla_A^0(\phi V)  + \dfrac{1}{2}\phi\circ \mB^N(A,V)
$$
so that we can conclude according to theorem~\ref{difference-tensor} that
\begin{equation}\label{eq-phi-B=B-phi}
\phi\circ \mB^N= \phi^*\mB-\Phi,
\end{equation}
which can also  be verified directly using the expressions of $\mB^N$, $\mB$ and $\Phi$.
\paragraph{If $H/K$ is (locally) symmetric.} In this case, we have $T^c=\Phi$ (see (\ref{hme}), or (\ref{lift-T}) and (\ref{lift-Omega})). Moreover, according to (\ref{ex-nabla-p}), $\nabla^\pk$ lifts to $d + \theta_\kk$, so that we recover that $\nabla^\pk=\nabla^c$ in this case. Now, let us apply the equality $\nabla^v=\nabla^c$ in $\ver$ (theorem~\ref{difference-tensor})\footnote{We can also use directely theorem~\ref{coderivative}.}.\\
Let $f\colon  (L,b)\to N$ be a map then we have
$$
\tau^v(f)=\mrm{Tr}_b(\nabla^v d^vf)=* d^{\nabla^v}* d^v u = *\Ad F \left( d(*\alpha_\pk) + [\alpha_\kk\wedge (*\alpha_\pk)]\right).f
$$
where $F$ lifts $f$ in $G$ and $\alpha=F^{-1}dF$. Then $f$ is vertically harmonic \iif 
$$
 d(*\alpha_\pk) + [\alpha_\kk\wedge (*\alpha_\pk)]=0
$$
Moreover $f$ is flat ($f^*\Phi=0$) \iif it is vertically torsion free ($f^*T^c=0$) \iif 
$$
[\alpha_\mk\wedge\alpha_\mk]_\pk=0 \Longleftrightarrow  d\alpha_\pk + [\alpha_\kk\wedge\alpha_\pk]=0.
$$
\paragraph{$G/K$ is a (locally) $2k$-symmetric space} 
Let us suppose that there exists an order $2k$ automorphism $\tau\colon\g\to \g$ such that $K=G_0$ with $G_0$ such that $(G^\tau)^0\subset G_0\subset G^\tau$, and $(G^\sigma)^0\subset H\subset G^\sigma$ with $\sigma=\tau^2$ (see section~\ref{finitorderauto}). Then $H/K$ is (locally) symmetric (see section~\ref{finitorderauto}). We have the following identities (with the notation of section~\ref{finitorderauto})
$$
\mk=\oplus_{j=1}^{k-1} \mk_j\quad \text{and} \quad \kk=\g_0,\pk=\g_k.
$$
Then we have 
$$
\Omega_\pk=-\dfrac{1}{2}[\theta_\mk\wedge \theta_\mk]_\pk=-\dfrac{1}{2} 
\sum_{ \underset{i,j\in \Z_{2k}\setminus\{0,k\}}{i + j=k} } [\theta_j\wedge\theta_i],
$$
so that in particular
\begin{prop}\em
Let $(L,j)$ be a Riemann surface.
 If $f\colon (L,j)\to N$ satisfies the equations $\alpha_{-j}''=0$, $1\leq j\leq k-1$ then we have  $f^*\Phi=0$. In other words if $f\colon (L,j)\to N$ is horizontally holomorphic then it  is flat, that is to say $f$ is vertically torsion free or equivalently
\begin{equation}\label{torsionfree}
 d\alpha_k + [\alpha_0\wedge\alpha_k]=0.
\end{equation}
\end{prop}
\begin{thm}\label{evendeterhom}
In the even minimal determined elliptic integrable system $(\syst(k,\tau))$, the last equation $(S_k)$ is equivalent to 
$$
\left\{\begin{array}{ccccl} (\mrm{Re}(S_k)) &  \equiv & d\alpha_k + [\alpha_0\wedge\alpha_k]=0 & \Longleftrightarrow   & f \text{ is vertically torsion free i.e. } f \text{ is flat,}\\
(\mrm{Im}(S_k)) & \equiv &  d(*\alpha_k) + [\alpha_0\wedge(*\alpha_k)]=0 & \Longleftrightarrow   & f \text{ is vertically harmonic.}
\end{array}\right.
$$
In conclusion the even minimal determined elliptic system $(\syst(k,\tau))$ means that  the geometric map $f$ is horizontally holomorphic (which implies that $f$ is flat) and vertically harmonic.
\end{thm}
\begin{rmk}\em
The vertical torsion free equation (\ref{torsionfree}) is the projection on $\pk$ of the Maurer-Cartan  equation provided that we assume the horizontal holomorphicity $\alpha_{-j}''=0, 1\leq j \leq k-1$. In the same way, the equations $(S_j)$, $0\leq j\leq k-1$, of the elliptic system $(\syst(k,\tau))$, are the projections on the different spaces $\g_{-j}$, of the Maurer-Cartan equation, provided that we assume the horizontal holomorphicity.
\end{rmk}
\paragraph{Use of the canonical $2k$-structure $\mijo$.}\index{canonical!embedding|(}
Furthermore the morphism of bundle (over $M$) $\mijo\colon N\to \mZ_{2k,2}^{\alpha_0}(M,J_2)\subset \zdk^{\alpha_0}(M)$ defines a $2k$-structure on $\pi^*TM$ (still denoted by $\mijo$), which according to (\ref{go-gk}) allows to specify the subbundles $\kk_G$ and $\pk_G$ (under the embedding $\pi^*\hk_G\hookrightarrow\so(\pi^*TM)$)
\begin{eqnarray}\label{kg-pg}
\kk_G & = & \{ A\in \pi^*\hk_G| [A,\mijo]=0\}:=\so_{(+1)}(\pi^*TM,\mijo)\cap \pi^*\hk_G\\
\pk_G & = & \{ A\in \pi^*\hk_G| A\mijo +  \mijo  A\}:=\so_{(-1)}(\pi^*TM,\mijo)\cap \pi^*\hk_G
\end{eqnarray}
\begin{rmk}\em
The embedding $\hk_G\hookrightarrow\so(TM)$ is the $H$-equivariant extention of the map $a \in \hk \mapsto \adm a\in \so(\mk)\cong\so(T_{p_0}M)$, and in the same way $\mijo$ is the $H$-equivariant extention of the map $h.G_0\in H/G_0 \mapsto hJ_0h^{-1}\in\mZ(T_{p_0}M,J_0)$, so that the equations (\ref{kg-pg}) are obtained by $H$-equivariance from (\ref{go-gk}).
\end{rmk}
Let us now express the homogeneous fibre bundle tools $\phi$, $\Phi$, and $\nabla^c$ in terms of the embedding $\mijo$. To do not weigh the notation we will forget the index $J_0$ in $\mijo$, in the following theorem.
\begin{thm}\label{entermdeij0}
If $A,B\in TN$, $F\in\mal C(\pk_G)$ then 
\begin{description}
\item[(i)]  $\phi A=  - \dfrac{1}{2}\mai^{-1} \onabla{M}{0}{A} \mai$
\item[(ii)] $\Phi(A,B)= \dfrac{1}{2}\mai^{-1}[\mai,{\pi}^*R^{\onabla{M}{0}{}}(A,B)]$ where $R^{\onabla{M}{0}{}}$ is the curvature of $\onabla{M}{0}{}$.
\item[(iii)] $\nabla_A^{c} F=\onabla{N}{0}{A} F=\dfrac{1}{2}\mai^{-1}[\mai,\onabla{M}{0}{A} F]$. 
\end{description}
\end{thm}
\proof
Lifted  in $G$ (i) become
$$
\adm\theta_k= - \dfrac{1}{2}J_0^{-1}\left( d J_0 + \left[ \adm \theta_\hk, J_0 \right] \right), 
$$
 keeping in mind that under the embedding $\pi^*\hk_G\hookrightarrow\so(\pi^*TM)$, $\theta_k $ corresponds to $\adm\theta_k$. Moreover this equality holds since $\left[ \adm \theta_\hk, J_0 \right]= -2 J_0\adm \theta_k$. This proves (i).\\
For the following, let us keep in mind that $\hk_G$ is embedded into $\so(\pi^*TM)$.\\
Then lifted in $G$, (ii) becomes
$$
-\adm \left([\theta_\mk,\theta_\mk]_\pk \right) = \dfrac{1}{2}J_0^{-1}\left[ J_0, -\adm \left([\theta_\mk,\theta_\mk]_\hk\right)  \right] 
$$ 
which holds for the same reason as above.\\
Finally (iii) lifted in $G$ becomes
$$
 d + \theta_0 = -\dfrac{1}{2}J_0^{-1}\left[ J_0, d \tl F + \left[ \adm \theta_\hk, \tl F\right]\,   \right] 
$$ 
which holds since the right hand side is equal to $d \tl F + \left[ \adm \theta_0, \tl F\right]$ which lifts $ \onabla{N}{0}{}F$, because $\left[ \adm \theta_k, \tl F\right]$ commutes with $J_0$. We have denoted by $\tl F\colon G\to \g_k$ the lift of $F$. \comprf \hsq
\begin{thm}\label{verticalij0}
Let $s\in \mal C(\pi)$ and  $J=s^*\mijo\in\mal C(\mZ_{2k,2}^{\alpha_0}(M,J_2))$ be the corresponding $2k$-structure. Then
\begin{description}
\item[(i)] $I(d^{v} s)= -\dfrac{1}{2} J^{-1}\onabla{M}{0}{} J$. Thus $s$ is horizontal \iif $J$ is $\onabla{M}{0}{}$-parallel.
\item[(ii)] $I(\Pi^{v}(s)) = -\dfrac{1}{4}[J^{-1},(\onabla{M}{0}{})^2J ]$.\\[1mm]
Thus $s$ is superflat \iif $(\onabla{M}{0}{})^2J $ commutes with $J$.
\item[(iii)] $I(\tau^{v}(s))= \dfrac{1}{4}[J^{-1}, (\onabla{M}{0}{})^* \onabla{M}{0}{} J ]  $.\\[1mm]
Thus $s$ is a harmonic section \iif $(\onabla{M}{0}{})^* \onabla{M}{0}{} J $ commutes with $J$. 
\item[(iv)] $s^*\Phi=\dfrac{1}{2} J^{-1}[J,R^{\onabla{M}{0}{}}]$. 
\end{description}
These properties hold also for maps $f\in\mal C^\infty(L,N)$, $(L,b)$ being a Riemannian manifold: (i),(ii),(iii) without any change and (iv) becoming $f^*\Phi=\dfrac{1}{2}J^{-1}[J,u^*R^{\onabla{M}{0}{}}]$, with $u=\pi\circ f$.
\end{thm}
\proof
(i) and (iv) follows from theorem~\ref{entermdeij0} (i) and (ii) respectively. The assertions (iii) is a consequence of the assertion (ii). Finally (ii) follows from theorem~\ref{verticalij0}\,-(iii) and the fact that 
$$
\left[ \onabla{M}{0}{}J^{-1}, \onabla{M}{0}{}J  \right] = - \left[ J^{-1} \left( \onabla{M}{0}{} J\right) J^{-1},\onabla{M}{0}{}J\right] 
$$
commutes with $J$. \comprf \hsq \medskip\\
Now, we can conclude.
\begin{cory}\label{interpretnabla0}
Let $(L,j)$ be a Riemann surface, $f\colon L\to N$ a map and $J=f^*\mijo$ the corresponding map into $\mZ_{2k,2}^{\alpha_0}(M,J_2)$. Then $f$ is a geometric solution of the even minimal determined system $(\syst(k,\tau))$ \iif 
\begin{description}
\item[(i)] $J$ is an admissible twistor lift ($\Leftrightarrow$ f is horizontally holomorphic).
\item[(ii)] $J$ is vertically harmonic\footnote{The vertical harmonicity is \wrt the splitting defined by example~\ref{zdkjE}. See section~\ref{5.3.4} and (more precisely) theorem~\ref{verticalj=2}, for the above characterisation of vertical harmonicity in $\mZ_{2k,2}^{\alpha_0}(M,J_2)$.}:  $\left[(\onabla{M}{0}{})^* \onabla{M}{0}{} J,J\right]=0$ ($\Leftrightarrow$ f is vertically harmonic).
\end{description}
Moreover the first condition implies that $\left[ u^*R^{\onabla{M}{0}{}},J\right] =0$  i.e. that $J$ is a flat section in $(\End(u^*TM),u^*\onabla{M}{0}{})$ ($\Leftrightarrow$ f is  flat).\\
Furthermore $f$ is a primitive geometric solution (i.e. there exists $m\leq k$ such that $f$ is $m$-primitive, which is equivalent to say that $f$ is $k$-primitive) \iif
\begin{description}
\item[(i)] $J$ is an admissible twistor lift
\item[(ii)] $J$ is parallel: $\onabla{M}{0}{} J=0$  ($\Leftrightarrow$ f is horizontal).
\end{description}
\end{cory}
\index{canonical!connection, $G$-invariant|)}\index{canonical!embedding|)}

%%%%%%%%%%%%%%%%%%%%%%%%%%%%%%%%%%%%%%%%%%%%%%%%%%%%%%%%%%%%%%%%%%%%%%%%%%%%%%%%%%%%%%
%
\subsubsection{The twistor bundle of almost complex structures $\Sigma(E)$}\label{twistor2}
\index{twistor|(}
We give ourself the same ingredients as in example~\ref{sigmaE}. Let us suppose  that the vector bundle $E$ is oriented. Then the bundle of positive (resp. negative) orthogonal almost complex structure on $E$ (i.e. the component $\Sigma^\eps(E)$ of $\Sigma(E)$ with $\eps=\pm 1$), $\pi_\Sigma \colon \Sigma^\eps(E)\to M$ is a homogeneous fibre bundle. Indeed, we take $Q=\mal{SO}(E)$ the $SO(2n)$-bundle of positively oriented orthonormal frames of $E$, $H=SO(2n)$ and $K=U(n)$ (embedded in $SO(2n)$ via $A + i B\mapsto \begin{pmatrix} A & -B \\ B & A \end{pmatrix}$). $K $ is the subgroup of $SO(2n)$ which commutes with $J_0^\eps=\eps\begin{pmatrix} 0 & -\Id \\ \Id & 0\end{pmatrix}$. The involution $T=\Int J_0^\eps$ in $SO(2n)$ gives rise to the symmetric  space $H/K=\Sigma^\eps(\rdn)$, and to the following symmetric decomposition $\hk=\kk\oplus\pk$ with
\begin{eqnarray*}
\kk & = & \{A\in \so(2n)| [A,J_0^\eps]=0\}\\
\pk & = & \{A\in \so(2n)| AJ_0^\eps + J_0^\eps A=0\}.
\end{eqnarray*}
Concerning $\omega$, we take the $\so(2n)$-valued connection 1-form on $Q$ corresponding to the covariant derivative $\nabla$ in $E$: if $e=(e_1,\ldots,e_{2n})$ is a (local) moving frame of $E$ (i.e. a section of $Q$) then 
$$
\nabla(e_1,\ldots,e_{2n})=(e_1,\ldots,e_{2n})\omega(e;de).
$$
Now, let us consider the isomorphism of bundle:
$$
\mJ\colon e.U(n)\in \mal{SO}(E)/U(n)\overset{\cong}{\longmapsto} J \in \Sigma^\eps(E)|\, \mal{M}at_{e.U(n)}(J)=J_0^\eps.
$$
The isomorphism $\mJ$ defines a bijection between the set of section of $\pi\colon N\to M$\footnote{In all the section~\ref{examhomfibrbund}, as it was the case in all the section~\ref{Homfibrbund}, $N:=Q/K$.} and the set of complex structure of $E$ (sections of $\pi_\Sigma$): $s\in \mal C(N)\to J=\mJ\circ s\in \mal C(\Sigma^+(E))$.\\
The existence of a (positive) complex structure $J$ in $E$ -- i.e. a section of $\pi_\Sigma\colon\Sigma^+(E)\to M$--  is equivalent to the existence of an $U(n)$-reduction of the principal bundle $\mal{SO}(E)\to M$: $J$ defines a Hermitian structure on $E$ and then the $U(n)$-subbundle of unitary frames for this Hermitian structure, and vice versa.\\
The isomorphism of bundle over $M$, $\mJ\colon N\to \Sigma^+(E)$ defines tautologically a canonical complex structure on $\pi^*E\to N$ (which we still denote by $\mJ$)\footnote{and which is in fact nothing but $\mJ^* \mal I=\mal I\circ \mJ$, see example~\ref{sigmaE}.} $\mJ\colon N\to \Sigma^+(\pi^*E)$. Under this identification, let us specify the subbundles $\pk_Q$ and $\kk_Q$. First, we have $\hk_Q=\so(E)$\footnote{See remark~(\ref{mathcalJ}) (and more precisely equation~(\ref{hqsoe})) for the identification map.}, the bundle of skew-symmetric endomorphism of $E$ and then\footnote{Since $\pi^*E$ is canonically endowed with the complex structure $\mJ$, we need not to specify this latter in the notation $\so_{\pm}(\pi^*E)$, whereas $E_x$, for $x\in M$, could be endowed with any element $J_x\in \Sigma^+(E_x)$, this is why we must specify it in the notation $\so_{\pm}(E_x,J_x)$.}
\begin{eqnarray*}
(\kk_Q)_y  & = & \{ F\in \so(E_{\pi(y)})| [F,\mJ (y)]=0\}=:\so_+(E_{\pi(y)},\mJ(y))=:\so_+(\pi^*E)_y \\
(\pk_Q)_y & = & \{ F\in \so(E_{\pi(y)})| F\mJ (y) + \mJ(y) F=0\}=:\so_-(E_{\pi(y)},\mJ(y))=:\so_-(\pi^*E)_y.
\end{eqnarray*}
Then the decomposition following $\pi^*\hk_Q=\kk_Q\oplus\pk_Q$ of any element $F\in \pi^*\hk_Q=\so(\pi^*E)$ is given by 
$$
F = \dfrac{1}{2}\mJ\{F,\mJ\} + \dfrac{1}{2}\mJ[F,\mJ]
$$
where $\{\ ,\ \}$ is the anticommutator.
\begin{rmk}\label{mathcalJ}\em
The canonical complex structure $\mJ$ is a section of the associated bundle over $N$: $\Sigma^+(\pi^*E)=\pi^*(Q\times_H \Sigma^\eps (\rdn))=\pi^*Q\times_H \Sigma^\eps(\rdn)$, so that it can be lifted to a $H$-equivariant map $\tl\mJ\colon\pi^*Q\to \Sigma^\eps(\rdn)\subset\hk$, which is given by 
$$
 \tl \mJ\colon (e.K,e.h^{-1})\in \pi^*Q\longmapsto hJ_0^\eps h^{-1}\in\Sigma^\eps(\rdn).
$$
Remark that the restriction of $\tl \mJ$ to $Q\subset\pi^* Q$ is the constant map $J_0^\eps$ (the inclusion $Q\subset\pi^* Q$ is given by $e\mapsto (e.U(n),e)$), and that $\tl\mJ\colon\pi^*Q\to \Sigma^\eps(\rdn)$ is the $H$-equivariant extension of the $K$-equivariant constant map $J_0^\eps$ on $Q$. $\tl \mJ$ can also be given in terms of $\mJ$ by 
$$
\tl\mJ\colon(y;e)\in\pi^*Q\longmapsto \mal{M}at_e(\mJ(y))\in \Sigma^\eps(\rdn)\subset\hk.
$$
Furthermore, we have a canonical identification $N=Q\times_H H/K$ (via $[e,h.K]\mapsto (e.h).K$) and the identification depending on $J_0^\eps$: $H/K=\Sigma^\eps(\rdn)$ (via $h.K\mapsto hJ_0h^{-1}$) so that $N=Q\times_H\Sigma^\eps(\rdn)$ (via $e.K\mapsto [e,J_0]$).  Then under this last identification, 
$\mJ$ is the restriction to $N$ of the canonical identification
\begin{equation}\label{hqsoe}
\begin{array}{rcl}
\hk_Q:=Q\times_H\hk & \overset{\cong}{\longrightarrow} & \so(E)\\
{ [e,a]} & \longmapsto &  A |\ \mal{M}at_e(A)=a.
\end{array}
\end{equation}
Therefore
$$
\hk_Q=\mJ^*\so_+(\pi_{\Sigma}^*E) \quad   \text{ and } \quad \pk_Q=\mJ^*\so_-(\pi_{\Sigma}^*E),
$$
with the notations of example~\ref{sigmaE}.
\end{rmk}
Let us now express the homogeneous connection $\phi$, the curvature foms $\Phi$ and the canonical connection $\nabla^c$ in terms of $\mJ$ (following~\cite{cmw2}). 
\begin{thm}\label{entermdeJ}\cite[Prop.~4.1]{cmw2}
If  $A,B\in TN$, $F\in \mal C(\pk_Q)$ then:
\begin{description}
\item[(i)] $\phi A=\dfrac{1}{2}\mJ .\nabla_A\mJ$
\item[(ii)] $\Phi(A,B)=\dfrac{1}{2}\mJ[\pi^*R(A,B),\mJ]$, where $R$ is the curvature operator of the $\nabla$.
\item[(iii)] $\nabla_A^c F=\dfrac{1}{2}\mJ[\nabla_A F,\mJ]$
\end{description}
\end{thm}
\begin{thm}\label{vertical}\cite[Theorem~4.2]{cmw2}
Let $s\in \mal C(\pi)$ and  $J=s^*\mJ$ be the corresponding complex structure, and $\nabla^*\nabla = -\mrm{Tr}\nabla^2$, the rough Laplacian of $E$. Then
\begin{description}
\item[(i)] $I(d^v s)= \dfrac{1}{2}J.\nabla J=\dfrac{1}{4}[J,\nabla J]$. Thus $s$ is horizontal \iif $J$ is parallel.
\item[(ii)] $I(\Pi^v(s))=\dfrac{1}{4}[J,\nabla^2J]$. Thus $s$ is superflat \iif $\nabla^2J$ commutes with $J$.
\item[(iii)] $I(\tau^v(s))= -\dfrac{1}{4}[J,\nabla^*\nabla J]$. Thus $s$ is a harmonic section \iif $\nabla^*\nabla J$ commutes with $J$.
\item[(iv)] $s^*\Phi=\dfrac{1}{2}J[R,J]$.
\end{description}
These properties hold also for maps $f\in\mal C^\infty(L,N)$, $(L,b)$ being a Riemannian manifold: (i),(ii),(iii) without any change and (iv) becoming $f^*\Phi=\dfrac{1}{2}[u^*R,J]$, with $u=\pi\circ f$.
\end{thm}
From theorem~\ref{entermdeJ}-(i) (or theorem~\ref{vertical}-(i)) it follows that $d\mJ$ sends the decomposition $TN=\ver\oplus\hor$ onto the decomposition $T\Sigma^\eps(E)=\ver^\Sigma\oplus\hor^\Sigma$ coming from $\nabla$ (see example~\ref{sigmaE})
 so that we can consider $\pi_\Sigma\colon\Sigma^\eps(E)\to M$ as a homogeneous fibre bundle over $M$ with  structure group $H=SO(2n)$ and $K=U(n)$. Besides, since  the vertical and horizontal subbundles corresponds via $\mJ$, then we can conclude according to (\ref{metric}) and (\ref{metrich}) that \textbf {$\mal J$ is an isometry}. \\[1mm]
Moreover, we see that $s$ is vertically harmonic in $N$ \iif the rough Laplacian $\nabla^*\nabla J$ of $J$ in $\so(E)$ is vertical (i.e. in $\ver_J^\Sigma$, see example~\ref{sigmaE}) so that we recover the definition of  vertically harmonic twistor lifts used in \cite{ki3} and \cite{bk}. More precisely, via the isometry $\mJ$, \textbf{the vertical tension field of $s$} -- which is, let us recall it, defined using the Levi-Civita connection in $N$ which corresponds via the isometry $\mJ$ to the Levi-Civita connection in $\Sigma^+(E)$-- \textbf{is exactly the vertical part in $\pi_\Sigma^*\so(E)$ of the rough laplacian of $J$}:
$$
d\mJ(\tau^v(s))=\nabla_{\tau^v(s)}\mJ = -2\mJ \circ\phi(\tau^v(s))=\dfrac{1}{2}J[J,\nabla^*\nabla J]
$$ 
according to theorem~\ref{entermdeJ}-(i) and theorem~\ref{vertical}-(iii). Concretely, \textbf{to compute the vertical tension field in $\Sigma^+(E)$}, instead of using the (abstract) Levi-Civita connection, it is enough to \textbf{take the vertical part of the rough Laplacian} (which uses the concrete metric connection $\nabla$).
\subsubsection{The twistor bundle $\zdk(E)$ of a Riemannian vector bundle}\label{subsect-zdkE}
We give ourself  the same ingredients and notations as in example~\ref{zdkE}. Let us suppose that the vector bundle $E$ is oriented. Then the bundle $\pi_{\mal Z}\colon \zdk^\alpha(E)\to M$ is a homogeneous fibre bundle. Indeed, we take $Q=\mal{SO}(E)$, $H=SO(2n)$ and $K=\U_0(J_0^\alpha)$. Let us recall that the order $r$ automorphism $T=\Int J_0^\alpha$ in $SO(2n)$ gives rise to the $r$-symmetric space $H/K=\zdk^\alpha(\rdn)$, and to the following reductive decompostion $\hk=\kk\oplus \pk$ with 
$$
\kk=\so_0(J_0^\alpha) \quad \text{and} \quad \pk =\so_*(J_0^\alpha):= \left( \bigoplus_{j\in\Z/r\Z\setminus\{0\}} \negthickspace \negthickspace \so_j^\C(J_0^\alpha)\right) \bigcap\so(2n).
$$
Concerning $\omega$, we take the same as in the previous example. Now let us consider the isomorphism of bundle:
\begin{equation}\label{defjk}
\mJ\colon e.\U_0(J_0^\alpha)\in \mal{SO}(E)/\U_0(J_0^\alpha)\overset{\cong}{\longmapsto} J \in \zdk^\alpha(E)|\, \mal{M}at_{e}(J)=J_0^\alpha.
\end{equation}
The isomorphism $\mJ$ defines a bijection between the sections of $\pi\colon N\to M$ and the set of sections of $\pi_\mal{Z}\colon \zdk^\alpha(E)\to M$, $s\in \mal C(N)\mapsto J=\mJ\circ s\in \mal C(\zdk^\alpha(E))$.\\
The isomorphism of bundle  over $M$, $\mJ\colon N\to \zdk^\alpha(E)$ defines tautologically a canonical $2k$-structure\footnote{and which is in fact nothing but $\mJ^* \mal I=\mal I\circ \mJ$, see example~\ref{zdkE}.} on $\pi^*E\to N$ (still denoted by $\mJ$), $\mJ\colon N\to \zdk^\alpha(\pi^*E)$. Under this consideration, we therefore have $\hk=\so(E)$ and for all $y\in N$,
\begin{eqnarray*}
\kk_Q & = & \so_0(\pi^*E,\mJ)\\
\pk_Q & = & \so_*(\pi^*E,\mJ)=\left( \bigoplus_{j\in\Z/r\Z\setminus\{0\}} \negthickspace \negthickspace \so_j^\C(\pi^*E,\mJ)\right) \bigcap\so(\pi^*E).
\end{eqnarray*}
Since $\pi^*E$ is canonically endowed with $\mJ$, we will not specify it and use the notation $\so_j^\C(\pi^*E):=\so_j^\C(\pi^*E,\mJ)$.\\
Let us consider the surjective morphism of vector bundle
$$
\begin{array}{crcl}
\ad\mJ\colon & \pi^*\hk_Q=\so(\pi^*E) & \longrightarrow  & \maB_*(\pi^*E) =\mJ .\so_*(\pi^*E)=\mJ .\pk_Q\\
             & (J,A) & \longmapsto & \ad J (A)=[J, A]= J\sum_{j=1}^r(1-\omega_r^j)A_j
\end{array}
$$
where $A_j=[A]_{\so_j^\C(E_x)}$ is the $\so_j^\C(E_x)$-component of $A\in\so(E_x)$. The kernel of $\ad\mJ$ is $\kk_Q=\so_0(\pi^*E)$ so that $\ad\mJ$ induces an isomorphism from $\pk_Q$ onto $\mJ .\pk_Q$. We will set 
$$
(\ad\mJ)^{-1}=(\ad\mJ_{|\mJ .\pk_Q})^{-1}\oplus 0_{\mJ.\kk_Q}
$$
so that
\begin{eqnarray}
(\ad\mJ)^{-1}\circ \ad \mJ & = & \mrm{pr}_{\pk_Q}, \text{ the  projection on }\pk_Q \text{ along }\kk_Q, \text{  and} \label{proj-pq}\\
\ad \mJ \circ (\ad\mJ)^{-1} & = & \mrm{pr}_{\mJ .\pk_Q}, \text{ the  projection on } \mJ .\pk_Q \text{ along } \mJ .\kk_Q .\label{proj-jpq}
\end{eqnarray}
Let us remark that $\mJ.\pk_Q=\mJ^*\ver^\mal{ Z}$ is the (pullback by $\mJ$ of the) vertical space of $\pi_\mal{Z}$ (see example~\ref{zdkE}). More precisely the $\mJ$-pullback of the decomposition ${\ver^{SO(E)}}_{|N_\mZ}=\mal B_0(\pi_\mZ^* E)\oplus \mal B_*(\pi_\mZ^* E)$ (see example~\ref{zdkE}) is the decomposition $\mJ.\so(E)=\mJ.\kk_Q\oplus\mJ.\pk_Q$.\\[1mm]
Let us now express the homogeneous fibre bundle tools $\phi,\Phi$ and $\nabla^\pk$ in terms of $\mJ$. 
\begin{thm}\label{entermdeJ2}
If $A,B\in TN$, $F\in \mal C(\pk_Q)$ then
\begin{description}
\item[(i)] $\nabla \mJ= -\ad \mJ\circ\phi$ thus $\phi A=  -(\ad \mJ)^{-1}\nabla_A\mJ$
\item[(ii)] $\Phi(A,B)=(\ad \mJ)^{-1}[\mJ,\pi^*R(A,B)]$
\item[(iii)] $\nabla_A^\pk F=(\ad \mJ)^{-1}[\mJ,\nabla_A F]$
\end{description}
\end{thm}
\proof
We follow the proof  given in \cite{cmw2} of theorem~\ref{entermdeJ}.\\
(i) Let $D$ be the exterior covariant derivative for $\mak{gl}(2n)$-valued differential form on $\pi^*Q$. Let $\tl A\in TQ$ lifting $A\in TN$, then
$$
D\tl \mJ (\tl A)= d\mJ(\tl A) + [\omega(\tl A),\tl\mJ]= [\omega_\pk(\tl A), J_0^\alpha],
$$ 
since $\tl \mJ= J_0^\alpha$. Projection on $N$ yields%\footnote{Keep in mind that we have $N=Q/K=(\pi^* Q)/H$ and   $TN= TQ/\kk}= T(\pi^*Q)/\hk$.}
$$
\nabla_A\mJ = - \ad \mJ \circ \phi(A).
$$
(ii) First we have $\Omega_\pk= (\ad J_0^\alpha)^{-1} \left[ J_0^\alpha,\Omega\right]$, according to \eqref{proj-pq}. The right hand side is the restriction to $Q$ of the $\hk$-valued 2-form $(\ad \tl\mJ)^{-1}[\tl\mJ,\tl\Omega]$ on $\pi^*Q$, where $\tl \Omega$ is the pullback of $\Omega$. Since $\tl\Omega$ is the curvature of the pullback connection, on projection to $N$ we obtain $\Phi=(\ad\mJ)^{-1}[\mJ,\pi^*R]$ (by definition of  $\omega$, $\Omega$ is the lift in $Q$ of the curvature operator $R$ of $\nabla$).\\
(iii) $\nabla^\omega$ is the restriction of the tensor product connection in $E^*\otimes E= \End(E)$ to $\hk_Q=\so(E)$, and its $\pk_Q$-component follows from \eqref{proj-pq}. \comprf\hsq\medskip\\
Given $A,B$ 1-forms on $N$ with values in some Lie algebra, we define
$$
[A\odot B](X,Y)=[A(X), B(Y)] + [A(Y), B(X)], \quad \forall X,Y\in TN
$$ 
\begin{thm}\label{vertical2}
Let $s\in \mal C(\pi)$ and  $J=s^*\mJ$ be the corresponding $2k$-structure. Then
\begin{description}
\item[(i)] $I(d^v s)= -(\ad J)^{-1}\nabla J$. Thus $s$ is horizontal \iif $J$ is parallel.
\item[(ii)] $I(\Pi^v(s)) = -(\ad J)^{-1}\nabla^2J + \dfrac{1}{2}(\ad J)^{-1} [\nabla J \odot(\ad J)^{-1}\nabla J]$.\\[1mm]
Thus $s$ is superflat \iif $\nabla^2J - \dfrac{1}{2} [\nabla J \odot(\ad J)^{-1}\nabla J]$ commutes with $J$.
\item[(iii)] $I(\tau^v(s))= +(\ad J)^{-1}\nabla^*\nabla J  + \dfrac{}{}(\ad J)^{-1} \mrm{Tr}\left( [\nabla J ,(\ad J)^{-1}\nabla J]\right) $.\\[1mm]
Thus $s$ is a harmonic section \iif $\nabla^*\nabla J + \dfrac{}{}\mrm{Tr}\left( [\nabla J ,(\ad J)^{-1}\nabla J]\right)$ commutes with $J$. 
\item[(iv)] $s^*\Phi=(\ad J)^{-1}[J,R]$.
\end{description}
These properties hold also for maps $f\in\mal C^\infty(L,N)$, $(L,b)$ being a Riemannian manifold: (i),(ii),(iii) without any change and (iv) becoming $f^*\Phi=(\ad J)^{-1}[J,u^*R]$, with $u=\pi\circ f$.
\end{thm}
\proof  
(i) Take the  pullback by $f$  of theorem~\ref{entermdeJ2}-(i).\\
(iv) Take the  pullback by $f$  of theorem~\ref{entermdeJ2}-(ii). \\
(ii) According to theorems~\ref{tauv-tau} and \ref{phinablads}-(i), we have 
\begin{eqnarray*}
I(\Pi^v s) & = & \nabla^c (f^*\phi) + \dfrac{1}{2} f^*\phi^* \mB \\
 & = &  \nabla^c (f^*\phi) + \dfrac{1}{2}\left[  f^* \phi,f^*\phi\right]_\pk \quad  (\text{since  } H/K \text{ is naturally reductive}),\\ 
& = & \nabla^\pk (f^* \phi) - \dfrac{1}{2} \left[ f^*\phi ,f^*\phi\right ]_\pk  \quad (\text{according to \eqref{nabla-p}}).
\end{eqnarray*}
Moreover,
\begin{eqnarray*}
\nabla^\pk (f^*\phi) =\mrm{pr}_{\pk_Q} \left( \nabla^\omega (f^*\phi)\, \right) & = & \mrm{pr}_{\pk_Q} \left( \nabla \left( -(\ad J)^{-1} \nabla J\right) \,\right) \\
& = &   (\ad J)^{-1} \circ \ad(\nabla J)\circ (\ad J)^{-1} (\nabla J) - (\ad J)^{-1}\left( \nabla^2 J \right) \\
& = &   (\ad J)^{-1} \left( [\nabla J ,(\ad J)^{-1}\nabla J]\right)  - (\ad J)^{-1}\left( \nabla^2 J \right) 
\end{eqnarray*}
where in the second line we have used the lemma~\ref{nabla_ad_J-1} below.\\
Now, by an immediate computation using \eqref{proj-pq}, we obtain that
$$
\left[(\ad J)^{-1} \nabla J ,(\ad J)^{-1}\nabla J\right]_{\pk_Q} = \left [\nabla J \wedge(\ad J)^{-1}\nabla J\right] 
$$
This completes the proof of (ii).\\
Then (iii) follows immediately from (ii). \comprf\hsq
\begin{lemma}\label{nabla_ad_J-1}
Let $E, F$ be two vector bundles each one endowed with  a connection that we denote indisctinctely each by $\nabla$. Let $A\in E^*\otimes F$ be an linear morpism. Let us consider a  spliting $E=K_E \oplus \Bar E$ where $K_E=\ker A$. Then let us set 
$B=\left( A_{\Bar E}\right)^{-1} \oplus 0_{K_F}$, where $K_F$ is some complementary subbundle  to $\Bar F:= A(E)$. \\
Then  we have 
$$
\mrm{pr}_{\Bar F} \circ \nabla B \circ  \mrm{pr}_{\Bar E }= - B\circ \nabla A \circ B 
$$
\end{lemma}
\proof
Differentiating the equations $ A\circ  B=\mrm{pr}_{\bar F}$ and  $ \mrm{pr}_{\bar E} \circ B= B$, we obtain
\begin{eqnarray*}
\nabla A \circ  B +  A\circ \nabla B  & =  & \nabla \mrm{pr}_{\bar F}\\
\nabla  \mrm{pr}_{\bar E} \circ  B +  \mrm{pr}_{\bar E} \circ \nabla B  & =  & \nabla B
\end{eqnarray*}
so that multiplying the first equation   $B$  on the the left, and injecting it in the second equation, we obtain
\begin{eqnarray*}
\nabla B  &  =  &  \nabla  \mrm{pr}_{\bar E} \circ  B + B \circ \nabla \mrm{pr}_{\bar F}  - B\circ \nabla A \circ B\\
 &  =  &   \mrm{pr}_{K_E}\circ \nabla B + B \circ \nabla \mrm{pr}_{\bar F}  - B\circ \nabla A \circ B
\end{eqnarray*}
\comprf\hsq\medskip\\
As in \ref{twistor2} above,  we conclude from theorem~\ref{entermdeJ2}-(i)  that $d\mJ$ sends  the decomposition $TN=\ver\oplus \hor$ onto the decomposition $TN_\mZ=\ver^\mZ\oplus\hor^\mZ$ coming from $\nabla$ (see example~\ref{zdkE})  so that we can consider $\pi_\mZ\colon\zdk^\alpha(E)\to M$ as a homogeneous fibre bundle over $M$ with  structure groups $H=SO(2n)$ and $K=\U_0(J_0^\alpha)$. We will call this structure \emph{the homogeneous fibre bundle structure defined in $N_\mZ$ by  $\nabla$ (or by the Riemannian vector bundle $(E,\nabla)$)}.\\
Besides, since  the vertical and horizontal subbundles corresponds via $\mJ$, then we can conclude according to (\ref{metric}) and (\ref{metrich}) that \textbf {$\mal J$ is an isometry}.\\[1mm]
Moreover, the vertical tension field of $J$ in $N_\mZ=\zdk^{\alpha}$ is given by
\begin{multline*}
d\mJ(\tau^v(s))=\nabla_{\tau^v(s)}\mJ = -(\ad\mJ) \circ\phi(\tau^v(s))\\
  =-(\ad\mJ)\circ (\ad\mJ)^{-1}(\nabla^*\nabla J + \dfrac{}{}\mrm{Tr}\left( [\nabla J ,(\ad J)^{-1}\nabla J]\right))\\ 
 =-\left[ \nabla^*\nabla J + \dfrac{}{}\mrm{Tr}\left( [\nabla J ,(\ad J)^{-1}\nabla J]\right)\right]_{\ver^\mZ} 
\end{multline*} 
By taking $k=2$ in the two preceding theorems, we recover of course the results of the previous section: just remark that in this case, $\ad\mJ= 0_{\kk_Q}\oplus 2{L_{\mJ}}_{|\pk_Q}$, and that  $\nabla J$ anticommutes with $J$.
\begin{rmk}\label{Jj}\em
Let us consider the canonical identification
\begin{equation}\label{pi*identif}
\begin{array}{rcl}
H_Q\colon=Q\times_H H & \overset{\cong}{\longrightarrow} & SO(E)\\
{ [e,h]} & \longmapsto &  A |\ \mal{M}at_e(A)=h.
\end{array}
\end{equation}
then $\mJ$ is the restriction to $N\cong Q\times_H\zdk^\alpha(\rdn)$ (via $e.K\mapsto [e,J_0^\alpha]$) of (\ref{pi*identif}).\\
More generally, for $ j\in \Z$, we can consider $\mJ_j$ the restriction of (\ref{pi*identif}) to $Q/\U_{j-1}(J_0^\alpha)= 
%%(Q\times_{SO(2n)}SO(2n)/ \U_{j-1}(J_0^\alpha)=)
Q\times_{SO(2n)}(\zdk(\rdn))^j$ (via $e.\U_{j-1}(J_0^\alpha)\mapsto [e, (J_0^\alpha)^j]$):
\begin{equation}\label{isomJj}
\mJ_j\colon e.\U_{j-1}(J_0^\alpha)\in \mal{SO}(E)/\U_{j-1}(J_0^\alpha)\overset{\cong}{\longmapsto} J \in (\zdk^\alpha(E))^j|\, \mal{M}at_{e}(J)=(J_0^\alpha)^j.
\end{equation}
\end{rmk}
\begin{rmk}\label{suppl}\em
The previous study could have been done (without any change) for any component $\mal U_{2k}^\alpha(E)$. In particular, by replacing $J_0^\alpha$ by $(J_0^\alpha)^j$ in what precedes, we get the isomorphism~(\ref{isomJj})
$$
\mJ_j\colon \mal{SO}(E)/\U_{j-1}(J_0^\alpha)\overset{\cong}{\longmapsto} (\zdk^\alpha(E))^j=\mal U_{p}^{j\cdot\alpha}(E),
$$
where $p=\dfrac{2k}{(2k,j)}$, and by applying theorem~\ref{vertical2}, we see that a cross section  $s_j\colon M\to \mal{SO}(E)/\U_{j-1}(J_0^\alpha)$ is horizontal \iif the corresponding section $J_j=\mJ_j\circ s_j\colon M\to (\zdk^\alpha(E))^j$ is parallel: $\nabla J_j=0$.
\end{rmk}
\subsubsection{The Twistor subbundle $\mZ_{2k,j}^\alpha(E)$}\label{5.3.4}
We continue here the study of example~\ref{zdkjE}, $\pi_\mZ^j\colon \mZ_{2k,j}^\alpha(E,J_j)\to M$, and prove that it defines a homogeneous bundle fibre bundle. Let us recall that we have a bijection between the set of (global) sections $J_j$ in $(\zdk^\alpha(E))^j=\mal U_p^{j\cdot\alpha}(E)$ and the set of $\U_{j-1}(J_0^\alpha)$-reductions $\pi^j\colon Q^j \to M$ of $\mal{SO}(E)$, which is given by 
\begin{equation}\label{defQj}
Q^j=\Um_{j-1}^\alpha(E):=\{e \in \mal{SO}(E)| \mal Mat_e(J_j)=(J_0^\alpha)^j\}.
\end{equation}
Let us consider such a reduction $Q^j$ (defined by some $J_j$). Then $\pi^j\colon Q^j\to M$ is a principal bundle with structural group $H^j=\U_{j-1}(J_0^\alpha)$ and we take for the second structural group $K=\U_0(J_0^\alpha)$ as in the previous example. Let us recall that the order $j$  automorphism $T={\Int J_0^\alpha}_{|\U_{j-1}(J_0^\alpha)}$ gives rise to the $j$-symmetric space $H^j/K=\mZ_{2k,j}^\alpha(\rdn,(J_0^\alpha)^j)$, and to the following reductive decomposition $\hk^j=\kk\oplus\pk^j$ where 
\begin{eqnarray*}
\kk & = & \so_0(J_0^\alpha) \quad \text{and}\\
 \pk^j & = & \ul_{j-1}^*(J_0^\alpha):=\left( \oplus_{q=1}^{(r,j)-1} \so_{pq}^\C(J_0^\alpha)\right) \cap \so(\rdn)=\ul_{j-1}(J_0^\alpha)/\ul_0(J_0^\alpha)=\so_0((J_0^\alpha)^j)/\so_0(J_0^\alpha)\\
&   & \qquad \quad \quad\ \cong T_{J_0^\alpha}\mZ_{2k,j}^\alpha(\rdn,(J_0^\alpha)^j),
\end{eqnarray*}
the last identification is given by 
$$
A\in \oplus_{q=1}^{(r,j)-1} \so_{pq}^\C(J_0^\alpha)\longmapsto A\cdot J_0^\alpha =[A,J_0^\alpha]\in \oplus_{q=1}^{(r,j)-1}\mal B_{pq}^\C(J_0^\alpha)
$$
(see section~\ref{AdJj}).\\
For the connection form  on $Q^j$ we take 
$$
\omega^j:=\omega_{\hk_j|TQ^j}.
$$
 We set as usual $N^j=Q^j/K$ which  is a Homogeneous fibre bundle over $M$. Moreover the isomorphism of bundle (and isometry) $\mJ\colon N\to \zdk^\alpha(E)$ satisfies 
$$
\mJ(N^j)=\mZ_{2k,j}^\alpha(E,J_j)
$$
by definition of $\mJ$ and $\ Q^j$ (see (\ref{defjk}) and (\ref{defQj})), so that it induces an isomorphism of bundle from $N^j$ onto $N_\mZ^j$.\\
Let us denote by $TN^j=\ver^j \oplus \hor^j$ the splitting in terms of vertical and horizontal  subbundles given by $\omega^j$. Then denoting by $s_j$ the cross section in the associated bundle $Q/H^j=\mal{SO}(E)/\U_{j-1}(J_0^\alpha)$ defining the $H^j$-reduction $Q^j$ (i.e. $\mJ_j\circ s_j= J_j$)\footnote{See remark~\ref{Jj} and \ref{suppl}.}, according to section~\ref{Reductions}, we have the following equivalences
$$
\begin{array}{rcccc}
\omega \text{ is  reducible in  } Q^j \text{ (to } \omega^j) & \overset{\ref{Reductions}}{\Longleftrightarrow} &  s_j  \text{ is horizontal} & \overset{\ref{Reductions}}{\Longleftrightarrow}  & \hor^j=\hor_{|N^j}\\
  &  & \left.\text{\footnotesize Rmk~\ref{suppl}}\right\Updownarrow &  &  \Updownarrow \\
 &   &     \nabla J_j= 0  &  \overset{\text{ex.~\ref{zdkjE}}}{\Longleftrightarrow} &  \hor^{\mZ,j}={\hor^\mZ}_{|N_\mZ^{j}}
\end{array}
$$  
\begin{exam}\label{nabla0j_2}\index{canonical!connection, $G$-invariant}
Let $M= G/H$ be the $k$-symmetric space correponding to some $2k$-symmetric space $G/G_0$ (see section~\ref{evendetercase}), and take $(E,\nabla)=( TM,\onabla{M}{0}{})$, $j=2$ and $J_2$ given by lemma~\ref{lemma2}. Then we have
$$
\onabla{M}{0}{}J_2=0.
$$ 
Indeed   $\onabla{M}{0}{}J_2$ lifts  in $G$ into
$$
(d+\theta_\hk)J_0^2=dJ_0+ [\theta_\hk,J_0^2]=0
$$
(see lemma~\ref{lemma2}). Therefore we can conclude that in this case $\omega$ is reducible in $Q^2$ (to $\omega^2$).
\end{exam}                  
If $\omega$ is not reducible in $Q^j$ (to $\omega^j$), then according to (\ref{metric}) and (\ref{metrich}), $N^j\hookrightarrow N$ and $N_\mZ^{j}\hookrightarrow N_\mZ$ are not isometries, and thus we can not say directly that $\mJ$ induces an isometry from $N^j$ onto $N_\mZ^{j}$, even if as we will see below it is effectively the case. As above, the result of (\ref{metric}) and (\ref{metrich}), and $d\mJ(\ver^j)=\ver^{\mZ,j}$, is that: 
$ \mJ\colon N^j \to N_\mZ^{j}$ is an isometry \iif  $ d\mJ (\hor^j)= \hor^{\mZ, j}$.\\[1mm]
Now let us come back to the connection form $\omega^j\colon TQ^j\to \hk^j\subset \so(2n)$. It defines a metric covariant derivative $\rabla{j}{} $ in the associated vector bundle $E$. Then we have 
$$
\rabla{j}{} J^j=0 .
$$
Indeed $J_j$ lifts into the $H^j$-equivariant (constant) map $\tl\mJ_j\colon e\in Q^j \to (J_0^\alpha)^j\in (\zdk^\alpha(\rdn))^j\subset \mak{gl}_{2n}(\rdn)$ and $\rabla{j}{} J_j$ lifts into
$$
\Bar D^j \tl J^j =d\tl J_j+ [\omega^j,\tl J_j] = 0 + 0=0,
$$
since by definition $\hk^j=\ul_{j-1}(J_0^\alpha)=\so_0((J_0^\alpha)^j)$ commutes with $(J_0^\alpha)^j$.
\begin{rmk}\em
We can do the things more concretely by using a (local) moving frame $e$ in $Q^j$: $\rabla{j}{}$ is then caracterized by 
\begin{equation}\label{nablaje}
\rabla{j}{} (e_1,\ldots,e_{2n})=(e_1,\ldots,e_{2n}).\omega^j(e;de)
\end{equation}
Then by definition of $Q^j$ we have
\begin{equation}\label{defqj}
J_je=e.(J_0^\alpha)^j
\end{equation}
so that 
$$
(\rabla{j}{} J_j)e + J_j (\rabla{j}{} e) = e.\omega^j (J_0^\alpha)^j
$$
then using (\ref{defqj}) and (\ref{nablaje}), we obtain 
$$
(\rabla{j}{} J_j)e= e.\omega^j(J_0^\alpha)^j - J_j(e.\omega^j)= e.(\omega^j(J_0^\alpha)^j - (J_0^\alpha)^j \omega ^j)= 0
$$
since $\omega^j$ takes values in $\hk_j=\so_0((J_0^\alpha)^j)$.
\end{rmk}
In fact we can caracterize $\rabla{j}{}$ in the following more general way, which in particular generalizes a well-known result of Rawnsley \cite[p.~107]{rawnsley} about complex structures on vector bundles.
\begin{thm}\label{decdenabla}
Let $E$ be a Riemannian vector bundle as above. Let $p'\in \mathbb N^*$ and $J\in \mal C(\mal U_{p'}(E))$, then $\Ad J$ defines an automorphism of the linear bundle $\End(E)$ (over $\Id_M$), i.e. a section of $\End(\End(E))$. Then the metric covariant derivative $\nabla$ in $E$ admits an unique decomposition in the form: \footnote{As usual $r_{p'}$ is the order of $\Ad J$, i.e. $ r_{p'}=p'$ if  $p'$ is odd, and if $p'$ is even then $r=p'$ if $J^{\frac{p'}{2}}\neq -\Id$ and $r_{p'}=\dfrac{p'}{2}$ if $J^{\frac{p'}{2}}= -\Id$.}
\begin{equation}\label{decnabla}
\nabla= \overset{J}{\nabla^0} + \sum_{i=1}^{r_{p'}-1} A_i
\end{equation}
where $\overset{J}{\nabla^0}$ is a metric covariant  derivative for which 
$$
\overset{J}{\nabla^0}J=0
$$
and $A_i\in \mal C(T^*M\otimes\so_i^\C(E,J))$, i.e. $JA_iJ^{-1}=\omega_r^i A_i$ and $A_i\in \so(E)^\C$. \\
$\overset{J}{\nabla^0}$ will be called  the $J$-commuting component of $\nabla$, $A_*=\sum_{i=1}^{r_{p'}-1} A_i\in \mal C(T^*M\otimes\so_*(E,J))$ the $\so_*(E,J)$-component of $\nabla$, and $A_i$ the $\so_i^\C(E,J)$-component of $\nabla$.
\end{thm}
\textbf{Proof.} \emph{Uniqueness.} Let us suppose that (\ref{decnabla}) exists then we have 
$$
\nabla J=\sum_{i=1}^{r_{p'}-1} [A_i,J]
$$
so that 
$$
\sum_{i=1}^{r_{p'}-1} A_i=-(\ad J)^{-1}(\nabla J)
$$
(see section~\ref{twistor2}) which proves the uniqueness of $(A_i)_{1\leq i \leq r_{p'}-1}$ (these are determined by $\nabla$ and $J$, more precisely these are the components of $-(\ad J)^{-1}(\nabla J)$). Now  $\overset{J}{\nabla^0}=\nabla - \sum_{i=1}^{r_{p'}-1} A_i$ is also unique.\\[1mm]
\emph{Existence.} Let $\nabla^0$ be any metric covariant derivative commuting with $J$, that is to say $\nabla^0$ corresponds to a connection on the principal bundle of Hermitian  frames on $(E,\langle\ ,\ \rangle, J)$ (such a connection always exists, see \cite{KN}). Then consider
$$
A=\nabla -\nabla^0\in \mal C (T^*M\otimes\so(E))
$$
and let $A=\sum_{i=0}^{r_{p'}-1} A_i$ be the decomposition of $A$ following 
$\so(E,J)^\C=\oplus_{i=0}^{r_{p'}-1} \so_i^\C(E,J)$. Let us set
$$
\overset{J}{\nabla^0}=\nabla^0 + A_0
$$
then $\overset{J}{\nabla^0}$ is a $J$-commuting metric covariant derivative in $E$ and we have 
$$
\nabla= \overset{J}{\nabla^0} + \sum_{i=1}^{r_{p'}-1} A_i
$$
which proves the existence.\hfill$\square$\\[1.5mm]
Applying   this theorem to $J_j$, we obtain the following.
\begin{cory}
$\rabla{j}{}$ is the  $J_j$-commuting component of $\nabla$.
\end{cory}
\textbf{Proof.}
The $H$-equivariant lift of $\nabla$ is the following derivative on $Q$:\footnote{Remark that here $r_{p'}$ is the order of $\Ad (J_0^\alpha)^j$, so that $r_{p'}=\dfrac{r}{(r,j)}=p$.}
\begin{equation}\label{dplusomega}
d + \omega =(d + \omega_0) + \sum_{i=1}^{r_{p'} -1} \omega_i
\end{equation} 
where $ \omega_i=[\omega]_{\so_i^\C((J_0^\alpha)^j)}$, and in particular $\omega_0=\omega^j$. Then restricting (\ref{dplusomega}) to $Q^j$, and projecting on $M$, we obtain the decomposition~(\ref{decnabla}) of $\nabla$:
$$
\nabla=\rabla{j}{} + \sum_{i=1}^{r_{p'}-1} A_i
$$
that is to say $d+ \omega_0$ is the $H^j$-equivariant lift of $\overset{J_j}{\nabla^0}$, which is thus equal to  $\rabla{j}{}$, and $\omega_i$ is the $H^j$-equivariant lift of the $\so_i^\C(E,J_j)$-valued 1-form on $M$, $A_i$. This completes the proof.
\begin{rmk}\em
Moreover $\omega$ is reducible in $Q^j$ ($\nabla J_j=0$) \iif $\rabla{j}{}=\nabla$.
\end{rmk}
\begin{rmk}\em
Under the hypothesis of theorem~\ref{decdenabla} we have 
$$
\forall F\in \mal C(\mal A_0(E,J)),  \quad \overset{J}{\nabla^0} F=\mrm{pr}_{\mal A_0(E,J)}\circ \nabla F
$$
where $\mrm{pr}_{\mal A_0(E,J)}\colon \End(E)\to \mal A_0(E,J)$ is the orthogonal projection  (i.e. along $\mal A_*(E,J)$) so that in particular
$$
\forall F\in \mal C(\so_0(E,J)),  \quad \overset{J}{\nabla^0} F=\mrm{pr}_{\so_0(E,J)}\circ \nabla F
$$
where $\mrm{pr}_{\so_0(E,J)}\colon \so(E)\to \so_0(E,J)$ is the orthogonal projection. 
Indeed, 
$$
\nabla F= \overset{J}{\nabla^0} F + \sum_{i=1}^{r_{p'} -1} [A_i,F]
$$
and $J$ commutes with $\overset{J}{\nabla^0}$ and $F$ so with $\overset{J}{\nabla^0} F$ : $(\overset{J}{\nabla^0} F) . J = \overset{J}{\nabla^0}(F.J)- F\overset{J}{\nabla^0} J= \overset{J}{\nabla^0} (J.F)= J\overset{J}{\nabla^0}F$. Moreover $[A_i,F]\in [\mal A_i(J),\mal A_0(J)]\subset \mal A_i(J) $, so that we can conclude. $\centerdot$

\end{rmk}
The canonical $2k$-structure in $\pi^*E$, $\mJ \colon N\to \pi^*E$ induces by restriction a $2k$-structure in ${\pi^j}^*E$, still denoted by  $\mJ\colon N^j\to {\pi^j}^*E$.\\[1mm]
Now, let us specify the subbundles $\pk_{Q^j}^j$ and $\kk_{Q^j}$. First, we have $\hk_{Q^j}^j=\ul_{0}(E,J_j)$\footnote{ The restriction of the identification (\ref{hqsoe}), $\hk_Q\cong\so(E)$ to $\hk_{Q^j}^j$ gives rises to an identification $\hk_{Q^j}^j:=Q^j\times_{H^j} \hk^j\cong\so_0(E,J_j)$.}  and then 
\begin{eqnarray*}
\kk_{Q^j} & = & {\kk_Q}_{|Q^j}=\so_0({\pi^j}^*E,\mJ)\\
\pk_{Q^j}^j & = & \ul_{j-1}^*(E,\mJ)=\left(\oplus_{i\in p.\Z_r\setminus\{0\}}\so_i^\C({\pi^j}^*E,\mJ)\right)\bigcap\so({\pi^j}^*E). 
\end{eqnarray*}
The morphism of vector bundle $\ad \mJ\colon \so(\pi^*E)\to \mJ.\pk_Q$ induces a surjective morphism from ${\pi^j}^*\hk_{Q^j}^j=\ul_{j-1}({\pi^j}^*E,\mJ)$ onto $\mJ.\pk_{Q^j}^j$, with kernel $\kk_{Q^j}$:
$$
\begin{array}{crcl}
\ad \mJ\colon & {\pi^j}^*\hk_{Q^j}^j=\ul_{j-1}({\pi^j}^*E,\mJ) & \longrightarrow & \mJ.\ul_{j-1}^*({\pi^j}^*E,\mJ) 
=\mJ.\pk_{Q^j}^j\\
 & (J,A) & \longmapsto & \ad J(A)=[J,A]=J\displaystyle\sum_{i=1}^{(r,j)-1} (1-\omega_r^{ip})A_{ip}
\end{array}
$$
where $A_i=[A]_{\so_i^\C(E_x)}$.\\
As above, now we express the homogeneous fibre bundle tools $\phi^j$, $\Phi^j$ and $\nabla^{\pk^j}$ in terms of $\mJ$.
\begin{thm}\label{entermdeJj}
If $A,B\in TN^j$, $F\in\mal C(\pk_{Q^j}^j)$ then
\begin{description}
\item[(i)] $\rabla{j}{} \mJ= -\ad \mJ\circ\phi^j$ thus $\phi^j A=  -(\ad \mJ)^{-1}\rabla{j}{A}\mJ$
\item[(ii)] $\Phi^j(A,B)=(\ad \mJ)^{-1}\left[\mJ,{\pi^j}^*R^{\rabla{j}{}}(A,B)\right]$ where $R^{\rabla{j}{}}$ is the curvature of $\rabla{j}{}$.
\item[(iii)] $\nabla_A^{\pk^j} F=(\ad \mJ)^{-1}\left[\mJ,\rabla{j}{A} F\right]$
\end{description}
\end{thm}
In the following theorem, we use the notation of \ref{Reductions}. In particular, we denote by " $\cdot^\mrm{v}$ " instead of " $\cdot^v$ " the vertical component in $\ver^j\subset TN^j$.
\begin{thm}\label{verticalj}
Let $s\in \mal C(\pi^j)$ and  $J=s^*\mJ$ be the corresponding $2k$-structure. Then
\begin{description}
\item[(i)] $I(d^\mrm{v} s)= -(\ad J)^{-1}\rabla{j}{} J$. Thus $s$ is horizontal \iif $J$ is $\rabla{j}{}$-parallel.
\item[(ii)] $I(\Pi^\mrm{v}(s)) = -(\ad J)^{-1}(\rabla{j}{})^2J + \dfrac{1}{2}(\ad J)^{-1} \left[\rabla{j}{} J \odot(\ad J)^{-1}\rabla{j}{} J\right]$.\\[1mm]
Thus $s$ is superflat \iif $(\rabla{j}{})^2J - \dfrac{1}{2} \left[\rabla{j}{} J \odot(\ad J)^{-1}\rabla{j}{} J\right]$ commutes with $J$.
\item[(iii)] $I(\tau^\mrm{v}(s))= +(\ad J)^{-1} {\rabla{j}{}}^* \rabla{j}{} J  + \dfrac{}{}(\ad J)^{-1} \mrm{Tr}\left( \left[\rabla{j}{} J ,(\ad J)^{-1}\rabla{j}{} J\right]\right) $.\\[1mm]
Thus $s$ is a harmonic section \iif ${\rabla{j}{}}^* \rabla{j}{} J + \dfrac{}{}\mrm{Tr}\left( \left[\rabla{j}{} J ,(\ad J)^{-1}\rabla{j}{} J\right]\right)$ commutes with $J$. 
\item[(iv)] $s^*\Phi^j=(\ad J)^{-1}\left[J,R^{\rabla{j}{}}\right]$. 
\end{description}
These properties hold also for maps $f\in\mal C^\infty(L,N)$, $(L,b)$ being a Riemannian manifold: (i),(ii),(iii) without any change and (iv) becoming $f^*\Phi^j=(\ad J)^{-1}\left[ J, u^*R^{\rabla{j}{}}\right]$, with $u=\pi\circ f$.
\end{thm}
\proof Let us endow $E$ with $\rabla{j}{}$ and apply the theorems~\ref{entermdeJ2} and \ref{vertical2} (with the Riemannian vector bundle $(E,\rabla{j}{})$), then by restriction to $N^j$, we obtain theorems~\ref{entermdeJj} and \ref{verticalj}. Indeed,  in this case $\omega$ is reducible in $Q^j$ and then everything corresponds via the reduction $N^j\hookrightarrow N$, according to propositions~\ref{prop-tool-corresp} and \ref{prop-tau&Pi-corresp}. \comprf \hsq \medskip\\ 
As above, from theorem~\ref{entermdeJj}-(i), we conclude that $d\mJ$ sends  the decomposition $TN^j=\ver^j\oplus \hor^j$ onto the decomposition $TN_\mZ^j=\ver^{\mZ,j}\oplus\hor^{\mZ,j}$  (see example~\ref{zdkjE})  so that we can consider $\pi_\mZ^j\colon N_\mZ^j(E)\to M$ as a homogeneous fibre bundle over $M$ with  structure groups $H^j=\U_{j-1}(J_0^\alpha)$ and $K=\U_0(J_0^\alpha)$. We will call this structure \emph{the homogeneous fibre bundle structure defined by (the $J_j$-commuting part of) $\nabla$}.\\
Besides, since  the vertical and horizontal subbundles corresponds via $\mJ$, then we can conclude according to (\ref{metric}) and (\ref{metrich}) that $\mal J\colon N^j\to N_\mZ^{j}$ is an isometry.\\[1mm]
Moreover, the vertical tension field of $J$ in $N_\mZ^j=\mZ_{2k,2}^\alpha(E,J_j)$ is given by
$$
d\mJ(\tau^\mrm{v}(s))
 =-\left[ {\rabla{j}{}}^*\rabla{j}{} J + \dfrac{}{}\mrm{Tr}\left( \left[\rabla{j}{} J ,(\ad J)^{-1}\rabla{j}{} J\right]\right)\right]_{\ver^{\mZ,j}} 
$$ 

\begin{rmk}\label{nablac}\em
According to \ref{Reductions}, the canonical connection  in $\pk_{Q^j}^j\to N^j$ is the restriction of the canonical connection in $\pk_Q\to N$, to $\pk_{Q^j}^j$.
\end{rmk}
\begin{rmk}\label{everythingreducible}\em
If we endow $E$ with $\rabla{j}{}$ and apply the theorems~\ref{entermdeJ2} and \ref{vertical2} (with the Riemannian vector bundle $(E,\rabla{j}{})$), then by restriction to $N^j$, we obtain theorems~\ref{entermdeJj} and \ref{verticalj}.\\
In particular, superflatness and vertical harmonicity (for sections in $N^j$) are the same in $N^j$ and $N$. This is what happens in particular in example~\ref{nabla0j_2}.
\end{rmk}
%
%%%%%%%%%%%%%%%%%%%%%%%%%%%%%%%%%%%%%%%%%%%%%%%%
%
\paragraph{The particular case of $\mZ_{2k,2}^\alpha(E,J_2)$}
According to  theorem~\ref{embedding}, we will be especially interested by  this subcase in our interpretation of the even elliptic integrable system.
In this subcase the fibre $H^2/K=\mZ_{2k,2}(\rdn, (J_0^\alpha)^2)$ is symmetric so that we obtain simplifications (coming in particular from  the facts that $\nabla^c=\nabla^{\pk}$ and that any section $J\in \mal C(\pi_\mZ^2)$ anticommutes with $\rabla{2}{} J$) in theorems~\ref{entermdeJj}  and  \ref{verticalj} which   then take the same forms as theorems~\ref{entermdeJ} and \ref{vertical} about the twistor bundle $\Sigma^\eps(E)$, just by doing the change $\nabla\longleftrightarrow \rabla{2}{}$. Therefore the case $\mZ_{2k,2}^\alpha(E,J_2)$ is very similar to that of $\Sigma^\eps(E)$. \\
Before writing the simplified theorems for $j=2$, let us do some useful observations.\\[1mm]
First, we have\footnote{with notation defined in remark~\ref{nouvellenotation}.} 
\begin{eqnarray}\label{k2-p2}
\kk_{Q^2} & = & \so_0({\pi^2}^* E,\mJ)=\{ A \in \so({\pi^2}^* E) |\, [A,\mJ]=0\}=\so_{(+1)}({\pi^2}^* E,\mJ)\\
\pk_{Q^2}^2 & = & \so_{\frac{r}{2}} ({\pi^2}^* E,\mJ)=\{ A \in \so({\pi^2}^* E) |\, A.\mJ + \mJ .A=0\}=\so_{(-1)}({\pi^2}^* E,\mJ).\qquad 
\end{eqnarray}
Then $\ad \mJ$ induces  a surjective morphism from ${\pi^2}^*\hk_{Q^2}^2=\ul_1({\pi^2}^*E,\mJ)$ onto  $\mJ.\pk_{Q^2}^2=\mal B_{\frac{r}{2}} ({\pi^2}^* E,\mJ)$ with kernel $\kk_{Q^2}$
$$
\begin{array}{crcl}
\ad \mJ\colon & \ul_{1}({\pi^2}^*E,\mJ)= \so_{(+1)} ({\pi^2}^*E,\mJ) \oplus  \so_{(-1)} ({\pi^2}^*E,\mJ) & \longrightarrow & \dfrac{}{} \mal B_{(-1)}({\pi^2}^*E,\mJ) = \mJ.\so_{(-1)}({\pi^2}^*E,\mJ)  \\
    & (J, A_0 + A_1) & \longmapsto & \dfrac{}{}\ad J(A)=[J,A]= 2J A_1
\end{array}
$$
where we denote by $A_0 + A_1$ the decomposition following  $\so_{(+1)} ({\pi^2}^*E,\mJ) \oplus  \so_{(-1)}({\pi^2}^*E,\mJ)$ instead of $ A_0 + A_{\frac{r}{2}}$.

\begin{thm}\label{entermdeJj2}
If $A,B\in TN^2$, $F\in\mal C(\pk_{Q^2}^2)$ then 
\begin{description}
\item[(i)]  $\phi^2 A=  - \dfrac{1}{2}\mJ^{-1} \rabla{2}{} \mJ$
\item[(ii)] $\Phi^2(A,B)= \dfrac{1}{2}\mJ^{-1}\left[\mJ,{\pi^2}^*R^{\rabla{2}{}}(A,B)\right]$ where $R^{\rabla{2}{}}$ is the curvature of $\rabla{2}{}$.
\item[(iii)] $\nabla_A^{c} F=\dfrac{1}{2}\mJ^{-1}\left[\mJ,\rabla{2}{A} F\right]$ \footnote{see remark~\ref{nablac}}
\end{description}
\end{thm}
\begin{thm}\label{verticalj=2}
Let $s\in \mal C(\pi^2)$ and  $J=s^*\mJ$ be the corresponding $2k$-structure. Then
\begin{description}
\item[(i)] $I(d^\mrm{v} s)= -\dfrac{1}{2} J^{-1}\rabla{2}{A} J$. Thus $s$ is horizontal \iif $J$ is $\rabla{2}{}$-parallel.
\item[(ii)] $I(\Pi^\mrm{v}(s)) = -(\ad J)^{-1}\left(\rabla{2}{}\right)^2J=-\dfrac{1}{4}\left[J^{-1},\left(\rabla{2}{}\right)^2J\right] $.\\[1mm]
Thus $s$ is superflat \iif $\left(\rabla{2}{}\right)^2J $ commutes with $J$.
\item[(iii)] $I(\tau^\mrm{v}(s))= (\ad J)^{-1} {\rabla{2}{}}^* \rabla{2}{} J = \dfrac{1}{4}\left[J^{-1},{\rabla{2}{}}^* \rabla{2}{} J\right ] $.\\[1mm]
Thus $s$ is a harmonic section \iif ${\rabla{2}{}}^* \rabla{2}{} J $ commutes with $J$. 
\item[(iv)] $s^*\Phi^2=\dfrac{1}{2} J^{-1}\left[J,R^{\rabla{2}{}}\right]$. 
\end{description}
These properties hold also for maps $f\in\mal C^\infty(L,N)$, $(L,b)$ being a Riemannian manifold: (i),(ii),(iii) without any change and (iv) becoming $f^*\Phi^2=\dfrac{1}{2} J^{-1}\left[J,u^*R^{\rabla{2}{}}\right]$, with $u=\pi\circ f$.
\end{thm}
Let us add that the vertical tension field in $N_\mZ^2$ is given by
\begin{equation}\label{tauvj=2}
d\mJ(\tau^\mrm{v}(s))
 =-\left[ {\rabla{2}{}}^*\rabla{2}{} J \right]_{\ver^{\mZ,2}} =-\dfrac{1}{2}J\left[J^{-1}, {\rabla{2}{}}^*\rabla{2}{} J \right]
\end{equation}
\index{homogeneous fibre bundle|)}
\subsection{Geometric interpretation of the even minimal determined system}\label{geom-interpret}
\index{determined minimal@determined, minimal|(}\index{even case|(}
\subsubsection{The injective morphism of homogeneous fibre bundle $\mijo\colon G/G_0\hookrightarrow \mZ_{2k,2}^{\alpha_0}(G/H,J_2)$.}\label{geom-interpret-sub1}
\index{canonical!embedding|(}
Here, we want to ask ourself if the inclusion $\mijo\colon G/G_0\hookrightarrow \mZ_{2k,2}^{\alpha_0}(G/H,J_2)$ given by theorem~\ref{embedding} conserves the homogeneous fibre bundle structure, in particular: the vertical harmonicity is it conserved. We use the notations of \ref{homspacefibr} and \ref{5.3.4} (with $E=TM$, $\nabla$ a metric connection on $M$ and $j=2$)\footnote{That is to say the notations  for $\mZ_{2k,2}^{\alpha_0}(G/H,J_2)$ will have the subscript "2"  and these of  $G/G_0$ will not have  subscript according to \ref{homspacefibr} and \ref{5.3.4}}. First, we see that  $\mijo$ is obtained under quotient from the following injective morphism of bundle (which is an embedding if $G$ is closed in $\Is(M)$):
\begin{equation}\label{lesQ}
\begin{array}{crcl}
\mai_{e_0}\colon &  G & \hookrightarrow  & Q^2= \mal U_1^{\alpha_0}(G/H,J_2)\subset \mal{SO}(M)\\
                 & g  & \longmapsto  &  g\cdot e_0
\end{array}
\end{equation}
where $e_0\in \mal{SO}(T_{p_0}M)$ is such that $\mal{M}at_{e_0}(J_0)=J_0^{\alpha_0}$, and $J_0=\tm$. In other words $G\to M$ is a reduction of $\mal U_1^{\alpha_0}(G/H,J_2)\to M$ itself a reduction of $\mal{SO}(M)\to M$. \\
Further quotienting in (\ref{lesQ}) by $\U_0(J_0^\alpha)$  the target space and then by $G_0$  the domain, we obtain (by definition of $G_0=G^\tau\cap H$, see theorem~\ref{embedding}) the injective morphism of bundle 
$$
\mai_{\overline{e_0}}\colon g.G_0 \longmapsto (g\cdot e_0)\U_0(J_0^\alpha)\in \mal U_1^{\alpha_0}(M,J_2) / \U_0(J_0^\alpha) \subset \mal{SO}(M)/\U_0(J_0^\alpha)
$$ 
where $\overline{e_0}=e_0\U_0(J_0^\alpha)\in N^2$, and finally composing with $\mJ$ (in the target space) we obtain the map $\mijo$:
$$
g.G_0\longmapsto g\cdot (e_0\U_0(J_0^\alpha))\overset{\mJ}{\longmapsto} J=gJ_0 g^{-1}\in \mZ_{2k,2}^{\alpha_0}(M,J_2).
$$
Since $\mijo$ (resp. $\mai_{\overline{e}_0}$) is an injective morphism of bundle (and an immersion) $d\mijo$ (resp. $d\mai_{\overline{e}_0}$) induces an injective morphism of bundle from the vertical subbundle $\ver^{G/G_0}=[\g_k]$ into the vertical subbundle $\ver^{\mZ,2}$ (resp. $\ver^2$).\\
$\mijo$ is the restriction to $G/G_0$ of the inclusion map $\mai\colon \End (G/H)  \to M\times \End (\g)$ (see \ref{End}). Indeed, we have the inclusion depending on $J_0$: $g.G_0 \in G/G_0 \mapsto [g,J_0]\in G\times_H\End(\mk)=\End(G/H)$ which under the inclusion $\mai$ gives $g.G_0\in G/G_0\mapsto (g.x_0, \Ad g\circ\tm\circ\Ad g^{-1})\in M\times \End(\g)$ which is in nothing but $\mijo$ (as usual under the identification $TM=[\mk]$). Then under the inclusion $\hk_G\subset\so(TM)$, we have $\hk_G\subset\hk_{Q^2}^2=\ul_0(TM,J_2)$. Indeed, under the linear isotropy representation of $H$ in $T_{x_0}M$, we have $H\subset \U_0(T_{x_0}M,J_0^2)= \U_1(T_{x_0}M,J_0)$ so that $\hk\subset\ul_0(T_{p_0}M,J_0^2)$ and thus $\hk_G:=G\times_H\hk\subset \ul_0(TM,J_2)$. Moreover let  us remark that $\pi^2\circ\mai_{\overline e_0}=\pi$ so that $\pi^*\hk_G\subset {\pi^2}^*\hk_{Q^2}^2$ over $\mai_{\overline e_0}\colon N\to N^2$ (i.e. the inclusion is a morphism of bundle over $\mai_{\overline e_0}$).
\\
Furthermore, since $\Ad J_0$ leaves invariant $\hk\subset \ul_1(T_{x_0}M,J_0)$, the restriction to $\hk$ of the symmetric decomposition
$$
\ul_1(T_{x_0}M,J_0)=\so_{(+1)}(T_{x_0}M,J_0)\oplus\so_{(-1)}(T_{x_0}M,J_0)
$$
gives rise to the decomposition $\hk=\g_0\oplus\g_k$ according to (\ref{go-gk}), so that the symmetric decomposition given by $\Ad\mJ$ on ${\pi^2}^*\hk_{Q^2}^2=\ul_1({\pi^2}^*TM,\mJ)$, that is to say 
$$
\ul_1(TM,\mJ)=\so_{(+1)}({\pi^2}^*TM,\mJ)\oplus\so_{(-1)}({\pi^2}^*TM,\mJ)
$$
gives rise in the $\Ad \mJ$-invariant subspace $\pi^*\hk_G\subset {\pi^2}^*\hk_{Q^2}^2$ to the symmetric decomposition of $\Ad \mijo$ (restricted to $\pi^*\hk_G\subset \so(\pi^*TM)$)
$$
\pi^*\hk_G=\kk_G\oplus\pk_G
$$
according to (\ref{kg-pg}). In other words, the decomposition given by (\ref{kg-pg}) injects into the decomposition given by (\ref{k2-p2}) via the inclusion $\pi^*\hk_G\subset {\pi^2}^*\hk_{Q^2}^2$.\\[1mm]
Now let us interpret theorems~\ref{entermdeij0} and \ref{verticalij0} using the homogeneous fibre bundle structure in $\mZ_{2k}^{\alpha_0}(M,J_2)$ defined by the Riemannian vector bundle $(E,\nabla)=(TM,\onabla{M}{0}{})$ (in the sense of \ref{5.3.4}). We continue to use the same conventions for the notations in $N$ and $N^2$ (no subscript for $N$ and subscript 2 for $N^2$ and $N^{\mZ,2}$). Recall that we have $\mijo=\mJ\circ\mai_{\overline e_0}=\mai_{\overline e_0}^*\mJ$ and that $(\onabla{M}{0}{})^{[2]}=\onabla{M}{0}{}$. Then according to theorems~\ref{entermdeJj2} and \ref{verticalj=2},  theorems~\ref{entermdeij0} and \ref{verticalij0} implies
\begin{thm}\label{entermdenabla0} We have the following identities
\begin{description}
\item[(i)] $\phi=\mai_{\overline e_0}^*\phi^2$
\item[(ii)] $\Phi =\mai_{\overline e_0}^*\Phi^2$
\item[(iii)] $\nabla^c=\onabla{N}{0}{}|[\pk]=\mai_{\overline e_0}^*\nabla^{c,2}$, where $\nabla^{c,2}$ is the canonical connection in $\pk_{Q^2}^2$.
\end{description}
\end{thm}
\begin{thm}\label{verticalnabla0}
Let $s\in\mal C(\pi)$ and identify it (temporarily) with $s^*\mai_{\overline e_0}\in \mal C(\pi^2)$. Then under the inclusion $\mai_{\overline e_0}\colon N\to N^2$, we have:
\begin{description}
\item[(i)] $d^v s= d^{\mrm v, 2} s$
\item[(ii)] $\Pi^v s= \Pi^{\mrm v, 2}$  
\item[(iii)] $\tau^v s =\tau^{\mrm v,2} s$
\item[(iv)]  $s^*\Phi=s^*\Phi^2$  
\end{description}
These properties holds also, without any change, for maps $f\in \mal C^\infty (L,N)$, $(L,b)$ being a Riemannian manifold.
\end{thm}
Let us remark that since the connection form $\omega$, on $Q=\mal{SO}(TM)$ defined by $\onabla{M}{0}{}$ is reducible in $Q^2$, then in the previous theorems all the "quantities" in $N^2$ (right handside) can also be computed in $\mal{SO}(TM)/\U_0(J_0^{\alpha_0})\cong \zdk^{\alpha_0}(M)$, since "everything is reducible" in this case (see remark~\ref{everythingreducible}). \\[1.5mm]
Now, let us compute the vertical tension field of $J\colon L\to N_\mZ^2$ for the homogenous fibre bundle structure defined in $N_\mZ^2$ by $\onabla{M}{0}{}$: according to (\ref{tauvj=2}) we have
\begin{equation}
\tau^{\mrm v,2} (J)=-\dfrac{1}{2}J \left[ J^{-1}, (\onabla{M}{0}{})^*\onabla{M}{0}{} J\right]. 
\end{equation}
Then suppose that $J$ takes values in  $\mijo(G/G_0)$ i.e. $J=f^*\mijo$ for some $f\in \mal C^\infty(L,N)$, then according to theorem~\ref{verticalnabla0} (and $\mijo=\mJ\circ\mai_{\overline e_0}$) we have
$$
d\mijo (\tau^v(f))= d\mJ (\tau^{\mrm v,2} (\bar f))=\tau^{\mrm v,2}(J)
$$
where $\bar f=\mJ^{-1}\circ J$ i.e. $J={\bar f}^*\mJ$.\\
\textbf{The tension fields (and thus vertical harmonicity) correspond} (up to multiplicative constant) \textbf{via the different inclusions and identifications, in particular via} $\mijo\colon N\to N_\mZ^2$.
%
%%%%%%%%%%%%%%%%%%
%
\paragraph{The canonical embedding  is a reduction of homogeneous fibre bundles.}
In fact, we can recover theorems~\ref{entermdenabla0} and \ref{verticalnabla0} as well as theorems~\ref{entermdeij0}, \ref{verticalij0} by a more conceptual way: just remark that $\mijo\colon N\to N_\mZ^2$ is a reduction of homogeneous fibre bundles.
\begin{thm}\label{thm-inj-morph-bund} The injective morphism of bundle $\mijo\colon G/G_0\hookrightarrow \mZ_{2k,2}^{\alpha_0}(G/H,J_2)$ is a  reduction of homogeneous fibre bundles.
\end{thm}
\proof 
According to remark~\ref{rmk-can-con-restr}, the canonical connection in the $H$-principal bundle $G(G/H,H)$ corresponds (via the natural injection  $\mai_{e_0}$ defined by equation \eqref{lesQ}) to the canonical connection in the $G$-invariant $G$-structure $P'=G.e_0$, itself the restriction  to $P'$ of the canonical connection in the $G$-invariant $\U_1(J_0^{\alpha_0})$-structure $Q^2=\mal U_1^{\alpha_0}(G/H,J_2)$, itself the the restriction  to $Q^2$ of the canonical connection in  $\mal{SO}(M)$.
\comprf\hsq\medskip\\
Now, applying  propositions~\ref{prop-tool-corresp} and \ref{prop-tau&Pi-corresp} we obtain:
\begin{cory}
The theorems~\ref{entermdenabla0} and \ref{verticalnabla0} as well as theorems~\ref{entermdeij0}, \ref{verticalij0} are corollaries of theorems~\ref{entermdeJj2} and \ref{verticalj=2}.
\end{cory}
\paragraph{Using the Levi-Civita connection.}
In fact in what precedes we can replace the canonical connection in $M$, $\onabla{M}{0}{}$, by  (the $J_2$ commuting part of) the Levi-Civita connection in $M$.
\begin{prop}\label{prop-J2-comm-part}
The canonical linear connection on $M$ is the $J_2$-commuting component of the Levi-Civita connection $\onabla{M}{}{}$ on $M$:
$$
\onabla{M}{0}{}=\onabla{J_2}{0}{}=\onabla{M}{[2]}{}.
$$
\end{prop}
\proof  According to subsection~\ref{family}, we have 
$$
\onabla{M}{}{}=\onabla{M}{0}{} + \dfrac{1}{2}\left([\ ,\ ]_{[\mk]} +  \mU^M\right).
$$
Then it suffices to apply equation \eqref{eq-Ad-taum} in subsection~\ref{Action-of-Ad-taum}. \comprf \hsq\medskip\\
Now, we can rewrite the previous results by replacing the canonical connection by (the $J_2$ commuting part of) the Levi-Civita connection.
\begin{cory}
The homogeneous fibre bundle structures in $N_\mZ^2$ defined by the canonical affine coonection $\onabla{M}{0}{}$ and by (the $J_2$-commuting  part of) the Levi-Civita connection $\onabla{M}{}{}$, in $M$, are the same. Therefore theorems~\ref{entermdeij0}, \ref{verticalij0} and corollary~\ref{interpretnabla0} still hold if we replace $\onabla{M}{0}{}$ by $\onabla{M}{[2]}{}$. Moreover  theorems~\ref{entermdenabla0} and \ref{verticalnabla0} hold with the homogeneous fibre bundle structure defined in $N^2$ by the ($J_2$-commuting  part of) the Levi-Civita connection $\onabla{M}{}{}$.
\end{cory}
\paragraph{Some additionnal identities}
Let us conclude this subsection by some additionnal equalities.
\begin{prop}
The canonical linear connection on $M$ is the $J_2$-commuting component of the connections $\onabla{\mrm{met}}{t}{}$ on $M$:
$$
\onabla{M}{0}{}=\left( \onabla{\mrm{met}}{t}{}\right) ^{[2]}.
$$
\end{prop}
\proof 
Same proof as for proposition~\ref{prop-J2-comm-part}.  \hsq
\begin{prop}
Let $J_1\in \mal C(\mal U_{2k}^*(N))$ be the section defined by $\breve J_0=\tau_{|\mak n}^{-1}$ with, let us recall it, $\mak n=\pk\oplus\mk=\g_k\oplus\mk$, then
\begin{description}
\item[(i)]
The $J_1$-commuting component of the ($\pi$-pullback of the) canonical linear connection in $M$, $\pi^*\onabla{M}{0}{}$, is $\onabla{N}{0}{}$  the canonical linear connection in $N$. This latter is also the $J_1$-commuting component of the ($\pi$-pullback of the) Levi-Civita connection, and more generally of the connections $\onabla{\mrm{met,M}}{t}{}$.
\item[(ii)] The $J_1$-commuting component of the Levi-Civita connection in $N$, $\onabla{N}{}{}$, is $\onabla{N}{0}{}$.
\item[(iii)] More generally, the $J_1$-commuting component of $\onabla{\mrm{met,N}}{t}{}$ is $\onabla{N}{0}{}$.
\item[(iv)] Let $s\in \mal C(\pi)$ and $J=s^*\mijo$ the corresponding $2k$-structure on $M$, then $s^*\onabla{N}{0}{}$ is the $J_1$-commuting component of $\onabla{M}{0}{}$, and  also the $J_1$-commuting component of  (the $s$-pulback of) the Levi-Civita connection on $M$, $s^*\onabla{M}{}{}$; and more generally of (the $s$-pulback of)
the connections $\onabla{\mrm{met,N}}{t}{}$ .
\end{description}
\end{prop}
%
%We recover in particular \cite{bk} from (iv).
%
\proof
The first point of (i) becomes obvious when it is written in terms of the lift of the connections: the $\breve{J}_0$-commuting part of $d+\theta_\hk= d + \theta_0 + \theta_k$ is \ $ d + \theta_0 $ (whereas $\theta_k$ is its anticommuting part). For the second point of (i), (ii), (iii), and (iv) use the same method as for proposition~\ref{prop-J2-comm-part}.  \hsq

%%%%%%%%%%%%

\subsubsection{Conclusion}
Now\footnote{We still consider the same situation as in \ref{geom-interpret-sub1}} we can conclude:
\begin{thm}\index{primitive@primitive, case or system}
Let $(L,j)$ be a Riemann surface, $f\colon L\to N=G/G_0$ be a map and $J=f^*\mijo$ the corresponding map into $\mZ_{2k,2}^{\alpha_0}(M,J_2)$. $\bullet$ Then $f$ is a geometric solution of the even minimal determined system $(\syst(k,\tau))$ \iif
\begin{description}
\item[(i)] $J$ is an admissible twistor lift ($\Leftrightarrow$ $f$ is horizontally holomorphic)
\item[(ii)] $J$ is vertically harmonic in $\mZ_{2k,2}^{\alpha_0}(M,J_2)$ endowed with its homogeneous fibre bundle structure defined by the Levi-Civita connection, $\nabla$, in $M$: 
$$
\left[{\rabla{2}{}}^*\rabla{2}{} J,J  \right]=0,
$$ 
where $\rabla{2}{}$ is the $J_2$-commuting component of $\nabla$. ($\Leftrightarrow$ $f$ is vertically harmonic in $G/G_0$).
\end{description}
$\bullet$ Moreover the first condition implies that $J$ is flat in $\mZ_{2k,2}^{\alpha_0}(M,J_2)$:
$$
J^*\Phi^{\mZ,2}=\left[ u^* R^{\rabla{2}{}}, J \right]=0,
$$ 
where $\Phi^{\mZ,2}$ is the homogeneous curvature form in $\mZ_{2k,2}^{\alpha_0}(M,J_2)$, which means also that $J$ is a flat section in $\End(u^*TM, u^*\rabla{2}{})$. ($\Leftrightarrow$ $f$ is flat in the homogeneous fibre bundle $N\to M$).\\
$\bullet$ Furthermore $f$ is a primitive geometric solution (i.e. there exists $m\leq k$ such that $f$ is $m$-primitive, which is equivalent to say that $f$ is $k$-primitive) \iif
\begin{description}
\item[(i)] $J$ is an admissible twistor lift
\item[(ii)] $J$ is parallel: $\rabla{2}{} J=0$  ($\Leftrightarrow$ f is horizontal).
\end{description}
$\bullet$ Besides in the two previous characterizations, at each point (ii), the Levi-Civita connection can be replaced by any $G$-invariant metric connection $\nabla'$ whose the $J_2$-commuting component ${\nabla'}^{[2]}$ leaves invariant $\hk_G\subset \so(TM)$. This is the case in particular for the connections 
$$
\onabla{\mrm{met}}{t}{}=\onabla{M}{0}{} + t\left( [\ ,\ ]_{[\mk]} + \mU^M\right) ,\quad 0\leq t \leq 1,
$$ 
for which the $J_2$-commuting component is the canonical connection on $M$: $\onabla{M}{0}{}$.
\end{thm}
\index{vertically harmonic|)}\index{determined minimal@determined, minimal|)}\index{twistor|)}
\index{canonical!connection, $G$-invariant|)}
\index{horizontal subbundle|)}
\index{vertical subbundle|)}\index{even case|)}\index{vertical tension field|)}\index{canonical!embedding|)}

\subsection{Bibliographical remarks and summary of the results.}
Since our "generalized" twistor spaces have never been investigated before. A fortiori, it is the same for the study of their structures of homogeneous fibre bundles, as well as for the study of the preservation of these structures under the different inclusions $\mal Z_{2k}^\alpha (E,J_j) \to \zdk^\alpha (E)$. In this section, we refered to papers of C.M. Wood \cite{cmw2,cmw1} for the  definitions and properties of homogeneous fibre bundles and vertical harmonicity. Then we have applied these to the study of our "generalized" twistor spaces. Let us remark that the techniques that we have taken from \cite{cmw2} are in fact known by people of the "English school" like Rawnsley, Salomon, Burstall and al., \cite{rawnsley,Salomon,BuRaw,Obrian-rawnsley}). In particular the idea to use the tautological canonical complex structure $\mal J$ has been for example already used by  Rawnsley \cite{rawnsley}. However, C.M. Wood  investigated vertically harmonic sections in different types of bundle and in particular homogeneous fibre bundles \cite{cmw2}. In some sense, our study in  this section~\ref{vertically-harmonic} fits  in with the general spirit of the works on twistor methods done by several people of the "English school" (see more particulary \cite{rawnsley}).\smallskip\\
The originality of our study as well as its additional intrinsic difficulty come from the more complicated algebraic structure on the fibre which is a $2k$-symmetric space (instead of a symmetric space) (compare for example theorem~\ref{vertical} with theorems~\ref{vertical2} and \ref{verticalj}). This more complicated algebraic structure imposes to define and to use new linear connections $\rabla{j}{}$ when one considers the reductions $\mal Z_{2k}^\alpha (E,J_j) \to \zdk^\alpha (E)$. Moreover, the fact that the fibre is symmetric in the classical twistor bundle $\Sigma(E)$, implies that the different connections defined on this homogeneous fibre bundle can be expressed in a natural way in terms of the metric connection $\nabla$ on $E$. In the case of $\zdk^\alpha (E)$ and $\mal Z_{2k}^\alpha (E,J_j)$,  there are additional complicated terms which appear in these expressions of these different connections in terms of $\nabla$ (due to the structure of $2k$-symmetric space of the fibre).\smallskip\\
Moreover, another novelty is the fact that the fibre $\Sigma(\R^{2n})$ is a complex manifold, even more a Hermitian symmetric space, therefore we are dealing in this case with integrable geometry. In our case we are dealing with non integrable geometry since the structure defined in the fibre $\zdk(\rdn)$ are not parallel with respect to the Levi-Civita connection (there are parallel with respect to the canonical connection which is a metric connection with torsion). This fact will become very clear in the two next sections.\smallskip\\
Finally, a last novelty is our study of the injective morphism of homogeneous fibre bundle $\mijo\colon G/G_0\hookrightarrow \mZ_{2k,2}^{\alpha_0}(G/H,J_2)$, §\ref{geom-interpret-sub1}. We think that it is a nice application of the general theory. Moreover, in our opinion, the fact that the homogeneous fibre bundle structure is preserved (in particular, vertical harmonicity is preserved) is also a nice result. A far as we know, there is no such kind of example - of a non trivial embedding of a homogeneous space into a homogeneous fibre bundle preserving the structure - in the litterature about twistor methods.\smallskip\\
Again, let us make precise that all our results concerning the even minimal determined integrable system are completely new, like all the results concerning the elliptic integrable systems contained in this paper.

%%%%%%%%%%%%%%%%%%%%%%%%%%%%%%%%%%%%%%%%%%%%%%%%%%%%%%%%%%%%%%%%%%%%%%%%%%%%%%%%%%%%%%%%%%%%%%%%%
%                                                                                               %
%                                                                                               %
%                       Generalized harmonic maps                                                    %
%                                                                                               %
%                                                                                               %
%%%%%%%%%%%%%%%%%%%%%%%%%%%%%%%%%%%%%%%%%%%%%%%%%%%%%%%%%%%%%%%%%%%%%%%%%%%%%%%%%%%%%%%%%%%%%%%%%

\section{Generalized harmonic maps}\label{affineharmonic}
\subsection{Affine harmonic maps and holomorphically harmonic maps}
\index{harmonic map|(}
A map $u\colon M \to N$ between two Riemannian manifolds $(M,g)$ and $(N,h)$ is harmonic if it extremizes the energy functional
$$
E(u)=\frac{1}{2}\int_D |du|^2 d\mathrm{vol}_g
$$
for all compact subdomains $D\subset M$, where $|du|^2=\mrm{Tr}_g(u^*h)$. The associated Euler-Lagrange equation is
$\tau(u):=\mrm{Tr}_g(\nabla du)=0$, where $\nabla$ is the connection on $T^*M\otimes u^*TN$ induced by the Levi-Civita connections of $M$ and $N$.\\[1mm]
Now, we generalise this definition for maps from a Riemannian manifold into an affine manifold. We present two different ways to do that. The first one is the natural one (see also  \cite{higaki}) and concerns general affine manifolds whereas the second one concerns maps from Riemann surfaces into affine almost complex manifolds.
\subsubsection{Affine harmonic maps: general properties}\label{aff-harm-gener-prop}
\index{tension field}
\begin{defn}
Let $s\colon (M,g)\to (N,\nabla)$ be a smooth map from a Riemannian manifold $(M,g)$ into an affine manifold $(N,\nabla)$. We set
$$
\tau(s)=\mrm{Tr}_g(\nabla ds)=-\nabla^* ds= *d^\nabla * ds
$$
and we say that $s$ is affine harmonic with respect to $\nabla$ or $\nabla$-harmonic if $\tau(s)=0$.
\end{defn}
Now, let us consider the case where $(M,g)$ is a Riemannian surface  i.e. a Riemann surface $(L,j)$ with a Hermitian metric $g$. Then the action of the Hodge operator $*$ of $L$, is independent of the metric $g$ on 1-forms ($*\alpha=\alpha\circ j$), whereas in 2-forms (resp. 0-forms) it is multiplied by the factor $\lambda^2$ (resp. $\lambda^{-2}>0$) when the metric is multiplied by the factor $\lambda\in C^\infty(L,\R_+^*)$. Hence the tension field $\tau(f)=*d^\nabla(*df)$ is multiplied by $\lambda^2$, under this last transformation. In particular the affine harmonicity for maps $f\colon (L,j) \to (N,\nabla)$ does not depend on the hermitian metric $L$ but only on the conformal structure of $(L,j)$. Thus we have:

\index{strongly harmonic}
\begin{thm}\label{thmstrong}
Let $(L,j)$ be a Riemann surface and $f\colon (L,j)\to (N,\nabla)$ a smooth map. Let $TL^\C=T'L\oplus T''L$ be the decomposition of $TL^\C$ into the $(1,0)$ and $(0,1)$-parts, and $d=\partial + \bar\partial$ and $\nabla^{f^*(TN)}=\nabla' + \nabla''$ the corresponding splittings. Then we have 
$$
2\,\bar\partial^\nabla\partial f = d^\nabla df + i d^\nabla *df,
$$
moreover $d^\nabla df=f^*T$, where $T$ is the torsion of $\nabla$  and $d^\nabla *df=\tau(f) \mrm{vol}_g$ for any hermitian metric $g$ on $L$.
Therefore the following statements are equivalent:
\begin{description}
\item[(i)] $\nabla^{''}\partial f=0$
\item[(ii)] $\bar\partial^\nabla\partial f=0$
\item[(iii)] $\nabla_{\dl{}{\overline z}}\left(\dl{f}{z}\right)=0$, for any holomorphic local coordinate $z=x+iy$ (i.e. $(x,y)$ are conformal coordinates for any hermitian metric in $L$).
\item[(iv)] $f$ is $\nabla$-harmonic with respect to any hermitian metric in $L$ and  torsion free: $f^*T=0$ (i.e. $T(\dl{u}{x},\dl{u}{y})=0$ for any  conformal coordinates $(x,y)$).
\end{description}
We will say in this case that $f$ is strongly $\nabla$-harmonic.
\end{thm}
\begin{rmk}\em
We remark that the imaginary part (resp. the real part) of equation (ii) (resp. equation (iii)) is the $\nabla$-harmonic maps equation whereas its real part (resp. imaginary part) is the torsion free equation $f^*T=0$.\\
If $T=0$ or more generally $f^*T=0$, then $f$ is strongly $\nabla$-harmonic \iif it is $\nabla$-harmonic.
\end{rmk}
\index{harmonic map|)}
\subsubsection{Holomorphically harmonic maps}\label{holharmmap} \index{holomorphically harmonic|(}
In the case the target space $N$ is endowed with an almost complex structure $J$ then we have another way to generalise the definition of harmonicity to  maps from a Riemann surface into $N$.
\begin{defn}\label{holoharmdef}
Let $(L,j)$ be a Riemann surface and $(N,\nabla)$ be an affine manifold endowed with a  complex structure $J$. Let us  denote $TN^\C=T^{1,0}N\oplus T^{0,1}N$ the corresponding decomposition of $TN^\C$. We will say that $f\colon L\to N$ is holomorphically harmonic if 
$$
[\bar\partial^\nabla\partial f]^{1,0}=0.
$$
\end{defn}
\begin{prop}\label{holoharm}
Let $(L,j)$ be a Riemann surface and $(N,\nabla)$ be an affine manifold endowed with a  complex structure $J$.
Then $f$ is holomorphically harmonic \iif (for any hermitian metric $g$ in $M$)
$$
T_g(f) + J\tau_g(f)  =0
$$
where $T_g(f)=*(f^*T)=f^*T(e_1,e_2)$, with $(e_1,e_2)$ an orthonormal basis of $TL$, or equivalently
$$
 \tau_g(f)\mrm{vol}_g =J(f^*T ).
$$
Therefore $f$ is strongly harmonic \iif it is torsion free and holomorphically harmonic. In particular,
if $T=0$, or more generally $f^*T=0$, then $f$ is holomorphically harmonic \iif it is harmonic. Hence for torsion free connection $\nabla$ harmonicity and holomorphic harmonicity are the same.
\end{prop}
\textbf{Proof.}
Let $Z=X+iY\in TN^\C$ with $X,Y\in TN$, then since  $T^{1,0}N$ and $ T^{0,1}N$ are given respectively  by $\{V \mp iJV, V\in TN\}$, we deduce that  
$$
[Z]^{1,0}=0\Leftrightarrow X+JY=0 \quad  \text{and}  \quad [Z]^{0,1}=0\Leftrightarrow X-JY=0.
$$
 Now, let us apply that to the $TN^\C$-valued 2-form $\bar\partial^\nabla\partial f$, we obtain 
$$
[\bar\partial^\nabla\partial f]^{1,0}=0\Longleftrightarrow d^\nabla df + Jd^\nabla*df=0
$$
according to theorem~\ref{thmstrong}. This proves the first assertion. Then the assertion concerning strongly harmonicity follows from theorem~\ref{thmstrong}-(iv). This completes the proof.\hfill $\square$\\[1mm]
Let us remark that
\begin{prop}
In the same situation as in the previous proposition, let us suppose in addition that $\nabla J=0$. Then if a map $f\colon L\to N$ is holomorphic i.e. $df\circ j_L=J df$, then $f$ is anti-holomorphically harmonic (i.e. holomorphically harmonic with respect to $-J$).
\end{prop} 
\textbf{Proof.}
$f$ is holomorphic \iif $df(T^{1,0}L)\subset T^{1,0}N$ i.e. $[\partial f]^{0,1}=0$. Then we have
$$
[\bar\partial^\nabla\partial f]^{0,1}=\bar\partial^\nabla[\partial f]^{0,1}=0
$$
since $\nabla$ commutes with $J$. This completes the proof.\hfill $\square$\medskip\\
It can also be useful to observe the following.
\begin{prop}
Let $(N,J)$ be an almost complex manifold with an almost complex linear connection that we will denote by $\nabla^0$. Then let us define a family of connection 
$$
\nabla^t=\nabla^0 -tT^0, \quad 0\leq t\leq 1.
$$
Then a map $f\colon (L,j_L)\to (N,J)$ from a Riemann surface $L$ into the almost complex manifold $N$ is holomorphically harmonic w.r.t. $\nabla^1$ and $J$ \iif $f$ is holomorphically harmonic w.r.t. $\nabla^0$ and $-J$. We will say more simply that $f$ is $\nabla^1$-holomorphically harmonic \iif it is $\nabla^0$-anti-holomorphicallly harmonic.
\end{prop}
\paragraph{Holomorphic sections of complex vector bundles}
Now, we need to do some recalls about complex vector bundles that we will apply in the next paragraph to obtain an interpretation of the holomorphic harmonicity in terms of holomorphic 1-forms.\\[1mm]
Let $E\to M$ be a real vector bundle (over a manifold $M$) endowed with a complex structure $J\in \End(E)$. Then any frame in the form $(e_x^1,\ldots, e_x^r,Je_x^1,\ldots, Je_x^r)$ at some point $x\in M$ can be extended to a local frame $(e^1,\ldots,e^r,Je^1,\ldots,Je^r)$ in the neighbourhood of $x$. Then there exists an open covering $(U_\alpha)_{\alpha\in I}$ of $M$ and local trivialisations $\Phi_\alpha\colon (E_{|U_\alpha},J)\to U_\alpha\times (\C^r,i\Id)$ which are $\C$-linear isomorphisms ($\Phi_\alpha\circ J=i\Phi_\alpha$), or equivalently of which transition maps take values in the endomorphisms of $\C^r$: $\phi_{\alpha\beta}=\Phi_\beta\circ\Phi_\alpha^{-1}\colon U_\alpha\cap U_\beta \to GL(\C^r)$. Therefore $E$ is a complex vector bundle.%\\[1mm]
\begin{rmk}\em
Let us set $\hat\C=\R[J]$, then $\hat\C=\R[J]$ is a vector bundle over $M$ whose fibres are fields isomorphic to $\C=\R[i]$ and each fibre $E_x$ of $E$ is a $\hat\C_x$-vector space. Then $E^\C$ is endowed  with two different structures of vector bundle: one over the field $\C$ (the tautological one defined by the complexification of $E$) and another one "over the distribution of field $\hat\C$" (i.e. the one defined by $J$). In paticular, we have two different complex structures in $E^\C$.
\end{rmk}
Now, let us suppose that $E$ is endowed with a complex connection $\nabla$, i.e. a connection which commutes with $J$: $\nabla J=0$. Then  for all $X\in TM$, $\nabla_X\colon \mal C(E)\to \mal C(E)$ is $\C$-linear with respect to the complex vector space structure defined on $\mal C(E)$ by the complex vector bundle structure on $E$. Then we have two different ways to extend $\nabla$ to $TM^\C$.
\begin{enumerate}
\item 
 The canonical one: for any section $s\in \mal C(E^\C)$, we extend $\nabla s$ by $\C$-linearity to a linear morphism from $TM^\C$ to $E^\C$,
$$
\nabla_{iX}s=i\nabla_X s,\quad \forall X\in TM,s\in \mal C(E^\C)
$$
after, of course, having extended $\nabla$ to a connection on $E^\C$ by setting $ \nabla is=i\nabla s$, $\forall s\in \mal C(E)$. In conclusion, $\forall s\in \mal C(E^\C)$, $\nabla s\in \mal C(T^*M^\C\otimes E^\C)$. 
\item
By using the complex vector bundle structure of $E$ defined by $J$: for any $s\in \mal C(E)$, we extend $\nabla s$ by $\C$-linearity to a linear morphism from $TM^\C$ to $E$: 
$$
\widehat\nabla_{iX} s= J\widehat\nabla_X s, \quad  \forall X\in TM,s\in \mal C(E).
$$
\end{enumerate}
Let us remark that $\widehat\nabla$ depends on $J$, and since we use  the complex vector space structure defined by $J$, one needs that $\nabla$ and $J$ commute. One the other side the simple canonical complex extention defined in 1 (that we still denote by $\nabla$) is independant of $J$ and one needs not to do any additionnal hypothesis. Remark that the extention 1 is nothing but the extention $\widehat\nabla$ defined by the complex structure $i\Id_{E^\C}$ on $E$ (which commutes obviously with $\nabla$).\\[1mm] 
Now let us suppose that $M$ is an (almost) complex manifold with (almost) complex structure $j_M$. Then we have the splitting $TM^\C=T^{1,0}M\oplus T^{0,1}M$ defined by $j_M$ which gives rise respectively to the following decompositions of $\nabla$ and $\widehat\nabla$:
\begin{eqnarray*}
\widehat\nabla & = & \nabla^{(1,0)} + \nabla^{(0,1)}\\
\nabla & = & \nabla' + \nabla''.
\end{eqnarray*}
\begin{defn}
$\nabla^{(0,1)}$ is called the Cauchy-Riemann operator defined by $\nabla$ and $J$.
\end{defn}
More generally, let $\eta\in\mal C(T^*M\otimes E)$ be a 1-form on $M$ with values in $E$. Then we can extend it in two different ways by $\C$-linearity in $TM^\C$ by setting:
\begin{eqnarray*}
\eta^\C(X + iY) & = & \eta(X) + i\eta(Y),\quad \forall X,Y\in TM\\
\hat\eta(X + iY) & = & \eta(X) + J\eta(Y),\quad \forall X,Y\in TM.
\end{eqnarray*}
Remark that $\eta^\C\in \mal C(T^*M^\C\otimes E^\C)$ whereas  $\hat\eta\in\mal C(T^*M^\C\otimes E)$. As above we can decompose $\eta^\C$ and $\hat\eta$ according to the decomposition $TM^\C=T^{1,0}M\oplus T^{0,1}M$:
\begin{eqnarray}
\eta^\C & = &  \eta' + \eta'' \label{splitting-eta1}\\
\hat\eta & = & \eta^{(1,0)} + \eta^{(0,1)}\label{splitting-eta2}.
\end{eqnarray}
Then we have the following relations
\begin{lemma}
\begin{equation}\label{eta}
\begin{array}{cc}
\left[\eta'\right]^{1,0}=\eta^{(1,0)} -iJ\eta^{(1,0)} &  \left[\eta''\right]^{0,1}= \eta^{(1,0)} +iJ\eta^{(1,0)}\\
\left[\eta'\right]^{0,1}= \eta^{(0,1)} +iJ\eta^{(0,1)} &  \left[\eta''\right]^{1,0}= \eta^{(0,1)} -iJ\eta^{(0,1)}
\end{array}
\end{equation}
\end{lemma}
\textbf{Proof.}
Let $Z=X-ij_MX\in T^{1,0}M$ with $X\in TM$. Then 
\begin{eqnarray*}
\left[ \eta(Z)\right]^{1,0}=\left[\eta(X)-i\eta(j_MX)\right]^{1,0} & = & \eta(X)-iJ\eta(X) - i\left( \eta(j_MX)-iJ\eta(j_MX)\right)\\
& = & \eta(X)-J\eta(j_MX) -i J\left( \eta(X)-J\eta(j_MX)\right) \\
& = & \eta^{(1,0)}(Z) -iJ\eta^{(1,0)}(Z).
\end{eqnarray*}
This gives us $\left[\eta'\right]^{1,0}$. Then by apllying this to  $-J$, we obtain $\left[\eta'\right]^{0,1}$. Finally, the second column of (\ref{eta}) is obtained by $\C$-conjugaison from the first column. This completes the proof.\hfill $\square$\\[1mm] 
We can apply what precedes to the flat differentiation $d$. Let $(N,J)$ be an almost complex manifold and $s\colon M \to N$ a map. Then we consider the complex vector bundle $E=s^*TN$ over $M$. Then applying what precedes to the 1-form $\eta=ds$, we can consider the extensions $\widehat{ds}$ and $(ds)^\C$, which then allows us to define the following extension of $d$ to $TM^\C$: 
$$
\hat d s= \widehat{ds} \quad \text{and}\quad  d^\C s=(ds)^\C,   
$$
and \textbf{by abuse of notation\footnote{and to be coherent with the notation used until now, in the paper.} $d^\C$ will be still denoted by $d$}.
Then we can write the following decompositions
$$
\hat d =\hat\partial + \overline{\hat\partial} \quad \text{and}\quad  d= \partial + \bar\partial
$$
according to the decomposition $TM^\C=T^{1,0}M\oplus T^{0,1}M$.%\\[1mm]
\begin{defn}
$\overline{\hat\partial}$ is called the Cauchy-Riemann derivative defined by  $J$.
\end{defn}
Now let us come back to the general situation of a complex vector bundle $E$ over an almost complex manifold $(M,j_M)$,  endowed with a complex  connection $\nabla$. Let us set 
$$
\mH(M,E)=\{\eta\in T^*M\otimes E|  \eta\,j_M=J\eta \}. 
$$
Then $\mH(M,E)$ is a vector subbundle of the vector bundle $T^*M\otimes E$ and is naturally endowed with the complex structure defined by
\begin{equation}\label{def-I}
I(\eta)= \eta\, j_M=J \eta , \quad \forall\eta\in T^*M\otimes E,
\end{equation}
which makes  $\mH(M,E)$ being a complex vector bundle whose the set of sections is
$$
\mrm{Hom}((TM,j_M),(E,J))=\{\eta\in\mal C(T^*M\otimes E)| \eta\circ j_M=J\circ\eta \}.
$$
The sections of $\mH(M,E)$ can also be characterized by using the splittings (\ref{splitting-eta1}-\ref{splitting-eta2}):
\begin{lemma}
We have the following equivalences for 1-forms $\eta\in\mal C(T^*M\otimes E)$:
$$
\eta\circ j_M=J\circ\eta\ \Longleftrightarrow\ \eta'\in T^*M\otimes E^{1,0}\ \Longleftrightarrow\  \eta^{(0,1)}=0.
$$
\end{lemma}
Then we deduce in particular
\begin{lemma}\label{lemmamorphism}
Let $s\in \mal C(E)$, then we have the following equivalences:
$$
\nabla s\circ j_M =J\circ\nabla s \Longleftrightarrow \nabla^{(0,1)} s= 0 \Longleftrightarrow \left[ \nabla' s\right] ^{0,1}=0 \Longleftrightarrow \nabla'(s-iJs)=0.
$$
We will say that $s$ is a vertically holomorphic section.
\end{lemma}
In fact we can say more
\begin{lemma}\label{check-J}
Let us consider the splitting $TE=\hor\oplus\ver$ given by $\nabla$, where $\ver=\ker\pi=\pi^*E$ is the vertical subbundle and $\hor$ the horizontal one. Then let us define an  almost complex  structure $\check{J}$ on the manifold $E$  by setting
$$
\check{J}=\left((d\pi)^* j_M\right)_{|\hor} \oplus \pi^*J .
$$
Then a section $s\in \mal C(E)$ is $\check J$-holomorphic \iif it is vertically holomorphic. 
%Moreover $\check J$ is the unique complex structure on the manifold $E$ whose the holomorphic section are those %which are vertically holomorphic.
\end{lemma}
\textbf{Proof.}
It suffices to prove that any section $s\in\mal C(E)$ is horizontally holomorphic, i.e. satisfies the horizontal part of the equation $ds\circ j_M=J\circ ds$. We have $d\pi\circ (ds\circ j_M)=j_M$ since $s$ is a section.  In the other side we have $d\pi\circ(\check  J\circ ds)=j_M\circ d\pi\circ ds=j_M$, by definition of $\check{J}$ and using the fact $s$ is a section. In conclusion $d\pi\circ(ds\circ j_M)=d\pi\circ(J\circ ds)$. This completes the proof.\hfill $\square$\medskip\\
In the following, we will say that \emph{a section of a complex vector bundle $(E,J,\nabla)$ is holomorphic if it is $\check J$-holomorphic.}\medskip\\
Now, let us apply the two previous lemmas to the vector bundle $\mH(M,E)$. First, let us endow $M$ with an almost complex connection $\nabla^M$ (it means $\nabla^MJ=0$; such a connection always exists, see \cite{KN}). Then $T^*M\otimes E$ is naturally endowed with the connection  $\onabla{\otimes}{}{}$ defined by $\nabla^M$ and $\nabla$. Further, we denote by $\overline\nabla$ the restriction to $\mH(M,E)$ of $\onabla{\otimes}{}{}$. Then we remark that 
$\overline\nabla$ commutes with the complex structure $I$ (defined by (\ref{def-I})). Therefore, we can now apply the two previous lemmas to the complex vector bundle $ \left( \mH(M,E),I, \overline\nabla\right) $:
\begin{prop}\label{holomorphicsectionprop}
A section of $\mH(M,E)$, $\eta\in\mrm{Hom}((TM,j_M),(E,J))$, is holomorphic \iif it satisfies one of the following equivalent statements
\begin{description}
\item[(i)] $\overline\nabla^{(0,1)} \eta=0$
\item[(ii)] $\overline\nabla''\eta'=0$
\item[(iii)] $\left[ \overline\nabla''\eta \right] ^{1,0}=0$.
\end{description}
Moreover if $M$ is a Riemann surface\footnote{and $\nabla^M$ the unique torsion free complex connection in $M$, which coincides also with the Levi-Civita connection of any Hermitian metric on $M$.}, then it is also equivalent to 
\begin{description}
\item[(iv)]  $\overline{\hat\partial}^{\widehat\nabla}\hat\eta = \overline{\hat\partial}^{\widehat\nabla}\eta^{(1,0)}=0$, or
\item[(v)] $ \bar\partial^{\nabla}\eta'=0$
\end{description}
\end{prop}
Moreover, if $M$ is a complex manifold (i.e. $j_M$ is integrable) then we  choose for $\nabla^M$ the unique torsion free complex connection on $M$. Then we obtain the following result:
\begin{prop}\label{d-nabla-eta}
Let $l\in TM$ be a complex line in the tangent bundle of the complex manifold $M$. Then  for any section $\eta\in\mrm{Hom}((TM,j_M),(E,J))$ we have the following equality
$$
\overline\nabla^{(0,1)}\eta_{|l\times l}=d^\nabla\eta_{|l\times l}.
$$
Moreover if  $\eta$ is holomorphic then $d^\nabla\eta=0$. More particulary, if $M$ is a Riemann surface then we have the following equivalence
$$
\eta \text{ is holomorphic }\Longleftrightarrow d^\nabla\eta=0.
$$
\end{prop}
\begin{rmk}\em 
One could directly deduces the case of  a Riemann surface by using proposition~\ref{holomorphicsectionprop}.\\ Indeed, the first way to do that is to write $d^\nabla\eta=d^\nabla\eta'+ d^\nabla\eta''$. Then remark that $\eta'$ and $\eta''$ takes values in $E^{1,0}$ and $E^{0,1}$ respectively, according to lemma~\ref{lemmamorphism}. Therefore since $E^{1,0}$ and $E^{0,1}$ are $\nabla$-parallel, we can say that $d^\nabla\eta'$ and $d^\nabla\eta''$ take values resp. in $E^{1,0}$ and $E^{0,1}$ resp., so that $d^\nabla\eta=0\Leftrightarrow d^\nabla\eta'=0\Leftrightarrow d^\nabla\eta''=0$. Then if $M$ is a Riemann surface $d^\nabla\eta'=\bar\partial^\nabla\eta'$, and we conclude by using proposition~\ref{holomorphicsectionprop}.\\ 
The second way to do that is to use the $\hat\C$-linearity. Indeed, the extension to $TM^\C$ by $\hat\C$-linearity of
$d^\nabla\eta$ is $\widehat{d^\nabla\eta}=d^{\widehat\nabla}\hat\eta= d^{\widehat\nabla}\eta^{(1,0)}=0$, since $\eta^{(0,1)}=0$ (see lemma~\ref{lemmamorphism}). Then  if $M$ is a Riemann surface $d^{\widehat\nabla}\eta^{(1,0)}=
\overline{\hat\partial}^{\widehat\nabla}\eta^{(1,0)}$, and we conclude by using proposition~\ref{holomorphicsectionprop}-(iv).
\end{rmk}
\begin{rmk}\label{eta-beta} \em
Let us consider a 1-form $\beta\in\mal C(T^*M\otimes E)$, then we can associate to it 
$$
\eta=\beta -J\beta\circ j_M=\hat\beta\circ(\Id -i j_M)=\beta^{(1,0)}\circ(\Id - ij_M).
$$
By definition $\eta\in\mal C(\mH(M,E))$, i.e. $\eta\circ j_M=J\circ\eta$. Moreover, still suppposing that $M$ is complex and that $\nabla^M$ is the unique torsion free complex connection on $M$, we have
\begin{equation}\label{id-ijm}
\widehat{\overline\nabla}\hat\eta =\widehat{\overline\nabla}\beta^{(1,0)}\circ (\Id-ij_M)
\end{equation}
because $(\Id-ij_M)$ is $\nabla^M$-parallel. \\
Let us remark that since $j_M$ and the multiplication by $i$ coincide in $T^{1,0}M$, they define the same complex structure, which we will suppose $T^{1,0}M$ to be canonically endowed with. Then, since $\hat\beta$ is by definition a complex linear morphism from $TM^\C$ to $E$, $\beta^{(1,0)}$ is also a complex linear morphism from $T^{1,0}M$ to $E$. Hence $\beta^{(1,0)}$ is  a section of the complex vector bundle $T_{1,0}^*M\otimes_\C E$.
Therefore, from equation (\ref{id-ijm}), we deduce that $\eta$ is a holomorphic section of $\mH(M,E)$ \iif $\beta^{1,0}$ is a holomorphic section of $T_{1,0}^*M\otimes_\C E$. In particular if $M$ is a Riemann surface, we deduce from (\ref{id-ijm}), that
$$
\overline{\hat\partial}^{\widehat\nabla}\hat\eta=\overline{\hat\partial}^{\widehat\nabla}\beta^{(1,0)}\circ(\Id-ij_M).
$$
\end{rmk}
Now, we come back to our complex vector bundle $(E,J,\nabla)$ and  we recall a theorem (\cite{kobayashi}) which characterizes  when $\check J$ is integrable.
\begin{thm}
Let $(E,J,\nabla)\to (M,j_M)$ be a complex vector bundle over a complex manifold, with a complex connection $\nabla$. Then 
we will say that a holomorphic structure $\mal E$ is compatible with $\nabla$ (or that $\nabla$ is adapted to $\mal E$) if it is induced by the almost complex structure $\check J$ (defined by lemma~\ref{check-J}). In other words, a section $s\in \mal C(E)$ is holomorphic with respect to $\mal E$ \iif 
$$
\forall Z\in T^{1,0}M,\ \widehat\nabla_{\bar Z} s=0.
$$
An holomorphic structure $\mal E$ exists on $E$ \iif $\check J$ is integrable, and in this case $\mal E$ is unique.
Moreover $\check J$ is integrable \iif the $(0,2)$-component  of the curvature operator\footnote{i.e. the $(0,2)$-component of the extension $\hat R$ of $R$ to $\Lm^2 T^*\!M^\C$ by $\hat\C$-linearity.} $R$ of $\nabla$ vanishes.
\end{thm}
When $M$ is of dimension 2, then the $(0,2)$-component  of the curvature operator always vanishes so that $E$ always admits a holomorphic structure compatible with $\nabla$, that we will call, following \cite{BuRaw}, the Koszul-Malgrange holomorphic structure induced by $\nabla$. In the following, we suppose that  a complex vector bundle $(E,J,\nabla)$ over a Riemann surface is always endowed with its Koszul-Malgrange holomorphic structure.
\paragraph{Interpretation of the holomorphic harmonicity in terms of holomorphic 1-forms.}
\index{dd derivative@$\bar\partial\partial$-derivative}
Now we come back to the situation in the begining of \ref{holharmmap}. More precisely, we consider $(N,J)$ an almost complex manifold, with $\nabla$ an almost complex connection, $(L,j_L)$ a Riemann surface and $f\colon L \to N$ a map. Then we apply what precedes to the complex vector bundle $E=(f^*TN,f^*\nabla,f^*J)$ over $L$ (i.e. $L$ plays the role of $M$ and $f$ the one of $s$ with respect to the notation of the previous paragraphs). We obtain  a first theorem:
\begin{prop}\label{hat-df}
Let $f\colon (L,j_L) \to (N,J,\nabla)$ be a map from a Riemann surface into an almost complex affine manifold. Let us set 
$$
\eta=df -Jdf\circ j_M.
$$
Then $\eta$ is a section of $\mH(L,f^*TN)$, i.e. $\eta\circ j_M= J\circ\eta$. Moreover $f$ is holomorphically harmonic \iif $\eta$ is a holomorphic section of the complex vector bundle $\mH(L,f^*TN)$, i.e.
$$
\bar\partial^\nabla\eta'=0.
$$ 
\end{prop}
\textbf{Proof.}
We write 
$$
d^\nabla \eta=d^\nabla \left( df -Jdf\circ j_M\right) )=d^\nabla\left( df +*Jdf\right) = d^\nabla df + Jd^\nabla*df
$$
so that we can conclude according to proposition~\ref{d-nabla-eta} and proposition~\ref{holoharm}. This completes the proof.\hfill$\square$\medskip\\
Now, we can give a characterization of holomorphic harmonicity which looks like very closely to the one which holds for harmonic maps (\cite{BuRaw}):
\begin{thm}\label{thm-dbdf=0}
A map $f\colon (L,j_L) \to (N,J,\nabla)$  from a Riemann surface into an almost complex affine manifold, is holomorphically harmonic \iif 
\begin{equation}\label{C-laplacian-f=0}
\overline{\hat\partial}^{\widehat\nabla}\hat\partial f =0,
\end{equation}
i.e. $\hat\partial f$ is a holomorphic section of $T_{1,0}^*L\otimes_\C f^*TN$.
\end{thm}
\textbf{Proof.}
Apply remark~\ref{eta-beta} to $\beta=df$ and then use proposition~\ref{hat-df} to prove that $\hat\partial f$ is a holomorphic section and proposition~\ref{holomorphicsectionprop}-(iv) to prove the equation~(\ref{C-laplacian-f=0}).\hfill$\square$
\begin{rmk}\label{rmk-comp-dbdf} \em
Let us derive \eqref{C-laplacian-f=0} by a direct computation. Let be $Z=X + ij_L X\in T^{1,0} L$, then we have
\begin{eqnarray*}
\left( f^*\nabla\right)_{\bar Z}^{0,1} \hat\partial f(Z) &  = & \left( \nabla_{f(X)} + J \nabla_{f(j_L X)}\right) ( df(X) - Jdf(j_L X)) \\
 & = & \nabla_{f(X)} df(X) + \nabla_{f(j_L X)} df(j_L X) + J \left( \nabla_{f(j_L X)} df(X) - \nabla_{f(X)}df(j_L X)\right)\\
 &  =  &  d^\nabla * df (X,j_L X)  - J d^\nabla df (X,j_L X).
\end{eqnarray*}
\end{rmk}
\index{holomorphically harmonic|)}
\subsection{The sigma model with a Wess-Zumino term in Nearly K\"{a}hler manifolds}\index{Wess-Zumino|(}
Here we present an interpretation of the holomorphic harmonicity in terms of  a sigma model with a Wess-Zumino term. 
\subsubsection{Totally skew-symmetric torsion}\label{Total-Skew-Sym-tors}
\index{skew-symmetric torsion|(}
First, let us recall some useful properties about connections \footnote{See also \cite{srni} for a nice presentation about metric connections and their torsion}.
\begin{defn}
Two linear connection $\nabla$ and $\nabla'$ in a manifold $N$ are sayed to be geodesically equivalent if they have the same geodesics. A connection $\nabla$ on a Riemannian manifold $(N,h)$ is sayed to be geodesic-preserving if it is geodesically equivalent to the Levi-Civita connection $\nabla^h$ of $h$.
\end{defn}
\begin{prop}\label{prop-gene-A} Let $\nabla$ be a connection on a manifold $N$ and  $A\in \mal C(T^*N\otimes \End(TN))$. Then the connection 
$$
\nabla'=\nabla + A
$$ 
has the same geodesic as $\nabla$ \iif $A(\cdot, \cdot)$ is  skew-symmetric (as a bilinear map). In this case for any map $f\colon (M,g)\to N$, from a Riemannian manifold in to $N$, we have $\tau_g'(f)=\tau_g(f)$, where $\tau_g'(f)$ and $\tau_g(f)$ are the tension fields \wrt $\nabla$ and $\nabla'$ respectively. Moreover (still in this case), we have 
$$
T^{\nabla'}= T^\nabla + 2A.
$$
Now, let us suppose that $\nabla$ is metric \wrt some metric $h$ in $N$. Then $\nabla'$ is metric \iif $A$  takes values (as a 1-form) in the skew-symmetric endomorphisms of $TN$: $A\in \mal C(T^*N\otimes \so(TN))$. Therefore $\nabla'$ is metric and geodesically equivalent to $\nabla$ \iif $A$ is \emph{totally skew-symmetric} which means that the associated 3-linear map defined by $A^*(X,Y,Z)= \langle A(X,Y),Z\rangle$ is a 3-form on $N$.
\end{prop}
Now let us see how we can introduce the Levi-Civita connection starting from a given metric connection.
\begin{prop}\label{T-total-antisym} 
Let $N$ be a manifold endowed with some connection that we denote by $\nabla^0$ (for some reason that will appear clearly below). Let us set 
$$
\nabla^t=\nabla^0 -t T^0,\quad 0\leq t\leq 1,
$$
where $T^0=T^{\nabla^0}$ is the torsion of $\nabla^0$. Then we have
$$
T^t:=T^{\nabla^t}= -(2t-1) T^0.
$$
In particular, $\nabla^{\frac{1}{2}}$ is torsion free. Moreover all the connections $\nabla^t,\quad 0\leq t\leq 1$, are geodesically equivalent. \\
Now, let $h$ be a metric on $N$ which is $\nabla^0$-parallel. Then $\nabla^t$, $t\neq 0$, is metric \iif $T^0$ is totally skew-symmetric that is to say the 3-linear map defined by 
$$
(T^0)^*(X,Y,Z):=\langle T^0(X,Y),Z\rangle
$$
is a 3-form. In this case, $\nabla^{\frac{1}{2}}$ coincides with the Levi-Civita connection $\nabla^h$ of $h$.
\end{prop}
\begin{rmk}\label{rmk-strongharmequiv}\index{strongly harmonic}
We see that for a map $f\colon (M,g)\to N$, the strongly $\nabla^t$-harmonicities are all equivalent for $t\neq \frac{1}{2}$.
\end{rmk}
Conversely,
\begin{prop}\label{prop-5-9}
 Let $(N,h)$ be a Riemannian manifold, and let us denote by $\nabla^h$ its Levi-Civita connection. Then a metric connection $\nabla$ on $N$ is entirely determined  by its torsion $T$.  Moreover a metric connection  $\nabla$ on $N$ is  geodesic-preserving \iif its torsion $T$ is totally skew-symmetric. Then in this case we have
$$
\nabla = \nabla^h + \dfrac{1}{2}T.
$$
\end{prop}
\textbf{Proof.} For any metric connection $\nabla=\nabla^h + A$, we have 
\begin{eqnarray}
T(X,Y) & = & A(X,Y) -A(Y,X)\label{T-A}\\
2A^*(X,Y,Z) & = & T^*(X,Y,Z) + T^*(Z,X,Y) + T^*(Z,Y,X)\label{A-T},
\end{eqnarray} 
which proves the first assertion. 
Concerning the second assertion, we see (according to (\ref{T-A}-\ref{A-T})) that $A$ is totally skew-symmetric \iif $T$ is so, i.e., according to proposition~\ref{prop-gene-A},  $\nabla$  is geodesic preserving  \iif $T$ is totally skew-symmetric. Then in this case $T=2A$ i.e. $\nabla = \nabla^h + \dfrac{1}{2}T$. This completes the proof.\hfill$\square$
\begin{rmk}\label{T+U}\em
The second equation (\ref{A-T}) can be derived directly  from the first one (\ref{T-A}) (compute the right hand side of the second equation using the first equation which gives  $2A^*(X,Y,Z)$). But there is another way  (which will be useful in the following) to interpret this second equation. Indeed, first let us identify (via the metric $h$) in the following of this remark, each $TN$-valued bilinear form $B$ on $N$ with the corresponding trilinear form $B^*$. Then let us set $A := \frac{1}{2}( T + U)= \frac{1}{2}( T(X,Y,Z) + T(Z,X,Y) + T(Z,Y,X))$, where $ U(X,Y,Z)=\langle U(X,Y),Z\rangle =T(Z,X,Y) + T(Z,Y,X)$. We remark that $U$ is symmetric \wrt to the  variables $X,Y$, so that  the connection $\nabla  - A=\nabla -\frac{1}{2}(T+U)$  is torsion free. Moreover we see that $A(X,Y,Z)=\frac{1}{2}( T(X,Y,Z) + T(Z,X,Y) + T(Z,Y,X))$ is skew symmetric \wrt the two last variables $Y,Z$. Therefore $\nabla  - A$ is metric and thus this is the Levi-Civita connection $\nabla^h$: 
$$
\nabla^h=\nabla -\frac{1}{2}(T+U).
$$
Moreover, $T$ is totally skew-symmetric \iif the "natural reductivity term" $U=0$.\\
Furthermore, let us remark that the bijective correspondence between $T$ and $A$ comes simply from the isomorphism of vector bundle $T\in\Lm^2T^*N\otimes TN\longmapsto T + U \in T^*N\otimes\so(TN)$.
\end{rmk}
\subsubsection{The general case of an almost Hermitian manifold}\label{WZ}
\index{holomorphically harmonic|(}
Let $(E,J)$ be a complex vector space and let us set 
$$
\mrm{Bil}(E)=E^*\otimes E^*\otimes E \quad \text{and}\  \mal T (E)= (\Lm^2 E^*)\otimes E \subset \mrm{Bil}(E).
$$
and for $\eps,\eps'\in \Z_2$ we set 
$$
\mrm{Bil}^{\eps,\eps'}(E,J)=\{A\in\mrm{Bil}(E)| A(J\cdot,\cdot) = \eps JA,\  A(\cdot,J\cdot)=\eps'JA \}
$$
so that we have the decomposition 
\begin{equation}\label{dec-Bil}
\mrm{Bil}(E)=\oplus_{(\eps,\eps')\in \Z_2\times\Z_2}\mrm{Bil}^{\eps,\eps'}(E,J).
\end{equation}
Let us remark that for any $A\in\mrm{Bil}(E)$, its component $A^{\eps,\eps'}\in\mrm{Bil}^{\eps,\eps'}(E,J)$ is given by
\begin{equation}\label{component-eps,eps'}
A^{\eps,\eps'}(X,Y)=-\dfrac{1}{4} \left(\eps\eps' A(JX,JY) +\eps JA(JX,Y) + \eps' JA(X,JY) -A(X,Y)\right).
\end{equation}
Moreover we also have the decomposition 
$$
\mal T(E)= \mal T^{2,0} \oplus\mal T^{0,2} \oplus \mal T^{1,1},
$$
where $\mal T^{2,0}=(\Lm^{2,0}{E^*}^\C)\otimes_\C E=\mrm{Bil}^{++}(E,J)\cap \mal T(E)$, $\mal T^{2,0}=(\Lm^{0,2}{E^*}^\C)\otimes_\C E=\mrm{Bil}^{--}(E,J)\cap \mal T(E)$  and  $\mal T^{1,1}=(\Lm^{1,1}{E^*}^\C)\otimes_\C E=\left( \mrm{Bil}^{+-} + \mrm{Bil}^{-+}\right) (E,J)\cap \mal T(E)$.\\
Of course, these notation can be extented to the case $(E,J)$ is a complex vector bundle. In particular, we will use these for the tangent bundle $(TN,J)$ of an almost complex manifold, and will forget in this case the precision of the bundle in the notation and write for example simply $\mal T$ and $\mrm{Bil}$.\medskip\\
Given  an almost complex manifold $(N,J)$ with a connection $\nabla$ (that we do not suppose to be almost complex) then its torsion $T$ satisfies $T\in \mal T$ and\footnote{i.e. $T\in \mal C(\mal T)$} we can decompose it following (\ref{dec-Bil}): $T=T^{++} + T^{--} + T^{-+} + T^{+-}$. Then since $T$ is skew-symmetric, then so is $T^{++}$, $T^{--}$ and $T^{+-} + T^{-+}$. In other words, we have $T^{++}=T^{2,0}$, $T^{--}=T^{0,2}$ and $T^{+-} + T^{-+}=T^{1,1}$. In particular we have $T^{+-}(X,Y)=-T^{-+}(Y,X)$.
Now, let us see how this decomposition can have a geometric meaning.
\begin{lemma}\index{Nijenhuis tensor!|(}
Let $(N,J,\nabla)$ be an almost complex manifold with an almost complex connection $\nabla$.  Then  we have 
$$
N_J =4 T^{--}
$$ 
where $N_J$ denotes the torsion of $J$ i.e its Nijenhuis  tensor. 
\end{lemma}
\textbf{Proof.} According to \cite{KN}, Chap. IX, Prop. 3.6, the torsion  $N_J$ of $J$ can be expressed in terms of the torsion $T$ of the almost complex connection $\nabla$:
$$
-N_J(X,Y)=T(JX,JY)-JT(JX,Y) -JT(X,JY) -T(X,Y)
$$
which gives us $N_J=4T^{--}$. This completes the proof.\hfill$\square$\\
\begin{prop}\label{JanticomT}
 Let $(N,J,\nabla)$ be an almost complex manifold with an almost complex connection $\nabla$.  Then the following statements are equivalent. 
\begin{description}
\item[(i)] $J$ \emph{anticommutes with the torsion} $T$ of $\nabla$ : $T(X,JY)=-JT(X,Y)$.
\item[(ii)] $T=T^{--}$  i.e. $T\in \mal T^{2,0}$.
\item[(iii)] $T=\dfrac{1}{4}N_J$.
\end{description}
\end{prop}
\textbf{Proof.} (ii) $\Leftrightarrow$ (iii) follows from the previous lemma. Now, we have obviously (ii) $\Rightarrow$ (i). Conversely (i) implies that $T=T^{--} + T^{+-}$ but since $T$ is skew-symmetric this implies $T^{+-}=0$ and $T=T^{--}$. This completes the proof.\medskip\hfill$\square$\\
From now until the end of this section~\ref{WZ}, we consider $(N,J)$  an almost complex manifold with an almost complex linear connection $\nabla$ and a $\nabla$-parallel Hermitian metric $h$. Therefore $(N,J,h)$ is an almost Hermitian manifold with a Hermitian connection $\nabla$.
\begin{prop}\label{prop-H-closed}
Let $(N,J,h)$ be an almost Hermitian manifold with a Hermitian connection $\nabla$. Let us suppose that $J$ anticommutes with the torsion $T$ of $\nabla$. Let us suppose also that the torsion of $\nabla$ is totally skew-symmetric i.e. 
$$
T^*(X,Y,Z)=\langle T(X,Y),Z\rangle
$$
is a 3-form. Lastly, we suppose that the torsion is $\nabla$-parallel, i.e. $\nabla T^*=0$ which is equivalent to $\nabla T=0$.  Then  the trilinear map 
$$
H(X,Y,Z)=-T^*(X,Y,JZ)=\langle JT(X,Y),Z\rangle
$$ 
is 3-form  and is closed $dH=0$.
\end{prop}
\textbf{Proof.} Firstly, according to proposition~\ref{JanticomT}, we have $T^*(JX,Y,Z)=T^*(X,JY,Z)=T^*(X,Y,JZ)$, which prove that $H$ is a 3-form.\\
Let us compute the exterior differential of $H$ in terms of the connection (with torsion) $\nabla$:
\begin{eqnarray*}
dH(X_0,X_1,X_2,X_3) &  =  & \sum_{i=0}^3  (-1)^{i}\nabla_{X_i}^0 H(X_0,\ldots, \hat X_i,\ldots, X_3)\\
  &   &  - \sum_{0\leq i<j\leq 3} (-1)^{i+j} H(T(X_i,X_j), X_0,\ldots,\hat X_i,\ldots, \hat X_j,\ldots, X_3))\\
 & = & \underset{i,j,k}{\mal S} H(T(X_0,X_i), X_j, X_k) + H(T(X_i,X_j), X_0, X_k)
\end{eqnarray*}
where the last sum is on all the circular permutations of 1,2,3. Moreover, we have 
$$
\begin{array}{rcl}
H(T(X_0,X_i), X_j, X_k) & = & -H(X_k, X_j,T(X_0,X_i))\\
                        & = & \langle T(X_k, X_j),JT(X_0,X_i))\rangle\\
                        & = &  \langle T(X_0,X_i), J  T(X_j, X_k) \rangle\\
                        & = &  -H( X_0,X_i,   T(X_j, X_k))\\ 
                        & = &  -H( T(X_j, X_k),X_0, X_i)\\ 
\end{array}
$$
so that we can conclude that $dH=0$. This completes the proof.\hfill$\square$\\ 
\begin{thm}\label{WZW}
Let $(N,J,h)$ be an almost Hermitian manifold with a Hermitian connection $\nabla$. Then, under the 3 hyphothesis of the previous proposition ($T$ anticommutes with $J$, is totally skew-symmetric and $\nabla$-parallel), the equation for holomorphically harmonic maps $f\colon L\to N$ is the  equation of motion (i.e. the Euler-Lagrange equation) for the sigma model in $N$ with the Wess-Zumino term defined by the closed 3-form $H$.
The action functional is given by
$$
S(f)=E(f)  + S^{WZ}(f)= \dfrac{1}{2}\int_L|df|^2 d\mathrm{vol}_g + \int_B H,
$$
where $B$ is 3-submanifold (or indeed a 3-chain) in $N$ whose boundary is $f(L)$.
\end{thm}
\textbf{Proof.}  Since $dH=0$ we have
$$
\delta S^{WZ}=\int_B L_{\delta f}H=\int_Bd\imath_{\delta f}H=\int_{f(L)}\imath_{\delta f}H,
$$
therefore the Euler-Lagrange equation is
$$
-\tau_g(f) + JT_g(f)=0
$$
which is the equation for holomorphically harmonic maps \wrt $\nabla$, since $\nabla$ is geodesic preserving (see propositions~\ref{prop-gene-A} and \ref{prop-5-9}) ($g$ being as always a Hermitian metric on $L$). This completes the proof.\hfill$\square$\\ 
\index{Nijenhuis tensor!|)}
\subsubsection{The example of a 3-symmetric space}
\index{canonical!almost complex structure|(} \index{canonical!connection, $G$-invariant|(}\index{model@model case or system|(}
Let us suppose now that $N=G/G_0$ is a (locally) 3-symmetric space. We use the notations of subsection~\ref{2.1.2}. In particular, $N$ is endowed with its canonical almost complex structure $\undj$ defined by (\ref{def-odd}).
\begin{prop}\label{3-sym-prop}
The canonical connection $\nabla^0$ in $N$ commutes with the canonical almost complex structure $\undj$
$$
\nabla^0 \undj=0.
$$
Moreover,  $\undj$ anticommutes with the torsion $T^0$ of $\nabla^0$. Lastly, if $N$ is Riemannian, then 
$\nabla^0$ is metric and $(N,J,h)$ is almost Hermitian for any $G$-invariant metric $h$.\footnote{chosen according to our convention explained in subsection~\ref{2.1.2}: that is $\tau_\mk$ leaves invariant the inner product defining $h$.}\\
Furthermore, the torsion of $\nabla^0$ is totally skew-symmetric \iif $h$ is naturally reductive.
\end{prop}
Now, we can conclude 
\begin{thm}\label{WZ-3sym}
Let $N=G/G_0$ be a  (locally) 3-symmetric space, that we suppose to be Riemannian and naturally reductive, and endowed with its canonical almost complex structure $\undj$ and its canonical connection $\nabla^0$. Let $h$ be a $G$-invariant naturally reductive metric on $N$. Then the equation for holomorphically harmonic maps $f\colon L\to N$ is the   Euler-Lagrange equation for the sigma model in $N$ with the Wess-Zumino term defined by the closed 3-form $H$ defined by 
$$
H(X,Y,Z)=\langle \undj T(X,Y),Z\rangle
$$
where $T$ is the torsion of $\nabla^0$.

\end{thm}
\begin{rmk}\em
The hypothesis of natural reductivity is always satisfied if we allow us to use pseudo-Riemannian metrics and if $\g$ is semi-simple: the metric defined by the Killing form is then naturally reductive. Moreover, let us remark that \wrt  Riemannian metrics the natural reductivity is in fact an hypothesis of compactness: the Lie subgroup of $GL(\mk)$ generated by $\{[\ad_\mk(X)]_\mk, X\in\mk\}$ must be compact.
\end{rmk}
\index{canonical!almost complex structure|)} \index{canonical!connection, $G$-invariant|)}\index{model@model case or system|)}
\subsubsection{The good geometric context/setting}\label{good-geom-set-kahler}
In the previous variational interpretation given by theorem~\ref{WZW}, we need to make 3 hypothesis on the torsion of the almost Hermitian connection: $T$ anticommutes with $J$, is totally skew-symmetric and $\nabla$-parallel. Here, we want to understand what do these hypothesis mean geometrically and  what is the good geometric context in which these take place. It will turn out that the good geometric context is the one of \emph{Nearly K\"{a}hler manifold}.
\begin{defn}
An almost Hermitian manifold $(N,h,J)$ is called nearly K\"{a}hler if 
$$
(\nabla^h_X J)X=0,
$$
where $\nabla^h$ is the Levi-Civita connection of $h$.
\end{defn}
We can deduce immediately the following properties.
\begin{prop}\label{prop-equiv-nearly}\index{canonical!Hermitian connection|(}
Let $(N,h,J)$ be an almost Hermitian manifold. Let us consider its canonical  Hermitian connection
$$
\nabla^0:=\nabla^h -\dfrac{1}{2}J\nabla^h J,
$$
the torsion of which is denoted by $T^0$.
Then the following statements are equivalent:
\begin{description}
\item[(i)] $T^0(\cdot,J\cdot)=-JT^0(\cdot,\cdot)$ and $T^0$ is totally skew-symmetric.
\item[(ii)] $T^0=-J\nabla^h J$.
\item[(iii)] $\dfrac{1}{2}J\nabla^h J(\cdot,\cdot)$ is skew-symmetric.
\item[(iv)] $(N,h,J)$ is nearly K\"{a}hler.
\item[(v)] $\nabla^h _{JX} J= -J\nabla_X J$ and $T^0$ is totally skew-symmetric.
\end{description}
\end{prop}
\textbf{Proof.} First, we see that the implications (ii) $\Rightarrow$ (i) and (ii) $\Rightarrow$ (iii) $\Leftrightarrow$ (iv) are obvious. Then by definition of $\nabla^0$ we have 
\begin{equation}\label{T0antisym}
T^0(X,Y)=-\dfrac{1}{2}\left( J(\nabla_X^h J)Y -J(\nabla_Y^h J)X\right)
\end{equation}
which gives us the implication (iii) $\Rightarrow$ (ii). Furtermore, according to proposition~\ref{T-total-antisym}, if $T^0$ is totally skew-symmetric then we have $\nabla^h=\nabla^0 -\frac{1}{2}T^0$ which provides the implication (i) $\Rightarrow$ (ii). Finally the equivalence (i) $\Leftrightarrow$ (v) follows directly from (\ref{T0antisym}).
This completes the proof.\medskip\hfill$\square$\\ 
In particular, a nearly K\"{a}hler manifold endowed with its canonical Hermitian connection satisfies  2 of our 3 hypothesis on the torsion ($T^0$ anticommutes with $J$ and is totally skew-symmetric).\\
Conversely, we have
\begin{thm}\label{thm-equivalence-Nearly}
Let $(N,h,J,\nabla^0)$ be an almost Hermitian manifold with an almost Hermitian connection $\nabla$. If  the torsion $T^0$ of $\nabla^0$ anticommutes with $J$ and is  totally skew-symmetric then $(N,J,h)$ is nearly K\"{a}hler. Moreover, in this case, $\nabla^0$ is the canonical Hermitian connection.\\
Therefore the injective map
$$
(h,J)\longmapsto (h,J,\nabla^h -\dfrac{1}{2}J\nabla^h J)
$$
is in fact a bijection from the set of  nearly K\"ahler structures on $N$ into the set of almost Hermitian structures,  $(h,J,\nabla^0)$, with an almost Hermitian connection whose the torsion is totally skew-symmetric and anticommutes with $J$. 
\end{thm}
\begin{rmk} \em
In other words, in an almost Hermitian manifold there exists at most only one Hermitian connection with totally skew-symmetric and $J$-anticommuting torsion, and if this connection exists then it coincides with the canonical Hermitian connection and the almost Hermitian manifold is nearly K\"{a}hler.
\end{rmk}
Moreover, the third hypothesis (the torsion is parallel) is implied by the first two.
\begin{prop}\emph{[Kirichenko], \cite{Kirichenko77, AFS05}}
If $(N,h,J)$ is nearly K\"{a}hler then  the canonical Hermitian connection has a parallel torsion: $\nabla^0 T^0=0$.
\end{prop}
Now, we can reformulate our theorem~\ref{WZW} by using the right geometric context:
\begin{thm}\label{thm-WZW-Nearly}
Let $(N,h,J)$ be a nearly K\"{a}hler manifold then the equation of holomorphic harmonicity for maps $f\colon L\to N$ is exactly the Euler-Lagrange equation for the sigma model in $N$ with a Wess-Zumino term defined by the 3-form:
$$
H=\dfrac{1}{3}d\Omega_J
$$
where $\Omega_J=\langle J\cdot, \cdot\rangle$ is the K\"{a}hler form.
\end{thm}
\textbf{Proof of theorem~\ref{thm-equivalence-Nearly}.} If the torsion of $\nabla^0$ is totally skew-symmetric then we have $\nabla^h=\nabla^0 - \frac{1}{2}T^0$, so that if moreover $T^0$ anticommutes with $J$ then $- \frac{1}{2}T^0$ is the $J$-anticommuting part\footnote{In the sense of theorem~\ref{decdenabla}.} of $\nabla^h$ i.e. $- \frac{1}{2}T^0=\frac{1}{2}J\nabla^h J$ and $\nabla^0$ is the canonical Hermitian connection. Therefore, we can apply proposition~\ref{prop-equiv-nearly} which allows us to conclude. This completes the proof.\medskip\hfill$\square$\\
\textbf{Proof of theorem~\ref{thm-WZW-Nearly}}
With the notation of proposition~\ref{prop-H-closed} we have $H(X,Y,Z)=\langle JT^0(X,Y),Z\rangle = \langle(\nabla^h J)(X,Y), Z \rangle=\nabla^h\Omega_J(X,Y,Z)$ according to proposition~\ref{prop-equiv-nearly}, and since $\nabla^h \Omega_J$ is a 3-form ($(N,J,h)$ is nearly K\"{a}hler), we have $d\Omega_J=3\nabla^h \Omega_J$. This  completes the proof.\hfill$\square$%\\
\begin{rmk}\em
With this new setting, the  closeness of $H=\frac{1}{3}d\Omega_J$ is obvious.
\end{rmk}
\paragraph{Return to the example of a 3-symmetric space}\index{canonical!almost complex structure}
According to proposition~\ref{3-sym-prop}, a Riemannian (locally) 3-symmetric space $N=G/G_0$ is nearly K\"{a}hler \iif it is naturally reductive. Then, we can reformulate the theorem~\ref{WZ-3sym} as follows:
\begin{thm}\label{WZ-3sym'}
Let $N=G/G_0$ be a  (locally) 3-symmetric space, that we suppose to be Riemannian and naturally reductive, and endowed with its canonical almost complex structure $\undj$ and its canonical connection $\nabla^0$. Let $h$ be a $G$-invariant naturally reductive metric on $N$. Then the equation for holomorphically harmonic maps $f\colon L\to N$ is the   Euler-Lagrange equation for the sigma model in $N$ with the Wess-Zumino term defined by the closed 3-form $H=\dfrac{1}{3}d\Omega_{\undj}$.
\end{thm}\index{Wess-Zumino|)}
\index{canonical!Hermitian connection|)}\index{skew-symmetric torsion|)}
\subsubsection{$J$-twisted harmonic maps}
\begin{defn} Let $f\colon (M,g)\to N$ be a map from a Riemannian manifold $(M,g)$ to a manifold $N$. Let us suppose that the vector bundle $f^*TN$ is naturally endowed with some connection $\overline{\nabla}$. Then we will say that $f$ is roughly harmonic \wrt $\overline{\nabla}$ (or $\overline{\nabla}$-roughly harmonic) if 
$$
\mrm{Tr}_g(\overline{\nabla} df)=0.
$$
\end{defn}
This definition is useful in the case there exists a natural mapping which associates to each map $f\colon (M,g)\to N$ a connection in  the vector bundle $f^*TN$. For example, we have the following.
\begin{thm}
A map $f\colon (L,j_L) \to (N,J,\nabla)$ from  a Riemann surface into an almost complex manifold with a connection $\nabla$ is holomorphically harmonic \iif it is roughly harmonic \wrt
$$
\overline{\nabla}=f^*\nabla + \dfrac{1}{2}JT(df\circ j_L, .).
$$
\end{thm}
\begin{defn}
Let $(N,J)$ be an almost complex manifold with an arbitray connection $\nabla$. Then let us decompose it (in an unique way) as the sum of a $J$-commuting and respectively $J$-anticommuting part:
$
\nabla=\nabla^0+ A$,
where $\nabla^0 J=0$, and $A\in \mal C(T^*N\otimes \End(TN))$, $AJ=-JA$, i.e. $A=\frac{1}{2}J\nabla J$. Then to any map $f\colon L\to N$ let us  associate the connection
$$
\overline{\nabla} =f^*\nabla^0 - J A\circ j_L.
$$
We will say that $f\colon L\to N$ is \emph{$J$-twisted harmonic \wrt $\nabla$} if $f$ is roughly harmonic \wrt $\overline{\nabla}$.
\end{defn}
Now, we can conlude by the following  interpretation of holomorphic harmonicity.
\begin{thm}
Let $(N,J,h)$ be a Nearly K\"{a}hler manifold. Then a map $f\colon L\to N$ is holomorphically harmonic \wrt the canonical Hermitian connnection $\nabla^0$ \iif it is $J$-twisted harmonic \wrt $\nabla^h$.
\end{thm}
\index{holomorphically harmonic|)}
\subsection{The sigma model with a Wess-Zumino term in $\mal G_1$-manifolds}\label{WZ-g-1} \index{Wess-Zumino|(}
\subsubsection{$TN$-valued 2-forms}\label{TN-valued-2-forms}
\index{Nijenhuis tensor!|(}
Let  $(N,J,h)$ be  an almost Hermitian manifold.
In all the section~\ref{WZ-g-1}, each $TN$-valued 2-form on $N$, $B\in \mal C(\mal T)$, will be identified (via the metric $h$) with the corresponding trilinear form, skew-symmetric \wrt to the  2 first  arguments:
$$
B(X,Y,Z):=\langle B(X,Y), Z\rangle.
$$
In particular, the left multiplication by $J$ on $\mal C(\mal T)$ defined a multiplication on the set of corresponding trilinear forms $(JB)(X,Y,Z)=\langle JB(X,Y), Z\rangle$. Moreover, under this identification, the space $\Omega^3(N):=\mal C(\Lambda^3 T^*N)$ of 3-forms will be considered as a subspace of $\mal C(\mal T)$. We denote by $\mrm{Skew}$ the following surjective linear map from $\mal C(\mal T)$ onto $\Omega^3(N)$:
$$
\mrm{Skew}(B)(X,Y,Z) = B(X,Y,Z) + B(Y,Z,X) + B(Z,X,Y).
$$
Let us remark that $\frac{1}{3}\mrm{Skew}(B)$ is the skew-symmetric part of the trilinear form $B$ and $\frac{1}{3}\mrm{Skew}\colon \mal C(\mal T)\to \Omega^3(N)$ is a  projector (called the \emph{Bianci projector} in \cite {Gauduchon}).
To any trilinear form $\alpha\in\mal C(\otimes^3 T^*N)$ will be associated its $J$-twisted trilinear form 
$$
\alpha^c=-\alpha(J\cdot,J\cdot,J\cdot).
$$
In particular, if $\alpha=d\beta$, with $\beta\in \Omega^2(N):=\mal C(\Lambda^2 T^*N)$ then we set $d^c\beta:=\alpha^c$.\\
We will also use the following action of the complex structure $J$ on $C(\mal T)$: for any $B\in \mal C(\mal T)$
\begin{equation}\label{J.T-def}
J\cdot B:=-JB(J\cdot,J\cdot)=J(B^{++} + B^{--}) -J(B^{+-} + B^{-+})
\end{equation}
i.e. in terms of trilinear forms
$$
J\cdot B= B(J\cdot, J\cdot,J\cdot)= -B^c.
$$
Let us remark that $J\cdot(J\cdot B)=-B$\medskip.\\
Furthermore, let $(L,j_L)$ be a Riemann surface  and $B\in \mal C(\mal T)$, then for any map $f\colon L\to N$ and any Hermitian metric $g$ on $(L,j_L)$, we set 
\begin{equation}\label{def-B_g}
B_g(f)=* f^*B=B(f_*TL).
\end{equation}
We will use a second natural action of $J$ on $\mal C(\mal T)$ defined by: for any $B\in\mal C(\mal T)$,
$$
J{\circact} B = B(J\cdot,\cdot,\cdot) + B(\cdot,J\cdot,\cdot) + B(\cdot,\cdot,J\cdot)
$$
and in terms of the components $B^{\eps,\eps'}$:
$$
J{\circact} B = J B^{++} - 3JB^{--} -J(B^{+-} + B^{-+}).
$$
Moreover, by the aid of the two previous natural action, we can define a third action that will turn out to be the relevant one in the interpretation of the maximal odd determined system: for any $B\in\mal C(\mal T)$,
$$
J\star B= \dfrac{1}{2}\left( B(J\cdot, J\cdot,J\cdot) + B(J\cdot,\cdot,\cdot) + B(\cdot,J\cdot,\cdot) + B(\cdot,\cdot,J\cdot)\right) = \dfrac{1}{2}\left( J\cdot B + J\circact B \right), 
$$
and in terms of the components $B^{\eps,\eps'}$:
$$
J\star B=J B^{++} - JB^{--} -J(B^{+-} + B^{-+}).
$$
\begin{rmk}\em
Let us remark that all the three previous actions are independent of the metric $h$, as we can see it from the expressions in terms of the components $B^{\eps,\eps'}$, or more simply by writing their definitions using $TN$-valued 2-forms like in (\ref{J.T-def}): $J\circact B =B(J\cdot,\cdot) + B(\cdot,J\cdot) - JB(\cdot,\cdot)$  and 
$$
J\star B= -\dfrac{1}{2}J\left(B(J\cdot,J\cdot) +  JB(J\cdot,\cdot) + J B(\cdot,J\cdot) + B(\cdot,\cdot)\right).
$$ 
We remark that this last formula - up to the factor $-\frac{1}{2}J$ and to the signs - makes $J\star B$ look like to some kind of torsion tensor of $J$ \wrt $B$ (cf. the definition of $N_J$). \\
In particular, these three actions are defined in a general almost complex manifold $(N,J)$.
\end{rmk}
\begin{rmk}\em
The equations~(\ref{component-eps,eps'}) can be rewritten using the metric $h$ as follows
\begin{equation}\label{component-eps,eps'-h}
B^{\eps,\eps'}=-\dfrac{1}{4} \left(\eps\eps' B(J\cdot,J\cdot,\cdot) -\eps B(J\cdot,\cdot,J\cdot) - \eps' B(\cdot,J\cdot,J\cdot) -B(\cdot,\cdot,\cdot)\right).
\end{equation}
This leads us to define the following action
$$
J\circact_2 B= B(J\cdot,J\cdot,\cdot) + B(J\cdot,\cdot,J\cdot) + B(\cdot,J\cdot,J\cdot)=B-4B^{--}.
$$
\end{rmk}
It is also important to remark that $J\cdot B - J\circact B=4JB^{--}$ and $J\cdot B -J\star B= 2JB^{--}$, so that
\begin{prop}\label{J.T-J*T}
Let $T$ be the torsion of some Hermitian connection $\nabla$ on $(N,J,h)$, then we have
$$
J\cdot T - J\circact T=JN_J\quad \mrm{and}\quad J\cdot T -J\star T= \dfrac{1}{2}JN_J.
$$
\end{prop}
\index{Nijenhuis tensor!|)}
\subsubsection{Stringy Harmonic maps}\index{stringy harmonic|(}
We have seen two different ways to generalise the harmonicity to the case of affine target manifold. The first one is very natural and consist simply to write the harmonic map equation $\mrm{Tr}_g(\nabla df)=0$ for a linear connection $\nabla$. The second one concerns holomorphicaly harmonic maps (from a Riemann surface into an almost complex manifold) and was dicted to us by the geometric equation of the second elliptic integrable system associated to a 3-symmetric space (section~\ref{detodcase}, paragraph: The model case).
Furthermore, the preliminary study of the maximal determined system done in section~\ref{eq-maximal-deter},  leads us to introduce the following generalisation of harmonic maps (which will turn out to be a generalisation of holomorphically harmonic maps for  particular target spaces like nearly K\"{a}hler manifolds).
\begin{defn}\label{defn-j-stringy}
Let $(N,J)$ be an almost complex manifold with $\nabla$ a linear connection then we will say that a map $f\colon L\to N$ from a Riemann surface into $N$ is \emph{stringy harmonic} if it is solution of \emph{the harmonic map equation with a $JT$-term}:
$$
-\tau_g(f) + (J\cdot T)_g(f)=0.
$$
where $g$ is a Hermitian metric on $L$.
\end{defn}
We remark that  if $T$ anticommutes with $J$ then stringy harmoniciy coincides with holomorphic harmonicity (since in this case $J\cdot T= JT$). However, even though stringy harmonicty seems to be the more natural generalisation of (holomorphically) harmonic maps - in particular because of the property $J\cdot(J\cdot B)=-B$ which makes the first action look like very closely to the simple multiplication by $J$, which is not the case for the two other actions- it will turn out that the interpretation of the maximal odd determined system (see section~\ref{eq-maximal-deter}) will use the action $J\star T$ of $J$ on $T$. This leads us to the following modified definition.
\begin{defn}\label{defn-*j-stringy}
Let $(N,J)$ be an almost complex manifold with $\nabla$ a linear connection then we will say that a map $f\colon L\to N$ from a Riemann surface into $N$ is \emph{$\star$-stringy harmonic} if it is solution of \emph{the modified stringy harmonic maps equation}:
$$
-\tau_g(f) + (J\star T)_g(f)=0.
$$
where $g$ is a Hermitian metric on $L$.
\end{defn}
We remark that  if $T$ anticommutes with $J$ then $\star$-stringy harmoniciy coincides with antiholomorphic harmonicity (since in this case $J\star T= -JT$)\medskip.\\
Now, we will see that, under some hypothesis, the two previous definitions are in fact equivalent in the sense that there exists a new almost complex structure $J^\star$ such that $J\star T=J^\star\cdot T$.
\begin{prop}\label{prop-JstarJ}
Let $(N,J)$ be an almost complex manifold. Let us suppose that there exists a $J$-invariant decomposition $TN=E^+ \oplus E^-$ and some $B\in\mal C(\mal T)$ such that  we have 
\begin{equation}\label{conditions-B}
\forall \alpha,\alpha'\in \Z_2,\qquad B^{**}(E^\alpha,E^{\alpha'})\subset E^{\alpha\alpha'} \quad  \mrm{and}  \quad B^{--}(E^\alpha,E^{\alpha'})\subset E^{-\alpha\alpha'}.
\end{equation}
where $B^{**}:=B^{++} + B^{+-} + B^{-+}=B-B^{--}$.
Let us define $J^\star= J_{|E^+}\oplus -J_{|E^-}$, then we have
$$
J\star B=J^\star\cdot B.
$$
\end{prop}
\proof
Applying equation \eqref{component-eps,eps'} to $B^{**}$, we obtain
$$
B^{\eps,\eps'}(E^\alpha,E^{\alpha'})\subset E^{\alpha,\alpha'}, \qquad \forall (\eps,\eps') \in \Z_2^2\setminus\{(-,-)\}.
$$
Moreover, we have 
$$
B=\sum_{\eps,\eps'\in\Z^2,  \alpha,\alpha'\in\Z^2} B_{|E^\alpha\times E^{\alpha'}}^{\eps,\eps'}.
$$
Besides, since $B_{|E^\alpha\times E^{\alpha'}}^{\eps,\eps'}(J^\star\cdot\, ,\cdot) = \eps\alpha J B_{|E^\alpha\times E^{\alpha'}}^{\eps,\eps'}=\eps\alpha (\alpha\alpha') J^\star B_{|E^\alpha\times E^{\alpha'}}^{\eps,\eps'} $., it follows that $B_{|E^\alpha\times E^{\alpha'}}^{\eps,\eps'}$ is of type $(\eps \alpha',\eps'\alpha)$ \wrt to $J^\star$, if $(\eps,\eps')\neq \{(-,-)\}$. And $B_{|E^\alpha\times E^{\alpha'}}^{- -}$ is of type $(\alpha',\alpha)$. Therefore denoting by $\Bar B^{\eps,\eps'}$ the $(\eps,\eps')$-component of $B$ \wrt $J^\star$, we have 
$$
\Bar B^{\eps,\eps'}= \sum_{(\alpha',\alpha)\neq -(\eps,\eps')} B_{|E^\alpha\times E^{\alpha'}}^{\eps\alpha',\eps'\alpha}
+ B_{|E^{\eps'}\times E^\eps}^{--}.
$$
Finally, using this last equation  we check by computation that $J^\star \cdot B = J\star B$. \comprf \hsq\medskip\\
Therefore, this yields the following corollary.
\begin{cory}\label{cory-J*-J}
Let $(N,J)$ be an almost complex manifold with $\nabla$ a linear connection. Let us suppose that there exists a $J$-invariant decomposition $TN=E^+ \oplus E^-$ such that the torsion $T$ of $\nabla$ satisfies the conditions (\ref{conditions-B}). Let $J^\star= J_{|E^+}\oplus -J_{|E^-}$. Then the $\star$-stringy harmonicity with respect to $J$ is exactly the stringy harmonicity\index{stringy harmonic|)} with respect to $J^\star$.
\end{cory}
\begin{rmk}\label{rmk-hold-vectorbundle}\em
Let us remark that proposition~\ref{prop-JstarJ} is an algebraic identity, and holds in any complex vector space, or more generally any complex vector bundle $(E,J)$.
\end{rmk}

%
%%%%%%%%%%%%%%%%%%%%%%%%%%%%%%%%%%%%%%%%%%%%%%%%%%%%%%%%%%%%%%%%%

\subsubsection{Almost Hermitian $\mal G_1$-manifolds}\label{AlmostHermG_1}
\index{g1 manifold@$\mal G_1$-manifold|(}\index{Nijenhuis tensor!|(}\index{skew-symmetric torsion|(}
In all this subsection, we consider $(N,J,h)$ an almost Hermitian manifold with a Hermitian connection $\nabla$, whose the torsion is denoted by $T$.\medskip\\
We prove easily the following.
\begin{prop}
The components $B^{\eps\eps'}$ of an element $B\in\mal C(\mal T)$, considered as trilinear forms,  satisfy the following properties:
\begin{eqnarray*}
 B^{++} & \in  & \mal C\left( (\Lambda^{2,0}\otimes \Lambda^{0,1})\oplus(\Lambda^{0,2}\otimes \Lambda^{1,0}) \right)\\ 
B^{+-} & \in & \mal C\left( (\Lambda^{1,0}\otimes \Lambda^{0,1}\otimes\Lambda^{0,1}) \oplus  (\Lambda^{0,1}\otimes \Lambda^{1,0}\otimes\Lambda^{1,0})\right) \\
B^{-+} & \in & \mal C\left( (\Lambda^{1,0}\otimes \Lambda^{1,0}\otimes\Lambda^{0,1}) \oplus  (\Lambda^{0,1}\otimes \Lambda^{0,1}\otimes\Lambda^{1,0})\right) \\
B^{--} & \in & \mal C\left( (\Lambda^{2,0}\otimes \Lambda^{1,0})\oplus(\Lambda^{0,2}\otimes \Lambda^{0,1}) \right),
\end{eqnarray*}
where $\Lambda^{p,q}=\Lambda^{p,q} T^*N$.
\end{prop} 
The results of this is: 
\begin{cory}\label{cory-B-eps,eps'}
Let $B\in\Omega^3(N)\subset\mal C(\mal T)$ be a 3-form on $N$, then $B^{--}$ is also a 3-form and is of type $(3,0) + (0,3)$.
Moreover $B^{++} + B^{+-} + B^{-+}$ is a 3-form of type $(2,1) + (1,2)$ and we have the following relations:
\begin{eqnarray*}
B^{++}(X,Y,Z)= B^{+-}(Z,X,Y)\\
B^{-+}(X,Y,Z)=B^{+-}(Y,Z,X)
\end{eqnarray*} 
in other words $B^{**}:=B^{++} + B^{+-} + B^{-+}=\mrm{Skew}(B^{\eps\eps'})$, $\forall (\eps,\eps')\in \Z_2\setminus\{(-,-)\}$. \\
In particular, let us suppose that the torsion $T$ of the Hermitian connection $\nabla$ is totally skew-symmetric, then $T^{--}$ is also a 3-form and is of type $(3,0) + (0,3)$, and $T^{**}$ is a 3-form of type $(2,1) + (1,2)$. More particulary, the Nijenhuis tensor $N_J$ is totally skew-symmetric.
\end{cory}
\begin{cory}\label{cory-dOmega-J}
Let us suppose that the torsion $T$ of the Hermitian connection $\nabla$ is totally skew-symmetric, then
$$
d\Omega_J= 3JT^{--} + J(T^{+-} + T^{-+} -T^{++})=-J\circact T
$$
i.e. 
$$
d\Omega_J= JN_J -J\cdot T .
$$
\end{cory}
\textbf{Proof.}
Since $T$ is skew-symmetric, we have $\nabla=\nabla^h +\frac{1}{2}T$, so that $\nabla J=0$ implies
$$
\nabla^h\Omega_J=-\dfrac{1}{2}(T(\cdot, J\cdot,\cdot) + T(\cdot,\cdot,J\cdot))
$$
and therefore applying the operator $\mrm{Skew}$ to that and using the fact that $T$ is skew-symmetric we obtain
$$
d\Omega_J=-J\circact T.
$$
Then the last assertion follows from proposition~\ref{J.T-J*T}. This completes the proof.\medskip\hfill $\square$\\
Now, we can conclude  that
\begin{thm}\label{Unik}\index{connection!characteristic|(}
An almost  Hermitian manifold $(N,J,h)$ admits a Hermitian connection with totally skew-symmetric  torsion \iif the Nijenhuis tensor $N_J$ is itself totally skew-symmetric. In this case, the connection is unique and determined by its torsion which is given by 
$$
T=-d^c\Omega_J + N_J.
$$
\end{thm}
\textbf{Proof.} If such a Hermitian connection with skew-symmetric torsion exists, then according to corollary~\ref{cory-B-eps,eps'}, $N_J$ is itself skew-symmetric and moreover, according to corollary~\ref{cory-dOmega-J}, we have $d\Omega_J= JN_J -J\cdot T=J\cdot(N_J-T)$ therefore $T=J\cdot d\Omega_J + N_J$. This proves the uniqueness.\\
Conversely, let us suppose that $N_J$ is skew-symmetric and let $\nabla$ be the metric connection defined by the torsion $T=J\cdot d\Omega_J + N_J$, i.e. $\nabla=\nabla^h-\frac{1}{2}T$. We have to check that $\nabla J=0$. Let us recall (\cite{KN}, proposition~4.2) that we have
\begin{equation}\label{nablah-Omega-J}
2(\nabla_X^h\Omega_J)(Y,Z)=d\Omega_J(X,Y,Z) - d\Omega_J(X,JY,JZ) + N_J(Y,Z,JX).
\end{equation}
Applying $\mrm{Skew}$ to that, we obtain
\begin{multline*}
2d\Omega_J(X,Y,Z)=3d\Omega_J(X,Y,Z) - d\Omega_J(X,JY,JZ) - d\Omega_J(JX,Y,JZ) -d\Omega_J(JX,JY,Z)\\
 + N_J(Y,Z,JX) + N_J(X,Y,JZ) + N_J(Z,X,JY)
\end{multline*}
therefore 
\begin{eqnarray}
-4(d\Omega_J)^{0,2}(X,Y,Z) & = &  N_J(Y,Z,JX) + N_J(X,Y,JZ) + N_J(Z,X,JY)\label{dOmega--}\\
 & = & 3N_J(X,Y,JZ) \text{ (since $N_J$ is skew-symmetric (and of type $(0,2)$)}\nonumber.
\end{eqnarray}
Now, we can compute
\begin{eqnarray*}
\nabla\Omega_J & = & \left( \nabla^h + \dfrac{1}{2}(J\cdot d\Omega_J + N_J)\right) \Omega_J \\
& = &  \nabla^h\Omega_J + \dfrac{1}{2}(N_J(X,JY,Z) + N_J(X,Y,JZ) -d\Omega_J(JX,Y,JZ) - d\Omega_J(JX,JY,Z))\\ 
 &  =  & \dfrac{1}{2}(d\Omega_J(X,Y,Z) - d\Omega_J(X,Y,JZ) -d\Omega_J(JX,Y,JZ) - d\Omega_J(JX,JY,Z)\\
 &  &  + N_J(Y,Z,JX) + N_J(X,JY,Z) + N_J(X,Y,JZ))\\
 & = & \dfrac{1}{2}(3N_J(X,Y,JZ) + N_J(Y,Z,JX) + N_J(X,JY,Z) + N_J(X,Y,JZ)) = 0.
\end{eqnarray*}
This completes the proof.\hfill$\square$
\begin{rmk}\label{rmk-frie-ivan}\em 
This theorem  can also be deduced  from  a more general result  of Gauduchon \cite[Proposition~2]{Gauduchon} (see also section~\ref{sect-Gauduchon} below). Moreover, it has  already been  proved  by  Friedrich-Ivanov \cite[theorem 10.1]{Friedrich-ivanov} (but without writting completely the proof). To  our knowledge a completely written proof has never been given in the literature\footnote{See also \cite{srni} for the proof of the fact that the connection given in theorem~\ref{Unik} satisfies $\nabla J=0$.}. In fact, Friedrich and Ivanov present a unified approach to construct $G$-connections, for any $G$-structure on some manifold, with skew-symmetric torsion. For example, they characterize the class of $G_2$-structures and the class of contact metric structures  for which such a connection exists, and  prove the uniqueness of this connection.
\end{rmk} 
\begin{defn}\index{Gray-Hervella classification} 
The unique  Hermitian connection with skew-symmetric torsion is called the characteristic connection. According to the Gray-Hervella classification \emph{\cite{Gray-Hervella}} of almost Hermitian manifolds, $(N,J,h)$ admits a skew-symmetric Nijenhuis tensor \iif if is of class $\mal W_1\oplus \mal W_2 \oplus \mal W_3=:\mal G_1$ (see \emph{\cite{Gray-Hervella}}).
These manifolds are called $\mal G_1$-manifolds and according to the previous theorem they are exactly the almost Hermitian manifolds which admit a  characteristic connection.
\end{defn}
\begin{prop}\label{prop-dH_J=0}
Let us suppose that the almost Hermitian manifold $(N,J,h)$ is a $\mal G_1$-manifold. Let us suppose that its  characteristic connection $\nabla$ has a parallel torsion  $\nabla T=0$. Then the 3-form 
$$
H(X,Y,Z)=T(JX,JY,JZ)=\langle (J\cdot T)(X,Y),Z\rangle
$$
is closed $dH=0$.
\end{prop}
\textbf{Proof.}
Since, according to corollary~\ref{cory-dOmega-J}, we have $H=-d\Omega_J + JN_J$, we only have to prove that the 3-form $JN_J=4JT^{--}$ is closed. Moreover, since $J$ is $\nabla$-parallel, so is the decomposition $\mrm{Bil}(E)=\oplus_{\Z_2\times\Z_2}\mrm{Bil}^{\eps,\eps'}(E,J)$, so that if $T$ is $\nabla$-parallel then so are its components $T^{\eps\eps'}$. Then the $\nabla$-parallel 3-forms, $H$, $d\Omega_J$ and $JT^{--}$ satisfies the following formula for $\nabla$-parallel 3-forms $\alpha$:
\begin{eqnarray*}
d\alpha(X_0,X_1,X_2,X_3) &  =  & \sum_{i=0}^3  (-1)^{i}\nabla_{X_i} \alpha(X_0,\ldots, \hat X_i,\ldots, X_3)\\
  &   &  - \sum_{0\leq i<j\leq 3} (-1)^{i+j} \alpha(T(X_i,X_j), X_0,\ldots,\hat X_i,\ldots, \hat X_j,\ldots, X_3))\\
 & = & \underset{i,j,k}{\mak S} \alpha(T(X_0,X_i), X_j, X_k) + \alpha(T(X_i,X_j), X_0, X_k)\\
 & = & \underset{i,j,k}{\mak S} \alpha(T(X_0,X_i), X_j, X_k) + \underset{i,j,k}{\mak S} \alpha(T(X_j,X_k), X_0, X_i)\\
& = & \underset{X,Y,Z}{\mak S} \alpha(T(V,Z), X,Y) + \alpha(T(X,Y),V,Z)
\end{eqnarray*}
where we have set $(X,Y,V,Z)=(X_2,X_3,X_0,X_1)$. Then applying this formula to $d\Omega_J$, we obtain
\begin{eqnarray*}
0=-d(d\Omega_J) & = & -d\left(3JT^{--} + J(T^{+-} + T^{-+} -T^{++})\right)\\
 & = & \underset{X,Y,Z}{\mak S} \sum_{(X,Y)\rightleftarrows(V,Z)}\langle \left(3 T^{--} - T^{++} + T^{1,1}\right)(X,Y), J \left( T^{--} + T^{++} + T^{1,1}\right) (V,Z)\rangle
 \end{eqnarray*}
where $(X,Y)\rightleftarrows(V,Z)$ means that we sum on the set  $\{(X,Y,V,Z),(V,Z,X,Y)\}$. After a straightforward computation, we find
\begin{eqnarray}\label{ddOmega-J}
0=-d(d\Omega_J) & = & \underset{X,Y,Z}{\mak S} \sum_{(X,Y)\rightleftarrows(V,Z)} \left\lbrace 4\langle T^{--}(X,Y), JT^{++}(V,Z)\rangle \right. \\
&  & + 2 \langle T^{--}(X,Y), JT^{1,1}(V,Z)\rangle\nonumber\\
&  &  - 2  \left. \langle T^{++}(X,Y), JT^{1,1}(V,Z)\rangle.
\right\rbrace\nonumber
\end{eqnarray} 
Now, let us consider 4-linear forms on the variable $(X,Y,V,Z)\in TN^4$ and the associated decomposition $\otimes^4 T^*N^\C=\oplus_{\eps\in (\Z_2)^4} \Lambda^{\eps_1}\otimes \Lambda^{\eps_2}\otimes\Lambda^{\eps_3} \otimes\Lambda^{\eps_4}$, where $\Lambda^+=\Lambda^{1,0}T^*N$ and $\Lambda^-=\Lambda^{0,1}T^*N$. Then the term in the first line of (\ref{ddOmega-J}) is in $(\otimes^4\Lambda^+)\oplus (\otimes^4\Lambda^-)$ whereas the terms in the second and third lines are in $(\otimes^{+++-})\oplus (\otimes^{---+})\oplus (\otimes^{++-+})\oplus (\otimes^{--+-})$, where $\otimes^{\eps_1,\eps_2,\eps_3,\eps_4}= \Lambda^{\eps_1}\otimes \Lambda^{\eps_2}\otimes\Lambda^{\eps_3} \otimes\Lambda^{\eps_4}$, $\forall \eps\in(\Z_2)^4$. Hence the sums $\underset{X,Y,Z}{\mak S} \sum_{(X,Y)\rightleftarrows(V,Z)}$ of these terms are respectively in $\Lambda^{4,0} \oplus \Lambda^{0,4}$ (first line) and in $\Lambda^{3,1} \oplus \Lambda^{1,3}$ (second and third lines). Therefore we obtain (in particular) that the first line vanishes
$$
\underset{X,Y,Z}{\mak S} \sum_{(X,Y)\rightleftarrows(V,Z)}  4\langle T^{--}(X,Y), JT^{++}(V,Z)\rangle =0.
$$
Let us apply this to the computation of $d(JN_J)$:
\begin{eqnarray*}
d(JN_J) & = & \underset{X,Y,Z}{\mak S} \sum_{(X,Y)\rightleftarrows(V,Z)}  4\langle JT^{--}(X,Y), T^{**}(V,Z) \rangle\\ 
& = &  \underset{X,Y,Z}{\mak S} \sum_{(X,Y)\rightleftarrows(V,Z)}  4\langle JT^{--}(X,Y), T^{1,1}(V,Z) \rangle.
\end{eqnarray*}
We see that $\langle JT^{--}(X,Y), T^{+-}(V,Z)\rangle$ is in $(\otimes^{+++-}) \oplus (\otimes^{---+})$. But since $JT^{--}$ is a 3-form, we have
$$
\langle JT^{--}(X,Y), T^{+-}(V,Z)\rangle =-\langle JT^{--}(X,T^{+-}(V,Z)), Y\rangle
$$
and this second 4-linear form (in the variable $(X,Y,V,Z)$) is in $(\otimes^{++-+}) \oplus (\otimes^{--+-})$, which imposes that $\langle JT^{--}(X,Y), T^{+-}(V,Z)\rangle =0$, $\forall (X,Y,V,Z)\in TN^4$. We can prove the same result if we replace $T^{+-}$ by $T^{-+}$. Therefore $d(JN_J)=0$. This completes the proof.\medskip\hfill $\square$\\
Moreover, according to proposition~\ref{J.T-J*T},  we deduce
\begin{prop}
Let us suppose that the almost Hermitian manifold $(N,J,h)$ is a $\mal G_1$-manifold. Let us suppose that its  characteristic connection $\nabla$ has a parallel torsion  $\nabla T=0$. Then the 3-form 
$$
H^\star(X,Y,Z)=\langle (J\star T)(X,Y),Z\rangle
$$
is closed $dH^\star=0$.
\end{prop}
Now, we can conclude with the following variational interpretation of the stringy harmonicity.
\begin{thm}\label{thm-WZW-G1}\index{stringy harmonic}
Let us suppose that the almost Hermitian manifold $(N,J,h)$ is a $\mal G_1$-manifold. Let us suppose that its  characteristic connection $\nabla$ has a parallel torsion  $\nabla T=0$.\\
$\bullet$ Then the equation for stringy harmonic maps $f\colon L\to N$ is exactly the Euler-Lagrange equation for the sigma model in $N$ with a Wess-Zumino term defined by the closed 3-form
$$
H=-d\Omega_J + JN_J.
$$ 
$\bullet$ Moreover the equation for $\star$-stringy harmonic maps $f\colon L\to N$ is exactly the Euler-Lagrange equation for the sigma model in $N$ with a Wess-Zumino term defined by the closed 3-form
$$
H^\star=-d\Omega_J + \dfrac{1}{2}JN_J.
$$ 
\end{thm}
\begin{rmk} \em
We remark that the two previous sigma model differ by the Wess-Zumino term defined by the 3-form $\dfrac{1}{2}JN_J$.
\end{rmk}
\index{connection!characteristic|)}\index{g1 manifold@$\mal G_1$-manifold|)}\index{Nijenhuis tensor!|)} \index{skew-symmetric torsion|)}
\subsubsection{Characterization of Hermitian connections in terms of their torsion}\label{sect-Gauduchon}
In this subsection, we will give a result of Gauduchon \cite[Proposition~2]{Gauduchon} characterizing the Hermitian connections in terms of their torsion. We need to write it with our notations and inside our setting and to write one proof in such a way that it will appear clearly that this result holds as well for Riemannian $f$-structure (see theorem~\ref{thm-nablaOmega_F=0}) so that we will not have to reprove it (at least not entirely) in this more general context.
\begin{thm}\label{Gauduchon}
Let $(N,J,h)$ be a Hermitian manifold. Then a metric connection $\nabla$ is almost complex \iif the following statements hold
$$
N_J=4T^{0,2} \quad \mrm{and}\quad \mrm{Skew} (T^{2,0} - T^{1,1})=(d^c\Omega_J)^{**}.
$$
\end{thm}
\textbf{Proof.}
The metric connection $\nabla$ can be written in the form $\nabla=\nabla^h -\dfrac{1}{2}(T +U)$, where $U(X,Y,Z)=T(Z,X,Y) + T(Z,Y,X)$. Then $\nabla \Omega_J=0$ \iif 
\begin{equation}\label{nablahJ=AJ}
\nabla^h J=-\dfrac{1}{2}[ (T+U),J]=- \left( (T+U)^{--}J +  (T+U)^{+-}J\right)
\end{equation}
but\footnote{with $U(B)(X,Y,Z)=B(Z,X,Y) + B(Z,Y,X)$, for any $B\in\mal T$.} 
$$
U^{--}(X,Y,Z)=U(T^{--}) \quad \mrm{and}\quad U^{+-}(X,Y,Z) =T^{-+}(Z,X,Y) + T^{++}(Z,Y,X).
$$
so that, according to  (\ref{nablah-Omega-J})
\begin{eqnarray*}
(d\Omega_J)^{+-}=(\nabla^h \Omega_J)^{+-} & = & -T^{+-}(X,JY,Z) - T^{-+}(Z,X,JY) - T^{++}(Z,JY,X)\\
 & = & -\left(\mrm{Skew}(T)^{+-}(X,JY,Z) -2T^{++}(JY,Z,X)\right)\\
& = &  -J\cdot\left(\mrm{Skew}(T)^{+-}(X,Y,Z) -2T^{++}(Y,Z,X)\right).
\end{eqnarray*}
Therefore applying $\mrm{Skew}$,
\begin{eqnarray*}
J\cdot (d\Omega_J)^{**}= (J\cdot d\Omega_J)^{**} & = & \mrm{Skew}(T)^{**} -2\mrm{Skew}(T^{++}(Y,Z,X)) \\
 & = & \mrm{Skew}(T^{**} -2T^{++})\\
& = & \mrm{Skew}(-T^{2,0} + T^{1,1})
\end{eqnarray*}
Besides, taking the $(\eps,\eps')$-component of equation (\ref{nablahJ=AJ}) for $(\eps,\eps')=(++),(-+)$ instead of $(+-)$ would give the same result. Now, it remains to see what gives us the $(-,-)$-component of this equation. Equations (\ref{nablah-Omega-J}) and then (\ref{dOmega--}) yield
\begin{eqnarray*}
2(\nabla^h\Omega_J)^{0,2} & = & 2(d\Omega_J)^{0,2} + N_J(Y,Z,JX)= -\dfrac{1}{2} \left( N_J(Y,Z,JX) + N_J(JX,Y,Z) + N_J(Z,JX,Y) \right) + N_J(Y,Z,JX)\\
&  = &  -\dfrac{1}{2}\left(  N_J(JX,Y,Z) + N_J(Z,JX,Y) + N_J(Z,Y,JX)\right)
\end{eqnarray*}
so that the $(0,2)$-component of equation (\ref{nablahJ=AJ}) is written
$$
-\dfrac{1}{4}\left(  N_J(X,JY,Z) + N_J(Z,X,JY) + N_J(Z,JY,X)\right)=-(T^{--}(X,JY,Z) + T^{--}(Z,X,JY) + T^{--}(Z,JY,X)).
$$
Using the fact that the map $B\in \mal T\mapsto B + U(B)\in T^*N\otimes\so(TN)$ is bijective\footnote{See remark~\ref{T+U}.}, we obtain
$$
T^{--}=\dfrac{1}{4}N_J.
$$
This completes the proof.\hfill$\square$
%
%(\cdot,J\cdot,\cdot) + (T+U)(\cdot,\cdot,J\cdot)\right)
%
%
%
\subsubsection{The example of a naturally reductive homogeneous space}
\index{canonical!connection, $G$-invariant|(}
In this subsection we consider $G/K$ a reductive homogeneous space and  we denote by $\g=\kk\oplus \mk$ a reductive decomposition of the Lie algebra $\g$.
\begin{thm}\label{thm-Hom-G1}\index{connection!characteristic}
Let $N=G/K$ be a Riemannian naturally reductive homogeneous space. Then the canonical connection is a metric connection with skew-symmetric torsion (\wrt any naturally reductive $G$-invariant metric $h$). Let us suppose also that $N=G/K$ is endowed with some  $G$-invariant complex structure $J$ (i.e.  $\mk$ is endowed with  some $\Ad K$-invariant complex structure $J_0$). If moreover one can choose a naturally reductive $G$-invariant metric $h$ for which $J$ is orthogonal\footnote{which means that denoting by $G(\mk)$, the compact subgroup in $GL(\mk)$ generated by $\Lambda_\mk(\mk) :=\{[\ad_\mk(X)]_\mk, X\in\mk\}\subset\mak{gl}(\mk)$, and by $\langle G(\mk), J_0\rangle$ the closed subgroup generated by $G(\mk)$ and $J_0$, then $\langle G(\mk), J_0\rangle/G(\mk)$ is compact.}, then $(N,h,J)$ is an almost Hermitian $\mal G_1$-manifold and its characteristic connection coincides with the canonical connection. Therefore, in this case its  characteristic connection $\nabla$ has a parallel torsion  $\nabla T=0$. 
\end{thm}
\textbf{Proof.} The naturally reductivity means exactly that the torsion of the canonical connection is skew-symmetric. Then according to theorem~\ref{Unik}, we deduce that $(N,h,J)$ is $\mal G_1$-manifold. This completes the proof.\hfill$\square$%\\ 
\begin{rmk}\em
In particular, we see that the Nijenhuis tensor is skew-symmetric. We can recover that by saying that since the $G$-invariant complex structure is parallel with respect to the canonical connection, then $N_J=4T^{--}$ and moreover since $T$ is a 3-form so is its component $T^{--}$.
\end{rmk}
By definition of ($\star$-)stringy harmonicity and the expression of the torsion of $\nabla^0$ in terms of the Lie bracket, we have the following.
\begin{prop}\label{prop-stringy-hom}
Let $N=G/K$ be a Riemannian homogeneous manifold endowed with a $G$-invariant complex structure $J$. Let  $f\colon L\to N$ be a smooth map, $F\colon L\to G$ be a (local) lift of $f$ and $\alpha=F^{-1}.dF$ the corresponding Maurer-Cartan form. Then in terms of $\alpha$, the equation of stringy harmonicity (\wrt $\nabla^0$) is written
$$
d*\alpha_\mk + [\alpha_\kk\wedge *\alpha_\mk] - \dfrac{1}{2}J_0\left[J_0 \alpha _\mk\wedge J_0\alpha_\mk \right]_\mk =0
$$
whereas the equation of $\star$-stringy harmonicity is written:
$$
d*\alpha_\mk + [\alpha_\kk\wedge *\alpha_\mk] + \dfrac{1}{2}\left[J_0 \alpha _\mk\wedge\alpha_\mk \right]_\mk + \dfrac{1}{4}J_0\left( \left[J_0\alpha_\mk \wedge \,J_0 \alpha _\mk \right]_\mk + \left[ \alpha _\mk\wedge\alpha_\mk \right]_\mk\right) =0
$$
where $J_0$ is the complex structure on $\mk$ corresponding to  $J$.
\end{prop}
\index{g1 manifold@$\mal G_1$-manifold|)}
\subsubsection{Geometric interpretation of the maximal determined odd case.}
\index{canonical!almost complex structure|)} \index{determined maximal@determined, maximal|(}\index{odd case|(} \index{k symmetric space@$k'$-symmetric space|(}
In this subsection, we suppose that $N=G/K$ is a (locally) $(2k+1)$-symmetric space, and we use the notations  and the conventions of  \ref{finitorderauto}. In particular, we have defined the subspace $\mk_j\subset\mk$, for $1\leq j\leq k$ (see subsection~\ref{2.1.2}). We will set $\mk_{-j}=\mk_j$,  $\mk_{j+ 2k+1}=\mk_j$, and $\mk_0=\{ 0 \}$, so that $\mk_j$ is now defined for all $j\in \Z$.
%subsections~\ref{def-g-tau} and
%
%
\begin{prop}\label{T-eps-esp'-nabla0}
Let us suppose that $N=G/K$ is a (locally) $(2k+1)$-symmetric space endowed with its canonical almost complex structure $\undj$ and its canonical connection $\nabla^0$. Then the $(\eps,\eps')$-component of the torsion $T$ of 
$\nabla^0$ are given by
\begin{eqnarray*}
\widetilde{T^{++}(X,Y)} & = &  -\sum_{ \underset {1\leq i\leq j\leq k} { i+j \leq k}   } 
\left(1- \dfrac{\delta_{ij}}{2}\right)  \left(     [X_{\mk_i}, Y_{\mk_j}]_{\mk_{i+j}} +  [X_{\mk_j}, Y_{\mk_i}]_{\mk_{i+j}}\right)\\ 
\widetilde{T^{--}(X,Y)} & = &  -\sum_{ \underset {1\leq i\leq j\leq k} {i+j \geq k+1}   } 
\left(1- \dfrac{\delta_{ij}}{2}\right)  \left(     [X_{\mk_i}, Y_{\mk_j}]_{\mk_{i+j}} +  [X_{\mk_j}, Y_{\mk_i}]_{\mk_{i+j}}\right).\\
\widetilde{T^{+-}(X,Y)} & = &  -\sum_{1\leq i\leq j\leq k} [X_{\mk_j}, Y_{\mk_i}]_{\mk_{j-i}}\\
\widetilde{T^{-+}(X,Y)} & = &  -\sum_{1\leq i\leq j\leq k} [X_{\mk_i}, Y_{\mk_j}]_{\mk_{j-i}}.
\end{eqnarray*}
where $X,Y\in \mal C(TN)$ with lifts $X_\mk,Y_\mk \in \mal C^\infty(G,\mk)$, and $\widetilde{T^{\eps,\eps'}(X,Y)}\in \mal C^\infty(G,\mk)$ denotes the lift of $T^{\eps,\eps'}(X,Y)\in \mal C(TN)$. An other possibility is to consider that $X,Y\in TN$ are tangent vectors at some point $y\in N$ and that we have chosen $g\in G$ such that $g.G_0=y$ and  that we have set $X=\Ad g (X_\mk)$, and $Y=\Ad g (Y_\mk)$. Then the above equation, when written  in the form $\widetilde A = B_\mk$, means in fact that we have $A=\Ad g (B_\mk)$.
\end{prop}
\proof
This follows from the fact that the lift in $G$ of the torsion of $\nabla^0$ is given by $\tl T (X_\mk, Y_\mk)= -[ X_\mk, Y_\mk]_\mk$, from the fact that the commutation relations $[\g_i^\C,\g_j^\C] \subset \g_{i+j}^\C$ implies the following relations $[\mk_i,\mk_j]\subset \mk_{i+j}\oplus \mk_{i-j}$, and finally from the definition of $\undj$. \hsq 
\begin{thm}\label{oddmaxdeter-geom-interpret}
Let us suppose that $N=G/K$ is a (locally) $(2k+1)$-symmetric space endowed with its canonical almost complex structure $\undj$ and its canonical connection $\nabla^0$.  Then the associated maximal determined system, $\syst(2k,\tau)$ is \emph{the equation of $\star$-stringy harmonicity for the geometric map $f\colon L\to N$}:
${(\nabla^0)}^* df + (\undj\star T^0)(f)=0$. \footnote{Where we have removed the index "$g$" which specifies that the previous terms are computed with respect to some Hermitian metric $g$ on $L$.}\\
Moreover, if we consider now  that $N=G/K$ is endowed with the almost complex structure $\undj^\star:=\oplus_{j=1}^{k}(-1)^{j}{\undj}_{[\mk_j]}$, then this system is \emph{the equation of stringy harmonicity for the geometric map $f\colon L\to N$}: ${(\nabla^0)}^* df + (\undj^\star \cdot T^0)(f)=0$.\\
Now, Suppose also that $N=G/K$ is naturally reductive\footnote{And we choose a naturally reductive $G$-invariant metric $h$ for which $\taum$ and thus $\undj$ are orthogonal. See the Appendix, lemma~\ref{Gm-J_0-compact} for the existence of such an metric.}.
Therefore, the previous system is  exactly the Euler-Lagrange equation for the sigma model in $N$ with a Wess-Zumino term defined by the closed 3-form
$$
H^\star=-d\Omega_{\undj} + \dfrac{1}{2}\undj N_{\undj}.
$$ 
Moreover, if  $N=G/K$ is endowed with the almost complex structure $\undj^\star$, the previous system is  exactly the Euler-Lagrange equation for the sigma model in $N$ with a Wess-Zumino term defined by the closed 3-form
$$
H=-d\Omega_{\undj^\star} + \undj^\star N_{\undj^\star}.
$$  
\end{thm}
\proof
The first point follows from theorem~\ref{thm-max-deter-odd}, proposition~\ref{prop-stringy-hom}, and proposition~\ref{T-eps-esp'-nabla0}. Then the second point follows from corrolary~\ref{cory-J*-J}  and  proposition~\ref{T-eps-esp'-nabla0}. Let us make precise that  we apply corrolary~\ref{cory-J*-J} with the $J$-invariant decomposition $TN=E^+\oplus E^-$,  defined by $E^+=\oplus_{\underset{j\text{ even}}{j=1}}^k [\mk_j]$, and $E^-=\oplus_{\underset{j\text{ odd}}{j=1}}^k [\mk_j]$. Then according to proposition~\ref{T-eps-esp'-nabla0}, this decomposition  satisfies the hypothesis of corrolary~\ref{cory-J*-J}. Finally the two last points (Wess-Zumino formulations) follow then from theorem~\ref{thm-WZW-G1} and \ref{thm-Hom-G1}.  \comprf \hsq
\index{Wess-Zumino|)}\index{determined maximal@determined, maximal|)}
\index{canonical!almost complex structure|)} \index{canonical!connection, $G$-invariant|)}\index{odd case|)} \index{k symmetric space@$k'$-symmetric space|)}
\subsection{Stringy harmonicity  versus Holomorphic harmonicity}\index{holomorphically harmonic|(}
In this section, we will compare these two notions and prove that in the case of an almost complex affine manifold, these are equivalent. Indeed, the stringy harmonicity \wrt an almost complex connection $\nabla$ is equivalent to the holomorphic harmonicity \wrt a new almost complex connection $\nabla^\star$. In particular, according to subsection~\ref{holharmmap}, stringy harmonicity has an interpretation in terms of holomorphic 1-forms.
\begin{prop}\label{prop-stringy-star-bullet}
Let $(N,J,\nabla)$ be an almost complex affine manifold  and $(L,j_L)$ a Riemann surface. Let us set
\begin{eqnarray*}
\nabla^\bullet & = & \nabla - (T^{2,0} + T^{0,2})\\
\nabla^\star & = & \nabla - T^{2,0}
\end{eqnarray*}
Then \begin{description}
\item[$\bullet$] $f\colon L\to N$ is stringy harmonic \wrt $\nabla$ \iif $f$ is anti-holomorphically harmonic \wrt to $\nabla^\bullet$.
\item[$\bullet$] $f$ is $\star$-stringy harmonic \wrt $\nabla$ \iif $f$ is anti-holomorphically harmonic \wrt to $\nabla^\star$.
\end{description}
\end{prop}
\proof Since $T^{2,0} + T^{0,2}\in\mal T$ and $T^{2,0}\in\mal T$, then $\tau^{\nabla^\bullet}(f)=\tau^{\nabla^\star}(f)=\tau^\nabla(f)$, according to proposition~\ref{prop-gene-A}. Moreover, we have
\begin{eqnarray*}
T^{\nabla^\star} & = & T - 2T^{2,0}= - T^{2,0} + T^{0,2} + T^{1,1}\\
T^{\nabla^\bullet} & = & T - 2(T^{2,0} + T^{0,2}) = - (T^{2,0} + T^{0,2}) + T^{1,1}
\end{eqnarray*}
so that
$$
J T^{\nabla^\star} = -J\star T \quad \text{ and } \quad J T^{\nabla^\bullet} = -J\cdot T,
$$
where, of course, $T$, $T^{\nabla^\bullet} $ and $T^{\nabla^\star}$ are respectively the torsions of $\nabla$, $\nabla^\bullet$ and $\nabla^\star$. We  conclude by using  proposition~\ref{holoharm} and definitions~\ref{defn-j-stringy} and \ref{defn-*j-stringy}. \comprf\hsq\medskip\\
Now, since $\nabla^\star J=0$, and according to theorem~\ref{thm-dbdf=0}, we have
\begin{cory}\index{dd derivative@$\bar\partial\partial$-derivative}
Let  $f\colon (L,j_L) \to (N,J,\nabla)$  be a map from a Riemann surface into an almost complex affine manifold.   Then $f$ is $\star$-stringy harmonic \iif 
$$
\hat\partial^{\widehat\nabla^\star}\overline{\hat\partial} f =0,
$$
i.e. $\overline{\hat\partial} f$ is an anti-holomorphic section of $T_{1,0}^*L\otimes_\C f^*TN$, endowed with the holomorphic structure defined by $\nabla^\star$.
\end{cory}
\proof This follows from proposition~\ref{prop-stringy-star-bullet} above, and the theorem~\ref{thm-dbdf=0}.\hsq
\begin{rmk}\em
We can also check the previous corollary by direct computations as in remark~\ref{rmk-comp-dbdf}.
\end{rmk}
\begin{rmk}\em
Let us remark that, in general, the metricity of the connection is not preserved when one passes from $\nabla$ to $\nabla^\star$ (resp. $\nabla^\bullet$). Indeed, if $\nabla$ is metric then $\nabla^\star$ is metric \iif $T^{2,0}$ is a 3-form \iif $T^{2,0}=T^{**}=0$, and therefore $T=T^{0,2}$. Therefore, if $\nabla J=0$, then $(N,J,h)$ is nearly K\"ahler. 
\end{rmk}
Let us conclude this subsection by the following:
\begin{prop}
Let  $f\colon (L,j_L) \to (N,J,\nabla)$ be a map from  a Riemann surface into an almost complex manifold with a connection $\nabla$. Then \begin{description}
\item[$\bullet$]  
$f$ is stringy harmonic \iif it is roughly harmonic \wrt
$$
\overline{\nabla}=f^*\nabla + \dfrac{1}{2}(J\cdot T)(df\circ j_L, .).
$$
\item[$\bullet$] 
And $f$ is $\star$-stringy harmonic \iif it is roughly harmonic \wrt
$$
\overline{\nabla}^\star =f^*\nabla + \dfrac{1}{2}(J\star T)(df\circ j_L, .).
$$
\end{description}
\end{prop}
\index{holomorphically harmonic|)}

\subsection{Bibliographical remarks and summary of the results.}
$\bullet$ The notion of holomorphically harmonic maps and stringy harmonic maps are new notions which generalize the notion of harmonic maps. In some sense these "generalized harmonic" maps are what corresponds to harmonic maps when the Levi-Civita connection is replaced by a metric connection with torsion. This fact is related to analoguous facts in mathematical physics: in superstring theory \cite{Strominger}  and in the study of non linear sigma models \cite{braden}. Moreover the new PDE of stringy harmonic maps could  be studied from the point of view of analysis. An interesting problem would be to see if the properties of harmonic maps (existence, regularity etc..) can be generalized to stringy harmonic maps. In particular, one knows that the existence of a variational intepretation for a PDE allows to use the techniques and the methods of the calculus of variations.\smallskip\\
$\bullet$ Let us point out several results concerning the properties of stringy harmonic maps, obtained in the present section: the interpretation in terms of the vanishing of some $\bar\partial\partial$-derivative, the equivalence between $\star$-stringy harmonicity and stringy harmonicity \wrt a new almost complex structure (under some conditions), the interpretation in terms of $J$-twisted harmonic maps, the fact that $J$-holomorphic curves are stringy harmonic (if $\nabla J=0$), the equivalence between stringy harmonicity and holomorphic harmonicity \wrt a new connection.
\smallskip\\
$\bullet$  One of our main result is the variational interpretation of stringy harmonicity in terms of a sigma model with a Wess-Zumino. The key point is in fact our result that if the characteristic connection of a $\mal G_1$-manifold has a parallel torsion then the $3$-form $H=J\cdot T$ is closed. This then implies the variational interpretation of stringy harmonicity. But we think that the closure of this $3$-form is itself an important (and of course a completely new) result. Indeed the existence of a closed 3-form on some manifold could have several implications on the geometry of this manifold, especially in low dimensions. Especially since our closed 3-form $H=J\cdot T$ is $\nabla$-parallel, where $\nabla$ is the characteristic connection. This implies that $H$ is invariant by the holonomy group of $\nabla$.\\ 
Moreover, in string theory II, the mathematical model involves in particular a 3-form $H$ which when it is closed then the solutions are called strong solutions \cite{srni,Friedrich-ivanov}.\smallskip\\
$\bullet$ Moreover, the variational interpretation of stringy harmonicity gives rise to two new contributions. Indeed, from the point of view of geometric analysis, this result provides a new class of geometric variational problems taking place in some class of almost Hermitian manifolds. From the point of view of mathematical physics, this provides a  class of geometric structures in which it is possible to define a non linear sigma models with a Wess-Zumino term.\smallskip\\
$\bullet$ Furthermore, another novelty, which is also an implication of the closure of the 3-form $H$, takes place in reductive homogeneous spaces. In these spaces, in general the torsion $T$ of the canonical connection (which is a 3-form) is not closed. Therefore, our result provides in this context a closed 3-form which can be used to redo the work of Agricola~\cite{agricola} about homogeneous models in string theory, by using the closed 3-form $H$ instead of the torsion $T$ which is not closed.\smallskip \\
$\bullet$ Finally, we provide a new contribution to the field of (integrable) non linear sigma models. Indeed we give new examples of integrable two-dimensional non linear sigma models. These new examples take place in some homogeneous spaces, namely in $(2k+1)$-symmetric spaces, which are not symmetric spaces. At our knowledge, all the already known integrable two-dimensional non linear sigma models take place in symmetric spaces or (equivalently) in Lie groups.
\smallskip \\
$\bullet$ Let us conclude by repeating that theorem~\ref{Unik} has  already been  proved  by  Friedrich-Ivanov \cite{Friedrich-ivanov} but without writting completely the proof. Indeed they  applied their general method to construct $G$-connections with skew-symmetric torsion to the case of $G_2$-manifolds and contact manifolds. Then they mention that this method also works for almost hermitian manifolds, and they also give a sketch of direct proof. However, our proof is different and is a direct proof. To our knowledge a completely written proof has never been given in the literature.

%%%%%%%%%%%%%%%%%%%%%%%%%%%%%%%%%%%%%%%%%%%%%%%%%%%%%%%%%%%%%%%%%%%%%%%%%%%%%%%%%%%%%%%%%%%%%%%%%
%                                                                                               %
%                                                                                               %
%                       Generalized harmonic maps into $f$-manifolds.                                                   %
%                                                                                               %
%                                                                                               %
%%%%%%%%%%%%%%%%%%%%%%%%%%%%%%%%%%%%%%%%%%%%%%%%%%%%%%%%%%%%%%%%%%%%%%%%%%%%%%%%%%%%%%%%%%%%%%%%%

\section{Generalized Harmonic maps into $f$-manifolds.}\label{gene-Harm-f-structure}
\index{f structure@$f$-structure|(}\index{f manifold@$f$-manifold|(}
\subsection{$f$-structures: General definitions and properties.}\label{f-structures-generalities}
\subsubsection{$f$-structures, Nijenhuis tensor and natural action on the space of torsions $\mal T$.}\label{subsec-f-struct}
\index{horizontal subbundle|(}
\index{vertical subbundle|(}
%
%\paragraph{Definitions}
%
Let us consider $(N,F)$  \emph{an $f$-manifold}, i.e. a manifold endowed with an $f$-structure (see definition~\ref{f-structure}). Let us set $\hor=\im F$ and $\ver=\ker F$, then we have $TN=\hor\oplus \ver$. If we put $P=-F^2$, then $P$ is the projector on $\hor$ along $\ver$.  Moreover $PF=FP=F$ and $F^2P=-P$. In particular, $\Bar J:=F_{|\hor}$ is a complex structure in the vector bundle $\hor$.\\
Let us denote also by $q:=\Id -P$ the projector on $\ver$ along $\hor$. We denote by $X=X^\ver + X^\hor$,  or sometimes simply by $X=X^v + X^h$, the decomposition of any element $X\in TN$.\medskip\\
In all the section~\ref{gene-Harm-f-structure}, we will consider the bundles $\hor^*$ and $\ver^*$ as well as all their tensor products respectively, as subbundles of $T^*N$ and $\otimes^k T^*N$, $k\in \mathbb{N}^*$, respectively. For example, for any trilinear form $B\in\mal C(\otimes^3T^*N)$, we will consider $B_{|\hor\times\ver\times\hor}$  as an element of $\mal C(T^*N^3)$  by identifying it to $B(P\cdot,q\cdot,P\cdot)$.\\
Moreover, we will often identify a $k$-linear map with its expression in terms of the vectors $(X_1,\ldots,X_k)\in TN^k$. For example, given $B\in\mal C(\otimes^3T^*N)$, we will write " let $\beta\in \hor^*\times\hor^*\times\ver^*$ be defined by $\beta=B(Z^v,X^h,Y^h)$" instead of " let $\beta\in \hor^*\times\hor^*\times\ver^*$ be defined by $\beta(X^h,Y^h,Z^v)=B(Z^v,X^h,Y^h)$, for all $X,Y,Z\in TN$".\medskip \\
\textbf{Notations} We extend  the notations and definitions of section~\ref{WZ-g-1} and the begining of subsection~\ref{WZ}, concerning there the complex bundle $(TN,J)$ (defined by a complex manifold $(N,J)$) to the complex bundle $(\hor,\Bar J)$, defined in the present section by the $f$-manifold $(N,F)$. Then all the algebraic results  of section~\ref{WZ-g-1} -like corollary~\ref{cory-B-eps,eps'} - can be extended to the complex bundle $(\hor,\Bar J)$.\footnote{Or to the Hermitian bundle $(\hor,\Bar J,h_{|\hor})$ if $(N,F)$ is endowed with a (compatible) metric $h$; see  defintion~\ref{F&h-compatible} below for a precise definition of a compatible metric.} 
\begin{defn}\label{def-nijenhuis -F}
 The Nijenhuis tensor $N_F$ of $F$ is defined by 
$$
N_F(X,Y)=[FX,FY] -F[FX,Y] -F[X,FY] -P[X,Y],
$$
where $X,Y\in\mal C(TN)$.
\end{defn}
Then we obtain immediately (\cite{Ishihara-Yano})
\begin{prop}\label{identities-N_F}
 We have the following identities.
$$
\begin{array}{rcccc}
N_F(qX,qY) & = & -P[qX,qY]  & = &  PN_F(qX,qY)\\
qN_F(X,Y) & = & q[FX,FY] & = & qN_F(pX,pY) \\
N_F(qX,PY) & = & -F[qX,FY]-P[qX,PY] & & 
\end{array}
$$
so that 
$$
N_{|\ver\times\ver}=\mR_\ver \quad \text{and} \quad N^\ver =- \mR_\hor(F\cdot, F\cdot),
$$ 
where $\mR_\ver$ and $ \mR_\hor$ are the curvature of $\ver$ and $\hor$ respectively (in the sense of definition~\ref{CurvaturPfaffian}). In particular, $N^\ver(\ver,\ver)=N^\ver(\hor,\ver)=\{0\}$ i.e
$$
N(\hor,\ver)\subset \hor \quad \text{and} \quad N(\ver,\ver)\subset \hor.
$$
Moreover ${N_F}_{|\ver\times\hor}={N_F^\hor}_{|\ver\times\hor}$ satisfies the following property
$$
N_F(X^v,\Bar J Y^h)=-\Bar J N_F(X^v,Y^h)
$$
i.e. $N_F(X^v,\cdot)_{|\hor}$ anticommutes with $\Bar J$.
\end{prop}
\begin{defn}\label{defn-B-eps-eps'}
Let $(N,F)$ be an $f$-manifold. Then for any $B\in \mal T$, we set 
$$
B^{\eps,\eps'}(X,Y)=-\dfrac{1}{4} \left(\eps\eps' B(FX,FY) +\eps FB(FX,Y) + \eps' FB(X,FY) -B(X,Y)\right).
$$ 
and $B^{2,0}:=B^{++}$, $B^{1,1}:=B^{+,-} + B^{-+}$ and $B^{0,2}:=B^{--}$.
\end{defn}
Then, setting $\Bar B= B_{|\hor^2}^\hor$, we have $PB^{\eps,\eps'}= \Bar B^{\eps,\eps'}$ or in other words $B^{\eps,\eps'}=\Bar B^{\eps,\eps'}  -\dfrac{1}{4} (B^\ver(F\cdot,F\cdot) -B^\ver) $, where the components $\Bar B^{\eps,\eps'}$ are computed with respect to the complex structure $\Bar J$ on $\hor$.\medskip\\
As for the case of an almost complex structure (section~\ref{TN-valued-2-forms}), we can define natural actions of $F$ on elements $B\in\mal T$:
\begin{eqnarray*}
F\cdot B & := & B(F\cdot,F\cdot,F\cdot):=-B^c\\
F\circact B & := & B(F\cdot,\cdot) + B(\cdot,F\cdot) -FB(\cdot,\cdot)\\
F\bullet B &  = & F\cdot B + \dfrac{1}{2}F\circact (B- \Bar B)\\
F\star B &  = & \dfrac{1}{2}\left( F\cdot B + F\circact B\right) .
\end{eqnarray*}
It is then important to remark that $F\cdot B=\Bar J\cdot\Bar B$, so that 
\begin{eqnarray*}
 F\bullet B &  = & \bar J\cdot \Bar B + \dfrac{1}{2}F\circact (B- \Bar B)\\
 F\star B &  = & \bar J\star \Bar B + \dfrac{1}{2}F\circact (B- \Bar B).
\end{eqnarray*}
Moreover, it is also useful to remark that $F\bullet B - F \circacti B = 4\Bar J\cdot \Bar B^{--} - \frac{1}{2}F \circacti (B- \Bar B)$, and   $F\bullet B - F\star B= 2\Bar J\cdot \Bar B^{--}$.
\subsubsection{Introducing a linear connection.}\label{introd-linear-connect}
Now, we introduce  a linear connection and want to compare the vertical component of the torsion with the vertical torsion.\smallskip\\
We obtain immediately the two followings properties.
\begin{prop}
Let $(N,F,\nabla)$ be an affine manifold endowed with a parallel $f$-structure $(\nabla F=0$). Then the subbundles $\hor=\im F$ and $\ver=\ker F$ are $\nabla$-parallel.
\end{prop}
\begin{prop}
Let $(N,\nabla)$ be an affine manifold. Let us suppose that we have a $\nabla$-parallel splitting $TN=\hor\oplus \ver$, where $\hor,\ver$ inherit the names of horizontal and vertical subbundles respectively. Then the vertical torsion  coincides with the vertical component of the torsion:
$$
T^v=T^\ver,
$$
where we use notations of section~\ref{phitorsion} for $T^v$, and the  notations defined above (in the begining of \ref{subsec-f-struct}) for $T^\ver$.
\end{prop}
In a more general context we can relate $T^v$ and $T^\ver$ as follows.
\begin{prop}\label{relationTv-TV}
Let $(N,\nabla)$ be an affine manifold. Let us suppose that we have some splitting $TN=\ver\oplus\hor$, where $\hor,\ver$ inherit the names of horizontal and vertical subbundles respectively. Then the vertical torsion and the vertical part of the torsion satisfy the following relations
$$
T_{|\ver\times\ver}^v=T_{|\ver\times\ver}^\ver  \quad \text{and} \quad T_{|\hor\wedge\ver}^v = T_{|\hor\wedge\ver}^\ver + \sigma^v.
$$
where $\sigma^v$ is the restriction to $\hor\wedge\ver$of the $\ver$-valued 2-form $\nabla^v q (X,Y) - \nabla^v q(Y,X)$.
\end{prop}
\textbf{Notation.} In all the next of the present section, we denote $\Phi=\mR_\hor$ the curvature of $\hor$.
\begin{defn}
The term $\mal R= T_{|\hor\wedge\ver}^v$ will be called the reductivity term.
\end{defn}
\begin{prop}\label{dec-Tv} 
Furthermore, in the situation off proposition~\ref{relationTv-TV}, we have the following equality:
$$
T^v= \Phi \oplus \mal R \oplus T^\ver(q\cdot,q\cdot)
$$
\end{prop}
\textbf{Proof of proposition~\ref{relationTv-TV} and \ref{relationTv-TV}.} This is a straightforward computation. \hsq \medskip\\
Let $X^h,Y^h\in\mal C(\hor)$. Then for any $f\in \mal C^\infty(N)$, we  have $\nabla_{X^h}^v(fY^h) = f\nabla_{X^h}^v Y^h + (X^h\cdot f) (Y^h)^v= f \nabla_{X^h}^v Y^h$ so that $\nabla_{X^h}^v Y^h$  defines  a bilinear map from $\hor\times\hor$ into $\ver$. Let $\Psi$ be its skew-symmetric part: $\Psi(H_1,H_2) = \nabla_{H_1}^vH_2 -\nabla_{H_2}^v H_1$. Then we have $T^\ver(H_1,H_2)=\nabla_{H_1}^vH_2 -\nabla_{H_2}^v H_1 -[H_1,H_2]^v$
i.e.
$$
T_{|\hor\times\hor}^\ver= \Psi + \Phi.
$$
Therefore,
\begin{prop} The following relation holds
$$
T^\ver= (\Psi + \Phi ) \oplus (\mal R - \sigma^v) \oplus T^\ver(q\cdot,q\cdot).
$$
Therefore $T^\ver=T^v$ \iif $\Psi=0$ and $\sigma^v=0$, which happens in particular if $\hor$  is $\nabla$-parallel. 
\end{prop}
%
%%%%%%%%%%%%%%%%%%%%%%%%%%%%%%%%%%%%%%%%%%%%%%%%%%%%%%%%%%%%%%%%%%%%%%%%%%%%%%%%%%%%%%%%%%%%%%%
%
\subsection{The $f$-connections and their torsion.}\label{f-connection&torsion} 
\subsubsection{Definition, notations and first properties.}\label{Def-not-Firstprop}
Let us come back to the case of an $f$-manifold $(N,F)$.
\begin{defn}
A linear connection $\nabla$, on an $f$-manifold $(N,F)$, which preserves the $f$-structure, i.e  $\nabla F=0$, is called an $f$-connection.
\end{defn}
Then we obtain easily: 
\begin{prop}\label{prop-T&N_F}
Let  $(N,F,\nabla)$ be an $f$-manifold endowed with an affine  $f$-connection. Then the torsion $T$ satisfies the following identity
$$
T(FX,FY) - FT(FX,Y) -FT(X,FY) -PT(X,Y)=-N_F(X,Y)
$$
\end{prop}
Therefore, we deduce the following.
\begin{cory}
Setting $N_{\Bar J}={N_F^\hor}_{|\hor\times\hor}$, the torsion $T$ (of an $f$-connection $\nabla$ on an $f$-manifolds $(N,F)$) satisfies the following identities:
\begin{eqnarray*}
N_{\Bar J} & = & 4\left( T_{\hor\times\hor}^\hor\right) _{\Bar J}^{0,2}\\
T^\hor(\Bar J X^h,Y^v) -\Bar J T^\hor(X^h,Y^v) & = & -\Bar J N_F^\hor(X^h,Y^v)=-\left( [\Bar J X^h,Y^v]^\hor -\Bar J[X^h,Y^v]^\hor\right)\\
T^\hor(X^v,\Bar J Y^h) -\Bar J T^\hor(X^v,Y^h) & = & -\Bar J N_F^\hor(X^v,Y^h)=-\left( [ X^v,\Bar J Y^h]^\hor -\Bar J[X^v,Y^h]^\hor\right)\\
T^\ver(X^h,Y^h) & = & -N_F^\ver(FX^h,FY^h)  =  \Phi(X,Y)\\
T^\hor(X^v,Y^v) & = & N_F^\hor(X^v,Y^v)  =  \mR_\ver(X^v,Y^v) 
\end{eqnarray*}
where $X,Y\in \mal C(TN)$. Consequently, the following component of the torsion $T_{|\hor\times\hor}^\ver$, $T_{|\ver\times\ver}^\hor$, $\left( T_{|\hor\times\hor}^\hor\right)^{0,2}$ and $\left[ T_{|\ver\times\hor}^\hor , \Bar J\right] $ are independent of $\nabla$.
\end{cory}
%Definition of $F\cdot B$, $F\circact B$, etc...
%
%
\paragraph{Introducing a metric}
After having  introduced a metric $h$ on $N$,  we want to characterize the metric connections $\nabla$ which preserves $F$. More precisely, we want to find a necessary and sufficient condition on the torsion $T$ for $\nabla$ metric to be a  $f$-connection.  In a first time, we will  begin by characterizing the metric connections which preserve the decomposition $TN=\hor\oplus\ver$, then in a second time we will introduce  the additionnal condition that the induced connection on $\hor$ preserves the complex structure $\Bar J$. \\
Let us define some notations. In the following, since a metric is given we use the convention defined in section \ref{TN-valued-2-forms}: each $TN$-valued bilinear form on $N$, $B\in \mal C(T^*N\otimes T^*N\otimes TN)$, will be identified (via the metric $h$) with the corresponding trilinear form. 
Moreover, we denote by $\Omega_A$ the bilinear form associated (via the metric $h$) to an endomorphism $A\in \mal C(\End(TN))$:
$$
\Omega_A(X,Y)=\langle A(X),Y\rangle,\quad \forall X,Y\in TN.
$$
Then, under our convention, for any endomorphism $A\in \mal C(TN)$, $\nabla^h A$ is identified to $\nabla^h\Omega_A$. Moreover, we set  
$$
\mrm{Sym}(B)(X,Y)=B(X,Y) + B(Y,X),\quad \forall X,Y\in TN,
$$
for all $B\in\mal C(T^*N\otimes T^*N\otimes TN)$.\\
Furthermore, let $E_1,E_2,E_3$ be vector bundles over $N$, then we set also
$\mal S (E_1\times E_2\times E_3)= \underset{i,j,k}{\mal S} E_i\otimes E_j\otimes E_k$, where we do a direct sum on the  circular permutation of 1,2,3.\\ 
Finally, to avoid any risk of confusion of the index "$h$" denoting the metric in the  notation of the Levi-Civita connection $\nabla^h$, with the same index in the notation for the horizontal component $X^h$ of a vector $X\in TN$, we will denote in all this section~\ref{gene-Harm-f-structure} the Levi-Civita connection by $D$:
$$
D:=\nabla^h.
$$
\subsubsection{Characterization  of  metric connections preserving the splitting.}
\begin{thm}\label{nabla-q=0}
 Let $(N,h)$ be a Riemannian manifold with an orthogonal decomposition $TN=\hor\oplus\ver$. Then a metric connection $\nabla$ leaves invariant this decomposition (i.e. $\hor$ and $\ver$ are $\nabla$-parallel) \iif its torsion $T$ satisfies 
$$
T_{|\hor\times\hor\times\ver}=\Phi, \qquad T_{|\ver\times\ver\times\hor}=\mR_\ver
$$
and
\begin{eqnarray*}
\mrm{Sym}_{\hor\times\hor}\left(T_{|\ver\times\hor\times\hor} \right)(X^v,Y^h,Z^h) & = & \mrm{Sym}_{\hor\times\hor}\left( D{\Omega_q}_{|\hor\times\hor\times\ver} \right)(Y^h,Z^h,X^v)\\
\mrm{Sym}_{\ver\times\ver}\left(T_{|\hor\times\ver\times\ver} \right)(X^h,Y^v,Z^v) & = & \mrm{Sym}_{\ver\times\ver}\left( D{\Omega_q}_{|\ver\times\ver\times\hor} \right)(Y^v,Z^v,X^h)
\end{eqnarray*}
In particular,  the components $\mrm{Sym}_{\hor\times\hor}\left(T_{|\ver\times\hor\times\hor} \right)$ and $\mrm{Sym}_{\ver\times\ver}\left(T_{|\hor\times\ver\times\ver} \right)$ are independent of $\nabla$.
\end{thm}
\textbf{Proof.}
According to remark~\ref{T+U}, we have $D=\nabla-A=\nabla -\frac{1}{2}(T+U)$. Therefore $\nabla q=0$ \iif 
\begin{equation}\label{condition-onT}
Dq= -[A,q].
\end{equation}
But $\langle [A,q](X,Y),Z\rangle=A(X,Y^v,Z)-A(X,Y,Z^v)=A(X^h,Y^v,Z^h) +A(X^v,Y^v,Z^h) -A(X^h,Y^h,Z^v)-A(X^v,Y^h,Z^v)$, and according to the characterization of the Levi-Civita connection:
\begin{eqnarray*}
-2\langle (D_{X^h}q)Y^h,Z^v\rangle & =  & 2 \langle D_{X^h} Y^h,Z^v\rangle = -Z^v\cdot \langle X^h,Y^h\rangle + \langle [X^h,Y^h],Z^v\rangle + \langle [Z^v,X^h],Y^h\rangle + \langle X^h,[Z^v,Y^h]\rangle \\
& = & -\Phi(X^h,Y^h,Z^v) + \langle -D_{Z^v}X^h + [Z^v,X^h],Y^h\rangle + \langle -D_{Z^v}Y^h + [Z^v,Y^h],X^h\rangle\\
& = & -\Phi(X^h,Y^h,Z^v) - \langle D_{X^h}Z^v,Y^h\rangle - \langle D_{Y^h}Z^v,X^h\rangle \\
& = & -\Phi(X^h,Y^h,Z^v) -  D\Omega_q(X^h,Z^v,Y^h) -  D\Omega_q(Y^h,Z^v,X^h)\\
& =  & -\Phi(X^h,Y^h,Z^v) -  \mrm{Sym}_{\hor\times\hor}(D\Omega_q(X^h,Y^h,Z^v).
\end{eqnarray*}
In  the last line we have used the fact that $D\Omega_q$ is symmetric \wrt the two last variables (since $q$ is a symmetric projector). Therefore the condition~(\ref{condition-onT})  restricted to $\hor\times\hor\times\ver$ is written $(T+U)_{|\hor\times\hor\times\ver}= \Phi + \mrm{Sym}_{\hor\times\hor}({D\Omega_q}_{|\hor\times\hor\times\ver})$, that is to say by identifiying respectively the symmetric and skew-symmetric part (\wrt the two first variables) of the two hand sides of this equality respectively, we obtain
\begin{multline}\label{nabla-q=0-hhv}
T_{|\hor\times\hor\times\ver}=\Phi \quad \mrm{and} \quad \mrm{Sym}_{\hor\times\hor}\left(T_{|\ver\times\hor\times\hor} \right)(X^v,Y^h,Z^h)=: U_{|\hor\times\hor\times\ver}(Y^h,Z^h,X^v)\\ =\mrm{Sym}_{\hor\times\hor}\left({D\Omega_q}_{|\hor\times\hor\times\ver}\right)(Y^h,Z^h,X^v).
\end{multline}
Moreover, since $D\Omega_q$ is symmetric \wrt the two last variables and $A$ is skew-symmetric  \wrt these two last variables, we see that the restriction to $\hor\times\hor\times\ver$ and to $\hor\times\ver\times\hor$ of the condition~(\ref{condition-onT}) are in fact equivalent. Furthermore, we have $\langle( D_{X^v}q)Y^h, Z^h\rangle =0$ and $\langle[A,q](X^v,Y^h),Z^h\rangle=0$. Therefore, the restriction to $\mal S(\hor\times\hor\times\ver)$ of the condition~(\ref{condition-onT}) is equivalent to (\ref{nabla-q=0-hhv}).\\
Proceeding in the same way as above, we obtain that the restriction to $\mal S(\ver\times\ver\times\hor)$ of (\ref{condition-onT}) is equivalent to 
$$
T_{|\ver\times\ver\times\hor}=\mR_\ver \quad \mrm{and} \quad \mrm{Sym}_{\ver\times\ver}\left(T_{|\hor\times\ver\times\ver} \right)(X^h,Y^v,Z^v) =  \mrm{Sym}_{\ver\times\ver}\left( D{\Omega_q}_{|\hor\times\ver\times\ver} \right)(Y^v,Z^v,X^h).
$$
Finally, we have $\langle( D_{X^h}q)Y^h,Z^h\rangle=0=\langle[A,q](X^h,Y^h),Z^h\rangle$ and $\langle( D_{X^v}q)Y^v,Z^v\rangle=0=\langle[A,q](X^v,Y^v),Z^v\rangle$. This completes the proof.\hfill$\square$%\medskip\\
\paragraph{Skew-symmetric torsion.}\index{skew-symmetric torsion}
Now, let us see under which   condition on the Riemannian manifold, there exists a connection preserving the splitting and with skew-symmetric torsion. It will turn out that  the existence of  a connection preserving the splitting and of which the horizontal component of the torsion $T_{|\hor^3}$ is skew-symmetric does not imposes any condition on the Riemannian manifold $(N,h)$, but the skew-symmetry of the other components $T_{|\mal S(\hor\times\hor\times\ver)}$ and $T_{|\mal S(\ver\times\ver\times\hor)}$ imposes constraints on the Riemannian manifold.
\begin{cory}\label{caract-H^2V}
Let $(N,h)$ be a Riemannian manifold with an orthogonal decomposition $TN=\hor\oplus\ver$. Then, the following statements are equivalent
\begin{description}
\item[(i)] There exists a metric connection $\nabla$ leaving invariant the decomposition $TN=\hor\oplus\ver$, such that the following component of the torsion $T_{|\mal S(\hor\times\hor\times\ver)}$ is skew-symmetric, i.e. $T_{|\mal S(\hor\times\hor\times\ver)}=\mrm{Skew}(\Phi)$.
\item[(ii)] For any metric connection $\nabla$ leaving invariant the decomposition $TN=\hor\oplus\ver$, the component of the torsion $T_{|\mal S(\hor\times\hor\times\ver)}$ is skew-symmetric.
\item[(iii)] ${D\Omega_q}_{|\hor\times\hor\times\ver}$ is skew-symmetric \wrt the two first variables, i.e. ${D q}_{|\hor\times\hor}\in\mal C((\Lambda^2\hor^*)\otimes \ver)$, or equivalently $D\Omega_q(P\cdot,P\cdot,q\cdot)=\frac{1}{2}\Phi$.
\item[(iv)] ${D\Omega_q}_{|\hor\times\ver\times\hor}$ is skew-symmetric \wrt the first and third variables, i.e. $D\Omega_q(PX,qY,PZ)=-\frac{1}{2}\Phi(Z,X,Y)$.
\item[(v)]  $\mrm{Skew}\left(D{\Omega_q}_{|\mal S(\hor\times\hor\times\ver)}\right)=0$.
\end{description}
We will then say that $(N,q,h)$ is of type $\mal H^2\mal V$.
\end{cory}
\begin{cory}\label{caract-V^2H}
 Let $(N,h)$ be a Riemannian manifold with an orthogonal decomposition $TN=\hor\oplus\ver$. Then, the following statements are equivalent.
\begin{description}
\item[(i)] There exists a metric connection $\nabla$ leaving invariant the decomposition $TN=\hor\oplus\ver$, such that the following component of the torsion $T_{|\mal S(\ver\times\ver\times\hor)}$ is skew-symmetric, i.e. $T_{|\mal S(\ver\times\ver\times\hor)}=\mrm{Skew}(\mR_\ver)$.
\item[(ii)] For any metric connection $\nabla$ leaving invariant the decomposition $TN=\hor\oplus\ver$, the component of the torsion $T_{|\mal S(\ver\times\ver\times\hor)}$ is skew-symmetric.
\item[(iii)] ${D\Omega_q}_{|\ver\times\ver\times\hor}$ is skew-symmetric \wrt the two first variables, i.e. ${D q}_{|\ver\times\ver}\in\mal C((\Lambda^2\ver^*)\otimes \hor)$, or equivalently $D\Omega_q(q\cdot,q\cdot,P\cdot)=-\frac{1}{2}\mR_\ver$.
\item[(iv)] ${D\Omega_q}_{|\ver\times\hor\times\ver}$ is skew-symmetric \wrt the first and third variables, i.e. $D\Omega_q(qX,PY,qZ)=\frac{1}{2}\mR_\ver(Z,X,Y)$.
\item[(v)]  $\mrm{Skew}\left(D{\Omega_q}_{|\mal S(\ver\times\ver\times\hor)}\right)=0$.
\end{description}
We will  then say that $(N,q,h)$ is of type $\mal V^2\mal H$.
\end{cory}
\textbf{Proof of corollary~\ref{caract-H^2V} and \ref{caract-V^2H}.}
We have seen in the proof of theorem~\ref{nabla-q=0} that for any metric connection $\nabla$, we have $\nabla {\Omega_q}_{|\hor^3}=0$ and $\nabla {\Omega_q}_{|\ver^3}=0$, so that we have $\nabla q=0$ \iif ${\nabla q}_{|\mal S(\hor\times\hor\times\ver)}=0$ and ${\nabla q}_{|\mal S(\ver\times\ver\times\hor)}=0$, which is equivalent respectively to the conditions on $T_{|\mal S(\hor\times\hor\times\ver)}=0$ and $T_{|\mal S(\ver\times\ver\times\hor)}=0$ respectively, described by theorem~\ref{nabla-q=0}. In particular, we see that there always exists metric connection $\nabla$ leaving invariant the decomposition of $TN$, which provides  us the implication (ii) $\Rightarrow$ (i). Moreover, a necessary condition for (i) is $\mrm{Sym}_{\hor\times\hor}\left( D{\Omega_q}_{|\hor\times\hor\times\ver} \right)=0$, and respectively $\mrm{Sym}_{\ver\times\ver}\left( D{\Omega_q}_{|\ver\times\ver\times\hor} \right)=0$. Conversely if this condition is satisfied, then according to  theorem~\ref{nabla-q=0}, (ii) is also satisfied (since $T\in\mal T$, then the skew-symmetry of  $T_{|\mal S(\hor\times\hor\times\ver)}$ (resp. $T_{|\mal S(\ver\times\ver\times\hor)}$) is equivalent to $\mrm{Sym}_{\hor\times\hor}\left(T_{|\ver\times\hor\times\hor} \right)=0$, resp. $\mrm{Sym}_{\ver\times\ver}\left( T_{|\hor\times\ver\times\ver} \right)=0$). Therefore we have proved the sequence of implications (ii) $\Rightarrow$ (i) $\Rightarrow$ (iii) $\Rightarrow$ (ii), i.e. (i), (ii) and (iii) are equivalent. Concerning the equivalent reformulation of (i) and (iii), the former follows from theorem~\ref{nabla-q=0} and the latter from the fact that according to the proof of theorem~\ref{nabla-q=0}, the skew-symmetric part (\wrt the two first variables $X,Y$) of ${D\Omega_q}_{|\hor\times\hor\times\ver}$ (resp. ${D\Omega_q}_{|\ver\times\ver\times\hor}$) is $\frac{1}{2}\Phi$ (resp. $-\frac{1}{2}\mR_\ver$). The equivalence between (iii) and (iv) follows from the symmetry of $D\Omega_q$ \wrt the two last variables. Finally, using the computation done in the proof of compute theorem~\ref{nabla-q=0}, we  can compute
\begin{eqnarray*}
 \mrm{Skew}\left(D{\Omega_q}_{|\mal S(\hor\times\hor\times\ver)}\right) & = & \dfrac{1}{2}\left(\mrm{Skew}(\Phi) + \mrm{Skew}(U(P\cdot,P\cdot,q\cdot))\right) -\dfrac{1}{2}\left(\mrm{Skew}(\Phi) - \mrm{Skew}(U(P\cdot,P\cdot,q\cdot))\right)\\
& = & \mrm{Skew}(U(P\cdot,P\cdot,q\cdot))
\end{eqnarray*}
and idem for $\mal S(\ver\times\ver\times\hor)$. This completes the proof.\hfill $\square$
\begin{rmk}\label{rmk-extend-red}\em
According to the previous proof, we see that if $(N,q,h)$ is of type $\mal H^2\mal V$, then for any extension $T\in\mal C(\mal T)$ of the skew-symmetric trilinear form $\mrm{Skew}(\Phi)\in\mal C\left( \mal S(\hor^*\times\hor^*\times\ver^*) \right) $, the corresponding metric connection $\nabla$ satisfies ${\nabla q}_{|\mal S(\hor\times\hor\times\ver)} =0$. In the same way, if $(N,q,h)$ is of type $\mal V^2\mal H$, then for any extension $T\in\mal C(\mal T)$ of the skew-symmetric trilinear form $\mrm{Skew}(\mR_\ver)\in\mal C\left( \mal S(\ver^*\times\ver^*\times\hor^*) \right) $, the corresponding metric connection $\nabla$ satisfies ${\nabla q}_{|\mal S(\ver\times\ver\times\hor)} =0$.
\end{rmk}
\begin{defn}\index{reductive $f$-manifold}
We will say that the (orthogonal) decomposition on the Riemannian manifold $(N,h)$ is \textbf{reductive} (\wrt the metric $h$) or that $(N,q,h)$ is \textbf{reductive} if $(N,q,h)$ is of type $\mal H^2\mal V$ and of type $\mal V^2\mal H$. This is equivalent to say that there exists a metric connection $\nabla$ leaving invariant the decomposition $TN=\hor\oplus\ver$, and with skew-symmetric torsion.
\end{defn}
\begin{prop}
$(N,q,h)$ is reductive \iif the trilinear map
$$
\mrm{Skew}({D\Omega_q}_{|\hor\times\hor\times\ver}) \oplus \mrm{Skew}({D\Omega_q}_{|\ver\times\ver\times\hor}) 
$$
is skew-symmetric.
\end{prop}
\textbf{Proof.} An element $\alpha\in\mal C(\hor^*\otimes\hor^*\otimes\ver^*)$ satisfies $\mrm{Skew}(\alpha)\in\mal C(\Lambda^2\hor^*)\wedge\ver^*)$ \iif $\alpha(X^h,Y^h,Z^v)$ is skew-symmetric \wrt $(X^h,Y^h)$. \hfill $\square$
\begin{prop}\label{D^v-nabla^c}
Let $(N,q,h)$ be a reductive Riemannian manifold. Let us suppose that is given some metric connection $\nabla^c$ on $\ver$.  Then there exists a metric connection $\nabla$ on $N$ preserving the decomposition $TN=\ver\oplus\hor$, with skew-symmetric torsion, and which coincides with $\nabla^c$ on $\ver$ \iif 
$$
(D^v-\nabla^c)(X^h,Y^v,Z^v)=-\dfrac{1}{2}\mR_\ver(Y^v,Z^v,X^h) \quad  \text{and} \quad(D^v-\nabla^c)_{|\ver^3}\in \mal C(\Lambda^3\ver^*).
$$
\end{prop}
\textbf{Proof.} Let us suppose that such a metric connection $\nabla$ exists. Then we have $D=\nabla -\frac{1}{2}T$ and thus $(D^v-\nabla^c)_{|\ver^3}\in \mal C(\Lambda^3\ver^*)$ and $(D^v-\nabla^c)(X^h,Y^v,Z^v)=-\dfrac{1}{2}T(X^h,Y^v,Z^v)=-\dfrac{1}{2}T(Y^v,Z^v,X^h)=-\dfrac{1}{2}\mR_\ver(Y^v,Z^v,X^h)$, since $T$ is skew-symmetric and $\nabla$ leaves invariant the decomposition of $TN$.\\
Conversely, if $\nabla^c$ satisfies the above conditions , then let $\beta\in\mal C(\Lambda^3\hor^*)$ and $\alpha= (D^v-\nabla^c)_{|\ver^3}\in \mal C(\Lambda^3\ver^*)$. Let us consider the 3-form 
$$
T=\beta \oplus \mrm{Skew}(\Phi) \oplus \mrm{Skew}(\mR_\ver)\oplus \alpha,
$$
as well the corresponding metric connection $\nabla=D +\frac{1}{2}T$. Then, since $(N,q,h)$ is reductive, according to theorem~\ref{nabla-q=0}, $\nabla$ preserves the decomposition of $TN$. Moreover, by definition of $\nabla$, we have
$(\nabla-\nabla^c)_{TN\times\ver\times\ver}=0$. This completes the proof.\hfill$\square$%\medskip\\
\paragraph{A useful additionnal property.}
Let us add the following characterization of the type $\mal \ver^2\mal\hor$ in terms of the vertical torsion of the Levi-Civita connection.
\begin{prop}\label{V^2H-Tv}
Let $(N,h)$ be a Riemannian manifold with an orthogonal decomposition $TN=\hor\oplus\ver$. Let $T^v$ be the vertical torsion  of the Levi-Civita connection. Then, $(N,q,h)$ is of type $\mal \ver^2\mal\hor$ \iif 
$$
\mrm{Sym}_{\ver\times\ver}(T_{|\hor\times\ver\times\ver}^v)=0.
$$
\end{prop}
\textbf{Proof}
Let $H\in\hor$, $V_1,V_2\in\ver$. Then 
\begin{eqnarray*}
T^v(H,V_1,V_2) + T^v(H,V_2,V_1) & = &
\langle D_H^v V_1 -[H,V_1]^v,V_2\rangle + \langle D_H^v V_2 -[H,V_2]^v,V_1\rangle\\
& = & \langle D_H V_1 -[H,V_1],V_2\rangle + \langle D_H V_2 -[H,V_2],V_1\rangle\\
& = & \langle D_{V_1}H,V_2\rangle + \langle D_{V_2}H,V_1\rangle \\
& = & \mrm{Sym}_{\ver\times\ver}(DP)(V_1,H,V_2)= \mrm{Sym}_{\ver\times\ver}(DP)(V_1,V_2,H)\\
& =  & -\mrm{Sym}_{\ver\times\ver}(Dq)(V_1,V_2,H).
\end{eqnarray*}
Then we completes the proof by applying corollary~\ref{caract-V^2H}\medskip.\hfill $\square$ 

\subsubsection{Characterization of metric  $f$-connections. Existence of a characteristic connection.}\label{Character-metic-f-connect}
\index{connection!characteristic|(}\index{global type $\mal G_1$|(}
Now, let us come back to the case of an $f$-manifolds. Then the condition $\nabla F=0$ is equivalent to the fact that $\nabla$ leaves invariant the decomposition $TN=\hor\oplus \ver$ and moreover $\nabla^\hor \Bar J=0$, where $\nabla^\hor$ is the connection induced by $\nabla$ on $\hor$.  Heuristically, we have to add to the conditions of theorem~\ref{nabla-q=0} those of theorem~\ref{Gauduchon} as well as the condition $\nabla_{\ver}^\hor \Bar J=0$.
\begin{defn}\label{F&h-compatible}
We will say that an $f$-structure $F$ and  a metric $h$ on a manifold $N$ are compatible if $\hor\perp\ver$ and if $\Bar J$ is an orthogonal complex structure on $\hor$ endowed with the metric induced by $h$. This is equivalent to say that $F$ is skew-symmetric \wrt the metric $h$: $F\in\so(TN)$, or equivalently that $I=\Bar J\oplus \Id_\ver$ is orthogonal: $I^*h=h$. We will then say that  $(N,F,h)$ is a metric $f$-manifold.
\end{defn}
\begin{thm}\label{thm-nablaOmega_F=0}
Let $(N,F,h)$ be a metric $f$-manifold. Then a metric  connection $\nabla$ preserves the $f$-structure $F$ \iif all the following statements hold:
$$
\begin{array}{lcl}
{\nabla \Omega_F}_{|\hor^3}=0 & \Longleftrightarrow  & N_{\Bar J}=4\left(T_{|\hor^3} \right)^{0,2}\text{ and } \mrm{Skew}\left( \left( T_{|\hor^3} \right)^{2,0} - \left( T_{|\hor^3}\right)^{1,1} \right) = (d^c{\Omega_F}_{|\hor^3})^{**}\\
{\nabla \Omega_F}_{|\mal S(\hor\times\hor\times \ver)} = 0 & \Longleftrightarrow  & 
\begin{array}{lcl} 
\left\lbrace \begin{array}{l} \negthickspace\negthickspace \boxed{ \begin{array}{l} {\nabla \Omega_F}_{|\hor\times\hor\times \ver}=0\\
{\nabla \Omega_F}_{|\hor\times\ver\times \hor}=0 \end{array}} \\
{\nabla \Omega_F}_{|\ver\times\hor\times\hor }=0\end{array}\right.
                      &   \Longleftrightarrow   &   \left\lbrace \begin{array}{l}
\boxed{{\nabla \Omega_q}_{|\mal S(\hor\times\hor\times\ver)}=0}\\
  \begin{array}{c} T(X^v,FY^h,Z^h) + T(X^v, Y^h,FZ^h)=\\ N(X^v, Y^h, FZ^h)\end{array}
\end{array}\right.
\end{array}  \\
{\nabla \Omega_F}_{|\mal S(\ver\times\ver\times \hor)} = 0 & \Longleftrightarrow  &  {\nabla \Omega_q}_{|\mal S(\ver\times\ver\times \hor)} = 0
\end{array}
$$
\end{thm}
\textbf{Proof.} We first notice that $\nabla\Omega_F(X^h,Y^h,Z^v)=\langle (\nabla_{X^h}F)Y^h,Z^v \rangle=
\langle \nabla_{X^h}(FY^h),Z^v \rangle= -\nabla\Omega_q(X^h,FY^h,Z^v)$. Therefore  ${\nabla\Omega_F}_{|\hor\times\hor\times \ver}=0 \Leftrightarrow {\nabla\Omega_q}_{|\hor\times\hor\times \ver}=0$. In the same way, 
${\nabla\Omega_F}_{|\hor\times\ver\times \hor}=0 \Leftrightarrow {\nabla\Omega_q}_{|\hor\times\ver\times \hor}=0$. Therefore, according to the proof of corollary~\ref{caract-H^2V}, we have proved the equivalence between the two boxes. \\
Now, let us compute
\begin{multline*}
T(X^v,FY^h,Z^h) + T(X^v,Y^h,FZ^h)  =  \langle \nabla_{X^v} (FY^h) - \nabla_{FY^h}X^v -[X^v,FY^h], Z^h\rangle \\
+ \langle \nabla_{X^v} Y^h - \nabla_{Y^h} X^v -[X^v,Y^h], FZ^h\rangle\\
  =  \nabla\Omega_F(X^v,Y^h,Z^h) - \nabla\Omega_q(FY^h,X^v,Z^h) 
   - \nabla\Omega_q(Y^h,X^v,FZ^h) +  N_F(X^v,Y^h,FZ^h)
\end{multline*}
so that if ${\nabla \Omega_q}_{|\mal S(\hor\times\hor\times\ver)}=0$ then we obtain the equivalence $({\nabla \Omega_F}_{|\ver\times\hor\times\hor }=0) \Leftrightarrow (T(X^v,FY^h,Z^h) + T(X^v, Y^h,FZ^h)= N(X^v, Y^h, FZ^h)$. \\
We have also to compute that $\nabla\Omega_F(X^v,Y^v,Z^h)=\langle (\nabla_{X^v}F)Y^v,Z^h\rangle =-\langle F\nabla_{X^v} Y^v, Z^h\rangle = \nabla q(X^v,Y^v,FZ^h)$,  also that $\nabla\Omega_F(X^v,Y^h,Z^v) = -\nabla q(X^v,FY^h,Z^v)$, and that $\nabla\Omega_F(X^h,Y^v,Z^v)=0$. Therefore, according to the proof of corollary~\ref{caract-V^2H}, we have proved the equivalence $({\nabla \Omega_F}_{|\mal S(\ver\times\ver\times \hor)} = 0 ) \Leftrightarrow    ({\nabla \Omega_q}_{|\mal S(\ver\times\ver\times \hor)} = 0)$.\\
Furthermore, we can prove (by a straightforward but a little bit long computation) the following formula %(see section~\ref{appendix-f-structure}) 
generalizing equation~(\ref{nablah-Omega-J}):
\begin{multline}\label{nabla-Omega-F}
2D\Omega_F(X,Y,Z)  = 
d\Omega_F(X,Y,Z) - d\Omega_F(X,FY,FZ) + N(Y,Z,FX) \\
+ \Phi(Y,FZ,X) + \Phi(FY,Z,X) - \Phi(X,FY,Z) -\Phi(FZ,X,Y) \\
+ \mrm{Sym}_{\ver\times\ver}({D \Omega_F}_{|\ver\times\hor\times\ver }) + \mrm{Sym}_{\ver\times\ver}({D \Omega_F}_{|\ver\times\ver\times\hor}).
\end{multline}
We deduce from this, that 
$$
2D\Omega_F(X^h,Y^h,Z^h)  = 
d\Omega_F(X^h,Y^h,Z^h) - d\Omega_F(X^h,FY^h,FZ^h) + N(Y^h,Z^h,FX^h).
$$ 
Moreover, $\nabla\Omega_{|\hor^3}=0$ \iif 
$$
{D\Omega_F}_{|\hor^3}=-\dfrac{1}{2}(T+U)_{|\hor^3}.
$$
Therefore, we can proceed as in the proof of theorem~\ref{Gauduchon} to prove the  equivalence concerning the restriction to $\hor^3$.
This completes the proof.
\medskip\hfill$\square$%\\
\begin{rmk}\label{rmk-admissible}\em
A priori, we could  think that we can deduce from this theorem that in a metric $f$-manifold, it could not exist any metric $f$-connection. Indeed, we see that the condition $T(X^v,FY^h,Z^h) + T(X^v,Y^h,FZ^h)= N_F(X^v,Y^h,FZ^h)$ is compatible with the condition ${\nabla \Omega_q}_{|\mal S(\hor\times\hor\times\ver)}=0$ \iif
\begin{equation}\label{eq-admissible}
\mrm{Sym}_{\hor\times\hor}(N_{F|\ver\times\hor\times\hor})= 2 \left[ \mrm{Sym}_{\hor\times\hor} ({D \Omega_q}_{|\hor\times\hor\times\ver})\right]^{(2,0) + (0,2)}(Y^h,Z^h,X^v) 
\end{equation}
according to theorem~\ref{nabla-q=0}. In other words, the existence of some metric connection $\nabla$ such that ${\nabla \Omega_F}_{|\mal S(\ver\times\hor\times\hor) }=0$ is possible \iif   equation~(\ref{eq-admissible}) holds. Nevertheless, the existence of some metric connection $\nabla$ such that ${\nabla F}_{|\hor^3}=0$ and ${\nabla F}_{|\mal S(\hor\times\ver\times\ver)}=0$ always holds without condition. But in fact, equation~(\ref{eq-admissible}) also always holds.
\end{rmk}
\begin{lemma}\label{lemma-admissible}
In a metric $f$-manifold $(N,F,h)$, the following identity holds   
$$
\mrm{Sym}_{\hor\times\hor}(N_{F|\ver\times\hor\times\hor})= 2 \left[ \mrm{Sym}_{\hor\times\hor} ({D \Omega_q}_{|\hor\times\hor\times\ver})\right]^{(2,0) + (0,2)}(Y^h,Z^h,X^v) .
$$
Any  metric $f$-manifold $(N,F,h)$ admits a metric $f$-connection.
\end{lemma}
\proof This can be checked directly by computation (use the characterization of the Levi-Civita connection, as in the proof of theorem~\ref{nabla-q=0}). \hsq 
\begin{rmk}\em
In a metric $f$-manifold, there are several constraints on the component $T_{|\ver\times\hor\times\hor}$ of the torsion of some metric $f$-connection  $\nabla$: the components
$\left( T_{|\ver\times\hor\times\hor}\right)^{(2,0)+(0,2)}$ and $\mrm{Sym}_{\hor\times\hor}\left( \left[ T_{|\ver\times\hor\times\hor} \right]^{1,1} \right) $ are determined. The only degree of freedom is on $ \mrm{Skew}_{\hor\times\hor}\left( \left[T_{|\ver\times\hor\times\hor} \right]^{1,1}\right) $ which can be chosen freely. Therefore, this is also the only degree of freedom of $T_{|\mal S(\ver\times\hor\times\hor) }$, according to theorem~\ref{nabla-q=0}.
\end{rmk}
\begin{rmk}\em
A compatible $f$-structure $F$ on a Riemannian manifold $(N,h)$ defines a reduction of the orthogonal frame bundle of $TN$ (which is of course a $O(n)$-bundle) to some $U(2p)\times O(n-2p)$-bundle, where $n=\dim N$ and $2p$ is the rank of $F$. Therefore a metric $f$-connection is nothing but a connection on this $U(2p)\times O(n-2p)$-bundle, which always exists. This provides a new proof of lemma~\ref{lemma-admissible}.
\end{rmk}
\begin{rmk}\label{N-f-skew} \em
 If $(N,F,h)$ is of type $\hor^2\ver$, then   $N_{F|\ver\times\hor\times\hor}$ is skew-symmetric \wrt the two last variables. 
\end{rmk}
%
%%%%%%%%%%%%%%%%%%%%%%%%%%%%%
%
\paragraph{Skew-symmetric torsion.}\index{skew-symmetric torsion}
Further, we are interested by metric $f$-connections with skew-symmetric torsion. As we have done above we have to study first the condition of skew-symmetry on each component of the torsion and then to group all the obtained conditions to obtain a global condition on the metric $f$-manifold for the existence of metric $f$-connection with skew-symmetric torsion.
\begin{defn}\label{def-Nijenhuis-extended}\index{Nijenhuis tensor!extended}
Let $(N,F,h)$ be a metric $f$-manifold. We define the \textbf{extended Nijenhuis tensor} $\tl N_F$ as the $TN$-valued 2-form on $N$ (whose corresponding trilinear map is) defined by
$$
\tl N_F:= N_F + \Phi + \mR_\ver(Z^v,X^v,Y^h) + \mR_\ver(Y^v,Z^v,X^h).
$$
We remark that $\tl N_{F|\mal S (\ver\times\ver\times\hor)}=\mrm{Skew}(\mR_\ver)$ is always skew-symmetric.
\end{defn}
\begin{prop}\label{prop-N-f-hhv}
 Let $(N,F,h)$ be a metric $f$-manifold. Then the following statements are equivalent.
\begin{description}
\item[(i)] There exists a metric $f$-connection $\nabla$ (satisfying then $\nabla F=0$) with a torsion $T$ such that ${T^{0,2}}_{|\mal S(\hor\times\hor\times \ver)}$ is skew-symmetric.
\item[(ii)] There exists a metric connection $\nabla$, satisfying ${\nabla F}_{|\mal S(\hor\times\hor\times \ver)}=0,$ with a torsion $T$ such that ${T^{0,2}}_{|\mal S(\hor\times\hor\times \ver)}$ is skew-symmetric.
\item[(iii)] $N_F(FY^h,Z^h,X^v) + N_F(Y^h,FZ^h,X^v)= N_F(X^v, Y^h, FZ^h)$.
\item[(iv)] The extented Nijenhuis tensor $\tl N_F$ satisfies: $\tl N_{F|\mal S(\hor\times\hor\times \ver)} $ is skew-symmetric.

\end{description}
\end{prop}
\textbf{Proof.}
Since the condition "${T^{0,2}}_{|\mal S(\hor\times\hor\times \ver)}$ is skew-symmetric" concerns only the subspace $\mal S(\hor\times\hor\times \ver)$, then (i) and (ii) are equivalent. Moreover,  we have
$$ 
\tl N_{F|\mal S(\hor\times\hor\times \ver)} = N_{F|\mal S(\hor\times\hor\times \ver)} - N_F(F\cdot,F\cdot, q\cdot)
$$
(according to the definition of $\tl N_F$ and proposition~\ref{identities-N_F}), and this equality gives easily the equivalence between (iii) and (iv). Now, it remains to prove the equivalence between (i) and (iv).\\
Let us first recall that
$$
-4 (T^{0,2}):= T(F\cdot,F\cdot,q\cdot) + T(F\cdot,q\cdot, F\cdot) + T(q\cdot,F\cdot,F\cdot) - T
%_{|\mal S(\hor\times\hor\times \ver)}.
$$
Let us suppose (i). Then we have according to proposition~\ref{prop-T&N_F}, that
$T(F\cdot,F\cdot,q\cdot) + T(F\cdot,q\cdot, F\cdot) + T(q\cdot,F\cdot,F\cdot) - T(\cdot,\cdot,P\cdot)=-N_F$ 
and thus
$$
-4T^{0,2} = -N_F - T(\cdot,\cdot,q\cdot)= -N_F -\Phi -T_{|\ver^3} -T_{|(\ver\wedge\hor)\otimes\ver}
$$
and hence
$-4T_{|\mal S(\hor\times\hor\times \ver)}^{0,2}= - {\tl N}_{F|\mal S(\hor\times\hor\times \ver)}$, which implies (iv).\\
Conversely, let us suppose (iv). Then we have to construct a trilinear map $T_{|\mal S(\hor\times\hor\times \ver)}$, on $\mal S(\hor\times\hor\times \ver)$, such that \medskip\\
$(C_1)$:  $-4T_{|\mal S(\hor\times\hor\times \ver)}^{0,2}$ is a 3-form, and\\
$(C_2)$:  the following conditions of theorem~\ref{nabla-q=0} hold

$$
T_{|\hor\times\hor\times \ver}=\Phi, \quad \mrm{Sym}_{\hor\times\hor}\left(T_{|\ver\times\hor\times\hor} \right) =  \mrm{Sym}_{\hor\times\hor}\left( D{\Omega_q}_{|\ver\times\hor\times\hor} \right)(Y^h,Z^h,Z^v).
$$
Indeed, for any trilinear map $T_{|\mal S(\hor\times\hor\times \ver)}\in\mal C(\mal T)$, let us set $T_{|\ver\times\hor\times\hor}=S+A$, where $S,A$ are resp. symmetric and skew-symmetric \wrt the two last variables. Then we can prove easily that a necessary  condition on $S,A$  so that $T_{|\mal S(\hor\times\hor\times \ver)}$ satisfies  the required conditions $(C_1)$-$(C_2)$ above, is
\begin{subequations}\label{eq-A-S}
\begin{eqnarray}
S & = & \mrm{Sym}_{\hor\times\hor}\left( D{\Omega_q}_{|\hor\times\hor\times\ver} \right)(Y^h,Z^h,X^v) \\
S^{(2,0)+(0,2)} & = & 0 \label{eq-A-S-b}\\
A^{(2,0)+(0,2)}(X^v,Y^h,Z^h) & = & \Phi^{(2,0)+(0,2)}(Y^h,Z^h,X^v)
\end{eqnarray}
\end{subequations}
We then see that the two equations on $S$ are compatible according to (iv) and lemma~\ref{lemma-admissible}, and defines then uniquely $S$. Moreover the equation on $A$ implies  $ T(X^v,FY^h,Z^h) + T(X^v, Y^h,FZ^h)=\\ N(X^v, Y^h, FZ^h)$, according to (iv) (written in the form of (iii)) and \eqref{eq-A-S-b}. Therefore, we see that, according to theorem~\ref{thm-nablaOmega_F=0}, the conditions $(C_1)$-$(C_2)$ on $T_{|\mal S(\hor\times\hor\times \ver)}$ implies that ${\nabla F}_{|\mal S(\hor\times\hor\times \ver)}=0$, for any  metric connection $\nabla$ whose the torsion $T\in \mal C(\mal T)$ is an extension of $T_{|\mal S(\hor\times\hor\times \ver)}$. \\
Now, we can verify easily that any trilinear map $T_{|\mal S(\hor\times\hor\times \ver)}\in\mal T$ defined by  $T_{|\hor\times\hor\times \ver}=\Phi$ and $T_{|\ver\times\hor\times\hor}=S+A$ with $A,S$ satisfying \eqref{eq-A-S}, satisfies the required conditions $(C_1)$-$(C_2)$. We have then proved (ii). This completes the proof.\hfill$\square$ 
\begin{prop}\label{H^2V-redu-typ-G_1}
Let $(N,F,h)$ be a metric $f$-manifold. Then the following statements are equivalent.
\begin{description}
\item[(i)] There exists a metric $f$-connection $\nabla$ (satisfying then $\nabla F=0$) with a torsion $T$ such that
$T_{|\mal S(\hor\times\hor\times \ver)}$ is skew-symmetric.
\item[(ii)] There exists a metric connection $\nabla$, satisfying ${\nabla F}_{|\mal S(\hor\times\hor\times \ver)}=0,$ with a torsion $T$ such that $T_{|\mal S(\hor\times\hor\times \ver)}$ is skew-symmetric.
\item[(iii)] $(N,q,h)$ is of type $\mal H^2\mal V$, and $\tl N_{F|\mal S(\hor\times\hor\times\ver)}$ is skew-symmetric.
\end{description} 
Furthermore, under these statements, for any such connection satisfying (i) or (ii), then $T_{|\mal S(\hor\times\hor\times \ver)}$ is unique (i.e. uniquely determined by the metric $f$-manifold $(N,F,h)$) and equal to $\mrm{Skew}(\Phi)$. Conversely, any extension $T\in\mal C(\mal T)$ of this unique skew-symmetric trilinear form $\mrm{Skew}(\Phi)$ defines a metric connection $\nabla$, satisfying ${\nabla F}_{|\mal S(\hor\times\hor\times \ver)}=0$. 
\end{prop}
\textbf{Proof.} (i) and (ii) are equivalent for the same reason as in the proof of the previous proposition. Moreover, (i) and (iii) are equivalent according to theorem~\ref{thm-nablaOmega_F=0} and corollary~\ref{caract-H^2V}. Moreover, by skew-symmetry and theorem~\ref{nabla-q=0}, we have $T_{|\mal S(\hor\times\hor\times\ver)}=\mrm{Skew}(\Phi)$. This completes the proof.\medskip\hfill$\square$\\
We are led to the following definition.
\begin{defn}
We will say that a  metric $f$-manifold $(N,F,h)$ is  reductive if  $(N,q,h)$ is reductive, where $q$ is defined by $F$.\\
We will say that a  metric $f$-manifold $(N,F,h)$ is \emph{ reductively of  type $\mal G_1$} if $\tl N_{F|\mal S(\hor\times\hor\times\ver)}$ is skew-symmetric.
\end{defn}
Now, let us turn ourself on the horizontal component $\hor^3$ of $TN^3$. 
\begin{defn}
We will say that a metric $f$-manifold is \emph{horizontally of type $\mal G_1$} or that it is \emph{of horizontal type $\mal G_1$} if one the following equivalent statements holds.
\begin{description}
\item[(i)] The \emph{horizontal Nijenhuis tensor} $N_{\Bar J}$ is skew-symmetric.
\item[(ii)] There exists a metric $f$-connection $\nabla$, such that $(T_{|\hor^3})^{0,2}$ is skew-symmetric.
\item[(iii)] There exists a metric connection $\nabla$, satisfying ${\nabla F}_{|\hor^3}=0$,  such that $(T_{|\hor^3})^{0,2}$ is skew-symmetric.
\end{description}
\end{defn}
\begin{prop}\label{prop-horG_1}
Let $(N,F,h)$ be a metric $f$-manifold. Then the following statements are equivalent.
\begin{description}
\item[(i)] $(N,F,h)$ is horizontally of type $\mal G_1$.
\item[(ii)] There exists a metric $f$-connection $\nabla$, such that $T_{|\hor^3}$ is skew-symmetric.
\item[(iii)] There exists a metric connection $\nabla$, satisfying ${\nabla F}_{|\hor^3}=0$,  such that $T_{|\hor^3}$ is skew-symmetric.
\end{description}
In this case, for any  such connection satisfying (i) or (ii), then $T_{|\hor^3}$ is unique (i.e. uniquely determined by the metric $f$-manifold $(N,F,h)$). Conversely any extension $T\in\mal C(\mal T)$ of this unique skew-symmetric trilinear form $T_{|\hor^3}$ defines a metric connection $\nabla$, satisfying ${\nabla F}_{|\hor^3}=0$.
\end{prop}
\textbf{Proof.}
(ii) and (iii) are equivalent for the same reason as above. Furthermore, according to theorem~\ref{thm-nablaOmega_F=0}, if ${\nabla F}_{|\hor^3}=0$ and $T_{|\hor^3}$ is skew-symmetric, then $N_{\Bar J}=4(T_{|\hor^3})^{0,2}$ is also skew-symmetric (according to corollary~\ref{cory-B-eps,eps'} applied to the Hermitian bundle $(\hor,\Bar J,h_\hor)$) and moreover $-(T_{|\hor^3})^{**}=(d^c\Omega_{F|\hor^3})^{**}$, which proves the uniqueness of $T_{|\hor^3}$.\\
Conversely, if (i) is satisfied, then let $T\in\mal C(\mal T)$ such that  $(T_{|\hor^3})^{0,2}=\frac{1}{4}N_{\Bar J}$, $ (T_{|\hor^3})^{**}=-(d^c\Omega_{F|\hor^3})^{**}$, and the other components being arbitrary. Then $T_{|\hor^3}$ is skew-symmetric, and the corresponding metric connection $\nabla$ satisfies ${\nabla\Omega_F}_{|\hor^3}=0$, according to theorem~\ref{thm-nablaOmega_F=0}. This completes the proof.\medskip\hfill$\square$\\
Now, let us regroup the previous results to conclude.
\begin{defn}
A metric $f$-manifold  $(N,F,h)$ with skew-symmetric extended Nijenhuis tensor $\tl N_F$ will be sayed   \textbf{of global type $\mal G_1$} or \text{globally of type $\mal G_1$}.
\end{defn}
\begin{thm}\label{characteristic-equiv-G1} 
A metric $f$-manifold $(N,F,h)$ admits a metric $f$-connection $\nabla$ with skew-symmetric torsion \iif it is reductive and of global type $\mal G_1$ . Moreover, in this case, for any $\alpha\in \mal C(\Lambda^3\ver^*)$, there exists a unique metric connection $\nabla$ with skew-symmetric torsion such that $T_{|\Lambda^3\ver}=\alpha$. This unique connection is given by 
$$
T=(-d^c\Omega_F + N_{F|\hor^3}) + \mrm{Skew}(\Phi) + \mrm{Skew}(\mR_\ver) + \alpha.
$$
\end{thm}
\textbf{Proof.}
The first assertion follows from propositions~\ref{prop-horG_1}, \ref{H^2V-redu-typ-G_1}, theorem~\ref{thm-nablaOmega_F=0},  corollary~\ref{caract-V^2H} and remark~\ref{rmk-extend-red}. Then in this case, $T_{|\hor^3}$ is entirely determined, according to proposition~\ref{prop-horG_1}. Moreover, by skew-symmetry and theorem~\ref{nabla-q=0}, we have $T_{|\mal S(\hor\times\hor\times\ver)}=\mrm{Skew}(\Phi)$  and  $T_{|\mal S(\ver\times\ver\times\hor)}=\mrm{Skew}(\mR_\ver)$. Now, let us determine $T_{|\hor^3}$. Since, $\nabla= D + \frac{1}{2}T$, the equation  $\nabla F=0$ can be written 
$$
D\Omega_F + \dfrac{1}{2}(T(\cdot,F\cdot,\cdot) + T(\cdot,\cdot,F\cdot))=0
$$
so that 
$$
d\Omega_F=-F\circact T
$$
and thus $d{\Omega_F}_{|\hor^3}=-\Bar J\circact T_{|\hor^3}= 4\Bar J \cdot (T_{|\hor^3})^{--} - \Bar J\cdot T_{|\hor^3}=\Bar J \cdot N_{\Bar J} - \Bar J\cdot T_{|\hor^3}$. This completes the proof.\hfill$\square$
\begin{rmk}\em
In the particular case where $\ver$ is a line bundle, $(N,F,h)$ is an almost contact metric manifold. Then in this particular case, the theorem~\ref{characteristic-equiv-G1} above allows us to recover the result of \cite[theorem~{8.2}]{Friedrich-ivanov}.
\end{rmk}
\begin{defn}
On a metric $f$-manifold $(N,F,h)$, a metric $f$-connection $\nabla$ with skew-symmetric torsion will be called  characteristic connection. 
\end{defn}
\begin{cory}\label{characteristic-Nabla-c}
 Let $(N,F,h)$ be a reductive metric $f$-manifold of global type $\mal G_1$. Let us suppose that is given some metric connection $\nabla^c$ on the vertical subbundle $\ver$. There exists a metric $f$-connection $\nabla$ on $N$ with skew symmetric torsion, which coincides with $\nabla^c$ on $\ver$ \iif  
$$
(D^v-\nabla^c)(X^h,Y^v,Z^v)=-\dfrac{1}{2}\mR_\ver(Y^v,Z^v,X^h) \quad  \text{and} \quad(D^v-\nabla^c)_{|\ver^3}\in \mal C(\Lambda^3\ver^*).
$$
%_{|\hor\times\ver\times\ver}
In this case this connection $\nabla$ is unique and will be called the characteristic connection extending or defined by $\nabla^c$.
\end{cory}
\begin{rmk}\label{rmk-nabla-c-set}\em
In other words, in a reductive metric $f$-manifold of global type $\mal G_1$, the set of metric connection $\nabla^c$ on the vertical subbundle $\ver$ which can be extended  to a characteristic connection, is the affine space
$$
D^v -\dfrac{1}{2}\mR_\ver(Y^v,Z^v,X^h) + \mal C(\Lambda^3\ver^*).
$$
\end{rmk}
\textbf{Proof of corollary~\ref{characteristic-Nabla-c}.} This follows immediately from the theorem~\ref{characteristic-equiv-G1}, the proposition~\ref{D^v-nabla^c} and the  theorem~\ref{thm-nablaOmega_F=0}.\hfill$\square$
\begin{prop}\label{prop-F-bullet-T}
Let $(N,F,h)$ be a reductive metric $f$-manifold of  global type $\mal G_1$. Let $\nabla$ be some characteristic connection on $N$. Then the we have
\begin{eqnarray*}
d\Omega_F=-F \circact T & = & F\cdot N_F-F\cdot T-F \circact(T_{|\mal S(\hor\times\hor\times\ver)}+ T_{|\mal S (\ver\times\ver\times\hor)})\\
 &  = &   F\cdot N_F  -F\bullet T  -  \frac{1}{2}F \circact(T_{|\mal S(\hor\times\hor\times\ver)}+ T_{|\mal S (\ver\times\ver\times\hor)}). 
\end{eqnarray*}
\end{prop}
\textbf{Proof.} We have seen (in the proof of theorem~\ref{characteristic-equiv-G1}) that  $d\Omega_F=-F\circact T$, and that $(-{F\circact} T)_{|\hor^3}=F\cdot N_F-F\cdot T$, moreover we have $(-F\circact T)_{|\ver^3}=0$, by definition of the action $F\circact$. This completes the proof.\hfill$\square$
\paragraph{A useful expresssion of the Nijenhuis tensor in terms of the Levi-Civita connection.}
${}$ \\
By a direct computation, we prove:
\begin{prop}\label{N_F-intermsof-D}
Let $(N,F,h)$ be a metric  $f$-manifold. Then we have
\begin{eqnarray*}
N_F(X,Y)  & = & \left(D_{FX}F\right) Y  - \left( D_{FY} F \right)X  - F\left( D_X F \right) Y  + F \left( D_Y F\right)X\\
 & = & \left(D_{FX}F\right) Y - \left( D_{FY} F \right)X  +  \left( D_X F \right) (FY) - \left( D_Y F\right)(FX) -  d^D q (X,Y).
\end{eqnarray*}
\end{prop}
\index{connection!characteristic|)}
\subsubsection{Precharacteristic and  paracharacteristic connections.}\label{Prechara-parachara-connect}
\index{connection!precharacteristic} \index{connection!paracharacteristic|(} 
Sometimes the condition of  global type $\mal G_1$ could be too much strong and it could happen that one needs the existence (and uniqueness up to the $\ver^3$-component of the torsion) of some characteristic connection by supposing weaker conditions on the metric $f$-manifold $(N,F,h)$.
\begin{defn}
Let $(N,F,h)$ be a  metric $f$-manifold of  horizontal type $\mal G_1$. Then any metric $f$-connection $\nabla$ with a skew-symmetric component $T_{|\hor^3}$ of its torsion, will be called a horizontal-characteristic connection.\\
Moreover, if we suppose that $(N,F,h)$ is also of type $\ver^2\hor$, then a metric $f$-connection $\nabla$ with  skew-symmetric components $T_{|\hor^3}$, $T_{|\mal S(\ver\times\ver\times\hor)}$ and $T_{|\ver^3}$ of its torsion, will be called a  precharacteristic connection.
\end{defn}
\begin{rmk} \em
Let us remark that in a  metric $f$-manifold of  horizontal type $\mal G_1$,    horizontal-charactersitic connections always exist,  and the component $T_{|\hor^3}$ of the torsion is unique. Moreover, if we suppose that $(N,F,h)$ is also of type $\ver^2\hor$, then  precharacteristic connections always exist and the components $T_{|\hor^3}$ and $T_{|\mal S(\ver\times\ver\times\hor)}$ are unique.
\end{rmk}
The following properties  will hold for the  horizontal curvature in all the examples of interest for us.
\begin{defn}\label{defn-pure}
Let $(N,F)$ be an $f$-manifold. Let $A\in\mal C(\hor^*\otimes\hor^*\otimes\ver)$, $B \in \mal C(\hor^*\otimes\ver^*\otimes\hor)$, and $C\in \mal C(\ver^*\otimes\hor^*\otimes\hor)$. Then we will say respectively that $A$, $B$ or $C$ is pure if respectively
\begin{description}
\item[(i)] $A(\bar JX,Y)=A(X,\Bar J Y)$, $\forall X,Y\in \hor$.
\item[(ii)] $B(\bar JX, Y)=-\Bar J B(X,Y)$, $\forall X\in \hor, Y\in\ver$.
\item[(iii)] $C(X, \bar JY)=-\Bar J B(X,Y)$, $\forall X\in\ver,Y\in \hor$. 
\end{description}
If $(N,F)$ is endowed with a compatible metric $h$, then this means that $A$, $B$ or $C$ considered as element of $\mal C(\hor^*\otimes\hor^*\otimes\ver)$, satisfies respectively $A^{1,1}=0$, $B^{1,1}=0$, $C^{1,1}=0$. Moreover, we will say that $A$, $B$ or $C$ resp. (considered as trilinear forms) is  skew-symmetric in $\hor\times\hor$ if resp. $A$ is skew-symmetric \wrt the 2 first variables, $B$ \wrt the first and third variables, and $C$ \wrt the 2 last variables.
\end{defn}
Let us remark that $N_{F|\ver\times\hor\times\hor}$ is pure by definition of $N_F$ (see proposition~\ref{identities-N_F}), and skew-symmetric in $\hor\times\hor$ if $(N,F,h)$ is of type $\hor^2\ver$ (see remark~\ref{N-f-skew}).
\begin{defn}
$\bullet$  $(N,F,h)$ will be called almost of type $\hor^2\ver$  if  one of the following equivalent statements holds
\begin{description}
\item[(i)] $\mrm{Sym}_{\hor\times\hor}\left( D{\Omega_q}_{|\hor\times\hor\times\ver} \right)$ is pure.
\item[(ii)] $\mrm{Sym}_{\hor\times\hor}\left( N_{F|\ver\times\hor\times\hor}\right)=2\,\mrm{Sym}_{\hor\times\hor}\left( D{\Omega_q}_{|\hor\times\hor\times\ver} \right)(Y^h,Z^h,X^v).$
\end{description}
$\bullet$ If moreover $(N,F,h)$ is of type $\ver^2\hor$, then we will say that it is almost reductive.%\medskip\\
%
%$\bullet$ We will say that $(N,F,h)$ is \textbf{almost of global type} $\mal G_1$ if it is horizontally of type $\mal G_1$ and almost of reductive type $\mal G_1$.
\end{defn}
\begin{rmk}\em
We remark that if $(N,F,h)$  is of type $\hor^2\ver$ then it is almost of type $\hor^2\ver$, and therefore
%is globally of type $\mal G_1$, then it is, in particular, almost of global type $\mal G_1$. Moreover 
if it is reductive then it is, in particular, almost reductive. 
\end{rmk}
\begin{thm}\label{thm-exist-uniq-parachara}
A metric $f$-manifold $(N,F,h)$ admits a precharacteristic connection $\nabla$ such that  the component $T_{|\ver\times\hor\times\hor}$ of its torsion is pure  \iif it is   almost reductive and   horizontally of type $\mal G_1$. Moreover, in this case, for any $\alpha\in \mal C(\Lambda^3\ver^*)$, there exists a unique precharacteristic connection $\nabla$ such that  the component $T_{|\ver\times\hor\times\hor}$ of its torsion is pure, and  such that $T_{|\Lambda^3\ver}=\alpha$. This unique connection is given by 
$$
T=(-d^c\Omega_F + N_{F|\hor^3}) + T_{|\mal S(\hor\times\hor\times\ver)} + \mrm{Skew}(\mR_\ver) + \alpha.
$$ 
where $T_{|\mal S(\hor\times\hor\times\ver)}=\Phi + \dfrac{1}{2}(N_F(X^v,Y^h,Z^h) - N_F(Y^v,X^h,Z^h) )$.\\
Moreover if we impose also the component $T_{|\ver\times\hor\times\hor}$ to be  skew-symmetric in $\hor\times\hor$, then this is possible \iif $(N,F,h)$  is  reductive and  of horizontal type $\mal G_1$.
\end{thm}
\proof
Concerning the components in $\hor^3$, $\mal S(\ver\times\ver\times\hor)$ and $\ver^3$, we can proceed as in the proof of theorem~\ref{characteristic-equiv-G1}. Concerning the component in $\mal S(\hor\times\hor\times\ver)$, use the identity $ T(X^v,FY^h,Z^h) + T(X^v, Y^h,FZ^h)= N(X^v, Y^h, FZ^h)$  from theorem~\ref{thm-nablaOmega_F=0} and the identity $ T_{|\hor\times\hor\times\ver}=\Phi$ from theorem~\ref{nabla-q=0}.  \comprf \hsq
\begin{defn}
On a  metric $f$-manifold $(N,F,h)$, a precharacteristic connection $\nabla$ such that   the component $T_{|\ver\times\hor\times\hor}$ of the torsion is pure, will be called a \textbf{paracharacteristic} connection.
\end{defn}
\begin{prop}
Let $(N,F,h)$ be a reductive metric $f$-manifold of global type $\mal G_1$. Then the paracharacteristic connection (defined by some $\alpha\in\mal C(\Lambda^3\ver^*)$) coincide with  the characteristic connection (defined by the same $\alpha\in\mal C(\Lambda^3\ver^*)$) \iif the horizontal curvature $\Phi$ is pure.
\end{prop}
\proof
According to theorems~\ref{characteristic-equiv-G1}  and \ref{thm-exist-uniq-parachara}, we only have to prove that the components $T_{|\mal S(\hor\times\hor\times\ver)}$ of the torsions of the two connections coincide if $\Phi$ is pure. But, since $\Phi$ is pure,  the component $T_{|\mal S(\hor\times\hor\times\ver)}=\mrm{Skew}(\Phi)$ of the characteristic connection is such that $T_{|\ver\times\hor\times\hor} $ is pure, and therefore by uniqueness of the paracharacteristic connection, the two connections coincide. Conversely, if they coincide, then $\mrm{Skew}(\Phi)_{|\ver\times\hor\times\hor}$ is pure, i.e. $\Phi$ is pure. \comprf \hsq

\paragraph{The Linear representations of the curvatures.} Using the metric $h$, we have  canonical isomorphisms $\Lambda^2\hor^*\cong \so(\hor)$ and $\Lambda^2\ver^*\cong \so(\ver)$. Let us denote by $\rho\in\so(\hor)\otimes\ver^*$ and $\sigma\in \so(\ver)\otimes\hor^*$ respectively, the elements   corresponding to $\Phi$ and $\mR_\ver$ respectively  under these isomorphisms:
\begin{equation}\label{rho-def}
\langle \rho(V)H_1,H_2\rangle =\langle\Phi(H_1,H_2),V\rangle, \qquad  H_1,H_2\in\hor, V\in\ver,
\end{equation}
and
\begin{equation}\label{sigma-def} 
\langle \sigma(H)V_1,V_2\rangle =\langle \mR_\ver(V_1,V_2),H\rangle, \qquad  V_1,V_2\in\ver, H\in\hor.
\end{equation}
To do not weigh down the notation, we denote $\Bar J=\Bar J$. Let us introduce \emph{the horizontal curvature operator}\,:
$$
\Bar R(X,Y)Z=\rho(\Phi(X,Y))Z.
$$
as well as \emph{its derivation term}
$$
\Bar{\mrm  A}(X,Y)=\Bar R(\Bar J X,Y) + \Bar R(X,\Bar J Y) - [\Bar J, \Bar R(X,Y)].
$$
(which vanishes \iif $\Bar J$ is a derivation of $\Bar R$ hence its name). We denote by $\Phi=\Phi^{(+)} + \Phi^{(-)}$ the splitting of $\Phi$ according to the eigenspace decomposition of the endomorphism of $ \Lambda^2\hor^*\otimes\ver $ defined by $B\mapsto B(\Bar J\cdot,\Bar J\cdot)$, i.e.
$$
\Phi^{(\eps)}(\Bar J\cdot,\Bar J\cdot)=\eps \Phi^{(\eps)}.
$$
In other words, $\Phi^{(+)}$ is the $(1,1)$-type part of $\Phi$ whereas $\Phi^{(-)}$ is the part of $\Phi$ of type $(2,0) + (0,2)$. Under the isomorphism $\Lambda^2\hor^*\cong \so(\hor)$, to this, corresponds  the decomposition 
$\rho=\rho^+ + \rho^-$, where $\Bar J\rho^\eps\Bar J^{-1}=\eps \rho^\eps$, according to the splitting $\so(\hor)=\so_+(\hor)\oplus\so_-(\hor)$ of $ \so(\hor)$ following its $\Bar J$-commuting and $\Bar J$-anticommuting parts. Then, this being done, we can define the corresponding curvature operator and  antiderivation terms:
$$
\Bar R^{(\eps)}=\rho^\eps(\Phi^{(\eps)}) \quad \text{and} \quad \Bar{\mrm  A}^{(\eps)}(X,Y)=\Bar R^{(\eps)}(\Bar J X,Y) + \Bar R^{(\eps)}(X,\Bar J Y) - [\Bar J, \Bar R^{(\eps)}(X,Y)].
$$
We remark that $\Bar{\mrm  A}^{(+)}=0$ and $\Bar{\mrm  A}^{(-)}(X,Y)=2\Bar R^{(-)}(X,\Bar JY)-2\Bar J\circ\Bar R^{-}(X,Y)$.
\begin{rmk}\label{rmk-H2v-skew-Nf} \em
$\bullet$ Let us remark that $\Phi$ is pure \iif  $\rho$ anticommute with $\Bar J$: $\rho(V)\Bar J =-\Bar J\rho(V)$, $\forall V\in\ver$.\\
$\bullet$ Furthermore, the condition that $\tilde N_{F|\mal S(\hor\times\hor\times\ver)}$ is skew-symmetric can be expressed in terms of the linear representations introduced above. Indeed, according to proposition~\ref{prop-N-f-hhv}-(iii), this condition can be written: $\Phi(H_1,\Bar J H_2,V) + \Phi(\Bar J H_1, H_2,V)=N_F(V,H_1,\Bar J H_2)$, for all $H_1,H_2\in \hor$ and $V\in\ver$. This last equation is then equivalent to: 
$$
[\rho(V),\Bar J]=-\Bar J N_F(V),  \quad \forall V\in \ver,
$$
where we have set $\langle N_F(V)H_1,H_2\rangle := N_F(V,H_1, H_2)$. This means that $\rho^-(V) = 2 N_F(V)$, $\forall V\in \ver$. We will see that this condition could happen to be  very strong, for example in the twistor bundles (see \ref{variational-twistor}).
\end{rmk}
\index{global type $\mal G_1$|)}
%
%
%%%%%%%%%%%%%%%%%%%%%%%%%%%%%%%%%%%%%%%%%%%%%%%%%%%%%%%%%%%%%%%%%%%%%%%%%%%%%%%
%
%
\subsubsection{Reductions of $f$-manifolds}\label{subsec-Reduct-f-manif}
\begin{defn}
A $q$-manifold $(N,q)$ is a manifold endowed with a linear projector $q\in\End(TN)$, i.e. with a splitting $TN=\ver\oplus\hor$. We will say  that $q$ is a $q$-structure. Moreover, we will say that $(N,q)$ is associated to the $f$-manifold $(N,F)$ if $q$ is associated to $F$: i.e $-F^2=\Id -q$, or equivalently $\ker F=\ver$ and $\im F=\hor$.\\
A metric $q$-manifold $(N,q,h)$ is a $q$-manifold endowed with a metric $h$ for which $q$ is orthogonal, or equivalently the splitting $TN=\ver\oplus\hor$ is orthogonal.
\end{defn}
\begin{defn}
A map $i_\mv\colon  (N^\mv,q^\mv) \to (N,q)$ from a $q$-manifold to another one is called a \textbf{$q$-immersion} if it is an immersion which satisfies 
\begin{equation}\label{ivq}
(i_\mv)^* q=q^\mv.
\end{equation}
Then $(N^\mv,q^\mv)$ is \textbf{$q$-immersed} in $(N,q)$. If $i_\mv$ is injective we will say that $(N^\mv,q^\mv)$ is a \textbf{$q$-submanifold} of the $q$-manifold $(N,q)$, and that $i_\mv$ is a \textbf{$q$-imbedding}.
Moreover we will say that $i_v\colon (N^\mv,q^\mv)\to (N,q) $ is a \textbf{reduction}  of $q$-manifolds, or that $(N^\mv,q^\mv)$ is a reduction of $(N,q)$, if $i_\mv$ is a $q$-imbedding  and the vector bundles $\hor^\mv$ and $\hor$ have the same rank that is to say
\begin{equation}\label{hor-v=hor}
\hor_{|N^\mv}=\hor^\mv.
\end{equation}
A map $i_\mv\colon (N^\mv, q^\mv, h^\mv)\to (N, q, h)$ between two metric $q$-manifolds is called a \textbf{$q$-isometry} if it is a $q$-immersion  and an isometry. If $i_\mv$ is injective we will say that $(N^\mv,q^\mv)$ is a metric \textbf{$q$-submanifold} of the metric $q$-manifold $(N,q)$, and that $i_\mv$ is a isometric \textbf{$q$-imbedding}.  If $i_\mv$ is moreover a reduction of  $q$-manifolds,  we will  say that $i_\mv$ is a \textbf{reduction} of metric $q$-manifolds.\\
A  map $i_v\colon (N^\mv,F^\mv)\to (N,F)$  from a  $f$-manifold $(N^\mv,F^\mv)$  into another one $(N,F)$  is  called a \textbf{$f$-immersion} if it is an immersion which satisfies
\begin{equation}\label{ivF}
(i_\mv) ^*F=F^\mv.
\end{equation}
When $i_\mv$ is also injective, we will say that $(N^\mv,F^\mv)$ is a  $f$-submanifold of $(N,F)$ and that $i_\mv$ is an $f$-imbedding. Moreover we will  then say that $i_v\colon (N^\mv,F^\mv)\to (N,F) $ is a reduction of $f$-manifolds or that $(N^\mv,F^\mv)$ is a reduction of $(N,F)$ if the vector bundles $\hor^\mv=\im F^\mv$ and $\hor=\im F$ have the same rank that is to say $\hor_{|N^\mv}=\hor^\mv.$\\
A map $i_\mv\colon (N^\mv, F^\mv, h^\mv)\to (N, F, h)$ between two metric $f$-manifolds is called a \textbf{$f$-isometry} if it is an $f$-immersion  and an isometry. If $i_\mv$ is injective we will say that $(N^\mv,F^\mv)$ is a metric \textbf{$f$-submanifold} of the metric $f$-manifold $(N,F)$, and that $i_\mv$ is a isometric \textbf{$f$-imbedding}.  If $i_\mv$ is moreover a reduction of  $f$-manifolds,  we will  say that $i_\mv$ is a \textbf{reduction} of metric $f$-manifolds.
\end{defn}
\begin{rmk}\em
\eqref{ivq} means  that $i_\mv$ sends the splitting $TN^\mv=\hor^\mv\oplus\ver^\mv$ into the splitting $TN=\hor\oplus\ver$. In other words,  the former is the trace on $TN^\mv$ of the latter. It can be written also in the  following form (which suggests the definition of holomorphicity)
$$
di_\mv\circ q^\mv= q\circ di_\mv.
$$
\eqref{ivF} implies \eqref{ivq}, so that a  $f$-immersion is in particular a $q$-immersion. Moreover a reduction of (metric) $f$-manifolds is a reduction of the corresponding (metric) $q$-manifolds. Remark also that \eqref{ivF} means simply that $i_\mv$ is a \emph{$f$-holomorphic map} from $(N^\mv,F^\mv)$ to $(N,F)$: $di_\mv\circ F^\mv =F\circ di_\mv$.
\end{rmk}
\begin{lemma}\label{lemma-N_F}
$\bullet$ Let $i_\mv \colon (N^\mv,q^\mv)\to (N,q)$ be a $q$-immersion. Then we have
$$
(i_\mv)^* \Phi=\Phi^\mv \quad  \text{and} \quad  (i_\mv)^* \mR_\ver=\mR_{\ver^\mv}
$$
where $\Phi^\mv$ is of course the curvature of the horizontal subbundle $\hor^\mv$.\medskip\\
$\bullet$ Let $i_\mv \colon (N^\mv,F^\mv)\to (N,F)$ be an $f$-immersion. Then we have
$$
(i_\mv)^* N_F=N_{F^\mv}
$$
\end{lemma}
%%%
%
\begin{defn}
Let $i_\mv\colon (N^\mv,h^\mv)\to (N,h)$ be a isometry between two Riemannian manifolds.  Let $\nabla$ be a linear  connection on $N$. The projection on $N^v$ of $\nabla$ is the connection $\nabla^\mv$ on $N^v$ defined by 
$$
\langle \nabla_X^\mv Y,Z\rangle = \langle \nabla_X \tl Y,Z\rangle \qquad \forall X,Z\in TN^\mv, Y\in \mal C(TN^\mv)
$$
where $\tl Y\in \mal C(TN)$ is some extension of $Y$.\smallskip\\
Let $(E,h)\to N$ be a Riemannian vector bundle endowed with a linear connection $\nabla$. Let $N^\mv\subset N$ be an immersed submanifold and  $E^\mv\to N^\mv$ be a vector subbundle of $E_{|N^\mv}$. Then we define  analogously the projection of $\nabla$ on $E^\mv$.
\end{defn}
\begin{lemma}
If $\nabla^\mv$ is the projection of $\nabla$ then their respective torsion 3-forms are related by 
$$
\langle T^{\nabla^\mv}(\cdot,\cdot),\,\cdot\,\rangle = (i_\mv)^*\langle T^{\nabla}(\cdot,\cdot),\,\cdot\,\rangle  
$$
\end{lemma}
\begin{lemma}\label{domegaq}
$\bullet$ Let  $i_\mv\colon (N^\mv,h^\mv,q^\mv)\to (N,h,q)$ be a  $q$-isometry. \\
\textbf{\em (i}) Then the Levi-Civita connection  $D^\mv$ of $(N,h^\mv)$ is the projection on $N^v$ of the Levi-Civita connection of $(N,h)$.\\
\textbf{\em (ii)} The following identity holds:
$$
D^\mv \Omega_{q^\mv}=(i_\mv)^*(D \Omega_q)
$$
$\bullet$ Let  $i_\mv\colon (N^\mv,h^\mv,F^\mv)\to (N,h,F)$ be a  $f$-isometry. Then $(i_\mv)^*\Omega_F=\Omega_{F^\mv}$.
\end{lemma}
\begin{thm}\label{thm-proj-chara}
Let  $i_\mv\colon (N^\mv,h^\mv,F^\mv)\to (N,h,F)$ be an $f$-isometry.\\
$\bullet$ Then if $(N,F,h)$ is horizontally of type $\mal G_1$, resp. of global type $\mal G_1$, resp. almost reductive, resp. reductive, then so is $(N^\mv,F^\mv, h^\mv)$.\\
$\bullet$ Moreover, if $(N,F,h)$ is horizontally of type $\mal G_1$ and almost reductive, then the projection on $N^\mv$ of its paracharacteristic connection (defined by some $\alpha\in \mal C(\Lambda^3\ver^*)$)  coincides with the paracharacteristic connection of $(N^\mv,F^\mv,h^\mv)$ (defined by  $(i_\mv)^*\alpha\in \mal C(\Lambda^3(\ver^\mv)^*)$). In the same way, if $(N,F,h)$ is reductive and of global type $\mal G_1$, then the projection on $N^\mv$ of its characteristic connection (defined by some $\alpha\in \mal C(\Lambda^3\ver^*)$)  coincides with the characteristic connection of $(N^\mv,F^\mv,h^\mv)$ (defined by  $(i_\mv)^*\alpha$).
\end{thm}
\textbf{Proof.} 
The first assertion follows from lemmas~\ref{lemma-N_F} and \ref{domegaq}. The second assertion follows from 
%the fact that the unique horizonal 3-form $\Bar T\in\mal C(\Lambda^3\hor^*)$ such that the torsion $T$ of any  metric $f$-connection $\nabla$  with a  skew-symmetric component $T_{\hor^3}$, satisfies $T_{\hor^3}=\Bar T$ (see  proposition~\ref{prop-horG_1}), is given by $d\Omega_{F|\hor^3}=\Bar J\circact \Bar T$, i.e. $\Bar T = -d^c\Omega_F + N_{F|\hor^3}$.
theorems~\ref{thm-exist-uniq-parachara} and \ref{characteristic-equiv-G1}, and lemmas~\ref{lemma-N_F} and \ref{domegaq}. 
This completes the proof.\hfill$\square$%\medskip\\
\begin{defn}
Let $E\to N$ be a vector bundle endowed with a linear connection $\nabla$. Let $N^\mv\subset N$ be an immersed submanifold and  $E^\mv\to N^\mv$ be a vector subbundle of $E_{|N^\mv}$. We will way that $\nabla$ is reducible in $E^\mv$ if $E^\mv$ is $\nabla$-parallel, that is to say parallel \wrt the induced connection on $E_{|N^\mv}$.  Then the induced connection on $E^\mv$ will be called the reduction of $\nabla$ in $E^\mv$.
\end{defn}
\begin{rmk}\em
If $E=TN$ and $E^\mv=TN^\mv$ then we recover the usual definition of reducibility of linear connection.
\end{rmk}
\begin{lemma}\label{torsion-red-related}
If $i_\mv\colon L\to N$ is an immersion and $\nabla^\mv$ the reduction of $\nabla$ then their torsion 2-forms are related by 
$$
T^{\nabla^\mv}=(i_\mv)^* T^\nabla
$$
\end{lemma}
\begin{prop}
Let $i_\mv\colon (N^\mv,h^\mv,F^\mv)\to (N,h,F)$ be a  reduction of metric $f$-manifolds. Let us suppose that the vertical component $D^v$ of the Levi-Civita connection of $(N,h)$is reducible in $\ver^\mv$.
Then the Levi-Civita connection of $(N,h)$ is reducible in $N^\mv$ and its reduction  is the Levi-Civita connection  $D^\mv$ of $(N,h^\mv)$.
\end{prop}
\begin{defn}\label{def-con-ver}
Let $(N,q,h)$ be a metric  $q$-manifold. We define the  affine space of \emph{compatible vertical connection} by 
$$
\mrm{Con}(\ver) = \left\lbrace D^v - \dfrac{1}{2}\mR_\ver(Y^v,Z^v,X^h) + \mal C(\Lambda^3\ver^*)\right\rbrace .
$$
If $(N,q,h)$ is associated to a metric $f$-manifold $(N,F,h)$ which is  horizontally of type $\mal G_1$ and almost reductive (respectively reductive and of global type $\mal G_1$) then this set corresponds to the set of metric connection $\nabla^c$ on the vertical subbundle $\ver$ which can be extended  to a paracharacteristic connection (respectively a characteristic connection), (see remark~\ref{rmk-nabla-c-set}).
\end{defn}
\begin{rmk}\em
If $\mR_\ver=0$, then  $D^v\in \con(\ver)$.
\end{rmk}
\begin{defn}\index{complete reduction}
Let  $i_\mv\colon (N^\mv,h^\mv,F^\mv)\to (N,h,F)$ be a reduction of metric $f$-manifolds. We will  say that $i_\mv\colon (N^\mv,h^\mv,F^\mv)\to (N,h,F)$ is a \textbf{complete reduction} of metric $f$-manifold if there exists $\nabla^c\in \con(\ver)$ which is reducible in  $\ver^\mv$.\\
More generally, we define analogously a complete reduction of metric $q$-submersions.
\end{defn}
\begin{prop}\label{con-ver-red}
Let $i_\mv\colon (N^\mv,q^\mv,h^\mv)\to (N,q,h)$ be a reduction of metric $q$-manifolds.
If $\nabla^c\in\con(\ver)$ is reducible in $\ver^\mv$, then its reduction is also in $\con(\ver^\mv)$.
\end{prop}
\proof 
According to lemmas~\ref{lemma-N_F} and \ref{domegaq} and definition~\ref{def-con-ver}, the projection of $\con(\ver)$ on $\ver^\mv$ is $\con(\ver^\mv)$. This completes the proof. \hsq
\begin{prop}\label{Tr-nabla-c-con-ver}
Let $(N,q,h)$ be a metric $q$-manifold, and $f\colon L\to N$ be  a map from a Riemann surface in $N$. Then $\mrm{Tr}_g(\nabla^c d^v f)$ does not depend on the choice of $\nabla^c\in \con(\ver)$.
\end{prop}
\proof
The difference between two elements of $\con(\ver)$ is in $\mal C(\Lambda^3\ver^*)$.
\begin{prop}\label{para-chara-red}
Let $i_\mv\colon (N^\mv,h^\mv,F^\mv)\to (N,F,h)$ be a complete reduction of metric $f$-manifolds. Let $\nabla^c\in\con(\ver)$ reducible in $\ver^\mv$, and $\nabla^{c,\mv}$ its reduction. Then if $(N,F,h)$ is horizontally of type $\mal G_1$ and almost reductive (respectively reductive and of global type $\mal G_1$), its paracharacteristic connection (respectively its characteristic connection) extending $\nabla^{c}$ is reducible to the paracharacteristic connection (respectively the characteristic connection) in $(N^\mv, h^\mv, F^\mv)$ extending $\nabla^{c,\mv}$.
\end{prop}
\proof 
This follows from theorem~\ref{thm-proj-chara}, equation \eqref{hor-v=hor}, and proposition~\ref{con-ver-red}. \comprf \hsq \medskip\\
Therefore, according to proposition~\ref{Tr-nabla-c-con-ver}, we deduce that
\begin{prop}\label{prop-tensionfield=}
Let $i_\mv\colon (N^\mv,h^\mv,F^\mv)\to (N,h,F)$ be a complete reduction of metric $f$-manifolds.\\
Let us suppose that $(N,F,h)$ is horizontally of type $\mal G_1$ and almost reductive. Let $f\colon L\to N^\mv$ be a map.\\
$\bullet$ The tension field of $f$ \wrt to one paracharactersitic connection $\nabla$  of $(N,h,F)$  coincides with the  tension field of $f$ \wrt to any paracharactersitic connection of $(N^\mv,h^\mv,F^\mv)$. \\
$\bullet$ In the same way, if $(N,F,h)$ is reductive and of global type $\mal G_1$, then the tension field of $f$ \wrt to one charactersitic connection $\nabla$  of $(N,h,F)$  coincides with the  tension field of $f$ \wrt to any charactersitic connection of $(N^\mv,h^\mv,F^\mv)$.
\end{prop}
\subsection{$f$-connections on fibre bundles}\label{f-connect-bundl}
Here, we consider the case where the vertical subbundle is the tangent space of the fibres of a fibration (or more generally a submersion) $\pi\colon N\to M$, i.e. $\ver=\ker d\pi$. Let us first remark that in this case $\mR_\ver=0$, which leads to immediate simplifications in the preceding results (obtained in section~\ref{f-connection&torsion}).\\
\textbf{Convention} In all the next of  the paper, all the submersions $\pi\colon N\to M$ are supposed to be surjective (i.e. the open set $\pi(N)$ coincides with $M$).
\subsubsection{Riemannian submersion and metric $f$-manifolds of global type $\mal G_1$.}\label{Riem-sub-Global-G_1}
\begin{prop}\label{prop-hor-Levi}
Let $\pi\colon (N,h)\mapsto (M,g)$ be a Riemannian submersion, over which we consider the natural orthogonal decomposition: $TN=\ver\oplus\hor$,  where $\ver=\ker d\pi$ and $\hor=\ver^\perp$.  Denote by $D$ and $D^g$ respectively the Levi-Civita connections of $(N,h)$ and $(M,g)$, respectively. Let $\widetilde{D^g}$ be the connection on $\hor$ defined by the lift of $D^g$:
$\widetilde{D_A^g}B = (d\pi_{|\hor})^{-1}(D_A^g\pi_*(B))\in \hor$ for all $A,B\in \mal C(TN)$.\\
Let us suppose that $(N,q,h)$ is of type $\mal V^2\mal H$. Then the horizontal component of the Levi-Civita connection in $N$ is related to $\widetilde{D^g}$ by the following formula:
$$
\langle D_A B, H \rangle= \langle\widetilde{D_A^g}B, H \rangle + \dfrac{1}{2}\left(\Phi(A,H,B^v) + \Phi(B,H,A^v) \right) 
$$
forall $A,B\in \mal C(TN)$ and $H\in\mal C(\hor)$.
\end{prop}
\textbf{Proof.} 
Let us set 
$$
S^M(A,B)=\pi_*(D_A B)- D_A^g(\pi_* B),\quad \forall A,B\in \mal C(TN).
$$
Then it is easy to see that $S^M$ is in fact a tensor, i.e. $S^M\in\mal C(T^*N\otimes T^*N\otimes \pi^*TM)$. Let $A,B\in TN$ and $H\in \hor$, and let us extend these to vector fields, denoted by the same notations, such that the horizontal components of these extension  are projectible: there exist vector fields $\Bar A,\Bar B, \bar H$ on $M$ such that 
$$
\pi_* A= \Bar A\circ\pi, \pi_* B= \Bar B\circ\pi, \pi_* H= \Bar H\circ\pi.
$$ 
Using  the fact that $h_{|\hor\times\hor}=\pi^*g$, the characterization of Levi-Civita yields:
\begin{eqnarray*}
2\langle S^M(A,B),\pi_* H\rangle & = & 2 h(D_A B, H) - 2 g(D_{\Bar A}^g(\Bar B), \Bar H)\circ\pi\\
   & = &  A.h(B^\hor,H) + B.h(A^\hor,H) - H.h(A,B) + h([A,B],H) + h([H,A],B) + h([H,B],A)\\
 &  & - \left( \Bar A .g(\Bar B,\Bar H) + \Bar B .g(\Bar A,\Bar H) - \Bar H .g(\Bar A,\Bar B)  + g([\Bar A,\Bar B],H) + g([\Bar H, \Bar A], \Bar B) + g([ \Bar H, \Bar B], \Bar A)\right) \circ\pi\\
& = & - H\cdot\langle qA,qB \rangle + \langle q[H,A],qB\rangle + \langle q[H,B], qA\rangle\\
& = & -\langle D_H^v (qA), qB\rangle - \langle qA, D_H^v (qB)\rangle + \langle q[H,A],qB\rangle + \langle q[H,B], qA\rangle\\
& = &  \langle T^v(A,H),B^v\rangle +  \langle T^v(B,H),A^v\rangle\\
& = & -\mrm{Sym}_{\ver\times\ver}(T_{|\hor\times\ver\times\ver}^v)(H,A^v,B^v) + \Phi(A^h,H,B^v) + \Phi(B^h,H,A^v)
\end{eqnarray*}
where we have used in the last line $T_{|\hor\times\hor}^v=\Phi$ (see proposition~\ref{dec-Tv}). We conclude by using proposition~\ref{V^2H-Tv}. This completes the proof.\hfill$\square$%\\
\begin{rmk}\label{rmk-nabla_ver}\em
By a direct computation using the characterization of Levi-Civita, we can prove
$$
\langle D_{X^v} Y^h, Z^h\rangle= \dfrac{1}{2}\Phi(Y^h,Z^h,X^v) + \dfrac{1}{2}\mrm{Sym}_{\hor\times\hor}\left( \nabla\Omega_{q|\hor\times\ver\times\hor}\right) (Z^h,X^v,Y^h) + \langle [X^v,Y^h],Z^h \rangle
$$
Then according to the previous proposition, we obtain
$$
\langle \widetilde{D_{X^v}^g} Y^h, Z^h\rangle = \dfrac{1}{2}\mrm{Sym}_{\hor\times\hor}\left( \nabla\Omega_{q|\hor\times\ver\times\hor}\right) (Z^h,X^v,Y^h) + \langle [X^v,Y^h],Z^h \rangle
$$
Moreover, this equation still holds if we replace  $D^g$ by any metric connection $\nabla^M$ on $M$, since then 
$\widetilde {(D^g -\nabla^M)}$ is a horizontal trilinear form.
\end{rmk}
\begin{prop}\label{prop-ver-Levi}
Let $\pi\colon (N,h)\mapsto (M,g)$ be a Riemannian submersion, with the same notations and definitions as in the previous proposition. Let us suppose that some metric connection $\nabla^c$ on $\ver$ is given, and denote by $T^c$ its vertical connection.\\
Then the vertical component of the Levi-Civita connection on $N$ is given by
$$
\langle D_A B, V \rangle= \langle\nabla_A^c B^v, V\rangle + \dfrac{1}{2}\left(\mB(A^v,B^v,V) - \Phi(A^h,B^h,V)  - \mrm{R}_a^c(A^h,B^v,V) -\mrm{R}_s^c(B^h,V,A^v)\right) 
$$
where
$$
\begin{array}{c}
\mB(V_1,V_2,V_3)= -T^c(V_1,V_2,V_3) - \mU^c(V_1,V_2,V_3)= -T^c(V_1,V_2,V_3) + T^c(V_1,V_3, V_2) + T^c (V_2,V_3,V_1)\\
\mrm{R}_a^c = \mrm{Skew}_{\ver\times\ver}(T_{|\hor\times\ver\times\ver}^c) \quad \mrm{and} \quad \mrm{R}_s^c= \mrm{Sym}_{\ver\times\ver}(T_{|\hor\times\ver\times\ver}^c).
\end{array}
$$
\end{prop}
\textbf{Proof.}  
Let us set $S^\ver(A,B)=q(D_A B) - \nabla_A^c(qB)$ for all $A,B\in \mal C(TN)$, this defines a element  $S^\ver\in \mal C(T^*N\otimes  T^*N\otimes \ver$.  Then let $A,B\in TN$ and $V\in \ver$, and let us extend these to vector fields, denoted by the same notations, such that the horizontal components of these extension  are projectible: there exist vector fields $\Bar A,\Bar B$ on $M$ such that 
$$
\pi_* A= \Bar A\circ\pi, \pi_* B= \Bar B\circ\pi, \pi_* V=0.
$$
Then using the characterization of the Levi-Civita connection, we have 
\begin{eqnarray*}
2\langle S^\ver(A,B),V\rangle & = & 2\langle D_AB,V\rangle - 2\langle \nabla_A^c B^v,V\rangle\\
 & = & A\cdot\langle B^v, V\rangle + B\cdot\langle A^v,V\rangle - V\cdot\langle A^v, B^v\rangle  \\
& - & \langle A^v,[B^v,V]^v\rangle + \langle B^v,[V,A]^v\rangle + \langle V,[A,B]^v\rangle -2\langle \nabla_A^c B^v, V\rangle\\
&  - &  V\cdot g(\pi_* A, \pi_*B) - g(\pi_*A, \pi_*[B,V]) + g(\pi_*B,\pi_*[V,A]) \\
& = & \langle \nabla_A^c B^v + \nabla_B^c A^v - 2\nabla_A^c B^v + [A,B]^v,V\rangle + \langle B^v, \nabla_A^c V - \nabla_V^c A^v + [V,A]^v\rangle\\
& &  + \langle \nabla_B^c V -\nabla_V^c B^v - [B,V]^v, A^v\rangle\\
& = & \langle T^c(B,A),V\rangle + \langle T^c(A,V),B^v\rangle + \langle T^c(B,V),A^v\rangle.
\end{eqnarray*}
Then we complete the proof by using proposition~\ref{Tc-Phi}.\hfill$\square$
\begin{rmk}\label{rmk(D-nablac)}\em
The results of that is that 
\begin{eqnarray*}
2(D-\nabla^c)_{|\hor\times\ver\times\ver} & = & \mrm{R}_a^c\\
2(D-\nabla^c)_{|\ver\times\hor\times\ver} & = & \mrm{R}_s^c(B^h,V,A^v)
\end{eqnarray*}
so that  $(D-\nabla^c)_{|\mal S(\ver\times\hor)}=0$ \iif the reductivity term $T^c_{|\ver\times\hor}$ vanishes, and in this case the vertical component of the Levi-Civita connection $D$ is given by
$$
\langle D_A B, V \rangle= \langle\nabla_A^c B^v, V\rangle + \dfrac{1}{2}\left(\mB(A^v,B^v,V) - \Phi(A^h,B^h,V)\right) 
$$ 
and we recover the corresponding relation in homogeneous fibre bundles of theorem~\ref {difference-tensor}, since in a homogeneous fibre bundle we have $T^c_{|\ver\times\hor}=0$. Moreover, from these considerations and proposition~\ref{V^2H-Tv}, we deduce the following.
\end{rmk}
\begin{cory}
The symmetric component of the reductivity term, $\mrm{Sym}_{\ver\times\ver}(T_{|\hor\times\ver\times\ver}^c)$, is independent of $\nabla^c$ and vanishes \iif $(N,q,h)$ is of type $\mal\ver^2\mal\hor$. Moreover, the reductivity term $T^c_{|\ver\times\hor}$ vanishes \iif $(N,q,h)$ is of type $\mal\ver^2\mal\hor$ and $(D-\nabla^c)_{|\hor\times\ver\times\ver}=0$. We will then say that $(N,q)$ is $\nabla^c$-reductive. In particular, if we take $\nabla^c=D^v$, the restriction to $\ver$ of the vertical component of Levi-Civita, then the $D^v$-reductivity means that $(N,q,h)$ is of type $\mal\ver^2\mal\hor$.
\end{cory}
Applying propositions~\ref{prop-hor-Levi}, \ref{prop-ver-Levi}, we obtain:
\begin{prop}\label{prop-Levi-dec}
Let $\pi\colon (N,h)\mapsto (M,g)$ be a Riemannian submersion. Let us suppose that there exists some metric connection $\nabla^c$ on $\ver$  for which $(N,q)$ is $\nabla^c$-reductive. Then, \wrt the decomposition $TN=\hor\oplus\ver$, the Levi-Civita connection in $N$ admits the following  decomposition:
$$
D_A =\begin{pmatrix}
\nabla_A^c & 0 \\ 0 & \widetilde{D_A^g} \end{pmatrix} 
+
\dfrac{1}{2}
\begin{pmatrix}
\mB(A^v,\cdot) & -\Phi(A^h,\cdot)\\ \rho(\cdot)A^h &  \rho(A^v) \end{pmatrix}.
$$
where, let us recall it, $\rho$ is defined by \eqref{rho-def}. 
%$\langle\Phi(A)^t(B^v), H\rangle=\langle B^v, \Phi(A,H)\rangle$.
In particular, this decomposition holds for homogeneous fibre bundles.
\end{prop}
\begin{rmk}\em
We see that according to this decomposition of Levi-Civita, we have $Dq_{|\hor\times\hor\times\ver}=\frac{1}{2}\Phi$, so that the $\nabla^c$-reductivity implies that $(N,q,h)$ is of type $\hor^2\ver$. In particular let us take $\nabla^c=D^v$ restricted to $\ver$, then the $D^v$-reductivity means that  $(N,q,h)$ is of type $\ver^2\hor$, and thus we see  that the type $\ver^2\hor$ implies the type $\hor^2\ver$. 
\end{rmk}
\begin{cory}\label{cory-reductivity}\index{homogeneous fibre bundle|(}
Let $\pi\colon (N,h)\mapsto (M,g)$ be a Riemannian submersion, endowed with its canonical orthogonal splitting $TN=\ver\oplus\hor$. Then if $(N,q,h)$ is of type $\ver^2\hor$ then it is also of type $\hor^2\ver$ and thus it is reductive. In particular, if $\ver$ can be endowed with a metric connection $\nabla^c$ with a vanishing reductivity term $T_{|\ver\times\hor}^c$, then $(N,q,h)$ is reductive. In particular, a homogeneous fibre bundle is reductive.
\end{cory}
\begin{cory}\label{cory-equiv-chara-Nablac}
Let $\pi\colon (N,h)\mapsto (M,g)$ be a Riemannian submersion, endowed with its canonical orthogonal splitting $TN=\ver\oplus\hor$. Let us suppose that $\hor$ is endowed with an orthogonal complex structure, that is to say $N$ is endowed with a  metric $f$-structure compatible\footnote{i.e. $\ker F=\ver$ and $\im F=\hor$.} with the previous splitting.\\ 
 Let us suppose that there exists some metric connection $\nabla^c$ on $\ver$  for which $(N,q)$ is $\nabla^c$-reductive, and that $T_{|\ver^3}^c$ is skew-symmetric. Then the following statements are equivalent
\begin{description}
\item[(i)] There exists a characteristic connection on $(N,F,h)$.
\item[(ii)] $(N,F,h)$ is of global type $\mal G_1$.
\item[(iii)] The canonical connection $\nabla^c$ can be extended to a characteristic connection.
\item[(iv)] There exists a Hermitian connection $\nabla^\hor$ on $\hor$ such that $\nabla:=\nabla^c\oplus\nabla^\hor$ has a skew-symmetric torsion.
\end{description}
In particular, these equivalences hold when $\pi\colon (N,h)\mapsto (M,g)$ is a homogeneous fibre bundle with a naturally reductive fibre $H/K$.
\end{cory}
\begin{rmk}\em
In other words, if $(N,q,h)$ is of type $\ver^2\hor$, then the  existence of a characteristic connection is equivalent to the global type $\mal G_1$, and in this case, the set of metric connections $\nabla^c$ on the vertical subbundle $\ver$ which can be extended  to a characteristic connection, is the affine space
$$
D^v  + \mal C(\Lambda^3\ver^*).
$$
\end{rmk}
\textbf{Proof of corollary~\ref{cory-equiv-chara-Nablac}.}
Since, according to corollary~\ref{cory-reductivity} , the $\nabla^c$-reductivity implies the reductivity, then (i) $\Leftrightarrow$ (ii) according to theorem~\ref{characteristic-equiv-G1}. Moreover the equivalence (ii) $\Leftrightarrow$ (iii) follows from corollary~\ref{characteristic-Nabla-c} and remark~\ref{rmk(D-nablac)}. Finally, the equivalence (iii) $\Leftrightarrow$ (iv) is obvious. This completes the proof.\hfill $\square$\medskip\\
We can rewrite the  proposition~\ref{cory-equiv-chara-Nablac} for paracharacteristic connections. 
\begin{cory}\label{cory-equiv-para-Nablac}
Let $\pi\colon (N,h)\mapsto (M,g)$ be a Riemannian submersion, endowed with its canonical orthogonal splitting $TN=\ver\oplus\hor$. Let us suppose that $\hor$ is endowed with an orthogonal complex structure, that is to say $N$ is endowed with a  metric $f$-structure compatible\footnote{i.e. $\ker F=\ver$ and $\im F=\hor$.} with the previous splitting.\\ 
 Let us suppose that there exists some metric connection $\nabla^c$ on $\ver$  for which $(N,q)$ is $\nabla^c$-reductive, and that $T_{|\ver^3}^c$ is skew-symmetric. Then the following statements are equivalent
\begin{description}
\item[(i)] There exists a paracharacteristic connection on $(N,F,h)$.
\item[(ii)] $(N,F,h)$ is of horizontal type $\mal G_1$.
\item[(iii)] The canonical connection $\nabla^c$ can be extended to a paracharacteristic connection.
%\item[(iv)] There exists a Hermitian connection $\nabla^\hor$ on $\hor$ such that $\nabla:=\nabla^c\oplus\nabla^\hor$ has a skew-symmetric torsion.
\end{description}
In particular, these equivalences hold when $\pi\colon (N,h)\mapsto (M,g)$ is a homogeneous fibre bundle with a naturally reductive fibre $H/K$.
\end{cory}
\index{homogeneous fibre bundle|)}
\subsubsection{Reductions of $f$-submersions}
\begin{defn}
We will say that $\pi\colon (N,q,h)\to (M,g)$ is a metric $q$-submersion if $\pi\colon (N,h)\to (M,g)$ is a Riemannian submersion and if $(N,q,h)$ is the metric $q$-manifold defined by the orthogonal spitting 
$TN=\hor\oplus\ver$, with $\ver=\ker d\pi$ and $\hor=\ver^\perp$.
\end{defn}
\begin{defn} $\bullet$ A map $\pi\colon (N,F)\to  M$ from an $f$-manifold $(N,F)$ to a manifold $M$ is called a \textbf{$f$-submersion} if it is a submersion and $\ker F =\ker d\pi$.\\
$\bullet$ A map $\pi\colon (N,F,h)\to (M,g)$ from a metric $f$-manifold to a Riemannian manifold is called a \textbf{metric $f$-submersion} if $\pi\colon (N,h)\to (M,g)$ is a Riemannian submersion and if $\pi\colon (N,F)\to M$ is an $f$-submersion.
\end{defn}
\begin{rmk}\em
Since $F$ is compatible with $h$, a metric $f$-submersion is also a metric $q$-submersion: the splitting defined by $F$ coincides with the orthogonal splitting defined by the Riemannian submersion.
\end{rmk}
\begin{rmk}\em
A metric $f$-submersion $\pi\colon (N,F,h)\to (M,g)$ is a metric $q$-submersion $\pi\colon (N,h)\to (M,g)$ with an orthogonal complex structure $\Bar J\in \mal C(\Sigma (\pi^*TM))\simeq \mal C(\Sigma(\hor))$.
\end{rmk}
\begin{defn}
Let $\pi_1\colon (N_1,F_1)\to M_1$ and $\pi_2\colon (N_2,F_2)\to M_2$ be two $f$-submersions. Then a map $\Psi\colon (N_1,F_1)\to (N_2,F_2)$ is a morphism of $f$-submersions if it satisfies the following conditions
\begin{description}
\item[(i)] $\Psi$ is $f$-holomorphic: $d\Psi\circ F_1 = F_2\circ d\Psi$ 
\item[(ii)] $\Psi$ is a morphism of submersion: there exists a map $\psi\colon M_1\to M_2$ such that $\pi_2\circ \Psi=\psi\circ\pi_1$.
\end{description}
If $M_1=M_2$ and $\psi=\Id$, we will say that $\Psi$ is a morphism of $f$-submersion over $M$.
\end{defn}
\begin{defn}
$\bullet$ A morphism of metric $f$-submersion is a morphism of $f$-submersion which is also an isometry.\smallskip\\
$\bullet$ A reduction of  $f$-submersion is a morphism of  $f$-submersion which is also  an injective immersion.\smallskip\\
$\bullet$ A reduction of metric $f$-submersion is a morphism of metric $f$-submersion which is injective. We will then say that it is a complete reduction of metric $f$-submersion it is a complete reduction of metric $f$-manifolds.
\end{defn}
\begin{defn}\label{defn-Hom-fibr-f-bund}\index{homogeneous fibre bundle|(}
$\bullet$ We will say that $\pi\colon (N,F)\to M$ is an $f$-fibration if it is an $f$-submersion and a fibration.\smallskip\\
$\bullet$ We will say that $\pi\colon (N,F,h)\to (M,g)$ is a homogeneous fibre $f$-bundle if $\pi\colon(N,h)\to (M,g)$ is a homogeneous fibre bundle and $\pi\colon (N,F,h)\to (M,g)$ is a metric $f$-submersion.
\end{defn}
\begin{prop}\label{prop-comp-red-Hbf}\index{complete reduction|(}
Let $\pi\colon N\to M$ be a homogeneous fibre bundle. Then a homogeneous fibre subbundle $\pi\colon N^\mv\to M$, of $\pi\colon N\to M$ defines a reduction of homogeneous fibre bundle\footnote{See definition~\ref{def-red-Hfb}} \iif the inclusion $i_\mv\colon N^\mv\to N$ is a reduction of $q$-manifolds. Moreover, a metric reduction of homogenous fibre bundle is always a complete reduction.
\end{prop}
\proof 
We have seen in subsection~\ref{Reductions} that $\pi\colon N^\mv\to M$ is a reduction of homogeneous fibre bundle \iif the equation \eqref{hor-v=hor} holds. This proves the first assertion. The second one follows from \eqref{egality}, lemma~\ref{lemma-ivdeB}, and the fact that according to subsection~\ref{Reductions}, the connection $\nabla^c$  leaves invariant the subbundle $\ver^\mv$, so that $D^v$ is reducible in $\ver^\mv$. \comprf \hsq\medskip\\
Therefore, we obtain, more particularly,
\begin{prop}
In the  same situation as in proposition~\ref{prop-comp-red-Hbf}, let us suppose that the homogeneous fibre bundle $\pi\colon N\to M$ admits an orthogonal complex structure $\Bar J\in \mal C(\Sigma(\pi^* TM))$ defining  then a metric $f$-structure and therefore a structure of homogeneous fibre $f$-bundle. Then the homogeneous fibre subbundle $\pi^\mv\colon N^\mv\to M$ inherits also a structure of homogeneous fibre $f$-subbundle defined by $\Bar J^\mv=(i_\mv)^* \Bar J\in \mal C(\Sigma\left( (\pi^\mv)^* TM)\right)$. If this subbundle defines a metric reduction of homogeneous fibre bundle, then it defines also a complete reduction of homogeneous fibre $f$-subbundle.
\end{prop}
\begin{defn}
Let $\pi\colon (N,F,h)\mapsto (M,g)$ be a homogeneous fibre $f$-bundle with a naturally reductive fibre $H/K$. Suppose that\footnote{Keep in mind that  a homogeneous fibre bundle is always reductive, according to corollary~\ref{cory-reductivity}.} $(N,F,h)$ is horizontally of type $\mal G_1$ (respectively of global type $\mal G_1$). Then we call the canonical paracharacteristic (respectively characteristic) connection of the homogeneous fibre $f$-bundle $(N,F,h)$ its paracharacteristic (respectively characteristic) connection extending its vertical canonical connection\footnote{See section~\ref{geometry}.} $\nabla^c$.
\end{defn}
Proposition~\ref{para-chara-red} yields the following.
\begin{prop}\label{hom-fib-bund-para-red}
Let $\pi\colon (N,h)\to (M,g)$ be a homogeneous fibre bundle with a naturally reductive fibre $H/K$, and  $\pi^\mv\colon N^\mv\to M$ a  homogeneous fibre subbundle such that  the inclusion $i_\mv\colon N^\mv\to N$ is a metric reduction of homogeneous fibre bundle\footnote{in the sens of definition~\ref{def-red-Hfb}}. Suppose also that  is given an orthogonal complex structure $\Bar J\in \mal C(\Sigma(\pi^* TM))$ defining  then  structures of homogeneous fibre $f$-bundles: $\ \pi\colon (N,F,h) \to (M,g)$ and $\pi^\mv\colon (N^\mv,F^\mv,h^\mv) \to (M,g)$. Suppose that $(N,F,h)$ is  horizontally of type $\mal G_1$ (of global type $\mal G_1$) then the canonical paracharacteristic (respectively characteristic) connection of $(N,F,h)$ is reducible in $(N^\mv,F^\mv,h^\mv)$ to the canonical paracharacteristic (respectively characteristic) connection of $(N^\mv,F^\mv,h^\mv)$.
\end{prop}
\index{homogeneous fibre bundle|)} \index{complete reduction|)}
%%%%%%%%%%%%%%%%%%%%%%%%%%%%%%%%%%%%%%%%%%%%%%%%%%%%%%%%%%%%%%%%%%%%%
%

\subsubsection{Horizontally K\"{a}hler $f$-manifolds and  horizontally projectible $f$-submersions.}
\index{horizontally K\"{a}hler}
\begin{defn} 
Let $(N,F,h)$ be a metric $f$-manifold. We will say that $(N,F,h)$ is 
\begin{description}
\item[$\bullet$] horizontally Hermitian if  $N_{\Bar J}=0$.
\item[$\bullet$] horizontally K\"ahler if $DF_{|\hor^3}=0$.
\end{description}
\end{defn}
\begin{lemma}\label{lemma-obvious}
Let $\pi\colon (N,F,h)\to (M,g)$ be a metric $f$-submersion.  Let us suppose that $(N,F,h)$ is reductive. Then it is horizontally K\"ahler  \iif
$$
D_H^g \Bar J=0, \forall H\in\hor.
$$
where, as usual $D^g$ denotes the Levi-Civita connections of $g$.
\end{lemma}
\proof
This follows from proposition~\ref{prop-hor-Levi}. \hsq
\begin{thm}\label{thm-Hor-Kahler}
Let $\pi\colon (N,F,h)\to (M,g)$ be a metric $f$-submersion. We suppose  that $(N,F,h)$ is reductive. 
$\bullet$ Then we have
$$
N(X^v,Y^h,Z^h)=\langle (D_{X^v}^g \Bar J) \Bar J Y^h, Z^h\rangle
$$
where $X,Y,Z\in TN$. \\
$\bullet$ Let us suppose moreover that $(N,F,h)$ is horizontally K\"{a}hler. Then $(N,F,h)$ is horizontally Hermitian: $N_{\Bar J}=0$.  In particular, $(N,F,h)$ is horizontally  of type $\mal G_1$. Therefore $(N,F,h)$  admits a  paracharacteristic connection (unique up to $\mal C(\Lm^3 \ver^*)$).\\
Furthermore, let $\Bar T\in\mal C(\Lambda^3\hor^*)$ be the unique horizontal 3-form, such that the torsion $T$ of any metric $f$-connection $\nabla$ in $N$ with a  skew-symmetric component $T_{\hor^3}$, satisfies $T_{\hor^3}=\Bar T$ (see proposition~\ref{prop-horG_1}). Then we have $\Bar T=0$, so that for any such connection $\nabla$, we have $\langle \nabla_{X^h}Y^h, Z^h\rangle=\langle\widetilde{D^g}_{X^h}Y^h, Z^h\rangle$, $\forall X^h,Y^h,Z^h\in\mal C(\hor)$.
\end{thm}
\proof
Concerning the first assertion, it can be proved by applying propositions~\ref{N_F-intermsof-D} and \ref{prop-Levi-dec}. Or more simply by using directly remark~\ref{rmk-nabla_ver}.\\
Further, according to proposition~\ref{N_F-intermsof-D} , $DF_{|\hor^3}=0$ implies $N_{F|\hor^3}=0$. Moreover, $DF_{|\hor^3}=0$ implies also $d\Omega_{F|\hor^3}=0$ which implies $\Bar T=0$ since $d\Omega_{F|\hor^3}=\Bar J{\circact} {\Bar T}$.  \comprf\hsq
\begin{defn}\index{canonical!connection \wrt an $f$-submersion}
Let $\pi \colon (N,F,h)\to (M,g)$ be a metric $f$-submersion. Let us suppose that $(N,F,h)$ is a reductive and horizontally $\mal G_1$. A metric connection $\ovr\nabla$ will be sayed to be the \textbf{canonical connection} on $M$ (\wrt the $f$-submersion) if  its torsion $\ovr T$ satisfies the equation $\widetilde{\ovr T}=T_{|\hor^3}$ where $T$ is the torsion of any paracharacteristic connection (see proposition~\ref{prop-horG_1}) and $\widetilde{\ovr T}$ is the lift in $\hor$ of $\ovr T$. If such a connection exists, we will say that the metric $f$-submersion is \textbf{horizontally projectible}.
\end{defn}
\begin{rmk} \em
The canonical connection is unique when it exists (since it is metric and its torsion is given). Moreover it has a skew-symmetric torsion. The metric $f$-submersion is \textbf{horizontally projectible} \iif  the horizontal 3-form $T_{|\hor^3}$ is projectible to a 3-form on $M$.
\end{rmk}
\begin{prop}\label{nablaJ=0-Hor-G1}
Let $\pi \colon (N,F,h)\to (M,g)$ be a metric $f$-submersion. Let us suppose that there exists a metric connection $\ovr\nabla$ geodesically equivalent to the Levi-Civita connection, i.e. with a skew-symmetric torsion, such that
$$
\ovr \nabla_H \Bar J=0, \quad\forall H\in \hor.
$$
Then $(N,F,h)$ is horizontally of type $\mal G_1$, and horizontally projectible, $\ovr\nabla$ being then the canonical connection of $M$.
\end{prop}
\proof
Any metric connection $\nabla$ on $N$ which satisfies 
$$
\langle \nabla_{H_1} H_2, H_3\rangle = \langle \widetilde {\ovr\nabla}_{H_1} H_2, H_3\rangle,  
$$
satisfies $\nabla F_{|\hor^3}=0 $ and $T_{|\hor^3}= \widetilde{\ovr T}$. Then proposition~\ref{prop-horG_1} allows us to conclude. \comprf \hsq \medskip\\
In particular, we have

\begin{cory}
Let $\pi \colon (N,F,h)\to (M,g)$ be a metric $f$-submersion. Let us suppose that $(N,F,h)$ is reductive and horizontally k\"{a}hler. Then it is horizontally projectible and the canonical connection on $M$ is its Levi-Civita connection $D^g$.
\end{cory}
\proof This follows immediately from lemma~\ref{lemma-obvious}  and proposition~\ref{nablaJ=0-Hor-G1}. More directly, according to theorem~\ref{thm-Hor-Kahler}, we have $T_{|\hor^3}=0$ and therefore this 3-form is projectible. \comprf\hsq

%%%%%%%%%%%%%%%%%%%%%%%%%%%%%%%%%%%%%%%%%%

\subsubsection{The example of a naturally reductive homogeneous space}\label{example-homogenenous-G1-red}
\index{canonical!connection, $G$-invariant|(}
\begin{prop}\label{thm-Hom-G1-red}
Let $N=G/K$ be a Riemannian homogeneous space and $\g=\kk\oplus\nk$ a reductive decomposition of $\g$. Let us suppose that $\nk$ admits an $\Ad K$-invariant decomposition $\nk =\mk\oplus \pk$, defining then a splitting $TN=\hor\oplus\ver$, where $\hor=[\mk]$ and $\ver=[\pk]$. Let us suppose that there exists on $\mk$ an $\Ad K$-invariant complex structure $\Bar J_0$, defining  then  an $f$-structure $F$ on $N$. Then for any $G$-invariant metric $h$ for which $\Bar J_0$  and the decomposition $\nk =\mk\oplus \pk$ are orthogonal (such a metric always exists),  $(N,F,h)$ is a metric $f$-manifold.\\
Furthermore, let us suppose that $N=G/K$ is naturally reductive, and that one can choose a naturally reductive metric $h$ as above\footnote{which means that denoting by $G(\nk)$, the compact subgroup in $GL(\nk)$ generated by $\Lambda_\nk(\nk) :=\{[\ad_\nk(X)]_\nk, X\in\nk\}\subset\mak{gl}(\nk)$, and by $\langle G(\nk), I_0\rangle$ the closed subgroup generated by $G(\nk)$ and $I_0:=\Bar J_0\oplus -\Id_\pk$, then $\langle G(\nk), I_0\rangle/G(\nk)$ is compact.}, then $(N,F,h)$ is reductive and of global type $\mal G_1$. Moreover, the canonical connection $\nabla^0$ is a characteristic connection. Therefore, in this case this  characteristic connection $\nabla^0$ has a parallel torsion  $\nabla^0 T=0$.
\end{prop}
\textbf{Proof.} The naturally reductivity means exactly that the torsion of the canonical connection is skew-symmetric. \comprf\hsq
\begin{prop}\label{hom-curvature}
Let us consider the situation described by (the 2 first sentences of) proposition~\ref{thm-Hom-G1-red}. Then the horizontal curvature and its linear representation, and the horizontal curvature operator are given (in terms of their lifts in $G$) by
\begin{eqnarray*}
\widetilde{\Phi( X,Y)} & = & - \left[ X_\mk, Y_\mk \right]_{\pk} \\
\widetilde{\rho(V)} & = & -\adm V_\pk\\
\widetilde{\Bar R(X,Y)Z} & =  & \adm \left[ X_\mk, Y_\mk \right]_{\pk} 
\end{eqnarray*}
In the second equality, we have supposed that $G/K$ is endowed with a naturally reductive metric.
\end{prop}
\begin{defn}\label{defn-hinducesg}
Let $N=G/K$ be a  Riemannian homogeneous space and $H\supset K$  a subgroup of $G$ such that $M=G/H$ is Riemannian.  Then we will say that a  $G$-invariant  metric $h$ on $N$ induces a  $G$-invariant metric $g$ on $M$ if  $\pi\colon (G/K,h)\to (G/H,g)$ is a Riemannian submersion and  therefore a homogeneous fibre bundle. We will then say that $h$ is projectible on $M$ and that $g$ is induced by $h$. 
\end{defn}
\begin{rmk}\label{rmk-hinducesg}\em
If $\pi\colon (G/K,h)\to (G/H,g)$ is a Riemannian submersion, then we are in the situation described by \ref{homspacefibr}. Indeed if $\mk$ is an $\Ad H$-invariant complement of $\hk$ in $\g$  and $\pk$ is an $\Ad H$-invariant complement of $\kk$ in $\hk$, then $\nk=\pk\oplus\mk$ is an $\Ad K$-invariant complement of $\kk$ in $\g$. Moreover, the $\Ad K$-invariant inner product $\langle\cdot, \cdot \rangle_\nk : = h_{y_0}$, where $y_0=1.K$, defining the $G$-invariant metric $h$, satisfies the equation \eqref{eq-met-p+m} i.e. its restriction $\langle \cdot,\cdot \rangle_\mk$ to $\mk$ is $\Ad H$-invariant. Conversely, if there exists such $\pk$ and $\mk$ such that $\langle \cdot,\cdot \rangle_\mk := h_{y_0|\mk\times\mk}$ is $\Ad H$-invariant then $h$ is projectible on $M$. 
\end{rmk} 
\begin{prop}\label{prop-2ksym-G1&pure}
Let us suppose that $N=G/K$ is a (locally) $2k$-symmetric space endowed with its canonical $f$-structure $F$ and its canonical connection $\nabla^0$. Let us suppose that $N=G/K$ is naturally reductive.\\ 
$\bullet$ Then there exists a $G$-invariant naturally reductive metric $h$ on $N$ such that $\tau_{|\nk}$ is orthogonal  and thus which is compatible with $F$. Then $(N,F,h)$ is a reductive metric $f$-manifold  of global type $\mal G_1$, and its horizontal curvature $\Phi$ is pure.\\
$\bullet$ Furthermore, any naturally reductive metric $h$ chosen as above, induces a $G$-invariant Riemannian metric $g$ on the $k$-symmetric space $M=G/H$, and $\pi\colon (G/K,h)\to (G/H,g)$ is a homogeneous fibre $f$-bundle.
\end{prop}
\proof
The existence of $h$ follows from lemma~\ref{Gn-F-compact}, in the Appendix. The fact that $(N,F,h)$ is a reductive metric $f$-manifold  of global type $\mal G_1$ follows from proposition~\ref{thm-Hom-G1-red}. Moreover, $\Phi$ is pure because of proposition~\ref{hom-curvature} and since
$$
[X_\mk,Y_\mk]_\pk=[X_\mk,Y_\mk]_{\g_k}=\sum_{1\leq |j| \leq k-1} [X_j,Y_{-j+k}], \quad \forall X_\mk,Y_\mk\in \mk,
$$
 and by definition of $\undj$.\\
The second point of the proposition follows from proposition~\ref{prop-nat-red-induce-appendix}, in the Appendix, definition~\ref{defn-hinducesg} and definition~\ref{defn-Hom-fibr-f-bund}. \comprf \hsq
\index{canonical!connection, $G$-invariant|)}\index{connection!paracharacteristic|)}
\subsubsection{The example of the twistor space $\mal Z_{2k}^\alpha(M)$.}\label{example-zdk}
\index{twistor|(}
Let $(M,g)$ be a (even dimensional) Riemannian manifold endowed with a metric connection $\ovr\nabla$. Then we consider the homogeneous fibre bundle $\pi\colon (\mal Z_{2k}^\alpha(M),h)\to (M,g)$  defined by the Riemannian vector bundle $(TM,g,\ovr\nabla$). (See  sections~\ref{preliminaires-examples} (example~\ref{zdkE}) and \ref{subsect-zdkE}.)
%and studied in sections~\ref{preliminaires-examples} (example~\ref{zdkE}) and \ref{subsect-zdkE}, where we take as the Riemannian vector bundle $E$, the tangent bundle $TM$ endowed with $g$ and $\ovr\nabla$.
Let us consider also its canonical  $2k$-structure $\mJ\in\mal C(\mal Z_{2k}^\alpha(\pi^*TM))$, to which corresponds the orthogonal complex structure $\und\mJ\in\mal C \left( \Sigma(\pi^*TM)\right) \cong \mal C \left( \Sigma(\hor)\right) $. This complex structure defines then a metric $f$-structure $\mal F$ on $\mal Z_{2k}^\alpha(M)$. Then the twistor bundle $(\mal Z_{2k}^\alpha(M),\mal F,h)$ is a reductive metric $f$-manifold, more precisely a homogeneous fibre $f$-bundle.

\begin{prop}\label{nabla-undJ=0}
Let $(M,g)$ be a Riemannian manifold endowed with a metric connection $\ovr\nabla$. Then we have
$$
\ovr\nabla_H \und\mJ=0, \quad\forall H\in \hor.
$$
\end{prop}
\proof
We have the following sequence of morphisms of bundle over $M$ which are $SO(2n)$-invariant projections on the fibres:
$$
\begin{CD}
Q=\mal {SO}(TM) @> \pi_{\mal Z} >> Q/\U_0(J_0^\alpha) =\zdk^\alpha(M) @>\mal P >> \Sigma (M)=Q/\U(\undj_0^\alpha)
\end{CD}
$$
In particular,  the structure of homogeneous fibre bundle of $N_\Sigma=\Sigma (M)$ (defined by $\ovr\nabla$) is  the image by the projection $\mal P$ of the structure of homogeneous fibre bundle of $N_{\mal Z}=\zdk^\alpha(M)$ (defined by $\ovr\nabla$):
$$
\begin{CD}
TQ = \hor^Q\oplus\ver^Q @> \left(\pi_{\mal Z}\right)_*  >> TN_{\mal Z} =\hor^\mZ \oplus \ver^\mZ @> (\mal P)_* >> TN_\Sigma = \hor^\Sigma\oplus\ver^\Sigma.
\end{CD}
$$
In particular we have $\ovr\nabla \und\mJ = \mal P _*(\ovr\nabla \mJ)$ and therefore
$$
\ovr\nabla_H \und\mJ=0, \quad\forall H\in \hor,
$$
according to theorem~\ref{entermdeJ2} (i). \comprf \hsq 

\begin{thm}\label{thm-Phi-pure-Zdk}
Let $(M,g)$ be a Riemannian manifold endowed with a metric connection $\ovr\nabla$. 
Let us suppose that $\ovr\nabla$ is geodesically equivalent to the Levi-Civita connection, i.e. it admits a skew-symmetric torsion. Then the twistor bundle $(\mal Z_{2k}^\alpha(M),\mal F,h)$ is horizontally of type $\mal G_1$, and horizontally projectible, $\ovr\nabla$ being the canonical connection of $M$. 
\end{thm}
\proof
This follows immediately from proposition~\ref{nabla-undJ=0} and proposition~\ref{nablaJ=0-Hor-G1}.\hsq
\begin{cory}
Let $(M,g)$ be a Riemannian manifold endowed with a metric connection $\ovr\nabla$. If $\ovr\nabla=D^g$ is the Levi-Civita connection, then $(\mal Z_{2k}^\alpha(M),\mal F,h)$ is horizontally K\"ahler. 
\end{cory}
\proof
This follows immediately from proposition~\ref{nabla-undJ=0} and  lemma~\ref{lemma-obvious}. \hsq
\subsubsection{The example of the twistor space $\mal Z_{2k,j}^\alpha(M,J_j)$.}\label{example-zdkj}
Let $(M,g)$ be a Riemannian manifold endowed with a metric connection $\ovr\nabla$ and $J_j$ a section of $\left( \mal Z_{2k}^\alpha(M)\right)^j$. Then we consider the  the homogeneous fibre bundle $\pi\colon (\mal Z_{2k,j}^\alpha(M,J_j),h)\to (M,g)$ defined by the Riemannian vector bundle $(TM,g,\overline\nabla^{[j]}$). We denote by $\mal F_j$ the canonical $f$-structure defined as above in \ref{example-zdk}. Then the twistor bundle $(\mal Z_{2k,j}^\alpha(M,J_j),\mal F_j,h)$ is a reductive metric $f$-manifold, more precisely a homogeneous fibre $f$-bundle.
\begin{thm}
Let $(M,g)$ be a Riemannian manifold endowed with a metric connection $\ovr\nabla$ and $J_j$ a section of $\left( \mal Z_{2k}^\alpha(M)\right)^j$. Let us suppose that $\ovr\nabla^{[j]}$ is geodesically equivalent to the Levi-Civita connection, i.e. it admits a skew-symmetric torsion.\\
$\bullet$  Then the twistor bundle $(\mal Z_{2k,j}^\alpha(M,J_j),\mal F_j,h)$ is horizontally  of type $\mal G_1$, and horizontally projectible, $\ovr\nabla^{[j]}$ being the canonical connection of $M$.\medskip\\
$\bullet$ Furthermore, if $\ovr\nabla J_j=0$ so that $ (\mal Z_{2k,j}^\alpha(M,J_j),\mal F_j,h)\to (M,g)$ is a reduction of $(\mal Z_{2k}(M),\mal F,h)\to (M,g)$, the  paracharacteristic connection of $(\mal Z_{2k}(M),\mal F,h)$ is reducible in  $(\mal Z_{2k,j}^\alpha(M,J_j),\mal F_j,h)$ and its reduction is the paracharacteristic connection of $(\mal Z_{2k,j}^\alpha(M,J_j),\mal F_j,h)$. 
\end{thm}
\proof 
For the first point, proceed as in \ref{example-zdk}. For the second point, apply proposition~\ref{hom-fib-bund-para-red}. \hsq 

\begin{cory}
Let $(M,g)$ be a Riemannian manifold endowed with a metric connection $\ovr\nabla$ and $J_j$ a section of $\left( \mal Z_{2k}^\alpha(M)\right)^j$.  If in particular $\ovr\nabla=D^g$ is the Levi-Civita connection, then  $(\mal Z_{2k,j}^\alpha(M,J_j),\mal F_j,h)$ horizontally K\"ahler. 
\end{cory}
\proof Proceed as in \ref{example-zdk}. \hsq

%%%
%%

\subsubsection{The reduction of homogeneous fibre bundle $\mijo\colon G/G_0\hookrightarrow \mZ_{2k,2}^{\alpha_0}(G/H,J_2)$.}\label{subsub-red-mijo}
\index{canonical!connection, $G$-invariant|(}\index{canonical!embedding|(}\index{connection!paracharacteristic|(}
We consider the situation described by \ref{geom-interpret-sub1} and  use the same notations.\\
We have seen  in \ref{geom-interpret-sub1} that $\mijo\colon G/G_0\longrightarrow\mZ_{2k,2}^{\alpha_0}(G/H,J_2)$ is a reduction of homogeneous fibre bundles (theorem~\ref{thm-inj-morph-bund}). Moreover by definition of $\mijo$ and $\mJ$, this is also an $f$-immersion. Therefore, since $(G/G_0,F,h)$ and\footnote{We denote the metric in $\mZ_{2k,2}^{\alpha_0}(G/H,J_2)$ by $h_2$ to differentiate it from the metric in $G/G_0$. See the convention of notations in \ref{geom-interpret-sub1}.} $(\mZ_{2k,2}^{\alpha_0}(G/H,J_2),\mal F_2,h_2)$ are both metric $f$-manifolds, it follows that $\mijo$ is a reduction of homogeneous fibre $f$-bundle. Moreover, we suppose that the $G$-invariant metric $g$ on $G/H$, induced\footnote{In the sens of definition~\ref{defn-hinducesg}, see remark~\ref{rmk-hinducesg}, and also \ref{homspacefibr}.} by the metric $h$ on $G/G_0$, is naturally reductive. This is the case in particular if $h$ is naturally reductive (see proposition~\ref{prop-2ksym-G1&pure}).\\ 
Therefore according to proposition~\ref{hom-fib-bund-para-red}, we have
\begin{prop}\label{red-mijo-canonical}
The canonical paracharacteristic connection of $\mZ_{2k,2}^{\alpha_0}(G/H,J_2)$ is reducible in $G/G_0$ and its reduction is the canonical connection of $G/G_0$.
\end{prop}
\begin{rmk}\em
To apply proposition~\ref{hom-fib-bund-para-red}, we need a priori a metric reduction of homogeneous fibre bundle, i.e. the inner product $\langle \cdot,\cdot\rangle_\pk$ on the fibre $H/G_0$ of the fibration $G/G_0\to G/H$, must be in the form described by remark~\ref{rmk-metric-inducedbyTr}. That is to say the  inner product on $\nk$ defining $h$ on $G/G_0$ is on the form $\langle\cdot,\cdot\rangle_\mk + \langle \cdot,\cdot\rangle_\pk$ with $\langle\cdot,\cdot\rangle_\mk$ $\Ad H$-invariant and naturally reductive (and invariant by $\taum$), and $\langle \cdot,\cdot\rangle_\pk$ given by remark~\ref{rmk-metric-inducedbyTr}. But in fact, proposition~\ref{red-mijo-canonical} still holds for any $G$-invariant metric $h$ which induces  a naturally reductive metric $g$ on $M$: in other words $\langle \cdot,\cdot\rangle_\pk$ can be any $\Ad K$-invariant inner product on $\pk$. Indeed, it suffices to remark that $\nabla^0$ is always the paracharacteristic connection of $(G/G_0,F,h)$ for any $h$ inducing a naturally reductive metric $g$ on $G/H$. Indeed, we have $(\adm V )\undj_0 = -\undj_0(\adm V)$, $\forall V\in\pk=\g_k$, by definition of $\undj_0$. Therefore, the component $T_{|\ver\times\hor\times\hor}$ of $\nabla^0$ is pure so that $\nabla^0$ is always the paracharacteristic connection of $(G/G_0,F,h)$  and proposition~\ref{red-mijo-canonical} holds in general for all possible choice of $\langle \cdot,\cdot\rangle_\pk$.
\end{rmk}
\index{twistor|)}\index{canonical!connection, $G$-invariant|)}\index{canonical!embedding|)}
%
%
%
%%%%%%%%%%%%%%%%%%%%%%%%%%%%%%%%%%%%%%%%%%%%%%%%%%%%%%%%%%%%%%%%%%%%%%%%%%%
%
%
\subsection{Stringy Harmonic maps in $f$-manifolds.}\label{Stringy-Harm-f-manif}
\index{stringy harmonic|(}
\subsubsection{Definitions}
We have defined the notion of stringy harmonic maps in the context of almost complex manifolds (endowed with a linear connection) and we have seen that it corresponds to a  generalisation of harmonic maps. Now, we will extend this notion of stringy harmonicity to  $f$-manifolds endowed with a linear connection.
Indeed, the preliminary study of the maximal determined system done in section~\ref{eq-maximal-deter},  leads us to introduce the following generalisation of (stringy) harmonic maps.
\begin{defn}
Let  $(N,F)$ be an $f$-manifold with $\nabla$ a linear connection. Then we will say that a map $f\colon L\to N$ from a Riemann surface into $N$ is \emph{stringy harmonic} if it is solution of \emph{the stringy harmonic maps equation}:
$$
-\tau_g(f) + (F\bullet T)_g(f)=0
$$
\end{defn}
\begin{defn}
Let $(N,F)$ be an $f$-manifold with $\nabla$ a linear connection. Then we will say that a map $f\colon L\to N$ from a Riemann surface into $N$ is \emph{$\star$-stringy harmonic} if it is solution of \emph{the modified stringy harmonic maps equation}:
$$
-\tau_g(f) + (F\star T)_g(f)=0.
$$
\end{defn}
\begin{prop}
Let $(N,F)$ be an $f$-manifold with $\nabla$ a linear $f$-connection. Let us suppose that $\mR_\ver=0$. Then, setting $\Bar T=T_{\hor^2}^\hor$,  the equation of stringy harmonicity, \wrt $\nabla$, for maps $f\colon L\to N$ is written
$$
\left\lbrace \begin{array}{l}
\displaystyle 
-*\tau^v(f) + \dfrac{1}{2} \Phi(df\wedge \Bar J df)=0\\
\displaystyle 
-*\tau^h(f) + \dfrac{1}{2}(\Bar J\cdot \Bar T)(df\wedge df) +  \dfrac{1}{2} N_F(d^v f\wedge\Bar J d^h f)=0
\end{array}
\right.
$$
whereas the equation of $\star$-stringy harmonicity is written:
$$
\left\lbrace \begin{array}{l}
\displaystyle
-*\tau^v(f) + \dfrac{1}{2} \Phi( df\wedge \Bar J df)=0\\
\displaystyle 
-*\tau^h(f) + \dfrac{1}{2}(\Bar J\star \Bar T)(df\wedge df) +  \dfrac{1}{2} N_F(d^v f\wedge\Bar J d^h f)=0

\end{array}
\right.
$$
where $\tau^v(f)$, $\tau^h(f)$ are respectively the vertical and horizontal component of the tension field \wrt $\nabla$ (which coincide also with the vertical and horizontal tension fields respectively, since $\nabla$ preserves the splitting $TN=\ver\oplus\hor$). Moreover, as usual, $*$ is the Hodge operator  for some Hermitian metric on $L$.
\end{prop}
\proof
Use theorem~\ref{thm-nablaOmega_F=0}. \hsq

\begin{rmk}\em
We see that stringy harmonicity in $(N,F)$ can be viewed  as a kind of coupling between some equation of $\Bar J$-stringy harmonicity and the equation of vertical harmonicity. We will come back to this in details elsewhere \cite{ki7,ki8}.
\end{rmk}
\begin{prop}\label{prop-same-stringy-harm-if-red}\index{global type $\mal G_1$|(}\index{connection!characteristic|(}
Let $i_\mv\colon (N^\mv,h^\mv,F^\mv)\to (N,F,h)$ be a complete reduction of metric $f$-manifolds. Let $\nabla^c\in\con(\ver)$ reducible in $\ver^\mv$, and $\nabla^{c,\mv}$ its reduction. Let us suppose that $(N,F,h)$ is horizontally of type $\mal G_1$ and almost reductive. Let $f\colon L\to N^\mv$ be a map.\\
$\bullet$ The \emph{($\star$\,-)}stringy harmonicity of $f$ \wrt to the paracharactersitic connection $\nabla$  of $(N,h,F)$ extending $\nabla^{c}$, is equivalent to the \emph{($\star$\,-)}stringy harmonicity of $f$ \wrt to the paracharactersitic connection of $(N^\mv,h^\mv,F^\mv)$ extending $\nabla^{c,\mv}$.\\
$\bullet$ In the same way, if $(N,F,h)$ is reductive and of global type $\mal G_1$, then the \emph{($\star$\,-)}stringy harmonicity of $f$ \wrt to the charactersitic connection $\nabla$  of $(N,h,F)$ extending $\nabla^{c}$, is equivalent to the \emph{($\star$\,-)}stringy harmonicity of $f$ \wrt to the charactersitic connection of $(N^\mv,h^\mv,F^\mv)$ extending $\nabla^{c,\mv}$.\smallskip\\
In particular, these properties hold when $i_\mv\colon (N^\mv,h^\mv,F^\mv)\to (N,F,h)$ is a reduction of homogeneous fibre bundle as described in proposition~\ref{hom-fib-bund-para-red}.
\end{prop}
\proof
This results immediately from proposition~\ref{para-chara-red} and lemma~\ref{torsion-red-related}. \comprf \hsq
\begin{prop}
Let $i_\mv\colon (N^\mv,h^\mv,F^\mv)\to (N,h,F)$ be a complete reduction of metric $f$-manifolds.\\
Let us suppose that $(N,F,h)$ is horizontally of type $\mal G_1$ and almost reductive. Let $f\colon L\to N^\mv$ be a map.\\
$\bullet$ The \emph{($\star$\,-)}stringy harmonicity of $f$ \wrt to one paracharactersitic connection $\nabla$  of $(N,h,F)$ is equivalent to the \emph{($\star$\,-)}stringy harmonicity of $f$ \wrt to any paracharactersitic connection of $(N^\mv,h^\mv,F^\mv)$.\\
$\bullet$ In the same way, if $(N,F,h)$ is reductive and of global type $\mal G_1$, then the \emph{($\star$\,-)}stringy harmonicity of $f$ \wrt to one charactersitic connection $\nabla$  of $(N,h,F)$ is equivalent to the \emph{($\star$\,-)}stringy harmonicity of $f$ \wrt to any charactersitic connection of $(N^\mv,h^\mv,F^\mv)$.
\end{prop}
\proof
We deduce from proposition~\ref{prop-tensionfield=}, theorem~\ref{thm-proj-chara} and equation \eqref{hor-v=hor}, that the terms in $\nabla^* df + (F\bullet T)_g(f)$ which could be no reducible in $N^\mv$ are only those corresponding to the vertical terms $(F\bullet T)_{\hor\times\ver\times\ver}$ and $(F\bullet T)_{\ver\times\hor\times\ver}$ equal to $\mR_\ver(Y^v,Z^v,\bar J X^h)$ and $\mR_\ver(Z^v,X^v,\Bar J Y^h)$ respectively. These terms do not depend on the choice of $\nabla$ (and $\nabla^\mv$ for the corresponding terms in $N^\mv$) and from proposition~\ref{prop-same-stringy-harm-if-red}, we know that for a convenient choice of $\nabla$ and $\nabla^\mv$, these terms are reductible, i.e. if $B$ (resp. $B^\mv$) is one of these terms then $B_g^\mv(f)=B_g(f)$. \comprf \hsq
\index{connection!paracharacteristic|)}

\subsubsection{The closeness of the 3-forms $F\bullet T$ and $F\star T$.}\label{dH=0}
Let us see under which conditions on a reductive metric $f$-manifold of global type $\mal G_1$, the 3-forms $F\bullet T$ and $F\star T$ defined by  one characteristic connection  are closed.
\begin{prop}\label{different3-form}
Let $(N,F,h)$ be a reductive metric $f$-manifold of  global type $\mal G_1$.  Let $\nabla$ be some characteristic connection on $N$. Let us consider the following statements:
\begin{description}
\item[(i)] The  3-forms $H=F\bullet T$ is closed.
\item[(ii)] The  3-forms $H^\star=F\star T$ is closed.
\item[(iii)] The horizontal 3-form $F\cdot N_F$ is closed.
\item[(iv)] The 3-form $F\cdot N_F -\dfrac{1}{2}F\circact \left(\mrm{Skew}(\Phi) + \mrm{Skew}(\mR_\ver)\right)$ is closed. 
\item[(v)] The 3-form $\dfrac{1}{2} \left( F\cdot N_F -F\circact \left(\mrm{Skew}(\Phi) + \mrm{Skew}(\mR_\ver)\right)\right) $ is closed. 
\item[(vi)] The 3-form $ F\circact \left(\mrm{Skew}(\Phi) + \mrm{Skew}(\mR_\ver)\right)$ is closed.   
\end{description}
Then we have the following equivalences: (i) $\Leftrightarrow$ (iv) and  (ii) $\Leftrightarrow$ (v). Moreover any couple of statements in $\{ (i), (ii), (iii)\}$ imply the third one. Any couple of statements in $\{ (iii),(iv),(v),(vi)\}$ imply
the two others. In other words identifying (i) with (iv) and (ii) with (v) (since they are respectively equivalent) in the set of six elements $\{(i),\ldots, (vi)\}$, then in the obtained set of four elements, any couple of two statements is equivalent to the  couple of the two others.
\end{prop}
\textbf{Proof.}
This follows from the following identities: 
\begin{eqnarray*}
\frac{1}{2}F\cdot N_F & = & F\bullet T - F\star T\\
F\bullet T & = &  F\cdot N_F - \dfrac{1}{2}F\circact \left(\mrm{Skew}(\Phi) + \mrm{Skew}(\mR_\ver)\right) -d\Omega_F \quad \text{(proposition~\ref{prop-F-bullet-T})} \\
F\star T & = & F\bullet T - \dfrac{1}{2}F\cdot N_F  = \dfrac{1}{2}F\cdot N_F - \dfrac{1}{2}F\circact \left(\mrm{Skew}(\Phi) + \mrm{Skew}(\mR_\ver)\right) -d\Omega_F.
\end{eqnarray*}
\hfill$\square$
%
%F\cdot N_F +  F\circact T - \dfrac{1}{2}F \circact (T - T_{|\hor^3})=
%
\begin{defn}
Let $(N,F,h)$ be a reductive metric $f$-manifold of  global type $\mal G_1$. We will say that $(N,F,h)$ \textbf{has a closed stringy structure} if the 3-forms $F\cdot N_F$ and $ F{\circact} \left(\mrm{Skew}(\Phi) + \mrm{Skew}(\mR_\ver)\right)$ are closed. This is equivalent to say that  the two 3-forms $H$ and $H^\star$ are closed.
\end{defn}
\paragraph{Characteristic connections with  parallel torsion.} 
In this paragraph, $(N,F,h)$ is a reductive metric $f$-manifold of global type $\mal G_1$. Let us suppose that one of its  characteristic connections, $\nabla$, has a parallel torsion $\nabla T=0$.\\
Now, let us compute all the components in each space $(\Lambda^p\hor^*)\wedge(\Lambda^q\ver^*)$, $p+q=4$, of the exterior differential $d\alpha$ of the different 3-forms $\alpha$ we have considered in proposition~\ref{different3-form}.
\begin{lemma}\label{d-alpha}
Let $\alpha$ be a $\nabla$-parallel 3-form. Then
$$
d\alpha(X,Y,V,Z)=\underset{X,Y,Z}{\mak S} \alpha(T(V,Z), X,Y) + \alpha(T(X,Y),V,Z)
$$
In particular, if $\alpha$ is horizontal then 
\begin{eqnarray*}
d\alpha_{|\Lambda^4\ver} & = & 0   \\ 
d\alpha_{|(\Lambda^3\ver)\wedge\hor} &= & 0  \\
d\alpha_{|(\Lambda^2\ver)\wedge\hor^2}(V_0,V_1,H_2,H_3) & =  & \alpha(\mR_\ver(V_0,V_1),H_2,H_3) \\
 d\alpha_{|\ver\wedge(\Lambda^3\hor)}(V_0,H_1,H_2,H_3) & = &  \underset{i,j,k}{\mak S} \alpha(T^\hor(V_0,H_i), H_j,H_k)
\end{eqnarray*}
\end{lemma}
\textbf{Notation} We set $\mal S(\hor,\ver)=\mal S(\hor\times\hor\times\ver) \oplus \mal S(\hor\times\ver\times\ver)$. Then for any $B\in\mal T$, we set $\mathring B=B_{|\mal S(\hor,\ver)}$. Let us remark that we have 
$$
-F\circact \, \mathring T= d\Omega_{F|\mal S(\hor,\ver)} \quad \text{and } -\Bar J\circact \Bar T= d\Omega_{F|\hor^3}.
$$
\begin{lemma}\label{lemma d=d}
We have the following identity:
$$
-d(\Bar J\circact \Bar T)= d(F\circact \,\mathring T)
$$

\end{lemma}
\textbf{Proof.} We have 
$$
0= -d (d\Omega_F)= d (\Bar J\circact \Bar T) + d(F\circact \, \mathring T).
$$
This completes the proof.\hfill $\square$
\begin{prop}\label{prop-Hor4}\index{cyclic derivation of the horizontal curvature}
The following identity holds:
$$
d(F\cdot N_F)_{|\hor^4} = d \left( F\circact (\mrm{Skew}(\Phi) + \mrm{Skew}(\mR_\ver) )\, \right)_{|\hor^4} 
$$
Therefore the following statements are equivalent
$$
\begin{array}{lcl}
\textbf{\em (i) }\ dH_{|\hor^4}=0  &   &   \textbf{\em (iv) }\   \mak S_{X,Y,Z} \,\Bar{\mrm A}(X,Y)Z=0  \\
\textbf{\em (ii) }\  dH_{|\hor^4}^\star=0   &   &  \textbf{\em (v) }\  \mak S_{X,Y,Z} \,\Bar{\mrm A}^{(-)}(X,Y)Z=0 .\\ 
 \textbf{\em (iii) }\  d(F\cdot N_F)_{|\hor^4}=0   &  & 
\end{array}
$$
We will  say that $\Bar J$ is a cyclic derivation of  the horizontal curvature when (iv) holds. $\centerdot $
\end{prop}
\textbf{Proof.}
We consider 4-forms on  $\hor$ and denote  by $(X,Y,V,Z)$ the variable in $\hor^4$. Proceeding as in the proof of proposition~\ref{prop-dH_J=0}, we obtain according to lemma~\ref{d-alpha}
\begin{multline*}
d (\Bar J\circact \Bar T) = \underset{X,Y,Z}{\mak S} (\Bar J\circact \Bar T)\left( X,Y, \Bar T (V,Z)\right) +  (\Bar J\circact \Bar T)\left( V,Z, \Bar T (X,Y)\right) \\
=\underset{X,Y,Z}{\mak S} \sum_{(X,Y)\rightleftarrows(V,Z)}\langle (\Bar J\circact \Bar T)(X,Y), \Bar T (V,Z)\rangle\\
= \underset{X,Y,Z}{\mak S} \sum_{(X,Y) \rightleftarrows(V,Z)}  4\langle  {\Bar T}^{--}(X,Y), \Bar J{\Bar T}^{++} (V,Z)\rangle + 2 \langle  {\Bar T}^{--}(X,Y), \Bar J{\Bar T}^{1,1} (V,Z)\rangle - 2\langle  {\Bar T}^{++}(X,Y), \Bar J{\Bar T}^{1,1} (V,Z)\rangle. 
\end{multline*}
Moreover, using the same arguments as in the end of the proof of  proposition~\ref{prop-dH_J=0}, we can conclude that the two last terms inside the sum vanish. More precisely, using the fact  that $\Bar J \Bar T^{--}$ is a 3-form, we can write  the second term inside the sum in two different forms which have different types in $\otimes^4\hor^*$, therefore, this terms vanishes. Idem for the third term, by using the  corollary~\ref{cory-B-eps,eps'} applied to $\Bar T^{++}$. Therefore
\begin{equation}\label{4,0+0,4}
d (\Bar J\circact \Bar T) = \underset{X,Y,Z}{\mak S} \sum_{(X,Y) \rightleftarrows(V,Z)}  4\langle  {\Bar T}^{--}(X,Y), \Bar J{\Bar T}^{++} (V,Z)\rangle.
\end{equation}
Furthermore, we also have according to lemma~\ref{d-alpha}
\begin{eqnarray*}
d(F\cdot N_F) & = & \underset{X,Y,Z}{\mak S} (F\cdot N_F)\left( X,Y, \Bar T (V,Z)\right)  + (F\cdot N_F)\left( V,Z, \Bar T (X,Y)\right) \\
 & = & \underset{X,Y,Z}{\mak S} \sum_{(X,Y)\rightleftarrows(V,Z)} 4\langle  \Bar J T^{--}(X,Y),\Bar T(V,Z)\rangle\\
 & = & \underset{X,Y,Z}{\mak S} \sum_{(X,Y)\rightleftarrows(V,Z)} 4\langle  \Bar J T^{--}(X,Y),{\Bar T}^{**}(V,Z)\rangle\\
 & = & \underset{X,Y,Z}{\mak S} \sum_{(X,Y)\rightleftarrows(V,Z)} 4\langle  \Bar J T^{--}(X,Y),{\Bar T}^{++}(V,Z)\rangle
\end{eqnarray*}
Hence, we obtain
$$
 d(F\cdot N_F) = - d (\Bar J\circact \Bar T). 
$$
This proves the first assertion according to lemma~\ref{lemma d=d}. Moreover, using lemma~\ref{d-alpha}, we compute that 
\begin{equation}\label{A-def-4form}
d(F\circact\, \mathring T)_{|\hor^4}   =   \underset{X,Y,Z}{\mak S}  \langle (\Bar J \circact\,\Phi) (X,Y), \Phi(W,Z)\rangle +  \langle (\Bar J \circact\,\Phi) (W,Z), \Phi(X,Y)\rangle.
\end{equation}
Furthermore, introducing the linear representation $\rho\colon\ver\to\so(\hor)$, we have
\begin{eqnarray}
\langle (\Bar J \circact\,\Phi) (X,Y), \Phi(W,Z)\rangle +  \langle (\Bar J \circact\,\Phi) (W,Z), \Phi(X,Y)\rangle &  =  &    \langle \rho\left( \Bar J \circact\,\Phi (X,Y)\right)W, Z\rangle -  \langle [J,\rho(\Phi(X,Y))] W,Z \rangle \nonumber\\
 &  =  &  \langle \Bar{\mrm A}(X,Y)W,Z\rangle = - \langle \Bar{\mrm A}(X,Y)Z,W\rangle\label{4,0 + 0,4}
\end{eqnarray}
Finally, according to \eqref{4,0+0,4}, we see that the 4-form $d (\Bar J\circact \Bar T)$ is of type $(4,0) +(0,4)$, and thus 
 $$
\underset{X,Y,Z}{\mak S} \langle \Bar{\mrm A}(X,Y)Z,W\rangle =\underset{X,Y,Z}{\mak S}\langle \Bar{\mrm A}^{(-)}(X,Y)Z,W\rangle.
$$
This completes the proof.\hfill$\square$\bigskip\\
\textbf{Notation} Given $B\in\mal C (\Lambda^2 T^*N\otimes\hor)$, we denote simply  $\im B=\{B(X,Y)\in \hor, X,Y\in TN\}\subset \hor$. In particular, we have $\im N_{\Bar J}=N_{\Bar J}(\hor,\hor)\subset\hor$ and $\im\mR_\ver=\mR_\ver(\ver,\ver)\subset\hor$.\\
 Moreover, we will also use the notations $\ker C=\{X\in\hor|\, C(X,\cdot)=0\}$ and $\mrm{Supp}(C)=(\ker C)^\perp$, for any  $C\in\mal C (\Lambda^2 \hor^*\otimes TN)$. 
\begin{prop}
Let us supppose that $\mR_\ver=0$. Then the following identities hold:
$$
\begin{array}{lcl}
\textbf{\em (i) }\  dH_{|\hor^2\times\ver^2}=0 &  &   \textbf{\em (iii) }\   d(F\cdot N_F)_{|\hor^2\times\ver^2} =0 \\
\textbf{\em (ii) }\  dH_{|\hor^2\times\ver^2}^\star=0  &   & \textbf{\em (iv) }\    d \left( F\circact (\mrm{Skew}(\Phi) + \mrm{Skew}(\mR_\ver) )\, \right)_{|\hor^2\times\ver^2}=0.                           
 \end{array}
$$
\end{prop}
\textbf{Proof.}It results from lemma~\ref{d-alpha} applied to the  horizontal 3-forms $F\cdot N_F$ and $\Bar J{\circact}\, \Bar T$, and then from lemma~\ref{lemma d=d}.\\
%
%
%\begin{rmk}\em
%Let us remark that if $\mR_\ver=0$, in particular if we are in presence with a fibre bundle, then (iii) holds.
%\end{rmk}
%%
%
\begin{prop}
Let us supppose that $\mR_\ver=0$ and that the horizontal curvature $\Phi$ is pure. Then the following statements are equivalent
$$
\begin{array}{lcl}
\textbf{\em (i) }\  dH_{|\hor^3\times\ver}=0 &  &  \textbf{\em (iii) }\   d(F\cdot N_F)_{|\hor^3\times\ver} =0 \\
\textbf{\em (ii) }\  dH_{|\hor^3\times\ver}^\star=0  &   & \textbf{\em (iv) } \   \mak S_{X,Y,Z} \, \Bar J N_{\Bar J}(\rho(V)X,Y,Z)=0.
\end{array}
$$
Moreover these later are also equivalent to the following equivalent statements
\begin{description}
\item[(iv)] $N_{\Bar J}(X,Y,\rho(V)Z)=0$.
\item[(v)] $  N_{\Bar J}( \hor,\hor)\ \perp  \rho(\ver)(\hor)$, or equivalently $\ker N_{\Bar J}\perp \ker\Phi $, i.e. $\mrm{Supp}(N_{\Bar J})\perp \mrm{Supp}(\Phi)$.
\end{description} 
We will  say that the 2-forms $N_{\Bar J}$ and $\Phi$ have orthogonal supports, when  (v) is satisfied.
\end{prop}
\textbf{Proof.} According to lemma~\ref{d-alpha}, we have
\begin{eqnarray}
d(\Bar J\circact\Bar T)(V_0,H_1,H_2,H_3) & = & \underset{i,j,k}{\mak S} (\Bar J\circact\Bar T)(\rho(V_0)H_i,H_j,H_k)\nonumber\\
 & = &  \underset{i,j,k}{\mak S} \langle \Bar J\left(-3\Bar T^{--} + \Bar T^{++} - \Bar T^{1,1}\right) (H_i,H_j),\rho(V_0)H_k\rangle \nonumber\\
& = &  \underset{i,j,k}{\mak S} \langle \Bar J\left(-3\Bar T^{--} + \Bar T^{++} - \Bar T^{1,1}\right) (H_i,H_j),\rho^{(-)}(V_0)H_k\rangle\label{dBarJ}
\end{eqnarray}
In the same way, we have 
$$
d(F\cdot N_F)(V_0,H_1,H_2,H_3)= \underset{i,j,k}{\mak S} \langle (4\Bar J \Bar T^{--})(H_i,H_j),\rho^{(-)}(V_0)H_k\rangle.
$$
Moreover, for any $\nabla$-parallel 3-form $\alpha\in\mal C\left( (\Lambda^2\hor^*)\wedge\ver^* \right)$, we have
\begin{eqnarray*}
d\alpha((V_0,H_1,H_2,H_3) & = & \underset{i,j,k}{\mak S} \alpha(T(V_0,H_i),H_j,H_k) + \alpha( T(H_i,H_j), V_0, H_k) \\
 & = & \underset{i,j,k}{\mak S} \alpha(T^\ver(V_0,H_i),H_j,H_k) + \alpha( \Bar T (H_i,H_j), V_0, H_k)\\
& = & \underset{i,j,k}{\mak S} \alpha(\Bar T (H_i,H_j), V_0, H_k)
\end{eqnarray*}
since $R^\ver=0$ and hence $T^\ver(V_0,H_i)=0$. Therefore let us apply this:
\begin{multline}\label{dBarJTrond}
d(F\circact \,\mathring T)(V_0,H_1,H_2,H_3)   =  d (F\circact \, \mrm{Skew}(\Phi))(V_0,H_1,H_2,H_3) 
 =  \underset{i,j,k}{\mak S} (F\circact \, \mrm{Skew}(\Phi))(H_k, \Bar T (H_i,H_j), V_0)\\
 =  \underset{i,j,k}{\mak S} (\Bar J\circact \, \Phi)(H_k, \Bar T (H_i,H_j), V_0)
 = \underset{i,j,k}{\mak S} \langle \left[\rho(V_0),\Bar J \,\right] H_k, \Bar T (H_i,H_j)\rangle\\
 = 2\underset{i,j,k}{\mak S} \langle \rho^{(-)}(V_0)H_k, \Bar J\Bar T (H_i,H_j)\rangle
\end{multline}
Then applying lemma~\ref{lemma d=d} to \eqref{dBarJ} and \eqref{dBarJTrond} , we obtain 
$$
\underset{i,j,k}{\mak S} \langle \Bar J\left( -\Bar T^{--} + \Bar T^{1,1} + 3 \Bar T^{++} \right) (H_i,H_j),\rho^{(-)}(V_0)H_k \rangle=0
$$
and then taking respectively the $(3,0)+(0,3)$-type and the $(2,1)+(1,2)$-type part \wrt to the horizontal variables $(H_1,H_2,H_3)\in \hor^3$, of this equation we obtain
\begin{gather}
\langle \Bar J \Bar T^{++}  (H_i,H_j),\rho^{(-)}(V_0)H_k \rangle=0\\
\langle \Bar J\left( -\Bar T^{--} + \Bar T^{1,1}\right) (H_i,H_j),\rho^{(-)}(V_0)H_k \rangle=0
\end{gather}
This proves the equivalence between the four first statements.\medskip\\
The equivalence (v) $\Leftrightarrow$ (vi) is obvious. It remains to prove the equivalence (iv) $\Leftrightarrow$  (v). But this results immediately from the fact that $N_{\Bar J}=4\Bar T^{--}$ is a 3-form of type $(3,0)+(0,3)$ and $\rho=\rho^{(-)}$.\\
This completes the proof. \hfill$\square$\medskip\\
Finally,
\begin{lemma}
The components in $\Lambda^4\ver^*$ and $(\Lambda^3\ver^*)\wedge \hor^*$ of the following exterior derivatives: $dH$, $dH^\star$, $d (F\cdot N_F)$ and $d  \left( F{\circact} (\mrm{Skew}(\Phi) + \mrm{Skew}(\mR_\ver) )\,\right) $, 
 vanishes.
\end{lemma}
\proof
Apply lemma~\ref{d-alpha} to the horizontal 3-forms $\Bar J \circact \Bar T$ and $F\cdot N_F$. Then use lemma~\ref{lemma d=d}. \comprf \hsq \medskip\\
Let us summarize:
\begin{thm}\label{thm-stringy-closed}
Let $(N,F,h)$ be a reductive metric $f$-manifold of global type $\mal G_1$. Let us suppose that one of its  characteristic connections, $\nabla$,  has a parallel torsion $\nabla T=0$. Let us supppose that $\mR_\ver=0$ and that the horizontal curvature $\Phi$ is pure.\medskip\\
$\bullet$ Then $(N,F,h)$ has a closed stringy structure \iif the horizontal 3-form $F\cdot N_F$ is closed.\medskip\\
$\bullet$ Moreover, this happens  \iif the horizontal complex structure $\Bar J$ is a cyclic permutation of the horizontal curvature,  and  the 2-forms $N_{\Bar J}$ and $\Phi$ have orthogonal supports.
\end{thm}

\paragraph{The particular case of Horizontally Hermitian $f$-manifolds.} 
\index{horizontally K\"{a}hler}
In horizontally Hermitian $f$-manifolds, we have $F\cdot N_F=0$ and therefore $F\cdot N_F=0$  is a closed 3-form. Therefore we can apply theorem~\ref{thm-stringy-closed}.
\begin{cory}
Let $\pi\colon (N,F,h)\to (M,g)$ be a  Riemannian submersion from a  metric $f$-manifold into a Riemannian manifold.  Let us suppose that $(N,F,h)$ is reductive and horizontally Hermitian, so that it is in particular  horizontally of type $\mal G_1$. Let us suppose also that one of its  characteristic connections, $\nabla$,  has a parallel torsion $\nabla T=0$ and that the horizontal curvature $\Phi$ is pure. Then $(N,F,h)$ has a closed stringy structure. 
\end{cory}
\subsubsection{The sigma model with a Wess-Zumino term in reductive metric $f$-manifold of global type $\mal G_1$.} \index{Wess-Zumino|(}
Now, we can conclude with the following variational interpretation of the stringy harmonicity.
\begin{thm}\label{thm-WZW-G1-red}
Let $(N,F,h)$ be a reductive  metric $f$-manifold of global type $\mal G_1$. Let us suppose that  $(N,F,h)$ is stringy closed. Let $\nabla$ be one characteristic connection.\\
$\bullet$ Then the equation for stringy harmonic maps (\wrt $\nabla$)  $f\colon L\to N$ is exactly the Euler-Lagrange equation for the sigma model in $N$ with a Wess-Zumino term defined by the closed 3-form
$$
H=-d\Omega_F + F\cdot N_F  - \dfrac{1}{2}F\circact \left(\mrm{Skew}(\Phi)\right). 
$$ 
$\bullet$ Moreover the equation for $\star$-stringy harmonic maps $f\colon L\to N$ is exactly the Euler-Lagrange equation for the sigma model in $N$ with a Wess-Zumino term defined by the closed 3-form
$$
H^\star=-d\Omega_F + \dfrac{1}{2}F\cdot N_F - \dfrac{1}{2}F\circact \left(\mrm{Skew}(\Phi)\right). 
$$ 
\end{thm}
\index{connection!characteristic|)}
%
%%%%%%%%%%%%%%%%%%%%%%%%%%%%%%%%%%%

\subsubsection{The example of a naturally reductive homogeneous space}\label{example-homogenenous}
\index{canonical!connection, $G$-invariant|(}
By definition of ($\star$-)stringy harmonicity and the expression of the torsion of $\nabla^0$ in terms of the Lie bracket, we have the following.
\begin{prop}\label{prop-stringy-hom-even}
Let $N=G/K$ be a Riemannian homogeneous manifold endowed with a $G$-invariant complex structure $F$. Let $\g=\kk\oplus\nk$ be a reductive decomposition of $\g$, and $\nk =\mk\oplus \pk$ the $\Ad K$-invariant decomposition defined by $F$. Let  $f\colon L\to N$ be a smooth map, $U\colon L\to G$ be a (local) lift of $f$ and $\alpha=U^{-1}.dU$ the corresponding Maurer-Cartan form. Then in terms of $\alpha$, the equation of stringy harmonicity (\wrt $\nabla^0$) is written
$$
\left\lbrace \begin{array}{l}
\displaystyle d*\alpha_\pk + [\alpha_\kk\wedge *\alpha_\pk] + \dfrac{1}{2}[\Bar J_0 \alpha_\mk\wedge\alpha_\mk]_{\pk}=0\\
d*\alpha_\mk + [\alpha_\kk\wedge *\alpha_\mk] - \dfrac{1}{2}\Bar J_0\left[\Bar J_0 \alpha _\mk\wedge \Bar J_0\alpha_\mk \right]_\mk + \dfrac{1}{2}\left(  [\alpha_\pk\wedge \Bar J_0 \alpha_\mk ] - \Bar J_0 [\alpha_\pk\wedge \alpha_\mk ]\right) =0
\end{array}\right.
$$
whereas the equation of $\star$-stringy harmonicity is written:
$$
\left\lbrace \begin{array}{l}
\displaystyle d*\alpha_\pk + [\alpha_\kk\wedge *\alpha_\pk] + \dfrac{1}{2}[\Bar J_0 \alpha_\mk\wedge\alpha_\mk]_{\pk}=0\\
\displaystyle d*\alpha_\mk + [\alpha_\kk\wedge *\alpha_\mk] + \dfrac{1}{2}[\Bar J_0 \alpha _\mk\wedge\alpha_\mk]_\mk  + \dfrac{1}{4}J_0\left( \left[\Bar J_0\alpha_\mk \wedge \,\Bar J_0 \alpha _\mk \right]_\mk + \left[ \alpha _\mk\wedge\alpha_\mk \right]_\mk\right)\\
 \hfill + \dfrac{1}{2}\left(  [\alpha_\pk\wedge \Bar J_0 \alpha_\mk ] - \Bar J_0 [\alpha_\pk\wedge \alpha_\mk ]\right) =0. 
\end{array}\right.
$$
where $\Bar J_0$ is the complex structure on $\mk$ corresponding to  $F$.
\end{prop}
\subsubsection{Geometric interpretation of the maximal determined even case.}
\index{canonical!$f$-structure|(}\index{even case|(} \index{k symmetric space@$k'$-symmetric space|(}
In this subsection, we suppose that $N=G/K$ is a (locally) $2k$-symmetric space, and we use the notations  and the conventions of  \ref{finitorderauto}.
\begin{thm}\label{thm-max-deter-even-stringy}
Let us suppose that $N=G/K$ is a (locally) $2k$-symmetric space endowed with its canonical $f$-structure\footnote{defined by \eqref{def-even}} $F$ and its canonical connection $\nabla^0$.  Then the associated maximal determined system, $\syst(2k-1,\tau)$ is the equation of $\star$-stringy harmonicity for the geometric map $f\colon L\to N$: ${(\nabla^0)}^* df + (F\star T)(f)=0$.\footnote{Where we have removed the index "$g$" which specifies that the previous terms are computed with respect to some Hermitian metric $g$ on $L$.}\\
Moreover, if we consider now  that $N=G/K$ is endowed with the $f$-structure $F^\star:=\oplus_{j=1}^{k-1}(-1)^{j}{F}_{[\mk_j]}\oplus 0_{[\g_k]}$, then this system is the equation of stringy harmonicity for the geometric map $f\colon L\to N$: ${(\nabla^0)}^* df + (F^\star \cdot T)(f)=0$.
\end{thm}
\proof
The first point follows from theorem~\ref{thm-max-deter-even}, proposition~\ref{prop-stringy-hom-even} and the fact that the component $\Bar T= T_{|\hor^2}^\hor$ of the torsion $T$ of  $\nabla^0$ satisfies the same kind of identities as in proposition~\ref{T-eps-esp'-nabla0}. The second point follows from proposition~\ref{prop-JstarJ} and remark~\ref{rmk-hold-vectorbundle} applied to the complex vector bundle $(\hor,\Bar J)$, the 2-form $\Bar T\in \mal T(\hor)$ and the $\Bar J$-invariant decomposition $TN=E^+\oplus E^-$,  defined by $E^+=\oplus_{\underset{j\text{ even}}{j=1}}^{k -1}[\mk_j]$, and $E^-=\oplus_{\underset{j\text{ odd}}{j=1}}^{k -1}[\mk_j]$. \comprf\hsq 
\begin{thm}\label{thm-geom-interpret}
Let us suppose that $N=G/K$ is a (locally) $2k$-symmetric space endowed with its canonical $f$-structure $F$ and its canonical connection $\nabla^0$. Let us suppose that $N=G/K$ is naturally reductive and we choose a naturally reductive $G$-invariant metric $h$ for which $\tau_{|\nk}$ is orthogonal\footnote{See proposition~\ref{prop-2ksym-G1&pure}.} and thus which is compatible with $F$. Then $(N,F,h)$ is a reductive metric $f$-manifold  of global type $\mal G_1$, and its horizontal curvature $\Phi$ is pure. Moreover, 
$(N,F,h)$ has a closed stringy structure.\\
Therefore, the  associated maximal determined system, $\syst(2k-1,\tau)$, is  exactly the Euler-Lagrange equation for the sigma model in $N$ with a Wess-Zumino term defined by the closed 3-form
$$
H^\star=-d\Omega_{F} + \dfrac{1}{2}F\cdot N_{F} - \dfrac{1}{2}F\circact \left(\mrm{Skew}(\Phi)\right) .
$$ 
Moreover, if we consider now  that $N=G/K$ is endowed with the $f$-structure $F^\star$, then the previous system is  exactly the Euler-Lagrange equation for the sigma model in $N$ with a Wess-Zumino term defined by the closed 3-form
$$
H=-d\Omega_{F^\star} + F^\star \cdot N_{F^\star} - \dfrac{1}{2}F\circact \left(\mrm{Skew}(\Phi)\right) .
$$  
\end{thm}
\index{Wess-Zumino|)}
The theorem will follow from the two following lemmas.
\begin{lemma}\label{proof-cyclic}
In the metric $f$-manifold $(N,F,h)$ of theorem~\ref{thm-geom-interpret},  $\Bar J$ is a cyclic derivation of the horizontal curvature.
\end{lemma}
\textbf{Proof.} According to subsection~\ref{example-homogenenous-G1-red}, we have 
\begin{eqnarray}
\widetilde{\Bar{\mrm A}(X,Y)Z} & = & \left[\, [\Bar J_0 X_\mk,Y_\mk]_{\pk}, Z_\mk \right]  +  \left[\, [ X_\mk, \Bar J_0 Y_\mk]_{\pk}, Z_\mk \right] - \Bar J_0 \left[\, [ X_\mk,Y_\mk]_{\pk}, Z_\mk \right] +  \left[\, [ X_\mk,Y_\mk]_{\pk},\Bar J_0 Z_\mk \right]\nonumber\\
 & = & 2 \left( \left[\, [\Bar J_0 X_\mk,Y_\mk]_{\pk}, Z_\mk \right] +  \left[\, [ X_\mk,Y_\mk]_{\pk},\Bar J_0 Z_\mk \right] \right)\label{equation-tilde-A} 
\end{eqnarray}
In  this equation, $X,Y,Z\in \mal C(\hor)$ with lifts $X_\mk,Y_\mk,Z_\mk\in \mal C^\infty(G,\mk)$, and $\widetilde{\Bar{\mrm A}(X,Y)Z}\in \mal C^\infty(G,\mk)$ denotes the lift of $\Bar{\mrm A}(X,Y)Z\in \mal C(\hor)$. An other possibility is to consider that $X,Y,Z\in \hor$ are horizontal vectors at some point $y\in N$ and that we have chosen $g\in G$ such that $g.G_0=y$ and  that we have set $X=\Ad g (X_\mk)$, and idem for $Y$ and $Z$. Then the equation~\eqref{equation-tilde-A}, when written  in the form $\widetilde A = B_\mk$, means in fact that we have $A=\Ad g (B_\mk)$.\\
 In the following, we will use this kind of notation without recalling these precisions.\\
Moreover, since $\Phi$ is pure (proposition~\ref{prop-2ksym-G1&pure}), we have $[ \Bar J_0 X_\mk, Y_\mk]_{\pk} = [  X_\mk, \Bar J_0 Y_\mk]_{\pk}$ and $\left(\ad_{\mk}\g_k\right)\Bar J_0 = - \Bar J_0 \left(\ad_{\mk}\g_k\right)$, then we obtain
\begin{eqnarray*} 
\underset{X,Y,Z}{\mak S}  \widetilde{\Bar{\mrm A}(X,Y)Z}  & = &  2i\underset{X,Y,Z}{\mak S} \left( \sum_{1\leq |j|,|l|\leq k-1}  \negthickspace   s(-j) \left[ [X_j, Y_{-j + k}], Z_l \right] \ +  \sum_{1\leq |j|,|l|\leq k-1} \negthickspace   s(-l) \left[ [X_j, Y_{-j + k}], Z_l \right]\right) \\
 & = &  2i\underset{X,Y,Z}{\mak S} \sum_{1\leq |l|\leq k-1 } s(-l) \sum_{ \underset{jl\,>0}{ 1\leq |j|\leq k-1} }   \left[ [X_j, Y_{-j + k}], Z_l \right]
\end{eqnarray*}
where $s(n)$ is the sign of $n\in \Z$.
Furthermore, $\alpha:=\underset{X,Y,Z}{\mak S}  \widetilde{\Bar{\mrm A}(X,Y)Z}$ defines a element of $\Lambda^3\g^*\otimes\g$, which in fact corresponds to a  4-form of type $(4,0) + (0,4)$ according to the proof of proposition~\ref{prop-Hor4} (see in particular, equations \eqref{A-def-4form} and \eqref{4,0 + 0,4}).
Let us consider the component in $\left(\g_j^*\wedge \g_{-j+k}^* \wedge \g_l^*\right)\otimes\g$ of this the element of  this element of $\Lambda^3\g^*\otimes\g$. We remark first  that this component is in fact in $\left(\g_j^*\wedge \g_{-j+k}^* \wedge \g_l^*\right)\otimes\g_{k+l}$. Let us set 
\begin{eqnarray*}
I^+ & = & \left \{ \, \{ j,-j + k, l \},\ 1\leq j,l\leq k-1 \right\} \\
I^- & = & \left \{ \, \{ j,-j + k, l \},\  -(k-1)\leq j,l\leq -1  \right\}=- I^+\\
I & = & I^+ \sqcup I^- .
\end{eqnarray*}
Then, denoting by $\ovr{k+l}$ the representant modulo $2k$ of $k+l$ in $\{-(k-1),\ldots, -1\}\cup\{ 1,\ldots, k-1\}$, we have  $\{j,\ovr{k+l}, l\}\notin I$ because  $\ovr{k+l}$ has a opposite sign to $j$ and $l$. Therefore, the component in  $\left(\g_j^*\wedge \g_{k+l}^* \wedge \g_l^*\right)\otimes\g_{-j + k}$ of our element $\alpha$ vanishes, which since $\alpha$ is a 4-form implies that $\alpha=0$. This completes the proof.\hfill$\square$
\begin{lemma}\label{proof-ortho-suppor}
In the metric $f$-manifold $(N,F,h)$ of theorem~\ref{thm-geom-interpret},  the 2-forms $N_{\Bar J}$ and $\Phi$ have orthogonal supports.
\end{lemma}
\textbf{Proof.} Denoting by $T$ the torsion of the canonical connection $\nabla^0$, we have
$$
 \widetilde{T^{--}(X,Y)}=-\sum_{ \underset {1\leq i\leq j\leq k-1} {k+1\leq i+j \leq 2(k-1)}   } 
\left(1- \dfrac{\delta_{ij}}{2}\right)  \left(     [X_{\mk_i}, Y_{\mk_j}]_{\mk_{i+j}} +  [X_{\mk_j}, Y_{\mk_i}]_{\mk_{i+j}}\right). 
$$
Therefore
$$
\langle T^{--}(X,Y),\rho(V)Z\rangle = -\sum_{ \underset {1\leq i\leq j\leq k-1} {k+1\leq i+j \leq 2(k-1)}   }  \left(1- \dfrac{\delta_{ij}}{2}\right) \langle [X_{\mk_i}, Y_{\mk_j}]_{\mk_{i+j}}, [Z_{\mk_{i+j-k}},V_k]\rangle
$$
but we have
\begin{eqnarray*}
\langle [X_{\mk_i}, Y_{\mk_j}]_{\mk_{i+j}}, [Z_{\mk_{i+j-k}},V_k]\rangle & =  &  
\langle [X_{\mk_i}, Y_{\mk_j}]_{i+j}, [Z_{i+j-k},V_k]\rangle +  \langle [X_{\mk_i}, Y_{\mk_j}]_{-(i+j)}, [Z_{-(i+j)-k},V_k]\rangle \\
 & = & \langle \left[  [X_{\mk_i}, Y_{\mk_j}]_{i+j}, Z_{i+j-k} \right] ,V_k \rangle +  \langle\left[  [X_{\mk_i}, Y_{\mk_j}]_{-(i+j)}, Z_{-(i+j)-k} \right]  , V_k\rangle 
\end{eqnarray*}
and these two scalar products vanishes, because $2(i+j) - k\neq k \mod 2k$ and  $-2(i+j) + k\neq k \mod 2k$. Indeed 
$2(i+j)\neq 0 \mod 2k$, since $k+1 \leq i+j \leq 2(k-1)$. This completes the proof.\hfill$\square$\medskip\\
\textbf{Proof of theorem~\ref{thm-geom-interpret}}
It follows from proposition~\ref{prop-2ksym-G1&pure}, lemmas~\ref{proof-cyclic} and \ref{proof-ortho-suppor}, theorems~\ref{thm-stringy-closed} and \ref{thm-WZW-G1-red}. \comprf \hsq 
\begin{rmk}\em
Let us remark that a 4-symmetric space, endowed with its canonical $f$-structure $F$ and  a  naturally reductive metric is horizontally K\"{a}hler.
\end{rmk}
%
%%%%%%%%%%%%%%%%%%%%%%%%%%%%%

\subsubsection{Twistorial Geometric interpretation of the maximal determined even case.}\index{twistor|(}
We come back to the situation of \ref{subsub-red-mijo}. We use also the notation of definition~\ref{def-twistor-lift}. Then according to proposition~\ref{red-mijo-canonical} \footnote{Or we can use proposition~\ref{prop-same-stringy-harm-if-red} as well} and theorem~\ref{thm-max-deter-even-stringy}, we obtain
\index{connection!paracharacteristic}
\begin{thm}
Let $(L,j)$ be a Riemann surface, $f\colon L\to N=G/G_0$ be a map and $J=f^*\mijo$ the corresponding map into $\mZ_{2k,2}^{\alpha_0}(M,J_2)$. Then $f$ is a geometric solution of the even maximal determined system $(\syst(2k-1,\tau))$ \iif $J\colon L\to \mZ_{2k,2}^{\alpha_0}(M,J_2)$  is $\star$-stringy harmonic \wrt the canonical paracharacteristic connection in the twistor space $\mZ_{2k,2}^{\alpha_0}(M,J_2)$.\smallskip\\
Moreover, let us consider now  that $N=G/G_0$ is endowed with the $f$-structure $F^\star:=\oplus_{j=1}^{k-1}(-1)^{j}{F}_{[\mk_j]}\oplus 0_{[\g_k]}$, and that $\mZ_{2k,2}^{\alpha_0}(M,J_2)$ is endowed with the $f$-structure $\mal F^\star:=\oplus_{j=1}^{k-1}(-1)^{j}{\und \mJ}_{\mk_j(\mJ)}\oplus 0_{\ver^{\mZ,2}}$, where $\mk_j(\mJ)$ is defined as in 
definition~\ref{def-twistor-lift}, and $\ver^{\mZ,2}$ is the vertical space of $\mZ_{2k,2}^{\alpha_0}(M,J_2)$.\\
Then $f$ is a geometric solution of this system \iif $J\colon L\to \mZ_{2k,2}^{\alpha_0}(M,J_2)$ is stringy harmonic \wrt the canonical paracharacteristic connection in the twistor space $\mZ_{2k,2}^{\alpha_0}(M,J_2)$.
\end{thm}\index{stringy harmonic|)}
\index{canonical!$f$-structure|)}\index{even case|)}\index{k symmetric space@$k'$-symmetric space|)}
\subsubsection{About the variational interpretation in the Twistor spaces.}\label{variational-twistor}
We have seen that the twistor spaces $(\zdk^\alpha(M),\mal F,h)$, and in particular $(\Sigma(M),\mal F,h)$, endowed with their structure of homogeneous fibre $f$-bundle defined in \ref{example-zdk}, are reductive and  horizontally of type $\mal G_1$. Are they in general of global type $\mal G_1$? More generally, which are their subbundles which are of global type $\mal G_1$? We give below the answer to these questions which are studied in \cite{ki-twistor} where the following results are proved. \\
Let us remark that since the twistor spaces and their subbundles are horizontally of type $\mal G_1$, the question is in fact to know if these spaces satisfy the condition that $\tilde N_{F|\mal S(\hor\times\hor\times\ver)}$ is skew-symmetric, or equivalently according to remark~\ref{rmk-H2v-skew-Nf}, the condition
$$
[\rho(V),\Bar J]=-\Bar J N_F(V),  \quad \forall V\in \ver.
$$
This condition implies strong conditions on the curvature of the metric connection $\ovr\nabla$ on $M$, as we can see from the expression of $\Phi$ in terms of the curvature of $\ovr \nabla$ (theorem~\ref{entermdeJ2}-(ii)). Now, let us present the  results obtained after some investigations (see \cite{ki-twistor} for more details).

\begin{lemma}\cite{ki-twistor}
Let $\nabla$ be a metric connection on the Riemannian manifold $(M,g)$. Suppose that the sectional curvature $k(P)$ of $\nabla$ depends only on the point $x\in M$ and not on the plan $P\in \Lambda^2 T_x M$, and that $k$ does not vanish. Then we have
$$
\left\{\begin{array}{l}
R(X,Y)Z= k \left(\langle Z,Y\rangle X - \langle Z,X\rangle Y\right)\\
T(X,Y)=\dfrac{1}{2k}\left( (Y\cdot k) X - (X\cdot k) Y\right) .
\end{array}\right.
$$
In particular, $\nabla$ is geodesically equivalent to Levi-Civita \iif it coincides with Levi-Civita. Moreover,  $k$ is constant \iif $\nabla$ coincides with Levi-Civita. 
\end{lemma}
\begin{thm}\cite{ki-twistor}
Let $(M,g)$ be a Riemannian manifold endowed with a metric connection $\nabla$.
Let us consider the homogeneous fibre $f$-bundle $(\Sigma(M),h,\mal F)$ defined in \ref{example-zdk}. Then  $(\Sigma(M),h,\mal F)$ is globally of type $\mal G_1$ \iif $\nabla=\nabla^g$ and $(M,g)$ has a constant sectional curvature. Therefore $M$ is a locally symmetric space and $\Sigma(M)$ is locally 4-symmetric.
\end{thm}
\begin{thm}\cite{ki-twistor}
Let $(M,g)$ be a Riemannian manifold endowed with a metric connection $\nabla$.
Let us consider the homogeneous fibre $f$-bundle $(\zdk^\alpha(M),\mal F,h)$ defined in \ref{example-zdk}. Then it is globally of type $\mal G_1$ \iif $\nabla=\nabla^g$ and $(M,g)$ has a constant sectional curvature. Therefore $M$ is a locally symmetric space and $\Sigma(M)$ is locally $2k$-symmetric.
\end{thm}
\begin{thm}\cite{ki-twistor}
Let $(M,g)$ be a Riemannian manifold endowed with a metric connection $\nabla$ with parallel torsion $\nabla T=0$. Let $(N,F,h)$ be a subbundle of $(\zdk^\alpha(M),\mal F,h)$. Then $N$ is globally of type $\mal G_1$ \iif $M$ is a locally homogeneous space, i.e. its universal covering can be written $\tl M= G/H$, and $\nabla$ coincides with  $\onabla{M}{0}{}$ the canonical connection of $M$. Moreover $N$ is a locally homogeneous space, i.e. $\tl N=G/K$, and the characteristic connection coincides with the canonical connection of $N$. 
\end{thm}
\index{twistor|)}\index{canonical!connection, $G$-invariant|)}
\index{f structure@$f$-structure|)}\index{global type $\mal G_1$|)}\index{f manifold@$f$-manifold|)}
\index{horizontal subbundle|)}
\index{vertical subbundle|)}

\subsection{Bibliographical remarks and summary of the results.}
All together, in this section, we have investigated an unexplored domain (at least for us and from the specific geometric point of view under consideration), by passing from the study of almost Hermitian structure to the one of metric $f$-manifolds . Indeed, firstly we had to find the good notions to consider: e.g. which action of $F$ on $T$ one should choose (to define for example the new notion of stringy harmonicity). Secondly, we had to find the good geometric datas to express all the geometric objects that we use. Indeed, in the almost Hermitian case: $\Omega_J$ and $N_J$ was enough, but in the case of metric $f$-manifolds: in addition to $\Omega_F$ and $N_F$, we need also the curvatures $\Phi:=\mrm R_\hor$ and $\mrm R_\ver$, but sometimes also the Levi-Civita derivative of the projection $q$. Concretely, the previous familiar study done in the case of almost Hermitian manifolds corresponds to what happens here for the horizontal part $\hor$ of the tangent space. This horizontal part uses only horizontal components $\hor^3$ (of the several geometric objects in presence which are most of the time $TN$-valued 2-forms on $N$). To this horizontal component $\hor^3$, we have to add the vertical one $\ver^3$ which in fact plays a negligible role, but also and mostly the coupling terms $\mal S(\hor\times\hor\times\ver)$ and  $\mal S(\hor\times\ver\times\ver)$. Frequently, we have to innovate in order to invent the pertinent new notions which will allow the global study of $(N,F,h)$ taking account of all these coupling components.\smallskip\\
Since 1960's, of great concern are $f$-structures introduced by K. Yano,  generalizing almost complex and almost contact structures. In turn, metric $f$-structures  are important objects of study in \emph{generalized Hermitian geometry},  an area of modern differential geometry developped since the middle of 1980's in particular by the Russian school (see Kirichenko \cite{Kirichenko-a,Kirichenko-b,Kirichenko-c}). In particular the notions of K\"{a}hler $f$-structure, Hermitian $f$-structure, $\mal G_1 f$-structure, nearly K\"{a}hler $f$-structure etc.. However no one of these notions (generalizing the corresponding notions on almost Hermitian manifolds) corresponds to our new notions of \emph{reductive}, \emph{of global type} $\mal G_1$, \emph{horizontaly K\"{a}hler} manifolds and so on, defined  in the present section. It seems therefore that there are different ways to generalize familiar geometric notions in almost Hermitian geometry to notions in metric $f$-geometry,  and that our own way is totally new. Remark that each of these  generalizations contains the particular case of $2k$-symmetric spaces (see \cite{Balashchenko1, Balashchenko2,Balashchenko3}).\smallskip\\
We would like to  mention some related works. Harmonic maps into metric $f$-manifolds have been studied by many mathematicians (in several contexts and for several kind of harmonic maps): Rawnsley \cite{rawnsley}, Black \cite{Black}, Urakawa \cite{Urakawa}, Bryant, Gherghe, Ianus and Pastore \cite{Bryant, Gherghe, Ianus}.\medskip\\
Now, let us summarize the main novelties of the present section.\smallskip\\
In  subsection~\ref{f-connection&torsion}, we provided a general studies of metric $f$-connections on metric $f$-manifolds.  We first characterized metric connections preserving the splitting $TN=\hor\oplus\ver$, and  then we characterize the Riemannian manifolds endowed with an orthogonal splitting, which admit such a connection with skew-symmetric torsion. We called them \emph{reductive} $f$-manifolds. Secondly, we characterized metric $f$-connnections in terms of their torsions (theorem~\ref{thm-nablaOmega_F=0}) which generalizes to metric $f$-structures the theorem of Gauduchon \cite[Proposition~2]{Gauduchon} (see theorem~\ref{Gauduchon}) about almost Hermitian structures. Then we characterize metric $f$-manifolds which admit metric $f$-connections with skew-symmetric torsion (in several steps: one step for each component: $\hor^3$, $\mal S (\hor\times\hor\times\ver)$ and $\mal S(\hor\times\ver\times\ver)$) which are  our so-called \emph{reductive manifold of global type} $\mal G_1$. This is the content of theorem~\ref{characteristic-equiv-G1} which is in our opinium an important theorem (and of course completely new).  In the particular case where $\ver$ is a line bundle, $(N,F,h)$ is an almost contact metric manifold and  this theorem contains the result of \cite[theorem~{8.2}]{Friedrich-ivanov}. \medskip\\
In subsection~\ref{f-connect-bundl}, we specialize our study to the case of fibre bundles. We give several geometric conditions which characterize on fibre bundles the geometric notions defined in subsection~\ref{f-connection&torsion} (reductivity, global type $\mal G_1$, horizontal global type $\mal G_1$, Horizontally K\"ahler, etc..). Then we apply these results to the examples of the twistor spaces $\mal Z_{2k}^\alpha(M)$ and $\mal Z_{2k,j}^\alpha(M,J_j)$ (and more easily to the example of $2k$-symmetric spaces). In particular, these twistor spaces provide classes of examples of manifolds which satisfy  the previous geometric properties. Finally, we prove that the previous geometric structures are preserved under the canonical embedding $\mijo\colon G/G_0\hookrightarrow \mZ_{2k,2}^{\alpha_0}(G/H,J_2)$.\medskip\\
In subsection~\ref{Stringy-Harm-f-manif}, we define the new notion of stringy harmonic maps \wrt an $f$-structure. Then we give some sufficient conditions so that this stringy harmonicity admits a variational interpretation in terms of a sigma model with a Wess-Zumino term. There are several ways to interpret this result, which leads to several novelties. Indeed, from the point of view of geometric analysis, this result provides a new class of geometric variational problems taking place in some class of metric $f$-manifolds. From the point of view of mathematical physic this provides a new class of geometric structures in which it is possible to define a non linear sigma models with a Wess-Zumino term.\smallskip\\
We conclude this subsection by the interpretation of the maximal determined even elliptic integrable systems in terms of stringy harmonic maps. We give this interpretation in two different settings: in the homogeneous $2k$-symmetric space and in the twistor space by using the canonical embedding.\smallskip\\
A consequence of this interpretation is a new contribution to the field of (integrable) non linear sigma models. Indeed we give new examples of integrable two-dimensional non linear sigma models. These new examples take place in some homogeneous spaces, namely $2k$-symmetric spaces, which are not symmetric spaces. Again, at our knowledge, all the already known integrable two-dimensional non linear sigma models take place in symmetric spaces or (equivalently) in Lie groups.\\
Several other consequences using this interpretation in the setting of the twistor spaces will be given in forthcoming papers \cite{ki-twistor,ki7,ki8}.

%%%%%%%%%%%%%%%%%%%%%%%%%%%%%%%%%%%%%%%%%%%%%%%%%%%%%%%%%%%%%%%%%%%%%%%%%%%%%%%%%%%%%%%%%%%%%%%%%
%                                                                                               %
%                                                                                               %
%                       Generalized harmonic maps  into homogeneous space                                                  %
%                                                                                               %
%                                                                                               %
%%%%%%%%%%%%%%%%%%%%%%%%%%%%%%%%%%%%%%%%%%%%%%%%%%%%%%%%%%%%%%%%%%%%%%%%%%%%%%%%%%%%%%%%%%%%%%%%%

\section{Generalized harmonic maps into reductive homogeneous spaces}\label{affineharmhom}
\index{canonical!connection, $G$-invariant|(}
%%%
%
This section has 3 objectives.\medskip\\
\textbf{1.} We want  to specify clearly and concretely the geometric meaning of each component of the zero curvature equation on a 1-form $\alpha_\lm$ taking values in a finite dimensional vector subspace of a loop Lie algebra. Indeed, we would like to interpret each component of this equation with a geometric property, as we have already done, for example for the even miminal determined system, in section~\ref{homspacefibr}. But here, we would like to keep  our study in the Lie algebra setting and in the corresponding homogeneous space for the geometric interpretation. Moreover, we would like to see how one can obtain naively some of the results obtained in sections~\ref{vertically-harmonic}, \ref{affineharmonic}, and \ref{gene-Harm-f-structure}. In other words, we would like  to continue and complete, in the light of the studies done in sections~\ref{vertically-harmonic}, \ref{affineharmonic}, and \ref{gene-Harm-f-structure}, the preliminary and naive study begun in section~\ref{melliptic}, keeping the same naive spirit. Then, over and above recovering by a different way,  precising, or completing results of the previous sections,  as well as outlining results already implicitly contained in the previous more general studies but which can be specified in more particular contexts, we also prove some new results. For example, we prove that the strongly harmonicity in a reductive homogeneous space has a formulation in terms of a zero curvature equation (that we call elliptic integrable system associated to a reductive homogeneous space). Moreover, we complete the study of the odd minimal determined system begun in section~\ref{detodcase}.\smallskip\\
\textbf{2.} Furthermore, we study the notions of vertical harmonicity, strongly vertical harmonicity, torsion freedom in the general context of homogeneous spaces fibrations $G/K\to G/H$. We obtain new results concerning how these notions can be related together and to harmonicity. We also give some conditions under which these notions have a formulation in terms of a zero curvature equation.\\ 
In the context of homogeneous spaces endowed with an invariant Pfaffian system, we study the notions of strongly vertically harmonic maps, vertically holomorphically harmonic maps and such. We also prove some relations between these notions. Again, we give some conditions under which these notions have a formulation in terms of a zero curvature equation.
\smallskip\\
\textbf{3.} Moreover, we want to specify the geometric interpretation of the intermediate determined systems ($m_{k'}< m < k'-1$). We already know that these are the equations of stringy harmonic maps which satisfies some additonnal holomorphicity conditions. We would like to make this more precise. \smallskip\\
\textbf{3.} Finally, we conclude with some remarks about the twistorial interpretation.

%%%%%%%%%%%%%%%%%%%%%%%%%%%%%%%%%%%%%%%%%%%%%%%%%%%%%%%%%%%%%%%%%%%%%%%%%%%%%%%%%%%%%%%%%%%%%%%%%
%
\subsection{Affine harmonic maps into reductive homogeneous spaces.}\label{Aff-harm-red-hom}
\index{harmonic map|(}\index{strongly harmonic|(}
Let $N=G/K$ be a reductive homogeneous space and $\g=\kk\oplus \mk$ a reductive decomposition of the Lie algebra $\g$. We use the notations of section~\ref{1} (applied to $N=G/K$ instead of $M=G/H$).
\begin{thm}\label{thmharmonic}
Let $(L,j)$ be a Riemann surface and $f\colon (L,j)\to N$ be a smooth map, let $F\colon L\to G$ be a (local) lift of $u$ and $\alpha=F^{-1}.dF$. Then the following statements are equivalent:
\begin{description}
\item [(i)] $f$ is $\nabla^t$-harmonic for one $t\in [0,1]$.
\item[(ii)] $f$ is $\nabla^t$-harmonic for every $t\in [0,1]$.
\item[(iii)] $d*\alpha_\mk +[\alpha_\kk\wedge *\alpha_\mk]=0$.
\item[(iv)] $\mrm{Im}\left( \bar\partial\alpha_\mk' + [\alpha_\kk''\wedge\alpha_\mk'] + t[\alpha_\mk''\wedge\alpha_\mk']_\mk\right) =0$, $\forall t\in [0,1]$.
\end{description}
In fact, the tension field $\tau^t(f)$ of $f$ with respect to $\nabla^t$ is independent of $t\in [0,1]$.
\end{thm}
\begin{thm}\label{thmstrongharmonic}
In the same situation as above, the following statements are equivalent:
\begin{description}
\item[(i)] $f$ is strongly $\nabla^t$-harmonic for one $t\in [0,1]\setminus\lbrace \frac{1}{2}\rbrace $.
\item[(ii)] $f$ is strongly $\nabla^t$-harmonic for every $t\in [0,1]\setminus\lbrace \frac{1}{2}\rbrace $.
 \item[(iii)] $\bar\partial\alpha_\mk' + [\alpha_\kk''\wedge\alpha_\mk''] + t[\alpha_\mk''\wedge\alpha_\mk']_\mk=0$, $\forall t\in [0,1]\setminus\left\lbrace \frac{1}{2}\right\rbrace$.
\item[(iv)] $f$ is $\nabla^t$-harmonic for one $t\in [0,1]$ and  $[\alpha_\mk\wedge\alpha_\mk]_\mk=0$.
\item[(iv)] $d\alpha_\lm + \dfrac{1}{2}[\alpha_\lm\wedge\alpha_\lm]=0$, $\forall \lm\in S^1$, with $\alpha_\lm=\lm^{-1}\alpha_\mk' +\alpha_\kk + \lm\alpha_\mk''$.
\end{description}
Furthermore $f$ is strongly $\nabla^{\frac{1}{2}}$-harmonic \iif it is $\nabla^{\frac{1}{2}}$-harmonic: indeed $\nabla^{\frac{1}{2}}$ is  torsion free.
\end{thm}
\textbf{Proof of theorem~\ref{thmharmonic}} The tension field $\tau^t(f)$ of $f$ with respect to $\nabla^t$ is given by 
\begin{eqnarray*}
\tau^t(f)=*d^{\nabla^t}*df & = & *\Ad F\left(  d*\alpha_\mk + [\alpha_\kk\wedge *\alpha_\mk] + t[\alpha_\mk\wedge*\alpha_\mk]_\mk\right)  \\
 & = & *\Ad F\left(  d*\alpha_\mk + [\alpha_\kk\wedge *\alpha_\mk]\right) 
\end{eqnarray*}
(see section~\ref{assocovarder} (especially equation (\ref{nabla-0})) and  section~\ref{family}). This proves the equivalence between (i), (ii) and (iii). Then we conclude by remarking that $2\,\mrm{Im}\left(\bar \partial\alpha_\mk' + [\alpha_\kk''\wedge\alpha_\mk']\right) = d*\alpha_\mk + [\alpha_\kk\wedge *\alpha_\mk]$ and that $[\alpha_\mk'\wedge\alpha_\mk'']_\mk= \dfrac{1}{2}[\alpha_\mk\wedge\alpha_\mk]_\mk$ is real. This completes the proof.\hfill $\square$\\[1mm]
\textbf{Proof of theorem~\ref{thmstrongharmonic}}
We have for all $t\in [0,1]$
\begin{equation}\label{nablat''}
\bar\partial^{\nabla^t}\partial f=\Ad F \left( \bar\partial\alpha_\mk' + [\alpha_\kk''\wedge\alpha_\mk''] + t[\alpha_\mk''\wedge\alpha_\mk']_\mk\right) 
\end{equation}
so that the $\nabla^t$-strongly harmonicity of $f$ is written: 
$$
 \bar\partial\alpha_\mk' + [\alpha_\kk''\wedge\alpha_\mk''] + t[\alpha_\mk''\wedge\alpha_\mk']_\mk=0.\qquad (S_\mk^t)
$$
Then the imaginary part of $\bar\partial^{\nabla^t}\partial f=0$ gives us the $\nabla^t$-harmonicity whereas the real part gives us
$$
d\alpha_\mk + [\alpha_\kk\wedge\alpha_\mk] + t[\alpha_\mk\wedge\alpha_\mk]_\mk=0 \qquad (\mrm{Re}(t))
$$
which is nothing but the lift of the torsion free equation: $f^*T^t=0$, where $T^t=T^{\nabla^t}$. Moreover the projection on $\mk$ of the Maurer-Cartan equation (on $\alpha$) gives us the structure equation
$$
d\alpha_\mk + [\alpha_\kk\wedge\alpha_\mk] + \dfrac{1}{2}[\alpha_\mk\wedge\alpha_\mk]_\mk=0 \qquad  [\mrm{MC}]_\mk
$$
which is nothing but  $(\mrm{Re}(\frac{1}{2}))$ (so that we recover that  $T^{\frac{1}{2}}=0$) but (since it can be written $((\mrm{Re}(0)) +\frac{1}{2}[\alpha_\mk\wedge\alpha_\mk]_\mk=0$) it is also  the lift  of (the $f$-pullback of) the equation expressing the canonical torsion $T^0$ in terms of the Lie bracket (see  theorem~\ref{T-R} or equation~(\ref{Tt})) :
\begin{equation}\label{torsionT^0}
T^0 + [\ ,\ ]_{[\mk]}=0.
\end{equation}
which combining with  the fact that the left hand side of $(\mrm{Re}(t))$ is the lift of $f^*T^t$, gives us back $T^t=(2t-1)[\ ,\ ]_{[\mk]}$ (see (\ref{Tt})).\\
Hence according to (\ref{nablat''}) and $[\mrm{MC}]_\mk$ the strongly harmonicity for one $t\neq\dfrac{1}{2}$ is equivalent to the harmonicity (imaginary part) and $[\alpha_\mk\wedge\alpha_\mk]_\mk=0$ (real part $(\mrm{Re}(t))$ combining with $[\mrm{MC}]_\mk$). We can also simply say that $f$ is strongly harmonic \iif $f$ is harmonic and torsion free i.e. $[\alpha_\mk\wedge\alpha_\mk]_\mk=0$ according to (\ref{Tt}). This proves the equivalence between (i), (ii), and (iii). Now, let us decompose the curvature of $\alpha_\lm$, with respect to powers of $\lm$:
\begin{eqnarray*}
d\alpha_\lm + \dfrac{1}{2}[\alpha_\lm \wedge\alpha_\lm] & = & \lm^{-1}\left( d\alpha_\mk' + [\alpha_\kk\wedge \alpha_\mk']\right) \\
 & + & ( d\alpha_\kk + \dfrac{1}{2}[\alpha_\kk\wedge\alpha_\kk] + \dfrac{1}{2}[\alpha_\mk\wedge\alpha_\mk]_\kk )   + \dfrac{1}{2}[\alpha_\mk\wedge\alpha_\mk]_\mk\\
 &   &  \lm\left( d\alpha_\mk'' + [\alpha_\kk\wedge\alpha_\mk'']\right) 
\end{eqnarray*}
hence using the fact that $\alpha_\lm$ is real (i.e. $\g$-valued) 
$$
d\alpha_\lm + \dfrac{1}{2}[\alpha_\lm \wedge\alpha_\lm]=0 \Leftrightarrow \left\lbrace \begin{array}{l} d\alpha_\mk' + [\alpha_\kk\wedge\alpha_\mk']=0 \quad (S_\mk^0)\\
\frac{}{} [\mrm{MC}]_\kk\\
 {}[\alpha_\mk\wedge\alpha_\mk]_\mk =0  \end{array}\right. \Leftrightarrow  \left\lbrace \begin{array}{l} (S_\mk^0)\\ \frac{}{}[\mrm{MC}]\end{array}\right.,
$$ 
In the last equivalence, we use the fact that $(S_\mk^0) + \overline{(S_\mk^0)}$ is the equation $d\alpha_\mk + [\alpha_\kk\wedge\alpha_\mk]=0$ which combined with $[\mrm{MC} ]_\mk$ (above) gives us $[\alpha_\mk\wedge\alpha_\mk]_\mk =0 $. 
Thus the zero curvature equation on $\alpha_\lm$ is equivalent to the strongly $\nabla^0$-harmonicity, i.e. the strongly $\nabla^t$-harmonicity for all $t\in [0,1]\setminus\{\frac{1}{2}\}$. Finally  the last assertion is obvious. This completes the proof.\hfill $\square$\\[1mm]
We are led naturally to the following definitions.
\begin{defn}
We will say that $f\colon L\to G/K$ is  torsion free if $ f^*T^t=0$ for  $t\in [0,1]\setminus\left\lbrace \frac{1}{2}\right\rbrace $ (this equation does not depend on $t$).
\end{defn}
\begin{defn} In the situation described by theorem~\ref{thmstrongharmonic}-(iv), we will say that the $\g$-valued 1-form on $L$,  $\alpha$, is solution of the the first elliptic system associated to the reductive homogeneous space $G/K$, and that the corresponding geometric map $f$ is a geometric solution of this system.
\end{defn}
\begin{rmk}\em
The independence of the $\nabla^t$-tension field  \wrt to $t$ (and therefore the equivalence between the statements (i) and (ii) in theorem~\ref{thmharmonic}), are in fact a consequence of the geodesic equivalence of the connections $\nabla^t$ (see propositions~\ref{prop-gene-A} and \ref{T-total-antisym}). Moreover, the equivalence of the $\nabla^t$-strongly harmonicities, for $t\neq\frac{1}{2}$, in theorem~\ref{thmstrongharmonic}\,-(i,\,ii), is a particular case of remark~\ref{rmk-strongharmequiv}.
\end{rmk}
\begin{rmk}\label{nobody-strongly}\em
Our so-called "strongly harmonic" maps in homogeneous space already appeared in the literature . It seems that they appeared for the first time in  \cite{12} (see laos \cite{higaki}). Moreover, in \cite{dik}, a concrete example of these maps is given: Gauss maps of CMC surfaces in $H^3$. \\
However, it seems to us very strange  that nobody has remarked before that the torsion free condition $[\alpha_\mk\wedge\alpha_\mk]_\mk =0 $  is already contained in the zero curvature equation $d\alpha_\lm + \dfrac{1}{2}[\alpha_\lm\wedge\alpha_\lm]=0$, $\forall \lm\in S^1$. Indeed everybody in the literature (in particular \cite{12}) adds the  torsion free condition to the zero curvature equation to characterise the torsion free harmonic maps - i.e. strongly harmonic maps -, whereas the strongly harmonicity is exactly and simply equivalent to the zero curvature equation.
\end{rmk}
\subsubsection*{Affine harmonic maps into symmetric spaces}
Now, if we suppose in particular that  $N$ is (locally) symmetric, i.e. $[\mk,\mk]\subset\kk$, then all the connections $\nabla^t$, $0\leq t\leq 1$, coincide. Moreover, if $N$ is also Riemannian then these are equal to the Levi-Civita connection. Therefore we obtain:
\begin{cory}
The first elliptic integrable system associated to a (locally) symmetric space $N=G/K$ is  the equation for $\nabla^0$-harmonic maps $f\colon L\to N$. If $N$ is Riemannian this means that it is the equation for harmonic maps $f\colon L\to N$ (with respect to Levi-Civita in $N$).
\end{cory}
\index{harmonic map|)}
\subsection{Affine/holomorphically harmonic maps into 3-symmetric spaces}
\index{holomorphically harmonic|(}\index{determined maximal@determined, maximal|(}\index{model@model case or system|(}
Let us suppose now that $N=G/G_0$ is a (locally) 3-symmetric space. We use the notations of section~\ref{melliptic}. $N$ is endowed with its canonical almost complex structure $\undj$ defined by (\ref{def-odd}). We continue here the study begun in  \ref{detodcase} concerning the lowest order determined odd system.
\begin{thm}\label{thm3sym}
Let $(L,j)$ be a Riemann surface and $f\colon L\to N$ a smooth map. Let $F\colon L\to G$ be a (local) lift of $f$ and $\alpha=F^{-1}.dF$. Then the following statement are equivalent
\begin{description}
\item[(i)] $\bar\partial\alpha_{-1}' + [\alpha_0''\wedge\alpha_{-1}'] +
[\alpha_1''\wedge\alpha_1']=0 \qquad (S_1)$
\item[(ii)] $\bar\partial\alpha_{1}' + [\alpha_0''\wedge\alpha_{1}']=0 \qquad (S_2)$
\item[(iii)] $f$ is holomorphically $\nabla^1$-harmonic: $\left[ \bar\partial^{\nabla^1}\partial f\right]^{1,0}=0$.
\item[(iv)] $f$ is anti-holomorphically $\nabla^0$-harmonic: $\left[ \bar\partial^{\nabla^0}\partial f\right]^{0,1}=0$.
\item[(v)] $f$ is a geometric solution of the second elliptic integrable system associated to the (locally) 3-symmetric space $G/G_0$:
$$
d\alpha_\lm + \dfrac{1}{2}[\alpha_\lm \wedge\alpha_\lm] =0,\quad \forall \lm\in S^1,
$$
where $\alpha_\lm =\lm^{-2}\alpha_1' + \lm^{-1}\alpha_{-1}' + \alpha_0 + \lm\alpha_1'' + \lm^2\alpha_{-1}''$.
\end{description}
\end{thm}
\textbf{Proof.} The equivalences (i) $\Leftrightarrow$ (ii) $\Leftrightarrow$ (v) have been proved in \ref{detodcase}.
To prove  (i) $\Leftrightarrow$ (iii): just take the $(1,0)$-component in $TN^\C$ of (\ref{nablat''}) for $t=1$. Idem for (ii) $\Leftrightarrow$ (iv). This completes the proof.\hfill $\square$%\medskip\\
\begin{rmk}\em
In fact, the equivalences (iv) $\Leftrightarrow$ (v) and  (iii) $\Leftrightarrow$ (iv) has also been already derived in \ref{detodcase} (theorem~\ref{conclusion-odd-mindeter}), so that this theorem has already been completely proved there. 
\end{rmk}
Now, additionning theorems~\ref{thm3sym}, \ref{thmstrongharmonic}, \ref{thmstrong} and proposition~\ref{holoharm}, we obtain
\begin{cory}
The following  statements are equivalent
\begin{description}
\item[(i)] $f$ is strongly $\nabla^t$-harmonic for one $t\in [0,1]\setminus\left\lbrace \frac{1}{2}\right\rbrace $.
\item[(ii)] $f$ is $\nabla^t$-harmonic for one $t\in [0,1]$ and torsion free.
\item[(iii)] $f$ is holomorphically $\nabla^t$-harmonic for one $t\in [0,1]$ and torsion free.
\item[(iv)] $f$ is a geometric solution of the first elliptic system associated to the reductive homogeneous space $G/G_0$.
\item[(v)] $f$ is a geometric solution of the second elliptic system associated to the 3-symmetric space $G/G_0$, and moreover $ [\alpha_1\wedge\alpha_1] =0$.
\item[(vi)] $f$ is in the same time a geometric solution of the determined odd elliptic systems $(\syst(2,\tau))$ and $(\syst(2,\tau^{-1}))$.
\end{description}
\end{cory}
Now, let us apply theorem~\ref{thm-WZW-Nearly} to the equivalence (iv) $\Leftrightarrow$ (v) of theorem~\ref{thm3sym}.
\begin{thm}
The second elliptic integrable system associated to a naturally reductive 3-symmetric space $N=G/G_0$ is the Euler-Lagrange equation  for the sigma model in $N$ with the Wess-Zumino term defined by the closed 3-form $H=-\dfrac{1}{3}d\Omega_{\undj}$, where $\undj$ is the canonical almost complex structure.
\end{thm}
\index{holomorphically harmonic|)}\index{determined maximal@determined, maximal|)}\index{model@model case or system|)}\index{strongly harmonic|)}
%
%%%%%%%%%%%%%%%%%%%%%%%%%%%%%%%%%%%%%%%%%%%%%%%%%%%%%%%%%%%%%%%%%%%%%%%%%%%%%%%%%%%%%%%%%%%%%%%%%%%%%%%%%%%%%
%   (Affine) vertically (holomorphically) harmonic maps         %%%%%%%%%%%%%%%%%%%%%%%%%%%%%%%%%%%%%%%%%%%%%%%
%%%%%%%%%%%%%%%%%%%%%%%%%%%%%%%%%%%%%%%%%%%%%%%%%%%%%%%%%%%%%%%%%%%%%%%%%%%%%%%%%%%%%%%%%%%%%%%%%%%%%%%%%%%%%
%
\subsection{(Affine) vertically (holomorphically) harmonic maps}\index{vertically harmonic|(}\index{vertical tension field|(}\index{strongly vertically harmonic|(}\index{vertically holomorphically harmonic|(}
\subsubsection{Affine vertically harmonic maps: general properties}
Here we generalise the definition of vertical harmonicity for maps from a Riemannian surface into an affine manifold.
\begin{defn}\label{def-aff-vert-harm}
Let $(N,\nabla)$ be an affine manifold. Let us suppose that we have a splitting $TN=\ver\oplus \hor$. In other words $N$ is endowed with a Pfaffian system (the vertical subbundle $\ver$) and with a connection on this Pfaffian system. Let $f\colon (L,b)\to N$ be a smooth map from a Riemannian manifold $(L,b)$ into $N$. Then we set
$$ 
\tau^v(f)=\mrm{Tr}_b(\nabla^v d^v f) =*d^{\nabla^v}*d^v f,
$$
where $\nabla^v d^v f$ is the vertical component of the covariant derivative of $df$ with respect to the  connection on $T^*L\otimes f^*TN$ induced by the Levi-Civita connection of $L$ and the linear connection $\nabla$. We will say that $f$ is affine vertically harmonic with respect to $\nabla$ or $\nabla$-vertically harmonic if $\tau^v(f)=0$.
\end{defn}
\begin{thm}\label{strongly-vertical-harm}
Let $(L,j)$ be a Riemann surface and $f\colon (L,j)\to (N,\nabla)$ a smooth map. Then we have 
$$
2 \bar\partial^{\nabla^v} \partial^v f = d^{\nabla^v} d^v f + i d^{\nabla^v} *d^v f,
$$
moreover $d^{\nabla^v} d^v f=f^*T^v$, where $T^v$ is the vertical torsion (see \ref{phitorsion})  and $d^{\nabla^v} *d^v f=\tau^v(f) \mrm{vol}_b$ for any hermitian metric $b$ in $L$.
Therefore the following statements are equivalent:
\begin{description}
\item[(i)] $\left( \nabla''\right)^v \partial^v f=0$.
\item[(ii)] $\bar\partial^{\nabla^v}\partial^v f=0$.
\item[(iii)] $\nabla_{\dl{}{\overline z}}^v\left(\dl{^v f}{z}\right)=0$, for any holomorphic local coordinate $z=x+iy$ (i.e. $(x,y)$ are conformal coordinates for any hermitian metric in $L$).
\item[(iv)] $f$ is $\nabla^v$-vertically harmonic with respect to any hermitian metric in $L$ and  vertically torsion free: $f^*T^v=0$ (i.e. $T^v(\dl{f}{x},\dl{f}{y})=0$ for any  conformal coordinates $x,y$).
\end{description}
We will say in this case that $f$ is strongly $\nabla$-vertically harmonic.
\end{thm}
\begin{rmk}\em
In the previous characterisation of strongly vertical harmonicity, we only need a connection $\nabla^v$ in $\ver$, and $\nabla$ is useless.
\end{rmk}
\subsubsection{Affine vertically holomorphically harmonic maps}
Here we generalize the notion of holomorphic harmonicity by introducing a new notion of vertical holomorphic harmonicity (in the same way that the vertical harmonicity generalizes the harmonicity).
\begin{defn}
Let $(N,\nabla)$ be an affine manifold with a splitting $TN=\ver\oplus \hor$ as in definition~\ref{def-aff-vert-harm}. Let us suppose that the subbundle $\ver$ admits a complex structure $J^\ver$. Let us denote by $\ver^\C=\ver^{1,0}\oplus\ver^{0,1}$ the splitting induced by the complex  structure $J^\ver$. Then we will say that  a map $f\colon (L,j_L)\to N$ from a Riemann surface into $N$, is vertically holomorphicaly harmonic with respect to $\nabla$ or $\nabla$-vert. hol. harmonic if
$$
\left[ \bar\partial^{\nabla^v}\partial^v f\right]^{1,0}=0.
$$ 
\end{defn}
\begin{thm}\label{thm-vert-hol-harm-charact}
Let  $(L,j_L)$ be a Riemann surface and $(N,\nabla)$ be an affine manifold with a splitting $TN=\ver\oplus \hor$ and a complex structure $J^\ver$ on $\ver$ as in the previous definition. Then $f\colon L\to N$ is $\nabla$-vert. hol. harmonic \iif 
$$
f^*T^v + J^vd^{\nabla^v}*d^vf=0.
$$
\end{thm}
\textbf{Proof.} The same as the one of proposition~\ref{holoharm}.\index{vertically harmonic|)}\hfill $\square$ \medskip\\
\index{vertically holomorphically harmonic|)}
%
%
%
%%%%%%%%%%%%%%%%%%%%%%%%%%%%%%%%%%%%%%%%%%%%%%%%%%%%%%%%%%%%%%%%%%%%%%%%%
%   Affine vertically harmonic maps into reductive homogeneous space
%%%%%%%%%%%%%%%%%%%%%%%%%%%%%%%%%%%%%%%%%%%%%%%%%%%%%%%%%%%%%%%%%%%%%%%%%%
%
\subsection{Affine vertically harmonic maps into reductive homogeneous space}\label{Affvertharmredhom}
\index{vertically harmonic|(}
Let $G$ be a Lie group, and $K\subset H\subset G$ subgroups of $G$ such that $M=G/H$ and $H/K$ are reductive. We use the notations of \ref{homspacefibr} (but we do not suppose a priori that the reductive homogeneous spaces are Riemannian).
\begin{thm}\label{thmvertiharmonic}
Let $(L,j)$ be a Riemann surface and $f\colon L \to N=G/K$ be a smooth map, $F\colon L\to G$ a (local) lift of $f$ and $\alpha=F^{-1}.dF$. Then the following statements are equivalent:
\begin{description}
\item[(i)] $f$ is $\nabla^t$-vertically harmonic for one $t\in [0,1]$.
\item[(ii)] $f$ is $\nabla^t$-vertically harmonic for all $t\in [0,1]$.
\item[(iii)] $d*\alpha_\pk + [\alpha_\kk\wedge *\alpha_\pk]=0$.
\item[(iv)] $\mrm{Im}\left( \bar\partial\alpha_\pk' + [\alpha_\kk''\wedge\alpha_\pk'] + t [\alpha_\pk''\wedge \alpha_\pk']_\pk\right) =0$, $\forall t\in [0,1]$.
\end{description}
The vertical tension field $\tau^{t,v}(f)$, with respect of $\nabla^t$, is independent of $ t\in [0,1]$.
\end{thm}
\textbf{Proof.} 
Setting $\nabla^{t,v}=(\nabla^t)^v$, we have 
\begin{eqnarray*}
\tau^{t,v}(f) =*d^{\nabla^{t,v}} *d^v f & = & * \Ad F (d*\alpha_\pk + [\alpha_\kk\wedge *\alpha_\pk] + t[\alpha_\nk\wedge *\alpha_\pk]_\pk)\\
 & = & * \Ad F (d*\alpha_\pk + [\alpha_\kk\wedge *\alpha_\pk] + t[\alpha_\pk\wedge *\alpha_\pk]_\pk)\\
& & \text{ since } [\mk,\pk]\subset [\mk,\hk]\subset \mk\\
& = & * \Ad F (d*\alpha_\pk + [\alpha_\kk\wedge *\alpha_\pk]).
\end{eqnarray*}
This gives us the equivalences (i) $\Leftrightarrow$ (ii) $\Leftrightarrow$ (iii) as well as the last assertion of the theorem.
Moreover let us compute the complex second derivative:
$$
\bar\partial^{\nabla^{t,v}}\partial^v f= * \Ad F (\Bar \partial \alpha_\pk' + [\alpha_\kk''\wedge \alpha_\pk'] + [\alpha_\pk''\wedge\alpha_\pk']_\pk).
$$
Then the equivalence (ii) $\Leftrightarrow$ (iv) follows from theorem~\ref{strongly-vertical-harm}. For this equivalence, we could also  remark that $[\alpha_\pk''\wedge\alpha_\pk']_\pk=\dfrac{1}{2}[\alpha_\pk\wedge\alpha_\pk]_\pk$ is in the real subspace $\pk$ and that $2\,\mrm{Im}\left( \bar\partial\alpha_\pk' + [\alpha_\kk''\wedge\alpha_\pk']\right)=d*\alpha_\pk + [\alpha_\kk\wedge *\alpha_\pk]$.
This completes the proof.\hfill $\square$\medskip\\
According to theorem~\ref{thmharmonic}, we deduce the following.
\begin{cory}\label{vertdetau=tauv},
In the situation of the previous theorem, if $f\colon L\to N$ is harmonic (\wrt some $\nabla^t$) then it is also vertically harmonic. More generally, the vertical tension field is the vertical component of the tension field.
\end{cory}
\begin{rmk}\em
In the case $N=G/K$ is endowed with a naturally reductive metric, then $\nabla^\frac{1}{2}$ coincides with the Levi-civita connection and therefore the previous corollary is nothing but a particular case of theorem~\ref{m,n}-(iii).
\end{rmk}
Now, let $f\colon L\to N$ be an arbitrary map from a Riemann surface into $N$.
Then the $f$-pullback of the vertical torsion with respect ot $\nabla^t$ is 
\begin{eqnarray*}
f^*T^{t,v}= d^{\nabla^{t,v}} d^v f & = & \Ad F\left( d\alpha_\pk + [\alpha_\kk\wedge\alpha_\pk] + t[\alpha_\nk\wedge\alpha_\pk]_\pk \right) \\
 & = & \Ad F\left( d\alpha_\pk + [\alpha_\kk\wedge\alpha_\pk] + t[\alpha_\pk\wedge\alpha_\pk]_\pk \right) \\
 & = & f^*\left(  T^{0,v} + t[\phi\wedge\phi]_{[\pk]}\right) 
\end{eqnarray*}
where $\phi\colon TN\to [\pk]$ is the projection on the vertical subbundle along the horizontal subbundle $[\mk]$. Therefore 
\begin{equation}\label{eq-Tt,v}
T^{t,v}=T^{0,v} + t[\phi\wedge\phi]_{[\pk]}.
\end{equation}
Moreover, recall that, according to section~\ref{homspacefibr}, the projection on $[\pk]$ of the Maurer-Cartan equation gives us  the homogeneous structure equation (see equations (\ref{lift-T}), (\ref{formulePhi}) and footnote \ref{HSE}) 
$$
T^{0,v}=\Phi -\dfrac{1}{2}[\phi\wedge\phi]_{[\pk]}
$$
where $\Phi=-\dfrac{1}{2}[\psi\wedge\psi]_{[\pk]}$ is the homogeneous curvature form and $\psi\colon TN\to [\mk]$ is the projection on $[\mk]$ along $[\pk]$. Then we have 
\begin{equation}\label{Ttvreductive}
T^{t,v}  =   \Phi + \left(t-\dfrac{1}{2}\right)[\phi\wedge\phi]_{[\pk]}.
\end{equation}
Therefore
\begin{thm}\label{the2affineconditions}
Let us consider  the same situation as in theorem~\ref{thmvertiharmonic}\medskip.\\
$\bullet$ If f is flat then the strongly $\nabla^t$-vertical harmonicity and  the freedom from torsion, for $f$, do not depend on $t$, if $t\in [0,1]\setminus \left\lbrace \frac{1}{2}\right\rbrace $.\\
Moreover $T^{\frac{1}{2},v}=\Phi$ so that (if $f$ is flat) strongly vertical harmonicity and vertical harmonicity with respect to $\nabla^{\frac{^1}{2}}$ are equivalent.\\[1.5mm]
$\bullet$ If $H/K$ is locally symmetric, i.e. $[\pk,\pk]\subset \kk$, then $\forall t\in [0,1]$, $T^{t,v}=\Phi$.\\
In particular, the $\nabla^t$-vertical torsion does not depend on $t\in [0,1]$, and thus neither does strongly harmonicity.
\end{thm}
Now, we can conclude
\begin{cory}\label{nablatsyst}\index{determined minimal@determined, minimal|(}
Let us suppose now that $N=G/K$ is a (locally) $2k$-symmetric space and that $M=G/H$ is the corresponding (locally) $k$-symmetric space. Then the even minimal determined system $(\syst(k,\tau))$ associated to $N$ means that the geometric map $f\colon L\to N$ is horizontally holomorphic and vertically harmonic with respect to any linear connection $\nabla^t$, $0\leq t\leq 1$. Moreover the horizontal holomorphicity implies the flatness of $f$ and thus its freedom from vertical torsion  (with respect to any connection $\nabla^t$, $0\leq t\leq 1$). More precisely the (last) equation $(S_k)$ of the system means
$$
\bar\partial^{\nabla^{t,v}}\partial^v f=0
$$
i.e. that $f$ is strongly $\nabla^t$-vertically harmonic, so that its real part means that $f$ is vertically torsion free and its imaginary part that $f$ is vertically harmonic.
\end{cory}
\subsubsection*{The Riemannian case}
Now, let us suppose that the reductive homogeneous space $M=G/H$ is Riemannian, and then so is $N=G/K$. In other words, we are in the situation described by \ref{homspacefibr}. Let us consider the metric connections in $N$:
$$
\nablamet{}{t}=\nabla^0 + t\mB^N, \qquad 0\leq t\leq 1
$$
with $\mB^N=[\ ,\ ]_{[\nk]} + \mU^N$ and $\mU^N$ defined by equation (\ref{defofU}).\\
For any $\Ad K$-invariant subspace $\mak l\subset \nk$, we will denote by $\mU^{\mak l}\colon \mak l \times \mak l\to \mak l$ the bilinear symmetric map defined by
$$
\langle \mU^{\mak l}(X,Y),Z\rangle = \langle [Z,X]_{\mak l},Y\rangle + \langle X,[Z,Y]_{\mak l}\rangle\quad \forall X,Y\in \mak l,
$$ 
and by $\mU^{[\mak l]}$ its extension to the subbundle $[\mak l]\subset TN$. Then we have in particular
$$
\mU^N=\mU^{[\nk]} \quad \text{and}\quad \mU=\mU^{[\pk]}
$$
where, let us recall it,  $\mU$ is defined  by (\ref{defofmU}).\\
Now, let us project the  definition equation of $\nablamet{}{t}$ in the vertical subbundle: we obtain $\forall V\in \mal C(TN)$,
$$
\phi(\nablamet{}{t} V)=\nabla^0 \phi V + t\phi\circ\mB^N(\cdot, V).
$$
Moreover, according to \ref{homspacefibr}, equation \eqref{eq-phi-B=B-phi}, we have $\phi\circ\mB^N=\phi^*\mB -\Phi$
so that
$$
\phi(\nablamet{}{t} V)=\nabla^0 \phi V + t(\phi^*\mB -\Phi)(\cdot, V)
$$
and in particular $\forall V\in \mal C(\ver)$, $\forall A \in TN$,
\begin{equation}\label{eq-phi-nablamet-t}
\phi(\nablamet{A}{t} V)= \nabla_A^0 V + t\left( [\phi A, V]_{[\pk]} + \mU^{[\pk]}(\phi A, V)\right) .
\end{equation}
Then according to theorem~\ref{coderivative} and remark~\ref{naturreductiv}, and denoting by $\nablamet{}{t,v}$ the vertical component of $\nablamet{}{t}$,  it follows:
\begin{thm}\label{thm-coincide-even} 
\begin{description}
\item[$\bullet$] If $H/K$ is naturally reductive, then the connections defined by the restriction to $\ver$ of $\nablamet{}{t,v}$, $0\leq t\leq 1$, are all $\phi$-equivalent. Therefore the vertical harmonicity, with respect to $\nablamet{}{t}$, is the same for all $0\leq t\leq 1$.
\item[$\bullet$] If $H/K$ is locally symmetric, then all the $\nablamet{}{t,v}$, $0\leq t\leq 1$, coincide in $\ver$. In particular the strongly vertical harmonicity coincides for all the connections $\nablamet{}{t}$, $0\leq t\leq 1$.
\end{description}
\end{thm}
\textbf{The vertical torsion of $\nablamet{}{t}$.}
We have seen in \ref{family} that the torsion of  $\nablamet{}{t}$ is the same as the one  of $\nabla^t$. Now let us see what happens for the vertical torsion. The vertical torsion with respect to $\nablamet{}{t,v}$ is given, according to \eqref{eq-phi-nablamet-t}, by 
\begin{eqnarray*}
\overset{\mrm met}{T}{}^{t,v}= d^{\nablamet{}{t,v}}\phi & = & d^{\nabla^0} \phi + t [\phi\wedge\phi]_{[\pk]}+ \mU^{[\pk]}(\phi\wedge\phi)\\
 & = & T^{0,v} + t [\phi\wedge\phi]_{[\pk]} + \mU^{[\pk]}(\phi\wedge\phi).
\end{eqnarray*}
But $ \mU^{[\pk]}(\phi\wedge\phi)=0$ because $ \mU^{[\pk]}$ is symmetric, so that, according to \eqref{eq-Tt,v},
\begin{equation}\label{sameverticaltorsion}
\overset{\mrm met}{T}{}^{t,v}= T^{t,v}.
\end{equation} 
\begin{rmk}\em
We see that the value $t=\frac{1}{2}$, i.e. the Levi-Civita connection, plays a special role according to theorem~\ref{the2affineconditions} and equation (\ref{sameverticaltorsion}). Indeed for the Levi-Civita connection, we always have $\overset{\mrm met}{T}{}^{\frac{1}{2},v}=\Phi$, so that if $f$ is flat, the strongly harmonicity and  the vertical harmonicity are equivalent.\\
However, if $H/K$ is (locally) symmetric, then we have $\forall t\in [0,1]$, $\overset{\mrm met}{T}{}^{t,v}= T^{t,v}=\Phi$, and we have even more, since all the connections $\nablamet{}{t,v}$ coincides on $\ver$. Therefore the special role played by the Levi-Civita connection is shared, in this case, with all the other connections $\nablamet{}{t}$.
\end{rmk}
\begin{rmk}\em
As already mentioned in section~\ref{homspacefibr}, the equation \eqref{eq-phi-B=B-phi}, i.e. $\phi\circ\mB^N=\phi^*\mB -\Phi$, can be  obtained directly by computation (without using the general theorem~\ref{difference-tensor}, as in section~\ref{homspacefibr}). Let us do this. Lifting this equation we then have to prove that $\forall X,Y\in \g$,
$$
[ X_\nk, Y_\nk]_\pk  + \left[ \mU^\nk(X_\nk,Y_\nk)\right]_\pk  = [ X_\pk, Y_\pk]_\pk  + \mU^\pk(X_\pk,Y_\pk) + [X_\mk, Y_\mk]_\pk, 
$$ 
since $[\mk,\pk]_\pk=\{0\}$. It then suffices  to prove that   $ \left[ \mU^\nk(X_\nk,Y_\mk)\right]_\pk=0$.
This equality holds. Indeed, we have $\forall Z\in \g$,
$$
\langle \mU^\nk(X_\nk,Y_\mk), Z_\pk \rangle = \langle [Z_\pk,X_\nk]_\nk,Y_\mk\rangle  +  \langle  X_\nk, [Z_\pk,Y_\mk]_\nk \rangle = \langle [Z_\pk,X_\mk] ,Y_\mk\rangle  +  \langle  X_\mk, [Z_\pk,Y_\mk] \rangle =0
$$
since $\adm \pk\subset\adm \hk\subset \so(\mk)$.
\end{rmk}
\textbf{The metric geometric interpretation of the even minimal determined system.}
Now, according to theorem~\ref{thm-coincide-even}, we can conclude by rewriting corollary~\ref{nablatsyst} in terms of the metric connection $\nablamet{}{t}$ instead of the linear connection $\nabla^t$.
\begin{cory}
Let us suppose now that $N=G/K$ is a (locally) $2k$-symmetric space and that $M=G/H$ is the corresponding (locally) $k$-symmetric space. Then the even minimal determined system $(\syst(k,\tau))$ associated to $N$ means that the geometric map $f\colon L\to N$ is horizontally holomorphic and vertically harmonic with respect to any metric connection $\nablamet{}{t}$, $0\leq t\leq 1$. Moreover the horizontal holomorphicity implies the flatness of $f$ and thus its freedom from vertical torsion  (with respect to any connection $\nablamet{}{t}$, $0\leq t\leq 1$). More precisely the (last) equation $(S_k)$ of the system means
$$
\bar\partial^{\nablamet{}{t,v}}\partial^v f=0
$$
i.e. that $f$ is strongly $\nablamet{}{t}$-vertically harmonic, so that its real part means that $f$ is vertically torsion free and its imaginary part that $f$ is vertically harmonic.
\end{cory}
\begin{rmk}\em
In particular for $t=\frac{1}{2}$, we recover theorem~\ref{evendeterhom}.
\end{rmk}
\index{vertically harmonic|)}\index{determined minimal@determined, minimal|)}
%
%
%%%%%%%%%%%%%%%%%%%%%%%%%%%%%%%%%%%%%%%%%%%%%%%%%%%%%%%%%%%%%%%%%%%%%%%%%%%%%%%%%%%%%%%
%
%
\subsection{Harmonicity vs vertically harmonicity}\label{sect-harm-vs-vertharm}
\index{vertically harmonic|(}\index{harmonic map|(}\index{strongly harmonic|(}
In this subsection, we endow all the reductive homogeneous spaces with their connections $\nabla^t$. The harmonicity, the strongly and vertical harmonicity are therefore considered \wrt to these connections. These notions are independent of $t$, with a  exception at $t=\frac{1}{2}$ for strongly harmonicity, see \ref{Aff-harm-red-hom} and \ref{Affvertharmredhom}). When we will say "strongly harmonic" without other precisions, we will mean "\wrt to any $\nabla^t$, $t\neq \frac{1}{2}$". Let us recall that if a naturally reductive invariant metric is given, then $\nabla^{\frac{1}{2}}$ coincides with Levi-Civita. We use the notations of the previous sections: in particular, $L$ is a Riemann surface. Moreover, we keep in mind the results of \ref{Aff-harm-red-hom} and \ref{Affvertharmredhom}.
%, given a map $f\colon L \to G/K$, we denote by $\alpha$ an associated $\g$-valued Maurer-Cartan form.
%
%
\begin{thm}\label{thm-equiv-harm&vertharm-t=0}
Let  $N=G/K$ be a (locally) $2k$-symmetric space. Let $f\colon L\to N=G/K$ be horizontally holomorphic. Then the following statements are equivalent:
\begin{description}
\item[(i)]  $f$ is harmonic,  
\item[(ii)] $f$ is vertically harmonic  and torsion free,
\item[(iii)] $f$  is strongly harmonic. 
\end{description}
\end{thm}
\proof 
$\bullet$ Let us prove: (i) $\Rightarrow$ (ii) and (i) $\Leftrightarrow$ (iii). According to corollary~\ref{vertdetau=tauv}, the harmonicity implies the vertical harmonicity. 
We have seen that the horizontal holomorphicity implies the flatness so that $[\alpha_\nk\wedge\alpha_\nk ]_{\g_k}=[\alpha_\mk\wedge\alpha_\mk]_{\g_k}=0$. 
Moreover we know that the $\nabla^{\frac{1}{2}}$-harmonicity coincides with the strongly $\nabla^{\frac{1}{2}}$-harmonicity so that we have ${\bar\partial}^{\nabla^{\frac{1}{2}}}\partial f=0$. Writing this equation in terms of $\alpha$ and projecting it on $\mk^{0,1}=\sum_1^{k-1} \g_j^\C$, yields
$$
\bar\partial\alpha_j' + [\alpha_0''\wedge \alpha_j'] + \frac{1}{2}[\alpha_\nk''\wedge\alpha_\nk']_{\g_j}=0,\quad 1\leq j\leq k-1.
$$
Now,  the horizontal holomorphicity $\alpha_j'=0$, $1\leq j\leq k-1$,  reduces this equation to $[\alpha_\nk''\wedge\alpha_\nk']_{\g_j}=0$, $1\leq j\leq k-1$, i.e. $[\alpha_\nk''\wedge\alpha_\nk']_\mk=0$. Therefore we obtain finally $[\alpha_\nk''\wedge\alpha_\nk']_\nk=0$ i.e. $f$ is torsion free.  Combining the torsion freedom with  the strongly $\nabla^{\frac{1}{2}}$-harmonicity, we obtain that $f$ is strongly  harmonic \wrt any connection $\nabla^t $.\\
$\bullet$  Now we prove: (ii) $\Rightarrow$ (iii). In terms of the Maurer-Cartan form $\alpha$, the torsion freedom of $f$ is written $[\alpha_\nk\wedge\alpha_\nk]_\nk=0$. Moreover the horizontal holomorphicity implies the flatness: $[\alpha_\mk\wedge\alpha_\mk]_\pk=0$ so that $[\alpha_\nk\wedge\alpha_\nk]_\pk=0$. Finally we have $[\alpha_\nk\wedge\alpha_\nk]_\mk=0$ and therefore, the projection on $\mk$ of the Maurer-Cartan equation becomes 
$$
d\alpha_\mk + [\alpha_0\wedge\alpha_\mk]=0.
$$
Then projecting it on $\mk^{1,0}$ and using the horizontal holomorphicity, we obtain
$$
\bar\partial\alpha_\mk' + [\alpha_0''\wedge\alpha_\mk']=0,
$$
which is nothing but the horizontal part of the strongly $\nabla^0$-harmonicity for $f$. The vertical part comes from the vertical harmonicity and the fact that the flatness implies the vertical torsion freedom (see theorem~\ref{the2affineconditions} or corollary~\ref{nablatsyst}). \comprf\hsq\medskip\\
From the expression of the tension field $\tau(f)$ computed in section~\ref{Aff-harm-red-hom}, we deduce the following. 
\begin{prop}
Let $G$ be a Lie group, and $K\subset H\subset G$ subgroups of $G$ such that $M=G/H$ and $H/K$ are reductive. Let $f\colon L \to N=G/K$ be a map and $u=\pi\circ f$ its projection on $M=G/H$.  Then the horizontal component of the tension field of $f$ is  related to the tension field of $u$  by:
$$
\tau^h(f)=\tau(u) + *[d^v f\wedge *d^hf].
$$
In the case $N$ is endowed with a naturally reductive invariant metric $h$, inducing then a  naturally reductive invariant metric\footnote{See proposition~\ref{prop-nat-red-induce-appendix}} $g$ on $M$, and $N$ and $M$ are resp. endowed with their Levi-Civita connections, this relation is nothing but a particular case of theorem~\ref{thm-LtoNvertHarm-vs-Harm}-(ii). In particular, strongly flatness  of $f$ means $[d^v f\wedge *d^hf]=0$. If $f$ is strongly flat and vertically harmonic, then its harmonicity is equivalent to the harmonicity of its projection $u$ on $M$.
\end{prop}
\textbf{Convention.} We will continue to say that $f$ is \emph{strongly flat} when $[d^v f\wedge *d^hf]=0$, even if $N$ and $M$ are not endowed with  metrics.
\begin{rmk}\em
We have $-[V,H]=T(V,H)=\rho(V)H$, $\forall V\in \ver, H\in\hor$, where $T$ is the torsion of the canonical connection on $N$ and $\rho\colon \ver\to \so(\hor)$ is the linear representation of the curvature $\Phi$ defined by \eqref{rho-def}.\\
Suppose that $N=G/K$ is loc. $2k$-symmetric and endowed with its canonical $f$-structure, then using the fact $(\adm V)\undj_0=-\undj_0(\adm V)$, $\forall V\in \g_k$, we obtain that $N_F(V,H)=-2[V,H]$, $\forall V\in\ver, H\in\hor$. Therefore $f$ is strongly flat \iif $\mrm{Tr}\left( f^*N_{F|\ver\times\hor}\right) =0$.
\end{rmk}
\begin{prop}
Let  $N=G/K$ be a (locally) $2k$-symmetric space and $f\colon L\to N=G/K$ a map. We suppose that $f\colon L\to N=G/K$ is horizontally holomorphic and that its projection $u=\pi\circ f\colon L\to M=G/H$, on the associated $k$-symmetric space, is torsion free. Then the following statements are equivalent:
\begin{description}
\item[(i)] $f$ is harmonic.
\item[(ii)] $u$ is harmonic and $f$ is vertically harmonic.
\item[(iii)] $f$ is strongly harmonic.
\item[(iv)] $u$ is strongly harmonic and $f$ is strongly vertically harmonic.
\end{description}
\end{prop}
\proof According to theorem~\ref{thm-equiv-harm&vertharm-t=0}, we have (i) $\Leftrightarrow $ (iii). Moreover, since $u$ is torsion free, its harmonicity is equivalent to its strongly harmonicity. Furthermore, since $f$ is flat because horizontally holomorphic then it is vertically torsion free, according to theorem~\ref{the2affineconditions}. Therefore (ii) $\Leftrightarrow $ (iv). Now, we have $\dfrac{1}{2}[\alpha_\nk\wedge\alpha_\nk]_\mk = [\alpha_\pk\wedge\alpha_\mk] +  \dfrac{1}{2}[\alpha_\mk\wedge\alpha_\mk]_\mk=[\alpha_\mk\wedge\alpha_\pk] $ because $u$ is torsion free. Besides,  (i) implies that $f$ is torsion free, according to theorem~\ref{thm-equiv-harm&vertharm-t=0}, hence $[\alpha_\pk\wedge \alpha_\mk]=0$ so that owing to the horizontal holomorphicity $[\alpha_\pk\wedge *\alpha_\mk]= [\alpha_\pk\wedge -\undj_0\alpha_\mk]=\undj_0[\alpha_\pk\wedge \alpha_\mk]=0$. Therefore $f$ is strongly flat which gives us: (i) $\Rightarrow$ (ii). Conversely, let us suppose (ii) and therefore (iv). Then writting the equation of strongly harmonicity for $u$, projecting it on $\mk^{0,1}$, and using the horizontal holomorphicity of $f$ (same method as in the proof of theorem~\ref{thm-equiv-harm&vertharm-t=0}), we obtain $[\alpha_\pk''\wedge\alpha_{\mk^{1,0}}']=0$ and hence $f$ is strongly flat so that (ii) $\Rightarrow$ (i). \comprf\hsq\medskip\\
In particular, if $k=2$ then $M=G/H$ is loc. symmetric: $[\mk,\mk]_\mk=\{0\}$, so that we recover:
\begin{cory}\emph{\cite{ki3}}
Let  $N=G/K$ be a (locally) $4$-symmetric space. Let $f\colon L\to N=G/K$ be a map and  $u\colon L\to M=G/H$ its projection on the associated symmetric space. We suppose that $f\colon L\to N=G/K$ is horizontally holomorphic. Then 
the following statements are equivalent:
\begin{description}
\item[(i)] $f$ is harmonic.
\item[(ii)] $u=\pi\circ f$ is harmonic and $f$ is vertically harmonic.
\end{description}
\end{cory}
\index{vertically harmonic|)}\index{harmonic map|)}\index{strongly harmonic|)}

%\begin{cory}
% endowed with a $G$-invariant naturally reductive metric. 
%
%
%
%
%%%%%%%%%%%%%%%%%%%%%%%%%%%%%%%%%%%%%%%%%%%%%%%%%%%%%%%%%%%%%%%%%%%%%%%%%%%%%%%%%%%%%%%%%%%%%
%   Affine vertically (holomorphically) harmonic maps into reductive homogeneous space      %
%%%%%%%%%%%%%%%%%%%%%%%%%%%%%%%%%%%%%%%%%%%%%%%%%%%%%%%%%%%%%%%%%%%%%%%%%%%%%%%%%%%%%%%%%%%%%
%
\subsection{(Affine) vertically (holomorphically) harmonic maps into reductive homogeneous space with an invariant Pfaffian structure}\label{subsec-odd-vert-hol-harm-Pfaff}
\index{vertically harmonic|(}\index{vertically holomorphically harmonic|(}
Let $N=G/K$ be a reductive homogeneous space and $\g=\kk\oplus\mk$ a reductive decomposition of $\g$. Let us suppose that $\mk$ admits an $\Ad K$-invariant decomposition
$$
\mk =\mk'\oplus \pk.
$$
Then $\pk$ defines a vertical subbundle $\ver=[\pk]$ and $\mk'$ an horizontal subundle $\hor=[\mk']$ giving a splitting $TN=\hor\oplus\ver$.\\
 The curvature of the horizontal distribution $\hor$ is given by 
$$
R^\hor=-[\psi,\psi]_{[\pk]}=-\dfrac{1}{2}[\psi\wedge\psi]_{[\pk]}
$$
where $\psi\colon TN\to [\mk']$ is the projection on $[\mk']$ along $[\pk]$. We will set 
$$
\Phi:=R^\hor.
$$
The vertical torsion of the linear connection $\nabla^t$ is given by $T^{t,v}=d^{(\nabla^t)^v} \phi$ and lifts into
\begin{eqnarray}
\widetilde{T^{t,v}} & = &  d\theta_\pk + [\theta_\kk\wedge\theta_\pk] + t[\theta_\mk\wedge\theta_\pk]_\pk \label{liftTtv1}\\
   & = &  d\theta_\pk + [\theta_\kk\wedge\theta_\pk] + t[\theta_{\mk'}\wedge\theta_\pk]_\pk + t[\theta_{\pk}\wedge\theta_\pk]_\pk  \label{liftTtv2} .
\end{eqnarray}
On the other hand, the projection on $\pk$ of the Maurer-Cartan equation gives 
$$
d\theta_\pk + [\theta_\kk\wedge\theta_\pk] + \dfrac{1}{2}[\theta_{\mk'}\wedge\theta_{\mk'}]_\pk +  [\theta_{\mk'}\wedge \theta_\pk]_\pk + \dfrac{1}{2}[\theta_{\pk}\wedge \theta_\pk]_\pk=0
$$
so that (\ref{liftTtv2}) can be written 
$$
\widetilde{T^{t,v}}=-\dfrac{1}{2}[\theta_{\mk'}\wedge\theta_{\mk'}]_\pk + (t-1)[\theta_{\mk'}\wedge \theta_\pk]_\pk + \left( t-\dfrac{1}{2}\right)[\theta_{\pk}\wedge \theta_\pk]_\pk
$$
which projected in $N$ becomes 
$$
T^{t,v}= \Phi + (t-1)[\psi\wedge\phi]_{[\pk]} + \left( t-\frac{1}{2}\right) [\phi\wedge\phi]_{[\pk]}.
$$
We remark that the values $t=\frac{1}{2},1$ play special roles. In particular:
\begin{description}
\item[$\bullet$] If $[\mk',\pk]_\pk=\{0\}$ then we have $T^{\frac{1}{2},v}=\Phi$. More generally we recover equation~(\ref{Ttvreductive}) and  the results of theorem~\ref{the2affineconditions} (by taking the  following values in the notations $\mk:=\nk$ and $\mk':=\mk$).
\item[$\bullet$] If $[\pk,\pk]_\pk=\{0\}$ then we have $T^{1,v}= \Phi$.
\end{description}
Now, if the two conditions are satisfied, $[\mk',\pk]_\pk=[\pk,\pk]_\pk=\{0\}$, then we have
$$
\forall t \in [0,1],\quad T^{t,v}=\Phi.
$$
Now let $f\colon (L,j_L)\to N$ be a map from a Riemann surface into $N$.
Let us compute the vertical tension field $\tau^{t,v}(f)$ of $f$ with respect to $\nabla^t$ (and some Hermitian metric $b$ in $L$). In order to do that, let $F\colon L\to G$ be a lift  of $f$ and $\alpha=F^{-1}.dF$. Then we have
\begin{eqnarray*}
\tau^{t,v}(f) =  *d^{\nabla^{t,v}} * d^v f & = &  *\Ad F\left( d *\alpha_\pk + [\alpha_\kk\wedge *\alpha_\pk] + t [\alpha_\mk\wedge *\alpha_\pk]_\pk \right)\\
& = & *\Ad F\left( d *\alpha_\pk + [\alpha_\kk\wedge *\alpha_\pk] + t [\alpha_{\mk'}\wedge *\alpha_\pk]_\pk + t [\alpha_\pk\wedge *\alpha_\pk]_\pk\right)\\
& = & \tau^{0,v}(f) + t * [f^*\psi \wedge *(f^*\phi)]_{[\pk]} \\
 & = & \tau^{0,v}(f) + t \mrm{Tr}_b\left( [f^*\psi,f^*\phi]_{[\pk]}\right) 
\end{eqnarray*}
Now, let us consider the $\Ad K$-invariant vector subspace
$$
\mk_* =\left\{X\in\mk'|\, [X,\pk]_\pk=\{0\}\right\}
$$
and let $\mk_1$ be  an $\Ad K$-invariant complement\footnote{Such an $\Ad K$-invariant complement always exists if $\mk'$ admits an $\Ad K$-invariant Pseudo-Euclidean inner product non degenerated on $\mk_*$. For example if $\g$ is semisimple, take the restriction to $\mk'$ of the Killing form and then $\mk_1=([\pk,\pk]_{\mk'})^\perp$.} of $\mk_*$ in $\mk'$
$$
\mk'=\mk_*\oplus\mk_1.
$$
Then we can rewrite the $\nabla^t$-vertical torsion in the form 
\begin{equation}\label{Tt,v-eq}
T^{t,v}= \Phi + (t-1) [\psi_1\wedge\phi]_{[\pk]} + \left(t-\dfrac{1}{2} \right) [\phi\wedge\phi]_{[\pk]}
\end{equation}
and the $\nabla^t$-vertical tension field (of $f$) in the form:
\begin{equation}\label{tautv-eq}
\tau^{t,v}(f)=\tau^{0,v}(f) + t *[f^*\psi_1\wedge *(f^*\phi)]_{[\pk]}
\end{equation}
where $\psi_1\colon TN\to [\mk_1]$ is the projection on $[\mk_1]$ along $[\mk_*]\oplus [\pk]$ i.e. the $[\mk_1]$-component of $\psi$.
\begin{defn}
Let us suppose that $N=G/K$ admits a $G$-invariant almost complex structure $\undj$ which leaves invariant the decomposition $TN=\ver\oplus \hor$, that is to say the vector space $\mk$ admits an $\Ad K$-invariant almost complex structure $\undj_0$ leaving invariant the decomposition $\mk=\mk'\oplus\pk$. Then we will say that $\undj$ anticommutes with the reductivity term $[\psi,\phi]_{[\pk]}$ if 
$$
\undj[\psi,\phi]_{[\pk]}=-[\undj\psi,\phi]_{[\pk]}=-[\psi,\undj\phi]_{[\pk]}
$$
\end{defn}
If $\undj$ anticommutes with the reductivity term then $\mk_*$ is $\undj$-invariant so that it admits a $\undj$-invariant complement $\mk_1$ in $\mk'$.\\
We obtain immediately the following characterisation of the anticommutation of $\undj$ with the reductivity term.
\begin{prop}\label{prop-J-anticom-red}
Let us  that $N=G/K$ is endowed with a $G$-invariant almost complex structure $\undj$ leaving invariant the decomposition $TN=\hor\oplus\ver$. For any $\undj_0$-invariant $\Ad K$-invariant subspace $\mak l\subset\mk$, let us denote by $\mak l^{\pm}$ respectively the $\pm i$-eigenspace of $\undj_{|\mak l}$.
Then $\undj$ anticommutes with the reducivity term $[\psi,\phi]_{[\pk]}$ \iif
$$
\left[  {\mk'}^\pm,\pk^\pm\right]_\pk\subset \pk^\mp   \quad \text{and} \quad \left[ {\mk'}^\pm,\pk^\mp\right]_\pk\subset \{0\}.
$$
In particular, if $\mk_*=\left\{X\in\mk'|\, [X,\pk]_\pk=\{0\}\right\}$ admits a $\undj$-invariant $\Ad K$-invariant complement $\mk_1$, then these conditions are equivalent to 
$$
\left[  \mk_1^\pm,\pk^\pm\right]_\pk\subset \pk^\mp   \quad \text{and} \quad \left[ \mk_1^\pm,\pk^\mp\right]_\pk\subset \{0\}.
$$
\end{prop}
The following  theorem presents some relations between  vert. hol. harmonicity, vertical harmonicity, flatness and torsion freedom. 
\begin{thm}\label{vert-hol-harm-equiv-vert-harm}
Let us suppose that $N=G/K$ is endowed with a $G$-invariant almost complex structure $\undj$ leaving invariant the decomposition $TN=\hor\oplus\ver$ and which anticommutes with the reductivity term $[\psi,\phi]_{[\pk]}$. Let $\mk_1$ be a $\undj$-invariant $\Ad K$-invariant complement in $\mk'$ of $\mk_*=\left\{X\in\mk'|\, [X,\pk]_\pk=\{0\}\right\}$. Let $f\colon (L,j_L)\to N$ be a map from a Riemann surface into $N$, $F\colon L\to G$ a (local) lift of $f$ and $\alpha=F^{-1}.dF$.\\
$\bullet$  Then if $f$ is flat, $f^*\Phi=0$, and $[\mk_1]$-holomorphic then the following statements are equivalent:
\begin{description}
\item[(i)] $f$ is vert. hol. harmonic w.r.t. $\nabla^1$ and $\undj$ : $\left[ \bar\partial^{\nabla^{1,v}}\partial^v f\right]^{1,0}=0$.
\item[(ii)] $f$ is vert. hol. harmonic \wrt  $\nabla^0$ and $-\undj$ : $\left[ \bar\partial^{\nabla^{0,v}} \partial^v f\right]^{0,1}=0$. 
\end{description}
Moreover if $[\pk,\pk]_{\pk}=\{0\}$, then these are also equivalent to
\begin{description}
\item[(iii)] $f$ is vertically harmonic \wrt  $\nabla^1$.
\end{description}
$\bullet$ Furthermore, if $[\pk,\pk]_{\pk}=\{0\}$ and $f$ is flat, then  $f$ is $\nabla^1$-torsion free so that $\nabla^1$-vertical harmonicity is equivalent to strongly $\nabla^1$-vertical harmonicity.
\end{thm}
\textbf{Proof.}
The $\nabla^1$-vertical holomorphic harmonicity is written 
\begin{equation}\label{nabla1-vert-hol-harm}
f^*T^{1,v} + \undj \, d^{\nabla^{1,v}}\negthickspace *d^v f=0
\end{equation}
but $f^*T^{1,v}=\dfrac{1}{2}f^*[\phi\wedge\phi]_{[\pk]}$ and $f^*T^{0,v}=-[f^*\psi\wedge f^*\phi]_{[\pk]} - \dfrac{1}{2}f^*[\phi\wedge\phi]_{[\pk]}$ whereas 
$$
d^{\nabla^{1,v}}*d^v f=d^{\nabla^{0,v}}*d^v f + [f^*\psi\wedge (*f^*\phi)]_{[\pk]}
$$
so that 
\begin{equation}\label{f*(1,v)}
f^*T^{1,v} + \,\undj\, d^{\nabla^{1,v}}\negthickspace *d^v f= \dfrac{1}{2}[f^*\phi\wedge f^*\phi]_{[\pk]}\, + \, \undj\, d^{\nabla^{0,v}}\negthickspace *d^v f \, + \,\undj[f^*\psi\wedge (*f^*\phi)]_{[\pk]}.
\end{equation}
Now let us use the fact that $\undj$ anticommutes with $[\psi,\phi]_{[\pk]}$:
\begin{eqnarray*}
\undj\left[f^*\psi\wedge *(f^*\phi)\right]_{[\pk]} & = & -\left[\undj (f^*\psi_1) \wedge *(f^*\phi)\right]_{[\pk]} \\
& = & \left [*(f^*\psi_1)\wedge *(f^*\phi)\right]  \quad \text{because } f \text{ is } [\mk_1]\text{-holomorphic},\\
& =  & \left[f^*\psi_1 \wedge f^*\phi\right]= \left[f^*\psi \wedge f^*\phi\right].
\end{eqnarray*}
Therefore, injecting this in (\ref{f*(1,v)}), we obtain
$$
f^*T^{1,v} + \,\undj\, d^{\nabla^{1,v}}\negthickspace *d^v f= -\left( f^*T^{0,v} - \,\undj\, d^{\nabla^{0,v}}\negthickspace *d^v f\right). 
$$
This proves the equivalence (i) $\Leftrightarrow$ (ii). Now, if we suppose that $[\pk,\pk]_{\pk}=\{0\}$, then $f^* T^{1,v}=0$. Therefore the $\nabla^1$-vertical holomorphic harmonicity (\ref{nabla1-vert-hol-harm}) is equivalent to the $\nabla^1$-vertical harmonicity $d^{\nabla^{1,v}} *d^v f=0$. This completes the proof.\hfill $\square$ \medskip\\
Now, let us see how the vertical holomorphic harmonicity is written in terms of the Maurer-Cartan form $\alpha$ of a  lift $F$ of $f\colon L\to N$. We obtain immediately
\begin{prop}\label{vert-hol-harm-MC}
Let us suppose that $N=G/K$ is endowed with a $G$-invariant almost complex structure $\undj$ leaving invariant the decomposition $TN=\hor\oplus\ver$. Then $f\colon L\to N$ is vert. hol. harmonic \wrt $\nabla^0$ and $-\undj$ \iif 
$$
\bar\partial \alpha_{\pk^+}' + \left[\alpha_\kk'' \wedge \alpha_{\pk^+}'\right]=0.
$$
Moreover, if $\undj$ anticommutes with the reductivity term $[\psi,\phi]_{[\pk]}$ and $\mk_*$ admits an  $\Ad K$-invariant $\undj$-invariant complement $\mk_1$ in $\mk'$, then a $[\mk_1]$-holomorphic map $f\colon L\to N$ is vert. hol. harmonic \wrt $\nabla^1$ and $\undj$ \iif 
$$
\bar\partial \alpha_{\pk^+}' + \left[\alpha_\kk'' \wedge \alpha_{\pk^+}'\right] + \left[\alpha_{\mk_1^-}'' \wedge \alpha_{\pk^-}'\right]_{\pk} + \left[ \alpha_\pk ''\wedge \alpha_\pk'\right]_{\pk^+}=0 .  
$$  
\end{prop}
\index{determined minimal@determined, minimal|(}
Now, let us suppose that $N=G/K$ is a (locally) $(2k+1)$-symmetric space, then the $\Ad K$-invariant decomposition $\mk=\mk'\oplus\pk$ is given by $\pk=\mk_k$ and $\mk'=\oplus_{j+1}^{k-1} \mk_j$ with the notations of \ref{2.1.2}. Moreover according to the commutation relations $[\g_i^\C,\g_j^\C]\subset\g_{i+j}^\C$, we have 
$$
\mk_*=\left\{ X\in \mk', [X,\pk]_\pk=\{0\}\right\} =\oplus_{j=2}^{k-1} \mk_j
$$
so that $\mk_1$ is an $\Ad G_0$-invariant supplement to $\mk_*$. Moreover $N$ is endowed naturally with its canonical almost complex structure $\undj$ defined in \ref{2.1.2}, which leaves invariant all the $\mk_j$ and thus the subspaces $\mk_1$, $\mk_*$, $\pk$. Furthermore, using once again the commutation relations, one can see that $\undj$ anticommutes with the reductivity term $[\psi,\phi]_{[\pk]}$. Finally, let us remark that $[\pk,\pk]_\pk =[\mk_k,\mk_k]_\pk=\{0\}$. Now, the theorem~\ref{vert-hol-harm-equiv-vert-harm} can be applied.
\begin{cory}\label{7.5}
Let us suppose that $N=G/K$ is a (locally) $(2k+1)$-symmetric space endowed with its canonical almost complex structure $\undj$, and with the $\undj$-invariant splittings $TN=[\mk']\oplus [\pk]$ and $[\mk']=[\mk_1]\oplus [\mk_*]$. Let $ f\colon L\to N$ be a map, $F\colon L\to G$ a lift of $f$ and $\alpha=F^{-1}.dF$. Then if $f$ is flat, $f^*\Phi=0$, and $[\mk_1]$-holomorphic then $f$ is $\nabla^1$-vertically torsion free $f^*T^{1,v}=0$, and the following statements are equivalent
\begin{description}
\item[(i)] $f$ is vert. hol. harmonic \wrt $\nabla^1$ and $\undj$ : $\left[ \bar\partial^{\nabla^{1,v}}\partial^v f\right]^{1,0}=0$.
\item[(ii)] $f$ is vert. hol. harmonic \wrt  $\nabla^0$ and $-\undj$ : $\left[ \bar\partial^{\nabla^{0,v}} \partial^v f\right]^{0,1}=0$. 
\item[(iii)] $f$ is vertically harmonic \wrt  $\nabla^1$.
\item[(iv)] $f$ is strongly vertically harmonic \wrt  $\nabla^1$.
\end{description}
\end{cory}
Now, let us apply proposition~\ref{vert-hol-harm-MC}.
\begin{prop}\label{MC-vert-hol-harm}
Let us suppose that $N=G/K$ is a (locally) $(2k+1)$-symmetric space. Then $f\colon L\to N$ is vert. hol. harmonic \wrt $\nabla^0$ and $-\undj$ \iif 
$$
\bar\partial \alpha_{k}' + \left[\alpha_0'' \wedge \alpha_{k}'\right]=0.
$$
Moreover, if $f\colon L\to N$ is flat and $[\mk_1]$-holomorphic, then it is vert. hol. harmonic \wrt $\nabla^1$ and $\undj$ \iif
$$
\bar\partial\alpha_{-k}' + [\alpha_0''\wedge \alpha_{-k}'] + [\alpha_1''\wedge \alpha_k']=0.
$$  
\end{prop}
Furthermore, as in the even case (i.e. $N=G/K$ is (locally) $2k$-symmetric), the horizontal holomorphicity implies the flatness.
\begin{prop}\label{horhol--flat}
Let us suppose that $N=G/K$ is a (locally) $(2k+1)$-symmetric space. Then if $f\colon L\to N$ is horizontallly holomorphic (i.e. $[\mk']$-holomorphic) then $f$ is flat $f^*\Phi=0$.
\end{prop}
\proof It follows from the equation 
$$
\tl\Phi = -[\theta_{\mk'},\theta_{\mk'}]_\pk = -\sum_{\underset{1\leq |i|,|j|\leq k-1}{i+j=k}} [\theta_i,\theta_j]_{\g_k},
$$
and the fact that if $i+j=k$ and $1\leq |i|,|j|\leq k-1$, then $i$ and $j$ have the same sign. \comprf\hsq\medskip\\
Now let us conclude with the following geometric interpretation of the odd minimal determined system.
\begin{cory}\label{cory-2k+1-syst}
Let us suppose that $N=G/K$ is a (locally) $(2k+1)$-symmetric space. Then the odd minimal determined system $(\syst(k+1,\tau))$ associated to $N$ means that the geometric map $f\colon L\to N$ is horizontally holomorphic and vertically harmonic \wrt the linear connection $\nabla^1$. Moreover the horizontal holomorphicity implies the flatness of $f$ and its freedom from $\nabla^1$-vertical torsion,  $f^*T^{1,v}=0$.\\
More pecisely, the (last) equation $(S_{k+1})$ of the system (which lies in $\g_k$) means the vert. hol. harmonicity of $f$ \wrt $\nabla^0$ and $-\undj$
$$
\left[ \bar\partial^{\nabla^{0,v}} \partial^v f\right]^{0,1}=0
$$
whereas the equation $(S_k)$ (which lies in $\g_{-k}$) means  the vert. hol. harmonicity of $f$ \wrt  $\nabla^1$ and $\undj$  
$$
\left[ \bar\partial^{\nabla^{1,v}}\partial^v f\right]^{1,0}=0
$$
Moreover the sums $(S_k)+ (S_{k+1})$ (which lies in $\mk_k$) means (taking account of the $[\mk_1]$-holomorphicity $\alpha_{-1}''=0$) the strongly vertical harmonicity of $f$ \wrt $\nabla^1$:
$$
\bar\partial^{\nabla^{1,v}}\partial^v f=0,
$$ 
so that its real part means that $f$ is $\nabla^1$-vertically torsison free and its imaginary part that $f$ is $\nabla^1$-vertically harmonic.\\
All the other equations of the system, $(S_j)$, $0\leq j \leq k-1$ are (after having taken account of the horizontal holomorphicity $\alpha_{-j}''=0$, $1\leq j\leq k-1$) nothing but the projections on the subspace $\g_{-j}$, $1\leq j\leq k-1$ of the Maurer-Cartan equation (which means the existence of the geometric map $f$ corresponding to $\alpha$).
\end{cory}
\proof The first assertion has been  proved in section~\ref{detodcase}, theorem~\ref{conclusion-odd-mindeter}. Moreover, the $\nabla^1$-vertical harmonicity will follows from the next assertions, according to corollary~\ref
{7.5}.\\
The second assertion follows from proposition~\ref{horhol--flat} and theorem~\ref{vert-hol-harm-equiv-vert-harm}.
The third assertion follows from propositions \ref{prop-odd-mindeter} and \ref{MC-vert-hol-harm}. The fourth assertion has been proved in subsection~\ref{detodcase}, paragraph \emph{The strictly minimal determined case}. Finally, the fifth assertion follows from proposition~\ref{horhol--flat}. \comprf\hsq

\paragraph{Strongly vertical harmonicity \wrt $\nabla^t$.} Let us see what the strongly vertical harmonicity \wrt to $\nabla^t$, with $t\in [0,1]\setminus \{1\}$, means.\medskip\\
We have seen that the tension field of a map $f\colon L\to N=G/K$, with respect to $\nabla^t$, does not depend on $t\in [0,1]$ (see theorem~\ref{thmharmonic}). Let us set $\tau(f):=\tau^t(f)$.
\begin{prop}\label{prop-7.5}
Let $f\colon L\to N=G/K$ be a map.\\ 
$\bullet$ Then we have
$$
\tau^{0,v}(f)=[\tau(f)]^v=2\left[\bar\partial^{\nabla^{\frac{1}{2}}}\partial f\right]^v.
$$
$\bullet$  If $[\pk,\pk]_\pk=\{0\}$ and $f$ is flat, then $f^*T^{1,v}=0$ i.e. $\tau^{1,v}(f)=2\,\bar\partial^{\nabla^{1,v}}\partial f$.\medskip\\
$\bullet$ Moreover, under these hypothesis\footnote{i.e. the hypothesis $[\pk,\pk]=\{0\}$ and $f$ is flat.}, let  $\mk_1$ be an  $\Ad K$-invariant complement of $\mk_*$, in $\mk'$. Then,   the following statements are equivalent:
\begin{description}
\item[(i)] $f^*T^{t,v}=0$ for one $t\in [0,1]\setminus\{1\}$, 
\item[(ii)] $f^*[\psi_1\wedge \phi]_{[\pk]}=0$,
\item[(iii)] $\tau^{t,v}(f)=\tau^{0,v}(f)$ for one $t\in [0,1]\setminus\{0\}$.
\end{description}
\end{prop}
\proof In the first point, the first equality comes from the fact that $\nabla^0$ leaves invariant the splitting $TN=\ver\oplus\hor$, since then $\left( \nabla^0\right)^v d^v f = \left[ \nabla^0 df\right]^v$. The second equality follows from theorem~\ref{thmstrong} and the fact that $T^{ \frac{1}{2} }=0$.\\
The second point follows from theorem~\ref{vert-hol-harm-equiv-vert-harm}. The third point follows from the equations \eqref{Tt,v-eq} and \eqref{tautv-eq}. \comprf\hsq
\begin{cory}\label{cory-7.7}
Let us suppose that $N=G/K$ is endowed with a $G$-invariant almost complex structure $\undj$ leaving invariant the decomposition $TN=\hor\oplus\ver$ and which anticommutes with the reductivity term $[\psi,\phi]_{[\pk]}$. We also suppose that there exists   a $\undj$-invariant $\Ad K$-invariant complement $\mk_1$ in $\mk'$ of $\mk_*$.
Lastly, we suppose that $[\pk,\pk]_\pk=\{0\}$.\\
 Let $f\colon L\to N=G/K$ be a map which is flat. Let $F\colon L\to N$ be a lift of $f$ and $\alpha=F^{-1}.dF$.
Then the following statements are equivalent
\begin{description}
\item[(i)] $\bar\partial^{\nabla^{t,v}}\partial^v f=0$ for one $t\in [0,1]\setminus \{1\}$,
\item[(ii)] $\left[\bar\partial^{\nabla^t}\partial f\right]^v=0$ for one $t\in [0,1]\setminus\{\frac{1}{2}\}$,
\item[(iii)] $\tau^{1,v}(f)=0$ and $f^*[\psi_1\wedge \phi]_{[\pk]}=0$,
\item [(iv)] $[\tau(f)]^v=0$ and $f^*[\psi_1\wedge \phi]_{[\pk]}=0$,
\item[(v)]  $\left[ \bar\partial^{\nabla^{t,v}}\partial^v f\right] ^{1,0}=0$ for one $t\in [0,1]$ and $f^*[\psi_1\wedge \phi]_{[\pk]}=0$,
\item[(vi)] $\left[ \bar\partial^{\nabla^{1-t,v}}\partial^v f\right] ^{0,1}=0$ for one $t\in [0,1]$ and $f^*[\psi_1\wedge \phi]_{[\pk]}=0$,
\item[(vii)] $f$ is a geometric solution of the \emph{first elliptic integrable system associated to the $\ad \kk$-invariant decomposition} $\g=(\kk\oplus\mk')\oplus \pk$, i.e. the 1-form $\beta_\lm=\lm^{-1}\alpha_\pk' + (\alpha_\kk+ \alpha_{\mk'}) + \lm\alpha_\pk''$ satisfies the zero curvature equation
$$
d\beta_\lm + \dfrac{1}{2}[\beta_\lm\wedge \beta_\lm]=0 \quad \forall \lm\in \C^*.
$$
\end{description}
\end{cory}
\proof The statement (i) is equivalent to the statement:  " $f^*[\psi_1\wedge \phi]_{[\pk]}=0$ and $\tau^{t,v}(f)=0$, $\forall t\in [0,1]$", according to proposition~\ref{prop-7.5} and theorem~\ref{strongly-vertical-harm}. Moreover, since $\widetilde{ \left[ T^t\right]^v}= (2t-1)[\theta_\mk,\theta_\mk]_\pk =(2t-1)\left( [\theta_{\mk'}\wedge\theta_{\mk'}]_\pk  + \frac{1}{2}[\theta_{\mk_1}\wedge\theta_\pk]_\pk\right) $, we have 
$$
\left[ T^t\right]^v =(2t-1) \left( \Phi +  [\psi_1\wedge\phi]_{[\pk]} \right).
$$ 
Therefore, (ii) $\Leftrightarrow$  $(\, [\psi_1\wedge \phi]_{[\pk]}=0$ and $[\tau(f)]^v=0 \,)$, which  is nothing but (iv). Moreover, according to theorem~\ref{thm-vert-hol-harm-charact} and equations \eqref{Tt,v-eq} and \eqref{tautv-eq}, we have (v) $\Leftrightarrow$ (vi) $\Leftrightarrow$ ( $f^*[\psi_1\wedge \phi]_{[\pk]}=0$ and $\tau^{t,v}(f)=0$, $\forall t\in [0,1]$). Moreover, proposition~\ref{prop-7.5} gives us the equivalences 
$$
(iii) \Leftrightarrow (iv) \Leftrightarrow (\, f^*[\psi_1\wedge \phi]_{[\pk]}=0 \text{ and  } \tau^{t,v}(f)=0,\  \forall t\in [0,1] \,). 
$$
We have then proved the equivalence between the six first assertions.
Finally, we have
\begin{equation}\label{system-(vii)}
(vii)  \Leftrightarrow \left\lbrace 
\begin{array}{l} \begin{array}{rcl} \bar\partial \alpha_\pk' + [\alpha_\kk\wedge\alpha_\pk'] & = & 0\\
 {} [\alpha_{\mk'}\wedge\alpha_\pk'] & = &  0 \end{array}\\
d(\alpha_\kk + \alpha_{\mk'})  + [\alpha_\pk'\wedge\alpha_\pk''] + \dfrac{1}{2}[(\alpha_\kk + \alpha_{\mk'})\wedge (\alpha_\kk + \alpha_{\mk'})]=0
\end{array}\right.
\end{equation}
Moreover, we have $[\alpha_\pk'\wedge\alpha_\pk'']_\pk=0$, since $[\pk,\pk]_\pk=\{0\}$. Furthermore, $[(\alpha_\kk + \alpha_{\mk'})\wedge (\alpha_\kk + \alpha_{\mk'})]_\pk=0$ because $[\mk',\pk]_\pk=0$ and the flatness of $f$ means that $[\alpha_{\mk'}\wedge\alpha_{\mk'}]_\pk=0$.  Therefore, the last equation in \eqref{system-(vii)} is nothing but $[\mrm{MC}]_{\kk\oplus\mk'}$, the projection of the Maurer-Cartan equation on $\kk\oplus\mk'$. Hence, we conclude that
$$
(vii) \Leftrightarrow  (\, \tau^{0,v}(f)=0 \text{ and  }  f^*[\psi_1\wedge \phi]_{[\pk]}=0  \,)
$$
\comprf\hsq\medskip\\
Now, we come back to the case of a (locally) $(2k+1)$-symmetric space.
\begin{cory}\label{cory-2k+1}
Let us suppose that $N=G/K$ is a (locally) $(2k+1)$-symmetric space. Let $f\colon L\to N=G/K$ be a map,  $F\colon L\to N$ be a lift of $f$ and $\alpha=F^{-1}.dF$. Then the following statements are equivalent:
\begin{description}
\item[(i)] $f$ is horizontally holomorphic and strongly $\nabla^t$-vertically harmonic for one $t\in [0,1]\setminus \{1\}$,
\item[(ii)] $f$ is a geometric solution of $(\syst(k+1,\tau))$ and $[\alpha_1\wedge\alpha_k]=0$,
\item[(iii)] $f$ is horizontally holomorphic and is a geometric solution of the first elliptic integrable system associated to the $\ad \g_0$-invariant decomposition $\g=(\g_0\oplus\mk')\oplus \mk_k$.
\end{description}
\end{cory}
\proof According to corollary~\ref{cory-7.7} and corollary~\ref{cory-2k+1-syst}, we have (i) $\Leftrightarrow$ (iii). Moreover, according to corollaries~\ref{cory-2k+1-syst} and \ref{cory-7.7}, the equivalent statements (i) and (iii) are also equivalent to (ii). \hsq
\index{vertically harmonic|)}\index{determined minimal@determined, minimal|)}\index{vertical tension field|)}
\index{strongly vertically harmonic|)}\index{vertically holomorphically harmonic|)}
%
%
%%%%%%%%%%%%%%%%%%%%%%%%%%%%%%%%%%%%%%%%%%%%%%%%%%%%%%%%%%%%%%%%%%%%%%%%%%%%%%%%%%%%%%%%%%%%%%%%%
%
%
%
\subsection{The intermediate determined systems.}\label{subsec-hom-intermediate}
\subsubsection{The odd case.}\label{subsec-intermediate-oddcase}
\index{odd case|(}
Here, we consider  a (locally) $(2k+1)$-symmetric space, $N=G/K$ , endowed with its canonical almost complex structure $\undj$ and its canonical connection $\nabla^0$. Let $ k+1 \leq m\leq 2k$ be an integer and set $\und m= 2k - m$. Then we consider the $G$-invariant splitting $TN= \hor^{ m}\oplus \ver^{ m}$, where
$$
\hor^{ m}=\oplus_{j=1}^{\und m} [\mk_j] \quad \text{and} \quad  \ver^{ m}=\oplus_{j=\und m + 1}^{k} [\mk_j].
$$
Let us remark that this splitting is also $\undj$-invariant and $\nabla^0$-parallel. \\
According to proposition~\ref{hol-hor-condition} and theorem~\ref{oddmaxdeter-geom-interpret}, we have the following.
\begin{prop}\label{prop-intermediate-stringy}
Let  $N=G/K$ be a (locally) $(2k+1)$-symmetric space endowed with its canonical almost complex structure $\undj$ and its canonical connection $\nabla^0$.  Let $m$ be an integer such that $ k+1 \leq m\leq 2k$. Then the associated  determined system, $\syst(m,\tau)$  means that the geometric map $f\colon L\to N$ is $\star$-stringy harmonic and 
$\hor^m$-holomorphic. 
\end{prop}
We will need the following definition.
\begin{defn}
Let $(N,J)$ be an almost complex  manifold endowed with a linear connection $\nabla$. Let $E\subset TN$ be some subbundle. We will say that  a map $f\colon L\to N$ from a Riemann surface into $N$ is  $E$-$\star$-stringy harmonic if 
$$
\left[-\tau_g(f) + (J\star T)_g(f)\right]^{E} =0.
$$
where $g$ is a Hermitian metric on $L$. If $E$ inherits the name of vertical or horizontal subbundle then we will say that $f$ is vertically, resp. horizontally $\star$-stringy harmonic. 
\end{defn}  
Then the proposition~\ref{prop-intermediate-syst} gives us:
\begin{thm}\label{thm-intermed-deter-stringy}
Let  $N=G/K$ be a (locally) $(2k+1)$-symmetric space endowed with its canonical almost complex structure $\undj$ and its canonical connection $\nabla^0$.  Let $m$ be an integer such that $ k+1 \leq m\leq 2k$. Then the associated  determined system, $\syst(m,\tau)$  means that the geometric map $f\colon L\to N$ is $\hor^m$-holomorphic and vertically $\star$-stringy harmonic. 
\end{thm}
\proof Let us come back to section~\ref{eq-maximal-deter}, then the equation $(E_\mk)$, in theorem~\ref{thm-max-deter-odd},  means that $f$ is $\star$-stringy harmonic, according to theorem~\ref{oddmaxdeter-geom-interpret} and its proof. Therefore, the projection $[(E)]_{\pk}$ in proposition~\ref{prop-intermediate-syst} means that $f$ is vertically  $\star$-stringy harmonic. Then the proposition~\ref{prop-intermediate-syst} allows us to conclude. \comprf \hsq\medskip\\
In fact this theorem is a particular case of a more general result:
\begin{prop}
Let $(N,J,\nabla)$ be an almost complex affine manifold, with $TN=\ovr\hor\oplus\ovr\ver$ a $J$-invariant, $\nabla$-parallel splitting. Let  $T$ be the torsion of $\nabla$. Let us suppose that 
\begin{equation}\label{hyp-T++}
T^{++} (\ovr\ver,\ovr\ver)\subset \ovr\ver \quad \text{and} \quad  T^{++}(\ovr\hor,\ovr\ver)\subset\ovr\ver.
\end{equation}
Let $(L,j)$ be a Riemann surface. Then any $\ovr\hor$-holomorphic map $f\colon L\to N$ is horizontally $\star$-stringy harmonic. Therefore a horizontally holomorphic map $f\colon L\to N$  is  $\star$-stringy harmonic \iif it is vertically $\star$-stringy harmonic.
\end{prop}
\proof
Using the fact that the splitting is $\nabla$-parallel and the horizontal holomorphicity, we have
$$
[\tau^\nabla(f)]^{\ovr\hor} = d^\nabla *[df]^{\ovr\hor}= - J d^\nabla [df]^{\ovr\hor}= - J T^{\ovr\hor}(df,df).
$$
Moreover, recall that $J\star T=J(T^{++} - T^{--} - T^{1,1})$. We then have to prove that $f^* (T^{++})^{\ovr\hor}=0$, which will imply the horizontal  $\star$-stringy harmonicity of $f$. Moreover, remark that we have 
$$
T^{\eps,\eps}(T^{1,0} N, T^{0,1} N) = 0, \quad \forall \eps\in\Z_2.
$$
Therefore, since $f$ is horizontally holomorphic, we have $T^{\eps,\eps}([df]^{\ovr\hor},[df]^{\ovr\hor})=0$. Furthermore, according to the hypothesis \eqref{hyp-T++} on $T^{++}$, we conclude that $f^* (T^{++})^{\ovr\hor}=0$. \comprf \hsq\medskip\\
It is not difficult to see that a (locally) $(2k+1)$-symmetric space endowed with its canonical almost complex structure $\undj$, its canonical connection $\nabla^0$ and the splitting $TN=\hor^m\oplus\ver^m$, satisfies the hypothesis of the previous proposition.
\index{odd case|)}

%
%%%%%%%%%%%%%%%%%%%%%%%%%%%%%%%%%%%%%%%%%%%%%%%%%%%%%%%%%%%%%%%%%%%%%%%%%%%%%%%%
%
\subsubsection{The even case.}\index{f structure@$f$-structure}\index{even case|(}
Here, we consider  a (locally) $2k$-symmetric space, $N=G/K$ , endowed with its canonical $f$-structure $F$ and its canonical connection $\nabla^0$. Let $ k \leq m\leq 2k-1$ be an integer and set $\und m= 2k -1- m$. Then we consider the $G$-invariant splitting $TN= \hor^{ m}\oplus \ver^{ m}$, where
$$
\hor^{ m}=\oplus_{j=1}^{\und m} [\mk_j] \quad \text{and} \quad  \ver^{ m}=\oplus_{j=\und m + 1}^{k} [\mk_j],
$$
and $\mk_k=\g_k$. This splitting is also $\undj$-invariant and $\nabla^0$-parallel. \\
According to proposition~\ref{hol-hor-condition} and theorem~\ref{thm-max-deter-even-stringy}, we have the following.
\begin{prop}
Let  $N=G/K$ be a (locally) $2k$-symmetric space endowed with its canonical $f$-structure $F$ and its canonical connection $\nabla^0$.  Let $m$ be an integer such that $ k \leq m\leq 2k-1$. Then the associated  determined system, $\syst(m,\tau)$  means that the geometric map $f\colon L\to N$ is $\star$-stringy harmonic and 
$\hor^m$-holomorphic. 
\end{prop}
We will need the following definition.
\begin{defn}
Let $(N,F)$ be an $f$-manifold endowed with a linear connection $\nabla$, and $E\subset TN$ some subbundle.
We will say that  a map $f\colon L\to N$ from a Riemann surface into $N$ is  $E$-$\star$-stringy harmonic if 
$$
\left[-\tau_g(f) + (F\star T)_g(f)\right]^{E} =0.
$$
where $g$ is a Hermitian metric on $L$. If $E$ inherits the name of vertical or horizontal subbundle then we will say that $f$ is vertically, resp. horizontally $\star$-stringy harmonic. 
\end{defn}  
Then the proposition~\ref{prop-intermediate-syst} gives us:
\begin{thm}
Let  $N=G/K$ be a (locally) $2k$-symmetric space endowed with its canonical $f$-structure $F$ and its canonical connection $\nabla^0$.  Let $m$ be an integer such that $ k \leq m\leq 2k-1$. Then the associated  determined system, $\syst(m,\tau)$  means that the geometric map $f\colon L\to N$ is $\hor^m$-holomorphic and vertically $\star$-stringy harmonic. 
\end{thm}
\proof It is analogous to the one of theorem~\ref{thm-intermed-deter-stringy}. \hsq\medskip\\
In fact this theorem is a particular case of a more general result:
\begin{prop}
Let $(N,F,\nabla)$ be an $f$-manifold endowed with an $f$-connection. As usual, we set $\hor=\im F$ and $\ver=\ker F$. Suppose  that   $\mR_{\ver}=0$. Moreover, let  $TN=\ovr\hor\oplus\ovr\ver$ be an $F$-invariant, $\nabla$-parallel splitting, such that $\ovr\hor\subset\hor$ and $\ovr\ver\supset\ver$.  Let us set $\ovr T:=T_{\hor^2}^{\hor}$, where $T$ is the torsion of $\nabla$. Let us suppose that 
$$
\ovr T^{++} (\ovr\ver\cap \hor,\ovr\ver\cap \hor)\subset \ovr\ver\cap \hor \quad \text{and} \quad  \ovr T^{++}(\ovr\hor,\ovr\ver\cap \hor)\subset\ovr\ver\cap \hor.
$$
Suppose also that $T_{\hor\times\ver}^\hor$ is pure\footnote{i.e. $T_{\hor\times\ver}^\hor$ anticommutes with $\Bar J$, see definition~\ref{defn-pure}.}.\\
Let $(L,j)$ be a Riemann surface. Then any $\ovr\hor$-holomorphic map $f\colon L\to N$ is horizontally $\star$-stringy harmonic. Therefore a horizontally holomorphic map $f\colon L\to N$  is  $\star$-stringy harmonic \iif it is vertically $\star$-stringy harmonic.
\end{prop}
\proof It is analogous to the one of theorem~\ref{thm-intermed-deter-stringy}. \hsq\medskip\\
It is not difficult to see that a (locally) $2k$-symmetric space endowed with its canonical $f$-structure $F$, its canonical connection $\nabla^0$ and the splitting $TN=\hor^m\oplus\ver^m$, satisfies the hypothesis of the previous proposition. 
\index{even case|)}

% Let us set $\ovr T:=T_{\ovr\hor^2}^{\ovr\hor}$, where $T$ is the torsion of $\nabla$. Let us suppose that $ \ovr %T^{++} (\ovr\ver,\ovr\ver)\subset \ovr\ver$ and $\ovr T^{++}(\ovr\hor,\ovr\ver)\subset\ovr\ver$.
%
%
%

\subsubsection{Sigma model with a Wess-Zumino term.}
We have seen that the maximal determined system has an interpretation in terms of a sigma models with a Wess-Zumino term defined by a 3-form $H$. In fact, more generally, let $m_{k'}\leq m \leq k' -1$,  and let us consider the splitting $TN=\hor^m \oplus\ver^m$ defined above. Then one can prove that  any $m$-th determined system is the Euler-Lagrange equation \wrt vertical variations (i.e. in $\ver^m$) of the following functional
$$
E^{\bar v}(f)= \dfrac{1}{2}\int_L |d^{\bar v } f|^2 d\mrm{vol}_g  + \int_B H^{\bar v }
$$
where $d^{\bar v } f= [df]^{\ver^m}$, $H^{\bar v }=H-\ovr H= H_{|\mal S(\ver^m,\hor^m)}$, $\ovr H=H_{|(\hor^m)^3}$, and $B$ is a 3-submanifold of $N$ with boundary $\partial B= f(L)$. \\
We will come back to this in \cite{ki6}.

%%%%%%
%%%%%%

\subsection{Some remarks about the twistorial interpretation.}
\subsubsection{The even case.}\index{even case|(}
We have seen that in the even case, each geometric property of the geometric map $f\colon L\to N=G/G_0$ like horizontal holomorphicity,  vertical harmonicity, stringy harmonicity is naturally translated into the same geometric property of the twistor lift $J\colon L\to \mZ_{2k,2}^\alpha(M,J_2)$  in the twistor space $\mZ_{2k,2}^\alpha(M,J_2)$. Moreover, the very particular structure of homogeneous fibre $f$-bundle of $N=G/G_0$ can be realised as a subbundle of the universal homogeneous fibre $f$-bundle $\mZ_{2k,2}^\alpha(M,J_2)$. "Universal" means that we can define it for any Riemannian manifold $(M,g)$ endowed with a global $k$-structure $J_2$ (and a metric connection $\ovr\nabla$).  Remark that since $\mZ_{2k,2}^\alpha(M,J_2)$ is a complete reduction of $\zdk(M)$ (\wrt to the structures of  homogeneous fibre $f$-bundles defined by the canonical connection of $M$) then we can also as well embedd $N$ in  $\zdk(M)$ which is more universal because defined for any Riemannian manifold $(M,g)$ (endowed with a metric connection $\ovr\nabla$). Then all the geometric properties below are also preserved under this embedding and hold in $\zdk(M)$. However, it is perhaps better to use the universal homogeneous fibre $f$-bundle containing $N$ which looks like the most to $N$ and which is also the smallest one (a kind of "universal homogeneous fibre $f$-bundle closure"). This exists and this is $\mZ_{2k,2}^\alpha(M,J_2)$. For example, $\mZ_{2k,2}^\alpha(M,J_2)$ has a symmetric fibre like $N$, whereas $\zdk(M)$ has a $2k$-symmetric fibre.
\index{even case|)} 

\subsubsection{The odd case.}\index{odd case}
In the odd case the use of the twistor space $\mZ_{2k+1}(N)$ is less pertinent than in the even case. Indeed in the even case, we had some particular fibration that twistor space allows to realise more universally as some bundle of endormophisms over $M$. Here we do not have this problem of fibration and therefore do not need a priori the twistor space. In the odd case, we have a canonical section\footnote{Defined by lemma~\ref{lemma-def-J_1}} $J_1\colon G/G_0\to \mZ_{2k+1}(G/G_0)$, which allows to duplicate each geometric property satisfied by the geometric map $f\colon L\to N$ into 2 "identical" properties in each subbundle $\hor$ and $\ver$ of the tangent bundle of the twistor space. \\
However, the twistor space $\mZ_{2k+1}(N)$ is still universal since it can be defined for any Riemannian manifold $(N,h)$ (endowed with some metric connection $\ovr\nabla$). Moreover, under the embedding defined by the canonical section $J_1\colon G/G_0\to \mZ_{2k+1}(G/G_0)$, any geometric property like $\undj$-holomorphicity, $F^{[m]}$-holomorphicity, vertical harmonicity, stringy harmonicity and so on is preserved and holds for the twistor lift $J=f^*J_1 \colon L\to \mZ_{2k+1}(N)$  in the twistor space.\medskip\\ 
We have to specify that the structure of homogeneous fibre bundle of $\mZ_{2k+1}(N)$ is defined with respect to\footnote{Recall that $\nabla^0$ is metric \wrt any $G$-invariant metric, whereas $\nabla^t$, for $t\neq 0$ is metric when $N$ is endowed with a naturally reductive metric} $\nabla^t$. In particular, this connection defines a splitting $TN=\hor\oplus\ver$. Then, we can define a canonical almost complex structure $\check J:=\left((d\pi)^* \undj\right)_{|\hor} \oplus \pi^*\und{(\Ad J_1)}$. That is to say on the horizontal subbundle we take complex structure defined by (the lift of) $\undj$ and moreover the fibre  is $(2k+1)$-symmetric and is therefore endowed with a canonical almost complex structure, that we endow  the vertical subbundle with.\\
Then, with respect  to this almost complex structure, $J_1\colon (G/G_0,\undj)\to (\mZ_{2k+1}(G/G_0), \check J)$ is holomorphic.\medskip\\
We will come back to that in \cite{ki6}.
\index{canonical!connection, $G$-invariant|)}

\subsection{Bibliographical remarks and summary of the results.}
A summary of the results has already been given in the introduction of the present section~\ref{affineharmhom}. All the results listed their are new. \\
Let us mention some related works about harmonic maps into homogeneous space: Higaki \cite{higaki}, Burstall-Pedit \cite{12}, Black \cite{Black}.  In \cite{higaki}, strongly harmonic maps are defined (under a different name) but nobody has remarked that these has a formulation in terms of a zero curvature equation (see remark~\ref{nobody-strongly}).

%

%%%%%
%%%%%%%%%
%%%
\section{Appendix}
\subsection{Vertical harmonicity}
\begin{thm}\label{same-vertical-Harm}
Let us consider the situation described by example~\ref{example1} and suppose that $\pi\colon N\to M$ is a Riemannian submersion and $u\colon L\to M$ is an isometry. Then $f\colon L\to N$ is vertically harmonic \iif the corresponding section $\tl f\colon L\to u^*N$ is a harmonic section. Furthermore $f\colon L \to N$ is harmonic \iif $\tl f\colon L\to u^*N$ is harmonic and $[\tau(f)]_{u_*(TL)^\perp}^\hor=0$ i.e. the component of the tension field in the subspace of $\hor$ corresponding by the isometry $d\pi_{|\hor}$ to the normal bundle $u_*(TL)^\perp$ in $TM$, vanishes, or equivalently $\left[ d\pi(\tau(f))\right]_{u_*(TL)^\perp}=0$.
\end{thm}
\textbf{Proof.} The Levi-Civita in $u^*N$ is the orthogonal projection ofthe Levi-Civita connection in $L\times N$, on the tangent bundle $T(u^*N)$. Let us determine this orthonormal projection. First let us express clearly what is the tangent subbundle $T(u^*N)$ in $T(L\times  N)$.
$$
T_{(l,n)} (u^*N)=\{(\xi,\eta)\in T_{(l,n)}L\times N| du(\xi)=d\pi(\eta)\}. 
$$
Let us do some identifications. First an usual one: consider that $TL$ is a subbundle of $TM_{|L}$ (and forget the "$u_*$" in $u_*(TL)$), secondly: we consider that $\pi^*TM=\hor$, identifying these by the isometry $d\pi_{|\hor}$, so that we will write $\hor_{|\pi^{-1}(L)}=\pi^* TL\oplus\pi^* TL^\perp$, where $TL^\perp$ is the normal bundle of $L$ in $M$. Moreover, for any $\eta\in TN_{|\pi^{-1}(L)}$ let us write its decomposition following $TN_{\pi^{-1}(L)}=\pi^*TL \oplus \pi^*TL^\perp\oplus\ver_{\pi^{-1}(L)}$ as 
$$
\eta=\eta_{TL}^\hor + \eta_{TL^\perp}^\hor + \eta^\ver.
$$
Then under the previous identifications, we have
\begin{eqnarray*}
T_{(l,n)} (u^*N) & = & \{(\xi,\eta)\in T_l L\times T_n N|\eta_{TL}^\hor=\xi, \eta_{TL^\perp}^\hor=0\}\\
                 & = & \{(\xi,\xi + \eta^\ver), \xi\in T_lL, \eta^\ver\in\ver_n\}.
\end{eqnarray*} 
This gives us a splitting $T(u^*N)=\ver^{u^*N}\oplus \hor^{u^*N}$ where $\forall (l,n)\in u^*N$,
$$
\ver_{(l,n)}^{u^*N}=\{0\}\times\ver_n \quad \text{and} \quad \hor_{(l,n)}^{u^*N}=T_l L\times \hor_n \cap T_{(l,n)} (u^*N)=\{(\xi,\xi)\in T_lL\times T_lL\}.
$$
Let us determine the orthogonal of the tangent space $T(u^*N)$:
\begin{multline*}
(\alpha,\beta)\in \left( T_{(l,n)} (u^*N)\right)^\perp  \Longleftrightarrow   \\
\begin{array}{rcl}
\forall (\xi,\eta)\in T_{(l,n)} (u^*N), \quad 0 & =  & \langle (\xi,\eta), (\alpha,\beta)\rangle\\
  & = & \langle \xi,\alpha\rangle + \langle\eta ,\beta\rangle \\
 & = &  \langle \xi,\alpha\rangle + \langle \xi,\beta_{TL}^\hor\rangle + \langle 0 ,\beta_{TL^\perp}^\hor\rangle + \langle \eta^\ver,\beta^\ver\rangle\\
& = & \langle \xi,\alpha + \beta_{TL}^\hor\rangle + \langle \eta^\ver,\beta^\ver\rangle
\end{array}
\\
  \Longleftrightarrow    (\alpha + \beta_{TL}^\hor, \beta^\ver)=0.
\end{multline*}
Therefore
$$
\left( T_{(l,n)} (u^*N)\right)^\perp= \{(-\beta_{TL}^\hor,\beta),\beta\in \hor_n\}.
$$
Decomposing each $(a,b)\in T(L\times N)_{|u^*N}$ following the decomposition $T(L\times N)_{|u^*N}=T(u^*N)\oplus T(u^*N)^\perp$: $(a,b)=(\xi,\eta) + (\alpha,\beta)$, then we obtain
$$
\left\{
\begin{array}{l}
a=\eta_{TL}^\hor -\beta_{TL}^\hor\\
b=(\eta_{TL}^\hor +\beta_{TL}^\hor) + \beta_{TL}^\hor + \eta^\ver
\end{array}
\right.
$$
so that this decomposition is therefore given by 
$$
(a,b)=\left( \dfrac{a + b_{TL}^\hor}{2}, a + b_{TL}^\hor + b^\ver\right)  + \left( -\dfrac{(b_{TL}^\hor - a)}{2}, \dfrac{(b_{TL}^\hor - a)}{2} + b_{TL^\perp}^\hor\right). 
$$
Now, let us come back to our fonction $f\colon L\to N$ and  the corresponding section $\tl f\colon (L,b)\to u^*N$. Then let us compute
$$
\begin{array}{rcl}
\onabla{u^*N}{v}{}d^v \tl f =\onabla{u^*N}{v}{}(dl,df)^{\ver^{u^*\!N}} & = & \onabla{u^*N}{v}{}(0,d^v f)
= \left( \left[ \nabla(0,d^v f)\right]_{T(u^*N)} \right)^{\ver^{u^*\!N}}\\ & = & \left( \left[ (0, \nabla d^v f)\right]_{T(u^*N)} \right)^{\ver^{u^*\!N}}\\
 & = &  \left( \dfrac{1}{2}(\nabla d^v f)_{TL}^\hor, \dfrac{1}{2}(\nabla d^v f)_{TL}^\hor + \nabla^v d^v f \right)^{\ver^{u^*\!N}}\\
& = & (0,\nabla^v d^v f)
\end{array}
$$
Finally, we have proved 
\begin{equation}\label{superflat}
\onabla{u^*N}{v}{}d^v \tl f = \nabla^v d^v f
\end{equation}
and by taking the trace, we obtain the first assertion of the theorem.\\
Now, in the same way we obtain 
\begin{equation}\label{totallygeodesic}
\onabla{u^*N}{}{}d \tl f=\left( \dfrac{1}{2}(\nabla d f)_{TL}^\hor, \dfrac{1}{2}(\nabla d f)_{TL}^\hor + \nabla^v d f \right)
\end{equation}
so that $\tl f\colon N\to u^*N$ is harmonic \iif $[\tau(f)]_{TL}^\hor=0$ and $[\tau(f)]^\ver=0$. Therefore $f\colon L\to N$ is harmonic \iif $\tl f\colon N\to u^*N$ is harmonic and $[\tau(f)]_{TL^\perp}^\hor=0$. This completes the proof.\hfill$\square$\\[1mm]
From the proof of theorem~\ref{same-vertical-Harm} (more precisely from (\ref{superflat}) and (\ref{totallygeodesic})), we obtain:
\begin{thm}\label{same2form}
Let us consider the situation described by theorem~\ref{same-vertical-Harm}. Then $f\colon L\to N$ is superflat \iif the corresponding section $\tl f\colon L\to u^*N$ is superflat. Furthermore $f\colon L\to N$ is totally geodesic \iif $\tl f\colon L\to u^*N$ is totally geodesic and $\left[ \nabla df\right] _{TL^\perp}^\hor=0$ (i.e. $\left[ d\pi(\nabla df)\right] _{TL^\perp}=0$).
\end{thm}
\begin{rmk}\em
The metric defined in example~\ref{example1} in $u^*N$ (and thus in theorems~\ref{same-vertical-Harm} and \ref{same2form}, i.e. the metric induced by the product metric, is given by
\begin{equation}\label{product-metric}
|(\xi,\eta)|^2=2|\xi|^2 + |\eta^\ver|^2
\end{equation}
whereas, when $\pi\colon N\to M$ is a Homogeneous fibre bundle, the metric in $u^*N$, considered as an Homogeneous fibre bundle, is defined in \ref{Homfibrbund} by equation (\ref{metrich}) and is given by 
\begin{equation}\label{hom-metric}
|(\xi,\eta)|^2=|\xi|^2 + |\eta^\ver|^2.
\end{equation}
However, theorems~\ref{same-vertical-Harm} and \ref{same2form} hold, of course, also with the metric (\ref{hom-metric}). Indeed, first remark that the theorems hold if we multiply the product metric in $L\times N$ by a constant factor. Then just apply these theorems with the same $(M,g)$ (and thus the same $(L,u^*g)$), $N$ endowed with the new metric $|\cdot|^2_\hor + 2|\cdot|_\ver^2$ (the old one being $|\cdot|^2_\hor + |\cdot|_\ver^2$) and endow $L\times N$ with $\frac{1}{2}$ times the product metric, then the induced metric on $u^*N$ is (\ref{hom-metric}): $\dfrac{1}{2}(
|\xi|^2 +(|\xi|^2 + 2|\eta^\ver|^2))=|\xi|^2 + |\eta^\ver|^2$.
\end{rmk}
\subsection{$G$-invariant metrics}\label{G-invariant-metric}
\subsubsection{About the natural reductivity.}\label{appendix-natred}
\begin{lemma}\label{lemma-nat-red-H-riem-appendix}
Let $N=G/K$ be a naturally reductive Riemannian homogeneous space. Let $H\supset K$ be a subgroup of $G$ such that $M=G/H$ is reductive and $TM$ admits a $G$-invariant lift in $TN$ \wrt the $G$-invariant projection $\pi\colon N\to M$. Then any naturally reductive metric $h$ induces a Riemannian metric on $\tl M=G/H^0$, in particular $G/H^0$ is Riemannian. 
\end{lemma}
\proof 
A $G$-invariant metric $h$ on $N=G/K$ is naturally reductive \iif the torsion $T$ of the canonical connection is a 3-form \wrt $h$. In particular, the identity
\begin{equation}\label{eq-natureductiv-nk}
\langle [Z,X]_\nk, Y \rangle = -\langle X, [Z,Y]_\nk\rangle, \quad \forall X,Y,Z\in\nk
\end{equation}
characterising the natural reductivity holds, where $\g=\kk\oplus\nk$ is an $\Ad K$-invariant decomposition. Moreover according to the hypothesis on  $G/H$, there exists an $\Ad H$-invariant decomposition $\g=\hk\oplus\mk$ such that $\mk\subset\nk$. Then, the subspace $\pk:=\hk\cap \nk$ is $\Ad K$-invariant and we have $\hk=\kk\oplus\pk$ (because $\dim(\hk\cap\nk)= \dim \hk + \dim\nk - \dim\g=\dim\hk -\dim\kk$).   Then, let us apply \eqref{eq-natureductiv-nk} for $X_\mk,Y_\mk\in\mk$ and $Z=V\in\pk$, we obtain: $\langle [V,X_\mk], Y_\mk \rangle = -\langle X_\mk, [V,Y_\mk]\rangle$, i.e. $\adm \hk \subset \so(\mk,h_{|\mk})$ and therefore $\Adm H^0\subset SO(\mk,h_{|\mk})$, or in other words $h_{|\mk}$ is $\Ad H^0$-invariant. \comprf \hsq\medskip\\
In fact, we can do much better by using a result of Kostant \cite{Kost56}.
\begin{defn}\emph{\cite{Kost56}}
Let $\g=\kk\oplus\nk$, $\kk\cap\nk=\{0\}$, $[\kk,\nk]\subset \nk$. We will say that a  inner product B on $\nk$ is stricly invariant if $\ad_\nk\kk\subset \so(\nk,B)$ and it satisfies \eqref{eq-natureductiv-nk}; i.e. $[\ad_\nk X]_\nk$ is skew-symmetric for all $X\in\g$. 
\end{defn}
\begin{thm}\label{kostant56-thm}\emph{\cite{Kost56}}
Let $\g=\kk\oplus\nk$, $\kk\cap\nk=\{0\}$, $[\kk,\nk]\subset \nk$ and $(\g,\kk)$ effective. Let $B$ be a stricly invariant inner product on $\nk$. Let $\g(\nk)=\nk + [\nk,\nk]$, $\kk_1=\g(\nk)\cap\kk$, so that the ideal $\g(\nk)=\kk_1\oplus\nk$. There exists one and only one invariant symmetric bilinear form $B^*$ on $\g(\nk)$  extending $B$ and such that $B^*(\kk_1,\nk)=0$. Moreover, $B^*$ is nonsingular on $\g(\nk)$ and hence on $\kk_1$.
\end{thm}
\begin{defn}\emph{\cite{Kost56}}  We will say, in the situation of the previous theorem that $\nk$ is \emph{pervasive} in $\g$ if $\g(\nk)=\g$. This is also equivalent to $[\nk,\nk]_\kk=\kk$.
\end{defn}
According to this theorem, we see that  if $B$ is stricly invariant on $\nk$, and $\g(\nk)=\g$, then it is automatically $\Ad K$-invariant (and not only $\Ad K^0$-invariant) for any subgroup $K\subset G$ with Lie algebra $G$ (remember that $G$ is always supposed connected, according to our convention). Therefore the fact for $\pk$ to be \emph{natural} (\cite{Kost56}) \wrt $G/K$ is a purely Lie algebra concept. Recall that $\pk$ is natural means that $\pk$ is $\Ad K$-invariant and admits an $\Ad K$-invariant naturally reductive inner product. In other words $G/K$ is Riemannian and $\nabla^{ \frac{1}{2}}$ coincides with the Levi-civita connection of some $G$-invariant metric (or equivalently $G/K$ is Riemannian and the torsion of the canonical connection is totally skew-symmetric \wrt some $G$-invariant metric) and $\pk= T_o G/K$.\\
Moreover, in the situation of lemma~\ref{lemma-nat-red-H-riem-appendix}, we can take $\mk=\hk^\perp$ -\wrt $B^*$ -  which is then $\Ad H$-invariant, as well as $B_{|\mk\times\mk}=\langle\cdot,\cdot\rangle_\mk$  and therefore $G/H$ is Riemannian.
\begin{prop}\label{prop-nat-red-induce-appendix}
Let $N=G/K$ be a naturally reductive Riemannian homogeneous space. Let $H\supset K$ be a subgroup of $G$. Then $M=G/H$ is Riemannian. Moreover, any $G$-invariant naturally reductive metric $h$ on $N$  induces a $G$-invariant metric $g$ on $M$ such that  $\pi\colon (G/K,h)\to (G/H,g)$ is a Riemannian submersion and  therefore a homogeneous fibre bundle. 
\end{prop} 
\proof
First remark that, if $\g=\g(\nk)$, then this follows immediately from theorem~\ref{kostant56-thm}. Moreover, according to \cite[Chap. I]{kowalski},  the  connected normal subgroup $G_1\subset G$ generated by the ideal $\g(\nk)$ is acting transitively on $N=G/K$. Moreover, since of course the projection $\pi\colon G/K\to G/H$ is $G$-invariant: $g.\pi(y)=\pi(g.y)$, $\forall g\in G, y\in N$, we see that $G_1 .\pi(y_0)=\pi(G_1.y_0)=\pi(N)=M$ so that $G_1$ acts transitively on $M=G/H$. This completes the proof. \hsq\medskip\\
Moreover, let us remark the following interesting fact.
\begin{prop}
Let $\g$ be a real Lie algebra and $\tau\colon \g \to \g$ be an automorphism. Let us  assume that $\tau$ defines in $\g$ a $\tau$-invariant reductive decomposition: $\g=\g_0\oplus \nk$ with $\nk=\im(\Id-\tau)$. Let us suppose that there exists a stricly invariant inner product $B$ on $\nk$, then we have $\g(\nk)=\g$.
\end{prop} 
\proof 
According to \cite[Corollary~4]{Kost56}, there exists an ideal $\mak a\subset\g$ complementary to the ideal $\g(\nk)$ : $\g=\g(\nk)\oplus\mak a$. Moreover, the restriction to $\mak a$ of the projections $X\mapsto -X_\nk$ and $X\mapsto X_\hk$ are isomorphisms. Furthermore, $\g(\nk)$ is clearly invariant by $\tau$ and idem for $\mak a$ (this results from the exact definition of $\mak a$ in \cite[Corollary~4]{Kost56}). Therefore, we then have $\mak a=\mak a\cap\kk\oplus \mak a\cap\nk$, which  implies that $\mak a=0$, since the restriction to $\mak a$ of the projections are isomorphisms. \comprf\hsq\bigskip\\
\emph{In the following subsections, we use the notations of section~\ref{finitorderauto}.}
\subsubsection{Existence of an $\Ad H$-invariant inner product on $\mk$ for which $\tm$ is an isometry.}
\begin{lemma}\label{lemma2.7-gener}
Let $G$ be a connected Lie group and $\sigma\colon G\to G$ an automorphism of order $p$. Then $G^\sigma$ has finite many connected components. 
\end{lemma}
\proof This is proved in \cite[Lemma~2.7]{an}, for $p=2$, i.e. $\sigma$ is an involution. The proof in \cite{an} can be generalized without any difficulty to the case of an order $p$ automorphism. \hsq
\begin{prop}\label{AdHsurAdH0finite}
Let $\sigma\colon\g\to\g$ be an automorphism of order $p$, of a real Lie algebra $\g$ with trivial center. Let $G$ be a connected Lie group with Lie algebra $\g$, then $\Ad G^\sigma/\Ad (G^\sigma)^0$ is finite.
\end{prop}
\proof
According to remark~\ref{rmk-G0&G'}, it suffices to apply  lemma~\ref{lemma2.7-gener} to the adjoint group $G'=\Ad G$. \hsq
\begin{thm}\label{AdmH,taum-compact}
Let $\tau\colon\g\to\g$ be an automorphism of order $2k$, of a real Lie algebra $\g$ with trivial center. Let $G$ be a connected Lie group with Lie algebra $\g$ and $H$ a subgroup such that $(G^\sigma)^0\subset H\subset G^\sigma$, where $\sigma=\tau^2$. If $\Adm H$ is compact, then the subgroup generated by $\Adm H$ and $\tm$ is compact.
\end{thm}
\proof Follow the proof of \cite[theorem~21]{ki3}   by using the lemma~\ref{lemma2.7-gener} and prop.~\ref{AdHsurAdH0finite}  above  instead of \cite[Lemma~2.7]{an}. \hsq

\subsubsection{Existence of a naturally reductive metric for which  $\mbox{\und J}$ is an isometry, resp. $F$ is metric.}

\paragraph{The odd case}
Let $\tau\colon \g\to \g$ be an automorphism of order $2k+1$.
We have 
$$
\taum [\adm(X)]_\mk \taum^{-1} = [\adm(\taum X)]_\mk, \quad \forall X\in \mk.
$$
Therefore 
\begin{lemma}\label{Gm-J_0-compact}
Let  $N=G/K$ be a (locally) $(2k+1)$-symmetric space endowed with its canonical almost complex structure $\undj$. Suppose also that $N=G/K$ is naturally reductive. Denoting by  $G(\mk)$ the compact subgroup in $GL(\mk)$ generated by $\Lambda_\mk(\mk) :=\{[\ad_\mk(X)]_\mk, X\in\mk\}\subset\mak{gl}(\mk)$, then $\langle G(\mk), \taum\rangle$, the subgroup generated by $G(\mk)$ and $\taum$,  is compact. Therefore for any metric $h$ on $\mk$ invariant by $\langle G(\mk), \taum\rangle$, then $h$ is naturally reductive, and $\undj_0$ is orthogonal. 
\end{lemma}

\paragraph{The even case}
Let $\tau\colon \g\to \g$ be an automorphism of order $2k$.
We have 
$$
\tau_{\nk} [\ad_\nk(X)]_\nk \tau_{\nk}^{-1} = [\ad_{\nk}(\tau_{\nk} X)]_\nk, \quad \forall X\in \nk.
$$
Therefore 
\begin{lemma}\label{Gn-F-compact}
Let  $N=G/K$ be a (locally) $2k$-symmetric space endowed with its canonical $f$-structure $F$. Suppose also that $N=G/K$ is naturally reductive. Denoting by  $G(\nk)$ the compact subgroup in $GL(\nk)$ generated by $\Lambda_\nk(\nk) :=\{[\ad_\nk(X)]_\nk, X\in\nk\}\subset\mak{gl}(\nk)$, then $\langle G(\nk), \tau_{\nk}\rangle$, the subgroup generated by $G(\nk)$ and $\tau_{\nk}$,  is compact. Therefore for any metric $h$ on $\nk$ invariant by $\langle G(\nk), \tau_{\nk}\rangle$, then $h$ is naturally reductive, and $I_0=\undj_0\oplus -\Id_{\g_k}$ is orthogonal, i.e. $F$ is metric. 
\end{lemma}

\index{m g tau system@$m$-th $(\g,\tau)$-system|see{$m$-elliptic integrable system}}

\printindex

\section*{List of symbols.}
We give here a list of symbols used in the present paper divided by sections. For each symbols, are given its meaning and if necessary the place where it is defined in the paper.  Let us make precise that a list of conventional notations is given at §\ref{generalities}. They are not recalled here.
\subsection*{Section~\ref{1}}
\begin{supertabular}{llr}
$G$ & real Lie group &  \\
$H$ & a closed subgroup of $G$ &  \\
$\g$ &  a real Lie algebra or the Lie algebra of $G$ & \\
$\hk$ & a subalgebra of $\g$ or the Lie algebra of $H$ & \\
$\mk$ & $\Ad H$-invariant complement of $\hk$ in $\g$: $\g=\hk\oplus\mk$ & §\ref{reductivehomspaces}\\ 
$[\mak l]$ & subbundle of $(G/H)\times\g\cong G\times_H \g$ defined by the $\Ad H$-invariant subspace $\mak l\subset \g$ & §\ref{reductivehomspaces}\\
$ [ g , \xi ]$ & element of $G\times_H \mk$ defined by $(g,\xi)\in G\times \mk$ & §\ref{reductivehomspaces}\\
$\theta$ & Maurer-Cartan form of $G$ & §\ref{reductivehomspaces}\\  
$\beta$ & Maurer-Cartan form of $G/H$ & §\ref{reductivehomspaces}\\
$\nabla^0$ & $G$-invariant canonical connection on $M=G/H$ & §\ref{assocovarder}\\ 
$\nabla^t$  &   $\nabla^0 + t [ \, , \, ]_{[\mk]}$ &   §\ref{family}\\ 
$\mU^M$ & natural reductivity term &  §\ref{family}\\
$\overset{\mrm{met}}{\nabla^t}$ &  $\nabla^0 + t\left([\ ,\ ]_{[\mk]} +  \mU^M\right)$ &  §\ref{family}\\
$\nabla^{\mrm{L.C.}}$ & Levi-Civita connection $=\overset{\mrm{met}}{\nabla^\frac{1}{2}}$ & §\ref{family}\\
$\alpha$  &  $U^{-1}.dU$, where $U$ is a local section of $\pi\colon G\to G/H$ &  §\ref{assocovarder} \\
\end{supertabular}

\subsection*{Section~\ref{melliptic}}

\begin{supertabular}{llr}
$\g$ &  a real Lie algebra &  §\ref{def-g-tau}\\
$\tau\colon \g\to\g$ &  automorphism of $\g$ & §\ref{def-g-tau} \\
$G$ & Lie group with Lie algebra $\g$ & §\ref{def-g-tau} \\
$G^\tau$ &  subgroup of $G$ fixed by $\tau$ &  §\ref{def-g-tau} \\
$\g_0=\g^\tau$  & Lie subalgebra fixed by $\tau$ & §\ref{def-g-tau} \\
$G_0$ & closed subgroup of $G$ s.t. $(G^\tau)^0\subset G_0\subset G^\tau$ & §\ref{def-g-tau} \\
$\omega_{k'}$ & $k'$-th primitive root of unity of $\tau\colon\g\to\g$ & §\ref{finitorderauto}.0 \\
$\g_j^\C$ &  $\omega_{k'}^j$-eigenspace of $\tau\colon\g\to\g$, the automorphism of order $k'$ &  §\ref{finitorderauto}.0 \\
$\g_k$ & real $(-1)$-eigenspace of $\tau$, i.e. $ \g_k \otimes \C=\g_k^\C$  &  §\ref{evendetercase}\\
$\mk_j $ & the real subspace in $\g$ s.t. $\mk_j^\C=\g_j^\C\oplus\g_{-j}^\C$, for $j\neq 0,k$ & §\ref{finitorderauto}\\
$\nk$ & the real susbspace in $\g$ s.t. $\nk^\C=\bigoplus_{j\in \Z_k'\setminus \{0\} }\g_j^\C $ &  §\ref{finitorderauto}\\ 
$\sigma$  &  $\tau^2$ & §\ref{evendetercase}\\
$G_0$ &  closed subgroup of $G$ s.t. $(G^\tau)^0\subset G_0\subset G^\tau$ &  §\ref{finitorderauto} \\
$H$ & closed subgroup of $G$ s.t. $(G^\sigma)^0\subset H\subset G^\sigma$   &  §\ref{evendetercase}\\
$N$  &  $G/G_0$  &    §\ref{finitorderauto} \\
$M$ &  $G/H$   &  §\ref{evendetercase}\\
$\undj$ & canonical almost complex structure on a $(2k+1)$-symmetric space &  §\ref{2.1.2}\\
$L$ & Riemann surface & \\
$u=(u_0,\ldots,u_m)$  & \!\!\! \begin{tabular}{l} \!\!\!$(1,0)$-type 1-form on $L$ with values in $\prod_{j=0}^m\g_{-j}^\C$, unknown of   \\
the $m$-th $(\g,\tau)$-system \end{tabular}  &  §\ref{def-subsec}\\
$\alpha_\lm$  &   $\sum_{j=0}^m \lm^{-j}u_j + \lm^j \bar u_j $   &   §\ref{def-subsec}\\
$(S_j)$ &   the component in $\g_{-j}$ of  (Syst($m,\g,\tau$))  &  §\ref{def-subsec}\\
\!\!\begin{tabular}{l}  (Syst($m,\g,\tau$))\\ (Syst($m,\tau$))\\ (Syst($m$))\\ (Syst)\end{tabular}
  &  $m$-th elliptic integrable system defined by $\tau$  &  §\ref{def-subsec}\\
$m_{k'}$  &   $0$ if $k'=1$, and $\displaystyle \left[\dfrac{k'+1}{2}\right]$ if $k'>1$  &  §\ref{def-subsec}\\
$\alpha$ & $\alpha_{\lm=1}$ & §\ref{def-subsec}\\
$\alpha_j$ & $[\alpha]_{\g_j}$  &  §\ref{def-subsec}\\
$f\colon L\to G/G_0$ & geometric solution of $(\syst)$ &   §\ref{geom-sol}\\
$U\colon L\to G$ & lift of $\alpha$ i.e. $U^{-1}.dU = \alpha$ &  §\ref{geom-sol} \\
$\Lm\gtau$  &     $\{\eta_\bullet\colon S^1\to \g|\,\eta_{\omega\lm}=\tau(\eta_\lm),\forall \lm \in S^1\}$ & §\ref{geom-sol}\\
$\Lm\gtc$  &     $\{\eta_\bullet\colon S^1\to \g^\C|\,\eta_{\omega\lm}=\tau(\eta_\lm),\forall \lm \in S^1\}$ & §\ref{geom-sol}\\
$\Lm_m\gtau$ &     $\{\eta_\bullet\in\Lm\gtau|\,\eta_\lm=\sum_{|j|\leq m}\lm^j \hat{\eta}_j\}$  &  §\ref{geom-sol} \\
$\Lm^-\gtc$   &     $\{\eta_\bullet\in\Lm\gtau^\C|\,\eta_\lm=\sum_{j\leq 0} \lm^j \hat{\eta}_j\}$  &  §\ref{geom-sol}\\
$\Lm_*^{\pm}\gtau$  &  $\{\eta_\bullet\in\Lm\gtau^\C|\,\eta_\lm=\sum_{j \gtrless 0} \lm^j \hat{\eta}_j\}$ &  §\ref{geom-sol}\\
$\Lm_*\gtau$   &  $\{\eta_\bullet\in\Lm\gtau| \eta_\lm=\sum_{j \neq 0} \lm^j \hat{\eta}_j\}\cong \Lm\gtau/\g_0$ & §\ref{geom-sol}\\
$\Lm\Gt$  &   $\{U_{\bullet} \colon S^1\to G | U_{\omega\lm}=\tau\left(U_\lm\right)\}$  &  §\ref{geom-sol}\\
$\mathcal S(m)$ &  the space of solutions $\alpha_\lm$ of  (Syst($m$))  &  §\ref{geom-sol}\\
$\nabla^v$ & vertical part of $\nabla$ &  §\ref{determined}\\
$\mrm{Tr}_g$ &  trace \wrt $g$  &   §\ref{determined}\\
$\tau^v(f)$ &  vertical tension  field of $f$  &   §\ref{determined}\\
$F^{[m]}$  &  $G$-invariant $f$-structure on $N=G/G_0$ $2k$-symmetric  &  eq. \eqref{def-even}\\
$F$ &  $F^{[k-1]}$, the canonical $G$-invariant structure on $N=G/G_0$  &  eq. \eqref{def-even}\\
$\und m$  &  $k' -1 - m$  &  prop.~\ref{prop-intermediate-syst}\\
\end{supertabular}

\subsection*{Section~\ref{isometry-twistor}}

\begin{supertabular}{llr}
$\mal{U}_p(E)$  &  $\{A\in SO(E),\, A^p=\Id,\, A^i\neq \Id \text{ if } 1\leq i < p\}$ & §\ref{isometry-twistor}.0\\
$\mal{U}_p^*(E)$  &  $\{A\in\mal{U}_p(E)|1\notin\mrm{Spect}(A)\}$  &  §\ref{isometry-twistor}.0\\
$\mal{U}_p^{**}(E)$   &  $\{A\in\mal{U}_p(E)|\pm 1\notin\mrm{Spect}(A)\}$  &  §\ref{isometry-twistor}.0\\
$\mal{Z}_{2k}(E)$   &  $\mal{U}_{2k}^{**}(E)$ &  §\ref{isometry-twistor}.0\\
$\mal{Z}_{2k+1}(E)$  &  $\mal{U}_{2k+1}^{*}(E)=\mal{U}_{2k+1}^{**}(E)$  &  §\ref{isometry-twistor}.0\\
$E_A(\lm)$   &  $\ker(A-\lm\Id)$  &  §\ref{order 2k}.0 \\
$\omega_{r}$   &   $e^{2i\pi/r}$  &   \\
$\mk_j$    &   the real subspace s.t. $\mk_j^\C= E_A(\omega_{2k}^j)\oplus E_A(\omega_{2k}^{-j})$ for $j\geq 0$  &  §\ref{order 2k}.0\\
$p_j$      &  $\frac{\dim\mk_j}{2}$  &    §\ref{setconnectcomp}\\
$\zdk^\alpha(\rdn)$  &   connected component of $\zdk(\rdn)$  &  §\ref{setconnectcomp}\\
$\zdk^0(\rdn)$ &    $ \left\{ A\in\zdk(\rdn)| A^k = -\Id\right\}=\bigsqcup_{\{\alpha|\forall j, p_{2j}=0\}}\zdk^\alpha$   &   §\ref{setconnectcomp}\\
$\zdk^*(\rdn)$    &   $ \left\{ A\in\zdk(\rdn)| A^k \neq -\Id\right\}=\bigsqcup_{\{\alpha|\exists j, p_{2j}\neq 0\}}\zdk^\alpha$     &    §\ref{setconnectcomp}\\
$r$        &    order of $\Ad J$, for $J\in\zdk(\rdn)$ i.e. $r=\begin{cases} 2k & \text{ in } \zdk^*(\rdn)\\
 k & \text{ in } \zdk^0(\rdn)\end{cases}$ &  §\ref{setconnectcomp}\\
$\displaystyle\mal A_j^\C(J) \dfrac{}{}$  &   $\displaystyle \dfrac{}{} \ker(\Ad J -\omega_r^j\Id)$   &  §\ref{AdJ} \\
$\displaystyle \dfrac{}{}\mal B_j^\C(J)$    &   $\displaystyle \dfrac{}{}\mal A_j^\C(J)\cap(J.\so(2n))^\C$  &   §\ref{AdJ}\\
$\displaystyle \dfrac{}{}\so_j^\C(J)$     &    $\displaystyle \dfrac{}{}\mal A_j^\C(J)\cap\so(2n)^\C $    & §\ref{AdJ}  \\
$\U_{j-1}(J)$   &      $\{ g\in SO(2n) | gJ^j g^{-1}= J^j\}$  &  §\ref{AdJ}  \\
$\mal Z_{2k,j}^\alpha(\rdn,J^j)$   &   $\{J'\in\zdk^\alpha(\rdn)|(J')^j=J^j\}$  &  §\ref{AdJ} \\
$\ul_{j-1}(J)$  &   $\lie (\U_{j-1}(J))=\so_0(J^j)=\left( \oplus_{q=0}^{(r,j)-1}\so_{qp}^\C(J)\right)\cap \so(2n)$  &  §\ref{AdJj}\\
$j\cdot \alpha$  &  image of the component $\alpha\in\pi_0(\maU_p(E))$ by the map $J\mapsto J^j$ & §\ref{effect-power-map}\\
$J_0^\alpha$ &  fixed element of $\zdk(\rdn)$ defined by eq. \eqref{eq-def-J0alpha} &  §\ref{effect-power-map}\\ 
$\mal{SO}(M)$  &   $SO(2n)$-bundle of positively oriented orthonormal frames on $M$  &  §\ref{3twistorspacesandreduction}\\
$J_j$ &    global section of $\left( \zdk^\alpha(M)\right)^j$  &   §\ref{3twistorspacesandreduction}\\
$\Um_{j-1}^\alpha(M)$   & \begin{tabular}{l}  $\U_{j-1}(J_0^\alpha)$-reduction of $\mal{SO}(M)$ defined by $J_j$ i.e. \\  $\{e=(e_1,\ldots,e_{2n})\in \mal{SO}(M)|\,\mal Mat_e(J_j)=(J_0^\alpha)^j\}$ \end{tabular}  &  §\ref{3twistorspacesandreduction}\\
$J_2$ & section of $\left(\zdk(G/H)\right)^2$ defined by $J_0^2=\taum^2$, Lemma~\ref{lemma2}  &  §\ref{subsub-canonical-embedding} \\
$\mak I_{J_0}$  &  canonical embedding $ G/G_0  \hookrightarrow  \mal Z_{2k,2}^{\alpha_0}(G/H, J_2)$, Th.~\ref{embedding}   &  §\ref{subsub-canonical-embedding}\\
$\undj$ &  complex structure defined by $J\in \zdk(E)$, Def.~\ref{defn-j-undj}  &   §\ref{subsub-twistorlift}\\
$\undj$ &  complex structure defined by $J\in \mal Z_{2k+1}(E)$, Def.~\ref{def-undj-odd}  &   §\ref{subsub-twistorlift}\\
$\mk_j(J)$ &  real subspace s.t. $\mk_j(J)^\C = \ker( J -\omega_{2k}^{-j}\Id)\oplus \ker( J -\omega_{2k}^{j}\Id)$ & §\ref{subsub-twistorlift}\\
\end{supertabular}

\subsection*{Section~\ref{vertically-harmonic}}

\begin{supertabular}{llr}
$\mrm{Tr}_g$ &  trace \wrt $g$  &   §\ref{verticalenergyfunctional}\\
$\tau^v(f)$ &  vertical tension  field of  $f$, $\mrm{Tr}_g(\nabla^v d^v f)$  &   §\ref{verticalenergyfunctional}\\
$E^v(u)$  &  $\displaystyle\dfrac{1}{2}\int_M |d^v u|^2 d\mrm{vol}_g$  &  §\ref{verticalenergyfunctional}\\
$\ver$  &  the horizontal subbundle $\ker d\pi$  & §\ref{verticalenergyfunctional}\\
$\hor$  &  the horizontal subbundle $\ver^\perp$ &  §\ref{verticalenergyfunctional}\\
$E$  &  $p\colon (E,\nabla, \langle\cdot,\cdot\rangle)\mapsto (M,g)$ a Riem. vector bundle of rank $2n$  & §\ref{preliminaires-examples}\\
$\so_{\pm}(E_x,J)$  &  $\{ A\in \so(E_x)|\, AJ =\pm JA\}$ &  §\ref{preliminaires-examples}\\
$\so_*(E_x,J)$   &    $\displaystyle\left( \bigoplus_{j\in\Z/r\Z\setminus\{0\}} \negthickspace \negthickspace \so_j^\C(E_x,J)\right) \bigcap\so(E_x)$  &   §\ref{preliminaires-examples}\\
$\ver^\Sigma$, $\hor^\Sigma$ &  vertical and horizontal subbundles on $\Sigma(E)$  &  §\ref{preliminaires-examples} \\
$\ver^\mZ$,  $\hor^\mZ$   &   vertical and horizontal subbundles on $\Sigma(E)$  &  §\ref{preliminaires-examples} \\
$N_\Sigma$, $N_\mZ$ &   resp. $\Sigma(E)$, $\zdk^\alpha(E)$  &  §\ref{preliminaires-examples}\\
$\maB_*(E_x,J)$   &  $\left( \bigoplus_{j\in\Z/r\Z\setminus\{0\}}\negthickspace \negthickspace \mal B_j^\C(E_x,J)\right)  \bigcap\End(E_x)=\ver_J^{\mal Z}$  & §\ref{preliminaires-examples} \\
$\mal I$  & canonical $2k$-structure on $\pi_\mZ^* E$.  & §\ref{preliminaires-examples} \\
$N_\mal{Z}^j$   &   $\mal Z_{2k,j}^\alpha(E,J_j)$  & §\ref{preliminaires-examples} \\
$\ver^{\mZ,j}$,  $\hor^{\mZ,j}$  &  vertical and horizontal subbundles on $\mal Z_{2k,j}^\alpha(E,J_j)$ &  §\ref{preliminaires-examples}\\
$T^\Psi$ &  $\Psi$-torsion  &  §\ref{phitorsion}\\
$T^v$  &  vertical torsion & §\ref{phitorsion} \\
$\mR_\hor$  &  curvature of the subbundle $\hor$ &  §\ref{phitorsion}\\
$Q$  &  principal $H$-bundle over $M$ &  §\ref{geometry}\\
$\omega$ &  $\hk$-valued connection form on $Q$  & §\ref{geometry} \\
$N$  &   Homogeneous fibre bundle $Q/K$ over $M$ with fiber $H/K$  & §\ref{geometry} \\
$\pk$  &  $\Ad K$-invariant summand of $\kk$ in $\hk$: $\hk=\kk\oplus\pk$  & §\ref{geometry} \\
$\pk_Q$ &  $Q\times_K\pk\to N$  &  §\ref{geometry} \\
$I$  &   canonical isomorphism $I\colon \ver\to\pk_Q$  &  §\ref{geometry}\\
$g$  &  Riemannian metric on $M$ & §\ref{geometry} \\
$h$  &  Riemannian metric on  $N$ defined by eq. \eqref{metrich} & §\ref{geometry}\\
$\phi$  &  Homogeneous connection form, i.e. $I\circ\mrm{pr}_\ver$, with $\mrm{pr}_\ver$ projection on $\ver$ & §\ref{geometry}\\
$\Phi$  & Homogenenous curvature form, i.e. $I\circ\mR_\hor$ & §\ref{geometry}\\
$T^c$    &  $\phi$-torsion of $\nabla^c$, i.e. the vertical torsion defined by $\nabla^c$ & §\ref{geometry}\\
$\mU$    &  $\langle \mU(a,b), c\rangle =\langle[c,a]_\pk,b\rangle + \langle a, [c, b]_\pk\rangle$  &  §\ref{geometry} \\
$\mB$    &  $\mU+ [\ ,\ ]_\pk$  &  §\ref{geometry} \\
$\mJ$  &  canonical complex structure on $\pi^*E\to N$  &  §\ref{twistor2}\\
$\mJ$  &  canonical $2k$-structure on $\pi^*E\to N$  &   §\ref{subsect-zdkE}  \\
$\rabla{j}{}$   &  linear metric connection in $E$ defined by $\omega^j:=\omega_{\hk_j|TQ^j}$ &  §\ref{5.3.4} \\
\end{supertabular}

\subsection*{Section~\ref{affineharmonic}}

\begin{supertabular}{llr}
$\tau(f)$  &  tension field $\mrm{Tr}_g(\nabla df)$ &  §\ref{aff-harm-gener-prop}\\
$(L,j)$  &   Riemann surface  &  \\
$T$  &  Torsion of a connection $\nabla$  &  §\ref{aff-harm-gener-prop}\\
$T_g(f)$  &   $*_g (f^*T)$   &   §\ref{holharmmap}\\
$\widehat\nabla$, $\hat d$  &  \begin{tabular}{l} $\C$-linear extension of $\nabla$ (resp. $d$) from $TM^\C$ to $(E,J)$,\\ where $\nabla$ is a connection on $(E,J)$  \end{tabular}   &  §\ref{holharmmap}\\
$\hat\partial$  &  $\hat d =\hat\partial + \overline{\hat\partial}$ &   §\ref{holharmmap}\\
$ U(X,Y)$  &  $\langle U(X,Y),Z\rangle =T(Z,X,Y) + T(Z,Y,X)$  &   §\ref{Total-Skew-Sym-tors}\\
$\nabla^h$   &   Levi-Civita connection of $(N,h)$  &  §\ref{Total-Skew-Sym-tors}\\
$\mrm{Bil}(E)$   &   $E^*\otimes E^*\otimes E$  &   §\ref{WZ}\\
$\mal T (E)$   &   $(\Lm^2 E^*)\otimes E \subset \mrm{Bil}(E)$  & §\ref{WZ} \\
$\mrm{Bil}^{\eps,\eps'}(E,J)$ &  $\{A\in\mrm{Bil}(E)| A(J\cdot,\cdot) = \eps JA,\  A(\cdot,J\cdot)=\eps'JA \}$ & §\ref{WZ}\\
$\mal T$, $\mrm{Bil}$  &   $\mal T(TN)$, $\mrm{Bil}(TN)$   &    §\ref{WZ} \\
$T^{\eps,\eps'}$  &   component of $T$ in  $\mrm{Bil}^{\eps,\eps'}(E,J)$  &  §\ref{WZ}\\
$N_J$   &   Nijenhuis tensor  of  $J$   &   §\ref{WZ}\\
$\Omega_J$  &  K\"{a}hler form of $J$   &    §\ref{good-geom-set-kahler}\\
$\mrm{Skew}(B)(X,Y,Z)$   &  $ B(X,Y,Z) + B(Y,Z,X) + B(Z,X,Y)$  &  §\ref{TN-valued-2-forms}\\
$B^c$ &   $-B(J\cdot,J\cdot,J\cdot)$  &   §\ref{TN-valued-2-forms}\\
$d^c\beta$  &  $(d\beta)^c$    &        §\ref{TN-valued-2-forms}\\
$J\cdot B$  &   $-JB(J\cdot,J\cdot)$\,  i.e. $B(J\cdot,J\cdot,J\cdot)=-B^c$ in terms of trilinear forms &   §\ref{TN-valued-2-forms}\\
$B_g(f)$   &   $* f^*B=B(f_*TL)$   &  §\ref{TN-valued-2-forms} \\
$J {\circact} B$   &  $B(J\cdot,\cdot,\cdot) + B(\cdot,J\cdot,\cdot) + B(\cdot,\cdot,J\cdot)$   &    §\ref{TN-valued-2-forms}  \\
$J\star B$    &   $\dfrac{1}{2}\left( B(J\cdot, J\cdot,J\cdot) + B(J\cdot,\cdot,\cdot) + B(\cdot,J\cdot,\cdot) + B(\cdot,\cdot,J\cdot)\right) = \dfrac{1}{2}\left( J\cdot B + J\circact B \right) $  &  §\ref{TN-valued-2-forms}\\ 
$\mal G_1$  &  class $\mal W_1\oplus \mal W_2 \oplus \mal W_3$ in the Gray-Hervella classification &  §\ref{AlmostHermG_1}\\
$H$  &   3-form $J\cdot T$   &   §\ref{AlmostHermG_1}\\
$H^\star$  & 3-form $J\star T$    &   §\ref{AlmostHermG_1}\\
\end{supertabular}

\subsection*{Section~\ref{gene-Harm-f-structure}}

\begin{supertabular}{llr}
$F$  &  $f$-structure  &  §\ref{subsec-f-struct} \\
$\hor$,  $\ver$  &  $\hor=\mrm{Im} F$,  $\ver=\ker F$   &   §\ref{subsec-f-struct}\\
$P$  &   $-F^2$  &   §\ref{subsec-f-struct}\\
$q$  &    $\Id - P$ &   §\ref{subsec-f-struct}\\
$\Bar J$   &   $F_{|\hor}$   &  §\ref{subsec-f-struct}\\
$N_F$   &   Nijenhuis tensor of $F$  &  §\ref{subsec-f-struct}\\
$\Bar B$   &  $ B_{|\hor^2}^\hor$, for any $B\in \mal T$ &  §\ref{subsec-f-struct}\\
$F\cdot B$    &   $B(F\cdot,F\cdot,F\cdot):=-B^c$   &  §\ref{subsec-f-struct} \\
$F\circact B$    &   $B(F\cdot,\cdot) + B(\cdot,F\cdot) -FB(\cdot,\cdot)$  &  §\ref{subsec-f-struct}\\
$F\bullet B$   &  $F\cdot B + \dfrac{1}{2}F\circact (B- \Bar B)$   &  §\ref{subsec-f-struct}\\
$F\star B$   &   $\dfrac{1}{2}\left( F\cdot B + F\circact B\right)$   &   §\ref{subsec-f-struct}\\
$\Phi$  &  $\mR_\hor$ the curvature of $\hor$  &   §\ref{introd-linear-connect} \\
$N_{\Bar J}$  &   $N_{F|\hor\times\hor}^\hor$  &    §\ref{Def-not-Firstprop}\\
$\Omega_A(X,Y)$  &   $\langle A(X),Y\rangle$   &  §\ref{Def-not-Firstprop}\\
$\mrm{Sym}(B)(X,Y)$  &   $B(X,Y) + B(Y,X),\quad \forall X,Y\in TN$  &   §\ref{Def-not-Firstprop}\\
$\mal S (E_1\times E_2\times E_3)$   &  $ \underset{i,j,k}{\mal S} E_i\otimes E_j\otimes E_k$   &   §\ref{Def-not-Firstprop} \\
$D$   &   $\nabla^h$ the Levi-Civita connection of $h$  &  §\ref{Def-not-Firstprop}\\
$ U(X,Y)$  &  $\langle U(X,Y),Z\rangle =T(Z,X,Y) + T(Z,Y,X)$  & §\ref{Total-Skew-Sym-tors} \\
$\tl N_F$  &  extended Nijenhuis tensor, Def.~\ref{def-Nijenhuis-extended}  &   §\ref{Character-metic-f-connect}\\
$\rho(V)$    &   $\langle \rho(V)H_1,H_2\rangle =\langle\Phi(H_1,H_2),V\rangle$,   $H_1,H_2\in\hor, V\in\ver$, & §\ref{Prechara-parachara-connect}\ \\
$\sigma(H)$  &  $\langle \sigma(H)V_1,V_2\rangle =\langle \mR_\ver(V_1,V_2),H\rangle$,   $V_1,V_2\in\ver, H\in\hor$  &  §\ref{Prechara-parachara-connect}\\
$\Bar R(X,Y)$   &  $\rho(\Phi(X,Y))$  &  §\ref{Prechara-parachara-connect}\\
$\Bar{\mrm  A}(X,Y)$    &   $\Bar R(\Bar J X,Y) + \Bar R(X,\Bar J Y) - [\Bar J, \Bar R(X,Y)]$   &  §\ref{Prechara-parachara-connect}\\
$\Phi^{(\eps)}$  &  $\Phi^{(\eps)}(\Bar J\cdot,\Bar J\cdot)=\eps \Phi^{(\eps)}$  &  §\ref{Prechara-parachara-connect} \\
$\Bar R^{(\eps)}$   &   $\rho^\eps(\Phi^{(\eps)}) $   &  §\ref{Prechara-parachara-connect} \\
$\Bar{\mrm  A}^{(\eps)}(X,Y)$   &  $\Bar R^{(\eps)}(\Bar J X,Y) + \Bar R^{(\eps)}(X,\Bar J Y) - [\Bar J, \Bar R^{(\eps)}(X,Y)]$ &   §\ref{Prechara-parachara-connect} \\
$\mrm{Con}(\ver)$   &    $\left\lbrace D^v - \dfrac{1}{2}\mR_\ver(Y^v,Z^v,X^h) + \mal C(\Lambda^3\ver^*)\right\rbrace$, Def.~\ref{def-con-ver}   &   §\ref{subsec-Reduct-f-manif}\\
$D^g$  &   Levi-Civita connection of  $(M,g)$   &  §\ref{Riem-sub-Global-G_1}\\
$\widetilde{D^g}$  & connection in $\hor$ defined by the lift of $D^g$  &  §\ref{Riem-sub-Global-G_1}\\
$\mal F$  &  canonical metric $f$-structure on $\zdk^\alpha(M)$  & §\ref{example-zdk} \\
$\mal F_j$  &  canonical metric $f$-structure on $\mal Z_{2k,j}^\alpha(M)$  & §\ref{example-zdkj} \\ 
$H$  &  $F\cdot T$  &  §\ref{dH=0} \\
$H^\star$   &   $  F\star T$   &  §\ref{dH=0}  \\
$\mal S(\hor,\ver)$  &   $\mal S(\hor\times\hor\times\ver) \oplus \mal S(\hor\times\ver\times\ver)$   &  §\ref{dH=0}  \\
$\mathring B$     &   $B_{|\mal S(\hor,\ver)}$  &  §\ref{dH=0} \\ 
$\mrm{Supp}(C)$  &    $(\ker C)^\perp$, $\forall C\in \mal C (\Lambda^2 \hor^*\otimes TN)$  &  §\ref{dH=0}  \\  
\end{supertabular}

\subsection*{Section~\ref{affineharmhom}}

\begin{supertabular}{llr}
$\tau^{t,v}$ & vertical tension field \wrt $\nabla^t$ &  §\ref{Affvertharmredhom}\\
$T^{t,v}$ &  vertical torsion \wrt $\nabla^t$  &   §\ref{Affvertharmredhom}\\
$\mU^{\mak l}$  &  $\langle \mU^{\mak l}(X,Y),Z\rangle = \langle [Z,X]_{\mak l},Y\rangle + \langle X,[Z,Y]_{\mak l}\rangle\quad \forall X,Y\in \mak l$,  &  §\ref{Affvertharmredhom}\\
$\mk'$  &   $\Ad K$-invariant complement of $\pk$ in $\mk$ &  §\ref{subsec-odd-vert-hol-harm-Pfaff}\\
$\psi$   &   $\psi\colon TN \to [\mk']$ projection on $[\mk']$ along $[\pk]$  &  §\ref{subsec-odd-vert-hol-harm-Pfaff}\\
$\mk_*$   &   $\left\{X\in\mk'|\, [X,\pk]_\pk=\{0\}\right\}$   &  §\ref{subsec-odd-vert-hol-harm-Pfaff}\\
$\mk_1$ &  $\undj$-invariant $\Ad K$-invariant complement of $\mk_*$ in $\mk'$  &  §\ref{subsec-odd-vert-hol-harm-Pfaff}\\
$\mak l^\pm$ &   $\pm i$-eingenspace of $\undj_{|\mak l}$  &  §\ref{subsec-odd-vert-hol-harm-Pfaff}\\
$\psi_1$ &  $[\mk_1]$-component of $\psi$ &  §\ref{subsec-odd-vert-hol-harm-Pfaff}\\
$\und m $   &    $k' - 1 - m$  &   §\ref{subsec-hom-intermediate}\\
$\hor^m$  &   $\oplus_{j=1}^{\und m} [\mk_j]$  &  §\ref{subsec-hom-intermediate}\\
$\ver^{ m}$  &   $\oplus_{j=\und m + 1}^{[ \frac{k'}{2}] } [\mk_j]$ &  §\ref{subsec-hom-intermediate} \\
\end{supertabular}

\end{document}